\documentclass{article}

\usepackage{longtable}
\usepackage{booktabs}
\def\tightlist{}

\usepackage{charter}
\usepackage{fancyhdr}
\usepackage[margin=1.6in]{geometry}

\usepackage{amsmath,amssymb}
\usepackage[colorlinks=true]{hyperref}
\usepackage{xcolor}
\definecolor{myurlcolor}{rgb}{0.6,0,0}
\definecolor{mycitecolor}{rgb}{0,0,0.8}
\definecolor{myrefcolor}{rgb}{0,0,0.8}
\hypersetup{linkcolor=myrefcolor,citecolor=mycitecolor,urlcolor=myurlcolor}

\renewcommand{\texttt}[1]{%
  \begingroup
  \ttfamily
  \begingroup\lccode`~=`/\lowercase{\endgroup\def~}{/\discretionary{}{}{}}%
  \begingroup\lccode`~=`[\lowercase{\endgroup\def~}{[\discretionary{}{}{}}%
  \begingroup\lccode`~=`.\lowercase{\endgroup\def~}{.\discretionary{}{}{}}%
  \catcode`/=\active\catcode`[=\active\catcode`.=\active
  \scantokens{#1\noexpand}%
  \endgroup
}

\usepackage{titlesec}

\renewcommand{\thesection}{Week~\arabic{section}}
\titleformat{\section}[display]{\normalfont}{\Large\bfseries\thesection}{1em}{\large\normalfont}

\usepackage{titletoc}
\titlecontents{section}[0em]{\normalfont}{\bfseries\thecontentslabel\hspace{1em}\normalfont\small}{}{\titlerule*[0.3pc]{.}\small\itshape\thecontentspage}[\vspace{0.5em}]

\usepackage[toc]{multitoc}

\title{This Week's Finds in Mathematical Physics (101--150)}
\author{John Baez}
\date{April 9, 1997 to June 18, 2000}

\usepackage{tikz}
\usetikzlibrary{knots}
\usetikzlibrary{arrows}
\usetikzlibrary{decorations.pathmorphing}

\usepackage{etoolbox}
\AtBeginEnvironment{quote}{\itshape}

\usepackage{tikz-cd}

\usepackage{graphicx}
\usepackage[export]{adjustbox}

\usepackage[main]{embedall}

\begin{document}

\begin{titlepage}
  \begin{center}
    {\Huge\textbf{This Week's Finds in}}
  \\[0.7em]{\Huge\textbf{Mathematical Physics}}
  \\[1em]{\huge\textit{Weeks 101 to 150}}
  \\[4em]{\LARGE \textit{April 9, 1997} to \textit{June 18, 2000}}
  \\[4em]{\huge by John Baez}
  \\[0.5em]{\Large{Typeset by Tim Hosgood}}

\vskip 23em

 \begin{tikzpicture}
\scalebox{2.4}{
    \node (mr) at (2,2) {$(a (b c)) d$};
    \node (t) at (0,3.5) {$((a b) c) d$};
    \node (ml) at (-2,2) {$(a b) (c d)$};
    \node (bl) at (-1.3,0) {$a (b (c d))$};
    \node (br) at (1.3,0) {$a ((b c) d)$};
    \draw[double,double equal sign distance,-implies] (t) to (ml);
    \draw[double,double equal sign distance,-implies] (ml) to (bl);
    \draw[double,double equal sign distance,-implies] (t) to (mr);
    \draw[double,double equal sign distance,-implies] (mr) to (br);
    \draw[double,double equal sign distance,-implies] (br) to (bl);
}
  \end{tikzpicture}
  \end{center}

\end{titlepage}

These are issues 51 to 100 of \emph{This Week's Finds of Mathematical Physics}.
This series has sometimes been called the world's first blog, though it was originally 
posted on a ``usenet newsgroup'' called \texttt{sci.physics.research} --- a form of communication that predated the world-wide web.  I began writing this series as a way to talk about papers in mathematics and physics, and continued doing this in issues 101--150, focusing strongly on topics connected to quantum gravity, topoogical quantum field theory, and \(n\)-categories.  However, I digressed into topics ranging from biology to the fiction of Greg Egan to the game of Go.  I also explained some topics in homotopy theory in this series of mini-articles:

\begin{itemize}
\item \hyperlink{homotopy_A}{\textbf{A.}} Presheaf categories.
\item \hyperlink{homotopy_B}{\textbf{B.}} The category of simplices, $\Delta$.
\item \hyperlink{homotopy_C}{\textbf{C.}} Simplicial sets.
\item \hyperlink{homotopy_D}{\textbf{D}.} Simplicial objects.
\item \hyperlink{homotopy_E}{\textbf{E}.} Geometric realization.
\item \hyperlink{homotopy_F}{\textbf{F}.} Singular simplicial set.
\item \hyperlink{homotopy_G}{\textbf{G}.} Chain complexes.
\item \hyperlink{homotopy_H}{\textbf{H}.} The chain complex of a simplicial abelian
group.
\item \hyperlink{homotopy_I}{\textbf{I}.} Singular homology.
\item \hyperlink{homotopy_J}{\textbf{J}.} The nerve of a category.
\item \hyperlink{homotopy_K}{\textbf{K}.} The classifying space of a category.
\item \hyperlink{homotopy_L}{\textbf{L}.} $\Delta$ as the free monoidal category on 
a monoid object.
\item \hyperlink{homotopy_M}{\textbf{M}.} Simplicial objects from adjunctions.
\item \hyperlink{homotopy_N}{\textbf{N}.} The loop space of a topological space.
\item \hyperlink{homotopy_O}{\textbf{O}.} The group completion of a topological monoid.
\end{itemize}

Tim Hosgood kindly typeset all 300 issues of \emph{This Week's Finds} in 2020. They will be released in six installments of 50 issues each.  I have edited the issues here to make the style a bit more uniform and also to change some references to preprints, technical reports, etc.\ into more useful arXiv links.  This accounts for some anachronisms where references in an issue appeared only after that issue was written.

There are undoubtedly many typos and other mistakes.
If you find any, please contact me and I will try to fix them.

\tableofcontents

\hypertarget{week101}{%
\section{April 9, 1997}\label{week101}}

Darwinian evolution through natural selection is an incredibly powerful
way to explain the emergence of complex organized structures. However,
it is not the \emph{only} important process that naturally gives rise to
complex structures. Maybe when we study biology we should also look for
other ways that order can spontaneously arise.

After all, there are plenty of complex structures in the nonbiological
world. When it snows, we see lots of beautiful snowflakes with similar
but different hexagonal structures. Do we conclude that snowflakes
\emph{evolved} to be hexagonal through natural selection? No.

But wait! Maybe in some sense a hexagonal snowflake is ``more fit'' in
certain weather conditions. Perhaps this shape is more efficient at
getting water molecules to adhere to it than other shapes. We can think
of different snowflakes as engaged in ``competition'' for water
molecules, and the ones that grow fastest as the ``winners''. In fact,
the exact shapes of snowflakes in a snowstorm depend crucially on the
temperature, humidity and so on\ldots{} so who the ``winners'' are
depends on the environment, just as in Darwinian evolution!

A biologist will reply: fine, but this is still not ``Darwinian
evolution''. For Darwinian evolution in the strict sense, we require
that there be a ``lineage''. Darwinian evolution applies only to
entities that reproduce and pass some of their traits down to
descendants. The idea is that over the course of many generations,
traits that aid reproduction will accumulate, while traits that hinder
it will be weeded out. Snowflakes don't have kids. A one-shot
competition for resources, followed by melting into oblivion the next
day, is not what Darwinian evolution is about.

Okay, okay, so it's not Darwinian evolution. But it's still interesting.
It's showing us that Darwinian evolution is just \emph{one} of various
ways that order can arise. So we shouldn't study Darwinian evolution in
isolation. We should study \emph{all} the ways that systems
spontaneously generated complex patterns, and see how they relate. If we
do that, perhaps we'll see a bunch of interesting relationships between
physics and chemistry and biology. Also, maybe we'll get a better handle
on how life arose in the first place\ldots{} that curious transition
from chemistry to biology.

If I wasn't so hooked on quantum gravity I would love to work on this
stuff. It's obviously cool, and obviously a lot more \emph{practical}
than quantum gravity. The origin of complexity a very hot topic these
days. But alas, I am just an old-fashioned guy in love with simplicity.
Whenever I see a new journal come out with a title like ``Complex
Systems'' or ``Journal of Complexity'' or ``Santa Fe Institute Studies
in the Science of Complexity'', I heave a wistful sigh and dream of
starting a journal entitled ``Simplicity''.

Actually, the fun lies in the interplay between complexity and
simplicity: how complex phenomena can arise from simple laws, and
sometimes obey new simple laws of their own. I like to hang out on the
simple end of things, but that doesn't stop me from enjoying the new
work on complexity. At one point I got a big kick out of Manfred Eigen's
work on ``hypercycles'' --- systems of chemicals that catalyze each
others formation. (You may remember Eigen as the discoverer of the
``Eigenvalue''\ldots{} in which case I pity you.) Presumably life
started as some sort of hypercycle, so the mathematical study of the
competition between hypercycles may shed some light on why there is only
one genetic code. There is a lot of nice math of this type in:

\begin{enumerate}
\def\labelenumi{\arabic{enumi})}
\tightlist
\item
  Manfred Eigen, \emph{The Hypercycle, a Principle of Natural
  Self-Organization}, Springer, Berlin, 1979.
\end{enumerate}
\noindent
Another name that comes up in this context is Ilya Prigogine, mainly for
his work on non-equilibrium thermodynamics and the spontaneous formation
of patterns in dissipative systems. The following are just a few of his
many books:

\begin{enumerate}
\def\labelenumi{\arabic{enumi})}
\setcounter{enumi}{1}
\item
  G. Nicolis and I. Prigogine, \emph{Self-Organization in Nonequilibrium
  Systems: from Dissipative Structures to Order Through Fluctuations},
  Wiley, New York, 1977.

  Ilya Prigogine, \emph{From Being to Becoming: Time and Complexity in
  the Physical Sciences}, W. H. Freeman, San Francisco, 1980.

  Ilya Prigogine, \emph{Introduction to Thermodynamics of Irreversible
  Processes}, Interscience Publishers, New York, 1967.
\end{enumerate}

A bit more recently, the work of Stuart Kauffman has dominated the
subject. It's really he who has pushed for the unified study of the
whole gamut of methods of spontaneous generation of order, particularly
in the context of biological systems. He's written two books. The
latter, in particular, includes a lot of math problems just
\emph{waiting} to be tackled by good mathematicians and physicists.

\begin{enumerate}
\def\labelenumi{\arabic{enumi})}
\setcounter{enumi}{2}
\item
  Stuart A. Kauffman, \emph{At Home in the Universe: the Search for Laws
  of Self-Organization and Complexity}, Oxford U.\ Press, Oxford, 1995.

  Stuart A. Kauffman, \emph{The Origins of Order: Self-Organization and
  Selection in Evolution}, Oxford U.\ Press, Oxford 1993.
\end{enumerate}

If non-Darwinian forms of spontaneous pattern-formation can be important
in biology, can Darwinian evolution be important in non-biological
contexts? Well, as I mentioned in \protect\hyperlink{week31}{``Week
31''} and \protect\hyperlink{week33}{``Week 33''}, the physicist Lee
Smolin has an interesting hypothesis about how the laws of nature may
have evolved to their present point by natural selection. The idea is
that black holes beget new ``baby universes'' with laws similar but not
necessarily quite the same as their ancestors. Now this is extremely
speculative, but it has the saving virtue of making a lot of testable
predictions: it predicts that all the constants of nature are tuned so
as to maximize black hole production. Smolin has just come out with a
book on this, which also happens to be a good place to learn about his
work on quantum gravity:

\begin{enumerate}
\def\labelenumi{\arabic{enumi})}
\setcounter{enumi}{3}
\tightlist
\item
  Lee Smolin, \emph{The Life of the Cosmos}, Crown Press, 1997.
\end{enumerate}
\noindent
Interestingly, Stuart Kauffman and Lee Smolin have teamed up to write a
paper on the problem of time in quantum gravity:

\begin{enumerate}
\def\labelenumi{\arabic{enumi})}
\setcounter{enumi}{4}
\tightlist
\item
  Stuart Kauffman and Lee Smolin, ``A possible solution to the problem
  of time in quantum cosmology'', available as
  \href{https://arxiv.org/abs/gr-qc/9703026}{\texttt{gr-qc/9703026}}.
\end{enumerate}

Right now you can also read this paper on John Brockman's website called
``Edge''. This website features all sorts of fun interviews and
discussions. For example, if you look now you'll find an intelligent
interview with my favorite living musician, Brian Eno. More to the
point, a discussion of Kauffman and Smolin's paper is happening there
now. As a long-time fan of USENET newsgroups and other electronic forms
of chitchat, I'm really pleased to see how Brockman has set up a kind of
modern-day version of the French salon.

\begin{enumerate}
\def\labelenumi{\arabic{enumi})}
\setcounter{enumi}{5}
\tightlist
\item
  Edge: \url{http://www.edge.org}
\end{enumerate}

Okay. Now\ldots{} what's even more fashionable, trendy, and close to the
cutting edge than complexity theory? You guessed it: homotopy theory!
Currently known only to hippest of the hip, this is bound to hit the
bigtime as soon as they figure out how to make flashy color graphics
illustrating the Adams spectral sequence.

Last week I went to the Workshop on Higher Category Theory and Physics
at Northwestern University, and also, before that, part of a conference
on homotopy theory they had there. Actually these two subjects are
closely related: homotopy theory is a highly algebraic way of studying
the topology of spaces of various dimensions, and lots of what we
understand about ``higher dimensional algebra'' comes from homotopy
theory. So it was a nice combination.

Lots of the homotopy theory was over my head, alas, but what I
understood I enjoyed. It may seem sort of odd, but the main thing I got
out of the homotopy theory conference was an explanation of why the
number 24 is so important in string theory! In bosonic string theory
spacetime needs to be 26-dimensional, but subtracting 2 dimensions
for the surface of the string itself we get 24, and it turns out that
it's really the special properties of the number 24 that make all the
magic happen.

I began to delve into these mysteries in
\protect\hyperlink{week95}{``Week 95''}. There, however, I was mainly
reporting on very fancy stuff that I barely understand, stuff that seems
like a pile of complicated coincidences. Now, I am glad to report, I am
beginning to understand the real essence of this 24 business. It turns
out that the significance of the number 24 is woven very deeply into the
basic fabric of mathematics. To put it rather mysteriously, it turns out
that every integer has some subtle ``hidden symmetries''. These
symmetries have symmetries of their own, and in turn \textbf{these}
symmetries have symmetries of \textbf{their} own --- of which there are
exactly 24.

Hmm, mysterious. Let me put it another way. It probably won't be obvious
why this is another way of saying the same thing, but it has the
advantage of being more concrete. Suppose that the integer \(n\) is
sufficiently large --- 4 or more will do. Then there are 24 essentially
different ways to wrap an \((n+3)\)-dimensional sphere around an
\(n\)-dimensional sphere. More precisely still, given two continuous
functions from an \((n+3)\)-sphere to an \(n\)-sphere, let's say that
they lie in the same ``homotopy class'' if you can continuously deform
one into another. Then when \(n\) is \(4\) or more, it turns out that
there are exactly 24 such homotopy classes.

Now that I have all the ordinary mortals confused and all the homotopy
theorists snickering at me for making such a big deal out of something
everyone knows, I should probably go back and explain what the heck I'm
getting at, and why it has to do with string theory. But I'm getting
worn out, and your attention is probably flagging, so I'll do this next
time. I'll say a bit about homotopy theory, stable homotopy theory, the
sphere spectrum, and why Andre Joyal says we should call the sphere
spectrum the ``integers'' (thus explaining my mysterious remark above).

\begin{center}\rule{0.5\linewidth}{0.5pt}\end{center}

\begin{quote}
\emph{Deep, deep infinity! Quietness. To dream away from the tensions of
daily living; to sail over a calm sea at the prow of a ship, toward a
horizon that always recedes; to stare at the passing waves and listen to
their monotonous soft murmur; to dream away into unconsciousness\ldots.}

--- Maurits Escher
\end{quote}

\hypertarget{week102}{%
\section{April 21, 1997}\label{week102}}

In \protect\hyperlink{week101}{``Week 101''} I claimed to have figured
out the real reason for the importance of the number 24 in string
theory. Now I'm not so sure --- some pieces of the puzzle that I thought
would fit together don't seem to be fitting. Maybe if I explain what I
know so far, some experts will hand me some of the missing pieces, or
tell me where the ones I have go.

Most of the puzzle pieces came from a talk at a conference on homotopy
theory that I went to:

\begin{enumerate}
\def\labelenumi{\arabic{enumi})}
\tightlist
\item
  Ulrike Tillmann, ``The moduli space of Riemann surfaces --- a homotopy
  theory approach'', talk at Northwestern University Algebraic Topology
  Conference, March 27, 1997
\end{enumerate}
\noindent
However, some conversations with Andre Joyal during this conference
really helped turn my attention towards what might be going on here.

Let's start by recalling some stuff about homotopy groups of spheres.
There are often lots of topologically different ways of wrapping an
\(m\)-dimensional sphere around a \(k\)-dimensional sphere. For example,
if \(m = k = 1\), we're talking about the ways of wrapping a circle
around a circle. These are classified by an integer called the ``winding
number''. We can make this concrete by thinking of the circle as the
unit circle in the complex plane. Take your favorite integer and call it
\(n\). Then the function \[f(z) = z^n\] maps the unit circle (the
complex numbers with \(|z| = 1\)) to itself. If \(n\) is positive, this
function wraps the unit circle around itself \(n\) times in the
counterclockwise direction. If \(n\) is negative, the circle gets
wrapped around in the other direction. If \(n\) is zero, \(f(z) = 1\),
so we have a constant function --- no ``wrapping around'' at all!

It turns out that any continuous function from the circle to itself can
be continuously deformed to exactly one of these functions
\(f(z) = z^n\). Homotopy theory is all about such continuous
deformations. In the jargon of homotopy theory, we say two functions
from some space to some other space are ``homotopic'' if we can
continuously deform the first function to the second. Another way of
putting it is that the two functions lie in the same ``homotopy class''.
Speaking of jargon, real topologists never say ``continuous function'':
instead, they say ``map''. So, using this jargon: we know the homotopy
class of a map from the circle to itself if we know its winding number.

Now: what happens if we go to higher dimensions? What are all the
homotopy classes of maps from the \(m\)-dimensional sphere to the
\(k\)-dimensional sphere? Spheres are pretty simple spaces, so one might
at first guess there is some simple answer to this question for all
\(m\) and \(k\).

Unfortunately, it's far from simple. In fact, nobody knows the answer
for all \(m\) and \(k\)! People \emph{do} know the answer for zillions
of particular values of \(m\) and \(k\). But there is no simple pattern
to it: instead, there is an incredibly complicated and beautiful weave
of subtle patterns, which we have not gotten to the bottom of\ldots{}
and perhaps never will.

To get a little feel for this, let's bring in some standard notation:
folks use \(\pi_m(X)\) to denote the set of homotopy classes of maps
from an \(m\)-dimensional sphere to the space \(X\). When \(m > 0\),
this set is actually a group, called the ``\(m\)th homotopy group'' of
\(X\). These groups are of major importance in algebraic topology.

So, what we are talking about is \(\pi_m(S^k)\): the set of all homotopy
classes of ways of wrapping an \(m\)-sphere around an \(k\)-sphere. I
already implicitly said that \[\pi_1(S^1) = \mathbb{Z}\] where
\(\mathbb{Z}\) stands for the integers, since the winding number is an
integer. The same thing happens if we go up a dimension:
\[\pi_2(S^2) = \mathbb{Z}\] In other words: you can wrap a 2-sphere (an
ordinary sphere) \(n\) times around itself for any integer \(n\). How?
Well, say we use spherical coordinates and describe a point on the
sphere using its angle \(\varphi\) from the north pole, together with
the angle \(\theta\) saying how far east it is from Greenwich. Then the
map \[f(\varphi,\theta) = (\varphi, n\theta)\] does the job. Any map
from \(S^2\) to itself is homotopic to exactly one of these.

The same basic idea works in any higher dimension, too:
\[\pi_k(S^k) = \mathbb{Z} \quad\text{for any}\quad k\geqslant 1\]

In other words, there is always an integer \(n\) that plays the role of
the ``winding number'' of a map from the \(k\)-sphere to itself ---
though only uncouth physicists call it the ``winding number'';
mathematicians call it the ``degree''.

So far, so good. Now, what about mapping a sphere to another sphere of
\emph{higher} dimension? This is nice and simple:
\[\pi_m(S^k) =\{0\} \quad\text{whenever}\quad m < k\]

The \(\{0\}\) there is just a standard way to write a set with only one
element, which we call ``zero''. So what we mean is that there's only
\emph{one} homotopy class of ways to map a sphere to a sphere of higher
dimension. There is always enough ``room'' to wiggle around one map
until it looks like another.

What about mapping a sphere to another sphere of \emph{lower} dimension?
Here is where the trouble starts! --- or the fun, depending on your
attitude towards complexity. For example, there is only one homotopy
class of maps from a 2-sphere to a circle: \[\pi_2(S^1) = \{0\}\] There
is just no way a 2-sphere can get interestingly ``stuck'' on the
``hole'' of the circle. This may seem obvious. But it's not really quite
as obvious at it seems, because if we move up one dimension, we have:
\[\pi_3(S^2) = \mathbb{Z}\] This came as a big shock when Heinz Hopf
first discovered it in the 1930's; before then, people had no idea how
sneaky homotopy groups were!

There is a beautiful way to compute an integer called the ``Hopf
invariant'' that keeps track of the homotopy class of a map from the
3-sphere to the 2-sphere. There are lots of nice ways to compute it, but
alas, I only have time to briefly sketch one! Suppose that the map
\(f\colon S^3\to S^2\) is smooth (otherwise we can always smooth it up).
Then most points \(p\) in \(S^2\) have the property that the points
\(x\) in \(S^3\) with \(f(x) = p\) form a ``link'': a bunch of knots in
\(S^3\). If we take two different points in \(S^2\) with this property,
we get two links. From these two links we can compute an integer called
the ``linking number'': for example, we can just draw these two links
and count the times one crosses over or under the other (with
appropriate plus or minus signs for each crossing). This number turns
out not to depend on how we picked the two points! Moreover, it only
depends on the homotopy class of \(f\). It's called the Hopf invariant
of \(f\).

Moving up one dimension, it turns out that \[\pi_4(S^3) = \mathbb{Z}/2\]
Here \(\mathbb{Z}/2\) is the group with two elements, usually written
\(0\) and \(1\), with addition \(\mod 2\). Why only two homotopy classes
of maps from \(S^4\) to \(S^3\)? Well, you can compute something like
the Hopf invariant for these maps, exactly as we did before, but the
thing is, links in 4 dimensions are easy to unlink. You can unlink
something like \[
  \begin{tikzpicture}
    \begin{knot}
      \strand[thick] (0,0)
        to [out=down,in=up] (0.75,-1.5)
        to [out=down,in=up] (0,-3);
      \strand[thick] (0.75,0)
        to [out=down,in=up] (0,-1.5)
        to [out=down,in=up] (0.75,-3);
      \flipcrossings{2}
    \end{knot}
  \end{tikzpicture}
\] and make it look like \[
  \begin{tikzpicture}
    \begin{knot}
      \strand[thick] (0,0)
        to (0,-3);
      \strand[thick] (0.75,0)
        to (0.75,-3);
    \end{knot}
  \end{tikzpicture}
\] so the linking number in 4 dimensions is only defined \(\mod 2\).
Thus the ``Hopf invariant'' is only defined \(\mod 2\).

The exact same thing happens in higher dimensions, too, so in fact we
have:
\[\pi_{k+1}(S^k) = \mathbb{Z}/2 \quad\text{for any}\quad k \geqslant 3.\]
This illustrates an important general fact: when the dimensions get high
enough, there's more room to wiggle things around, and as we keep
jacking up the dimension, homotopy groups simplify a bit and settle down
after a while. This is the idea behind ``stable homotopy theory''.

Let's look at some more examples. We have

\begin{itemize}
\tightlist
\item
  \(\pi_3(S^1) = \{0\}\)
\item
  \(\pi_4(S^2) = \mathbb{Z}/2\)
\item
  \(\pi_5(S^3) = \mathbb{Z}/2\)
\item
  \(\pi_6(S^4) = \mathbb{Z}/2\)
\end{itemize}

and so on:

\[\pi_{k+2}(S^k) = \mathbb{Z}/2 \quad\text{for any}\quad k\geqslant 2\]

Sadly, I do \emph{not} understand why this is true. How do you wrap a
4-sphere around a 2-sphere in an interesting way? Dunno.

(Thanks to Dan Christensen, an answer appears at the end of this post.)

Plunging on undeterred, we have:

\begin{itemize}
\tightlist
\item
  \(\pi_4(S^1) = \{0\}\)
\item
  \(\pi_5(S^2) = \mathbb{Z}/2\)
\item
  \(\pi_6(S^3) = \mathbb{Z}/12\)
\item
  \(\pi_7(S^4) = \mathbb{Z}\oplus\mathbb{Z}/12\)
\item
  \(\pi_8(S^5) = \mathbb{Z}/24\)
\item
  \(\pi_9(S^6) = \mathbb{Z}/24\)
\end{itemize}

and so on:

\[\pi_{k+3}(S^k) = \mathbb{Z}/24 \quad\text{for any}\quad k\geqslant 5.\]

Here is where the magic number 24 comes in! What the above says is that
if \(k\) is large enough, there are exactly 24 different homotopy class
of maps from an \((k+3)\)-sphere to an \(k\)-sphere!

Now I should explain what this has to do with string theory. But first I
should say more about the homotopy groups of spheres. There are some
simple patterns worth knowing about. First,
\[\pi_m(S^1) = \{0\} \quad\text{for any}\quad m\geqslant 2.\] Second,
there is a nice formula for when the homotopy groups settle down as we
jack up the dimension:
\[\mbox{$\pi_{k+n}(S^k)$ is independent of $k$ as long as $k\geqslant n+2$.}\]
The homotopy groups can stabilize sooner, as we saw for \(n = 2\), but
never later, and often they stabilize right at \(k = n+2\). There is a
simple reason for this. We saw that \(\pi_{k+1}(S^k)\) stabilized at
\(k = 3\) because it's easy to unlink links in 4 or more dimensions.
Similarly, \(\pi_{k+n}(S^k)\) must stabilize by the time \(k = n+2\),
because it's easy to untie knotted \(n\)-dimensional surfaces in
\(2n+2\) or more dimensions!

For more on stable homotopy groups of spheres, try:

\begin{enumerate}
\def\labelenumi{\arabic{enumi})}
\setcounter{enumi}{1}
\item
  Douglas C. Ravenel, \emph{Complex Cobordism and Stable Homotopy Groups
  of Spheres}, Academic Press, New York, 1986.

  Douglas C. Ravenel, \emph{Nilpotence and Periodicity in Stable
  Homotopy Theory}, Princeton U.\ Press, Princeton, 1992.
\end{enumerate}

Ravenel also spoke at this conference and is a real expert on stable
homotopy groups of spheres. Unfortunately his talk was too high-powered
for me. The 2nd book above is a bit more forgiving to the amateur, but
the first one has lots of nice tables of stable homotopy groups of
spheres.

The relationship between homotopy groups of spheres and higher-
dimensional knot theory is a wonderful thing. James Dolan and I are
learning a lot about \(n\)-categories by pondering it. When I spoke to
him at the conference at Northwestern, it became clear that Andre Joyal
had also thought about it very deeply. Joyal is famous for his work
relating category theory, combinatorics and topology, and his way of
thinking about the homotopy groups of spheres reflects these interests.
He said a very fascinating thing; he said ``really we should call the
sphere spectrum the `true integers'\,''. I would like to explain
this\ldots{} but here things get a bit technical, and I am afraid they
will get a lot more technical when I get around to the string theory
stuff.

What's the ``sphere spectrum''? Well, roughly it's just the list of
spheres \(S^0\), \(S^1\), \(S^2\), \ldots, but the word ``spectrum''
refers to the way all these spaces are all related, all aspects of one
big thing.

Here's a nice way to think of it. Start with the integers. Normally we
think of these as just a set, or actually a group, since we can add
them. But if we avoid the sin of mistaking isomorphism for equality we
can think of them as a category.

I already began to explain this in my parable about the shepherd in
\protect\hyperlink{week99}{``Week 99''}. The shepherd started with the
category of finite sets and ``decategorified'' it to obtain the set of
natural numbers, associating to each finite set a natural number, its
number of elements. Taking disjoint unions of sets corresponds to
addition, the empty set corresponds to zero, and so on.

Okay. What are the \emph{integers} the decategorification of?

Well, we can imagine finite sets that have both ``positive'' and
``negative'' elements. The ``number of elements'' of such a set will be
the number of positive elements minus the number of negative elements.
This is a bit weird if we're talking about sheep, but perhaps not so
weird if we talk about positrons and electrons, which can annihilate
each other. (In \protect\hyperlink{week92}{``Week 92''} I explain what
I'm hinting at here: the relation between antiparticles and
adjunctions.)

Topologists prefer to speak of ``positively and negatively oriented
points''. We can draw a set of positively and negatively oriented points
like this: \[
  \begin{tikzpicture}
    \node at (0,0) {$-$};
    \node at (1,0) {$+$};
    \node at (2,0) {$+$};
    \node at (3,0) {$+$};
    \node at (4,0) {$+$};
    \node at (5,0) {$-$};
    \node at (6,0) {$-$};
  \end{tikzpicture}
\] We can add them by setting them side by side. But how do the
positively and negatively oriented points cancel? Well, remember, we're
trying to get a category! If finite lists of positively and negatively
oriented points are our objects, what are our morphisms? How about
tangles, like this: \[
  \begin{tikzpicture}
    \node[label=above:{$-$}] at (0,0) {};
    \node[label=above:{$+$}] at (1,0) {};
    \node[label=above:{$+$}] at (2,0) {};
    \node[label=above:{$+$}] at (3,0) {};
    \node[label=above:{$+$}] at (4,0) {};
    \node[label=above:{$-$}] at (5,0) {};
    \node[label=above:{$-$}] at (6,0) {};
    \begin{knot}
      \strand[thick] (0,0)
        to [out=down,in=down,looseness=2] (1,0);
      \strand[thick] (0,-3)
        to [out=up,in=up,looseness=2] (1,-3);
      \strand[thick] (2,0) to (2,-3);
      \strand[thick] (3,0)
        to [out=down,in=down,looseness=2] (6,0);
      \strand[thick] (4,0)
        to [out=down,in=down,looseness=2] (5,0);
    \end{knot}
    \node[label=below:{$+$}] at (0,-3) {};
    \node[label=below:{$-$}] at (1,-3) {};
    \node[label=below:{$+$}] at (2,-3) {};
  \end{tikzpicture}
\] These let us cancel or create positive and negative points in pairs.
Voila! The categorified integers! Just as the integers form a monoid
under addition, these form a monoidal category (see
\protect\hyperlink{week89}{``Week 89''} for these concepts). The
monoidal structure here is disjoint union, which we can denote with a
plus sign if we like. Similarly, we can write the empty set as 0.

Above it looks like I'm drawing a \(1\)-dimensional tangle in
\(2\)-dimensional space. To understand the ``commutativity'' of the
categorified integers we should work with \(1\)-dimensional tangles in
higher-dimensional space. If we consider them in \(3\)-dimensional
space, we have room to switch things around: \[
  \begin{tikzpicture}
    \begin{knot}[clip width=7]
      \node[label=above:{$+$}] at (0,0) {};
      \node[label=above:{$+$}] at (1,0) {};
      \begin{knot}
        \strand[thick] (1,0)
          to [out=down,in=up] (0,-2);
        \strand[thick] (0,0)
          to [out=down,in=up] (1,-2);
      \end{knot}
      \node[label=below:{$+$}] at (0,-2) {};
      \node[label=below:{$+$}] at (1,-2) {};
    \end{knot}
  \end{tikzpicture}
\] This gets us commutativity, as I explained in
\protect\hyperlink{week100}{``Week 100''}. Technically speaking, we get
a ``braided'' monoidal category. However, there are two different ways
to switch things around; for example, in addition to the above way there
is \[
  \begin{tikzpicture}
    \begin{knot}[clip width=7]
      \node[label=above:{$+$}] at (0,0) {};
      \node[label=above:{$+$}] at (1,0) {};
      \begin{knot}
        \strand[thick] (0,0)
          to [out=down,in=up] (1,-2);
        \strand[thick] (1,0)
          to [out=down,in=up] (0,-2);
      \end{knot}
      \node[label=below:{$+$}] at (0,-2) {};
      \node[label=below:{$+$}] at (1,-2) {};
    \end{knot}
  \end{tikzpicture}
\] To get rid of this problem (if you consider it a problem) we can work
with \(1\)-dimensional tangles in \(4\)-dimensional space, where we can
deform the first way of switching things to the second. We get a
``symmetric'' monoidal category. Working in higher dimensions doesn't
change anything: things have stabilized.

If we impose the extra condition that the morphisms \[
  \begin{tikzpicture}
    \begin{knot}
      \strand[thick] (0,0)
        to [out=up,in=up,looseness=2] (1,0);
    \end{knot}
    \node[label=below:{$+$}] at (0,0) {};
    \node[label=below:{$-$}] at (1,0) {};
  \end{tikzpicture}
\] and \[
  \begin{tikzpicture}
    \begin{knot}
      \strand[thick] (0,0)
        to [out=down,in=down,looseness=2] (1,0);
    \end{knot}
    \node[label=above:{$+$}] at (0,0) {};
    \node[label=above:{$-$}] at (1,0) {};
  \end{tikzpicture}
\] are inverses, as are \[
  \begin{tikzpicture}
    \begin{knot}
      \strand[thick] (0,0)
        to [out=up,in=up,looseness=2] (1,0);
    \end{knot}
    \node[label=below:{$-$}] at (0,0) {};
    \node[label=below:{$+$}] at (1,0) {};
  \end{tikzpicture}
\] and \[
  \begin{tikzpicture}
    \begin{knot}
      \strand[thick] (0,0)
        to [out=down,in=down,looseness=2] (1,0);
    \end{knot}
    \node[label=above:{$-$}] at (0,0) {};
    \node[label=above:{$+$}] at (1,0) {};
  \end{tikzpicture}
\] then all morphisms become invertible, so we have not just a monoidal
category but a monoidal groupoid --- a groupoid being a category with
all morphisms invertible (see \protect\hyperlink{week74}{``Week 74''}).
In fact, not only are morphisms invertible, so are all objects! By this
I mean that every object \(x\) has an object \(-x\) such that \(x + -x\)
and \(-x + x\) are isomorphic to \(0\). For example, if \(x\) is the
positively oriented point, \(-x\) is the negatively oriented point, and
vice versa. So we have not just a monoidal groupoid but a ``groupal
groupoid''. (I have adopted this charming terminology from James Dolan.)

Very nice. We seem to have avoided the sin of decategorification, and
are no longer treating the integers as a mere \emph{set} (or group, or
commutative group). We are treating them as a \emph{category} (or
groupal groupoid, or braided groupal groupoid, or symmetric groupal
groupoid).

On the other hand, it's a bit odd to say that \[
  \begin{tikzpicture}
    \begin{knot}
      \strand[thick] (0,0)
        to [out=up,in=up,looseness=2] (1,0);
    \end{knot}
    \node[label=below:{$+$}] at (0,0) {};
    \node[label=below:{$-$}] at (1,0) {};
  \end{tikzpicture}
\] and \[
  \begin{tikzpicture}
    \begin{knot}
      \strand[thick] (0,0)
        to [out=down,in=down,looseness=2] (1,0);
    \end{knot}
    \node[label=above:{$+$}] at (0,0) {};
    \node[label=above:{$-$}] at (1,0) {};
  \end{tikzpicture}
\] are inverses. This amounts to saying that the morphism: \[
  \begin{tikzpicture}
    \begin{knot}
      \strand[thick] (0,0)
        to [out=up,in=up,looseness=2] (1,0);
      \strand[thick] (0,0)
        to [out=down,in=down,looseness=2] (1,0);
    \end{knot}
    \node[fill=white] at (0,0) {$+$};
    \node[fill=white] at (1,0) {$-$};
  \end{tikzpicture}
\] is equal to the identity morphism from 0 to 0, which corresponds to
the empty picture: \[\phantom{.}\] Hmm. They sure don't \emph{look}
equal. We must be doing something wrong.

What are we doing wrong? We're committing the sin of decategorification:
treating isomorphisms as equations! We should treat the integers not as
a mere category, but as a \(2\)-category! See
\protect\hyperlink{week80}{``Week 80''} for the precise definition of
this concept; for now, it's enough to say that a \(2\)-category has
things called 2-morphisms going between morphisms. If we treat the
integers as a \(2\)-category, we can say there is a \(2\)-morphism going
from \[
  \begin{tikzpicture}
    \begin{knot}
      \strand[thick] (0,0)
        to [out=up,in=up,looseness=2] (1,0);
      \strand[thick] (0,0)
        to [out=down,in=down,looseness=2] (1,0);
    \end{knot}
    \node[fill=white] at (0,0) {$+$};
    \node[fill=white] at (1,0) {$-$};
  \end{tikzpicture}
\] to the identity morphism. This \(2\)-morphism has a nice geometrical
description in terms of a \(2\)-dimensional surface: the surface in 3d
space that's traced out as the above picture shrinks down to the empty
picture. It's hard to draw, but let me try: \[
  \begin{tikzpicture}
    \begin{scope}
      \begin{knot}
        \strand[thick] (0,0)
          to [out=up,in=up,looseness=2] (1,0);
        \strand[thick] (0,0)
          to [out=down,in=down,looseness=2] (1,0);
      \end{knot}
      \node[fill=white] at (0,0) {$+$};
      \node[fill=white] at (1,0) {$-$};
    \end{scope}
    \node at (1.75,0) {$\Longrightarrow$};
    \begin{scope}[shift={(2.5,0)},scale=0.8]
      \begin{knot}
        \strand[thick] (0,0)
          to [out=up,in=up,looseness=2] (1,0);
        \strand[thick] (0,0)
          to [out=down,in=down,looseness=2] (1,0);
      \end{knot}
      \node[fill=white] at (0,0) {\mbox{\scriptsize$+$}};
      \node[fill=white] at (1,0) {\mbox{\scriptsize$-$}};
    \end{scope}
    \node at (4,0) {$\Longrightarrow$};
    \begin{scope}[shift={(4.65,0)},scale=0.5]
      \begin{knot}
        \strand[thick] (0,0)
          to [out=up,in=up,looseness=2] (1,0);
        \strand[thick] (0,0)
          to [out=down,in=down,looseness=2] (1,0);
      \end{knot}
      \node[fill=white] at (0,0) {\mbox{\tiny$+$}};
      \node[fill=white] at (1,0) {\mbox{\tiny$-$}};
    \end{scope}
    \node at (5.75,0) {$\Longrightarrow$};
  \end{tikzpicture}
\] Okay, say we do this: treat the integers as a \(2\)-category. We
again are faced with a question: do we make all the \(2\)-morphisms
invertible? If we do, we get a ``\(2\)-groupoid'', or actually a
``groupal \(2\)-groupoid''. But again, to do so amounts to committing
the sin of decategorification. To avoid this sin, we should tread the
integers as a \(3\)-category. Etc, etc!

You may have noted how worlds of ever higher-dimensional topology are
automatically unfolding from our attempt to avoid the sin of
decategorification. This is what's so neat about \(n\)-categories. I
haven't told you how it all works, but let me summarize what's actually
happening here. Normally we treat the integers as the free group on one
generator, or else the free commutative group on one generator. But
groups and commutative groups are just the tip of the iceberg! The
following picture is similar to that in
\protect\hyperlink{week74}{``Week 74''}:
\vfill
\newpage
\begin{longtable}[]{@{}llll@{}}
\caption*{\(k\)-tuply groupal \(n\)-groupoids}\tabularnewline
\toprule
\begin{minipage}[b]{0.26\columnwidth}\raggedright
\strut
\end{minipage} & \begin{minipage}[b]{0.21\columnwidth}\raggedright
\(n=0\)\strut
\end{minipage} & \begin{minipage}[b]{0.21\columnwidth}\raggedright
\(n=1\)\strut
\end{minipage} & \begin{minipage}[b]{0.21\columnwidth}\raggedright
\(n=2\)\strut
\end{minipage}\tabularnewline
\midrule
\endfirsthead
\toprule
\begin{minipage}[b]{0.26\columnwidth}\raggedright
\strut
\end{minipage} & \begin{minipage}[b]{0.21\columnwidth}\raggedright
\(n=0\)\strut
\end{minipage} & \begin{minipage}[b]{0.21\columnwidth}\raggedright
\(n=1\)\strut
\end{minipage} & \begin{minipage}[b]{0.21\columnwidth}\raggedright
\(n=2\)\strut
\end{minipage}\tabularnewline
\midrule
\endhead
\begin{minipage}[t]{0.26\columnwidth}\raggedright
\(k=0\)\strut
\end{minipage} & \begin{minipage}[t]{0.21\columnwidth}\raggedright
sets\strut
\end{minipage} & \begin{minipage}[t]{0.21\columnwidth}\raggedright
groupoids\strut
\end{minipage} & \begin{minipage}[t]{0.21\columnwidth}\raggedright
\(2\)-groupoids\strut
\end{minipage}\tabularnewline
\begin{minipage}[t]{0.26\columnwidth}\raggedright
\strut
\end{minipage} & \begin{minipage}[t]{0.21\columnwidth}\raggedright
\strut
\end{minipage} & \begin{minipage}[t]{0.21\columnwidth}\raggedright
\strut
\end{minipage} & \begin{minipage}[t]{0.21\columnwidth}\raggedright
\strut
\end{minipage}\tabularnewline
\begin{minipage}[t]{0.26\columnwidth}\raggedright
\(k=1\)\strut
\end{minipage} & \begin{minipage}[t]{0.21\columnwidth}\raggedright
groups\strut
\end{minipage} & \begin{minipage}[t]{0.21\columnwidth}\raggedright
groupal groupoids\strut
\end{minipage} & \begin{minipage}[t]{0.21\columnwidth}\raggedright
groupal \(2\)-groupoids\strut
\end{minipage}\tabularnewline
\begin{minipage}[t]{0.26\columnwidth}\raggedright
\strut
\end{minipage} & \begin{minipage}[t]{0.21\columnwidth}\raggedright
\strut
\end{minipage} & \begin{minipage}[t]{0.21\columnwidth}\raggedright
\strut
\end{minipage} & \begin{minipage}[t]{0.21\columnwidth}\raggedright
\strut
\end{minipage}\tabularnewline
\begin{minipage}[t]{0.26\columnwidth}\raggedright
\(k=2\)\strut
\end{minipage} & \begin{minipage}[t]{0.21\columnwidth}\raggedright
commutative groups\strut
\end{minipage} & \begin{minipage}[t]{0.21\columnwidth}\raggedright
braided groupal groupoids\strut
\end{minipage} & \begin{minipage}[t]{0.21\columnwidth}\raggedright
braided groupal \(2\)-groupoids\strut
\end{minipage}\tabularnewline
\begin{minipage}[t]{0.26\columnwidth}\raggedright
\strut
\end{minipage} & \begin{minipage}[t]{0.21\columnwidth}\raggedright
\strut
\end{minipage} & \begin{minipage}[t]{0.21\columnwidth}\raggedright
\strut
\end{minipage} & \begin{minipage}[t]{0.21\columnwidth}\raggedright
\strut
\end{minipage}\tabularnewline
\begin{minipage}[t]{0.26\columnwidth}\raggedright
\(k=3\)\strut
\end{minipage} & \begin{minipage}[t]{0.21\columnwidth}\raggedright
`` "\strut
\end{minipage} & \begin{minipage}[t]{0.21\columnwidth}\raggedright
symmetric groupal groupoids\strut
\end{minipage} & \begin{minipage}[t]{0.21\columnwidth}\raggedright
weakly involutory groupal \(2\)-groupoids\strut
\end{minipage}\tabularnewline
\begin{minipage}[t]{0.26\columnwidth}\raggedright
\strut
\end{minipage} & \begin{minipage}[t]{0.21\columnwidth}\raggedright
\strut
\end{minipage} & \begin{minipage}[t]{0.21\columnwidth}\raggedright
\strut
\end{minipage} & \begin{minipage}[t]{0.21\columnwidth}\raggedright
\strut
\end{minipage}\tabularnewline
\begin{minipage}[t]{0.26\columnwidth}\raggedright
\(k=4\)\strut
\end{minipage} & \begin{minipage}[t]{0.21\columnwidth}\raggedright
`` "\strut
\end{minipage} & \begin{minipage}[t]{0.21\columnwidth}\raggedright
`` "\strut
\end{minipage} & \begin{minipage}[t]{0.21\columnwidth}\raggedright
strongly involutory groupal \(2\)-groupoids\strut
\end{minipage}\tabularnewline
\begin{minipage}[t]{0.26\columnwidth}\raggedright
\strut
\end{minipage} & \begin{minipage}[t]{0.21\columnwidth}\raggedright
\strut
\end{minipage} & \begin{minipage}[t]{0.21\columnwidth}\raggedright
\strut
\end{minipage} & \begin{minipage}[t]{0.21\columnwidth}\raggedright
\strut
\end{minipage}\tabularnewline
\begin{minipage}[t]{0.26\columnwidth}\raggedright
\(k=5\)\strut
\end{minipage} & \begin{minipage}[t]{0.21\columnwidth}\raggedright
`` "\strut
\end{minipage} & \begin{minipage}[t]{0.21\columnwidth}\raggedright
`` "\strut
\end{minipage} & \begin{minipage}[t]{0.21\columnwidth}\raggedright
`` "\strut
\end{minipage}\tabularnewline
\bottomrule
\end{longtable}

What are all these things? Well, an \(n\)-groupoid is an \(n\)-category
with all morphisms invertible, at least up to equivalence. An
\((k+n)\)-groupoid with only trivial \(j\)-morphisms for \(j < k\) can
be seen as a special sort of \(n\)-groupoid, which we call a
``\(k\)-tuply groupal \(n\)-groupoid''.

Part of Joyal's point was that we should really think of the integers as
the ``free \(k\)-tuply monoidal \(n\)-groupoid on one object''. Here the
idea is not to fix \(n\) and \(k\) once and for all --- this would only
prevent us from understanding the subtleties that show up when we
increase \(n\) and \(k\)! Instead, we should think of them as variable,
or perhaps consider the limit as they become large.

The other part of his point was that there's a correspondence between
\(n\)-groupoids and the information left in topological spaces when we
ignore all homotopy groups above dimension \(n\) --- so-called
``homotopy \(n\)-types''. Using this correspondence, the ``free
\(k\)-tuply monoidal \(n\)-groupoid on one object'' corresponds to the
homotopy \((k+n)\)-type of the \(k\)-sphere. Moreover, if we keep
jacking up \(k\), this stabilizes when \(k\geqslant n+2\). Actually, as
the dittos in the above chart suggest, it's a quite general fact that
the notion of \(k\)-tuply monoidal \(n\)-groupoid stabilizes for
\(k\geqslant n+2\).

Yet another point is that the pictures above explain the relation
between higher-dimensional knot theory and the homotopy groups of
spheres in a very vivid, direct way.

Okay. What about string theory and the magic number 24? Well, notice
that the pictures above started looking a bit like strings! Hmm\ldots.

Here's the idea, as far as I understand it. Presumably all but the hardy
have stopped reading this article by now, so I will pull out all the
stops. The string worldsheet is a Riemann surface so we'll need some
stuff about Riemann surfaces from \protect\hyperlink{week28}{``Week
28''}. Let \(M(g,n)\) be the moduli space of Riemann surfaces with genus
\(g\) and \(n\) punctures, and let \(G(g,n)\) be the corresponding
mapping class group. Since \(M(g,n)\) is the quotient of Teichm\"uller
space by \(G(g,n)\) and Teichm\"uller space is contractible, we have
\[M(g,n) = BG(g,n)\] where ``\(\mathcal{B}\)'' means ``classifying
space''. There's a natural inclusion \[G(g,n)\hookrightarrow G(g+1,n)\]
defined by sewing an torus with two punctures onto your genus-\(g\)
surface with \(n\) punctures, which increases the genus by 1. Let's
define \(G(\infty,n)\) to be direct limit as \(g\to\infty\), and let
\[M(\infty,n) = BG(\infty,n).\] Now it turns out \(M(\infty,1)\) has a
kind of product on it. The reason is that there are products
\[M(g,1)\times M(h,1)\to M(g+h,1)\] given sewing two surfaces together
with a 3-punctured sphere. Using this product we can define the group
completion \[M(\infty,1)^+\] and the result Tillmann stated which got me
so excited was that \[\pi_3(M(\infty,1)^+) = \mathbb{Z}/24 \oplus H\]
for some unknown group \(H\). Since this is really a result about the
mapping class groups of surfaces, it \emph{must} have something to do
with how conformal field theories always give projective representations
of these mapping class groups, with the ``phase ambiguity'' being of the
form \(\exp(2\pi ci/24)\), where \(c\) is the ``central charge''. No? I
just don't quite see why. Maybe someone will enlighten me.

Anyway, the way she proved this definitely ties right into the stuff
about stable homotopy groups of spheres. She used explicit maps between
the third stable homotopy group of spheres
\[\pi_{k+3}(S^k) = \mathbb{Z}/24 \quad\text{for}\quad k \geqslant 5\]
and \(\pi_3(M(\infty,1)^+)\)! And the way she got the map from the
latter to the former amounts to working with pictures I was drawing
above. To put it more precisely, in

\begin{enumerate}
\def\labelenumi{\arabic{enumi})}
\setcounter{enumi}{2}
\tightlist
\item
  ``Higher-dimensional algebra and topological quantum field theory'',
  by John Baez and James Dolan, \emph{Jour. Math. Phys.} \textbf{36}
  (1995), 6073--6105.  Also available as \href{https://arxiv.org/abs/q-alg/9503002}{\texttt{q-alg/9503002}}.
\end{enumerate}
\noindent
we argue that framed \(n\)-dimensional surfaces embedded in
\((n+k)\)-dimensions should be described by the free \(k\)-tuply
monoidal \(n\)-category with duals on one object. This should map down
to the free \(k\)-tuply groupal \(n\)-groupoid on one object, by the
usual yoga of ``freeness''. Taking \(n = 3\) and \(k\) sufficiently
large, we should obtain a homomorphism from the mapping class group of
any Riemann surface to the third stable homotopy group of spheres!
Presumably the idea is that in the limit of large genus this
homomorphism is onto!

Of course, Tillmann doesn't prove her result using the still-nascent
formalism of \(n\)-categories, but I think it will eventually be
possible. (Also, my rough argument applies to Riemann surfaces with no
punctures, while she considers those with one puncture, but various
things she said make me think this might not be such a big deal.) The
real puzzle is this: what does \[\pi_3(M(\infty,n)^+)\] have to do with
central extensions of \(G(g,n)\) for finite \(g\)? If I could figure
this out I'd be very happy.

\begin{center}\rule{0.5\linewidth}{0.5pt}\end{center}

\textbf{Addendum:} Dan Christensen answered a puzzle above. Here's how
to get a nontrivial element of \[\pi_4(S^2).\] Take the map
\(f\colon S^3\to S^2\) generating \[\pi_3(S^2)\] and compose it with the
map \(g\colon S^4\to S^3\) generating \[\pi_4(S^3)\] (which, by the way,
is obtained from \(f\) by ``suspension'') to obtain the desired map from
\(S^4\) to \(S^2\). This is an instance of a very general trick:
composing elements of homotopy groups of spheres to get new ones!

\begin{center}\rule{0.5\linewidth}{0.5pt}\end{center}

\begin{quote}
\emph{Think of one and minus one. Together they add up to zero, nothing,
nada, niente, right? Picture them together, then picture them
separating, peeling apart\ldots. Now you have something, you have two
somethings, where you once had nothing.}

--- John Updike, Roger's Version
\end{quote}

\[
  \begin{tikzpicture}
    \begin{knot}
      \strand[thick] (0,0)
        to [out=up,in=up,looseness=2] (1,0);
    \end{knot}
    \node[label=below:{$+$}] at (0,0) {};
    \node[label=below:{$-$}] at (1,0) {};
  \end{tikzpicture}
\]

\hypertarget{week103}{%
\section{April 26, 1997}\label{week103}}

As I segue over from the homotopy theory conference at Northwestern
University to the conference on higher-dimensional algebra and physics
that took place right after that, it's a good time to mention Ronnie
Brown's web page:

\begin{enumerate}
\def\labelenumi{\arabic{enumi})}
\tightlist
\item
  Ronald Brown, Higher-dimensional group theory,
  \href{https://web.archive.org/web/19970629093438/http://www.bangor.ac.uk/~mas010/hdaweb2.htm}{\texttt{https://web.archive.org/web/}}
\href{https://web.archive.org/web/19970629093438/http://www.bangor.ac.uk/~mas010/hdaweb2.htm}{\texttt{19970629093438/http://www.bangor.ac.uk/\textasciitilde{}mas010/hdaweb2.htm}}
\end{enumerate}
\noindent
Brown is the one who coined the phrase ``higher-dimensional algebra'',
and for many years he has been developing this subject, primarily as a
tool for doing homotopy theory. I wrote a bit about his ideas two years
ago, in \protect\hyperlink{week53}{``Week 53''}. A lot has happened in
higher-dimensional algebra since then, and the web page above is a good
place to get an overview of it. It includes a nice bibliography on the
subject. 

The Workshop on Higher Category Theory and Physics was exciting because
it pulled together a lot of people working on the interface between
these two subjects, many of whom had never before met. It was organized
by Ezra Getzler and Mikhail Kapranov. Getzler is probably most
well-known for his proof of the Atiyah--Singer index theorem. This
wonderful theorem captured the imagination of mathematical physicists
for many years starting in the 1960s. The reason is that it relates the
topology of manifolds to the the solutions of partial differential
equations on these manifolds, and thus ushered in a new age of
applications of topology to physics. In the 1980s, working with ideas
that Witten came up with, Getzler found a nice ``supersymmetric proof''
of the Atiyah--Singer theorem. Later Getzler turned to other things, such
as the use of ``operads'' (see \protect\hyperlink{week42}{``Week 42''})
to study conformal field theory (which shows up naturally in string
theory). Kapranov has also done a lot of work with operads and conformal
field theory, and many other things, but I first learned about him
through his paper with Voevodsky on ``braided monoidal
\(2\)-categories'' (see \protect\hyperlink{week4}{``Week 4''}). This got
me very excited since it turned me on to many of the main themes of
\(n\)-category theory.

Alas, my description of this fascinating conference will be terse and
dry in the extreme, since I am flying to Warsaw in 3 hours for a quantum
gravity workshop. I'll just mention a few papers that cover some of the
themes of this conference. Ross Street gave two talks on Batanin's
definition of weak \(n\)-categories (and even weak
\(\omega\)-categories), which one can get as follows:

\begin{enumerate}
\def\labelenumi{\arabic{enumi})}
\setcounter{enumi}{2}
\tightlist
\item
  Ross Street, \emph{The role of Michael Batanin's monoidal globular
  categories}, in \emph{Higher Category Theory}, eds.\ E. Getzler and 
  M. Kapranov, \emph{Contemp.\ Math.} \textbf{230}, AMS, Providence, Rhode     
   Island, 1998, pp.\ 99–116.
\end{enumerate}

Subsequently Batanin has written a more thorough paper on his
definition:

\begin{enumerate}
\def\labelenumi{\arabic{enumi})}
\setcounter{enumi}{3}
\tightlist
\item
  Michael Batanin, ``Monoidal globular categories as a natural
  environment for the theory of weak \(n\)-categories'', \emph{Adv.
  Math} \textbf{136} (1998), 39--103.  
\end{enumerate}

I gave a talk on Dolan's and my definition of weak \(n\)-categories,
which one can get as follows:

\begin{enumerate}
\def\labelenumi{\arabic{enumi})}
\setcounter{enumi}{4}
\tightlist
\item
  John Baez, ``An introduction to \(n\)-categories'', in \emph{7th Conference on 
  Category Theory and Computer Science,} eds.\ Eugenio Moggi and Giuseppe   
  Rosolini, Lecture Notes in Computer Science \textbf{1290}, Springer, Berlin, 1997. 
  Also available as
  \href{https://arxiv.org/abs/q-alg/9705009}{\texttt{q-alg/9705009}}.
\end{enumerate}
\noindent
Unfortunately Tamsamani was not there to present \emph{his} definition
of weak \(n\)-categories, but at least I have learned how to get his
papers electronically:

\begin{enumerate}
\def\labelenumi{\arabic{enumi})}
\setcounter{enumi}{5}
\item
  Zouhair Tamsamani, ``Sur des notions de \(\infty\)-categorie et
  \(\infty\)-groupoide non-strictes via des ensembles
  multi-simpliciaux''.  Also available as
  \href{https://arxiv.org/abs/alg-geom/9512006}{\texttt{alg-geom/9512006}}.

  Zouhair Tamsamani, ``Equivalence de la theorie homotopique des
  \(n\)-groupoides et celle des espaces topologiques \(n\)-tronques''.
  Also available as
  \href{https://arxiv.org/abs/alg-geom/9607010}{\texttt{alg-geom/9607010}}.
\end{enumerate}
\noindent
Also, Carlos Simpson has written an interesting paper using Tamsamani's
definition:

\begin{enumerate}
\def\labelenumi{\arabic{enumi})}
\setcounter{enumi}{6}
\tightlist
\item
  Carlos Simpson, ``A closed model structure for \(n\)-categories,
  internal Hom, \(n\)-stacks and generalized Seifert-Van Kampen''.
  Also available as
  \href{https://arxiv.org/abs/alg-geom/9704006}{\texttt{alg-geom/9704006}}.
\end{enumerate}

In a different but related direction, Masahico Saito discussed a paper
with Scott Carter and Joachim Rieger in which they come up with a nice
purely combinatorial description of all the ways to embed
\(2\)-dimensional surfaces in \(4\)-dimensional space:

\begin{enumerate}
\def\labelenumi{\arabic{enumi})}
\setcounter{enumi}{7}
\tightlist
\item
  J. Scott Carter, Joachim H. Rieger and Masahico Saito, ``A
  combinatorial description of knotted surfaces and their isotopies'',
   \emph{Adv. Math.} \textbf{127} (1997), 1--51.  
\end{enumerate}
\noindent
My student Laurel Langford has translated their work into \(n\)-category
theory and shown that ``unframed unoriented 2-tangles form the free
braided monoidal \(2\)-category on one unframed self-dual object'':

\begin{enumerate}
\def\labelenumi{\arabic{enumi})}
\setcounter{enumi}{8}
\tightlist
\item
  John Baez and Laurel Langford, ``2-Tangles'', \emph{Lett. Math. Phys.}
\textbf{43} (1998), 187--197.   (With many typos.)   Also available as
  \href{https://arxiv.org/abs/q-alg/9703033}{\texttt{q-alg/9703033}}.
\end{enumerate}
\noindent
This paper summarizes the results; the proofs will appear later.

While I was there, Carter also gave me a very nice paper he'd done with
Saito and Louis Kauffman. This paper discusses 4-manifolds and also
2-dimensional surfaces in \(3\)-dimensional space, again getting a
purely combinatorial description which is begging to be translated into
\(n\)-category theory:

\begin{enumerate}
\def\labelenumi{\arabic{enumi})}
\setcounter{enumi}{9}
\tightlist
\item
  J. Scott Carter, Louis H. Kauffman and Masahico Saito,
  ``Diagrammatics, singularities, and their algebraic interpretations'', in
   10th Brazilian Topology Meeting (Sao Carlos, 1996), \emph{Mat. Contemp.} Vol. 13,
   1997.    Draft version available as
 \href{https://citeseerx.ist.psu.edu/pdf/c1ea0a98e7d5a6bd9ad18a695da412ad5823d610}{https://} \break
\href{https://citeseerx.ist.psu.edu/pdf/c1ea0a98e7d5a6bd9ad18a695da412ad5823d610}{citeseerx.ist.psu.edu/pdf/c1ea0a98e7d5a6bd9ad18a695da412ad5823d610}
\end{enumerate}

I am sorry not to describe these papers in more detail, but I've been
painfully busy lately. (In fact, I am trying to figure out how to reform
my life to give myself more spare time. I think the key is to say ``no''
more often.)

Thanks to Justin Roberts for pointing out an error in
\protect\hyperlink{week102}{``Week 102''}. The phase ambiguity in
conformal field theories is not necessarily a 24th root of unity; it's
\(\exp(2\pi ic/24)\) where \(c\) is the central charge of the associated
Virasoro representation. This is a big hint as far as my puzzle goes.

Also I thank Dan Christensen for helping me understand \(\pi_4(S^2)\) in
a simpler way, and Scott Carter for a fascinating letter on the themes
of \protect\hyperlink{week102}{``Week 102''}. Alas, I have been too busy
to reply adequately to these nice emails!

Gotta run\ldots.

\hypertarget{week104}{%
\section{June 8, 1997}\label{week104}}

A couple of months ago I flew up to Corvallis, Oregon to an AMS meeting.
The AMS, in case you're unfamiliar with it, is the American Mathematical
Society. They have lots of regional meetings with special sessions on
various topics. One reason I went to this one is that there was a
special session on octonions, organized by Tevian Dray and Corinne
Manogue.

After the real numbers come the complex numbers, and after the complex
numbers come the quaternions, and after the quaternions come the
octonions, the most mysterious of all. The real numbers, complex
numbers, and quaternions have lots of applications to physics. What
about the octonions? Aren't they good for something too? This question
has been bugging me for a while now.

In fact, it bugs me so much that I decided to go to Corvallis to look
for clues. After all, in addition to Tevian Dray and Corinne Manogue ---
the former a mathematician, the latter a physicist, both deeply
interested in octonions --- a bunch of other octonion experts were going
to be there. One was my friend Geoffrey Dixon. I told you about him in
\protect\hyperlink{week59}{``Week 59''}. He wrote a book on the complex
numbers, quaternions and octonions and their role in physics. He has a
theory of physics in which these are related to electromagnetism, the
weak force, and the strong force, respectively. It's a bit far out, but
far from crazy! In fact, it's fascinating.

After writing about his book I got in touch with him in Cambridge,
Massachusetts. I found out that his other main interest in life, besides
the octonions, is the game Myst. This is probably not a coincidence. In
both the main question is ``What the heck is really going on here?'' He
has Myst all figured out, but he loves watching people play it, so he
got me to play it for a while. Someday I will buy a CD-ROM drive and
waste a few weeks on that game. Anyway, I got to know him back in the
summer of 1995, so it was nice to see him again in Corvallis.

Another octonion expert is Tony Smith. He too has a far-out but
fascinating theory of physics involving octonions! I wrote about his
stuff in \protect\hyperlink{week91}{``Week 91''}. I had never met him
before the Corvallis conference, but I instantly recognized him when I
met him, because there's a picture of him wearing a cowboy hat on his
homepage. It turns out he always wears that hat. He is a wonderful
repository of information concerning octonions and other interesting
things. He is also a very friendly and laid-back sort of guy. He lives
in Atlanta, Georgia.

I also met another octonion expert I hadn't known about, Tony Sudbery,
from York. (The original York, not the ``new'' one.) He gave a talk on
``The Exceptions that Prove the Rule''. The octonions are related to a
host of other mathematical structures in a very spooky way. In all sorts
of contexts, you can classify algebraic structures and get a nice
systematic infinite list together with a finite number of exceptions.
What's spooky is how the exceptions in one context turn out to be
related to the exceptions in some other context. These relationships are
complicated and mysterious in themselves. It's as if there were a hand
underneath the water and all we see is the fingers poking out here and
there. There seems to be some ``unified theory of exceptions'' waiting
to be discovered, and the octonions must have something to do with it. I
figure that to really understand what the octonions are good for, we
need to understand this ``unified theory of exceptions''.

Let's start by recalling what the octonions are!

I presume you know the real numbers. The complex numbers are things like
\[a + bi\] where \(a\) and \(b\) are real. We can multiply them using
the rule \[i^2 = -1\] They may seem mysterious when you first meet them,
but they lose their mystery when you see they are just a nice way of
keeping track of rotations in the plane.

Similarly, the quaternions are guys like \[a + bi + cj + dk\] which we
can multiply using the rules \[i^2 = j^2 = k^2 = -1\] and \[
  \begin{gathered}
    ij = k, jk = i, ki = j,
  \\ji = -k, kj = -i, ik = -j
  \end{gathered}
\] They aren't commutative, but they are still associative. Again they
may seem mysterious at first, but they lose their mystery when you see
that they are just a nice way of keeping track of rotations in 3 and 4
dimensions. Rotations in more than 2 dimensions don't commute in
general, so the quaternions had \emph{better} not commute. In fact,
Hamilton didn't invent the quaternions to study rotations --- his goal
was merely to cook up a ``division algebra'', where you could divide by
any nonzero element (see \protect\hyperlink{week82}{``Week 82''}).
However, after he discovered the quaternions, he used them to study
rotations and angular momentum. Nowadays people tend instead to use the
vector cross product, which was invented later by Gibbs. The reason is
that in the late 1800s there was a big battle between the fans of
quaternions and the fans of vectors, and the quaternion crowd lost. For
more on the history of this stuff, see:

\begin{enumerate}
\def\labelenumi{\arabic{enumi})}
\tightlist
\item
  Michael J. Crowe, \emph{A History of Vector Analysis}, Dover, Mineola,
  2011.
\end{enumerate}

Octonions were invented by Cayley later on in the 1800s. For these, we
start with \emph{seven} square roots of \(-1\), say \(e_1\) up to
\(e_7\). To learn how multiply these, draw the following diagram:
\[\includegraphics[max width=0.65\linewidth]{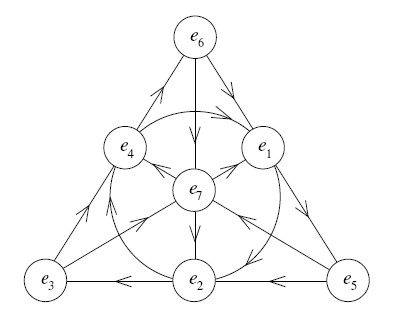}\] 
Draw
a triangle, draw a line from each vertex to the midpoint of the opposite
edge, and inscribe a circle in the triangle. Label the 7 points shown
with \(e_1\) through \(e_7\) --- it doesn't matter how, I've just drawn
my favorite way. Draw arrows on the edges of the triangle going around
clockwise, draw arrows on the circle also going around clockwise, and
draw arrows on the three lines pointing from each vertex of the triangle
to the midpoint of the opposite edge. Come on, \emph{do} it! I'm doing all this
work for you\ldots{} you should do some, too.

Okay. Now you have your very own octonion multiplication table. Notice
that there are six lines and a circle in your picture. Each one of these
gives us a copy of the quaternions inside the octonions. For example,
say you want to multiply \(e_6\) and \(e_7\). You notice that the the
vertical line says ``\(e_6\), \(e_7\), \(e_2\)'' on it as we follow the
arrow down. Thus, just as for \(i\), \(j\), and \(k\) in the
quaternions, we have \[
  \begin{gathered}
    e_6 e_7 =  e_2,   e_7 e_2 =  e_6,   e_2 e_6 =  e_7
  \\e_7 e_6 = -e_2,   e_2 e_7 = -e_6,   e_6 e_2 = -e_7
  \end{gathered}
\] 
So in particular we have \(e_6 e_7 = e_2\).

In case you lose your octonion table, don't worry: you don't really need
to remember the \emph{names} of those 7 square roots of \(-1\) and their
positions on the chart. You just need to remember the geometry of the
chart itself. Names are arbitrary and don't really matter, unless you're
talking to someone else, in which case you have to agree on them.

If you want to see spiffy high-tech octonion multiplication tables,
check out the following websites:

\begin{enumerate}
\def\labelenumi{\arabic{enumi})}
\setcounter{enumi}{1}
\item
  Tony Smith, \href{https://web.archive.org/web/19990203054838/http://galaxy.cau.edu/tsmith/TShome.html}{\texttt{https://web.archive.org/web/19990203054838/http://galaxy.}}  \href{https://web.archive.org/web/19990203054838/http://galaxy.cau.edu/tsmith/TShome.html}{\texttt{cau.edu/tsmith/TShome.html}}
\item
  Geoffrey Dixon, \url{http://www.7stones.com} 
\end{enumerate}

What's so great about the octonions? They are not commutative, and
worse, they are not even \emph{associative}. What's great about them is
that they form a division algebra, meaning you can divide by any nonzero
octonion. Better still, they form a ``normed'' division algebra. Just as
with the reals, complexes, and quaternions, we can define the norm of
the octonion
\[x = a_0 + a_1 e_1 + a_2 e_2 + a_3 e_3 + a_4 e_4 + a_5 e_5 + a_6 e_6 + a_7 e_7\]
to be
\[|x| = \sqrt{a_0^2 + a_1^2 + a_2^2 + a_3^2 + a_4^2 + a_5^2 + a_6^2 + a_7^2}.\]

What makes them a ``normed division algebra'' is that \[|xy| = |x||y|.\]
It's a wonderful fact about the world that the reals, complexes,
quaternions and octonions are the \emph{only} normed division algebras.
That's it!

However, the octonions remain mysterious, at least to me. They are
related to rotations in 7 and 8 dimensions, but not as simply as one
might hope. After all, rotations in \emph{any} number of dimensions are
still associative. Where is this nonassociative business coming from? I
don't really know. This question really bugs me.

A while ago, in \protect\hyperlink{week95}{``Week 95''}, I summarized a
paper by John Schwarz on supersymmetric Yang--Mills theory and why it
works best in dimensions 3, 4, 6, and 10. Basically, only in these
dimensions can you cook up spin-\(1/2\) particles that have as many
physical degrees of freedom as massless spin-\(1\) particles. I sort of
explained why. This in turn allows a symmetry between fermions and gauge
bosons. I didn't explain how \emph{this} works\ldots{} it seems pretty
tricky to me\ldots{} but anyway, it works.

So far, so good. But Schwarz wondered: is it a coincidence that the
numbers 3, 4, 6, and 10 are just two more than the numbers 1, 2, 4, and
8 --- the dimensions of the reals, complexes, quaternions, and
octonions?

Apparently not! The following papers explain what's going on:

\begin{enumerate}
\def\labelenumi{\arabic{enumi})}
\setcounter{enumi}{3}
\item
  Corinne A. Manogue and Joerg Schray, ``Finite Lorentz transformations
  automorphisms, and division algebras'', \emph{Jour. Math. Phys.}
  \textbf{34} (1993), 3746--3767.

  Corinne A. Manogue and Joerg Schray, ``Octonionic representations of
  Clifford algebras and triality'', available as
  \href{https://arxiv.org/abs/hep-th/9407179}{\texttt{hep-th/9407179}}.
\item
  Anthony Sudbery, ``Division algebras, (pseudo)orthogonal groups and
  spinors'', \emph{Jour. Phys. A} \textbf{17} (1984), 939--955.  

  Anthony Sudbery, ``Seven types of incongruity'', handwritten notes.
\end{enumerate}

Here's the basic idea. Let

\begin{itemize}
\tightlist
\item
  \(\mathbb{R}\) be the real numbers
\item
  \(\mathbb{C}\) be the complex numbers
\item
  \(\mathbb{H}\) be the quaternions
\item
  \(\mathbb{O}\) be the octonions.
\end{itemize}

Let \(\mathrm{SO}(n,1)\) denote the Lorentz group in \(n+1\) dimensions.
Roughly speaking, this is the symmetry group of \((n+1)\)-dimensional
Minkowski spacetime. Let \(\mathfrak{so}(n,1)\) be the corresponding Lie
algebra (see \protect\hyperlink{week63}{``Week 63''} for a lightning
introduction to Lie algebras). Then it turns out that:

\begin{itemize}
\tightlist
\item
  \(\mathfrak{sl}(2,\mathbb{R}) = \mathfrak{so}(2,1)\)
\item
  \(\mathfrak{sl}(2,\mathbb{C}) = \mathfrak{so}(3,1)\)
\item
  \(\mathfrak{sl}(2,\mathbb{H}) = \mathfrak{so}(5,1)\)
\item
  \(\mathfrak{sl}(2,\mathbb{O}) = \mathfrak{so}(9,1)\)
\end{itemize}

This relates reals, complexes, quaternions and octonions to the Lorentz
group in dimensions 3, 4, 6, and 10, and explains the ``coincidence''
noted by Schwarz! But it requires some explanation. Roughly speaking, if
\(\mathrm{SL}(2,\mathbb{K})\) is the group of \(2\times2\) matrices with
determinant \(1\) whose entries lie in the division algebra
\(\mathbb{K} = \mathbb{R}, \mathbb{C}, \mathbb{H}, \mathbb{O}\), then
\(\mathfrak{sl}(2,\mathbb{K})\) is defined to be the Lie algebra of this
group. This is simple enough for \(\mathbb{R}\) or \(\mathbb{C}\).
However, one needs to be careful when defining the determinant of a
\(2\times2\) quaternionic matrix, since quaternions don't commute. One
needs to be even more careful in the octonionic case. Since octonions
aren't even associative, it's far from obvious what the group
\(\mathrm{SL}(2,\mathbb{O})\) would be, so defining the Lie algebra
``\(\mathfrak{sl}(2,\mathbb{O})\)'' requires a certain amount of
finesse. For the details, read the papers.

As Corinne Manogue explained to me, this relation between the octonions
and Lorentz transformations in 10 dimensions suggests some interesting
ways to use octonions in \(10\)-dimensional physics. As we all know, the
10th dimension is where string theorists live. There is also a nice
relation to Geoffrey Dixon's theory. This theory relates the
electromagnetic force to the complex numbers, the weak force to the
quaternions, and the strong force to octonions. How? Well, the gauge
group of electromagnetism is \(\mathrm{U}(1)\), the unit complex
numbers. The gauge group of the weak force is \(\mathrm{SU}(2)\), the
unit quaternions. The gauge group of the strong force is
\(\mathrm{SU}(3)\)\ldots.

Alas, the group \(\mathrm{SU}(3)\) is \emph{not} the unit octonions. The
unit octonions do not form a group since they aren't associative.
\(\mathrm{SU}(3)\) is related to the octonions more indirectly. The
group of symmetries (or technically, ``automorphisms'') of the octonions
is the exceptional group \(\mathrm{G}_2\), which contains
\(\mathrm{SU}(3)\). To get \(\mathrm{SU}(3)\), we can take the subgroup
of \(\mathrm{G}_2\) that preserves a given unit imaginary
octonion\ldots{} say \(e_1\). This is how Dixon relates
\(\mathrm{SU}(3)\) to the octonions.

However, why should one unit imaginary octonion be different from the
rest? Some sort of ``symmetry breaking'', presumably? It seems a bit ad
hoc. However, as Manogue explained, there is a nice way to kill two
birds with one stone. If we pick a particular unit imaginary octonion,
we get a copy of the complex numbers sitting inside the octonions, so we
get a copy of \(\mathfrak{sl}(2,\mathbb{C})\) sitting inside
\(\mathfrak{sl}(2,\mathbb{O})\), so we get a copy of
\(\mathfrak{so}(3,1)\) sitting inside \(\mathfrak{so}(9,1)\)! In other
words, we get a particular copy of the good old \(4\)-dimensional
Lorentz group sitting inside the 10-dimensional Lorentz group. So fixing
a unit imaginary octonion not only breaks the octonion symmetry group
\(\mathrm{G}_2\) down to the strong force symmetry group
\(\mathrm{SU}(3)\), it might also get us from \(10\)-dimensional physics
down to \(4\)-dimensional physics.

Cool, no? There are obviously a lot of major issues involved in turning
this into a full-fledged theory, and they might not work out. The whole
idea could be completely misguided! But it takes guts to do physics, so
it's good that Tevian Dray and Corinne Manogue are bravely pursuing this
idea.

Upon learning that there is a deep relation between \(\mathbb{R}\),
\(\mathbb{C}\), \(\mathbb{H}\), \(\mathbb{O}\) and the Lorentz group in
dimensions 3, 4, 6, 10, one is naturally emboldened to take seriously a
few more ``coincidences''. For example, in
\protect\hyperlink{week82}{``Week 82''} I described the Clifford
algebras \(C_n\) --- i.e., the algebras generated by \(n\) anticommuting
square roots of \(-1\). These Clifford algebras are relevant to
\(n\)-dimensional \emph{Euclidean} geometry, as opposed to the Clifford
algebras relevant to \(n\)-dimensional \emph{Lorentzian} geometry, which
appeared in \protect\hyperlink{week93}{``Week 93''}. They go like this:

\begin{itemize}
\tightlist
\item
  \(C_0 = \mathbb{R}\)
\item
  \(C_1 = \mathbb{C}\)
\item
  \(C_2 = \mathbb{H}\)
\item
  \(C_3 = \mathbb{H}\oplus\mathbb{H}\)
\item
  \(C_4 = \mathbb{H}(2)\)
\item
  \(C_5 = \mathbb{C}(4)\)
\item
  \(C_6 = \mathbb{R}(8)\)
\item
  \(C_7 = \mathbb{R}(8)\oplus\mathbb{R}(8)\)
\item
  \(C_8 = \mathbb{R}(16)\)
\end{itemize}

where \(\mathbb{K}(n)\) stands for \(n\times n\) matrices with entries
taken from \(\mathbb{K} = \mathbb{R}, \mathbb{C}, \mathbb{H}\), and
``\(\oplus\)'' stands for ``direct sum''. Note that \(C_8\) is the same
as \(16\times16\) matrices with entries taken from \(C_0\). That's part
of a general pattern called ``Bott periodicity'': in general,
\(C_{n+8}\) is the same as \(16\times16\) matrices with entries taken
from \(C_n\).

Now consider the dimension of the smallest real representation of
\(C_n\). It's easy to work this out if you keep in mind that the
smallest representation of \(\mathbb{K}(n)\) or
\(\mathbb{K}(n)\oplus \mathbb{K}(n)\) is on \(\mathbb{K}^n\) --- the
vector space consisting of \(n\)-tuples of elements of \(\mathbb{K}\).
We get

The dimension of the smallest real representation:

\begin{itemize}
\tightlist
\item
  of \(C_0\) is 1
\item
  of \(C_1\) is 2
\item
  of \(C_2\) is 4
\item
  of \(C_3\) is 4
\item
  of \(C_4\) is 8
\item
  of \(C_5\) is 8
\item
  of \(C_6\) is 8
\item
  of \(C_7\) is 8
\item
  of \(C_8\) is 16
\end{itemize}

Note that it increases at \(n = 1, 2, 4, 8\). These are the dimensions
of \(\mathbb{R}\), \(\mathbb{C}\), \(\mathbb{H}\), and \(\mathbb{O}\).
Coincidence?

No! Indeed, \(C_n\) has a representation on a \(k\)-dimensional real
vector space if and only if the unit sphere in that vector space,
\(S^{k-1}\), admits \(n\) linearly independent smooth vector fields. So
the above table implies that:

\begin{itemize}
\tightlist
\item
  The sphere \(S^0\) admits 0 linearly independent vector fields.
\item
  The sphere \(S^1\) admits 1 linearly independent vector fields.
\item
  The sphere \(S^3\) admits 3 linearly independent vector fields.
\item
  The sphere \(S^7\) admits 7 linearly independent vector fields.
\end{itemize}

These spheres are the unit real numbers, the unit complex numbers, the
unit quaternions, and the unit octonions, respectively! If you know
about normed division algebras, it's obvious that these sphere admit the
maximum possible number of linear independent vector fields: you can
just take a basis of vectors at one point and ``left translate'' it to
get a bunch of linearly independent vector fields.

Now --- Bott periodicity has period 8, and the octonions have dimension
8. And as we've seen, both have a lot to do with Clifford algebras. So
maybe there is a deep relation between the octonions and Bott
periodicity. Could this be true? If so, it would be good news, because
while octonions are often seen as weird exceptional creatures, Bott
periodicity is bigtime, mainstream stuff!

And in fact it \emph{is} true. More on Bott periodicity and the
octonions coming up next Week.

\begin{center}\rule{0.5\linewidth}{0.5pt}\end{center}

\textbf{Addendum:} Robert Helling provided some interesting further
comments on supersymmetric gauge theories and the division algebras,
which I quote below. He recommends the following reference:

\begin{enumerate}
\def\labelenumi{\arabic{enumi})}
\setcounter{enumi}{5}
\tightlist
\item
  J. M. Evans, ``Supersymmetric Yang--Mills theories and division
  algebras'', \emph{Nucl. Phys.} \textbf{B298} (1988), 92--108.
\end{enumerate}

and he writes:

\begin{quote}
Let me add a technical remark that I extract from Green, Schwarz, and
Witten, Vol 1, Appendix 4A.

The appearance of dimensions 3,4,6, and 10 can most easily been seen
when one tries to write down a supersymmetric gauge theory in arbitrary
dimension. This means we're looking for a way to throw in some spinors
to the Lagrangian of a pure gauge theory: \[-\frac{1}{4} F^2\] in a way that the
new Lagrangian is invariant (up to a total derivative) under some
infinitesimal variations. These describe supersymmetry if their
commutator is a derivative (a generator of spacetime translations). As
usual, we parameterize this variation by a parameter \(\varepsilon\),
but now \(\varepsilon\) is a spinor.

From people that have been doing this for their whole life we learn that
the following Ansatz is common:
\[\delta A_m = \frac{i}{2} \overline{\varepsilon} \Gamma_m \psi\]
\[\delta \psi = -\frac{1}{4} F_{mn} \Gamma^{mn} \varepsilon\] 
Here \(A\) is the
connection, \(F\) its field strength and \(\psi\) a spinor of a type to
be determined. I suppressed group indices on all these fields. They are
all in the adjoint representation. \(\Gamma\) are the generators of the
Clifford algebra described by John Baez before.

For the Lagrangian we try the usual Yang--Mills term and add a minimally
coupled kinetic term for the fermions:
\[-\frac{1}{4} F^2 + \frac{ig}{2} \psi^\dagger \Gamma^m D_m \psi\] 
Here \(D_m\) is the
gauge covariant derivative and \(g\) is some number that we can tune to
to make this vanish under the above variations. When we vary the first
term we find \(g = 1\). In fact everything cancels without considering a
special dimension except for the term that is trilinear in \(\psi\) that
comes from varying the connection in the covariant derivative in the
fermionic term. This reads something like
\[f_{abc} \overline{\varepsilon} \Gamma_m \psi^a \psi^b \Gamma^m \psi^c\]
where I put in the group indices and the structure constants
\(f_{abc}\). This has to vanish for other reasons since there is no
other trilinear term in the fermions available. And indeed, after you've
written out the antisymmetry of \(f\) explicitly and take out the
spinors since this should vanish for all choices of \(\psi\) and
\(\varepsilon\). We are left with an expression that is only made of
gammas. And in fact, this expression exactly vanishes in dimensions 3,
4, 6, and 10 due to a Fierz identity. (Sorry, I don't have time to work
this out more explicitly.)

This is related to the division algebra as follows (as explained in the
papers pointed out by John Baez): Take for concreteness \(d = 10\). Here
we go to a light-cone frame by using coordinates
\[x^+ = x^0 + x^9 \quad\text{and}\quad x^- = x^0 - x^1.\] Then we write
the \(\Gamma_m\) as block matrices where \(\Gamma_+\) and \(\Gamma_-\)
have the \(+\)/\(-\) unit matrix as blocks and the others have
\(\gamma_i\) as blocks where \(\gamma_i\) are the \(\mathrm{SO}(8)\)
Dirac matrices (\(i=1,...,9\)). But they are intimately related to the
octonions. Remember there is triality in \(\mathrm{SO}(8)\) which means
that we can treat left-handed spinors, right-handed spinors and vectors
on an equal basis (see \protect\hyperlink{week61}{``Week 61''},
\protect\hyperlink{week90}{``Week 90''},
\protect\hyperlink{week91}{``Week 91''}). Now I write out all three
indices of \(\gamma_i\). Because of triality I can use \(i\),\(j\),\(k\)
for spinor, dotted spinor and vector indices. Then it is known that \[
  \gamma_{ijk} =
  \begin{cases}
    c_{ijk} &\mbox{for $i,j,k<8$;}
  \\\delta_{ij} &\mbox{for $k=8$ (and $ijk$ permuted);}
  \\0 &\mbox{for more than two of $ijk$ equal to $8$.}
  \end{cases}
\] is a representation of \(\mathrm{Cliff}(8)\) if \(c_{ijk}\) are the
structure constants of the octonions (i.e.~\(e_i e_j = c_{ijk} e_k\) for
the 7 roots of \(-1\) in the octonions).

When you plug this representation of the \(\Gamma\)'s in the above
mentioned \(\gamma\) expression you will will find that it vanishes due
to the antisymmetry of the associator \[[a,b,c] = a(bc) - (ab)c\]

in the division algebras. This is my understanding of the relation of
supersymmetry to the divison algebras.

Robert
\end{quote}

\hypertarget{week105}{%
\section{June 21, 1997}\label{week105}}

There are some spooky facts in mathematics that you'd never guess in a
million years\ldots{} only when someone carefully works them out do they
become clear. One of them is called ``Bott periodicity''.

A 0-dimensional manifold is pretty dull: just a bunch of points.
1-dimensional manifolds are not much more varied: the only possibilities
are the circle and the line, and things you get by taking a union of a
bunch of circles and lines. \(2\)-dimensional manifolds are more
interesting, but still pretty tame: you've got your n-holed tori, your
projective plane, your Klein bottle, variations on these with extra
handles, and some more related things if you allow your manifold to go
on forever, like the plane, or the plane with a bunch of handles added
(possibly infinitely many!), and so on\ldots. You can classify all these
things. \(3\)-dimensional manifolds are a lot more complicated: nobody
knows how to classify them. \(4\)-dimensional manifolds are a \emph{lot}
more complicated: you can \emph{prove} that it's \emph{impossible} to
classify them --- that's called Markov's Theorem.

Now, you probably wouldn't have guessed that a lot of things start
getting simpler when you get up around dimension 5. Not everything, just
some things. You still can't classify manifolds in these high
dimensions, but if you make a bunch of simplifying assumptions you sort
of can, in ways that don't work in lower dimensions. Weird, huh? But
that's another story. Bott periodicity is different. It says that when
you get up to 8 dimensions, a bunch of things are a whole lot like in 0
dimensions! And when you get up to dimension 9, a bunch of things are a
lot like they were in dimension 1. And so on - a bunch of stuff keeps
repeating with period 8 as you climb the ladder of dimensions.

(Actually, I have this kooky theory that perhaps part of the reason
topology reaches a certain peak of complexity in dimension 4 is that the
number 4 is halfway between 0 and 8, topology being simplest in
dimension 0. Maybe this is even why physics likes to be in 4 dimensions!
But this is a whole other crazy digression and I will restrain myself
here.)

Bott periodicity takes many guises, and I already described one in
\protect\hyperlink{week104}{``Week 104''}. Let's start with the real
numbers, and then throw in \(n\) square roots of \(-1\), say
\(e_1,\ldots,e_n\). Let's make them ``anticommute'', so
\[e_i e_j = - e_j e_i\] when \(i\) is different from \(j\). What we get
is called the ``Clifford algebra'' \(C_n\). For example, when \(n = 1\)
we get the complex numbers, which we call C. When \(n = 2\) we get the
quaternions, which we call H, for Hamilton. When \(n = 3\) we
get\ldots{} the octonions?? No, not the octonions, since we always
demand that multiplication be associative! We get the algebra consisting
of \emph{pairs} of quaternions! We call that
\(\mathbb{H}\oplus\mathbb{H}\). When \(n = 4\) we get the algebra
consisting of \(2\times2\) \emph{matrices} of quaternions! We call that
\(\mathbb{H}(2)\). And it goes on, like this:

\begin{itemize}
\tightlist
\item
  \(C_0 = \mathbb{R}\)
\item
  \(C_1 = \mathbb{C}\)
\item
  \(C_2 = \mathbb{H}\)
\item
  \(C_3 = \mathbb{H}\oplus\mathbb{H}\)
\item
  \(C_4 = \mathbb{H}(2)\)
\item
  \(C_5 = \mathbb{C}(4)\)
\item
  \(C_6 = \mathbb{R}(8)\)
\item
  \(C_7 = \mathbb{R}(8)\oplus\mathbb{R}(8)\)
\item
  \(C_8 = \mathbb{R}(16)\)
\end{itemize}

Note that by the time we get to \(n = 8\) we just have \(16\times16\)
matrices of real numbers. And that's how it keeps going: \(C_{n+8}\) is
just \(16\times16\) matrices of guys in \(C_n\)! That's Bott periodicity
in its simplest form.

Actually right now I'm in Vienna, at the Schroedinger Institute, and one
of the other people visiting is Andrzej Trautman, who gave a talk the
other day on ``Complex Structures in Physics'', where he mentioned a
nice way to remember the above table. Imagine the day is only 8 hours
long, and draw a clock with 8 hours. Then label it like this: \[
  \begin{tikzpicture}
    \draw (0,0) circle[radius=2.65cm];
    \node[label=below:{$\mathbb{R}$}] at (90:2.3) {0};
    \node[label=below left:{$\mathbb{C}$}] at (45:2.3) {1};
    \node[label=left:{$\mathbb{H}$}] at (0:2.3) {2};
    \node[label={[label distance=-2mm]above left:{$\mathbb{H}\oplus\mathbb{H}$}}] at (-45:2.3) {3};
    \node[label=above:{$\mathbb{H}$}] at (-90:2.3) {4};
    \node[label=above right:{$\mathbb{C}$}] at (-135:2.3) {5};
    \node[label=right:{$\mathbb{R}$}] at (180:2.3) {6};
    \node[label={[label distance=-2mm]below right:{$\mathbb{R}\oplus\mathbb{R}$}}] at (135:2.3) {7};
    \foreach \a in {0,45,90,135,180,-135,-90,-45}
      \draw (\a:2.5) to (\a:2.65);
  \end{tikzpicture}
\] The idea here is that as the dimension of space goes up, you go
around the clock. One nice thing about the clock is that it has a
reflection symmetry about the axis from 3 o'clock to 7 o'clock. To use
the clock, you need to know that the dimension of the Clifford algebra
doubles each time you go up a dimension. This lets you figure out, for
example, that the Clifford algebra in 4 dimensions is not really
\(\mathbb{H}\), but \(\mathbb{H}(2)\), since the latter has dimension
\(16 = 2^4\).

Now let's completely change the subject and talk about rotations in
infinite-dimensional space! What's a rotation in infinite-dimensional
space like? Well, let's start from the bottom and work our way up. You
can't really rotate in 0-dimensional space. In \(1\)-dimensional space
you can't really rotate, you can only \emph{reflect} things\ldots{} but
we will count reflections together with rotations, and say that the
operations of multiplying by \(1\) or \(-1\) count as ``rotations'' in
\(1\)-dimensional space. In \(2\)-dimensional space we describe
rotations by \(2\times2\) matrices like \[
  \left(
    \begin{array}{cc}
      \cos t & -\sin t
    \\\sin t & \cos t
    \end{array}
  \right)
\] and since we're generously including reflections, also matrices like
\[
  \left(
    \begin{array}{cc}
      \cos t & \sin t
    \\\sin t & -\cos t
    \end{array}
  \right)
\] These are just the matrices whose columns are orthonormal vectors. In
3-dimensional space we describe rotations by \(3\times3\) matrices whose
columns are orthonormal, and so on. In n-dimensional space we call the
set of \(n\times n\) matrices with orthonormal columns the ``orthogonal
group'' \(\mathrm{O}(n)\).

Note that we can think of a rotation in 2 dimensions \[
  \left(
    \begin{array}{cc}
      \cos t & -\sin t
    \\\sin t & \cos t
    \end{array}
  \right)
\] as being a rotation in 3 dimensions if we just stick one more row and
one column like this: \[
  \left(
    \begin{array}{ccc}
      \cos t & -\sin t & 0
    \\\sin t & \cos t & 0
    \\ 0 & 0 & 1
    \end{array}
  \right)
\] This is just a rotation around the z axis. Using the same trick we
can think of any rotation in \(n\) dimensions as a rotation in \(n+1\)
dimensions. So we can think of \(\mathrm{O}(0)\) as sitting inside
\(\mathrm{O}(1)\), which sits inside \(\mathrm{O}(2)\), which sits
inside \(\mathrm{O}(3)\), which sits inside \(\mathrm{O}(4)\), and so
on! Let's do that. Then let's just take the \emph{union} of all these
guys, and we get\ldots{} \(\mathrm{O}(\infty)\)! This is the group of
rotations, together with reflections, in infinite dimensions.

(Now if you know your math, or you read
\protect\hyperlink{week82}{``Week 82''}, you'll realize that I didn't
really change the subject, since the Clifford algebra \(C_n\) is really
just a handy way to study rotations in \(n\) dimensions. But never
mind.)

Now \(\mathrm{O}(\infty)\) is a very big group, but it elegantly
summarizes a lot of information about rotations in all dimensions, so
it's not surprising that topologists have studied it. One of the thing
topologists do when studying a space is to work out its ``homotopy
groups''. If you hand them a space \(X\), and choose a point \(x\) in
this space, they will work out all the topologically distinct ways you
can stick an \(n\)-dimensional sphere in this space, where we require
that the north pole of the sphere be at \(x\). This is what they are
paid to do. We call the set of all such ways the homotopy group
\(\pi_n(X)\). For a more precise description, try
\protect\hyperlink{week102}{``Week 102''} --- but this will do for now.

So, what are the homotopy groups of \(\mathrm{O}(\infty)\)? Well, they
start out looking like this:

\begin{longtable}[]{@{}ll@{}}
\toprule
\(n\) & \(\pi_n(\mathrm{O}(\infty))\)\tabularnewline
\midrule
\endhead
\(0\) & \(\mathbb{Z}/2\)\tabularnewline
\(1\) & \(\mathbb{Z}/2\)\tabularnewline
\(2\) & \(0\)\tabularnewline
\(3\) & \(\mathbb{Z}\)\tabularnewline
\(4\) & \(0\)\tabularnewline
\(5\) & \(0\)\tabularnewline
\(6\) & \(0\)\tabularnewline
\(7\) & \(\mathbb{Z}\)\tabularnewline
\bottomrule
\end{longtable}

And then they repeat, modulo 8. Bott periodicity strikes again!

But what do they mean?

Well, luckily Jim Dolan has thought about this a lot. Discussing it
repeatedly in the little cafe we tend to hang out at, we came up with
the following story. Most of it is known to various people already, but
it came as sort of a revelation to us.

The zeroth entry in the table is easy to understand. \(\pi_0\) keeps
track of how many connected components your space has. The rotation
group \(\mathrm{O}(\infty)\) has two connected components: the guys that
are rotations, and the guys that are rotations followed by a reflection.
So \(\pi_0\) of \(\mathrm{O}(\infty)\) is \(\mathbb{Z}/2\), the group
with two elements. Actually this is also true for \(\mathrm{O}(n)\)
whenever \(n\) is high enough, namely \(1\) or more. So the zeroth entry
is all about ``reflecting''.

The first entry is a bit subtler but very important in physics. It means
that there is a loop in \(\mathrm{O}(\infty)\) that you can't pull
tight, but if you go around that loop \emph{twice}, you trace out a loop
that you \emph{can} pull tight. In fact this is true for
\(\mathrm{O}(n)\) whenever \(n\) is \(3\) or more. This is how there can
be spin-\(1/2\) particles when space is \(3\)-dimensional or higher.
There are lots of nice tricks for seeing that this is true, which I hope
the reader already knows and loves. In short, the first entry is all
about ``rotating 360 degrees and not getting back to where you
started''.

The second entry is zero.

The third entry is even subtler but also very important in modern
physics. It means that the ways to stick a 3-sphere into
\(\mathrm{O}(\infty)\) are classified by the integers, \(\mathbb{Z}\).
Actually this is true for \(\mathrm{O}(n)\) whenever \(n\) is \(5\) or
more. It's even true for all sorts of other groups, like all ``compact
simple groups''. But can I summarize this entry in a snappy phrase like
the previous nonzero entries? Not really. Actually a lot of applications
of topology to quantum field theory rely on this \(\pi_3\) business. For
example, it's the key to stuff like ``instantons'' in Yang--Mills theory,
which are in turn crucial for understanding how the pion gets its mass.
It's also the basis of stuff like ``Chern--Simons theory'' and ``\(BF\)
theory''. Alas, all this takes a while to explain, but let's just say
the third entry is about ``topological field theory in 4 dimensions''.

The fourth entry is zero.

The fifth entry is zero.

The sixth entry is zero.

The seventh entry is probably the most mysterious of all. From one point
of view it is the subtlest, but from another point of view it is
perfectly trivial. If we think of it as being about \(\pi_7\) it's very
subtle: it says that the ways to stick a 7-sphere into
\(\mathrm{O}(\infty)\) are classified by the integers. (Actually this is
true for \(\mathrm{O}(n)\) whenever \(n\) is \(7\) or more.) But if we
keep Bott periodicity in mind, there is another way to think of it: we
can think of it as being about \(\pi_{-1}\), since \(7 = -1 \mod 8\).

But wait a minute! Since when can we talk about \(\pi_n\) when \(n\) is
\emph{negative?!} What's a -1-dimensional sphere, for example?

Well, the idea here is to use a trick. There is a space very related to
\(\mathrm{O}(\infty)\), called \(k\mathrm{O}\). As with
\(\mathrm{O}(\infty)\), the homotopy groups of this space repeat modulo
8. Moreover we have:
\[\pi_n(\mathrm{O}(\infty)) = \pi_{n+1}(k\mathrm{O})\] Combining these
facts, we see that the very subtle \(\pi_7\) of \(\mathrm{O}(\infty)\)
is nothing but the very unsubtle \(\pi_0\) of \(k\mathrm{O}\), which
just keeps track of how many connected components \(k\mathrm{O}\) has.

But what \emph{is} \(k\mathrm{O}\)?

Hmm. The answer is very important and interesting, but it would take a
while to explain, and I want to postpone doing it for a while, so I can
get to the punchline. Let me just say that when we work it all out, we
wind up seeing that the seventh entry in the table is all about
\emph{dimension}.

To summarize:

\begin{itemize}
\tightlist
\item
  \(\pi_0(\mathrm{O}(\infty)) = \mathbb{Z}/2\) is about
  \textbf{reflecting}
\item
  \(\pi_1(\mathrm{O}(\infty)) = \mathbb{Z}/2\) is about \textbf{rotating
  360 degrees}
\item
  \(\pi_3(\mathrm{O}(\infty)) = \mathbb{Z}\) is about
  \textbf{topological field theory in 4 dimensions}
\item
  \(\pi_7(\mathrm{O}(\infty)) = \mathbb{Z}\) is about \textbf{dimension}
\end{itemize}

But wait! What do those numbers 0, 1, 3, and 7 remind you of?

Well, after I stared at them for a few weeks, they started to remind me
of the numbers 1, 2, 4, and 8. And \emph{that} immediately reminded me
of the reals, the complexes, the quaternions, and the octonions!

And indeed, there is an obvious relationship. Let \(n\) be 1, 2, 4, or
8, and correspondingly let \(\mathbb{A}\) stand for either the reals
\(\mathbb{R}\), the complex numbers \(\mathbb{C}\), the quaternions
\(\mathbb{H}\), or the octonions \(\mathbb{O}\). These guys are
precisely all the ``normed division algebras'', meaning that the obvious
sort of absolute value satisfies \[|xy| = |x||y|.\] Thus if we take any
guy \(x\) in \(\mathbb{A}\) with \(|x| = 1\), the operation of
multiplying by \(x\) is length-preserving, so it's a reflection or
rotation in \(\mathbb{A}\). This gives us a function from the unit
sphere in \(\mathbb{A}\) to \(\mathrm{O}(n)\), or in other words from
the \((n-1)\)-sphere to \(\mathrm{O}(n)\). We thus get nice elements of
\[\pi_0(\mathrm{O}(1)), \quad\pi_1(\mathrm{O}(2)), \quad\pi_3(\mathrm{O}(4)), \quad\pi_7(\mathrm{O}(8))
\] which turn out to be precisely why these particular homotopy groups
of \(\mathrm{O}(\infty)\) are nontrivial.

So now we have the following fancier chart:

\begin{itemize}
\tightlist
\item
  \(\pi_0(\mathrm{O}(\infty))\) is about \textbf{reflecting} and the
  \textbf{real numbers}
\item
  \(\pi_1(\mathrm{O}(\infty))\) is about \textbf{rotating 360 degrees}
  and the \textbf{complex numbers}
\item
  \(\pi_3(\mathrm{O}(\infty))\) is about \textbf{topological field
  theory in 4 dimensions} and the \textbf{quaternions}
\item
  \(\pi_7(\mathrm{O}(\infty))\) is about \textbf{dimension} and the
  \textbf{octonions}
\end{itemize}

Now this is pretty weird. It's not so surprising that reflections and
the real numbers are related: after all, the only ``rotations'' in the
real line are the reflections. That's sort of what \(1\) and \(-1\) are
all about. It's also not so surprising that rotations by 360 degrees are
related to the complex numbers. That's sort of what the unit circle is
all about. While far more subtle, it's also not so surprising that
topological field theory in 4 dimensions is related to the quaternions.
The shocking part is that something so basic-sounding as ``dimension''
should be related to something so erudite-sounding as the ``octonions''!

But this is what Bott periodicity does, somehow: it wraps things around
so the most complicated thing is also the least complicated.

That's more or less the end of what I have to say, except for some
references and some remarks of a more technical nature.

Bott periodicity for \(\mathrm{O}(\infty)\) was first proved by Raoul
Bott in 1959. Bott is a wonderful explainer of mathematics and one of
the main driving forces behind applications of topology to physics, and
a lot of his papers have now been collected in book form:

\begin{enumerate}
\def\labelenumi{\arabic{enumi})}
\tightlist
\item
  \emph{The Collected Papers of Raoul Bott}, ed.~R. D. MacPherson. Vol.
  1: \emph{Topology and Lie Groups (the 1950s)}. Vol. 2:
  \emph{Differential Operators (the 1960s)}. Vol. 3: \emph{Foliations
  (the 1970s)}. Vol. 4: \emph{Mathematics Related to Physics (the
  1980s)}. Birkhauser, Boston, 1994, 2355 pages total.
\end{enumerate}
\noindent
A good paper on the relation between \(\mathrm{O}(\infty)\) and Clifford
algebras is:

\begin{enumerate}
\def\labelenumi{\arabic{enumi})}
\setcounter{enumi}{1}
\tightlist
\item
  M. F. Atiyah, R. Bott, and A. Shapiro, ``Clifford modules'',
  \emph{Topology} \textbf{3} (1964), 3-38.
\end{enumerate}
\noindent
For more stuff on division algebras and Bott periodicity try Dave
Rusin's web page, especially his answer to ``Q5. What's the question
with the answer \(n = 1\), \(2\), \(4\), or \(8\)?''

\begin{enumerate}
\def\labelenumi{\arabic{enumi})}
\setcounter{enumi}{2}
\tightlist
\item
  Dave Rusin, ``Binary products, algebras, and division rings'',
 \href{https://web.archive.org/web/20150511070342/http://www.math.niu.edu/~rusin/known-math/95/division.alg}{\texttt{https://web.archive.}}
 \href{https://web.archive.org/web/20150511070342/http://www.math.niu.edu/~rusin/known-math/95/division.alg}{\texttt{org/web/20150511070342/http://www.math.niu.edu/\textasciitilde{}rusin/known-math/95/}}
\href{https://web.archive.org/web/20150511070342/http://www.math.niu.edu/~rusin/known-math/95/division.alg}{\texttt{division.alg}}
\end{enumerate}

Let me briefly explain this \(k\mathrm{O}\) business. The space
\(k\mathrm{O}\) is related to a simpler space called
\(\mathcal{B}\mathrm{O}(\infty)\) by means of the equation
\[k\mathrm{O} = \mathcal{B}\mathrm{O}(\infty)\times\mathbb{Z},\] so let
me first describe \(\mathcal{B}\mathrm{O}(\infty)\). For any topological
group \(G\) you can cook up a space \(BG\) whose loop space is homotopy
equivalent to \(G\). In other words, the space of
(base-point-preserving) maps from \(S^1\) to \(BG\) is homotopy
equivalent to \(G\). It follows that \[\pi_n(G) = \pi_{n+1}(BG).\] This
space \(BG\) is called the classifying space of \(G\) because it has a
principal \(G\)-bundle over it, and given \emph{any} decent topological
space \(X\) (say a CW complex) you can get all principal \(G\)-bundles
over \(X\) (up to isomorphism) by taking a map \(f\colon X\to BG\) and
pulling back this principal \(G\)-bundle over \(BG\). Moreover,
homotopic maps to \(BG\) give isomorphic \(G\)-bundles over \(X\) this
way.

Now a principal \(\mathrm{O}(n)\)-bundle is basically the same thing as
an \(n\)-dimensional real vector bundle --- there are obvious ways to go
back and forth between these concepts. A principal
\(\mathrm{O}(\infty)\)-bundle is thus very much like a real vector
bundle of \emph{arbitrary} dimension, but where we don't care about
adding on arbitrarily many \(1\)-dimensional trivial bundles. If we take
the collection of isomorphism classes of real vector bundles over \(X\)
and decree two to be equivalent if they become isomorphic after adding
on trivial bundles, we get something called \(KX\), the ``real K-theory
of X''. It's not hard to see that this is a group. Taking what I've said
and working a bit, it follows that
\[KX = [X, \mathcal{B}\mathrm{O}(\infty)]\] where the right-hand side
means ``homotopy classes of maps from \(X\) to
\(\mathcal{B}\mathrm{O}(\infty)\)''. If we take \(X\) to be \(S^{n+1}\),
we see
\[KS^{n+1} = \pi_{n+1}(\mathcal{B}\mathrm{O}(\infty)) = \pi_n(\mathrm{O}(\infty))\]
It follows that we can get all elements of \(\pi_n\) of
\(\mathrm{O}(\infty)\) from real vector bundles over \(S^{n+1}\).

Of course, the above equations are true only for nonnegative \(n\),
since it doesn't make sense to talk about \(\pi_{-1}\) of a space.
However, to make Bott periodicity work out smoothly, it would be nice if
we could pretend that
\[KS^{-1} = \pi_0(\mathcal{B}\mathrm{O}(\infty)) = \pi_{-1}(\mathrm{O}(\infty)) = \pi_7(\mathrm{O}(\infty)) = \mathbb{Z}\]
Alas, the equations don't make sense, and
\(\mathcal{B}\mathrm{O}(\infty)\) is connected, so we don't have
\(\pi_0(\mathcal{B}\mathrm{O}(\infty)) = \mathbb{Z}\). However, we can
cook up a slightly improved space \(k\mathrm{O}\), which has
\[\pi_n(k\mathrm{O}) = \pi_n(\mathcal{B}\mathrm{O}(\infty))\] when
\(n > 0\), but also has \[\pi_0(k\mathrm{O}) = \mathbb{Z}\] as desired.
It's easy --- we just let
\[k\mathrm{O} = \mathcal{B}\mathrm{O}(\infty)\times\mathbb{Z}.\] So,
let's use this instead of \(\mathcal{B}\mathrm{O}(\infty)\) from now on.

Taking \(n = 0\), we can think of \(S^1\) as \(\mathbb{RP}^1\), the real
projective line, i.e.~the space of \(1\)-dimensional real subspaces of
\(\mathbb{R}^2\). This has a ``canonical line bundle'' over it, that is,
a \(1\)-dimensional real vector bundle which to each point of
\(\mathbb{RP}^1\) assigns the \(1\)-dimensional subspace of
\(\mathbb{R}^2\) that \emph{is} that point. This vector bundle over
\(S^1\) gives the generator of \(KS^1\), or in other words,
\(\pi_0(\mathrm{O}(\infty))\).

Taking \(n = 1\), we can think of \(S^2\) as the ``Riemann sphere'', or
in other words \(\mathbb{CP}^1\), the space of \(1\)-dimensional complex
subspaces of \(\mathbb{C}^2\). This too has a ``canonical line bundle''
over it, which is a 1-dimensional complex vector bundle, or
\(2\)-dimensional real vector bundle. This bundle over \(S^2\) gives the
generator of \(KS^2\), or in other words, \(\pi_1(\mathrm{O}(\infty))\).

Taking \(n = 3\), we can think of \(S^4\) as \(\mathbb{HP}^1\), the
space of \(1\)-dimensional quaternionic subspaces of \(\mathbb{H}^2\).
The ``canonical line bundle'' over this gives the generator of \(KS^4\),
or in other words, \(\pi_3(\mathrm{O}(\infty))\).

Taking \(n = 7\), we can think of \(S^8\) as \(\mathbb{OP}^1\), the
space of \(1\)-dimensional octonionic subspaces of \(\mathbb{O}^2\). The
``canonical line bundle'' over this gives the generator of \(KS^8\), or
in other words, \(\pi_7(\mathrm{O}(\infty))\).

By Bott periodicity,
\[\pi_7(\mathrm{O}(\infty)) = \pi_8(k\mathrm{O}) = \pi_0(k\mathrm{O})\]
so the canonical line bundle over \(\mathbb{OP}^1\) also defines an
element of \(\pi_0(k\mathrm{O})\). But
\[\pi_0(k\mathrm{O}) = [S^0,k\mathrm{O}] = KS^0\] and \(KS^0\) simply
records the \emph{difference in dimension} between the two fibers of a
vector bundle over \(S^0\), which can be any integer. This is why the
octonions are related to dimension.

If for any pointed space we define \[K^n(X) = K(S^n\wedge X)\] we get a
cohomology theory called K-theory, and it turns out that
\[K^{n+8}(X) = K(X)\] which is another way of stating Bott periodicity.
Now if \(\{*\}\) denotes a single point, \(K(\{*\})\) is a ring (this is
quite common for cohomology theories), and it is generated by elements
of degrees 1, 2, 4, and 8. The generator of degree 8 is just the
canonical line bundle over \(\mathbb{OP}^1\) and multiplication by this
generator gives a map \[K^n(\{*\})\to K^{n+8}(\{*\})\] which is an
isomorphism of groups --- namely, Bott periodicity! In this sense the
octonions are responsible for Bott periodicity.

\begin{center}\rule{0.5\linewidth}{0.5pt}\end{center}

\textbf{Addendum}: The Clifford algebra clock is even better than I
described above, because it lets you work out the fancier Clifford
algebras \(C_{p,q}\), which are generated by \(p\) square roots of
\(-1\) and \(q\) square roots of \(1\), which all anticommute with each
other. These Clifford algebras are good when you have \(p\) dimensions
of ``space'' and \(q\) dimensions of ``time'', and I described the
physically important case where \(q = 1\) in
\protect\hyperlink{week93}{``Week 93''}. To figure them out, you just
work out \(p - q \mod 8\), look at what the clock says for that hour,
and then take \(N\times N\) matrices of what you see, with \(N\) chosen
so that \(C_{p,q}\) gets the right dimension, namely \(2^{p+q}\). So say
you're a string theorist and you think there are 9 space dimensions and
1 time dimension. You say: ``Okay, \(9-1 = 8\), so I look and see what's
at 8 o'clock. Okay, that's \(\mathbb{R}\), the real numbers. But my
Clifford algebra \(C_{9,1}\) is supposed to have dimension
\(2^{9+1}=1024=32^2\), so my Clifford algebra must consist of
\(32\times32\) \emph{matrices} with real entries.''

By the way, it's not so easy to see that the canonical line bundle over
\(\mathbb{OP}^1\) is the generator of \(KS^8\) --- or equivalently, that
left multiplication by unit octonions defines a map from \(S^7\) into
\(\mathrm{SO}(8)\) corresponding to the generator of
\(\pi_7(\mathrm{O}(\infty))\). I claimed it's true above, but when
someone asked me why this was true, I realized I couldn't prove it! That
made me nervous. So I asked on \texttt{sci.math.research} if it was
really true, and I got this reply:

\begin{quote}
From: Linus Kramer  \hfill \break
Newsgroups: sci.math.research  \hfill \break
Subject: \(\pi_7(O)\) and octonions  \hfill \break
Date: Tue, 09 Nov 1999 12:44:33 +0100

John Baez asked if \(\pi_7(O)\) is generated by the (multiplication by)
unit octonions.

View this as a question in KO-theory: the claim is that \(H^8\)
generates the reduced real K-theory \(\tilde{K}\mathrm{O}(S^8)\) of the
8-sphere; the bundle \(H^8\) over \(S^8\) is obtained by the standard
glueing process along the equator \(S^7\), using the octonion
multiplication. So \(H^8\) is the octonion Hopf bundle. Its Thom space
is the projective Cayley plane \(\mathbb{OP}^2\). Using this and
Hirzebruch's signature theorem, one sees that the Pontrjagin class of
\(H^8\) is \(p_8(H^8)=6x\), for a generator \(x\) of the
\(8\)-dimensional integral cohomology of \(S^8\) {[}a reference for this
calulation is my paper ``The topology of smooth projective planes'',
\emph{Arch.\ Math.\ }\textbf{63} (1994){]}. We have a diagram
\[K\mathrm{O}(S^8) \xrightarrow{\mathrm{cplx}} K(S^8) \xrightarrow{\mathrm{ch}} H(S^8)\]
the left arrow is complexification, the second arrow is the Chern
character. In dimension 8, these maps form an isomorphism. Now
\(\mathrm{ch}(\mathrm{cplx}(H^8))=8+x\) (see the formula in the last
paragraph in Husemoller's \emph{Fibre Bundles}, the chapter on ``Bott
periodicity and integrality theorems''. The constant factor is
unimportant, so the answer is yes, \(\pi_7(O)\) is generated by the map
\(S^7\to\mathbb{O}\) which sends a unit octonion \(A\) to the map
\(l_A\colon x\mapsto Ax\) in \(\mathrm{SO}(8)\).

Linus Kramer
\end{quote}

More recently I got an email from Todd Trimble which cites another
reference to this fact:

\begin{quote}
From: Todd Trimble \hfill \break
Subject: Hopf bundles \hfill \break
To: John Baez \hfill \break 
Date: Fri, 25 Mar 2005 16:37:11 -0500 \hfill \break

John,

In the book \emph{Numbers} (GTM \textbf{123}), there is an article by
Hirzebruch where the Bott periodicity result is formulated as saying
that the generators of \(\tilde{K}\mathrm{O}(S^n)\) in the cases
\(n = 1, 2, 4, 8\) are given by \([\eta]-1\) where \(\eta\) is the Hopf
bundle corresponding to \(\mathbb{R}\), \(\mathbb{C}\), \(\mathbb{H}\),
\(\mathbb{O}\) and 1 is the trivial line bundle over these scalar
``fields'' (of real dimension 1, 2, 4, 8), and is 0 for
\(n = 3, 5, 6, 7\) {[}p.~294{]}. Also that the Bott periodicity
isomorphism\\
\[\tilde{K}\mathrm{O}(S^n) \to \tilde{K}\mathrm{O}(S^{n+8})\] is induced
by \([\eta(\mathbb{O})]-1\) {[}p.~295{]}. I know you are aware of this
already (courtesy of the response of Linus Kramers to your
sci.math.research query), but I thought you might find a
published reference, on the authority of no less than Hirzebruch,
handier (should you need it) than referring to a
sci.math.research exchange.

Unfortunately no proof is given. Hirzebruch says (p.~295),

\begin{quote}
Remark. Our formulation of the Bott periodicity theorem will be found,
in essentials, in {[}reference to Bott's Lectures on \(K(X)\), without
proofs{]}. A detailed proof within the framework of K-theory is given in
the textbook {[}reference to Karoubi's K-theory{]}. The reader will have
a certain amount of difficulty, however, in extracting the results used
here from Karoubi's formulation.
\end{quote}

Todd
\end{quote}

\begin{center}\rule{0.5\linewidth}{0.5pt}\end{center}

\begin{quote}
\emph{\ldots{} for geometry, you know, is the gate of science, and the
gate is so low and small that one can only enter it as a little child.}

--- William Clifford
\end{quote}

\hypertarget{week106}{%
\section{July 23, 1997}\label{week106}}

Well, it seems I want to talk one more time about octonions before
moving on to other stuff. I'm a bit afraid this obsession with octonions
will mislead the nonexperts, fooling them into thinking octonions are
more central to mainstream mathematical physics than they actually are.
I'm also worried that the experts will think I'm spend all my time
brooding about octonions when I should be working on practical stuff
like quantum gravity. But darn it, this is summer vacation! The only way
I'm going to keep on cranking out ``This Week's Finds'' is if I write
about whatever I feel like, no matter how frivolous. So here goes.

First of all, let's make sure everyone here knows what projective space
is. If you don't, I'd better explain it. This is honest mainstream stuff
that everyone should know, good nutritious mathematics, so I won't need
to feel too guilty about serving the extravagant octonionic dessert
which follows.

Start with \(\mathbb{R}^n\), good old \(n\)-dimensional Euclidean space.
We can imagine wanting to ``compactify'' this so that if you go sailing
off to infinity in some direction you'll come sailing back from the
other side like Magellan. There are different ways to do this. A
well-known one is to take \(\mathbb{R}^n\) and add on one extra ``point
at infinity'', obtaining the \(n\)-dimensional sphere \(S^n\). Here the
idea is that start anywhere in \(\mathbb{R}^n\) and start sailing in any
direction, you are sailing towards this ``point at infinity''.

But there is a sneakier way to compactify \(\mathbb{R}^n\), which gives
us not the \(n\)-dimensional sphere but ``projective \(n\)-space''. Here
we add on a lot of points, one for each line through the origin. Now
there are \emph{lots} of points at infinity, one for every direction!
The idea here is that if you start at the origin and start sailing along
any straight line, you are sailing towards the point at infinity
corresponding to that line. Sailing along any parallel line takes you
twoards the same point at infinity. It's a bit like a perspective
picture where different families of parallel lines converge to different
points on the horizon --- the points on the horizon being points at
infinity.

Projective \(n\)-space is also called \(\mathbb{RP}^n\). The
\(\mathbb{R}\) is for ``real'', since this is actually ``real projective
\(n\)-space''. Later we'll see what happens if we replace the real
numbers by the complex numbers, quaternions, or octonions.

There are some other ways to think about \(\mathbb{RP}^n\) that are
useful either for visualizing it or doing calculations. First a nice way
to visualize it. First take \(\mathbb{R}^n\) and squash it down so it's
just the ball of radius \(1\), or more precisely, the ``open ball''
consisting of all vectors of length less than \(1\). We can do this
using a coordinate transformation like:
\[x \mapsto x' =  \frac{x}{\sqrt{1+|x|^2}}\] Here \(x\) stands for a
vector in \(\mathbb{R}^n\) and \(|x|\) is its length. Dividing the
vector \(x\) by \(\sqrt{1 + |x|^2}\) gives us a vector \(x'\) whose
length never quite gets to \(1\), though it can get as close at it
likes. So we have squashed \(\mathbb{R}^n\) down to the open ball of
radius \(1\).

Now say you start at the origin in this squashed version of
\(\mathbb{R}^n\) and sail off in any direction in a straight line. Then
you are secretly heading towards the boundary of the open ball. So the
points an the boundary of the open ball are like ``points at infinity''.

We can now compactify \(\mathbb{R}^n\) by including these points at
infinity. In other words, we can work not with the open ball but with
the ``closed ball'' consisting of all vectors \(x\)' whose length is
less than or equal to \(1\).

However, to get projective \(n\)-space we also have to decree that
antipodal points \(x'\) and \(-x'\) with \(|x'| = 1\) are to be regarded
as the same. In other words, we need to ``identify each point on the
boundary of the closed ball with its antipodal point''. The reason is
that we said that when you sail off to infinity along a particular
straight line, you are approaching a particular point in projective
\(n\)-space. Implicit in this is that it doesn't matter which \emph{way}
you sail along that straight line. Either direction takes you towards
the same point in projective \(n\)-space!

This may seem weird: in this world, when the cowboy says ``he went
thataway'' and points at a particular point on the horizon, you gotta
remember that his finger points both ways, and the villian could equally
well have gone in the opposite direction. The reason this is good is
that it makes projective space into a kind of geometer's paradise: any
two lines in projective space intersect in a \emph{single} point. No
more annoying exceptions: even ``parallel'' lines intersect in a single
point, which just happens to be a point at infinity. This simplifies
life enormously.

Okay, so \(\mathbb{RP}^n\) is the space formed by taking a closed
\(n\)-dimensional ball and identifying pairs of antipodal points on its
boundary.

A more abstract way to think of \(\mathbb{RP}^n\), which is incredibly
useful in computations, is as the set of all lines through the origin in
\(\mathbb{R}^{n+1}\). Why is this the same thing? Well, let me
illustrate it in an example. What's the space of lines through the
origin in \(\mathbb{R}^3\)? To keep track of these lines, draw a sphere
around the origin. Each line through the origin intersects this sphere
in two points. Either one point is in the northern hemisphere and the
other is in the southern hemisphere, or both are on the equator. So we
can keep track of all our lines using points on the northern hemisphere
and the equator, but identifying antipodal points on the equator. This
is just the same as taking the closed 2-dimensional ball and identifying
antipodal points on the boundary! QED. The same argument works in higher
dimensions too.

Now that we know a point in \(\mathbb{RP}^n\) is just a line through the
origin in \(\mathbb{R}^{n+1}\), it's easy to put coordinates on
\(\mathbb{RP}^n\). There's one line through the origin passing through
any point in \(\mathbb{R}^{n+1}\), but if we multiply the coordinates
\((x_1,\ldots,x_{n+1})\) of this point by any nonzero number we get the
same line. Thus we can use a list of \(n+1\) real numbers to describe a
point in \(\mathbb{RP}^n\), with the proviso that we get the same point
in \(\mathbb{RP}^n\) if someone comes along and multiplies them all by
some nonzero number! These are called ``homogeneous coordinates''.

If you don't like the ambiguity of homogeneous coordinates, you can go
right ahead and divide all the coordinates by the real number \(x_1\),
getting \[(1, x_2/x_1, \ldots , x_{n+1}/x_1)\] which lets us describe a
point in \(\mathbb{RP}^n\) by n real numbers, as befits an
\(n\)-dimensional real manifold. Of course, this won't work if \(x_1\)
happens to be zero! But we can divide by \(x_2\) if \(x_2\) is nonzero,
and so on. \emph{One} of them has to be nonzero, so we can cover
\(\mathbb{RP}^n\) with \(n+1\) different coordinate patches
corresponding to the regions where different \(x_i\)'s are nonzero. It's
easy to change coordinates, too.

This makes everything very algebraic, which makes it easy to generalize
\(\mathbb{RP}^n\) by replacing the real numbers with other number
systems. For example, to define ``complex projective \(n\)-space'' or
\(\mathbb{CP}^n\), just replace the word ``real'' by the word
``complex'' in the last two paragraphs, and replace ``\(\mathbb{R}\)''
by ``\(\mathbb{C}\)''. \(\mathbb{CP}^n\) is even more of a geometer's
paradise than \(\mathbb{RP}^n\), because when you work with complex
numbers you can solve all polynomial equations. Also, now there's no big
difference between an ellipse and a hyperbola! This sort of thing is why
\(\mathbb{CP}^n\) is so widely used as a context for ``algebraic
geometry''.

We can go even further and replace the real numbers by the quaternions,
\(\mathbb{H}\), defining the ``quaternionic projective \(n\)-space''
\(\mathbb{HP}^n\). If we are careful about writing things in the right
order, it's no problem that the quaternions are noncommutative\ldots{}
we can still divide by any nonzero quaternion, so we can cover
\(\mathbb{HP}^n\) with n+1 different coordinate charts and freely change
coordinates as desired.

We can try to go even further and use the octonions, O. Can we define
``octonionic projective \(n\)-space'', \(\mathbb{OP}^n\)? Well, now
things get tricky! Remember, the octonions are nonassociative. There's
no problem defining \[\mathbb{OP}^1\]; we can cover it with two
coordinate charts, corresponding to homogeneous coordinates of the form
\[(x, 1)\] and \[(1, y),\] and we can change coordinates back and forth
with no problem. This amounts to taking \(\mathbb{O}\) and adding a
single point at infinity, getting the 8-dimensional sphere \(S^8\). This
is part of a pattern:

\begin{itemize}
\tightlist
\item
  \(\mathbb{RP}^1 = S^1\)
\item
  \(\mathbb{CP}^1 = S^2\)
\item
  \(\mathbb{HP}^1 = S^4\)
\item
  \(\mathbb{OP}^1 = S^8\)
\end{itemize}

I discussed the implications of this pattern for Bott periodicity in
\protect\hyperlink{week105}{``Week 105''}.

We can also define \(\mathbb{OP}^2\). Here we have 3 coordinate charts
corresponding to homogeneous coordinates of the form \[(1, y, z),\]
\[(x, 1, z),\] and \[(x, y, 1).\] We can change back and forth between
coordinate systems, but now we have to \emph{check} that if we start
with the first coordinate system, change to the second coordinate
system, and then change back to the first, we wind up where we started!
This is not obvious, since multiplication is not associative. But it
works, thanks to a couple of identities that are not automatic in the
nonassociative context, but hold for the octonions:
\[(xy)^{-1} = y^{-1} x^{-1}\] and \[(xy)y^{-1} = x.\] Checking these
equations is a good exercise for anyone who wants to understand the
octonions.

Now for the cool part: \(\mathbb{OP}^2\) is where it ends!

We can't define \(\mathbb{OP}^n\) for \(n\) greater than 2, because the
nonassociativity keeps us from being able to change coordinates a bunch
of times and get back where we started! You might hope that we could
weasel out of this somehow, but it seems that there is a real sense in
which the higher-dimensional octonionic projective spaces don't exist.

So we have a fascinating situation: an infinite tower of
\(\mathbb{RP}\)\^{}n's, an infinite tower of \(\mathbb{CP}^n\)'s, an
infinite tower of \(\mathbb{HP}^n\)'s, but an abortive tower of
\(\mathbb{OP}^n\)'s going only up to \(n = 2\) and then fizzling out.
This means that while all sorts of geometry and group theory relating to
the reals, complexes and quaternions fits into infinite systematic
patterns, the geometry and group theory relating to the octonions is
quirky and mysterious.

We often associate mathematics with ``classical'' beauty, patterns
continuing ad infinitum with the ineluctable logic of a composition by
some divine Bach. But when we study \(\mathbb{OP}^2\) and its
implications, we see that mathematics also has room for ``exceptional''
beauty, beauty that flares into being and quickly subsides into silence
like a piece by Webern. Are the fundamental laws of physics based on
``classical'' mathematics or ``exceptional'' mathematics? Since our
universe seems unique and special --- don't ask me how would we know if
it weren't --- Witten has suggested the latter. Indeed, it crops up a
lot in string theory. This is why I'm trying to learn about the
octonions: a lot of exceptional objects in mathematics are tied to them.

I already discussed this a bit in \protect\hyperlink{week64}{``Week
64''}, where I sketched how there are 3 infinite sequences of
``classical'' simple Lie groups corresponding to rotations in
\(\mathbb{R}^n\), \(\mathbb{C}^n\), and \(\mathbb{H}^n\), and 5
``exceptional'' simple Lie groups related to the octonions. After
studying it all a bit more, I can now go into some more detail.

In order of increasing dimension, the 5 exceptional Lie groups are
called \(\mathrm{G}_2\), \(\mathrm{F}_4\), \(\mathrm{E}_6\),
\(\mathrm{E}_7\), and \(\mathrm{E}_8\). The smallest, \(\mathrm{G}_2\),
is easy to understand in terms of the octonions: it's just the group of
symmetries of the octonions as an algebra. It's a marvelous fact that
all the bigger ones are related to \(\mathbb{OP}^2\). This was
discovered by people like Freudenthal and Tits and Vinberg, but a great
place to read about it is the following fascinating book:

\begin{enumerate}
\def\labelenumi{\arabic{enumi})}
\tightlist
\item
  Boris Rosenfeld, \emph{Geometry of Lie Groups}, Kluwer, Dordrecht,  
   1997.
\end{enumerate}

The space \(\mathbb{OP}^2\) has a natural metric on it, which allows us
to measure distances between points. This allows us to define a certain
symmetry group \(\mathbb{OP}^2\), the group of all its ``isometries'',
which are transformations preserving the metric. This symmetry group is
\(\mathrm{F}_4\)!

However, there is another bigger symmetry group of \(\mathbb{OP}^2\). As
in real projective \(n\)-space, the notion of a ``line'' makes sense in
\(\mathbb{OP}^2\). One has to be careful: these are octonionic
``lines'', which have 8 real dimensions. Nonetheless, this lets us
define the group of all ``collineations'' of \(\mathbb{OP}^2\), that is,
transformations that take lines to lines. This symmetry group is
\(\mathrm{E}_6\)! (Technically speaking, this is a ``noncompact real
form'' of \(\mathrm{E}_6\); the rest of the time I'll be talking about
compact real forms.)

To get up to \(\mathrm{E}_7\) and \(\mathrm{E}_8\), we need to take a
different viewpoint, which also gives us another way to get
\(\mathrm{E}_6\). The key here is that the tensor product of two
algebras is an algebra, so we can tensor the octonions with
\(\mathbb{R}\), \(\mathbb{C}\), \(\mathbb{H}\), or \(\mathbb{O}\) and
get various algebras:

\begin{itemize}
\tightlist
\item
  The algebra \((\mathbb{R}\otimes\mathbb{O})\) is just the octonions.
\item
  The algebra \((\mathbb{C}\otimes\mathbb{O})\) is called the
  ``bioctonions''.
\item
  The algebra \((\mathbb{H}\otimes\mathbb{O})\) is called the
  ``quateroctonions''.
\item
  The algebra \((\mathbb{O}\otimes\mathbb{O})\) is called the
  ``octooctonions''.
\end{itemize}

I'm not making this up: it's all in Rosenfeld's book! The poet Lisa
Raphals suggested calling the octooctonions the ``high-octane
octonions'', which I sort of like. But compared to Rosenfeld, I'm a
model of restraint: I won't even mention the dyoctonions, duoctonions,
split octonions, semioctonions, split semioctonions, \(1/4\)-octonions
or \(1/8\)-octonions --- for the definitions of these, you'll have to
read his book.

Apparently one can define projective planes for all of these algebras,
and all these projective planes have natural metrics on them, all of
them same general form. So each of these projective planes has a group
of isometries. And, lo and behold:

\begin{itemize}
\tightlist
\item
  The group of isometries of the octonionic projective plane is
  \(\mathrm{F}_4\).
\item
  The group of isometries of the bioctonionic projective plane is
  \(\mathrm{E}_6\).
\item
  The group of isometries of the quateroctonionic projective plane is
  \(\mathrm{E}_7\).
\item
  The group of isometries of the octooctonionic projective plane is
  \(\mathrm{E}_8\).
\end{itemize}

Now I still don't understand this as well as I'd like to --- I'm not
sure how to define projective planes for all these algebras (though I
have my guesses), and Rosenfeld is unfortunately a tad reticent on this
issue. But it looks like a cool way to systematize the study of the
expectional groups! That's what I want: a systematic theory of
exceptions.

I want to say a bit more about the above, but first let me note that
there are lots of other ways of thinking about the exceptional groups. A
great source of information about them is the following posthumously
published book by the great topologist Adams:

\begin{enumerate}
\def\labelenumi{\arabic{enumi})}
\setcounter{enumi}{1}
\tightlist
\item
  John Frank Adams, \emph{Lectures on Exceptional Lie Groups},
  eds.~Zafer Mahmud and Mamoru Mimura, U.\ Chicago Press,
  Chicago, 1996.
\end{enumerate}

He has a bit about octonionic constructions of \(\mathrm{G}_2\) and
\(\mathrm{F}_4\), but mostly he concentrates on constructions of the
exceptional groups using classical groups and spinors.

In \protect\hyperlink{week90}{``Week 90''} I explained Kostant's
constructions of \(\mathrm{F}_4\) and \(\mathrm{E}_8\) using spinors in
8 dimensions and triality --- which, as noted in
\protect\hyperlink{week61}{``Week 61''}, is just another way of talking
about the octonions. Unfortunately I don't yet see quite how this
relates to the above stuff, nor do I see how to get \(\mathrm{E}_6\) and
\(\mathrm{E}_7\) in a beautiful way using Kostant's setup.

There's also a neat construction of \(\mathrm{E}_8\) using spinors in 16
dimensions! Adams gives a nice explanation of this, and it's also
discussed in the classic tome on string theory:

\begin{enumerate}
\def\labelenumi{\arabic{enumi})}
\setcounter{enumi}{2}
\tightlist
\item
  Michael B. Green, John H. Schwarz, and Edward Witten,
  \emph{Superstring Theory}, two volumes, Cambridge U.\ Press, Cambridge,
  1987.
\end{enumerate}

The idea here is to take the direct sum of the Lie algebra
\(\mathfrak{so}(16)\) and its \(16\)-dimensional left-handed spinor
representation \(S_+\) to get the Lie algebra of \(\mathrm{E}_8\). The
bracket of two guys in \(\mathfrak{so}(16)\) is defined as usual, the
bracket of a guy in \(\mathfrak{so}(16)\) and a guy in \(S_+\) is
defined to be the result of acting on the latter by the former, and the
bracket of two guys in \(S_+\) is defined to be a guy in \(S_+\) by
dualizing the map \[\mathfrak{so}(16) \otimes S_+\to S_+\] to get a map
\[S_+ \otimes S_+\to \mathfrak{so}(16).\] This is a complete description
of the Lie algebra of \(\mathrm{E}_8\)!

Anyway, there are lots of different ways of thinking about exceptional
groups, and a challenge for the octonionic approach is to systematize
all these ways.

Now I want to wrap up by saying a bit about how the exceptional Jordan
algebra fits into the above story. Jordan algebras were invented as a
way to study the self-adjoint operators on a Hilbert space, which
represent observables in quantum mechanics. If you multiply two
self-adjoint operators \(A\) and \(B\) the result needn't be
self-adjoint, but the ``Jordan product'' \[A \circ B = (AB + BA)/2\] is
self-adjoint. This suggests seeing what identities the Jordan product
satisfies, cooking up a general notion of ``Jordan algebra'', seeing how
much quantum mechanics you can do with an arbitrary Jordan algebra of
observables, and classifying Jordan algebras if possible.

We can define a ``projection'' in a Jordan algebra to be an element
\(A\) with \(A \circ A = A\). If our Jordan algebra consists of
self-adjoint operators on the complex Hilbert space \(\mathbb{C}^n\), a
projection is a self-adjoint operator whose only eigenvalues are zero
and one. Physically speaking, this corresponds to a ``yes-or-no
question'' about our quantum system. Geometrically speaking, such an
operator is a projection onto some subspace of our Hilbert space. All
this stuff also works if we start with the real Hilbert space
\(\mathbb{R}^n\) or the quaternionic Hilbert space \(\mathbb{H}^n\).

In these special cases, one can define a ``minimal projection'' to be a
projection on a \(1\)-dimensional subspace of our Hilbert space.
Physically, minimal projections correspond to ``pure states'' --- states
of affairs in which the answer to some maximally informative question is
``yes'', like ``is the \(z\) component of the angular momentum of this
spin-\(1/2\) particle equal to \(1/2\)?'' Geometrically, the space of
minimal projections is just the space of ``lines'' in our Hilbert space.
This is either \(\mathbb{RP}^{n-1}\), or \(\mathbb{CP}^{n-1}\), or
\(\mathbb{HP}^{n-1}\), depending on whether we're working with the
reals, complexes or quaternions. So: the space of pure states of this
sort of quantum system is also a projective space! The relation between
quantum theory and ``projective geometry'' has been intensively explored
for many years. You can read about it in:

\begin{enumerate}
\def\labelenumi{\arabic{enumi})}
\setcounter{enumi}{3}
\tightlist
\item
  V. S. Varadarajan, \emph{Geometry of Quantum Theory}, Springer,
  Berlin, 1985.
\end{enumerate}

Most people do quantum mechanics with complex Hilbert spaces. Real
Hilbert spaces are apparently too boring, but some people have
considered the quaternionic case:

\begin{enumerate}
\def\labelenumi{\arabic{enumi})}
\setcounter{enumi}{4}
\tightlist
\item
  Stephen L. Adler, \emph{Quaternionic Quantum Mechanics and Quantum
  Fields}, Oxford U. Press, Oxford, 1995.
\end{enumerate}

If our Hilbert space is the complex Hilbert space \(\mathbb{C}^n\), its
group of symmetries is usually thought of as \(\mathrm{U}(n)\) --- the
group of \(n\times n\) unitary matrices. This group also acts as
symmetries on the Jordan algebra of self-adjoint \(n\times n\) complex
matrices, and also on the space \(\mathbb{CP}^{n-1}\).

Similarly, if we start with \(\mathbb{R}^n\), we get the group of
orthogonal \(n\times n\) matrices \(\mathrm{O}(n)\), which acts on the
Jordan algebra of real self-adjoint \(n\times n\) matrices and on
\(\mathbb{RP}^{n-1}\).

Likewise, if we start with \(\mathbb{H}^n\), we get the group
\(\mathrm{Sp}(n)\), which acts on the Jordan algebra of quaternionic
self-adjoint \(n\times n\) matrices and on \(\mathbb{HP}^{n-1}\).

This pretty much explains how the classical groups are related to
different flavors of quantum mechanics.

Now what about the octonions? Well, here we can only go up to \(n = 3\),
basically for the reasons explained before: the same stuff that keeps us
from defining octonionic projective spaces past a certain point keeps us
from getting Jordan algebras! The interesting case is the Jordan algebra
of \(3\times3\) self-adjoint octonionic matrices. This is called the
``exceptional Jordan algebra'', \(J\). The group of symmetries of this
is --- you guessed it, \(\mathrm{F}_4\). One can also define a ``minimal
projection'' in \(J\) and the space of these is \(\mathbb{OP}^2\).

Is it possible that octonionic quantum mechanics plays some role in
physics?  I don't know.

Anyway, here is my hunch about the bioctonionic, quateroctonionic, and
octooctonionic projective planes. I think to define them you should
probably tensor the exceptional Jordan algebra with \(\mathbb{C}\),
\(\mathbb{H}\), and \(\mathbb{O}\), respectively, and take the space of
minimal projections in the resulting algebra. Rosenfeld seems to suggest
this is the way to go. However, I'm vague about some important details,
and it bugs me, because the special identities I needed above to define
\(\mathbb{OP}^2\) are related to \(\mathbb{O}\) being an alternative
algebra, but \(\mathbb{C}\otimes\mathbb{O}\),
\(\mathbb{H}\otimes\mathbb{O}\) and \(\mathbb{O}\otimes\mathbb{O}\) are
not alternative.

I should add that in addition to octonionic projective geometry, one can
do octonionic hyperbolic geometry. One can read about this in Rosenfeld
and also in the following:

\begin{enumerate}
\def\labelenumi{\arabic{enumi})}
\setcounter{enumi}{5}
\tightlist
\item
  Daniel Allcock, ``Reflection groups on the octave hyperbolic plane'',
  \emph{Jour. Algebra} \textbf{213} (1999), 467--498.
\end{enumerate}

\begin{center}\rule{0.5\linewidth}{0.5pt}\end{center}

\textbf{Addenda:} Here's an email from David Broadhurst, followed by
various remarks.

\begin{quote}
John:

Shortly before his death I spent a charming afternoon with Paul Dirac.
Contrary to his reputation, he was most forthcoming.

Among many things, I recall this: Dirac explained that while trained as
an engineer and known as a physicist, his aesthetics were mathematical.
He said (as I can best recall, nearly 20 years on): At a young age, I
fell in love with projective geometry. I always wanted to use to use it
in physics, but never found a place for it.

Then someone told him that the difference between complex and
quaternionic QM had been characterized as the failure of theorem in
classical projective geometry.

Dirac's face beamed a lovely smile: Ah he said, it was just such a thing
that I hoped to do.

I was reminded of this when bactracking to your
\protect\hyperlink{week106}{``Week 106''}, today.

Best, David
\end{quote}
\noindent
The theorem that fails for quaternions but holds for \(\mathbb{R}\) and
\(\mathbb{C}\) is the ``Pappus theorem'', discussed in
\protect\hyperlink{week145}{``Week 145''}.

Next, a bit about \(\mathbb{OP}^n\). There are different senses in which
we can't define \(\mathbb{OP}^n\) for \(n\) greater than 2. One is that
if we try to define coordinates on \(\mathbb{OP}^n\) in a similar way to
how we did it for \(\mathbb{OP}^2\), nonassociativity keeps us from
being able to change coordinates a bunch of times and get back where we
started! It's definitely enlightening to see how the desired transition
functions \(g_{ij}\) fail to satisfy the necessary cocycle condition
\(g_{ij} g_{jk} = g_{ik}\) when we get up to \(\mathbb{OP}^3\), which
would require 4 charts.

But, a deeper way to think about this emerged in conversations I've had
with James Dolan. Stasheff invented a notion of ``\(A_\infty\) space'',
which is a pointed topological space with a product that is associative
up to homotopy which satisfies the pentagon identity up to\ldots{} etc.
Any \(A_\infty\) space \(G\) has a classifying space \(BG\) such that
\[\Omega(BG) \simeq G.\] In other words, \(BG\) is a pointed space such
that the space of loops based at this point is homotopy equivalent to
\(G\). One can form this space \(BG\) by the Milnor construction:
sticking in one 0-simplex, one \(1\)-simplex for every point of \(G\),
one \(2\)-simplex for every triple \((g,h,k)\) with \(gh = k\), one
\(3\)-simplex for every associator, and so on. If we do this where \(G\)
is the group of length-one elements of \(\mathbb{R}\)
(i.e.~\(\mathbb{Z}/2\)) we get \(\mathbb{RP}^\infty\), as we expect,
since \[\mathbb{RP}^\infty = B(\mathbb{Z}/2).\] Even better, at the
\(n\)th stage of the Milnor construction we get a space homeomorphic to
\(\mathbb{RP}^n\). Similarly, if we do this where \(G\) is the group of
length-one elements of \(\mathbb{C}\) or \(\mathbb{H}\) we get
\(\mathbb{CP}^\infty\) or \(\mathbb{HP}^\infty\). But if we take \(G\)
to be the units of \(\mathbb{O}\), which has a product but is not even
homotopy-associative, we get \(\mathbb{OP}^1 = S^7\) at the first step,
\(\mathbb{OP}^2\) at the second step, \ldots{} but there's no way to
perform the third step!

Next: here's a little more information on the octonionic, bioctonionic,
quateroctonionic and octooctonionic projective planes. Rosenfeld claims
that the groups of isometries of these planes are \(\mathrm{F}_4\),
\(\mathrm{E}_6\), \(\mathrm{E}_7\), and \(\mathrm{E}_8\), respectively.
The problem is, I can't quite understand how he constructs these spaces,
except for the plain octonionic case.

It appears that these spaces can also be constructed using the ideas in
Adams' book. Here's how it goes.

\begin{itemize}
\tightlist
\item
  The Lie algebra \(\mathrm{F}_4\) has a subalgebra of maximal rank
  isomorphic to \(\mathfrak{so}(9)\). The quotient space is
  \(16\)-dimensional --- twice the dimension of the octonions. It
  follows that the Lie group \(\mathrm{F}_4\) mod the subgroup generated
  by this subalgebra is a \(16\)-dimensional Riemannian manifold on
  which \(\mathrm{F}_4\) acts by isometries.
\item
  The Lie algebra \(\mathrm{E}_6\) has a subalgebra of maximal rank
  isomorphic to \(\mathfrak{so}(10)\oplus\mathfrak{u}(1)\). The quotient
  space is \(32\)-dimensional --- twice the dimension of the
  bioctonions. It follows that the Lie group \(\mathrm{E}_6\) mod the
  subgroup generated by this subalgebra is a \(32\)-dimensional
  Riemannian manifold on which \(\mathrm{E}_6\) acts by isometries.
\item
  The Lie algebra \(\mathrm{E}_7\) has a subalgebra of maximal rank
  isomorphic to \(\mathfrak{so}(12)\oplus\mathfrak{su}(2)\). The
  quotient space is 64-dimensional --- twice the dimension of the
  quateroctonions. It follows that the Lie group \(\mathrm{E}_6\) mod
  the subgroup generated by this subalgebra is a 64-dimensional
  Riemannian manifold on which \(\mathrm{E}_7\) acts by isometries.
\item
  The Lie algebra \(\mathrm{E}_8\) has a subalgebra of maximal rank
  isomorphic to \(\mathfrak{so}(16)\). The quotient space is
  128-dimensional --- twice the dimension of the octooctonions. It
  follows that the Lie group \(\mathrm{E}_6\) mod the subgroup generated
  by this subalgebra is a 128-dimensional Riemannian manifold on which
  \(\mathrm{E}_8\) acts by isometries.
\end{itemize}

According to:

\begin{enumerate}
\def\labelenumi{\arabic{enumi})}
\setcounter{enumi}{5}
\tightlist
\item
  Arthur L. Besse, \emph{Einstein Manifolds}, Springer, Berlin, 1987,
  pp.~313--316.
\end{enumerate}
\noindent
the above spaces deserve to be called the octonionic, bioctonionic, quateroctonionic 
and octooctonionic projective planes, respectively. However, I don't fully understand
the connection.

I thank Tony Smith for pointing out the reference to Besse (who, by the
way, is apparently a cousin of the famous Bourbaki). Thanks also go to
Allen Knutson for showing me a trick for finding the maximal rank
subalgebras of a simple Lie algebra.

Next, here's some more stuff about the biquaternions, bioctonions,
quaterquaternions, quateroctonions and octooctonions! I wrote this extra
stuff as part of a post to \texttt{sci.physics.research} on November 8,
1999\ldots.

\begin{quote}
One reason people like these algebras is that some of them --- the
associative ones --- are also Clifford algebras. I talked a bit about
Clifford algebras in \protect\hyperlink{week105}{``Week 105''}, but just
remember that we define the Clifford algebra \(C_{p,q}\) to be the
associative algebra you get by taking the real numbers and throwing in
\(p\) square roots of \(-1\) and \(q\) square roots of \(1\), all of
which anticommute with each other. This algebra is very important for
understanding spinors in spacetimes with \(p\) space and \(q\) time
dimensions. (It's also good for studying things in other dimensions, so
things can get a bit tricky, but I don't want to talk about that now.)

For example: if you just thrown in one square root of \(-1\) and no
square roots of \(1\), you get \(C_{1,0}\) --- the complex numbers!

Similarly, one reason people like the quaternions is because they are
\(C_{2,0}\). Start with the real numbers, throw in two square roots of
\(-1\) called \(I\) and \(J\), make sure they anticommute (\(IJ = -JI\))
and voila --- you've got the quaternions!

Similarly, one reason people like the biquaternions is because they are
\(C_{2,1}\). You take the quaternions and complexify them --- this
amounts to throwing in an extra number \(i\) that's a square root of
\(-1\) and \emph{commutes} with the quaternionic \(I\) and \(J\) --- and
you get an algebra which is also generated by \(I\), \(J\), and
\(K = iI\). Note that \(I\), \(J\), and \(K\) all anticommute, and \(K\)
is a square root of \(1\). Thus the biquaternions are \(C_{2,1}\)!

Similarly, one reason people like the quaterquaternions is because they
are \(C_{2,2}\). You take the quaternions and quaternionify them ---
this amounts to throwing in two square roots of \(-1\), say \(i\) and
\(j\), which anticommute but which \emph{commute} with the quaternionic
\(I\) and \(J\) --- and you get an algebra which is also generated by
\(I\), \(J\), \(K = iI\), and \(L = jI\). Note that \(I\), \(J\), \(K\),
and \(L\) all anticommute, and \(K\) and \(L\) are square roots of
\(1\). Thus the quaterquaternions are \(C_{2,2}\)!

Now, as soon as we thrown the octonions into the mix we don't get
Clifford algebras anymore, since octonions aren't associative, while
Clifford algebras are. However, there are still relationships to
Clifford algebras. For example, suppose we look at all the linear
transformations of the octonions generated by the left multiplication
operations \[x \mapsto ax\] This is an associative algebra, and it turns
out to be \emph{all} linear transformations of the octonions, regarded
as an \(8\)-dimensional real vector space. In short, it's just the
algebra of \(8\times8\) real matrices. And this is \(C_{6,0}\).

If you do the same trick for the bioctonions, quateroctonions and
octooctonions, you get other Clifford algebras\ldots{} but I'll leave
the question of \emph{which ones} as a puzzle for the reader. If you
need some help, look at the ``Footnote'' in
\protect\hyperlink{week105}{``Week 105''}.

Perhaps the fanciest example of this trick concerns the
biquateroctonions. Now actually, I've never heard anyone use this term
for the algebra \(\mathbb{C}\otimes\mathbb{H}\otimes\mathbb{O}\)! The
main person interested in this algebra is Geoffrey Dixon, and he just
calls it \(T\). But anyway, if we look at the algebra of linear
transformations of \(\mathbb{C}\otimes\mathbb{H}\otimes\mathbb{O}\)
generated by left multiplications, we get something isomorphic to the
algebra of \(16\times16\) complex matrices. And this in turn is
isomorphic to \(C_{9,0}\).

The biquateroctonions play an important role in Dixon's grand unified
theory of the electromagnetic, weak and strong forces. There are lots of
nice things about this theory --- for example, it gets the right
relationships between weak isospin and hypercharge for the fermions in
any one generation of the Standard Model (though, as in the Standard
Model, the existence of 3 generations needs to be put in ``by hand'').
It may or may not be right, but at least it comes within shooting
distance!

You can read a bit more about his work in
\protect\hyperlink{week59}{``Week 59''}.
\end{quote}

\begin{center}\rule{0.5\linewidth}{0.5pt}\end{center}

\begin{quote}
\emph{``Mainstream mathematics'' is a name given to mathematics that
more fittingly belongs on Sunset Boulevard.}

--- Gian-Carlo Rota, \emph{Indiscrete Thoughts}
\end{quote}

\hypertarget{week107}{%
\section{August 19, 1997}\label{week107}}

This summer I've been hanging out in Cambridge Massachusetts, working on
quantum gravity and also having some fun. Not so long ago I gave a talk
on cellular automata at Boston University, thanks to a kind invitation
from Bruce Boghosian, who is using cellular automata to model cool stuff
like emulsions:

\begin{enumerate}
\def\labelenumi{\arabic{enumi})}
\tightlist
\item
  Florian W. J. Weig, Peter V. Coveney, and Bruce M. Boghosian,
  ``Lattice-gas simulations of minority-phase domain growth in binary
  immiscible and ternary amphiphilic fluid'', available as
  \href{https://arxiv.org/abs/cond-mat/9705248}{\texttt{cond-mat/9705248}}.
\end{enumerate}
\noindent
As you add more and more of an amphiphilic molecule (e.g.~soap) to a
binary immiscible fluid (e.g.~oil and water), the boundary layer likes
to grow in area. This is why you wash your hands with soap. There are
various phases depending on the concentrations of the three substances
--- a ``spongy'' phase, a ``droplet phase'', and so on --- and it is
very hard to figure out what is going on quantitatively using analytical
methods.

Luckily, one can simulate this stuff using a cellular automaton!
Standard numerical methods for solving the Navier--Stokes equation tend
to outrun cellular automata when it comes to plain old hydrodynamics,
but with these fancy ``ternary amphiphilic fluids'', cellular automata
really seem to be the most practical way to study things - apart from
experiments, of course. This is very heartwarming to me, since like many
people I've been fond of cellular automata ever after learning of John
Conway's game of Life, and I've always hoped they could serve some
practical purpose.

I spoke about the thesis of my student James Gilliam and a paper we
wrote together:

\begin{enumerate}
\def\labelenumi{\arabic{enumi})}
\setcounter{enumi}{1}
\item
  James Gilliam, \emph{Lagrangian and Symplectic Techniques in Discrete
  Mechanics}, Ph.D.~thesis, Department of Mathematics, University of
  Riverside, 1996.  Available at  \hfill \break
\href{http://math.ucr.edu/home/baez/thesis_gilliam.pdf}{\texttt{http://math.ucr.edu/home/baez/thesis$\underline{\;}$gilliam.pdf}}

  John Baez and James Gilliam, An algebraic approach to discrete
  mechanics, \emph{Lett. Math. Phys.} \textbf{31} (1994), 205--212. Also
  available at \href{http://math.ucr.edu./home/baez/ca.pdf}{\texttt{http://math.ucr.edu./home/}} \href{http://math.ucr.edu./home/baez/ca.pdf}{\texttt{baez/ca.pdf}}
\end{enumerate}
\noindent
Here the idea was to set up as much as possible of the machinery of
classical mechanics in a purely discrete context, where time proceeds in
integer steps and the space of states is also discrete. The most famous
examples of this ``discrete mechanics'' are cellular automata, which are
the discrete analogs of classical field theories, but there are also
simpler systems more reminiscent of elementary classical mechanics, like
a particle moving on a line --- where in this case the ``line'' is the
integers rather than the real numbers. It turns out that with a little
skullduggery one can apply the techniques of calculus to some of these
situations, and do all sorts of stuff like prove a discrete version of
Noether's theorem --- the famous theorem which gives conserved
quantities from symmetries.

After giving this talk, I visited my friend Robert Kotiuga in the
Functorial Electromagnetic Analysis Lab in the Photonics Building at
Boston University. ``Photonics'' is the currently fashionable term for
certain aspects of optics, particularly quantum optics. As befits its
flashy name, the Photonics Building is brand new and full of gadgets
like a device that displays Maxwell's equations in moving lights when
you speak the words ``Maxwell's equations'' into an inconspicuous
microphone. (It also knows other tricks.) Robert told me about what he's
been doing lately with topology and finite-element methods for solving
magnetostatics problems - this blend of higbrow math and practical
engineering being the reason for the somewhat tongue-in-cheek name of
his office, inscribed soberly on a plaque outside the door.

Like the topologist Raoul Bott, Kotiuga started in electrical
engineering at McGill University, and gradually realized how much
topology there is lurking in electrical circuit theory and Maxwell's
equations. Apparently a paper of his was the first to cite Witten's
famous work on Chern--Simons theory - though presumably this is merely a
testament to the superiority of engineers over mathematicians and
physicists when it comes to rapid publication. In fluid dynamics, the
integral of the following quantity 
\[v\cdot \nabla \times v\]
(where \(v\) is the velocity vector field) is known as the ``helicity
functional''. Kotiuga been studying applications of the same
mathematical object in the context of magnetostatics, namely
\[A\cdot \nabla \times A\] where \(A\) is the magnetic vector
potential. It shows up in impedance tomography, for example. But in
quantum field theory, a generalization of this quantity to other forces
is known as the ``Chern--Simons functional'', and Witten's work on the
\(3\)-dimensional field theory having this as its Lagrangian turned out
to revolutionize knot theory. Personally, I'm mainly interested in the
applications to quantum gravity --- see
\protect\hyperlink{week56}{``Week 56''} for a bit about this. Here are
some papers Kotiuga has written on the helicity functional, or what we
mathematicians would call ``\(\mathrm{U}(1)\) Chern--Simons theory'':

\begin{enumerate}
\def\labelenumi{\arabic{enumi})}
\setcounter{enumi}{2}
\item
  P. R. Kotiuga, ``Metric dependent aspects of inverse problems and
  functionals based helicity'', \emph{Journal of Applied Physics},
  \textbf{70} (1993), 5437--5439.

  ``Analysis of finite element matrices arising from discretizations of
  helicity functionals'', \emph{Journal of Applied Physics}, \textbf{67}
  (1990), 5815--5817.

  ``Helicity functionals and metric invariance in three dimensions'',
  \emph{IEEE Transactions on Magnetics}, MAG-\textbf{25} (1989),
  2813--2815.

  ``Variational principles for three-dimensional magnetostatics based on
  helicity'', \emph{Journal of Applied Physics}, \textbf{63} (1988),
  3360--3362.
\end{enumerate}

Later Jon Doyle, a computer scientist at M.I.T. who had been to my talk,
invited me to a seminar at M.I.T. where I met Gerald Sussman, who with
Jack Wisdom has run the best long-term simulations of the solar system,
trying to settle the old question of whether the darn thing is stable!
It turns out that the system is afflicted with chaos and can only be
predicted with any certainty for about 4 million years\ldots{} though
their simulation went out to 100 million.

Here are some fun facts:

\begin{itemize}
\def\labelenumi{\arabic{enumi})}
\tightlist
\item
  They need to take general relativity into account even for the orbit
  of Jupiter, which precesses about one radian per billion years.
\item
  They take the asteroid belt into account only as modification of the
  sun's quadrupole moment (which they also use to model its oblateness).
\item
  The most worrisome thing about the whole simulation --- the most
  complicated and unpredictable aspect of the whole solar system in
  terms of its gravitational effects on everything else --- is the
  Earth-Moon system, with its big tidal effects. 
\item The sun loses one
  Earth mass per 100 million years due to radiation, and another quarter
  Earth mass due to solar wind. 
\item The first planet to go is Mercury! In
  their simulations, it eventually picks up energy through a resonance
  and drifts away.
\end{itemize}

For more, try:

\begin{enumerate}
\def\labelenumi{\arabic{enumi})}
\setcounter{enumi}{3}
\item
  Gerald Jay Sussman and Jack Wisdom, ``Chaotic evolution of the solar
  system'', \emph{Science}, \textbf{257}, 3 July 1992.

  Gerald Jay Sussman and Jack Wisdom, ``Numerical evidence that the
  motion of Pluto is chaotic'', \emph{Science}, \textbf{241}, 22 July
  1988.

  James Applegate, M. Douglas, Y. Gursel, Gerald Jay Sussman, Jack
  Wisdom, ``The outer solar system for 200 million years'',
  \emph{Astronomical Journal}, \textbf{92}, pp 176--194, July 1986,
  reprinted in \emph{Use of
  Supercomputers in Stellar Dynamics}, Lecture Notes in Physics \textbf{267}
   Springer, Berlin, 1986.

  James Applegate, M. Douglas, Y. Gursel, P Hunter, C. Seitz, Gerald Jay
  Sussman, ``A digital orrery'', in \emph{IEEE Transactions on
  Computers}, C-\textbf{34}, No.~9, pp.~822--831, September 1985,
  reprinted in Lecture Notes in Physics \textbf{267}, Springer, Berlin,
  1986.
\end{enumerate}

Meanwhile, I've also been trying to keep up with recent developments in
\(n\)-category theory. Some readers of ``This Week's Finds'' have
expressed frustration with how I keep tantalizing all of you with the
concept of \(n\)-category without ever quite defining it. The reason is
that it's a lot of work to write a nice exposition of this concept!

However, I eventually got around to taking a shot at it, so now you can
read this:

\begin{enumerate}
\def\labelenumi{\arabic{enumi})}
\setcounter{enumi}{4}
\tightlist
\item
  John Baez, ``Introduction to \(n\)-categories'', in
  \emph{7th Conference on Category Theory and Computer Science}, eds.~Eugenio
  Moggi and Giuseppe Rosolini, Lecture Notes in Computer Science
  vol.~\textbf{1290}, Springer, Berlin. Also available as
  \href{https://arxiv.org/abs/q-alg/9705009}{\texttt{q-alg/9705009}}.
\end{enumerate}
\noindent
There are different definitions of ``weak \(n\)-category'' out there now
and it will take a while of sorting through them to show a bunch are
equivalent and get the whole machinery running smoothly. In the above
paper I mainly talk about the definition that James Dolan and I came up
with. Here are some other new papers on this sort of thing\ldots{} I'll
just list them with abstracts.

\begin{enumerate}
\def\labelenumi{\arabic{enumi})}
\setcounter{enumi}{5}
\tightlist
\item
  Claudio Hermida, Michael Makkai and John Power, ``On weak higher
  dimensional categories'',  \emph{Jour.\ Pure Appl.\ Alg.} \textbf{154} 
 (2000), 221--246.  Also available at  \hfill \break \url{https://ncatlab.org/nlab/files/HermidaMakkaiPower01.pdf}
\end{enumerate}

\begin{quote}
Inspired by the concept of opetopic set introduced in a recent paper by
John C. Baez and James Dolan, we give a modified notion called
multitopic set. The name reflects the fact that, whereas the Baez/Dolan
concept is based on operads, the one in this paper is based on
multicategories. The concept of multicategory used here is a mild
generalization of the same-named notion introduced by Joachim Lambek in
1969. Opetopic sets and multitopic sets are both intended as vehicles
for concepts of weak higher dimensional category. Baez and Dolan define
weak \(n\)-categories as \((n+1)\)-dimensional opetopic sets satisfying
certain properties. The version intended here, multitopic
\(n\)-category, is similarly related to multitopic sets. Multitopic
\(n\)-categories are not described in the present paper; they are to
follow in a sequel. The present paper gives complete details of the
definitions and basic properties of the concepts involved with
multitopic sets. The category of multitopes, analogs of the opetopes of
Baez and Dolan, is presented in full, and it is shown that the category
of multitopic sets is equivalent to the category of set- valued functors
on the category of multitopes.
\end{quote}

\begin{enumerate}
\def\labelenumi{\arabic{enumi})}
\setcounter{enumi}{6}
\tightlist
\item
  Michael Batanin, Finitary monads on globular sets and notions of
  computad they generate.   (Apparently no longer available.)
\end{enumerate}

\begin{quote}
Consider a finitary monad on the category of globular sets. We prove
that the category of its algebras is isomorphic to the category of
algebras of an appropriate monad on the special category (of computads)
constructed from the data of the initial monad. In the case of the free
\(n\)-category monad this definition coincides with \({R}\). Street's
definition of \(n\)-computad. In the case of a monad generated by a
higher operad this allows us to define a pasting operation in a weak
\(n\)-category. It may be also considered as the first step toward the
proof of equivalence of the different definitions of weak
\(n\)-categories.
\end{quote}

\begin{enumerate}
\def\labelenumi{\arabic{enumi})}
\setcounter{enumi}{6}
\tightlist
\item
  Carlos Simpson, Limits in \(n\)-categories, 
  available as
  \href{https://arxiv.org/abs/alg-geom/9708010}{\texttt{alg-geom/9708010}}.
\end{enumerate}

\begin{quote}
We define notions of direct and inverse limits in an \(n\)-category. We
prove that the \((n+1)\)-category \(n\mathsf{CAT}'\) of fibrant
\(n\)-categories admits direct and inverse limits. At the end we
speculate (without proofs) on some applications of the notion of limit,
including homotopy fiber product and homotopy coproduct for
\(n\)-categories, the notion of n-stack, representable functors, and
finally on a somewhat different note, a notion of relative Malcev
completion of the higher homotopy at a representation of the fundamental
group.
\end{quote}

\begin{enumerate}
\def\labelenumi{\arabic{enumi})}
\setcounter{enumi}{7}
\tightlist
\item
  Sjoerd Crans, ``Generalized centers of braided and sylleptic monoidal
  \(2\)-categories'', \emph{Adv. Math.} \textbf{136} (1998), 183--223.
\end{enumerate}

\begin{quote}
Recent developments in higher-dimensional algebra due to Kapranov and
Voevodsky, Day and Street, and Baez and Neuchl include definitions of
braided, sylleptic and symmetric monoidal \(2\)-categories, and a center
construction for monoidal \(2\)-categories which gives a braided
monoidal \(2\)-category. I give generalized center constructions for
braided and sylleptic monoidal \(2\)-categories which give sylleptic and
symmetric monoidal \(2\)-categories respectively, and I correct some
errors in the original center construction for monoidal
\(2\)-categories.
\end{quote}

\begin{center}\rule{0.5\linewidth}{0.5pt}\end{center}

\begin{quote}
\emph{Time definitely repeats itself: that's its only job. }

--- Edward Dorn, \emph{Sirius in January}
\end{quote}

\hypertarget{week108}{%
\section{September 22, 1997}\label{week108}}

In the Weeks to come I want to talk about quantum gravity, and
especially the relation between general relativity and spinors, since
Barrett and Crane and I have some new papers out about how you can
describe ``quantum 4-geometries'' --- geometries of spacetime which have
a kind of quantum discreteness at the Planck scale --- starting from the
mathematics of spinors.

But first I want to say a bit about CTCS '97 --- a conference on
category theory and computer science organized by Eugenio Moggi and
Giuseppe Rosolini. It was so well-organized that they handed us the
conference proceedings when we arrived:

\begin{enumerate}
\def\labelenumi{\arabic{enumi})}
\tightlist
\item
  Eugenio Moggi and Giuseppe Rosolini, eds., \emph{Category Theory and
  Computer Science}, Lecture Notes in Computer Science \textbf{1290},
  Springer, Berlin, 1997.
\end{enumerate}

It was held in Santa Margherita Ligure, a picturesque little Italian
beach town near Genoa - the perfect place to spend all day in the
basement of a hotel listening to highly technical lectures. It's near
Portofino, famous for its big yachts full of rich tourists, but I didn't
get that far. The food was great, though, and it was nice to see the
lazy waves of the Mediterranean, so different from the oceans I know and
love. The vegetation was surprisingly similar to that in Riverside: lots
of palm trees and cacti.

I spoke about \(n\)-categories, with only the barest mention of their
possible relevance to computer science. But I was just the token
mathematical physicist in the crowd; most of the other participants were
pretty heavily into ``theoretical computer science'' --- a subject that
covers a lot of new-fangled aspects of what used to be called ``logic''.
What's neat is that I almost understood some of these talks, thanks to
the fact that category theory provides a highly general language for
talking about processes.

What's a computer, after all, but a physical process that simulates
fairly arbitrary processes --- including other physical processes? As we
simulate more and more physics with better and better computers based on
more and more physics, it seems almost inevitable that physics and
computer science will come to be seen as two ends of a more general
theory of processes. No?

A nice example of an analogy between theoretical computer science and
mathematical physics was provided by Gordon Plotkin (in the plane, on
the way back, when I forced him to explain his talk to me). Computer
scientists like to define functions recursively. For example, we can
define a function from the natural numbers to the natural numbers:
\[f\colon\mathbb{N}\to\mathbb{N}\] by its value at \(0\) together with a
rule to get its value at \(n+1\) from its value at \(n\): \[
  \begin{aligned}
    f(0) &= c
  \\f(n+1) &= g(f(n))
  \end{aligned}
\]

Similarly, physicists like to define functions by differential
equations. For example, we can define a function from the real numbers
to the real numbers: \[f\colon\mathbb{R}\to\mathbb{R}\] by its value at
\(0\) together with a rule to get its derivative from its value: \[
  \begin{aligned}
    f(0) &= c
  \\f'(t) &= g(f(t))
  \end{aligned}
\]

In both cases a question arises: how do we know we've really defined a
function? And in both cases, the answer involves a ``fixed-point
theorem''. In both cases, the equations above define the function \(f\)
\emph{in terms of itself}. We can write this using an equation of the
form: \[f = F(f)\] where \(F\) is some operator that takes functions to
functions. We say \(f\) is a ``fixed point for \(F\)'' if this holds. A
fixed-point theorem is something that says there exists a solution,
preferably unique, of this sort of equation.

But how do we describe this operator \(F\) more precisely in these
examples? In the case of the definition by recursion, here's how: for
any function \(f\colon\mathbb{N}\to\mathbb{N}\), we define the function
\(F(f)\colon\mathbb{N}\to\mathbb{N}\) by \[
  \begin{aligned}
    F(f)(0) &= c
  \\F(f)(n+1) &= g(f(n))
  \end{aligned}
\] The principle of mathematical induction says that any operator \(F\)
of this sort has a unique fixed point.

Similarly, we can formulate the differential equation above as a fixed
point problem by integrating both sides, obtaining:
\[f(t) = c + \int_0^t g(f(s)) ds\] which is an example of an ``integral
equation''. If we call the function on the right hand side \(F(f)\),
then this integral equation says \[f = F(f)\] In this case, ``Picard's
theorem on the local existence and uniqueness of solutions of ordinary
differential equations'' is what comes to our rescue and asserts the
existence of a unique fixed point.

You might wonder how Picard's theorem is proved. The basic idea of the
proof is very beautiful, because it \emph{takes advantage} of the
frightening circularity implicit in the equation \(f = F(f)\). I'll
sketch this idea, leaving out all the details.

So, how do we solve this equation? Let's see what we can do with it.
There's not much to do, actually, except substitute the left side into
the right side and get: \[f = F(F(f)).\] Hmm. Now what? Well, we can do
it again: \[f = F(F(F(f)))\] and again: \[f = F(F(F(f)))).\] Are we
having fun yet? It look like we're getting nowhere fast\ldots{} or even
worse, getting nowhere \emph{slowly}! Can we repeat this process so much
that the \(f\) on the right-hand side goes away, leaving us with the
solution we're after:
\[f = F(F(F(F(F(F(F(F(F(F(F(F(F(F(F(F(F(F(F(F(F(F(F(F(F(F(F(F(F(F(F(F\ldots \quad\text{?}\]

Well, actually, yes, if we're smart. What we do is this. We start by
\emph{guessing} the solution to our equation. How do we guess? Well, our
solution \(f\) should have \(f(0) = 0\), so just start with any function
with this property. Call it \(f_1\). Then we improve this initial guess
repeatedly by letting \[
  \begin{aligned}
    f_2 &= F(f_1)
  \\f_3 &= F(f_2)
  \\f_4 &= F(f_3)
  \end{aligned}
\] and so on. Now for the fun part: we show that these guesses get
closer and closer to each other\ldots{} so that they converge to some
function \(f\) with \(f = F(f)\)! Voila! With a little more work we can
show that no matter what our initial guess was, our subsequent guesses
approach the same function \(f\), so that the solution \(f\) is unique.

I'm glossing over some details, of course. To prove Picard's theorem we
need to assume the function \(g\) is reasonably nice (continuous isn't
nice enough, we need something like ``Lipschitz continuous''), and our
initial guess should be reasonably nice (continuous will do here). Also,
Picard's theorem only shows that there's a solution defined on some
finite time interval, not the whole real line. (This little twist is
distressing to Plotkin since it complicates the analogy with
mathematical induction. But there must be some slick way to save the
analogy; it's too cute not to be important!)

You can read about Picard's theorem and other related fixed-point
theorems in any decent book on analysis. Personally I'm fond of:

\begin{enumerate}
\def\labelenumi{\arabic{enumi})}
\setcounter{enumi}{1}
\tightlist
\item
  Michael Reed and Barry Simon, \emph{Methods of Modern Mathematical
  Physics}. Vol. 1: \emph{Functional Analysis}. Vol. 2: \emph{Fourier
  Analysis, Self-Adjointness}. Vol. 3: \emph{Scattering Theory}. Vol. 4:
  \emph{Analysis of Operators}. Academic Press, New York, 1980.
\end{enumerate}
\noindent
which is sort of the bible of analysis for mathematical physicists.

Now, it may seem a bit over-elaborate to reformulate the principle of
mathematical induction as a fixed point theorem. However, this way of
looking at recursion is the basis of a lot of theoretical computer
science. It applies not only to recursive definitions of functions but
also recursive definitions of ``types'' like those given in ``Backus-
Naur form'' --- a staple of computer science.

Let me take a simple example that Jim Dolan told me about. Suppose we
have some set of ``letters'' and we want to define the set of all
nonempty ``words'' built from these letters. For example, if our set of
letters was \(L = \{a,b,c\}\) then we would get an infinite set \(W\) of
words like \(a\), \(ca\), \(bb\), \(bca\), \(cbabba\), and on.

In Backus--Naur form we might express this as follows:

\begin{verbatim}
letter ::= a | b | c

word ::= <letter> | <word> <letter>
\end{verbatim}

In English the first line says ``a letter is either a, b, or c'', while
the second says ``a word is either a letter or a word followed by a
letter''. The second one is the interesting part because it's recursive.

In the language of category theory we could say the same thing as
follows. Let \(L\) be our set of letters. Given any set \(S\), let
\[F(S) = L + S \times L\] where \(+\) means disjoint union and
\(\times\) means Cartesian product. Then the set of ``words'' built from
the letters in \(L\) satisfies \(W = F(W)\), or in other words,
\[W = L + W \times L.\] This says ``a word is either a letter or an
ordered pair consisting of a word followed by a letter.'' In short, we
have a fixed point on our hands!

How do we solve this equation? Well, now I'm going to show you something
they never showed you when you first learned set theory. We just use the
usual rules of algebra: \[
  \begin{aligned}
    W &= L + W x L
  \\W - W x L &= L
  \\W x (1 - L) &= L
  \\W &= L/(1 - L)
  \end{aligned}
\] and then expand the answer out as a Taylor series, getting
\[W = L + L\times L + L\times L\times L + \ldots\] This says ``a word is
either a letter or an ordered pair of letters or an ordered triple of
letters or\ldots{}'' Black magic, but it works!

Now, you may wonder exactly what's going on --- when we're allowed to
subtract and divide sets and expand functions of sets in Taylor series
and all that. I'm not an expert on this, but one place to look is in
Joyal's work on ``analytic functors'' (functors that you can expand in
Taylor series):

\begin{enumerate}
\def\labelenumi{\arabic{enumi})}
\setcounter{enumi}{2}
\tightlist
\item
  Andre Joyal, ``Une théorie combinatoire des séries formelles'',
  \emph{Adv.\ Math.\ }\textbf{42} (1981), 1--82.
\end{enumerate}

Before I explain a little of the idea behind this black magic, let me do
another example. I already said that the principle of mathematical
induction could be thought of as guaranteeing the existence of certain
fixed points. But underlying this is something still more basic: the set
of natural numbers is also defined by a fixed point property! Suppose we
take our set of letters above to be set \(\{0\}\) which has only one
element. Then our set of words is \(\{0,00,000,0000,0000,\ldots\}\). We
can think of this as a funny way of writing the set of natural numbers,
so let's call it \(\mathbb{N}\). Also, let's follow von Neumann and
define \[1 = \{0\},\] which is sensible since it's a set with one element.
Then our fixed point equation says: \[\mathbb{N} = \mathbb{N} + 1\] This
is the basic fixed point property of the natural numbers.

At this point some of you may be squirming\ldots{} this stuff looks a
bit weird when you first see it. To make it more rigorous I need to
bring in some category theory, so I'll assume you've read
\protect\hyperlink{week73}{``Week 73''} and
\protect\hyperlink{week76}{``Week 76''} where I explained categories and
functors and isomorphisms.

If you've got a function \(F\colon S\to S\) from some set to itself, a
fixed point of \(F\) is just an element \(x\) for which \(F(x)\) is
\emph{equal} to \(x\). But now suppose we have a functor
\(F\colon \mathcal{C}\to\mathcal{C}\) from some category to itself.
What's a fixed point of this?

Well, we could define it as an object \(x\) of \(\mathcal{C}\) for which
\(F(x) = x\). But if you know a little category theory you'll know that
this sort of ``strict'' fixed point is very boring compared to a
``weak'' fixed point: an object \(x\) of \(\mathcal{C}\) equipped with
an \emph{isomorphism} \[f\colon F(x)\to x.\] Equality is dull,
isomorphism is interesting. It's also very interesting to consider a
more general notion: a ``lax'' fixed point, meaning an object x equipped
with just a \emph{morphism} \[f\colon F(x)\to x.\] Let's consider an
example. Take our category \(\mathcal{C}\) to be the category of sets.
And take our functor \(F\) to be the functor \[F(x) = x + 1\] by which
we mean ``disjoint union of the set \(x\) with the one-element set'' ---
I leave it to you to check that this is a functor. A lax fixed point of
\(F\) is thus a set \(x\) equipped with a function
\[f\colon x + 1\to x\] so the natural numbers
\(\mathbb{N} = \{0,00,000,\ldots\}\) is a lax fixed point in an obvious
way\ldots{} in fact a weak fixed point. So when I wrote
\(\mathbb{N} = \mathbb{N} + 1\) above, I was lying: they're not equal,
they're just isomorphic. Similarly with those other equations involving
sets.

Now, just as any function from a set to itself has a \emph{set} of fixed
points, any functor \(F\) from a category \(\mathcal{C}\) to itself has
a \emph{category} of lax fixed points. An object in this category is
just an object \(x\) of \(\mathcal{C}\) equipped with a morphism
\(f\colon F(x)\to x\), and a morphism from this object to some other
object \(g\colon F(y)\to y\) is just a commutative square: \[
  \begin{tikzcd}
    F(x) \rar["f"] \dar[swap,"F(h)"]
    & x \dar["h"]
  \\F(y) \rar[swap,"g"]
    & y
  \end{tikzcd}
\] In our example, the natural numbers is actually the ``initial'' lax
fixed point, meaning that in the category of lax fixed points there is
exactly one morphism from this object to any other.

So that's the real meaning of these funny recursive definitions in
Backus--Naur form: we have a functor \(F\) from some category like
\(\mathsf{Set}\) to itself, and we are defining an object by saying that
it's the initial lax fixed point of this functor. It's a souped-up
version of defining an element of a set as the unique fixed point of a
function!

I should warn you that category theorists and theoretical computer
scientists usually say ``algebra'' of a functor instead of ``lax fixed
point'' of a functor. Anyway, this gives a bit of a flavor of what those
folks talk about.

\begin{center}\rule{0.5\linewidth}{0.5pt}\end{center}

\textbf{Addendum:} Here's an interesting email that Doug Merritt sent me
after reading the above stuff:

\begin{quote}
A little web searching and discussion with Andras Kornai yields the
following info.

The original work on representing grammars as power series is

\begin{enumerate}
\def\labelenumi{\arabic{enumi})}
\setcounter{enumi}{3}
\tightlist
\item
  N. Chomsky and M. P. Schutzenberger, ``The algebraic theory of
  context-free languages'', in \emph{Computer Programming and Formal
  Systems}, North-Holland, Amsterdam, 1963.
\end{enumerate}

\ldots where Schutzenberger supplied the formal power series aspect,
basically just as the usual generating function trick.

The algebraic connection was developed through the 60's and 70's,
culminating in the work of Samuel Eilenberg, founder of category theory,
such as in

\begin{enumerate}
\def\labelenumi{\arabic{enumi})}
\setcounter{enumi}{4}
\tightlist
\item
  Samuel Eilenberg, Automata, \emph{Languages and Machines}, Academic
  Press, New York, 1974.
\end{enumerate}

A lot of the work in the area comes under the heading ``syntactic
semigroups'', which is fairly self-explanatory (and yields a lot of hits
when web surfing).

The question of expanding a grammar via synthetic division as usual
comes down to the question of whether it is represented as a complete
division algebra or not. Grammars are typically nonabelian, however in
order to use more powerful mathematical machinery, frequently
commutativity is often nonetheless assumed, and the grammar forced into
that Procrustean bed.

I happened across an interesting recent paper (actually a '94 PhD
thesis) that brings all the modern machinery to bear on this sort of
thing (e.g.~explaining how to represent grammars by power series via the
Lagrange Inversion Formula, and multi-non-terminal (multivariable)
grammars via the Generalized LIF), and that is even quite readable:

\begin{enumerate}
\def\labelenumi{\arabic{enumi})}
\setcounter{enumi}{5}
\tightlist
\item
  Ole Vilhelm Larsen, \emph{``Computing order-independent statistical
  characteristics of stochastic context-free languages''}, available at
  {\rm \href{https://web.archive.org/web/19970512050407/http://cwis.auc.dk/phd/fulltext/larsen/html/}{\texttt{https://web.}}} {\rm \href{https://web.archive.org/web/19970512050407/http://cwis.auc.dk/phd/fulltext/larsen/html/}{\texttt{archive.org/web/19970512050407/http://cwis.auc.dk/phd/full}}}  \hfill \break {\rm \href{https://web.archive.org/web/19970512050407/http://cwis.auc.dk/phd/fulltext/larsen/html/}{\texttt{text/larsen/html/}}} or 
 {\rm \href{https://web.archive.org/web/19971110034152/http://cwis.auc.dk/phd/fulltext/larsen/pdf/larsen.pdf}{\texttt{https://web.archive.org/web/19971110034}}}  {\rm \href{https://web.archive.org/web/19971110034152/http://cwis.auc.dk/phd/fulltext/larsen/pdf/larsen.pdf}{\texttt{152/http://cwis.auc.dk/phd/fulltext/larsen/pdf/larsen.pdf}}}
\end{enumerate}

You probably know all this better than I, but: As for fixed points, the
original theorem by Banach applies only to contractive mappings, but
beginning in '68 a flood of new theorems applying to various different
non-contractive situations began to appear, and research continues hot
and heavy. One danger of simply assuming fixed points is that there may
be orbits rather than attractive basins, which I alluded to briefly in
my {\rm \texttt{sci.math}} FAQ entry (which has become somewhat mangled over
the years) concerning the numeric solution of \(f(x) = x^x\) via direct
fixed point recurrence \((F(F(F(F(F...(\text{guess})..))))\). The orbits
cause oscillatory instability in some regions such that it becomes
appropriate to switch to a different technique.

Anyway that's merely to say that there are indeed spaces where one can't
just assume a fixed point theorem and that this can have practical
implications.

Hope that's of some interest.

Doug

\begin{verbatim}
Doug Merritt                       doug@netcom.com
Professional Wild-eyed Visionary   Member, Crusaders for a Better Tomorrow

Unicode Novis Cypherpunks Gutenberg Wavelets Conlang Logli Alife Anthro
Computational linguistics Fundamental physics Cogsci Egyptology GA TLAs
\end{verbatim}
\end{quote}

\hypertarget{week109}{%
\section{September 27, 1997}\label{week109}}

In the Weeks to come I want to talk about quantum gravity. A lot of cool
things have been happening in this subject lately. But I want to start
near the beginning\ldots.

In the 1960's, John Wheeler came up with the notion of ``spacetime
foam''. The idea was that at very short distance scales, quantum
fluctuations of the geometry of spacetime are so violent that the usual
picture of a smooth spacetime with a metric on it breaks down. Instead,
one should visualize spacetime as a ``foam'', something roughly like a
superposition of all possible topologies which only looks smooth and
placid on large enough length scales. His arguments for this were far
from rigorous; they were based on physical intuition. Electromagnetism
and all other fields exhibit quantum fluctuations - so gravity should
too. A wee bit of dimensional analysis suggests that these fluctuations
become significant on a length scale around the Planck length, which is
about 10\textsuperscript{-35} meters. This is very small, much smaller
than what we can probe now. Around this length scale, there's no reason
to suspect that ``perturbative quantum gravity'' should apply, where you
treat gravitational waves as tiny ripples on flat spacetime, quantize
these, and get a theory of ``gravitons''. Indeed, the
nonrenormalizability of quantum gravity suggests otherwise.

Wheeler didn't know what formalism to use to describe ``spacetime
foam'', but he was more concerned with building up a rough picture of
it. Since he is so eloquent, especially when he's giving handwaving
arguments, let me quote him here:

\begin{quote}
No point is more central than this, that empty space is not empty. It
is the seat of the most violent physics. The electromagnetic field
fluctuates. Virtual pairs of positive and negative electrons, in effect,
are constantly being created and annihilated, and likewise pairs of
\(\mu\) mesons, pairs of baryons, and pairs of other particles. All
these fluctuations coexist with the quantum fluctuations in the geometry
and topology of space. Are they additional to those geometrodynamic
zero-point disturbances, or are they, in some sense not now
well-understood, mere manifestations of them?
\end{quote}

That's from:

\begin{enumerate}
\def\labelenumi{\arabic{enumi})}
\tightlist
\item
  Charles Misner, Kip Thorne and John Wheeler, \emph{Gravitation},
  Freeman Press, 1973.
\end{enumerate}
\noindent
It's in the famous last chapter called ``Beyond the end of time''.
Strong stuff! This is what got me interested in quantum gravity in
college. Later I came to prefer less florid writing, and realized how
hard it was to turn gripping prose into actual theories\ldots{} but back
then I ate it up uncritically.

Part of Wheeler's vision was that ultimately physics is all about
geometry, and that particles might be manifestations of this geometry.
For example, electron-positron pairs might be ends of wormholes threaded
by electric field lines:

\begin{quote}
In conclusion, the vision of Riemann, Clifford and Einstein, of a
purely geometric basis for physics, today has come to a higher state of
development, and offers richer prospects --- and presents deeper
problems, than ever before. The quantum of action adds to this
geometrodynamics new features, of which the most striking is the
presence of fluctuations of the wormhole type throughout all space. If
there is any correspondence between this virtual foam-like structure and
the physical vacuum as it has come to be known through quantum
electrodynamics, then there seems to be no escape from identifying these
wormholes with `undressed electrons'. Completely different from these
`undressed electrons', according to all available evidence, are the
electrons and other particles of experimental physics. For these
particles the geometrodynamic picture suggests the model of collective
disturbances in a virtual foam-like vacuum, analogous to different kinds
of phonons or excitons in a solid.
\end{quote}
\noindent
That quote is from:

\begin{enumerate}
\def\labelenumi{\arabic{enumi})}
\setcounter{enumi}{1}
\tightlist
\item
  John Wheeler, \emph{Geometrodynamics}, Academic Press, New York, 1962.
\end{enumerate}

There are many problems with getting this wormhole picture of particles
to work. First, there was --- and is! --- no experimental evidence that
wormholes exist, virtual or otherwise. The main reason for believing in
virtual wormholes was the quantum-mechanical idea that ``whatever is not
forbidden is required''\ldots{} an idea which must be taken with a grain
of salt. Second, there was no mathematical model of ``spacetime foam''
or ``virtual wormholes''. It was just a vague notion.

However, Wheeler was mainly worried about two other problems. First, how
can we relate a space with a wormhole to one without? Since the two have
different topologies, there can't be any continuous way of going from
one to the other. In response to this problem, he suggested that the
description of spacetime in terms of a smooth manifold was not
fundamental, and that we really need some more other description, some
sort of ``pregeometry''. Second, what about the fact that electrons have
spin 1/2? This means that when you turn one around 360 degrees it
doesn't come back to the same state: it picks up a phase of -1. Only
when you turn it around twice does it come back to its original state!
This is nicely described using the mathematics of ``spinors'', but
\emph{not} so nicely described in terms of wormholes.

In his freewheeling, intuitive manner, Wheeler fastened on this second
problem as a crucial clue to the nature of ``pregeometry'':

``It is impossible to accept any description of elementary particles
that does not have a place for spin 1/2. What, then, has any purely
geometric description to offer in explanation of spin 1/2 in general?
More particularly and more importantly, what place is there in quantum
geometrodynamics for the neutrino - the only entity of half-integral
spin that is a pure field in its own right, in the sense that it has
zero rest mass and moves at the speed of light? No clear or satisfactory
answer is known to this question today. Unless and until an answer is
forthcoming, \emph{pure geometrodynamics must be judged deficient as a
basis of elementary particle physics}.''

Physics moves in indirect ways. Though Wheeler's words inspired many
students of relativity, progress on ``spacetime foam'' was quite slow.
It's not surprising, given the thorny problems and the lack of a precise
mathematical model. Quite a bit later, Hawking and others figured out
how to do calculations involving virtual wormholes, virtual black holes
and such using a technique called ``Euclidean quantum gravity''. Pushed
to its extremes, this leads to a theory of spacetime foam, though not
yet a rigorous one (see \protect\hyperlink{week67}{``Week 67''}).

But long before that, Newman, Penrose, and others started finding
interesting relationships between general relativity and the mathematics
of spin-\(1/2\) particles\ldots{} relationships that much later would
yield a theory of spacetime foam in which spinors play a crucial part!

The best place to read about spinorial techniques in general relativity
is probably:

\begin{enumerate}
\def\labelenumi{\arabic{enumi})}
\setcounter{enumi}{2}
\tightlist
\item
  Roger Penrose and Wolfgang Rindler, \emph{Spinors and Space-Time}.
  Vol. 1: \emph{Two-Spinor Calculus and Relativistic Fields}. Vol. 2:
  \emph{Spinor and Twistor Methods in Space-Time Geometry}. Cambridge
  U.\ Press, Cambridge, 1985--1986.
\end{enumerate}

There are roughly 3 main aspects to Penrose's work on spinors and
general relativity. The first is the ``spinor calculus'', described in
volume 1 of these books. By now this is a standard tool in relativity,
and you can find introductions to it in many textbooks, like
``Gravitation'' or Wald's more recent text:

\begin{enumerate}
\def\labelenumi{\arabic{enumi})}
\setcounter{enumi}{3}
\tightlist
\item
  Robert M. Wald, \emph{General Relativity}, U.\ Chicago
  Press, Chicago, 1984.
\end{enumerate}

The second is ``twistor theory'', described in volume 2. This is
mathematically more elaborate, and it includes an ambitious program to
reformulate the laws of physics in such a way that massless spin-\(1/2\)
particles, rather than points of spacetime, play the basic role.

The third is the theory of ``spin networks'', which was a very radical,
purely combinatorial approach to describing the geometry of space.
Penrose's inability to extend it to \emph{spacetime} is what made him
turn later to twistor theory. Probably the best explanation of Penrose's
original spin network ideas can be found in the thesis of one of his
students:

\begin{enumerate}
\def\labelenumi{\arabic{enumi})}
\setcounter{enumi}{4}
\tightlist
\item
  John Moussouris, \emph{Quantum Models of Space-Time Based on
  Recoupling Theory}, Ph.D. thesis, Department of Mathematics,
  University of Oxford, 1983.  Available at \href{https://ora.ox.ac.uk/objects/uuid:6ad25485-c6cb-4957-b129-5124bb2adc67/}{\texttt{https:}} \href{https://ora.ox.ac.uk/objects/uuid:6ad25485-c6cb-4957-b129-5124bb2adc67/}{\texttt{//ora.ox.ac.uk/objects/uuid:6ad25485-c6cb-4957-b129-5124bb2adc67/}}
\end{enumerate}

Here I want to talk about the spinor calculus, which is the most widely
used of these ideas. It's all about the rotation group in 3 dimensions
and the Lorentz group in 3+1 dimensions (by which we mean 3 space
dimensions and 1 time dimension). A lot of physics is based on these
groups. For general stuff about rotation groups and spinors in
\emph{any} dimension, see \protect\hyperlink{week61}{``Week 61''} and
\protect\hyperlink{week93}{``Week 93''}. Here I'll be concentrating on
stuff that only works when we start with *3* space dimensions.

Now I will turn up the math level a notch\ldots.

In the quantum mechanics of angular momentum, what matters is not the
representations of the rotation group \(\mathrm{SO}(3)\), but of its
double cover \(\mathrm{SU}(2)\). This group has one irreducible unitary
representation of each dimension \(d = 1, 2, 3,\ldots\). Physicists
prefer to call these the ``spin-\(j\)'' representations, where
\(j = 0, 1/2, 1, \ldots\). The relation is of course that
\(2j + 1 = d\).

The spin-\(0\) representation is the trivial representation. Physicists
call vectors in this representation ``scalars'', since they are just
complex numbers. Particles transforming in the spin-\(0\) representation
of \(\mathrm{SU}(2)\) are also called scalars. Examples include pions
and other mesons. The only \emph{fundamental} scalar particle in the
Standard Model is the Higgs boson --- hypothesized but still not seen.

The spin-\(1/2\) representation is the fundamental representation, in
which \(\mathrm{SU}(2)\) acts on \(\mathbb{C}^2\) in the obvious way.
Physicists call vectors in this representation ``spinors''. Examples of
spin-\(1/2\) particles include electrons, protons, neutrons, and
neutrinos. The fundamental spin-\(1/2\) particles in the Standard Model
are the leptons (electron, muon, tau and their corresponding neutrinos)
and quarks.

The spin-\(1\) representation comes from turning elements of
\(\mathrm{SU}(2)\) into \(3\times3\) matrices using the double cover
\(\mathrm{SU}(2)\to\mathrm{SO}(3)\). This is therefore also called the
``vector'' representation. The spin-\(1\) particles in the Standard
Model are the gauge fields: the photon, the W and Z, and the gluons.

Though you can certainly make composite particles of higher spin, like
hadrons and atomic nuclei, there are no fundamental particles of spin
greater than \(1\) in the Standard Model. But the Standard Model doesn't
cover gravity. In gravity, the spin-\(2\) representation is very
important. This comes from letting \(\mathrm{SO}(3)\), and thus
\(\mathrm{SU}(2)\), act on symmetric traceless \(3\times3\) matrices in
the obvious way (by conjugation). In perturbative quantum gravity,
gravitons are expected to be spin-\(2\) particles. Why is this? Well, a
cheap answer is that the metric on space is given by a symmetric
\(3\times3\) matrix. But this is not very satisfying\ldots{} I'll give a
better answer later.

Now, the systematic way to get all these representations is to build
them out of the spin-\(1/2\) representation. \(\mathrm{SU}(2)\) acts on
\(\mathbb{C}^2\) in an obvious way, and thus acts on the space of
polynomials on \(\mathbb{C}^2\). The space of homogeneous polynomials of
degree \(2j\) is thus a representation of \(\mathrm{SU}(2)\) in its own
right, called the spin-\(j\) representation. Since multiplication of
polynomials is commutative, in math lingo we say the spin-j
representation is the ``symmetrized tensor product'' of 2j copies of the
spin-\(1/2\) representation. This is the mathematical sense in which
spin-\(1/2\) is fundamental!

(In some sense, this means we can think of a spin-\(j\) particle as
built from \(2j\) indistinguishable spin-\(1/2\) bosons. But there is
something odd about this, since in physics we usually treat spin-\(1/2\)
particles as fermions and form \emph{antisymmetrized} tensor products of
them!)

Now let's go from space to spacetime, and consider the Lorentz group,
\(\mathrm{SO}(3,1)\). Again it's not really this group but its double
cover that matters in physics; its double cover is
\(\mathrm{SL}(2,\mathbb{C})\). Note that \(\mathrm{SL}(2,\mathbb{C})\)
has \(\mathrm{SU}(2)\) as a subgroup just as \(\mathrm{SO}(3,1)\) has
\(\mathrm{SO}(3)\) as a subgroup; everything fits together here, in a
very pretty way.

Now, while \(\mathrm{SU}(2)\) has only one \(2\)-dimensional irreducible
representation, \(\mathrm{SL}(2,\mathbb{C})\) has two, called the
left-handed and right-handed spinor representations. The ``left-handed''
one is the fundamental representation, in which
\(\mathrm{SL}(2,\mathbb{C})\) acts on \(\mathbb{C}^2\) in the obvious
way. The ``right-handed'' one is the conjugate of this, in which we take
the complex conjugate of the entries of our matrix before letting it act
on \(\mathbb{C}^2\) in the obvious way. These two representations become
equivalent when we restrict to \(\mathrm{SU}(2)\)\ldots{} but for
\(\mathrm{SL}(2,\mathbb{C})\) they're not! For example, when we study
particles as representations of \(\mathrm{SL}(2,\mathbb{C})\), it turns
out that neutrinos are left-handed, while antineutrinos are
right-handed.

All the irreducible representations of \(\mathrm{SL}(2,\mathbb{C})\) on
complex vector spaces can be built up from the left-handed and
right-handed spinor representations. Here's how: take the symmetrized
tensor product of \(2j\) copies of the left-handed spin representation
and tensor it with the symmetrized tensor product of \(2k\) copies of
the right-handed one. We call this the \((j,k)\) representation.

People in general relativity have a notation for all this stuff. They
write left-handed spinors as gadgets with one ``unprimed subscript'',
like this: \[v_A\] where \(A = 1,2\). Right-handed spinors are gadgets
with one ``primed subscript'', like: \[w_{A'}\]

where \(A' = 1,2\). As usual, fancier tensors have more subscripts. For
example, guys in the \((j,k)\) representation have \(j\) unprimed
subscripts and \(k\) primed ones, and don't change when we permute the
unprimed subscripts among themselves, or the primed ones among
themselves.

Now \(\mathrm{SO}(3,1)\) has an obvious representation on
\(\mathbb{R}^4\), called the ``vector'' representation for obvious
reasons. If we think of this as a representation of
\(\mathrm{SL}(2,\mathbb{C})\), it's the (1,1) representation. So when
Penrose writes a vector in 4 dimensions, he can do it either the old
way: \[v_a\] where \(a = 0,1,2,3\), or the new spinorial way:
\[v_{AA'}\] where \(A,A' = 1,2\).

Similarly, we can write \emph{any} tensor using spinors with twice as
many indices. This may not seem like a great step forward, but it
actually was\ldots{} because it lets us slice and dice concepts from
general relativity in interesting new ways.

For example, the Riemann curvature tensor describing the curvature of
spacetime is really important in relativity. It has 20 independent
components but it can split up into two parts, the Ricci tensor and Weyl
tensor, each of which have 10 independent components. Thanks to
Einstein's equation, the Ricci tensor at any point of spacetime is
determined by the matter there (or more precisely, by the flow of energy
and momentum through that point). In particular, the Ricci tensor is
zero in the vacuum. The Weyl tensor \[C_{abcd}\] describes aspects of
curvature like gravitational waves or tidal forces which can be nonzero
even in the vacuum. In spinorial notation this is written
\[C{AA'BB'CC'DD'}\] but we can also write it as
\[C_{AA'BB'CC'DD'} = \Phi_{ABCD} \varepsilon_{A'B'}\varepsilon_{C'D'} + \overline{\Phi_{ABCD} \varepsilon_{A'B'}\varepsilon_{C'D'}}\]
where \(\varepsilon\) is the matrix \[
  \left(
    \begin{array}{cc}
      0&1\\-1&0
    \end{array}
  \right)
\] and \(\Phi\) is the ``Weyl spinor''. The Weyl spinor is symmetric in
all its 4 indices so it lives in the (2,0) representation of
\(\mathrm{SL}(2,\mathbb{C})\). Note that this is a \(5\)-dimensional
complex representation, so the Weyl spinor has 10 real degrees of
freedom, just like the Weyl tensor --- but these degrees of freedom have
been encoded in a very efficient way! Even better, we see here why, in
perturbative quantum gravity, the graviton is a spin-\(2\) particle!

I'm only scratching the surface here, but the point is that spinorial
techniques are really handy all over general relativity. A great example
is Witten's spinorial proof of the positive energy theorem:

\begin{enumerate}
\def\labelenumi{\arabic{enumi})}
\setcounter{enumi}{5}
\tightlist
\item
  Edward Witten, ``A new proof of the positive energy theorem'',
  \emph{Comm. Math. Phys.} \textbf{80} (1981), 381--402.
\end{enumerate}

This says that for any spacetime that looks like flat Minkowski space
off at spatial infinity, but possibly has gravitational radiation and
matter in the middle, the ``ADM mass'' is greater than or equal to zero
as long as the matter satisfies the ``dominant energy condition'', which
says that the speed of energy flow is less than the speed of light.
What's the ADM mass? Well, basically the idea is this: if we go off
towards spatial infinity, where spacetime is almost flat and general
relativity effects aren't too big, we can imagine measuring the mass of
the stuff in the middle by seeing how fast a satellite would orbit it.
That's the ADM mass. If the satellite is \emph{attracted} by the stuff
in the middle, the ADM mass is positive. The proof of the positive
energy theorem was really complicated before Witten used spinors, which
let you write the ADM mass as an integral of an obviously nonnegative
quantity.

Next time I'll talk about spin networks and how they show up in recent
work on quantum gravity. We'll see that the idea of building up
everything from the spin-\(1/2\) representation of \(\mathrm{SU}(2)\)
assumes grandiose proportions: in this setup, \emph{space itself} is
built from spinors!

\begin{center}\rule{0.5\linewidth}{0.5pt}\end{center}

\begin{quote}
\emph{The universe is full of magical things, patiently waiting for our
wits to grow sharper.}

--- Eden Philpotts
\end{quote}

\hypertarget{week110}{%
\section{October 4, 1997}\label{week110}}

Last time I sketched Wheeler's vision of ``spacetime foam'', and his
intuition that a good theory of this would require taking spin-\(1/2\)
particles very seriously. Now I want to talk about Penrose's ``spin
networks''. These were an attempt to build a purely combinatorial
description of spacetime starting from the mathematics of spin-\(1/2\)
particles. He didn't get too far with this, which is why he moved on to
invent twistor theory. The problem was that spin networks gave an
interesting theory of \emph{space}, but not of spacetime. But recent
work on quantum gravity shows that you can get pretty far with spin
network technology. For example, you can compute the entropy of quantum
black holes. So spin networks are quite a flourishing business.

Okay. Building space from spin! How does it work?

Penrose's original spin networks were purely combinatorial gadgets:
graphs with edges labelled by numbers \(j = 0, 1/2, 1, 3/2,\ldots\)
These numbers stand for total angular momentum or ``spin''. He required
that three edges meet at each vertex, with the corresponding spins
\(j_1, j_2, j_3\) adding up to an integer and satisfying the triangle
inequalities \[|j_1 - j_2| \leqslant j_3 \leqslant j_1 + j_2.\] These
rules are motivated by the quantum mechanics of angular momentum: if we
combine a system with spin \(j_1\) and a system with spin \(j_2\), the
spin \(j_3\) of the combined system satisfies exactly these constraints.

In Penrose's setup, a spin network represents a quantum state of the
geometry of space. To justify this interpretation he did a lot of
computations using a special rule for computing a number from any spin
network, which is now called the ``Penrose evaluation'' or ``chromatic
evaluation''. In \protect\hyperlink{week22}{``Week 22''} I said how this
works when all the edges have spin 1, and described how this case is
related to the four-color theorem. The general case isn't much harder,
but it's a real pain to describe without lots of pictures, so I'll just
refer you to the original papers:

\begin{enumerate}
\def\labelenumi{\arabic{enumi})}
\item
  Roger Penrose, ``Angular momentum: an approach to combinatorial
  space-time'', in \emph{Quantum Theory and Beyond}, ed.~T. Bastin,
  Cambridge U.\ Press, Cambridge, 1971, pp.~151--180. Also available at
  \url{http://math.ucr.edu/home/baez/penrose/}

  Roger Penrose, ``Applications of negative dimensional tensors'', in
  \emph{Combinatorial Mathematics and its Applications}, ed.~D. Welsh,
  Academic Press, New York, 1971, pp.~221--244.  Also available at
  \url{http://math.ucr.edu/home/baez/penrose/}

  Roger Penrose, ``On the nature of quantum geometry'', in \emph{Magic
  Without Magic}, ed.~J. Klauder, Freeman, San Francisco, 1972,
  pp.~333--354.  Also available at
  \href{http://math.ucr.edu/home/baez/penrose/}{\texttt{http:}}  \href{http://math.ucr.edu/home/baez/penrose/}{\texttt{//math.ucr.edu/home/baez/penrose/}}

  Roger Penrose, ``Combinatorial quantum theory and quantized directions'',
  in \emph{Advances in Twistor Theory}, eds.~L. Hughston and R. Ward,
  Pitman Advanced Publishing Program, San Francisco, 1979, pp.~301--307.
   Also available at
   \href{http://math.ucr.edu/home/baez/penrose/}{\texttt{http://math.}}  \href{http://math.ucr.edu/home/baez/penrose/}{\texttt{ucr.edu/home/baez/penrose/}}
\end{enumerate}

It's easier to explain the \emph{physical meaning} of the Penrose
evaluation. Basically, the idea is this. In classical general
relativity, space is described by a \(3\)-dimensional manifold with a
Riemannian metric: a recipe for measuring distances and angles. In the
spin network approach to quantum gravity, the geometry of space is
instead described as a superposition of ``spin network states''. In
other words, spin networks form a basis of the Hilbert space of states
of quantum gravity, so we can write any state \(\psi\) as
\[\psi = \sum c_i \psi_i\] where \(\psi_i\) ranges over all spin
networks and the coefficients \(c_i\) are complex numbers. The simplest
state is the one corresponding to good old flat Euclidean space. In this
state, each coefficient \(c_i\) is just the Penrose evaluation of the
corresponding spin network \(\psi_i\).

Actually, this interpretation wasn't fully understood until later, when
Rovelli and Smolin showed how spin networks arise naturally in the
so-called ``loop representation'' of quantum gravity. They also came up
with a clearer picture of the way a spin network state corresponds to a
possible geometry of space. The basic picture is that spin network edges
represent flux tubes of area: an edge labelled with spin \(j\)
contributes an area proportional to \(\sqrt{j(j+1)}\) to any surface it
pierces.

The cool thing is that Rovelli and Smolin didn't postulate this, they
\emph{derived} it. Remember, in quantum theory, observables are given by
operators on the Hilbert space of states of the physical system in
question. You typically get these by ``quantizing'' the formulas for the
corresponding classical observables. So Rovelli and Smolin took the
usual formula for the area of a surface in a \(3\)-dimensional manifold
with a Riemannian metric and quantized it. Applying this operator to a
spin network state, they found the picture I just described: the area of
a surface is a sum of terms proportional to \(\sqrt{j(j+1)}\), one for
each spin network edge poking through it.

Of course, I'm oversimplifying both the physics and the history here.
The tale of spin networks and loop quantum gravity is rather long. I've
discussed it already in \protect\hyperlink{week55}{``Week 55''} and
\protect\hyperlink{week99}{``Week 99''}, but only sketchily. If you want
more details, try:

\begin{enumerate}
\def\labelenumi{\arabic{enumi})}
\setcounter{enumi}{1}
\tightlist
\item
  Carlo Rovelli, ``Loop quantum gravity'', \emph{Living Reviews in Relativity}
\textbf{1} (1998). Also available as
  \href{https://arxiv.org/abs/gr-qc/9710008}{\texttt{gr-qc/9710008}}.
\end{enumerate}
\noindent
The abstract gives a taste of what it's all about:

\begin{quote}
The problem of finding the quantum theory of the gravitational field,
and thus understanding what is quantum spacetime, is still open. One of
the most active of the current approaches is loop quantum gravity. Loop
quantum gravity is a mathematically well-defined, non-perturbative and
background independent quantization of general relativity, with its
conventional matter couplings. The research in loop quantum gravity
forms today a vast area, ranging from mathematical foundations to
physical applications. Among the most significant results obtained are:
(i) The computation of the physical spectra of geometrical quantities
such as area and volume; which yields quantitative predictions on
Planck-scale physics. (ii) A derivation of the Bekenstein--Hawking black
hole entropy formula. (iii) An intriguing physical picture of the
microstructure of quantum physical space, characterized by a
polymer-like Planck scale discreteness. This discreteness emerges
naturally from the quantum theory and provides a mathematically
well-defined realization of Wheeler's intuition of a spacetime ``foam''.
Longstanding open problems within the approach (lack of a scalar
product, overcompleteness of the loop basis, implementation of reality
conditions) have been fully solved. The weak part of the approach is the
treatment of the dynamics: at present there exist several proposals,
which are intensely debated. Here, I provide a general overview of
ideas, techniques, results and open problems of this candidate theory of
quantum gravity, and a guide to the relevant literature.
\end{quote}

You'll note from this abstract that the biggest problem in loop quantum
gravity is finding an adequate description of \emph{dynamics}. This is
partially because spin networks are better suited for describing space
than spacetime. For this reason, Rovelli, Reisenberger and I have been
trying to describe spacetime using ``spin foams'' --- sort of like soap
suds with all the bubbles having faces labelled by spins. Every slice of
a spin foam is a spin network.

But I'm getting ahead of myself! I should note that the spin networks
appearing in the loop representation are different from those Penrose
considered, in two important ways.

First, they can have more than 3 edges meeting at a vertex, and the
vertices must be labelled by ``intertwining operators'', or
``intertwiners'' for short. This is a concept coming from group
representation theory; as described in
\protect\hyperlink{week109}{``Week 109''}, what we've been calling
``spins'' are really irreducible representations of \(\mathrm{SU}(2)\).
If we orient the edges of a spin network, we should label each vertex
with an intertwiner from the tensor product of representations on the
``incoming'' edges to the tensor product of representations labelling
the ``outgoing'' edges. When 3 edges labelled by spins \(j_1, j_2, j_3\)
meet at a vertex, there is at most one intertwiner
\[f\colon j_1 \otimes j_2\to j_3,\] at least up to a scalar multiple.
The conditions I wrote down --- the triangle inequality and so on ---
are just the conditions for a nonzero intertwiner of this sort to exist.
That's why Penrose didn't label his vertices with intertwiners: he
considered the case where there's essentially just one way to do it!
When more edges meet at a vertex, there are more intertwiners, and this
extra information is physically very important. One sees this when one
works out the ``volume operators'' in quantum gravity. Just as the spins
on edges contribute \emph{area} to surfaces they pierce, the
intertwiners at vertices contribute \emph{volume} to regions containing
them!

Second, in loop quantum gravity the spin networks are \emph{embedded} in
some 3-dimensional manifold representing space. Penrose was being very
radical and considering ``abstract'' spin networks as a purely
combinatorial replacement for space, but in loop quantum gravity, one
traditionally starts with general relativity on some fixed spacetime and
quantizes that. Penrose's more radical approach may ultimately be the
right one in this respect. The approach where we take classical physics
and quantize it is very important, because we understand classical
physics better, and we have to start somewhere. Ultimately, however, the
world is quantum-mechanical, so it would be nice to have an approach to
space based purely on quantum-mechanical concepts. Also, treating spin
networks as fundamental seems like a better way to understand the
``quantum fluctuations in topology'' which I mentioned in
\protect\hyperlink{week109}{``Week 109''}. However, right now it's
probably best to hedge ones bets and work hard on both approaches.

Lately I've been very excited by a third, hybrid approach:

\begin{enumerate}
\def\labelenumi{\arabic{enumi})}
\setcounter{enumi}{3}
\tightlist
\item
  Andrea Barbieri, ``Quantum tetrahedra and simplicial spin networks'',
   available as
  \href{https://arxiv.org/abs/gr-qc/9707010}{\texttt{gr-qc/9707010}}.
\end{enumerate}
\noindent
Barbieri considers ``simplicial spin networks'': spin networks living in
a fixed \(3\)-dimensional manifold chopped up into tetrahedra. He only
considers spin networks dual to the triangulation, that is, spin
networks having one vertex in the middle of each tetrahedron and one
edge intersecting each triangular face.

In such a spin network there are 4 edges meeting at each vertex, and the
vertex is labelled with an intertwiner of the form
\[f\colon j_1 \otimes j_2\to j_3 \otimes j_4\] where \(j_1,\ldots,j_4\)
are the spins on these edges. If you know about the representation
theory of \(\mathrm{SU}(2)\), you know that \(j_1 \otimes j_2\) is a
direct sum of representations of spin \(j_5\), where \(j_5\) goes from
\(|j_1 - j_2|\) up to \(j_1 + j_2\) in integer steps. So we get a basis
of intertwining operators:
\[f\colon j_1 \otimes j_2\to j_3 \otimes j_4\] by picking one factoring
through each representation \(j_5\):
\[j_1 \otimes j_2\to j_5\to j_3 \otimes j_4\] where:

\begin{enumerate}
\def\labelenumi{\alph{enumi})}
\item
  \(j_1 + j_2 + j_5\) is an integer and
  \(|j_1 - j_2| \leqslant j_5 \leqslant j_1 + j_2\),
\item
  \(j_3 + j_4 + j_5\) is an integer and
  \(|j_3 - j_4| \leqslant j_5 \leqslant j_3 + j_4\).
\end{enumerate}

Using this, we get a basis of simplicial spin networks by labelling all
the edges \emph{and vertices} by spins satisfying the above conditions.
Dually, this amounts to labelling each tetrahedron and each triangle in
our manifold with a spin! Let's think of it this way.

Now focus on a particular simplicial spin network and a particular
tetrahedron. What do the spins \(j_1,\ldots,j_5\) say about the geometry
of the tetrahedron? By what I said earlier, the spins \(j_1,\ldots,j_4\)
describe the areas of the triangular faces: face number 1 has area
proportional to \(\sqrt{j_1(j_1+1)}\), and so on. What about \(j_5\)? It
also describes an area. Take the tetrahedron and hold it so that faces 1
and 2 are in front, while faces 3 and 4 are in back. Viewed this way,
the outline of the tetrahedron is a figure with four edges. The
midpoints of these four edges are the corners of a parallelogram, and
the area of this parallelogram is proportional to \(\sqrt{j_5(j_5+1)}\).
In other words, there is an area operator corresponding to this
parallelogram, and our spin network state is an eigenvector with
eigenvalue proportional to \(\sqrt{j_5(j_5+1)}\). Finally, there is also
a \emph{volume operator} corresponding to the tetrahedron, whose action
on our spin network state is given by a more complicated formula
involving the spins \(j_1,\ldots,j_5\).

Well, that either made sense or it didn't\ldots{} and I don't think
either of us want to stick around to find out which! What's the bottom
line, you ask? First, we're seeing how an ordinary tetrahedron is the
classical limit of a ``quantum tetrahedron'' whose faces have quantized
areas and whose volume is also quantized. Second, we're seeing how to
put together a bunch of these quantum tetrahedra to form a
\(3\)-dimensional manifold equipped with a ``quantum geometry'' ---
which can dually be seen as a spin network. Third, all this stuff fits
together in a truly elegant way, which suggests there is something good
about it. The relationship between spin networks and tetrahedra connects
the theory of spin networks with approaches to quantum gravity where one
chops up space into tetrahedra --- like the ``Regge calculus'' and
``dynamical triangulations'' approaches.

Next week I'll say a bit about using spin networks to study quantum
black holes. Later I'll talk about \emph{dynamics} and spin foams.

Meanwhile, I've been really lagging behind in describing new papers as
they show up\ldots{} so here are a few interesting ones:

\begin{enumerate}
\def\labelenumi{\arabic{enumi})}
\setcounter{enumi}{4}
\item
  Charles Nash, ``Topology and physics --- a historical essay'', in
  \emph{History of Topology}, edited by Ioan M. James, North-Holland, 
   Amsterdam, 1999, pp.\ 359--415. Also available as
  \href{https://arxiv.org/abs/hep-th/9709135}{\texttt{hep-th/9709135}}.
\item
  Luis Alvarez-Gaume and Frederic Zamora, ``Duality in quantum field
  theory (and string theory)'', available as
  \href{https://arxiv.org/abs/hep-th/9709180}{\texttt{hep-th/9709180}}.
\end{enumerate}

Quoting the abstract:

\begin{quote}
``These lectures give an introduction to duality in Quantum Field
Theory. We discuss the phases of gauge theories and the implications of
the electric-magnetic duality transformation to describe the mechanism
of confinement. We review the exact results of N=1 supersymmetric QCD
and the Seiberg--Witten solution of N=2 super Yang--Mills. Some of its
extensions to String Theory are also briefly discussed.''
\end{quote}

\begin{enumerate}
\def\labelenumi{\arabic{enumi})}
\setcounter{enumi}{6}
\tightlist
\item
  Richard E. Borcherds, ``What is a vertex algebra?'', available as
  \href{https://arxiv.org/abs/q-alg/9709033}{\texttt{q-alg/9709033}}.
\end{enumerate}

\begin{quote}
``These are the notes of an informal talk in Bonn describing how to
define an analogue of vertex algebras in higher dimensions.''
\end{quote}

\begin{enumerate}
\def\labelenumi{\arabic{enumi})}
\setcounter{enumi}{7}
\tightlist
\item
  J. M. F. Labastida and Carlos Lozano, ``Lectures in topological
  quantum field theory'', in \emph{AIP Conference Proceedings}
 \textbf{419}, American Institute of Physics, Woodbury, New York,
  1998, pp.\ 54--93. Also available as
  \href{https://arxiv.org/abs/hep-th/9709192}{\texttt{hep-th/9709192}}.
\end{enumerate}

\begin{quote}
``In these lectures we present a general introduction to topological
quantum field theories. These theories are discussed in the framework of
the Mathai--Quillen formalism and in the context of twisted N=2
supersymmetric theories. We discuss in detail the recent developments in
Donaldson--Witten theory obtained from the application of results based
on duality for N=2 supersymmetric Yang--Mills theories. This involves a
description of the computation of Donaldson invariants in terms of
Seiberg--Witten invariants. Generalizations of Donaldson--Witten theory
are reviewed, and the structure of the vacuum expectation values of
their observables is analyzed in the context of duality for the simplest
case.''
\end{quote}

\begin{enumerate}
\def\labelenumi{\arabic{enumi})}
\setcounter{enumi}{8}
\tightlist
\item
  Martin Markl, ``Simplex, associahedron, and cyclohedron'', 
  available as
  \href{https://arxiv.org/abs/alg-geom/9707009}{\texttt{alg-geom/9707009}}.
\end{enumerate}

\begin{quote}
``The aim of the paper is to give an `elementary' introduction to the
theory of modules over operads and discuss three prominent examples of
these objects --- simplex, associahedron (= the Stasheff polyhedron) and
cyclohedron (= the compactification of the space of configurations of
points on the circle).''
\end{quote}

\hypertarget{week111}{%
\section{October 24, 1997}\label{week111}}

This week I'll say a bit about black hole entropy, and next week I'll
say a bit about attempts to compute it using spin networks, as promised.
Be forewarned: all of this stuff should be taken with a grain of salt,
since there is no experimental evidence backing it up. Also, my little
``history'' of the subject is very amateur. (In particular, when I say
someone did something in such-and-such year, all I mean is that it was
published in that year.)

Why is the entropy of black holes so interesting? Mainly because it
serves as a testing ground for our understanding of quantum gravity. In
classical general relativity, any object that falls into a black hole
is, in some sense, lost and gone forever. Once it passes the ``event
horizon'', it can never get out again. This leads to a potential paradox
regarding the second law of thermodynamics, which claims that the total
entropy of the universe can never decrease. My office desk constantly
increases in entropy as it becomes more cluttered and dusty. I could
reduce its entropy with some work, dusting it and neatly stacking up the
papers and books, but in the process I would sweat and metabolize,
increasing my \emph{own} entropy even more --- so I don't bother. If,
however, I could simply dump my desk into a black hole, perhaps I could
weasel around the second law. True, the black hole would be more
massive, but nobody could tell if I'd dumped a clean desk or a dirty
desk into it, so in a sense, the entropy would be \emph{gone}!

Of course there are lots of potential objections to this method of
violating the second law. \emph{Anything} involving the second law of
thermodynamics is controversial, and the idea of violating it by
throwing entropy down black holes is especially so. The whole subject
might have remained a mere curiosity if it hadn't been for the work of
Hawking and Penrose.

In 1969, Penrose showed that, in principle, one could extract energy
from a rotating black hole:

\begin{enumerate}
\def\labelenumi{\arabic{enumi})}
\tightlist
\item
  Roger Penrose, ``Gravitational collapse: the role of general
  relativity'', \emph{Rev.~del Nuovo Cimento} \textbf{1}, (1969)
  272--276.
\end{enumerate}
\noindent
Basically, one can use the rotation of the black hole to speed up a rock
one throws past it, as long as the rock splits and one piece falls in
while the rock is in the ``ergosphere'' --- the region of spacetime
where the ``frame dragging'' is so strong that any object inside is
\emph{forced} to rotate along with it. This result led to a wave of
thought experiments and theorems involving black holes and
thermodynamics.

In 1971, Hawking proved the ``black hole area theorem'':

\begin{enumerate}
\def\labelenumi{\arabic{enumi})}
\setcounter{enumi}{1}
\tightlist
\item
  Stephen Hawking, ``Gravitational radiation from colliding black
  holes'', \emph{Phys. Rev.~Lett.} \textbf{26} (1971), 1344--1346.
\end{enumerate}
\noindent
The precise statement of this theorem is a bit technical, but loosely,
it says that under reasonable conditions --- e.g., no ``exotic matter''
with negative energy density or the like --- the total area of the event
horizons of any collection of black holes must always increase. This
result sets an upper limit on how much energy one can extract from a
rotating black hole, how much energy can be released in a black hole
collision, etc.

Now, this sounds curiously similar to the second law of thermodynamics,
with the area of the black hole playing the role of entropy! It turned
out to be just the beginning of an extensive analogy between black hole
physics and thermodynamics. I have a long way to go, so I will just
summarize this analogy in one cryptic chart:

\begin{longtable}[]{@{}rll@{}}
\toprule
& Black holes & Thermodynamics\tabularnewline
\midrule
\endhead
& black hole mass \(M\) & energy \(E\)\tabularnewline
& event horizon area \(A\) & entropy \(S\)\tabularnewline
& surface gravity \(K\) & temperature \(T\)\tabularnewline
\emph{First law:} & \(dM=KdA/8\pi+\mathrm{work}\) &
\(dE=TdS+\mathrm{work}\)\tabularnewline
\emph{Second law:} & \(A\) increases & \(S\) increases\tabularnewline
\emph{Third law:} & can't get \(K=0\) & can't get \(T=0\)\tabularnewline
\bottomrule
\end{longtable}

For a more thorough review by someone who really knows this stuff, try:

\begin{enumerate}
\def\labelenumi{\arabic{enumi})}
\setcounter{enumi}{2}
\tightlist
\item
  Robert M. Wald, ``Black holes and thermodynamics'', in \emph{Symposium
  on Black Holes and Relativistic Stars (in honor of S. Chandrasekhar)},
  December 14-15, 1996, available as
  \href{https://arxiv.org/abs/gr-qc/9702022}{\texttt{gr-qc/9702022}}.
\end{enumerate}

In 1973, Jacob Bekenstein suggested that we take this analogy really
seriously. In particular, he argued that black holes really do have
entropy proportional to their surface area. In other words, the total
entropy of the world is the entropy of all the matter \emph{plus} some
constant times the area of all the black holes:

\begin{enumerate}
\def\labelenumi{\arabic{enumi})}
\setcounter{enumi}{3}
\tightlist
\item
  Jacob Bekenstein, ``Black holes and entropy'', \emph{Phys. Rev.}
  \textbf{D7} (1973), 2333--2346.
\end{enumerate}
\noindent
This raises an obvious question --- what's the constant?? Also, in the
context of classical general relativity, there are serious problems with
this idea: you can cook up thought experiments where the total entropy
defined this way goes down, no matter what you say the constant is.

However, in 1975, Hawking showed that black holes aren't quite black!

\begin{enumerate}
\def\labelenumi{\arabic{enumi})}
\setcounter{enumi}{4}
\tightlist
\item
  Stephen Hawking, ``Particle creation by black holes'', \emph{Commun.
  Math. Phys.} \textbf{43} (1975), 199--220.
\end{enumerate}
\noindent
More precisely, using quantum field theory on curved spacetime, he
showed that a black hole should radiate photons thermally, with a
temperature \(T\) proportional to the surface gravity \(K\) at the event
horizon. It's important to note that this isn't a ``quantum gravity''
calculation; it's a semiclassical approximation. Gravity is treated
classically! One simply assumes spacetime has the ``Schwarzschild
metric'' corresponding to a black hole. Quantum mechanics enters only in
treating the electromagnetic field. The goal of everyone ever since has
been to reproduce Hawking's result using a full-fledged quantum gravity
calculation. The problem, of course, is to get a theory of quantum
gravity.

Anyway, here is Hawking's magic formula: \[T = K / 2\pi.\] Here I'm
working in units where \(\hbar = c = k = G = 1\), but it's important to
note that there is secretly an \(\hbar\) (Planck's constant) hiding in
this formula, so that it \emph{only makes sense quantum-mechanically}.
This is why Bekenstein's proposal was problematic at the purely
classical level.

This formula means we can take really seriously the analogy between
\(T\) and \(K\). We even know how to convert between the two! Of course,
we also know how to convert between the black hole mass \(M\) and energy
\(E\): \[E = M.\] Thus, using the first law (shown in the chart above),
we can convert between entropy and area as well: \[S = A/4.\] How could
we hope to get such a formula using a full-fledged quantum gravity
calculation? Well, in statistical mechanics, the entropy of any
macrostate of a system is the logarithm of the microstates corresponding
to that macrostate. A microstate is a complete precise description of
the system's state, while a macrostate is a rough description. For
example, if I tell you ``my desk is dusty and covered with papers'', I'm
specifying a macrostate. If there are \(N\) ways my desk could meet this
description (i.e., \(N\) microstates corresponding to that macrostate),
the entropy of my desk is \(\ln(N)\).

We expect, or at least hope, that a working quantum theory of gravity
will provide a statistical-mechanical way to calculate the entropy of a
black hole. In other words, we hope that specifying the area \(A\) of
the black hole horizon specifies the macrostate, and that there are
about \(N = \exp(A/4)\) microstates corresponding to this macrostate.

What are these microstates? Much ink has been spilt over this thorny
question, but one reasonable possibility is that they are \emph{states
of the geometry of the event horizon}. If we know its area, there are
still lots of geometries that the event horizon could have\ldots{} and
perhaps, for some reason, there are about \(\exp(A/4)\) of them! For
this to be true, we presumably need some theory of ``quantum geometry'',
so that the number of geometries is finite.

I presume you see what I'm leading up to: the idea of computing black
hole entropy using spin networks, which are designed precisely to
describe ``quantum geometries'', as sketched in
\protect\hyperlink{week55}{``Week 55''},
\protect\hyperlink{week99}{``Week 99''}, and
\protect\hyperlink{week110}{``Week 110''}. I'll talk about this next
week.

To be fair to other approaches, I should emphasize that string theorists
have their own rather different ideas about computing black hole entropy
using \emph{their} approach to quantum gravity. A nice review article
about this approach is:

\begin{enumerate}
\def\labelenumi{\arabic{enumi})}
\setcounter{enumi}{5}
\tightlist
\item
  Gary Horowitz, ``Quantum states of black holes'', available
  as \href{https://arxiv.org/abs/gr-qc/9704072}{\texttt{gr-qc/9704072}}.
\end{enumerate}
\noindent
Next time, however, I will only talk about the spin network (also known
as ``loop representation'') approach, because that's the one I
understand.

Okay, let me wrap up by listing a few interesting papers here and there
which are contributing to the entropy of my desk. It's 1:30 am and I'm
getting tired, so I'll just cop out and quote the abstracts. The first
one is a short readable essay explaining the limitations of quantum
field theory. The others are more technical.

\begin{enumerate}
\def\labelenumi{\arabic{enumi})}
\setcounter{enumi}{6}
\tightlist
\item
  Roman Jackiw, ``What is quantum field theory and why have some
  physicists abandoned it?'', available as
  \href{https://arxiv.org/abs/hep-th/9709212}{\texttt{hep-th/9709212}}.
\end{enumerate}

\begin{quote}
The present-day crisis in quantum field theory is described.
\end{quote}

\begin{enumerate}
\def\labelenumi{\arabic{enumi})}
\setcounter{enumi}{7}
\tightlist
\item
  Adel Bilal, ``M(atrix) theory: a pedagogical introduction'', \emph{Fortsch.\ 
  Phys.} \textbf{47} (1999), 5--28.  Also available as 
  \href{https://arxiv.org/abs/hep-th/9710136}{\texttt{hep-th/9710136}}.
\end{enumerate}

\begin{quote}
I attempt to give a pedagogical introduction to the matrix model of
M-theory as developed by Banks, Fischler, Shenker and Susskind (BFSS).
In the first lecture, I introduce and review the relevant aspects of
D-branes with the emergence of the matrix model action. The second
lecture deals with the appearance of eleven-dimensional supergravity and
M-theory in strongly coupled type IIA superstring theory. The third
lecture combines the material of the two previous ones to arrive at the
BFSS conjecture and explains the evidence presented by these authors.
The emphasis is not on most recent developments but on a hopefully
pedagogical presentation.
\end{quote}

Here's one on glueballs (for more on glueballs, see
\protect\hyperlink{week68}{``Week 68''}):

\begin{enumerate}
\def\labelenumi{\arabic{enumi})}
\setcounter{enumi}{8}
\tightlist
\item
  Gregory Gabadadze, ``Modeling the glueball spectrum by a closed
  bosonic membrane'', \emph{Phys.\ Rev.\ D} \textbf{58} (1998), 094015. 
  Also available as
  \href{https://arxiv.org/abs/hep-ph/9710402}{\texttt{hep-ph/9710402}}.
\end{enumerate}

\begin{quote}
We use an analogy between the Yang--Mills theory Hamiltonian and the
matrix model description of the closed bosonic membrane theory to
calculate the spectrum of glueballs in the large \(N_c\) limit. Some
features of the Yang--Mills theory vacuum, such as the screening of the
topological charge and vacuum topological susceptibility are discussed.
We show that the topological susceptibility has different properties
depending on whether it is calculated in the weak coupling or strong
coupling regimes of the theory. A mechanism of the formation of the
pseudoscalar glueball state within pure Yang--Mills theory is proposed
and studied.
\end{quote}

Fans of quaternions and octonions might like the following paper:

\begin{enumerate}
\def\labelenumi{\arabic{enumi})}
\setcounter{enumi}{9}
\tightlist
\item
  Jose M. Figueroa-O'Farrill, ``Gauge theory and the division
  algebras'', \emph{J.\ Geom.\ Phys.} \textbf{32} (1999), 
  227--240.  Also available as
  \href{https://arxiv.org/abs/hep-th/9710168}{\texttt{hep-th/9710168}}.
\end{enumerate}

\begin{quote}
We present a novel formulation of the instanton equations in
8-dimensional Yang--Mills theory. This formulation reveals these
equations as the last member of a series of gauge-theoretical equations
associated with the real division algebras, including flatness in
dimension 2 and (anti-)self-duality in 4. Using this formulation we
prove that (in flat space) these equations can be understood in terms of
moment maps on the space of connections and the moduli space of
solutions is obtained via a generalised symplectic quotient: a Kaehler
quotient in dimension 2, a hyperkaehler quotient in dimension 4 and an
octonionic Kaehler quotient in dimension 8. One can extend these
equations to curved space: whereas the \(2\)-dimensional equations make
sense on any surface, and the \(4\)-dimensional equations make sense on
an arbitrary oriented manifold, the \(8\)-dimensional equations only
make sense for manifolds whose holonomy is contained in
\(\mathrm{Spin}(7)\). The interpretation of the equations in terms of
moment maps further constraints the manifolds: the surface must be
orientable, the 4-manifold must be hyperkaehler and the 8-manifold must
be flat.
\end{quote}

\hypertarget{week112}{%
\section{November 3, 1997}\label{week112}}

This week I will talk about attempts to compute the entropy of a black
hole by counting its quantum states, using the spin network approach to
quantum gravity.

But first, before the going gets tough and readers start dropping like
flies, I should mention the following science fiction novel:

\begin{enumerate}
\def\labelenumi{\arabic{enumi})}
\tightlist
\item
  Greg Egan, \emph{Distress}, HarperCollins, 1995.
\end{enumerate}

I haven't been keeping up with science fiction too carefully lately, so
I'm not really the best judge. But as far as I can tell, Egan is one of
the few practitioners these days who bites off serious chunks of reality
--- who really tries to face up to the universe and its possibilies in
their full strangeness. Reality is outpacing our imagination so fast
that most attempts to imagine the future come across as miserably
unambitious. Many have a deliberately ``retro'' feel to them - space
operas set in Galactic empires suspiciously similar to ancient Rome,
cyberpunk stories set in dark urban environments borrowed straight from
film noire, complete with cynical voiceovers\ldots{} is science fiction
doomed to be an essentially \emph{nostalgic} form of literature?

Perhaps we are becoming too wise, having seen how our wildest
imaginations of the future always fall short of the reality, blindly
extrapolating the current trends while missing out on the really
interesting twists. But still, science fiction writers have to try to
imagine the unimaginable, right? If they don't, who will?

But how do we \emph{dare} imagine what things will be like in, say, a
century, or a millenium? Vernor Vinge gave apt expression to this
problem in his novel featuring the marooned survivors of a
``singularity'' at which the rate of technological advance became,
momentarily, \emph{infinite}, and most of civilization
inexplicably\ldots{} disappeared. Those who failed to catch the bus were
left wondering just where it went. Somewhere unimaginable, that's all
they know.

``Distress'' doesn't look too far ahead, just to 2053. Asexuality is
catching on bigtime\ldots{} as are the ``ultramale'' and ``ultrafemale''
options, for those who don't like this gender ambiguity business.
Voluntary Autists are playing around with eliminating empathy. And some
scary radical secessionists are redoing their genetic code entirely,
replacing good old ATCG by base pairs of their own devising. Fundamental
physics, thank god, has little new to offer in the way of new
technology. For decades, it's drifted off introspectively into more and
more abstract and mathematical theories, with few new experiments to
guide it. But this is the year of the Einstein Centenary Conference!
Nobel laureate Violet Masala will unveil her new work on a Theory of
Everything. And rumors have it that she may have finally cracked the
problem, and found --- yes, that's right --- the final, correct and true
theory of physics.

As science reporter Andrew Worth tries to bone up for his interviews
with Masala, he finds it's not so easy to follow the details of the
various ``All-Topology Models'' that have been proposed to explain the
10-dimensionality of spacetime in the Standard Unified Field Theory. In
one of the most realistic passages of imagined mathematical prose I've
ever seen in science fiction, he reads ``At least two conflicting
generalized measures can be applied to T, the space of all topological
spaces with countable basis. Perrini's measure {[}Perrini, 2012{]} and
Saupe's measure {[}Saupe, 2017{]} are both defined for all bounded
subsets of T, and are equivalent when restricted to M - the space of
n-dimensional paracompact Hausdorff manifolds - but they yield
contradictory results for sets of more exotic spaces. However, the
physical significance (if any) of this discrepancy remains
obscure\ldots.''

But, being a hardy soul and a good reporter, Worth is eventually able to
explain to us readers what's at stake here, and \emph{why} Masala's new
work has everyone abuzz. But that's really just the beginning. For in
addition to this respectable work on All-Topology Models, there is a lot
of somewhat cranky stuff going on in ``anthrocosmology'', involving
sophisticated and twisted offshoots of the anthropic principle. Some
argue that when the correct Theory of Everything is found, a kind of
cosmic self-referential feedback loop will be closed. And then there's
no telling \emph{what} will happen!

Well, I won't give away any more. It's fun: it made me want to run out
and do a lot more mathematical physics. And it raises a lot of deep
issues. At the end it gets a bit too ``action-packed'' for my taste, but
then, my idea of excitement is lying in bed thinking about
\(n\)-categories.

Now for the black holes.

In \protect\hyperlink{week111}{``Week 111''}, I left off with a puzzle.
In a quantum theory of gravity, the entropy of a black hole should be
the logarithm of the number of its microstates. This should be
proportional to the area of the event horizon. But what \emph{are} the
microstates? String theory has one answer to this, but I'll focus on the
loop representation of quantum gravity. This approach to quantum gravity
is very geometrical, which suggests thinking of the black hole
microstates as ``quantum geometries'' of the black hole event horizon.
But how are these related to the description of the geometry of the
surrounding space in terms of spin networks?

Starting in 1995, Smolin, Krasnov, and Rovelli proposed some answers to
these puzzles, which I have already mentioned in
\protect\hyperlink{week56}{``Week 56''},
\protect\hyperlink{week57}{``Week 57''}, and
\protect\hyperlink{week87}{``Week 87''}. The ideas I'm going to talk
about now are a further development of this earlier work, but instead of
presenting everything historically, I'll just present the picture as I
see it now. For more details, try the following paper:

\begin{enumerate}
\def\labelenumi{\arabic{enumi})}
\setcounter{enumi}{1}
\tightlist
\item
  Abhay Ashtekar, John Baez, Alejandro Corichi and Kirill Krasnov,
  ``Quantum geometry and black hole entropy'', to appear in \emph{Phys.
  Rev.~Lett.} \textbf{80} (1998), 904--907.   Also available as
  \href{https://arxiv.org/abs/gr-qc/9710007}{\texttt{gr-qc/9710007}}.
\end{enumerate}

This is a summary of what will eventually be a longer paper with two
parts, one on the ``black hole sector'' of classical general relativity,
and one on the quantization of this sector. Let me first say a bit about
the classical aspects, and then the quantum aspects.

One way to get a quantum theory of a black hole is to figure out what a
black hole is classically, get some phase space of classical states, and
then quantize \emph{that}. For this, we need some way of saying which
solutions of general relativity correspond to black holes. This is
actually not so easy. The characteristic property of a black hole is the
presence of an event horizon --- a surface such that once you pass it
you can never get back out without going faster than light. This makes
it tempting to find ``boundary conditions'' which say ``this surface is
an event horizon'', and use those to pick out solutions corresponding to
black holes.

But the event horizon is not a local concept. That is, you can't tell
just by looking at a small patch of spacetime if it has an event horizon
in it, since your ability to ``eventually get back out'' after crossing
a surface depends on what happens to the geometry of spacetime in the
future. This is bad, technically speaking. It's a royal pain to deal
with nonlocal boundary conditions, especially boundary conditions that
depend on \emph{solving the equations of motion to see what's going to
happen in the future just to see if the boundary conditions hold}.

Luckily, there is a purely local concept which is a reasonable
substitute for the concept of event horizon, namely the concept of
``outer marginally trapped surface''. This is a bit technical - and my
speciality is not this classical general relativity stuff, just the
quantum side of things, so I'm no expert on it! - but basically it works
like this.

First consider an ordinary sphere in ordinary flat space. Imagine light
being emitted outwards, the rays coming out normal to the surface of the
sphere. Clearly the cross-section of each little imagined circular ray
will \emph{expand} as it emanates outwards. This is measured
quantitatively in general relativity by a quantity called\ldots{} the
expansion parameter!

Now suppose your sphere surrounds a spherically symmetric black hole. If
the sphere is huge compared to the size of the black hole, the above
picture is still pretty accurate, since the light leaving the sphere is
very far from the black hole, and gravitational effects are small. But
now imagine shrinking the sphere, making its radius closer and closer to
the Schwarzschild radius (the radius of the event horizon). When the
sphere is just a little bigger than the Schwarzschild radius, the
expansion of light rays going out from the sphere is very small. This
might seem paradoxical - how can the outgoing light rays not expand? But
remember, spacetime is seriously warped near the event horizon, so your
usual flat spacetime intuitions no longer apply. As we approach the
event horizon itself, the expansion parameter goes to zero!

That's roughly the definition of an ``outer marginally trapped
surface''. A more mathematical but still rough definition is: ``an outer
marginally trapped surface is the boundary \(S\) of some region of space
such that the expansion of the outgoing family of null geodesics normal
to \(S\) is everywhere less than or equal to zero.''

We require that our space have some sphere \(S\) in it which is an outer
marginally trapped surface. We also require other boundary conditions to
hold on this surface. I won't explain them in detail. Instead, I'll say
two important extra features they have: they say the black hole is
nonrotating, and they disallow gravitational waves falling into \(S\).
The first condition here is a simplifying assumption: we are only
studying black holes of zero angular momentum in this paper! The second
condition is only meant to hold for the time during which we are
studying the black hole. It does not rule out gravitational waves far
from the black hole, waves that might \emph{eventually} hit the black
hole. These should not affect the entropy calculation.

Now, in addition to their physical significance, the boundary conditions
we use also have an interesting \emph{mathematical} meaning. Like most
other field theories, general relativity is defined by an action
principle, meaning roughly that one integrates some quantity called the
Lagrangian over spacetime to get an ``action'', and finds solutions of
the field equations by looking for minima of this action. But when one
studies field theories on a region with boundary, and imposes boundary
conditions, one often needs to ``add an extra boundary term to the
action'' --- some sort of integral over the boundary --- to get things
to work out right. There is a whole yoga of finding the right boundary
term to go along with the boundary conditions\ldots{} an arcane little
art\ldots{} just one of those things theoretical physicists do, that for
some reason never find their way into the popular press.

But in this case the boundary term is all-important, because
it's\ldots{}
\vskip 1em
\centerline{ \textbf{the Chern--Simons action}! }
\vskip 1em
(Yes, I can see people world-wide, peering into their screens, thinking
``Eh? Am I supposed to remember what that is? What's he getting so
excited about now?'' And a few cognoscenti thinking ``Oh, \emph{now} I
get it. All this fussing about boundary conditions was just an elaborate
ruse to get a topological quantum field theory on the event horizon!'')

So far we've been studying general relativity in honest
\(4\)-dimensional spacetime. Chern--Simons theory is a closely related
field theory one dimension down, in \(3\)-dimensional spacetime. As time
passes, the surface of the black hole traces out a \(3\)-dimensional
submanifold of our 4-dimensional spacetime. When we quantize our
classical theory of gravity with our chosen boundary conditions, the
Chern--Simons term will give rise to a ``Chern--Simons field theory''
living on the surface of the black hole. This field theory will describe
the geometry of the surface of the black hole, and how it changes as
time passes.

Well, let's not just talk about it, let's do it! We quantize our theory
using standard spin network techniques \emph{outside} the black hole,
and Chern--Simons theory \emph{on the event horizon}, and here is what we
get. States look like this. Outside the black hole, they are described
by spin networks (see \protect\hyperlink{week110}{``Week 110''}). The
spin network edges are labelled by spins \(j = 0, 1/2, 1,\ldots\). Spin
network edges can puncture the black hole surface, giving it area. Each
spin-\(j\) edge contributes an area proportional to \(\sqrt{j(j+1)}\).
The total area is the sum of these contributions.

Any choice of punctures labelled by spins determines a Hilbert space of
states for Chern--Simons theory. States in this space describe the
intrinsic curvature of the black hole surface. The curvature is zero
except at the punctures, so that \emph{classically}, near any puncture,
you can visualize the surface as a cone with its tip at the puncture.
The curvature is concentrated at the tip. At the \emph{quantum} level,
where the puncture is labelled with a spin \(j\), the curvature at the
puncture is described by a number \(j_z\) ranging from \(-j\) to \(j\)
in integer steps.

Now we ask the following question: ``given a black hole whose area is
within \(\varepsilon\) of \(A\), what is the logarithm of the number of
microstates compatible with this area?'' This should be the entropy of
the black hole. To figure it out, first we work out all the ways to
label punctures by spins \(j\) so that the total area comes within
\(\varepsilon\) of \(A\). For any way to do this, we then count the
allowed ways to pick numbers \(j_z\) describing the intrinsic curvature
of the black hole surface. Then we sum these up and take the logarithm.

That's roughly what we do, anyway, and for black holes much bigger than
the Planck scale we find that the entropy is proportional to the area.
How does this compare with the result of Bekenstein and Hawking,
described in \protect\hyperlink{week111}{``Week 111''}? Remember, they
computed that \[S = A/4\] where \(S\) is the entropy and \(A\) is the
area, measured in units where \(c = \hbar = G = k = 1\). What we get is
\[S = \frac{\ln 2}{4\pi\gamma\sqrt{3}} A\] To compare these results, you
need to know what that mysterious ``\(\gamma\)'' factor is in the second
equation! It's often called the Immirzi parameter, because many people
learned of it from the following paper:

\begin{enumerate}
\def\labelenumi{\arabic{enumi})}
\setcounter{enumi}{2}
\tightlist
\item
  Giorgio Immirzi, ``Quantum gravity and Regge calculus'',
  \emph{Nucl. Phys. Proc. Suppl.} \textbf{57} (1997) 65--72.  Also
  available as
  \href{https://arxiv.org/abs/gr-qc/9701052}{\texttt{gr-qc/9701052}}.
\end{enumerate}
\noindent
However, it was first discovered by Barbero:

\begin{enumerate}
\def\labelenumi{\arabic{enumi})}
\setcounter{enumi}{3}
\tightlist
\item
  Fernando Barbero, ``Real Ashtekar variables for Lorentzian signature
  space-times'', \emph{Phys. Rev.} \textbf{D51} (1995), 5507--5510.
  Also available as
  \href{https://arxiv.org/abs/gr-qc/9410014}{\texttt{gr-qc/9410014}}.
\end{enumerate}
\noindent
It's an annoying unavoidable arbitrary dimensionless parameter that
appears in the loop representation, which nobody had noticed before
Barbero. It's still rather mysterious. But it works a bit like this. In
ordinary quantum mechanics we turn the position \(q\) into an operator,
namely multiplication by \(x\), and also turn the momentum \(p\) into an
operator, namely \(-i(d/dx)\). The important thing is the canonical
commutation relations: \(pq-qp=-i\). But we could also get the canonical
commutation relations to hold by defining \[
  \begin{aligned}
    p &= -i \gamma \frac{d}{dx}
  \\q &= \frac{x}{\gamma}
  \end{aligned}
\] since the gammas cancel out! In this case, putting in a \(\gamma\)
factor doesn't affect the physics. One gets ``equivalent representations
of the canonical commutation relations''. In the loop representation,
however, the analogous trick \emph{does} affect the physics ---
different choices of the Immirzi parameter give different physics! For
more details try:

\begin{enumerate}
\def\labelenumi{\arabic{enumi})}
\setcounter{enumi}{3}
\tightlist
\item
  Carlo Rovelli and Thomas Thiemann, ``The Immirzi parameter in quantum
  general relativity'', \emph{Phys. Rev.} \textbf{D57} (1998), 1009--1014.
   Also available as
  \href{https://arxiv.org/abs/gr-qc/9705059}{\texttt{gr-qc/9705059}}.
\end{enumerate}

How does the Immirzi parameter affect the physics? It \emph{determines
the quantization of area}. You may notice how I keep saying ``each
spin-\(j\) edge of a spin network contributes an area proportional to
\(\sqrt{j(j+1)}\) to any surface it punctures''\ldots{} without ever
saying what the constant of proportionality is! Well, the constant is
\[8 \pi \gamma\] Until recently, everyone went around saying the
constant was \(1\). (As for the factor of \(8\pi\), I'm no good at these
things, but apparently at least some people were getting that wrong,
too!) Now Krasnov claims to have gotten these damned factors
straightened out once and for all:

\begin{enumerate}
\def\labelenumi{\arabic{enumi})}
\setcounter{enumi}{4}
\tightlist
\item
  Kirill Krasnov, ``On the constant that fixes the area spectrum in
  canonical quantum gravity'', \emph{Class. Quant. Grav.} \textbf{15}
  (1998), L1--L4.  Also available as
  \href{https://arxiv.org/abs/gr-qc/9709058}{\texttt{gr-qc/9709058}}.
\end{enumerate}

So: it seems we can't determine the constant of proportionality in the
entropy-area relation, because of this arbitrariness in the Immirzi
parameter. But we can, of course, use the Bekenstein--Hawking formula
together with our formula for black hole entropy to determine
\(\gamma\), obtaining 
\[\gamma = \frac{\ln 2 }{\pi\sqrt{3}}\] 
This may
seem like cheating, but right now it's the best we can do. All we can
say is this: we have a theory of the microstates of a black hole, which
predicts that entropy is proportional to area for largish black holes,
and which taken together with the Bekenstein--Hawking calculation allows
us to determine the Immirzi parameter.

What do the funny constants in the formula
\[S = \frac{\ln 2}{4\pi\gamma\sqrt{3}} A\] mean? It's actually simple.
The states that contribute most to the entropy of a black hole are those
where nearly all spin network edges puncturing its surface are labelled
by spin \(1/2\). Each spin-\(1/2\) puncture can have either
\(j_z = 1/2\) or \(j_z = -1/2\), so it contributes \(\ln(2)\) to the
entropy. On the other hand, each spin-\(1/2\) edge contributes
\(4\pi\gamma\sqrt{3}\) to the area of the black hole. Just to be
dramatic, we can call \(\ln 2\) the ``quantum of entropy'' since it's
the entropy (or information) contained in a single bit. Similarly, we
can call \(4\pi\gamma\sqrt{3}\) the ``quantum of area'' since it's the
area contributed by a spin-\(1/2\) edge. These terms are a bit
misleading since neither entropy nor area need come in \emph{integral}
multiples of this minimal amount. But anyway, we have
\[S = \frac{\text{quantum of entropy}}{\text{quantum of area}} A\] What
next? Well, one thing is to try to use these ideas to study Hawking
radiation. That's hard, because we don't understand \emph{Hamiltonians}
very well in quantum gravity, but Krasnov has made some progress\ldots.

\begin{enumerate}
\def\labelenumi{\arabic{enumi})}
\setcounter{enumi}{5}
\tightlist
\item
  Kirill Krasnov, ``Quantum geometry and thermal radiation from black
  holes'', \emph{Class.\ Quant.\ Grav.\ } \textbf{16} (1999),
   563--578.  Also available as
  \href{https://arxiv.org/abs/gr-qc/9710006}{\texttt{gr-qc/9710006}}.
\end{enumerate}
\noindent
Let me just quote the abstract:

\begin{quote}
A quantum mechanical description of black hole states proposed
recently within the approach known as loop quantum gravity is used to
study the radiation spectrum of a Schwarzschild black hole. We assume
the existence of a Hamiltonian operator causing transitions between
different quantum states of the black hole and use Fermi's golden rule
to find the emission line intensities. Under certain assumptions on the
Hamiltonian we find that, although the emission spectrum consists of
distinct lines, the curve enveloping the spectrum is close to the Planck
thermal distribution with temperature given by the thermodynamical
temperature of the black hole as defined by the derivative of the
entropy with respect to the black hole mass. We discuss possible
implications of this result for the issue of the Immirzi
\(\gamma\)-ambiguity in loop quantum gravity.
\end{quote}

This is interesting, because Bekenstein and Mukhanov have recently noted
that if the area of a quantum black hole is quantized in \emph{evenly
spaced steps}, there will be large deviations from the Planck
distribution of thermal radiation:

\begin{enumerate}
\def\labelenumi{\arabic{enumi})}
\setcounter{enumi}{6}
\tightlist
\item
  Jacob D. Bekenstein and V. F. Mukhanov, ``Spectroscopy of the quantum
  black hole'', \emph{Phys. Lett. B} \textbf{360} (1995), 7--12.  Also available as
  \href{https://arxiv.org/abs/gr-qc/9505012}{\texttt{gr-qc/9505012}}.
\end{enumerate}
\noindent
However, in the loop representation the area is not quantized in evenly
spaced steps: the area \(A\) can be any sum of quantities like
\(8\pi\gamma\sqrt{j(j+1)}\), and such sums become very densely packed
for large \(A\).

Let me conclude with a few technical comments about how Chern--Simons
theory shows up here. For a long time I've been studying the ``ladder of
dimensions'' relating field theories in dimensions 2, 3, and 4, in part
because this gives some clues as to how \(n\)-categories are related to
topological quantum field theory, and in part because it relates quantum
gravity in spacetime dimension 4, which is mysterious, to Chern--Simons
theory in spacetime dimension 3, which is well-understood. It's neat
that one can now use this ladder to study black hole entropy. It's worth
comparing Carlip's calculation of black hole entropy in spacetime
dimension 3 using a \(2\)-dimensional field theory (the
Wess-Zumino-Witten model) on the surface traced out by the black hole
event horizon --- see \protect\hyperlink{week41}{``Week 41''}. Both the
theories we use and those Carlip uses, are all part of the same big
ladder of theories! Something interesting is going on here.

But there's a twist in our calculation which really took me by surprise.
We do not use \(\mathrm{SU}(2)\) Chern--Simons theory on the black hole
surface, we use \(\mathrm{U}(1)\) Chern--Simons theory! The reason is
simple. The boundary conditions we use, which say the black hole surface
is ``marginally outer trapped'', also say that its extrinsic curvature
is zero. Thus the curvature tensor reduces, at the black hole surface,
to the intrinsic curvature. Curvature on a \(3\)-dimensional space is
\(\mathfrak{so}(3)\)-valued, but the intrinsic curvature on the surface
S is \(\mathfrak{so}(2)\)-valued. Since
\(\mathfrak{so}(3) = \mathfrak{su}(2)\), general relativity has a lot to
do with \(\mathrm{SU}(2)\) gauge theory. But since
\(\mathfrak{so}(2) = \mathfrak{u}(1)\), the field theory on the black
hole surface can be thought of as a \(\mathrm{U}(1)\) gauge theory.

(Experts will know that \(\mathrm{U}(1)\) is a subgroup of
\(\mathrm{SU}(2)\) and this is why we look at all values of \(j_z\)
going from \(-j\) to \(j\): we are decomposing representations of
\(\mathrm{SU}(2)\) into representations of this \(\mathrm{U}(1)\)
subgroup.)

Now \(\mathrm{U}(1)\) Chern--Simons theory is a lot less exciting than
\(\mathrm{SU}(2)\) Chern--Simons theory so mathematically this is a bit
of a disappointment. But \(\mathrm{U}(1)\) Chern--Simons theory is not
utterly boring. When we are studying \(\mathrm{U}(1)\) Chern--Simons
theory on a punctured surface, we are studying flat \(\mathrm{U}(1)\)
connections modulo gauge transformations. The space of these is called a
``Jacobian variety''. When we quantize \(\mathrm{U}(1)\) Chern--Simons
theory using geometric quantization, we are looking for holomorphic
sections of a certain line bundle on this Jacobian variety. These are
called ``theta functions''. Theta functions have been intensively
studied by string theorists and number theorists, who use them do all
sorts of wonderful things beyond my ken. All I know about theta
functions can be found in the beginning of the following two books:

\begin{enumerate}
\def\labelenumi{\arabic{enumi})}
\setcounter{enumi}{7}
\item
  Jun-ichi Igusa, \emph{Theta Functions}, Springer, Berlin, 1972.
\item
  David Mumford, \emph{Tata Lectures on Theta}, 3 volumes, Birkhauser,
  Boston, 1983--1991.
\end{enumerate}
\noindent
Theta functions are nice, so it's fun to see them describing states of a
quantum black hole!

\hypertarget{week113}{%
\section{November 26, 1997}\label{week113}}

This week I'd like to talk about ``spin foams''. People have already
thought a lot about using spin networks to describe the quantum geometry
of \(3\)-dimensional space at a given time. This is great for
kinematical aspects of quantum gravity, but not so good for dynamics.
For dynamics, it would be nice to have a description of the quantum
geometry of 4-dimensional \emph{spacetime}. That's where spin foams come
in.

If we use spin networks to describe the geometry of space at the Planck
scale, how might we describe spacetime? Well, space is supposed to be a
kind of slice of spacetime, so let's recall what a spin network is, and
see what it could be a slice of.

First of all, spin network is a graph: a bunch of vertices connected by
edges. What gives a graph when you slice it? Foam! Consider the soap
suds you get while washing the dishes. If we idealize it as a bunch of
2-dimensional surfaces meeting along edges, and imagine intersecting it
with a plane, we see that the result is typically a graph.

Topologists call this sort of space a ``2-dimensional complex''. It's a
generalization of a graph because we can form it by starting with a
bunch of ``vertices'', then connecting these with a bunch of ``edges''
to obtain a graph, and then taking a bunch of \(2\)-dimensional discs or
``faces'' and attaching them along their boundaries to this graph.
Mathematically, there's no reason to stop in dimension 2. Topologists
are interested in complexes of all dimensions. However,
\(2\)-dimensional complexes have been given special attention, thanks to
a number of famous unsolved problems involving them. A great place to
learn about them is:

\begin{enumerate}
\def\labelenumi{\arabic{enumi})}
\tightlist
\item
  C. Hog-Angeloni, W. Metzler, and A. Sieradski, \emph{Two-dimensional
  Homotopy and Combinatorial Group Theory}, London Mathematical Society
  Lecture Note Series \textbf{197}, Cambridge U.\ Press, Cambridge, 1993.
\end{enumerate}

However, a spin network is not \emph{merely} a graph: it's a graph with
edges labelled by irreducible representations of some symmetry group and
vertices labelled by intertwiners. If you don't know what this means,
don't panic! If we take our symmetry group to be \(\mathrm{SU}(2)\),
things simplify tremendously. If we take our graph to have 4 edges
meeting at every vertex, things simplify even more. In this case, all we
need to do is label each vertex and each edge with a number
\(j = 0, 1/2, 1, 3/2,\ldots\) called a ``spin''.

In this special case, we can get our spin network as a slice of a
2-dimensional complex with faces and edges labelled by spins. Such a
thing looks a bit like a foam of soap bubbles with edges and faces
labelled by spins --- hence the term ``spin foam''! More generally, a
spin foam is a \(2\)-dimensional complex with faces labelled by
irreducible representations of some group and edges labelled by
intertwining operators. When we slice a spin foam, each of its faces
gives a spin network edge, and each of its edges gives a spin network
vertex.

Actually, if you've ever looked carefully at soap suds, you'll know that
generically 3 faces meet along each edge. Spin foams like this are
important for quantum gravity in 3 spacetime dimensions. In 4 spacetime
dimensions it seems especially interesting to use spin foams of a
different sort, with 4 faces meeting along each edge. When we slice one
of these, we get a spin network with 4 edges meeting at each vertex.

What's so interesting about spin foams with 4 faces meeting along each
edge? Well, suppose we take a \(4\)-dimensional manifold representing
spacetime and triangulate it --- that is, chop it up into
\(4\)-simplices, the 4-dimensional analogs of tetrahedra. We get a mess
of \(4\)-simplices, which have tetrahedra as faces, which in turn have
triangles as faces.

Now we can form a spin foam with one vertex in the middle of each
4-simplex, one edge intersecting each tetrahedron, and one face
intersecting each triangle. This trick is called ``Poincar\'e duality'':
each \(d\)-dimensional piece of our spin foam intersects a
\((4-d)\)-dimensional piece of our triangulation. Since each tetrahedron
in our triangulated manifold has 4 triangular faces, our spin foam will
dually have 4 faces meeting at each edge. Since each \(4\)-simplex has 5
tetrahedra and 10 triangles, each spin foam vertex will have 5 edges and
10 faces meeting at it.

This seems to be a particularly interesting sort of spin foam for
quantum gravity in 4 dimensions: a spin foam dual to a triangulation of
spacetime. If we slice such a spin foam, we generically get a spin
network dual to a triangulation of space!

I discussed Barbieri's work on such spin networks in
\protect\hyperlink{week110}{``Week 110''}. A spin network like this has
a nice interpretation as a ``3-dimensional quantum geometry'', that is,
a quantum state of the geometry of space. Each spin network edge
labelled by spin \(j\) gives an area proportional to \(\sqrt{j(j+1)}\)
to the triangle it intersects. There's also a formula for the volume of
each tetrahedron, involving the spin on the corresponding spin network
vertex, together with the spins on the 4 spin network edges that meet
there.

It would be nice to have a similar geometrical interpretation of spin
foams dual to triangulations of spacetime. Some recent steps towards
this can be found in the following papers:

\begin{enumerate}
\def\labelenumi{\arabic{enumi})}
\setcounter{enumi}{1}
\item
  John Barrett and Louis Crane, ``Relativistic spin networks and quantum
  gravity'', \emph{Jour. Math. Phys.} \textbf{39} (1998), 3296--3302. Also available as
  \href{https://arxiv.org/abs/gr-qc/9709028}{\texttt{gr-qc/9709028}}.
\item
  John Baez, ``Spin foam models'', \emph{Class. Quant. Grav.}
\textbf{15} (1998), 1827-1858. Also available as
  \href{https://arxiv.org/abs/gr-qc/9709052}{\texttt{gr-qc/9709052}}.
\end{enumerate}

Perhaps I can be forgiven some personal history here. Michael
Reisenberger has been pushing the idea of spin foams (though not the
terminology) for quite a while\ldots{} see for example his paper:

\begin{enumerate}
\def\labelenumi{\arabic{enumi})}
\setcounter{enumi}{3}
\tightlist
\item
  Michael Reisenberger, ``Worldsheet formulations of gauge theories and
  gravity'', available as
  \href{https://arxiv.org/abs/gr-qc/9412035}{\texttt{gr-qc/9412035}}.
\end{enumerate}
\noindent
More recently, Carlo Rovelli and he gave a heuristic derivation of a
spin foam approach to quantum gravity starting with the so-called
canonical quantization approach:

\begin{enumerate}
\def\labelenumi{\arabic{enumi})}
\setcounter{enumi}{4}
\tightlist
\item
  Michael Reisenberger and Carlo Rovelli, ``\,`Sum over surfaces' form
  of loop quantum gravity', \emph{Phys. Rev.} \textbf{D56} (1997),
  3490--3508.  Also available as
  \href{https://arxiv.org/abs/gr-qc/9612035}{\texttt{gr-qc/9612035}}.
\end{enumerate}
\noindent
I started giving talks about spin foams in the spring of this year.
Following the ideas of Reisenberger and Rovelli, I was trying to
persuade everyone to think of spin foams as higher-dimensional analogs
of Feynman diagrams.

Mathematically, a Feynman diagram is just a graph with edges labelled by
representations of some group. But physically, a Feynman diagram
describes a \emph{process} in which a bunch of particles interact. Its
edges correspond to the worldlines traced out by some particles as time
passes, while its vertices represent interactions. Different quantum
field theories use Feynman diagrams with different kinds of vertices.
For any Feynman diagram in our theory, we want to compute a number
called an ``amplitude''. The absolute value squared of this amplitude
gives the probability that the process in question will occur.

We calculate this amplitude by computing a number for each for each edge
and each vertex and multiplying all these numbers together. The numbers
for edges are called ``propagators'' --- they describe the amplitude for
a particle to go from here to there. The numbers for vertices are called
``vertex amplitudes'' --- they describe the amplitude for various
interactions to occur.

Similarly, a spin foam is a \(2\)-dimensional complex with faces
labelled by representations and edges labelled by intertwiners. Each
spin foam face corresponds to the ``worldsheet'' traced out by a spin
network edge as time passes. So, in addition to thinking of a spin foam
as a ``4-dimensional quantum geometry'', we can think of it as a kind of
\emph{process}. The goal of the spin foam approach to quantum gravity is
to compute an amplitude for each spin foam. Following what we know about
Feynman diagrams, it seems reasonable to do it by computing a number for
each spin foam face, edge, and vertex, and then multiplying them all
together.

Quantum gravity is related to a simpler theory called \(BF\) theory,
which has a spin foam formulation known as the Crane--Yetter model ---
see \protect\hyperlink{week36}{``Week 36''},
\protect\hyperlink{week58}{``Week 58''}, and
\protect\hyperlink{week98}{``Week 98''}. There are various clues
suggesting that that the numbers for faces and edges --- the
``propagators'' --- should be computed in quantum gravity just as in the
Crane--Yetter model. The problem is the vertex amplitudes! The vertices
are crucial because these represent the interactions: the places where
something really ``happens''. The number we compute for a vertex
represents the amplitude for the corresponding interaction to occur.
Until we know this, we don't know the dynamics of our theory!

The ``spin foam vertex amplitudes for quantum gravity'' became my holy
grail. Whenever I gave a talk on this stuff I would go around afterwards
asking everyone if they could help me figure them out. I would lay out
all the clues I had and beg for assistance\ldots{} or at least a spark
of inspiration. In March I gave a talk a talk at Penn State proposing a
candidate for these vertex amplitudes --- a candidate I no longer
believe in. Afterwards Carlo Rovelli told me about his attempts to work
out something similar with Louis Crane and Lee Smolin\ldots{} attempts
that never quite got anywhere. We had a crack at it together but it
didn't quite gel. In May I asked John Barrett about the vertex
amplitudes at a conference in Warsaw. He said he couldn't guess them. I
couldn't get \emph{anyone} to guess an answer. In June, at a quantum
gravity workshop in Vienna, I asked Roger Penrose a bunch of questions
about spinors, hoping that this might be the key --- see
\protect\hyperlink{week109}{``Week 109''}. I learned a lot of
interesting stuff, but I didn't find the holy grail.

I kept on thinking. I started getting some promising ideas, and by the
summer I was hard at work on the problem, calculating furiously. I was
also writing a big fat paper about spin foams: the general formalism,
the relation to triangulations, the relationships to category theory,
and so on. I was very happy with it --- but I didn't want to finish it
until I knew the spin foam vertex amplitudes. That would be the icing on
the cake, I thought.

Then one weekend Louis Crane sent me email saying he and John Barrett
had written a paper proposing a model of quantum gravity. Aaargh! Had
they beat me to the holy grail? I frantically wrote up everything I had
while waiting for Monday, when their paper would appear on the preprint
server \texttt{gr-qc}. On Monday I downloaded it and yes, they had
beaten me. It was a skinny little paper and I absorbed it more or less
instantly. They didn't say a word about spin foams --- they were working
dually with triangulations --- but from my viewpoint, what they had done
was to give a formula for the spin foam vertex amplitudes.

Oh well. When you can't beat 'em, join 'em! I finished up my paper,
explaining how their formula fit in with what I'd written already, and
put it on the the preprint server the following weekend.

What did they do to get their formula? Well, the key trick was not to
use \(\mathrm{SU}(2)\) as the symmetry group, but instead use
\(\mathrm{SU}(2)\times\mathrm{SU}(2)\). This is the double cover of
\(\mathrm{SO}(4)\), the rotation group in 4 dimensions. Following the
idea behind Ashtekar's new variables for general relativity, I was only
using the ``left-handed half'' of this group, that is, one of the
\(\mathrm{SU}(2)\) factors. But the geometry of the \(4\)-simplex, and
its relation to quantum theory, is in some ways more easily understood
using the full \(\mathrm{SU}(2)\times\mathrm{SU}(2)\) symmetry group.

Not surprisingly, an irreducible representation of
\(\mathrm{SU}(2)\times\mathrm{SU}(2)\) is described by a pair of spins
\((j,k)\). The reason is that we can take the spin-\(j\) representation
of the ``left-handed'' \(\mathrm{SU}(2)\) and the spin-\(k\)
representation of the ``right-handed'' \(\mathrm{SU}(2)\) and tensor
them to get an irreducible representation of
\(\mathrm{SU}(2)\times\mathrm{SU}(2)\). If we use
\(\mathrm{SU}(2)\times\mathrm{SU}(2)\) as our group, our spin foams dual
to triangulations will thus have every face and every edge labelled by a
\emph{pair} of spins. However, Barrett and Crane's work suggests that
the only spin foams with nonzero amplitudes are those for which both
spins labelling a face or edge are equal! Thus in a way we are back down
to \(\mathrm{SU}(2)\) --- but we think of it all a bit differently.

I'm tempted to go into detail and explain exactly how the model works,
because it involves a lot of beautiful geometry. But it takes a while,
so I won't. First you need to really grok the phase space of all
possible geometries of the \(4\)-simplex. Then you need to quantize this
phase space, obtaining the ``Hilbert space of a quantum \(4\)-simplex''.
Then you need to note that there's a special linear functional on this
Hilbert space, called the ``Penrose evaluation'' - see
\protect\hyperlink{week110}{``Week 110''}. Putting all this together
gives the vertex amplitudes for quantum gravity\ldots{} we hope.

Anyway, back to my little personal story\ldots.

Though I'd been working on my paper before Barrett and Crane started,
and they finished before me, Michael Reisenberger started one even
earlier and finished even later! Indeed, he has been working on a spin
foam model of quantum gravity for several years now --- see
\protect\hyperlink{week86}{``Week 86''}. He took a purely left-handed
\(\mathrm{SU}(2)\) approach, a bit different what I'd been trying, but
closely related. He told lots of people about it, but unfortunately he's
very slow to publish.

When I heard Barrett and Crane were about to come out with their paper,
I emailed Reisenberger warning him of this. He doesn't like being
scooped any more than I do. Unfortunately I only had his old email
address in Canada, and now he's down in Uruguay, so he never got that
email. Thus Barrett and Crane's paper, followed by mine, came as a a big
shock to him! Luckily, this motivated him to hurry and come out with a
version of his paper:

\begin{enumerate}
\def\labelenumi{\arabic{enumi})}
\setcounter{enumi}{5}
\tightlist
\item
  Michael Reisenberger, ``A lattice worldsheet sum for 4-d Euclidean
  general relativity'', available as
  \href{https://arxiv.org/abs/gr-qc/9711052}{\texttt{gr-qc/9711052}}.
\end{enumerate}

Let me quote the abstract:

\begin{quote}
``A lattice model for four dimensional Euclidean quantum general
relativity is proposed for a simplicial spacetime. It is shown how this
model can be expressed in terms of a sum over worldsheets of spin
networks, and an interpretation of these worldsheets as spacetime
geometries is given, based on the geometry defined by spin networks in
canonical loop quantized GR. The spacetime geometry has a Planck scale
discreteness which arises''naturally" from the discrete spectrum of
spins of \(\mathrm{SU}(2)\) representations (and not from the use of a
spacetime lattice). The lattice model of the dynamics is a formal
quantization of the classical lattice model of {[}Reisenberger's paper
``A left-handed simplicial action for euclidean general relativity''{]},
which reproduces, in a continuum limit, Euclidean general relativity."
\end{quote}

To wrap up my little history, let me say what's been happening lately.
There is still a lot of puzzlement and mystery concerning these spin
foam models\ldots{} it's far from clear that they really work as hoped
for. We may be doing things a little bit wrong, or we may be on a
completely wrong track. The phase space of the \(4\)-simplex involves
some tricky constraint equations which could kill us if we're not
dealing with them right. Barbieri has suggested a modified version of
Barrett and Crane's approach which may overcome some problems with the
constraints:

\begin{enumerate}
\def\labelenumi{\arabic{enumi})}
\setcounter{enumi}{6}
\tightlist
\item
  Andrea Barbieri, ``Space of the vertices of relativistic spin
  networks'', available as
  \href{https://arxiv.org/abs/gr-qc/9709076}{\texttt{gr-qc/9709076}}.
\end{enumerate}
\noindent
John Barrett visited me last week and we made some progress on this
issue, but it's still very touchy.

Also, all the work cited above deals with so-called ``Euclidean''
quantum gravity --- that's why it uses the double cover of the rotation
group \(\mathrm{SO}(4)\). For ``Lorentzian'' quantum gravity we'd need
instead to use the double cover of the Lorentz group
\(\mathrm{SO}(3,1)\). This group is isomorphic to
\(\mathrm{SL}(2,\mathbb{C})\). As explained in
\protect\hyperlink{week109}{``Week 109''}, the finite-dimensional
irreducible representations of \(\mathrm{SL}(2,\mathbb{C})\) are also
described by pairs of spins, so the Lorentzian theory should be similar
to the Euclidean theory. However, most work so far has dealt with the
Euclidean case; this needs to be addressed.

Finally, Crane has written a bit more about the geometrical significance
of his work with Barrett:

\begin{enumerate}
\def\labelenumi{\arabic{enumi})}
\setcounter{enumi}{7}
\tightlist
\item
  Louis Crane, ``On the interpretation of relativistic spin networks and
  the balanced state sum'', available as
  \href{https://arxiv.org/abs/gr-qc/9710108}{\texttt{gr-qc/9710108}}.
\end{enumerate}

Next week I'll talk about other developments in the loop representation
of quantum gravity, some arising from Thiemann's work on the Hamiltonian
constraint. After that, I think I want to talk about something
completely different, like homotopy theory. Lately I've been trying to
make ``This Week's Finds'' very elementary and readable --- relatively
speaking, I mean --- but I'm getting in the mood for writing in a more
technical and incomprehensible manner, and homotopy theory is an ideal
subject for that sort of writing.

\hypertarget{week114}{%
\section{January 12, 1998}\label{week114}}

Classes have started! But I just flew back yesterday from the Joint
Mathematics Meetings in Baltimore --- the big annual conference
organized by the AMS, the MAA, SIAM, and other societies. Over 4000
mathematicians could be seen wandering in clumps about the glitzy harbor
area and surrounding crime-ridden slums, arguing about abstractions,
largely oblivious to the world around them. Everyone ate the obligatory
crab cakes for which Baltimore is justly famous. Some of us drank a bit
too much beer, too.

Witten gave a plenary talk on ``M-theory'', which was great fun even
though he didn't actually say what M-theory is. Steve Sawin and I ran a
session on quantum gravity and low-dimensional topology, so I'll say a
bit about what went on there. There was also a nice session on homotopy
theory in honor of J. Michael Boardman. I'll talk about that and various
other things next week.

A lot of the buzz in our session concerned the new ``spin foam''
approach to quantum gravity which I discussed in
\protect\hyperlink{week113}{``Week 113''}. The big questions are: how do
you test this approach without impractical computer simulations? Lee
Smolin's paper below suggests one way. Should you only sum over spin
foams that are dual to a particular triangulation of spacetime, or
should you sum over all spin foams that fit in a particular
\(4\)-dimensional spacetime manifold, or should you sum over \emph{all}
spin foams? There was a lot of argument about this. In addition to the
question of what is physically appropriate, there's the mathematical
problem of avoiding divergent infinite sums. Perhaps the sum required to
answer any truly physical question only involves finitely many spin
foams --- that's what I hope. Finally, should the time evolution
operators constructed using spin foams be thought of as describing true
time evolution, or merely the projection onto the kernel of the
Hamiltonian constraint? While it sounds a bit technical, this question
is crucial for the interpretation of the theory; it's part of what they
call ``the problem of time''.

Carlo Rovelli spoke about how spin foams arise in canonical quantum
gravity, while John Barrett and Louis Crane discussed them in the
context of discretized path integrals for quantum gravity, also known as
state sum models. As in the more traditional ``Regge calculus''
approach, these models start by chopping spacetime into simplices. The
biggest difference is that now \emph{areas of triangles} play a more
important role than lengths of edges. But Barrett, Crane and others are
starting to explore the relationships:

\begin{enumerate}
\def\labelenumi{\arabic{enumi})}
\item
  John W. Barrett, Martin Rocek and Ruth M. Williams, ``A note on area
  variables in Regge calculus'', \emph{Class.\ Quant.\ Grav.} \textbf{16} 
  (1999), 1373--1376.  Also available as
  \href{https://arxiv.org/abs/gr-qc/9710056}{\texttt{gr-qc/9710056}}.
\item
  Jarmo Makela, ``Variation of area variables in Regge calculus'', 
  \emph{Class.\ Quant.\ Grav.} \textbf{17} (2000), 4991--4998.  Also 
   available as
  \href{https://arxiv.org/abs/gr-qc/9801022}{\texttt{gr-qc/9801022}}.
\end{enumerate}

Also, there's been some progress on extracting Einstein's equation for
general relativity as a classical limit of the Barrett--Crane state sum
model. Let me quote the abstract of this paper:

\begin{enumerate}
\def\labelenumi{\arabic{enumi})}
\setcounter{enumi}{2}
\tightlist
\item
  Louis Crane and David N. Yetter, ``On the classical limit of the
  balanced state sum'', available as
  \href{https://arxiv.org/abs/gr-qc/9712087}{\texttt{gr-qc/9712087}}.
\end{enumerate}

\begin{quote}
The purpose of this note is to make several advances in the
interpretation of the balanced state sum model by Barrett and Crane in
{\rm \href{https://arxiv.org/abs/gr-qc/9709028}{\texttt{gr-qc/9709028}}} as a
quantum theory of gravity. First, we outline a shortcoming of the
definition of the model pointed out to us by Barrett and Baez in private
communication, and explain how to correct it. Second, we show that the
classical limit of our state sum reproduces the Einstein-Hilbert
lagrangian whenever the term in the state sum to which it is applied has
a geometrical interpretation. Next we outline a program to demonstrate
that the classical limit of the state sum is in fact dominated by terms
with geometrical meaning. This uses in an essential way the alteration
we have made to the model in order to fix the shortcoming discussed in
the first section. Finally, we make a brief discussion of the Minkowski
signature version of the model.
\end{quote}

Lee Smolin talked about his ideas for relating spin foam models to
string theory. He has a new paper on this, so I'll just quote the
abstract:

\begin{enumerate}
\def\labelenumi{\arabic{enumi})}
\setcounter{enumi}{3}
\tightlist
\item
  Lee Smolin, ``Strings as perturbations of evolving spin-networks'',    
  \emph{Nucl.\ Phys.\  Proc.\  Suppl.} \textbf{88} (2000), 103--113.  
  Also available as
  \href{https://arxiv.org/abs/hep-th/9801022}{\texttt{hep-th/9801022}}.
\end{enumerate}

\begin{quote}
A connection between non-perturbative formulations of quantum gravity
and perturbative string theory is exhibited, based on a formulation of
the non-perturbative dynamics due to Markopoulou. In this formulation
the dynamics of spin network states and their generalizations is
described in terms of histories which have discrete analogues of the
causal structure and many fingered time of Lorentzian spacetimes.
Perturbations of these histories turn out to be described in terms of
spin systems defined on \(2\)-dimensional timelike surfaces embedded in
the discrete spacetime. When the history has a classical limit which is
Minkowski spacetime, the action of the perturbation theory is given to
leading order by the spacetime area of the surface, as in bosonic string
theory. This map between a non-perturbative formulation of quantum
gravity and a 1+1 dimensional theory generalizes to a large class of
theories in which the group \(\mathrm{SU}(2)\) is extended to any
quantum group or supergroup. It is argued that a necessary condition for
the non-perturbative theory to have a good classical limit is that the
resulting 1+1 dimensional theory defines a consistent and stable
perturbative string theory.
\end{quote}

Fotini Markopolou spoke about her recent work with Smolin on formulating
spin foam models in a manifestly local, causal way.

\begin{enumerate}
\def\labelenumi{\arabic{enumi})}
\setcounter{enumi}{4}
\tightlist
\item
  Fotini Markopoulou and Lee Smolin, ``Quantum geometry with intrinsic
  local causality'', \emph{Phys.\ Rev.\ D} \textbf{58} (1998), 084032.  
  Also available as
  \href{https://arxiv.org/abs/gr-qc/9712067}{\texttt{gr-qc/9712067}}.
\end{enumerate}

\begin{quote}
The space of states and operators for a large class of background
independent theories of quantum spacetime dynamics is defined. The
\(\mathrm{SU}(2)\) spin networks of quantum general relativity are
replaced by labelled compact two-dimensional surfaces. The space of
states of the theory is the direct sum of the spaces of invariant
tensors of a quantum group \(G_q\) over all compact (finite genus)
oriented 2-surfaces. The dynamics is background independent and locally
causal. The dynamics constructs histories with discrete features of
spacetime geometry such as causal structure and multifingered time. For
\(\mathrm{SU}(2)\) the theory satisfies the Bekenstein bound and the
holographic hypothesis is recast in this formalism.
\end{quote}

The main technical idea in this paper is to work with ``thickened'' or
``framed'' spin networks, which amounts to replacing graphs by solid
handlebodies. One expects this ``framing'' business to be important for
quantum gravity with nonzero cosmological constant. This framing
business also appears in the \(q\)-deformed version of Barrett and
Crane's model and in my ``abstract'' version of their model, which
assumes no background spacetime manifold. Markopoulou and Smolin don't
specify a choice of dynamics; instead, they describe a \emph{class} of
theories which has my model as a special case, though their approach to
causality is better suited to Lorentzian theories, while mine is
Euclidean.

As I've often noted, spin foams are about spacetime geometry, or
dynamics, while spin networks are a way of describing the geometry of
space, or kinematics. Kinematics is always easier than dynamics, so the
spin network approach to the quantum geometry of space has been much
better worked out than the new spin foam stuff. Abhay Ashtekar gave an
overview of these kinematical issues in his talk on ``quantum Riemannian
geometry'', and Kirill Krasnov described how our understanding of these
already allows us to compute the entropy of black holes (see
\protect\hyperlink{week112}{``Week 112''}). Here it's worth mentioning
that the second part of Ashtekar's paper with Jerzy Lewandowski is
finally out:

\begin{enumerate}
\def\labelenumi{\arabic{enumi})}
\setcounter{enumi}{5}
\tightlist
\item
  Abhay Ashtekar and Jerzy Lewandowski, ``Quantum theory of geometry II:
  volume operators'', \emph{Adv.\ Theor.\ Math.\ Phys.\ } \textbf{1} (1998), 
  388--429.  Also available as
  \href{https://arxiv.org/abs/gr-qc/9711031}{\texttt{gr-qc/9711031}}.
\end{enumerate}

\begin{quote}
A functional calculus on the space of (generalized) connections was
recently introduced without any reference to a background metric. It is
used to continue the exploration of the quantum Riemannian geometry.
Operators corresponding to volume of three-dimensional regions are
regularized rigorously. It is shown that there are two natural
regularization schemes, each of which leads to a well-defined operator.
Both operators can be completely specified by giving their action on
states labelled by graphs. The two final results are closely related but
differ from one another in that one of the operators is sensitive to the
differential structure of graphs at their vertices while the second is
sensitive only to the topological characteristics. (The second operator
was first introduced by Rovelli and Smolin and De Pietri and Rovelli
using a somewhat different framework.) The difference between the two
operators can be attributed directly to the standard quantization
ambiguity. Underlying assumptions and subtleties of regularization
procedures are discussed in detail in both cases because volume
operators play an important role in the current discussions of quantum
dynamics.
\end{quote}

Before spin foam ideas came along, the basic strategy in the loop
representation of quantum gravity was to start with general relativity
on a smooth manifold and try to quantize it using the ``canonical
quantization'' approach. Here the most important and difficult thing is
to implement the ``Hamiltonian constraint'' as an operator on the
Hilbert space of kinematical states, so you can write down the
Wheeler-deWitt equation, which is, quite roughly speaking, the quantum
gravity analog of Schrodinger's equation. (For a summary of this
approach, try \protect\hyperlink{week43}{``Week 43''}.)

The most careful attempt to do this so far is the work of Thiemann:

\begin{enumerate}
\def\labelenumi{\arabic{enumi})}
\setcounter{enumi}{6}
\item
  Thomas Thiemann, ``Quantum spin dynamics (QSD)'', \emph{Class.\
  Quant.\ Grav.} \textbf{15} (1998), 839--873.  Also available
  as \href{https://arxiv.org/abs/gr-qc/9606089}{\texttt{gr-qc/9606089}}.

  ``Quantum spin dynamics (QSD) II'', \emph{Class.\ Quant.\ Grav.}
    \textbf{15} (1998), 875--905.    Also available as
  \href{https://arxiv.org/abs/gr-qc/9606090}{\texttt{gr-qc/9606090}}.

  ``QSD III: Quantum constraint algebra and physical scalar product in
  quantum general relativity'', \emph{Class.\ Quant.\  Grav.} \textbf{15}
  (1998), 1207--1247.  Also available as
  \href{https://arxiv.org/abs/gr-qc/9705017}{\texttt{gr-qc/9705017}}.

  ``QSD IV: 2+1 Euclidean quantum gravity as a model to test 3+1
  Lorentzian quantum gravity'', \emph{Class.\ Quant.\ Grav.} \textbf{15} 
   (1998), 1249--1280.  Also available as
  \href{https://arxiv.org/abs/gr-qc/9705018}{\texttt{gr-qc/9705018}}.

  ``QSD V: Quantum gravity as the natural regulator of matter quantum
  field theories'', \emph{Class.\ Quant.\ Grav.} \textbf{15} (1998), 
  1281--1314.  Also available as
  \href{https://arxiv.org/abs/gr-qc/9705019}{\texttt{gr-qc/9705019}}.

  ``QSD VI: Quantum Poincar\'e algebra and a quantum positivity of energy
  theorem for canonical quantum gravity'', \emph{Class.\ Quant.\ Grav.}
  \textbf{15} (1998), 1463--1485.   Also available as
  \href{https://arxiv.org/abs/gr-qc/9705020}{\texttt{gr-qc/9705020}}

  ``Kinematical Hilbert spaces for fermionic and Higgs quantum field
  theories'', \emph{Class.\ Quant.\ Grav.} \textbf{15} (1998), 
    1487--1512.  Also available as
  \href{https://arxiv.org/abs/gr-qc/9705021}{\texttt{gr-qc/9705021}}
\end{enumerate}

If everything worked as smoothly as possible, the Hamiltonian constraint
would satisfy nice commutation relations with the other constraints of
the theory, giving a representation of something called the ``Dirac
algebra''. However, as Don Marolf explained in his talk, this doesn't
really happen, at least in a large class of approaches including
Thiemann's:

\begin{enumerate}
\def\labelenumi{\arabic{enumi})}
\setcounter{enumi}{7}
\item
  Jerzy Lewandowski and Donald Marolf, ``Loop constraints: A habitat and
  their algebra'', \emph{Int.\ J.\ Mod.\ Phys.\ D} \textbf{7} (1998), 
  299--330.  Also available as
  \href{https://arxiv.org/abs/gr-qc/9710016}{\texttt{gr-qc/9710016}}.
\item
  Rodolfo Gambini, Jerzy Lewandowski, Donald Marolf, and Jorge Pullin,
  ``On the consistency of the constraint algebra in spin network quantum
  gravity'', \emph{Int.\ J.\ Mod.\ Phys.\ D} \textbf{7} (1998),  97--109.  
   Also available as
  \href{https://arxiv.org/abs/gr-qc/9710018}{\texttt{gr-qc/9710018}}.
\end{enumerate}

This is very worrisome\ldots{} as everything concerning quantum gravity
always is. Personally these results make me want to spend less time on
the Hamiltonian constraint, especially to the extent that it assumes a
the old picture of spacetime as a smooth manifold, and more time on
approaches that start with a discrete picture of spacetime. However, the
only way to make serious progress is for different people to push on
different fronts simultaneously.

There were a lot of other interesting talks, but since I'm concentrating
on quantum gravity here I won't describe the ones that were mainly about
topology. I'll wrap up by mentioning Steve Carlip's talk on spacetime
foam. He gave a nice illustration to how hard it is to ``sum over
topologies'' by arguing that this sum diverges for negative values of
the cosmological constant. He has a paper out on this:

\begin{enumerate}
\def\labelenumi{\arabic{enumi})}
\setcounter{enumi}{9}
\tightlist
\item
  Steven Carlip, ``Spacetime foam and the cosmological constant'',
  \emph{Phys. Rev.~Lett.} \textbf{79} (1997) 4071--4074.  Also
  available as
  \href{https://arxiv.org/abs/gr-qc/9708026}{\texttt{gr-qc/9708026}}.
\end{enumerate}

Again, I'll quote the abstract:

\begin{quote}
In the saddle point approximation, the Euclidean path integral for
quantum gravity closely resembles a thermodynamic partition function,
with the cosmological constant \(\Lambda\) playing the role of
temperature and the `density of topologies' acting as an effective
density of states. For \(\Lambda < 0\), the density of topologies grows
superexponentially, and the sum over topologies diverges. In
thermodynamics, such a divergence can signal the existence of a maximum
temperature. The same may be true in quantum gravity: the effective
cosmological constant may be driven to zero by a rapid rise in the
density of topologies.
\end{quote}

\hypertarget{week115}{%
\section{February 1, 1998}\label{week115}}

These days I've been trying to learn more homotopy theory. James Dolan
got me interested in it by explaining how it offers many important clues
to \(n\)-category theory. Ever since, we've been trying to understand
what the homotopy theorists have been up to for the last few decades.
Since trying to explain something is often the best way to learn it,
I'll talk about this stuff for several Weeks to come.

Before plunging in, though, I'd like mention yet another novel by Greg
Egan:

\begin{enumerate}
\def\labelenumi{\arabic{enumi})}
\tightlist
\item
  Greg Egan, \emph{Diaspora}, Orion Books, 1997.
\end{enumerate}

The main character of this book, Yatima, is a piece of software\ldots{}
and a mathematician. The tale begins in 2975 with ver birth as an
``orphan'', a citizen of the polis born of no parents, its mind seed
chosen randomly by the conceptory. Yatima learns mathematics in a
virtual landscape called the Truth Mines. To quote (with a few small
modifications):

\begin{quote}
The luminous object buried in the cavern floor broadcast the definition
of a topological space: a set of points, grouped into `open subsets'
which specified how the points were connected to one another --- without
appealing to notions like `distance' or `dimension'. Short of a raw set
with no structure at all, this was about as basic as you could get: the
common ancestor of virtually every entity worth of the name `space',
however exotic. A single tunnel led into the cavern, providing the link
to the necessary prior concepts, and half a dozen tunnels led out,
slanting gently `down' into the bedrock, pursuing various implications
of the definition. Suppose \(T\) is a topological space\ldots{} then
what follows? These routes were paved with small gemstones, each one
broadcasting an intermediate result on the way to a theorem.

Every tunnel in the Mines was built from the steps of a watertight
proof; every theorem, however deeply buried, could be traced back to
every one of its assumptions. And to pin down exactly what was meant by
a `proof', every field of mathematics used its own collection of formal
systems: sets of axioms, definitions, and rules of deduction, along with
the specialised vocabulary needed to state theorems and conjectures
precisely.

{[}\ldots.{]}

The library was full of the ways past miners had fleshed out the
theorems, and Yatima could have had those details grafted in alongside
the raw data, granting ver the archived understanding of thousands of
Konishi citizens who'd travelled this route before. The right
mind-grafts would have enabled ver effortlessly to catch up with all the
living miners who were pushing the coal face ever deeper in their own
inspired directions\ldots{} at the cost of making ver, mathematically
speaking, little more than a patchwork clone of them, capable only of
following in their shadow.

If ve ever wanted to be a miner in vis own right --- making and testing
vis own conjectures at the coal face, like Gauss and Euler, Riemann and
Levi-Civita, deRham and Cartan, Radiya and Blanca - then Yatima knew
there were no shortcuts, no alternatives to exploring the Mines first
hand. Ve couldn't hope to strike out in a fresh direction, a route no
one had ever chosen before, without a new take on the old results. Only
once ve'd constructed vis own map of the Mines --- idiosyncratically
crumpled and stained, adorned and annotated like no one else's --- could
ve begin to guess where the next rich vein of undiscovered truths lay
buried.
\end{quote}

The tale ends in a universe 267,904,176,383,054 duality transformations
away from ours, at the end of a long quest. What does Yatima do then?
Keep studying math! ``It would be a long, hard journey to the coal face,
but this time there would be no distractions.''

I won't give away any more of the plot. Suffice it to say that this is
hard science fiction --- readers in search of carefully drawn characters
may be disappointed, but those who enjoy virtual reality, wormholes, and
philosophy should have a rollicking good ride. I must admit to being
biased in its favor, since it refers to a textbook I wrote. A science
fiction writer who actually knows the Gauss-Bonnet theorem! We should be
very grateful.

Okay, enough fun --- it's time for homotopy theory. Actually homotopy
theory is \emph{tremendously} fun, but it takes quite a bit of
persistence to come anywhere close to the coal face. The original
problems motivating the subject are easy to state. Let's call a
topological space simply a ``space'', and call a continuous function
between these simply a ``map''. Two maps \(f,g\colon X\to Y\) are
``homotopic'' if one can be continuously deformed to the other, or in
other words, if there is a ``homotopy'' between them: a continuous
function \(F\colon[0,1]\times X\to Y\) with \[F(0,x) = f(x)\] and
\[F(1,x) = g(x).\] Also, two spaces \(X\) and \(Y\) are ``homotopy
equivalent'' if there are functions \(f\colon X\to Y\) and
\(g\colon Y\to X\) for which \(fg\) and \(gf\) are homotopic to the
identity. Thus, for example, a circle, an annulus, and a solid torus are
all homotopy equivalent.

Homotopy theorists want to classify spaces up to homotopy equivalence.
And given two spaces \(X\) and \(Y\), they want to understand the set
\([X,Y]\) of homotopy classes of maps from \(X\) to \(Y\). However,
these are very hard problems! To solve them, one needs high-powered
machinery.

There are two roughly two sides to homotopy theory: building machines,
and using them to do computations. Of course these are fundamentally
inseparable, but people usually tend to prefer either one or the other
activity. Since I am a mathematical physicist, always on the lookout for
more tools for my own work, I'm more interested in the nice shiny
machines homotopy theorists have built than in the terrifying uses to
which they are put.

What follows will strongly reflect this bias: I'll concentrate on a
bunch of elegant concepts lying on the interface between homotopy theory
and category theory. This realm could be called ``homotopical algebra''.
Ideas from this realm can be applied, not only to topology, but to many
other realms. Indeed, two of its most famous practitioners, James
Stasheff and Graeme Segal, have spent the last decade or so using it in
string theory! I'll eventually try to say a bit about how that works,
too.

Okay\ldots. now I'll start listing concepts and tools, starting with the
more fundamental ones and then working my way up. This will probably
only make sense if you've got plenty of that commodity known as
``mathematical sophistication''. So put on some Coltrane, make yourself
a cafe macchiato, kick back, and read on. If at any point you feel a
certain lack of sophistication, you might want to reread ``The Tale of
\(n\)-Categories'', starting with \protect\hyperlink{week73}{``Week
73''}, where a bunch of the basic terms are defined.

\begin{center}\rule{0.5\linewidth}{0.5pt}\end{center}

\hypertarget{homotopy_A}{\textbf{A.}}
\emph{Presheaf categories.} Given a category
\(\mathcal{C}\), a ``presheaf'' on \(\mathcal{C}\) is a contravariant
functor \(F\colon\mathcal{C}\to\mathsf{Set}\). The original example of
this is where \(\mathcal{C}\) is the category whose objects are open
subsets of a topological space \(X\), with a single morphism
\(f\colon U\to V\) whenever the open set \(U\) is contained in the open
set \(V\). For example, there is the presheaf of continuous real-valued
functions, for which \(F(U)\) is the set of all continuous real
functions on \(U\), and for any inclusion \(f\colon U\to V\),
\(F(f)\colon F(V)\to F(U)\) is the ``restriction'' map which assigns to
any continuous function on \(V\) its restriction to \(U\). This is a
great way of studying functions in all neighborhoods of \(X\) at once.

However, I'm bringing up this subject for a different reason, related to
a different kind of example. Suppose that \(\mathcal{C}\) is a category
whose objects are ``shapes'' of some kind, with morphisms
\(f\colon x\to y\) corresponding to ways the shape \(x\) can be included
as a ``piece'' of the shape \(y\). Then a presheaf on \(\mathcal{C}\)
can be thought of as a geometrical structure built by gluing together
these shapes along their common pieces.

For example, suppose we want to describe directed graphs as presheaves.
A directed graph is a bunch of vertices and edges, where the edges have
a direction specified. Since they are made of two ``shapes'', the vertex
and the edge, we'll cook up a little category \(\mathcal{C}\) with two
object, \(V\) and \(E\). There are two ways a vertex can be included as
a piece of an edge, either as its ``source'' or its ``target''. Our
category \(\mathcal{C}\), therefore, has two morphisms,
\(S\colon V\to E\) and \(T\colon V\to E\). These are the only morphisms
except for identity morphisms --- which correspond to how the edge is
part of itself, and the vertex is part of itself! Omitting identity
morphisms, our little category \(\mathcal{C}\) looks like this: \[
  \begin{tikzcd}
    V \rar[bend left=30,"S"] \rar[bend right=30,"T"]
    & E
  \end{tikzcd}
\] Now let's work out what a presheaf on \(\mathcal{C}\) is. It's a
contravariant functor \(F\colon\mathcal{C}\to\mathsf{Set}\). What does
this amount to? Well, it amounts to a set \(F(V)\) called the ``set of
vertices'', a set \(F(E)\) called the ``set of edges'', a function
\(F(S)\colon F(E)\to F(V)\) assigning to each edge its source, and a
function \(F(T)\colon F(E)\to F(V)\) assigning to each edge its target.
That's just a directed graph!

Note the role played by contravariance here: if a little shape \(V\) is
included as a piece of a big shape \(E\), our category gets a morphism
\(S\colon V\to E\), and then in our presheaf we get a function
\(F(S)\colon F(E)\to F(V)\) going the \emph{other way}, which describes
how each big shape has a bunch of little shapes as pieces.

Given any category \(\mathcal{C}\) there is actually a \emph{category}
of presheaves on \(\mathcal{C}\). Given presheaves
\(F,G\colon\mathcal{C}\to\mathsf{Set}\), a morphism \(M\) from \(F\) to
\(G\) is just a natural transformation \(M\colon F\Rightarrow G\). This
is beautifully efficient way of saying quite a lot. For example, if
\(\mathcal{C}\) is the little category described above, so that \(F\)
and \(G\) are directed graphs, a natural transformation
\(M\colon F\Rightarrow G\) is the same as:

\begin{itemize}
\tightlist
\item
  a map \(M(V)\) sending each vertex of the graph F to a vertex of the
  graph \(G\),
\end{itemize}

and

\begin{itemize}
\tightlist
\item
  a map \(M(E)\) sending each edge of the graph \(F\) to a edge of the
  graph \(G\),
\end{itemize}

such that

\begin{itemize}
\tightlist
\item
  \(M(V)\) of the source of any edge \(e\) of \(F\) equals the source of
  \(M(E)\) of \(e\),
\end{itemize}

and

\begin{itemize}
\tightlist
\item
  \(M(V)\) of the target of any edge \(e\) of \(F\) equals the target of
  \(M(E)\) of \(e\).
\end{itemize}
\noindent
Whew! Easier just to say \(M\) is a natural transformation between
functors!

For more on presheaves, try:

\begin{enumerate}
\def\labelenumi{\arabic{enumi})}
\setcounter{enumi}{1}
\tightlist
\item
  Saunders Mac Lane and Ieke Moerdijk, \emph{Sheaves in Geometry and
  Logic: a First Introduction to Topos Theory}, Springer, Berlin, 1992.
\end{enumerate}

\begin{center}\rule{0.5\linewidth}{0.5pt}\end{center}

\hypertarget{homotopy_B}{\textbf{B.}}
\emph{The category of simplices, $\Delta$.} This is a very important
example of a category whose objects are shapes --- namely, simplices ---
and whose morphisms correspond to the ways one shape is a piece of
another. The objects of \(\Delta\) are called \(1, 2, 3, \ldots\),
corresponding to the simplex with 1 vertex (the point), the simplex with
2 vertices (the interval), the simplex with 3 vertices (the triangle),
and so on. There are a bunch of ways for an lower-dimensional simplex to
be a face of a higher- dimensional simplex, which give morphisms in
\(\Delta\). More subtly, there are also a bunch of ways to map a
higher-dimensional simplex down into a lower-dimensional simplex, called
``degeneracies''. For example, we can map a tetrahedron down into a
triangle in a way that carries the vertices \(\{0,1,2,3\}\) of the
tetrahedron into the vertices \(\{0,1,2\}\) of the triangle as follows:
\[
  \begin{aligned}
    0&\mapsto0
  \\1&\mapsto0
  \\2&\mapsto1
  \\3&\mapsto2
  \end{aligned}
\] These degeneracies also give morphisms in \(\Delta\).

We could list all the morphisms and the rules for composing them
explicitly, but there is a much slicker way to describe them. Let's use
the old trick of thinking of the natural number \(n\) as being the
totally ordered \(n\)-element set \(\{0,1,2,\ldots,n-1\}\) of all
natural numbers less than \(n\). Thus for example we think of the object
\(4\) in \(\Delta\), the tetrahedron, as the totally ordered set
\(\{0,1,2,3\}\). These correspond to the 4 vertices of the tetrahedron.
Then the morphisms in \(\Delta\) are just all order-preserving maps
between these totally ordered sets. So for example there is a morphism
\(f\colon\{0,1,2,3\}\to\{0,1,2\}\) given by the order-preserving map
with \[
  \begin{aligned}
    f(0)&=0
  \\f(1)&=0
  \\f(2)&=1
  \\f(3)&=2
  \end{aligned}
\] The rule for composing morphisms is obvious: just compose the maps!
Slick, eh?

We can be slicker if we are willing to work with a category
\emph{equivalent} to \(\Delta\) (in the technical sense described in
\protect\hyperlink{week76}{``Week 76''}), namely, the category of
\emph{all} nonempty totally ordered sets, with order-preserving maps as
morphisms. This has a lot more objects than just \(\{0\}\), \(\{0,1\}\),
\(\{0,1,2\}\), etc., but all of its objects are isomorphic to one of
these. In category theory, equivalent categories are the same for all
practical purposes --- so we brazenly call this category \(\Delta\),
too. If we do so, we have following \emph{incredibly} slick description
of the category of simplices: it's just the category of finite nonempty
totally ordered sets!

If you are a true mathematician, you will wonder ``why not use the empty
set, too?'' Generally it's bad to leave out the empty set. It may seem
like ``nothing'', but ``nothing'' is usually very important. Here it
corresponds to the ``empty simplex'', with no vertices! Topologists
often leave this one out, but sometimes regret it later and put it back
in (the buzzword is ``augmentation''). True category theorists, like Mac
Lane, never leave it out. They define \(\Delta\) to be the category of
\emph{all} totally ordered finite sets. For a beautiful introduction to
this approach, try:

\begin{enumerate}
\def\labelenumi{\arabic{enumi})}
\setcounter{enumi}{2}
\tightlist
\item
  Saunders Mac Lane, \emph{Categories for the Working Mathematician},
  Springer, Berlin, 1988.
\end{enumerate}

\begin{center}\rule{0.5\linewidth}{0.5pt}\end{center}

\hypertarget{homotopy_C}{\textbf{C.}} 
\emph{Simplicial sets.} Now we put together the previous two
ideas: a ``simplicial set'' is a presheaf on the category of simplices!
In other words, it's a contravariant functor
\(F\colon\Delta\to\mathsf{Set}\). Geometrically, it's basically just a
bunch of simplices stuck together along their faces in an arbitrary way.
We can think of it as a kind of purely combinatorial version of a
``space''. That's one reason simplicial sets are so popular in topology:
they let us study spaces in a truly elegant algebraic context. We can
define all the things topologists love --- homology groups, homotopy
groups (see \protect\hyperlink{week102}{``Week 102''}), and so on ---
while never soiling our hands with open sets, continuous functions and
the like. To see how it's done, try:

\begin{enumerate}
\def\labelenumi{\arabic{enumi})}
\setcounter{enumi}{3}
\tightlist
\item
  J. Peter May, \emph{Simplicial Objects in Algebraic Topology}, Van
  Nostrand, Princeton, 1968.
\end{enumerate}

Of course, not everyone prefers the austere joys of algebra to the
earthy pleasures of geometry. Algebraic topologists thrill to
categories, functors and natural transformations, while geometric
topologists like drawing pictures of hideously deformed multi-holed
doughnuts in 4 dimensional space. It's all a matter of taste.
Personally, I like both!

\begin{center}\rule{0.5\linewidth}{0.5pt}\end{center}

\hypertarget{homotopy_D}{\textbf{D.}} 
\emph{Simplicial objects.} We can generalize the heck out of
the notion of ``simplicial set'' by replacing the category
\(\mathsf{Set}\) with any other category \(\mathcal{C}\). A ``simplical
object in \(\mathcal{C}\)'' is defined to be a contravariant functor
\(F\colon\Delta\to\mathcal{C}\). There's a category whose objects are
such functors and whose morphisms are natural transformations between
them.

So, for example, a ``simplicial abelian group'' is a simplicial object
in the category of abelian groups. Just as we may associate to any set
\(X\) the free abelian group on \(X\), we may associate to any
simplicial set \(X\) the free simplicial abelian group on \(X\). In
fact, it's more than analogy: the latter construction is a spinoff of
the former! There is a functor \[L\colon\mathsf{Set}\to\mathsf{Ab}\]
assigning to any set the free abelian group on that set (see
\protect\hyperlink{week77}{``Week 77''}). Given a simplicial set
\[X\colon\Delta\to\mathsf{Set}\] we may compose with \(L\) to obtain a
simplicial abelian group \[XL\colon\Delta\to\mathsf{Ab}\] (where I'm
writing composition in the funny order that I like to use). This is the
free simplicial abelian group on the simplicial set \(X\)!

Later I'll talk about how to compute the homology groups of a simplicial
abelian group. Combined with the above trick, this will give a very
elegant way to define the homology groups of a simplicial set. Homology
groups are a very popular sort of invariant in algebraic topology; we
will get them with an absolute minimum of sweat.

Just as a good firework show ends with lots of explosions going off
simultaneously, leaving the audience stunned, deafened, and content, I
should end with a blast of abstraction, just for the hell of it. Those
of you who remember my discussion of ``theories'' in
\protect\hyperlink{week53}{``Week 53''} can easily check that there is a
category called the ``theory of abelian groups''. This allows us to
define an ``abelian group object'' in any category with finite limits.
In particular, since the category of simplicial sets has finite limits
(any presheaf category has all limits), we can define an abelian group
object in the category of simplicial sets. And now for a very pretty
result: abelian group objects in the category of simplicial sets are the
same as simplicial abelian groups! In other words, an abstract ``abelian
group'' living in the world of simplicial sets is the same as an
abstract ``simplicial set'' living in the world of abelian groups. I'm
very fond of this kind of ``commutativity of abstraction''.

Finally, I should emphasize that all of this stuff was first explained
to me by James Dolan. I just want to make these explanations available
to everyone.

\hypertarget{week116}{%
\section{February 7, 1998}\label{week116}}

While general relativity and the Standard Model of particle physics are
very different in many ways, they have one important thing in common:
both are gauge theories. I will not attempt to explain what a gauge
theory is here. I just want to recommend the following nice book on the
early history of this subject:

\begin{enumerate}
\def\labelenumi{\arabic{enumi})}
\tightlist
\item
  Lochlainn O'Raifeartaigh, \emph{The Dawning of Gauge Theory},
  Princeton U. Press, Princeton, 1997.
\end{enumerate}

This contains the most important early papers on the subject, translated
into English, together with detailed and extremely intelligent
commentary. It starts with Hermann Weyl's 1918 paper ``Gravitation and
Electricity'', in which he proposed a unification of gravity and
electromagnetism. This theory was proven wrong by Einstein in a
one-paragraph remark which appears at the end of Weyl's paper ---
Einstein noticed it would predict atoms of variable size! --- but it
highlighted the common features of general relativity and Maxwell's
equations, which were later generalized to obtain the modern concept of
gauge theory.

It also contains Theodor Kaluza's 1921 paper ``On the Unification
Problem of Physics'' and Oskar Klein's 1926 paper ``Quantum Theory and
Five-Dimensional Relativity''. These began the trend, currently very
popular in string theory, of trying to unify forces by postulating
additional dimensions of spacetime. It's interesting how gauge theory
has historical roots in this seemingly more exotic notion. The original
Kaluza--Klein theory assumed a \(5\)-dimensional spacetime, with the
extra dimension curled into a small circle. Starting with
\(5\)-dimensional general relativity, and using the \(\mathrm{U}(1)\)
symmetry of the circle, they recovered \(4\)-dimensional general
relativity coupled to a \(\mathrm{U}(1)\) gauge theory --- namely,
Maxwell's equations. Unfortunately, their theory also predicted an
unobserved spin-0 particle, which was especially problematic back in the
days before mesons were discovered.

I wasn't familiar with another item in this book, Wolfgang Pauli's
letter to Abraham Pais entitled ``Meson-Nucleon Interactions and
Differential Geometry''. This theory, ``written down July 22--25 1953 in
order to see how it looks'', postulated 2 extra dimensions in the shape
of a small sphere. The letter begins, ``Split a \(6\)-dimensional space
into a \((4+2)\)-dimensional one.'' At the time, meson-nucleon
interactions were believe to have an \(\mathrm{SU}(2)\) symmetry
corresponding to conservation of ``isospin''. Pauli obtained a theory
with this symmetry group using the \(\mathrm{SU}(2)\) symmetry of the
sphere.

Apparently Pauli got a lot of his inspiration from Weyl's 1929 paper
``Electron and Gravitation'', also reprinted in this volume. This
masterpiece did all the following things: it introduced the concept of
2-component spinors (see \protect\hyperlink{week109}{``Week 109''}),
considered the possibility that the laws of physics violate parity and
time reversal symmetry, introduced the tetrad formulation of general
relativity, introduced the notion of a spinor connection, and explicitly
derived electromagnetism from the gauge principle! A famously critical
fellow, Pauli lambasted Weyl's ideas on parity and time reversal
violation --- which are now known to be correct. But even he conceded
the importance of deriving Maxwell's equations from the gauge principle,
saying ``Here I must admit your ability in Physics''. And he
incorporated many of the ideas into his 1953 letter.

An all-around good read for anyone seriously interested in the history
of physics! It's best if you already know some gauge theory.

Now let me continue the tour of homotopy theory I began last week. I was
talking about simplices. Simplices are amphibious creatures, easily
capable of living in two different worlds. On the one hand, we can think
of them as topological spaces, and on the other hand, as purely
algebraic gadgets: objects in the category of finite totally ordered
sets, which we call \(\Delta\). This gives simplices a special role as a
bridge between topology and algebra.

This week I'll begin describing how this works. Next time we'll get into
some of the cool spinoffs. I'll keep up the format of listing tools one
by one:

\begin{center}\rule{0.5\linewidth}{0.5pt}\end{center}

\hypertarget{homotopy_E}{\textbf{E.}} 
\emph{Geometric realization.} In
\protect\hyperlink{week115}{``Week 115''} I talked about simplicial
sets. A simplicial set is a presheaf on the category \(\Delta\).
Intuitively, it's a purely combinatorial way of describing a bunch of
abstract simplices glued together along their faces. We want a process
that turns such things into actual topological spaces, and also a
process that turns topological spaces back into simplicial sets.

Let's start with the first one. Given a simplicial set \(X\), we can
form a space \(|X|\) called the ``geometric realization'' of \(X\) by
gluing spaces shaped like simplices together in the pattern given by
\(X\). Given a morphism between simplicial sets there's an obvious
continuous map between their geometric realizations, so geometric
realization is actually a functor
\[|\cdot|\colon\mathsf{SimpSet}\to\mathsf{Top}\] from the category of
simplicial sets, \(\mathsf{SimpSet}\), to the category of topological
spaces, \(\mathsf{Top}\).

It's straightforward to fill in the details. But if we want to be slick,
we can define geometric realization using the magic of adjoint functors
--- see below.

\begin{center}\rule{0.5\linewidth}{0.5pt}\end{center}

\hypertarget{homotopy_F}{\textbf{F.}}
\emph{Singular simplicial set.} The basic idea here is that
given a topological space \(X\), its ``singular simplicial set''
\(\mathrm{Sing}(X)\) consists of all possible ways of mapping simplices
into \(X\). This gives a functor
\[\mathsf{Sing}\colon\mathsf{Top}\to\mathsf{SimpSet}.\] We make this
precise as follows.

By thinking of simplices as spaces in the obvious way, we can associate
a space to any object of \(\Delta\), and also a continuous map to any
morphism in \(\Delta\). Thus there's a functor
\[i\colon\Delta\to\mathsf{Top}.\] For any space \(X\) we define
\[\mathrm{Sing}(X)\colon\Delta\to\mathsf{Set}\] by
\[\mathrm{Sing}(X)(-) = \operatorname{Hom}(i(-),X)\] where the blank
slot indicates how \(\mathrm{Sing}(X)\) is waiting to eat a simplex and
spit out the set of all ways of mapping it --- thought of as a space!
--- into the space \(X\). The blank slot also indicates how
\(\mathrm{Sing}(X)\) is waiting to eat a \emph{morphism} between
simplices and spit out a \emph{function} between sets.

Having said what Sing does to \emph{spaces}, what does it do to
\emph{maps}? The same formula works: for any map \(f\colon X\to Y\)
between topological spaces, we define
\[\mathrm{Sing}(f)(-) = \operatorname{Hom}(i(-),f).\] It may take some
head scratching to understand this, but if you work it out, you'll see
it works out fine. If you feel like you are drowning under a tidal wave
of objects, morphisms, categories, and functors, don't worry! Medical
research has determined that people actually grown new neurons when
learning category theory.

In fact, even though it might not seem like it, I'm being incredibly
pedagogical and nurturing. If I were really trying to show off, I would
have compressed the last couple of paragraphs into the following one
line: \[\mathrm{Sing}(--)(-) = \operatorname{Hom}(i(-),--).\] where Sing
becomes a functor using the fact that for any category \(\mathcal{C}\)
there's a functor
\[\operatorname{Hom}\colon\mathcal{C}^{\mathrm{op}}\times\mathcal{C}\to\mathsf{Set}\]
where \(\mathcal{C}^{\mathrm{op}}\) denotes the opposite of
\(\mathcal{C}\), that is, \(\mathcal{C}\) with all its arrows turned
around. (See \protect\hyperlink{week78}{``Week 78''} for an explanation
of this.)

Or I could have said this: form the composite
\[\Delta^{\mathrm{op}}\times\mathsf{Top} \xrightarrow{i\times1} \mathsf{Top}^{\mathrm{op}}\times\mathsf{Top} \xrightarrow{\operatorname{Hom}} \mathsf{Set}\]
and dualize this to obtain
\[\mathrm{Sing}\colon\mathsf{Top}\to\mathsf{SimpSet}.\] These are all
different ways of saying the same thing. Forming the singular simplical
set of a space is not really an ``inverse'' to geometric realization,
since if we take a simplicial set \(X\), form its geometric realization,
and then form the singular simplicial set of that, we get something much
bigger than \(X\). However, if you think about it, there's an obvious
map from \(X\) into \(\mathrm{Sing}(|X|)\). Similarly, if we start with
a topological space \(X\), there's an obvious map from
\(|\mathrm{Sing}(X)|\) down to \(X\).

What this means is that \(\mathrm{Sing}\) is the right adjoint of
\(|\cdot|\), or in other words, \(|\cdot|\) is the left adjoint of
\(\mathrm{Sing}\). Thus if we want to be slick, we can just
\emph{define} geometric realization to be the left adjoint of
\(\mathrm{Sing}\). (See \protect\hyperlink{week77}{``Week
77''}-\protect\hyperlink{week79}{``Week 79''} for an exposition of
adjoint functors.)

\begin{center}\rule{0.5\linewidth}{0.5pt}\end{center}

\hypertarget{homotopy_G}{\textbf{G.}} 
\emph{Chain complexes.} Now gird yourself for some utterly
unmotivated definitions! If you've taken a basic course in algebraic
topology, you have probably learned about chain complexes already, and
if you haven't, you probably aren't reading this anymore --- so I'll
just plunge in.

A ``chain complex'' \(C_\bullet\) is a sequence of abelian groups and
``boundary'' homomorphisms like this:
\[C_0 \xleftarrow{d_1} C_1 \xleftarrow{d_2} C_2 \xleftarrow{d_3} C_3 \leftarrow \ldots\]
satisfying the magic equation \[d_i d_{i+1}x=0.\] This equation says
that the image of \(d_{i+1}\) is contained in the kernel of \(d_i\), so
we may define the ``homology groups'' to be the quotients
\[H_i(C_\bullet) = \mathrm{ker}(d_i)/\mathrm{im}(d_{i+1}).\] The study
of this stuff is called ``homological algebra''. You can read about it
in such magisterial tomes as:

\begin{enumerate}
\def\labelenumi{\arabic{enumi})}
\setcounter{enumi}{1}
\tightlist
\item
  Henri Cartan and Samuel Eilenberg, \emph{Homological Algebra},
  Princeton U.\ Press, 1956.
\end{enumerate}
\noindent
or
\begin{enumerate}
\def\labelenumi{\arabic{enumi})}
\setcounter{enumi}{2}
\tightlist
\item
  Saunders Mac Lane, \emph{Homology}, Springer, Berlin, 1995.
\end{enumerate}
\noindent
But if you want something a bit more user-friendly, try:

\begin{enumerate}
\def\labelenumi{\arabic{enumi})}
\setcounter{enumi}{3}
\tightlist
\item
  Joseph J. Rotman, \emph{An Introduction to Homological Algebra},
  Academic Press, New York, 1979.
\end{enumerate}

The main reason chain complexes are interesting is that they are similar
to topological spaces, but simpler. In ``singular homology theory'', we
use a certain functor to convert topological spaces into chain
complexes, thus reducing topology problems to simpler algebra problems.
This is usually one of the first things people study when they study
algebraic topology. In sections \hyperlink{homotopy_G}{G} and
\hyperlink{homotopy_H}{H}  below, I'll remind you how this goes.

Though singular homology is very useful, not everybody gets around to
learning the deep reason why! In fact, chain complexes are really just
another way of talking about a certain especially simple class of
topological spaces, called ``products of Eilenberg--Mac Lane spaces of
abelian groups''. In such spaces, topological phenomena in different
dimensions interact in a particularly trivial way. Singular homology
thus amounts to neglecting the subtler interactions between topology in
different dimensions. This is what makes it so easy to work with --- yet
ultimately so limited.

Before I keep rambling on I should describe the category of chain
complexes, which I'll call \(\mathsf{Chain}\). The objects are just
chain complexes, and given two of them, say \(C\) and \(C'\), a morphism
\(f\colon C\to C'\) is a sequence of group homomorphisms
\[f_i\colon C_i\to C'_i\] making the following big diagram commute: \[
  \begin{tikzcd}
    C_0 \dar[swap,"f_0"]
    & C_1 \lar[swap,"d_1"] \dar[swap,"f_1"]
    & C_2 \lar[swap,"d_2"] \dar[swap,"f_2"]
    & C_3 \lar[swap,"d_3"] \dar[swap,"f_3"]
    & \ldots \lar
  \\C'_0
    & C'_1 \lar[swap,"d'_1"]
    & C'_2 \lar[swap,"d'_2"]
    & C'_3 \lar[swap,"d'_3"]
    & \ldots \lar
  \end{tikzcd}
\]

The reason \(\mathsf{Chain}\) gets to be so much like the category
\(\mathsf{Top}\) of topological spaces is that we can define homotopies
between morphisms of chain complexes by copying the definition of
homotopies between continuous maps. First, there is a chain complex
called I that's analogous to the unit interval. It looks like this:
\[\mathbb{Z}\oplus\mathbb{Z} \xleftarrow{d_1} \mathbb{Z} \xleftarrow{d_2} 0 \xleftarrow{d_3} 0 \xleftarrow{d_4} \ldots\]
The only nonzero boundary homomorphism is \(d_1\), which is given by
\[d_1(x) = (x,-x)\] (Why? We take \(I_1 = \mathbb{Z}\) and
\(I_0 = \mathbb{Z}\oplus\mathbb{Z}\) because the interval is built out
of one \(1\)-dimensional thing, namely itself, and two 0-dimensional
things, namely its endpoints. We define \(d_1\) the way we do since the
boundary of an oriented interval consists of two points: its initial
endpoint, which is positively oriented, and its final endpoint, which is
negatively oriented. This remark is bound to be obscure to anyone who
hasn't already mastered the mystical analogies between algebra and
topology that underlie homology theory!)

There is a way to define a ``tensor product'' \(C\otimes C'\) of chain
complexes \(C\) and \(C'\), which is analogous to the product of
topological spaces. And there are morphisms
\[i,j\colon C\to I\otimes C\] analogous to the two maps from a space
into its product with the unit interval: \[
  \begin{gathered}
    i,j\colon X\to[0,1]\times X
  \\i(x) = (0,x), \quad j(x)=(1,x).
  \end{gathered}
\] Using these analogies we can define a ``chain homotopy'' between
chain complex morphisms \(f,g\colon C\to C'\) in a way that's completely
analogous to a homotopy between maps. Namely, it's a morphism
\(F\colon I\otimes C\to C'\) for which the composite
\[C\xrightarrow{i}I\otimes C\xrightarrow{F}C'\] equals \(f\), and the
composite \[C\xrightarrow{j}I\otimes C\xrightarrow{F}C'\] equals \(g\).
Here we are using the basic principle of category theory: when you've
got a good idea, write it out using commutative diagrams and then
generalize the bejeezus out of it!

The nice thing about all this is that a morphism of chain complexes
\(f\colon C\to C'\) gives rise to homomorphisms of homology groups,
\[H_n(f)\colon H_n(C)\to H_n(C').\] In fact, we've got a functor
\[H_n\colon\mathsf{Chain}\to\mathsf{Ab}.\]

And even better, if \(f\colon C\to C'\) and \(g\colon C\to C'\) are
chain homotopic, then \(H_n(f)\) and \(H_n(g)\) are equal. So we say:
``homology is homotopy-invariant''.

\begin{center}\rule{0.5\linewidth}{0.5pt}\end{center}

\hypertarget{homotopy_H}{\textbf{H.}} 
\emph{The chain complex of a simplicial abelian group}. Now
let me explain a cool way of getting chain complexes, which goes a long
way towards explaining why they're important. Recall from section
\hyperlink{homotopy_D}{D} in \protect\hyperlink{week115}{``Week 115''} 
that a simplicial abelian
group is a contravariant functor \(C\colon\Delta\to\mathsf{Ab}\). In
particular, it gives us an abelian group \(C_n\) for each object \(n\)
of \(\Delta\), and also ``face'' homomorphisms
\[\partial_0,\ldots\partial_{n-1}\colon C_n\to C_{n-1}\] coming from all
the ways the simplex with \((n-1)\) vertices can be a face of the
simplex with \(n\) vertices. We can thus can make \(C\) into a chain
complex by defining \(d_n\colon C_n\to C_{n-1}\) as follows:
\[d_n = \sum_i (-1)^i \partial_i.\] The thing to check is that
\[d_n d_{n+1} x = 0.\] The alternating signs make everything cancel out!
In the immortal words of the physicist John Wheeler, ``the boundary of a
boundary is zero''.

Unsurprisingly, this gives a functor from simplicial abelian groups to
chain complexes. Let's call it
\[\mathrm{Ch}\colon\mathsf{SimpAb}\to\mathsf{Chain}\] More surprisingly,
this is an equivalence of categories! I leave you to show this --- if
you give up, look at May's book cited in section \hyperlink{homotopy_C}{C} of
\protect\hyperlink{week115}{``Week 115''}. What this means is that
simplicial abelian groups are just another way of thinking about chain
complexes\ldots{} or vice versa. Thus, if I were being
ultra-sophisticated, I could have skipped the chain complexes and talked
only about simplicial abelian groups! This would have saved time, but
more people know about chain complexes, so I wanted to mention them.

\begin{center}\rule{0.5\linewidth}{0.5pt}\end{center}

\hypertarget{homotopy_I}{\textbf{I.}} 
\emph{Singular homology.} Okay, now that we have lots of
nice shiny machines, let's hook them up and see what happens! Take the
``singular simplicial set'' functor:
\[\mathrm{Sing}\colon\mathsf{Top}\to\mathsf{SimpSet},\] the ``free
simplicial abelian group on a simplicial set'' functor:
\[L\colon\mathsf{SimpSet}\to\mathsf{SimpAb},\] and the ``chain complex
of a simplicial abelian group'' functor:
\[\mathrm{Ch}\colon\mathsf{SimpAb}\to\mathsf{Chain},\] and compose them!
We get the ``singular chain complex'' functor
\[C\colon\mathsf{Top}\to\mathsf{Chain}\] that takes a topological space
and distills a chain complex out of it. We can then take the homology
groups of our chain complex and get the ``singular homology'' of our
space. Better yet, the functor \(C\colon\mathsf{Top}\to\mathsf{Chain}\)
takes homotopies between maps and sends them to homotopies between
morphisms of chain complexes! It follows that homotopic maps between
spaces give the same homomorphisms between the singular homology groups
of these spaces. Thus homotopy-equivalent spaces will have isomorphic
homology groups\ldots{} so we have gotten our hands on a nice tool for
studying spaces up to homotopy equivalence.

Now that we've got our hands on singular homology, we could easily spend
a long time using it to solve all sorts of interesting problems. I won't
go into that here; you can read about it in all sorts of textbooks,
like:

\begin{enumerate}
\def\labelenumi{\arabic{enumi})}
\setcounter{enumi}{4}
\tightlist
\item
  Marvin J. Greenberg, John R. Harper, \emph{Algebraic Topology: A First
  Course}, Benjamin/Cummings, Reading, Massachusetts, 1981.
\end{enumerate}

or

\begin{enumerate}
\def\labelenumi{\arabic{enumi})}
\setcounter{enumi}{5}
\tightlist
\item
  William S. Massey, \emph{Singular Homology Theory}, Springer,
  Berlin, 1980.
\end{enumerate}

which uses cubes rather than simplices.

What I'm trying to emphasize here is that singular homology is a
composite of functors that are interesting in their own right. I'll
explore their uses a bit more deeply next time.

\begin{center}\rule{0.5\linewidth}{0.5pt}\end{center}

\begin{quote}
\emph{At a very early age, I made an assumption that a successful
physicist only needs to know elementary mathematics. At a later time, to
my great regret, I realized that this assumption of mine was completely
wrong.}

--- Albert Einstein
\end{quote}

\hypertarget{week117}{%
\section{February 14, 1998}\label{week117}}

A true physicist loves matter in all its states. The phases we all
learned about in school - solid, liquid, and gas - are just the
beginning of the story! There lots of others: liquid crystal, plasma,
superfluid, and neutronium, for example. Today I want to say a little
about two more phases that people are trying to create: quark-gluon
plasma and strange quark matter. The first almost certainly exists; the
second is a matter of much discussion.

\begin{enumerate}
\def\labelenumi{\arabic{enumi})}
\tightlist
\item
  The E864 Collaboration, ``Search for charged strange quark matter
  produced in 11.5 A GeV/c Au + Pb collisions'', \emph{Phys. Rev.~Lett.}
  \textbf{79} (1997), 3612--3616.  Also available as
  \href{https://arxiv.org/abs/nucl-ex/9706004}{\texttt{nucl-ex/9706004}}.
\end{enumerate}

Last week I went to a talk on the search for strange quark matter by one
of these collaborators, Kenneth Barish. This talk was based on Barish's
work at the E864 experiment at the ``AGS'', the alternating gradient
synchrotron at Brookhaven National Laboratory in Long Island, New York.

What's ``strange quark matter''? Well, first remember from
\protect\hyperlink{week93}{``Week 93''} that in the Standard Model there
are bosonic particles that carry forces:

\begin{longtable}[]{@{}lll@{}}
\toprule
Electromagnetic force & Weak force & Strong force\tabularnewline
\midrule
\endhead
photon & W\textsubscript{+}, W\textsubscript{-}, Z & 8
gluons\tabularnewline
\bottomrule
\end{longtable}

and fermionic particles that constitute matter:

\begin{longtable}[]{@{}llll@{}}
\toprule
\textbf{Leptons} & & \textbf{Quarks} &\tabularnewline
\midrule
\endhead
electron & electron neutrino & down quark & up quark\tabularnewline
muon & muon neutrino & strange quark & charm quark\tabularnewline
tauon & tauon neutrino & bottom quark & top quark\tabularnewline
\bottomrule
\end{longtable}

(There is also the mysterious Higgs boson, which has not yet been seen.)

The quarks and leptons come in 3 generations each. The only quarks in
ordinary matter are the lightest two, those from the first generation:
the up and down. These are the constituents of protons and neutrons,
which are the only stable particles made of quarks. A proton consists of
two ups and a down held together by the strong force, while a neutron
consists of two downs and a up. The up has electric charge \(+2/3\),
while the down has electric charge \(-1/3\). They also interact via the
strong and weak forces.

The other quarks are more massive and decay via the weak interaction
into up and down quarks. Apart from that, however, they are quite
similar. There are lots of short-lived particles made of various
combinations of quarks. All the combinations we've seen so far are of
two basic sorts. There are ``baryons'', which consist of 3 quarks, and
``mesons'', which consist of a quark and an antiquark. Both of these
should be visualized roughly as a sort of bag with the quarks and a
bunch of gluons confined inside.

Why are they confined? Well, I sketched an explanation in
\protect\hyperlink{week94}{``Week 94''}, so you should read that for
more details. For now let's just say the strong force likes to ``stick
together'', so that energy is minimized if it stays concentrated in
small regions, rather than spreading all over the place, like the
electromagnetic field does. Indeed, the strong force may even do
something like this in the absence of quarks, forming short-lived
``glueballs'' consisting solely of gluons and virtual quark-antiquark
pairs. (For more on glueballs, see \protect\hyperlink{week68}{``Week
68''}.)

For reasons I don't really understand, the protons and neutrons in the
nucleus do not coalesce into one big bag of quarks. Even in a neutron
star, the quarks stay confined in their individual little bags. But
calculations suggest that at sufficiently high temperatures or
pressures, ``deconfinement'' should occur. Basically, under these
conditions the baryons and mesons either smash into each other so hard,
or get so severely squashed, that they burst open. The result should be
a soup of free quarks and gluons: a ``quark-gluon plasma''.

To get deconfinement to happen is not easy --- at low pressures, it's
expected to occur at a temperature of 2 trillion Kelvin! According to
the conventional wisdom in cosmology, the last time deconfinement was
prevalent was about 1 microsecond after the big bang! In the E864
experiment, they are accelerating gold nuclei to energies of 11.5 GeV
per nucleon and colliding them with a fixed target made of lead, which
is apparently \emph{not} enough energy to fully achieve deconfinement
--- they believe they are reaching temperatures of about 1 trillion
Kelvin. At CERN they are accelerating lead nuclei to 160 GeV per nuclei
and colliding them with a lead target. They may be getting signs of
deconfinement, but as Jim Carr explained in a recent post to
sci.physics, they're being very cautious about coming out and saying so.
By mid-1999, the folks at Brookhaven hope to get higher energies with
the Relativistic Heavy Ion Collider, which will collide two beams of
gold nuclei head-on at 100 GeV per nucleon\ldots{} see
\protect\hyperlink{week76}{``Week 76''} for more on this.

One of the hoped-for signs of deconfinement is ``strangeness
enhancement''. The lightest quark besides the up and down is the strange
quark, and in the high energies present in a quark gluon plasma, strange
quarks should be formed. Moreover, since Pauli exclusion principle
prevents two identical fermions from being in the same state, it can be
energetically favorable to have strange quarks around, since they can
occupy lower-energy states which are already packed with ups and downs.
They seem to be seeing strangeness enhancement at CERN:

\begin{enumerate}
\def\labelenumi{\arabic{enumi})}
\setcounter{enumi}{1}
\tightlist
\item
  J\"urgen Eschke, NA35 Collaboration, ``Strangeness enhancement in
  sulphur-nucleus collisions at 200 GeV/N''.  Also available as
  \href{https://arxiv.org/abs/hep-ph/9609242}{\texttt{hep-ph/9609242}}.
\end{enumerate}

As far as I can tell, people are just about as sure that deconfinement
occurs at high temperatures as they would be that tungsten boils at high
temperatures, even if they've never actually seen it happen. A more
speculative possibility is that as quark-gluon plasma cools down it
forms ``strange quark matter'' in the form of ``strangelets'': big bags
of up, down, and strange quarks. This is what they're looking for at
E864. Their experiment would only detect strangelets that live long
enough to get to the detector. When their experiment is running they get
\(10^6\) collisions per second. So far they've set an upper bound of
\(10^{-7}\) charged strangelets per collision, neutral strangelets being
harder to detect and rule out. For more on strangelets, try this:

\begin{enumerate}
\def\labelenumi{\arabic{enumi})}
\setcounter{enumi}{2}
\tightlist
\item
  E. P. Gilson and R. L. Jaffe, ``Very small strangelets'', \emph{Phys.
  Rev.~Lett.} \textbf{71} (1993) 332--335.  Also available as
  \href{https://arxiv.org/abs/hep-ph/9302270}{\texttt{hep-ph/9302270}}.
\end{enumerate}

Strange quark matter is also of interest in astrophysics. In 1984 Witten
wrote a paper proposing that in the limit of large quark number, strange
quark matter could be more stable than ordinary nuclear matter!

\begin{enumerate}
\def\labelenumi{\arabic{enumi})}
\setcounter{enumi}{3}
\tightlist
\item
  Edward Witten, ``Cosmic separation of phases'', \emph{Phys. Rev.}
  \textbf{D30} (1984) 272--285.
\end{enumerate}

More recently, a calculation of Farhi and Jaffe estimates that in the
limit of large quark number, the energy of strange quark matter is 301
MeV per quark, as compared with 310 Mev/quark for iron-56, which is the
most stable nucleus. This raises the possibility that under suitable
conditions, a neutron star could collapse to become a ``quark star'' or
``strange star''. Let me quote the abstract of the following paper:

\begin{enumerate}
\def\labelenumi{\arabic{enumi})}
\setcounter{enumi}{4}
\tightlist
\item
  Dany Page, ``Strange stars: Which is the ground state of QCD at finite
  baryon number?'', in \emph{High Energy Phenomenology} eds.~M. A. Perez
  and R. Huerta, World Scientific, Singapore, 1992, pp.~347--356.  Also
  available as
  \href{https://arxiv.org/abs/astro-ph/9602043}{\texttt{astro-ph/9602043}}.
\end{enumerate}

\begin{quote}
Witten's conjecture about strange quark matter (`Strange Matter') being
the ground state of QCD at finite baryon number is presented and stars
made of strange matter (`Strange Stars') are compared to neutron stars.
The only observable way in which a strange star differs from a neutron
star is in its early thermal history and a detailed study of strange
star cooling is reported and compared to neutron star cooling. One
concludes that future detection of thermal radiation from the compact
object produced in the core collapse of SN 1987A could present the first
evidence for strange matter.
\end{quote}

Here are a couple of books on the subject, which unfortunately I've not
been able to get ahold of:

\begin{enumerate}
\def\labelenumi{\arabic{enumi})}
\setcounter{enumi}{5}
\item
  \emph{Strange Quark Matter in Physics and Astrophysics: Proceedings of
  the International Workshop on Strange Quark Matter in Physics and
  Astrophysics}, ed.~Jes Madsen, North-Holland, Amsterdam, 1991.
\item
  \emph{International Symposium on Strangeness and Quark Matter},
  eds.~Georges Vassiliadis et al, World Scientific, Singapore, 1995.
\end{enumerate}

If anyone out there knows more about the latest theories of strange
quark matter, and can explain them in simple terms, I'd love to hear
about it.

Okay, enough of that.   Now, on with my tour of homotopy theory!

So far I've mainly been talking about simplicial sets. I described a
functor called ``geometric realization'' that turns a simplicial set
into a topological space, and another functor that turns a space into a
simplicial set called its ``singular simplicial set''. I also showed how
to turn a simplicial set into a simplicial abelian group, and how to
turn one of \emph{those} into a chain complex\ldots{} or vice versa.

As you can see, the key is to have lots of functors at your disposal, so
you can take a problem in any given context --- or more precisely, any
given category! --- and move it to other contexts where it may be easier
to solve. Eventually I want to talk about what all these categories
we're using have in common: they are all ``model categories''. Once we
understand that, we'll be able to see more deeply what's going on in all
the games we've been playing.

But first I want to describe a few more important tricks for turning
this into that. Recall from \protect\hyperlink{week115}{``Week 115''}
that there's a category \(\Delta\) whose objects \(0,1,2,\ldots\) are
the simplices, with \(n\) corresponding to the simplex with \(n\)
vertices --- the simplex with \(0\) vertices being the ``empty
simplex''. We can also define \(\Delta\) in a purely algebraic way as
the category of finite totally ordered sets, with \(n\) corresponding to
the totally ordered set \(\{0,1,\ldots,n-1\}\). The morphisms in
\(\Delta\) are then the order-preserving maps. Using this algebraic
definition we can do some cool stuff:

\begin{center}\rule{0.5\linewidth}{0.5pt}\end{center}

\hypertarget{homotopy_J}{\textbf{J.}} 
\emph{The nerve of a category.} This is a trick to turn a
category into a simplicial set. Given a category \(\mathcal{C}\), we
cook up the simplicial set \(\mathrm{Nerve}(\mathcal{C})\) as follows.
The 0-dimensional simplices of \(\mathrm{Nerve}(\mathcal{C})\) are just
the objects of \(\mathcal{C}\), which look like this: \[x\] The
\(1\)-simplices of \(\mathrm{Nerve}(\mathcal{C})\) are just the
morphisms, which look like this: \[x\xrightarrow{f}y\] The
\(2\)-simplices of \(\mathrm{Nerve}(\mathcal{C})\) are just the
commutative diagrams that look like this: \[
  \begin{tikzpicture}
    \node (x) at (0,0) {$x$};
    \node (y) at (1,1.7) {$y$};
    \node (z) at (2,0) {$z$};
    \draw[thick] (x) to node[fill=white]{$f$} (y);
    \draw[thick] (x) to node[fill=white]{$h$} (z);
    \draw[thick] (y) to node[fill=white]{$g$} (z);
  \end{tikzpicture}
\] where \(f\colon x\to y\), \(g\colon y\to z\), and \(h\colon x\to z\).
And so on. In general, the \(n\)-simplices of
\(\mathrm{Nerve}(\mathcal{C})\) are just the commutative diagrams in
\(\mathcal{C}\) that look like \(n\)-simplices!

When I first heard of this idea I cracked up. It seemed like an insane
sort of joke. Turning a category into a kind of geometrical object built
of simplices? What nerve! What use could this possibly be?

Well, for an application of this idea to computer science, see
\protect\hyperlink{week70}{``Week 70''}. We'll soon see lots of
applications within topology. But first, let me give a slick abstract
description of this ``nerve'' process that turns categories into
simplicial sets. It's really a functor
\[\mathrm{Nerve}\colon\mathsf{Cat}\to\mathsf{SimpSet}\] going from the
category of categories to the category of simplicial sets.

First, a remark on \(\mathsf{Cat}\). This has categories as objects and
functors as morphisms. Since the ``category of all categories'' is a bit
creepy, we really want the objects of \(\mathsf{Cat}\) to be all the
``small'' categories, i.e., those having a mere \emph{set} of objects.
This prevents Russell's paradox from raising its ugly head and
disturbing our fun and games.

Next, note that any partially ordered set can be thought of as a
category whose objects are just the elements of our set, and where we
say there's a single morphism from \(x\) to \(y\) if \(x\leqslant y\).
Composition of morphisms works out automatically, thanks to the
transitivity of ``less than or equal to''. We thus obtain a functor
\[i\colon\Delta\to\mathsf{Cat}\] taking each finite totally ordered set
to its corresponding category, and each order-preserving map to its
corresponding functor.

Now we can copy the trick we played in section \hyperlink{homotopy_F}{F} of
\protect\hyperlink{week116}{``Week 116''}. For any category
\(\mathcal{C}\) we define the simplicial set
\(\mathrm{Nerve}(\mathcal{C})\) by
\[\mathrm{Nerve}(\mathcal{C})(-) = \operatorname{Hom}(i(-),\mathcal{C})\]
Think about it! If you put the simplex \(n\) in the blank slot, we get
\(\operatorname{Hom}(i(n),\mathcal{C})\), which is the set of all
functors from that simplex, \emph{regarded as a category}, to the
category \(\mathcal{C}\). This is just the set of all diagrams in
\(\mathcal{C}\) shaped like the simplex \(n\), as desired!

We can say all this even more slickly as follows: take
\[\Delta^{\mathrm{op}}\times\mathsf{Cat} \xrightarrow{i\times1} \mathsf{Cat}^{\mathrm{op}}\times\mathsf{Cat} \xrightarrow{\operatorname{Hom}} \mathsf{Set}\]
and dualize it to obtain
\[\mathrm{Nerve}\colon\mathsf{Cat}\to\mathsf{SimpSet}.\]

I should also point out that topologists usually do this stuff with the
topologist's version of \(\Delta\), which does not include the ``empty
simplex''.

\begin{center}\rule{0.5\linewidth}{0.5pt}\end{center}

\hypertarget{homotopy_K}{\textbf{K.}} 
\emph{The classifying space of a category.} If compose our new
functor \[\mathrm{Nerve}\colon\mathsf{Cat}\to\mathsf{SimpSet}\] with the
``geometric realization'' functor
\[|\cdot|\colon\mathsf{SimpSet}\to\mathsf{Top}\] defined in section 
\hyperlink{homotopy_E}{E},
we get a way to turn a category into a space, called its ``classifying
space''. This trick was first used by Graeme Segal, the homotopy
theorist who later became the guru of conformal field theory. He
invented this trick in the following paper:

\begin{enumerate}
\def\labelenumi{\arabic{enumi})}
\setcounter{enumi}{7}
\tightlist
\item
  Graeme B. Segal, ``Classifying spaces and spectral sequences'',
  \emph{Publ. Math. Inst. des Haut. Etudes Scient.} \textbf{34} (1968),
  105--112.
\end{enumerate}

As it turns out, every reasonable space is the classifying space of some
category! More precisely, every space that's the geometric realization
of some simplicial set is homeomorphic to the classifying space of some
category. To see this, suppose the space \(X\) is the geometric
realization of the simplicial set \(S\). Take the set of all simplices
in \(S\) and partially order them by saying \(x\leqslant y\) if \(x\) is
a face of \(y\). Here by ``face'' I don't mean just mean a face of one
dimension less than that of \(y\); I'm allowing faces of any dimension
less than or equal to that of \(y\). We obtain a partially ordered set.
Now think of this as a category, \(\mathcal{C}\). Then
\(\mathrm{Nerve}(\mathcal{C})\) is the ``barycentric subdivision'' of
\(S\). In other words, it's a new simplicial set formed by chopping up
the simplices of \(S\) into smaller pieces by putting a new vertex in
the center of each one. It follows that the geometric realization of
\(\mathrm{Nerve}(\mathcal{C})\) is homeomorphic to that of \(S\).

There are lots of interesting special sorts of categories, like
groupoids, or monoids, or groups (see \protect\hyperlink{week74}{``Week
74''}). These give special cases of the ``classifying space''
construction, some of which were actually discovered before the general
case. I'll talk about some of these more next week, since they are very
important in topology.

Also sometimes people take categories that they happen to be interested
in, which may have no obvious relation to topology, and study them by
studying their classifying spaces. This gives surprising ways to apply
topology to all sorts of subjects. A good example is ``algebraic
K-theory'', where we start with some sort of category of modules over a
ring.

\begin{center}\rule{0.5\linewidth}{0.5pt}\end{center}

\hypertarget{homotopy_L}{\textbf{L.}}
\emph{\(\Delta\) as the free monoidal category on a monoid
object.} Recall that a ``monoid'' is a set with a product and a unit
element satisfying associativity and the left and right unit laws.
Categorifying this notion, we obtain the concept of a ``monoidal
category'': a category \(\mathcal{C}\) with a product and a unit object
satisfying the same laws. A nice example of a monoidal category is the
category \(\mathsf{Set}\) with its usual cartesian product, or the
category \(\mathsf{Vect}\) with its usual tensor product. We usually
call the product in a monoidal category the ``tensor product''.

Now, the ``microcosm principle'' says that algebraic gadgets often like
to live inside categorified versions of themselves. It's a bit like the
``homunculus theory'', where I have a little copy of myself sitting in
my head who looks out through my eyes and thinks all my thoughts for me.
But unlike that theory, it's true!

For example, we can define a ``monoid object'' in any monoidal category.
Given a monoidal category \(A\) with tensor product \(\otimes\) and unit
object \(1\), we define a monoid object \(a\) in \(A\) to be an object
equipped with a ``product'' \[m\colon a\otimes a\to a\] and a ``unit''
\[i\colon 1\to a\] which satisfy associativity and the left and right
unit laws (written out as commutative diagrams). A monoid object in
\(\mathsf{Set}\) is just a monoid, but a monoid object in
\(\mathsf{Vect}\) is an algebra, and I gave some very different examples
of monoid objects in \protect\hyperlink{week89}{``Week 89''}.

Now let's consider the ``free monoidal category on a monoid object''. In
other words, consider a monoidal category \(A\) with a monoid object
\(a\) in it, and assume that \(A\) has no objects and no morphisms, and
satisfies no equations, other than those required by the definitions of
``monoidal category'' and ``monoid object''.

Thus the only objects of \(A\) are the unit object together with \(a\)
and its tensor powers. Similarly, all the morphism of \(A\) are built up
by composing and tensoring the morphisms \(m\) and \(i\). So \(A\) looks
like this: \[
  \begin{tikzcd}[column sep=large]
    1
      \rar["i" description]
    & a
      \rar[shift left=5,"1\otimes i" description]
      \rar["i\otimes1" description]
    & a\otimes a
      \lar[shift left=5,"m" description]
      \rar[shift left=10,"1\otimes1\otimes i" description]
      \rar[shift left=5,"1\otimes i\otimes1" description]
      \rar["i\otimes1\otimes1" description]
    & a\otimes a\otimes a \quad \ldots
      \lar[shift left=5,"m\otimes1" description]
      \lar[shift left=10,"1\otimes m" description]
  \end{tikzcd}
\] Here I haven't drawn all the morphisms, just enough so that every
morphism in \(A\) is a composite of morphisms of this sort.

What is this category? It's just \(\Delta\)! The \(n\)th tensor power of
a corresponds to the simplex with \(n\) vertices. The morphisms going to
the right describe the ways the simplex with n vertices can be a face of
the simplex with \(n+1\) vertices. The morphisms going to the left
correspond to ``degeneracies'' --- ways of squashing a simplex with
\(n+1\) vertices down into one with \(n\) vertices.

So: in addition to its other descriptions, we can define \(\Delta\) as
the free monoidal category on a monoid object! Next time we'll see how
this is fundamental to homological algebra.

\hypertarget{week118}{%
\section{March 14, 1998}\label{week118}}

Like many people of a certain age, as a youth my interest in mathematics
and physics was fed by the Scientific American, especially Martin
Gardner's wonderful column. Since then the magazine seems to have gone
downhill. For me, the last straw was a silly article on the ``death of
proof'' in mathematics, written by someone wholly unfamiliar with the
subject. The author of that article later wrote a book proclaiming the
``end of science'', and went on to manage a successful chain of funeral
homes.

Recently, however, I was pleased to find a terse rebuttal of this
fin-de-siecle pessimism in an article appearing in --- none other than
Scientific American!

\begin{enumerate}
\def\labelenumi{\arabic{enumi})}
\tightlist
\item
  Michael J. Duff, ``The theory formerly known as strings'',
  \emph{Scientific American} \textbf{278} (February 1998), 64--69.
\end{enumerate}
\noindent
The article begins:

\begin{quote}
At a time when certain pundits are predicting the End of Science on the
grounds that all the important discoveries have already been made, it is
worth emphasizing that the two main pillars of 20th-century physics,
quantum mechanics and Einstein's general theory of relativity, are
mutually incompatible.
\end{quote}

To declare the end of science at this point, or even of particle physics
(the two are not the same!) would thus be ridiculously premature. It's
true that the quest for a unified theory of all the forces and particles
in nature is experiencing difficulties. On the one hand, particle
accelerators have become very expensive. On the other hand, it's truly
difficult to envision a consistent and elegant formalism subsuming both
general relativity and the Standard Model of particle physics - much
less one that makes new testable predictions. But hey, the course of
true love never did run smooth.

Duff's own vision certainly has its charms. He has long been advocating
the generalization of string theory to a theory of higher-dimensional
``membranes''. Nowadays people call these ``\(p\)-branes'' to keep track
of the dimension of the membrane: a \(0\)-brane is a point particle, a
\(1\)-brane is a string, a \(2\)-brane is a \(2\)-dimensional surface,
and so on.

For a long time, higher-dimensional membrane theories were unpopular,
even among string theorists, because the special tricks that eliminate
infinities in string theory don't seem to work for higher-dimensional
membranes. But lately membranes are all the rage: it seems they show up
in string theory whether or not you put them in from the start! In fact,
they seem to be the key to showing that the 5 different supersymmetric
string theories are really aspects of a single deeper theory ---
sometimes called ``M-theory''.

Now, I don't really understand this stuff at all, but I've been trying
to learn about it lately, so I'll say a bit anyway, in hopes that some
real experts will correct my mistakes. Much of what I'll say comes from
the following nice review article:

\begin{enumerate}
\def\labelenumi{\arabic{enumi})}
\setcounter{enumi}{1}
\tightlist
\item
  M. J. Duff, ``Supermembranes'', available as
  \href{https://arxiv.org/abs/hep-th/9611203}{\texttt{hep-th/9611203}}.
\end{enumerate}
\noindent
and also the bible of string theory:

\begin{enumerate}
\def\labelenumi{\arabic{enumi})}
\setcounter{enumi}{2}
\tightlist
\item
  Michael B. Green, John H. Schwarz, and Edward Witten,
  \emph{Superstring Theory}, two volumes, Cambridge U.\ Press, Cambridge,
  1987.
\end{enumerate}

Let's start with superstring theory. Here the ``super'' refers to the
fact that instead of just strings whose vibrational modes correspond to
bosonic particles, we have strings with extra degrees of freedom
corresponding to fermionic particles. We can actually think of the
superstring as a string wiggling around in a ``superspace'': a kind of
space with extra ``fermionic'' dimensions in addition to the usual
``bosonic'' ones. These extra dimensions are described by coordinates
that anticommute with each other, and commute with the usual bosonic
coordinates (which of course commute with each other). This amounts to
taking the boson/fermion distinction so seriously that we build it into
our description of spacetime from the very start! For more details on
the mathematics of superspace, try:

\begin{enumerate}
\def\labelenumi{\arabic{enumi})}
\setcounter{enumi}{3}
\tightlist
\item
  Bryce DeWitt, \emph{Supermanifolds}, Cambridge U. Press, Cambridge,
  1992.
\end{enumerate}

More deeply, ``super'' refers to ``supersymmetry'', a special kind of
symmetry transformation that mixes the bosonic and fermionic
coordinates. We speak of ``\(N = 1\) supersymmetry'' if there is one
fermionic coordinate for each bosonic coordinate, ``\(N = 2\)
supersymmetry'' if there are two, and so on.

Like all nice physical theories, we can in principle derive everything
about our theory of superstrings once we know the formula for the
\emph{action}. For bosonic strings, the action is very simple. As time
passes, a string traces out a \(2\)-dimensional surface in spacetime
called the ``string worldsheet''. The action is just the \emph{area} of
this worldsheet.

For superstring theory, we thus want a formula for the ``super-area'' of
a surface in superspace. And we need this to be invariant under
supersymmetry transformations. Suprisingly, this is only possible if
spacetime has dimension 3, 4, 6, or 10. More precisely, these are the
dimensions where \(N = 1\) supersymmetric string theory makes sense as a
\emph{classical} theory.

Note: these dimensions are just \(2\) more than the dimensions of the
four normed division algebras: the reals, complexes, quaternions and
octonions! This is no coincidence. Robert Helling recently posted a nice
article about this on \texttt{sci.physics.resarch}, which I have
appended to \protect\hyperlink{week104}{``Week 104''}. The basic idea is
that we can describe the vibrations of a string in \(n\)-dimensional
spacetime by a field on the string worldsheet with \(n-2\) components
corresponding to the \(n-2\) directions transverse to the worldsheet. To
get an action that's invariant under supersymmetry, we need some magical
cancellations to occur. It only works when we can think of this field as
taking values in one of the normed division algebras!

This is one of the curious things about superstring theory: the basic
idea is simple, but when you try to get it to work, you run into lots of
obstacles which only disappear in certain special circumstances ---
thanks to a mysterious conspiracy of beautiful mathematical facts. These
``conspiracies'' are probably just indications that we don't understand
the theory as deeply as we should. Right now I'm most interested in the
algebraic aspects of superstring theory --- and especially its
relationships to ``exceptional algebraic structures'' like the
octonions, the Lie group \(\mathrm{E}_8\), and so on. As I learn
superstring theory, I like keeping track of the various ways these
structures show up, like remembering the clues in a mystery novel.

Interestingly, the \emph{quantum} version of superstring theory is more
delicate than the classical version. When I last checked, it only makes
sense in dimension 10. Thus there's something inherently octonionic
about it! For more on this angle, see:

\begin{enumerate}
\def\labelenumi{\arabic{enumi})}
\setcounter{enumi}{4}
\item
  E. Corrigan and T. J. Hollowood, ``The exceptional Jordan algebra and
  the superstring'', \emph{Comm. Math. Phys.} \textbf{122} (1989),
  393.
\item
  E. Corrigan and T. J. Hollowood, ``A string construction of a
  commutative nonassociative algebra related to the exceptional Jordan
  algebra'', \emph{Phys. Lett.} \textbf{B203} (1988), 47.
\end{enumerate}
\noindent
and some more references I'll give later.

There are actually 5 variants of superstring theory in dimension 10, as
I explained in \protect\hyperlink{week72}{``Week 72''}:

\begin{enumerate}
\def\labelenumi{\arabic{enumi}.}
\tightlist
\item
  type I superstrings --- these are open strings, not closed loops.
\item
  type IIA superstrings --- closed strings where the left- and
  right-moving fermionic modes have opposite chiralities.
\item
  type IIB superstrings --- closed strings where the left- and
  right-moving fermionic modes have the same chirality.
\item
  \(\mathrm{E}_8\) heterotic superstrings --- closed strings where the
  left-moving modes are purely bosonic, with symmetry group
  \(\mathrm{E}_8\times\mathrm{E}_8\).
\item
  \(\mathrm{Spin}(32)/\mathbb{Z}_2\) heterotic superstrings --- closed
  strings where the left-moving modes are purely bosonic, with symmetry
  group \(\mathrm{Spin}(32)/\mathbb{Z}_2\)
\end{enumerate}

To get \(4\)-dimensional physics out of any of these, we need to think
of our \(10\)-dimensional spacetime as a bundle with a little
\(6\)-dimensional ``Calabi--Yau manifold'' sitting over each point of
good old 4-dimensional spacetime. But there's another step that's very
useful when trying to understand the implications of superstring theory
for ordinary particle physics. This is to look at the low-energy limit.
In this limit, only the lowest-energy vibrational modes of the string
contribute, each mode acting like a different kind of massless particle.
Thus in this limit superstring theory acts like an ordinary quantum
field theory.

What field theory do we get? This is a very important question. The
field theory looks simplest in \(10\)-dimensional Minkowski spacetime;
it gets more complicated when we curl up 6 of the dimensions and think
of it as a \(4\)-dimensional field theory, so let's just talk about the
simple situation.

No matter what superstring theory we start with, the low-energy limit
looks like some form of ``supergravity coupled to super-Yang--Mills
fields''. What's this? Well, supergravity is basically what we get when
we generalize Einstein's equations for general relativity to superspace.
Similarly, super-Yang--Mills theory is the supersymmetric version of the
Yang--Mills equations - which are important in particle physics because
they describe all the forces \emph{except} gravity. So superstring
theory has in it the seeds of general relativity and also the other
forces of nature --- or at least their supersymmetric analogues.

Like superstring theory, super-Yang--Mills theory only works in spacetime
dimensions 3, 4, 6, and 10. (See \protect\hyperlink{week93}{``Week 93''}
for more on this.) Different forms of supergravity make sense in
different dimensions, as explained in:

\begin{enumerate}
\def\labelenumi{\arabic{enumi})}
\setcounter{enumi}{6}
\tightlist
\item
  Y. Tanii, ``Introduction to supergravities in diverse dimensions'',
 "Introduction to supergravities in diverse dimensions'', \emph{Soryushiron 
  Kenkyu Electronics} \textbf{96} (1998), 315--351..  Also
   available as
  \href{https://arxiv.org/abs/hep-th/9802138}{\texttt{hep-th/9802138}}.
\end{enumerate}

In particular highest dimension in which supergravity makes sense is 11
dimensions (where one only has \(N = 1\) supergravity). Note that this
is one more than the favored dimension of superstring theory! This
puzzled people for a long time. Now it seems that M-theory is beginning
to resolve these puzzles. Another interesting discovery is that
11-dimensional supergravity is related to the exceptional Lie group
\(\mathrm{E}_8\):

\begin{enumerate}
\def\labelenumi{\arabic{enumi})}
\setcounter{enumi}{7}
\tightlist
\item
  Stephan Melosch and Hermann Nicolai, ``New canonical variables for
  \(d = 11\) supergravity'', \emph{Phys.\ Lett.\ B} \textbf{416} (1998),
  91--100.  Also available as
  \href{https://arxiv.org/abs/hep-th/9709227}{\texttt{hep-th/9709227}}.
\end{enumerate}

But I'm getting ahead of myself here! Right now I'm talking about the
low-energy limit of \(10\)-dimensional superstring theory. I said it
amounts to ``supergravity coupled to super-Yang--Mills fields'', and now
I'm attempting to flesh that out a bit. So: starting from \(N = 1\)
supergravity in 11 dimensions we can get a theory of supergravity in 10
dimensions simply by requiring that all the fields be constant in one
direction --- a trick called ``dimensional reduction''. This is called
``type IIA supergravity'', because it appears in the low-energy limit of
type IIA superstrings. There are also two other forms of supergravity in
10 dimensions: ``type IIB supergravity'', which appears in the
low-energy limit of type IIB superstrings, and a third form which
appears in the low-energy limit of the type I and heterotic
superstrings. These other two forms of supergravity are chiral --- that
is, they have a built-in ``handedness''.

Now let's turn to higher-dimensional supersymmetric membranes, or
``supermembranes''. Duff summarizes this subject in a chart he calls the
``brane scan''. This chart lists the known \emph{classical} theories of
supersymetric \(p\)-branes. Of course, a \(p\)-brane traces out a
(p+1)-dimensional surface as time passes, so from a spacetime point of
view it's p+1 which is more interesting. But anyway, here's Duff's chart
of which supersymmetric \(p\)-brane theories are possible in which
dimensions d of spacetime:

\begin{longtable}[]{@{}rccccccccccc@{}}
\toprule
\({}_d\diagdown^p\) & \(0\) & \(1\) & \(2\) & \(3\) & \(4\) & \(5\) &
\(6\) & \(7\) & \(8\) & \(9\) & \(10\)\tabularnewline
\midrule
\endhead
\(11\) & & & \checkmark & & & \checkmark & & & & ? &\tabularnewline
\(10\) & \checkmark & \checkmark & \checkmark & \checkmark & \checkmark
& \checkmark & \checkmark & \checkmark & \checkmark & \checkmark
&\tabularnewline
\(9\) & \checkmark & & & & \checkmark & & & & & &\tabularnewline
\(8\) & & & & \checkmark & & & & & & &\tabularnewline
\(7\) & & & \checkmark & & & \checkmark & & & & &\tabularnewline
\(6\) & \checkmark & \checkmark & \checkmark & \checkmark & \checkmark &
\checkmark & & & & &\tabularnewline
\(5\) & \checkmark & & \checkmark & & & & & & & &\tabularnewline
\(4\) & \checkmark & \checkmark & \checkmark & \checkmark & & & & & &
&\tabularnewline
\(3\) & \checkmark & \checkmark & \checkmark & & & & & & &
&\tabularnewline
\(2\) & \checkmark & & & & & & & & & &\tabularnewline
\(1\) & & & & & & & & & & &\tabularnewline
\bottomrule
\end{longtable}

We immediately notice some patterns. First, we see horizontal stripes in
dimensions 3, 4, 6, and 10: all the conceivable \(p\)-brane theories
exist in these dimensions. I don't know why this is true, but it must be
related to the fact that superstring and super-Yang--Mills theories make
sense in these dimensions. Second, there are four special \(p\)-brane
theories:

\begin{itemize}
\tightlist
\item
  the \(2\)-brane in dimension 4
\item
  the \(3\)-brane in dimension 6
\item
  the \(5\)-brane in dimension 10
\item
  the \(2\)-brane in dimension 11
\end{itemize}
\noindent
which are related to the real numbers, the complex numbers, the
quaternions and the octonions, respectively. Duff refers us to the
following papers for more information on this:

\begin{enumerate}
\def\labelenumi{\arabic{enumi})}
\setcounter{enumi}{8}
\item
  G. Sierra, ``An application of the theories of Jordan algebras and
  Freudenthal triple systems to particles and strings'', \emph{Class.
  Quant. Grav.} \textbf{4} (1987), 227.
\item
  J. M. Evans, ``Supersymmetric Yang--Mills theories and division
  algebras'', \emph{Nucl. Phys.} \textbf{B298} (1988), 92--108.
\end{enumerate}

From these four ``fundamental'' theories of \(p\)-branes in \(d\)
dimensions we can get theories of \((p-k)\)-branes in \(d-k\) dimensions
by dimensional reduction of both the spacetime and the \(p\)-brane. Thus
we see diagonal lines slanting down and to the left starting from these
``fundamental'' theories. Note that these diagonal lines include the
superstring theories in dimensions 3, 4, 6, and 10!

I'll wrap up by saying a bit about how M-theory, superstrings and
supergravity fit together. I've already said that:

\begin{enumerate}
\def\labelenumi{\arabic{enumi})}
\tightlist
\item
  type IIA supergravity in 10 dimensions is the dimensional reduction of
  11-dimensional supergravity; and
\item
  the type IIA superstring has typeIIA supergravity coupled to
  super-Yang--Mills fields as a low-energy limit.
\end{enumerate}

This suggests the presence of a theory in 11 dimensions that fills in
the question mark below: \[
  \begin{tikzpicture}
    \node (tl) at (0,0) {?};
    \node[text width=7em,align=center] (tr) at (5,0) {11-dimensional supergravity};
    \node[text width=7em,align=center] (bl) at (0,-3) {type IIA superstrings};
    \node[text width=7em,align=center] (br) at (5,-3) {type IIA supergravity in 10 dimensions};
    \draw[->] (tl) to node[label=above:{\scriptsize low-energy limit}]{} (tr);
    \draw[->] (bl) to node[label=above:{\scriptsize low-energy limit}]{} (br);
    \draw[->] (tl) to node[text width=4.2em,align=center,fill=white]{\scriptsize dimensional reduction} (bl);
    \draw[->] (tr) to node[text width=4.2em,align=center,fill=white]{\scriptsize dimensional reduction} (br);
  \end{tikzpicture}
\] This conjectured theory is called ``M-theory''. The actual details of
this theory are still rather mysterious, but not surprisingly, it's
related to the theory of supersymmetric \(2\)-branes in 11 dimensions
--- since upon dimensional reduction these give superstrings in 10
dimensions. More surprisingly, it's \emph{also} related to the theory of
\emph{5-branes} in 11 dimensions. The reason is that supergravity in 11
dimensions admits ``soliton'' solutions --- solutions that hold their
shape and don't disperse --- which are shaped like \(5\)-branes. These
are now believed to be yet another shadow of M-theory.

While the picture I'm sketching may seem baroque, it's really just a
small part of a much more elaborate picture that relates all 5
superstring theories to M-theory. But I think I'll stop here! Maybe
later when I know more I can fill in some more details. By the way, I
thank Dan Piponi for pointing out that Scientific American article.

For more on this business, check out the following review articles:

\begin{enumerate}
\def\labelenumi{\arabic{enumi})}
\setcounter{enumi}{10}
\item
  W. Lerche, ``Recent developments in string theory'', available as
  \href{https://arxiv.org/abs/hep-th/9710246}{\texttt{hep-th/9710246}}.
\item
  John Schwarz, ``The status of string theory'', available as
  \href{https://arxiv.org/abs/hep-th/9711029}{\texttt{hep-th/9711029}}.
\item
  M. J. Duff, ``M-theory (the theory formerly known as strings)'',
  \emph{Int.\ J.\ Mod.\ Phys.\ A} \textbf{11} (1996), 5623--5642.  Also
   available as 
  \href{https://arxiv.org/abs/hep-th/9608117}{\texttt{hep-th/9608117}}.
\end{enumerate}
\noindent
The first one is especially nice if you're interested in a nontechnical
survey; the other two are more detailed.

Okay. Now, back to my tour of homotopy theory! I had promised to talk
about classifying spaces of groups and monoids, but this post is getting
pretty long, so I'll only talk about something else I promised: the
foundations of homological algebra. So, remember:

As soon as we can squeeze a simplicial set out of something, we have all
sorts of methods for studying it. We can turn the simplicial set into a
space and then use all the methods of topology to study this space. Or
we can turn it into a chain complex and apply homology theory. So it's
very important to have tricks for turning all sorts of gadgets into
simplicial sets: groups, rings, algebras, Lie algebras, you name it! And
here's how\ldots.

\begin{center}\rule{0.5\linewidth}{0.5pt}\end{center}

\hypertarget{homotopy_M}{\textbf{M.}} 
\emph{Simplicial objects from adjunctions.} Remember from section
\hyperlink{homotopy_D}{D} of \protect\hyperlink{week115}{``Week 115''} 
that a ``simplicial object'' in some category is a contravariant functor from
\(\Delta\) to that category. In what follows, I'll take \(\Delta\) to be
the version of the category of simplices that contains the empty
simplex. Topologists don't usually do this, so what I'm calling a
``simplicial object'', they would call an ``augmented simplicial
object''. Oh well.

Concretely, a simplicial object in a category amounts to a bunch of
objects \(x_0, x_1, x_2,\ldots\) together with morphisms like this: \[
  \begin{tikzcd}[column sep=large]
    x_0
    & x_1
      \rar[shift left=5,"i_0" description]
      \lar["d_0" description]
    & x_2
      \rar[shift left=10,"i_1" description]
      \rar[shift left=5,"i_0" description]
      \lar["d_0" description]
      \lar[shift left=5,"d_1" description]
    & x_3\ldots
      \lar["d_0" description]
      \lar[shift left=5,"d_1" description]
      \lar[shift left=10,"d_2" description]
  \end{tikzcd}
\] The morphisms \(d_j\) are called ``face maps'' and the morphisms
\(i_j\) are called ``degeneracies''. They are required to satisfy some
equations which I won't bother writing down here, since you can figure
them out yourself if you read section \hyperlink{homotopy_B}{B} of
\protect\hyperlink{week114}{``Week 114''}.

Now, suppose we have an adjunction, that is, a pair of adjoint functors:
\[
  \begin{tikzcd}
    \mathcal{C} \rar[bend left=30,"L"]
    & \mathcal{D} \lar[bend left=30,"R"]
  \end{tikzcd}
\] This means we have natural transformations \[
  \begin{aligned}
    e&\colon LR\Rightarrow 1_{\mathcal{D}}
  \\i&\colon 1_{\mathcal{C}} \Rightarrow RL
  \end{aligned}
\] satisfying a couple of equations, which I again won't write down,
since I explained them in \protect\hyperlink{week79}{``Week 79''} and
\protect\hyperlink{week83}{``Week 83''}.

Then an object \(d\) in the category \(\mathcal{D}\) automatically gives
a simplicial object as follows: \[
  \begin{tikzcd}[column sep=huge]
    d
    & LR(d)
      \rar[shift left=5,"L\cdot i\cdot R" description]
      \lar["e" description]
    & LRLR(d)
      \rar[shift left=10,"LRL\cdot i\cdot R" description]
      \rar[shift left=5,"L\cdot i\cdot RLR" description]
      \lar["e\cdot LR" description]
      \lar[shift left=5,"LR\cdot e" description]
    & LRLRLR(d)
      \lar["e\cdot LRLR" description]
      \lar[shift left=5,"LR\cdot e\cdot LR" description]
      \lar[shift left=10,"LRLR\cdot e" description]
  \end{tikzcd}
\] where \(\cdot\) denotes horizontal composition of functors and
natural transformations.

For example, if \(\mathsf{AbGp}\) is the category of abelian groups, we
have an adjunction \[
  \begin{tikzcd}
    \mathsf{\phantom{|}Set} \rar[bend left=30,"L"]
    & \mathsf{AbGp} \lar[bend left=30,"R"]
  \end{tikzcd}
\] where \(L\) assigns to each set the free abelian group on that set, and \(R\)
assigns to each group its underlying set. Thus given an abelian group, the above
trick gives us a simplicial object in \(\mathsf{AbGp}\) --- or in other
words, a simplicial group. This has an underlying simplicial set, and
from this we can cook up a chain complex as in section 
\hyperlink{homotopy_H}{H} of
\protect\hyperlink{week116}{``Week 116''}. This lets us study groups
using homology theory! One can define the homology (and cohomology) of
lots other algebraic gadgets in exactly the same way.

Note: I didn't explain why the equations in the definition of adjoint
functors --- which I didn't write down --- imply the equations in the
definition of a simplicial object --- which I also didn't write down!

The point is, there's a more conceptual approach to understanding why
this stuff works. Remember from section \hyperlink{homotopy_K}{K} 
of last week that \(\Delta\)
is ``the free monoidal category on a monoid object''. This implies that
whenever we have a monoid object in a monoidal category \(\mathcal{M}\),
we get a monoidal functor \[F\colon\Delta\to\mathcal{M}.\] This gives a
functor \[G\colon\Delta^{\mathrm{op}}\to M^{\mathrm{op}}\] So: a monoid
object in \(\mathcal{M}\) gives a simplicial object in
\(\mathcal{M}^{\mathrm{op}}\).

Actually, if \(\mathcal{M}\) is a monoidal category,
\(\mathcal{M}^{\mathrm{op}}\) becomes one too, with the same tensor
product and unit object. So it's also true that a monoid object in
\(\mathcal{M}^{\mathrm{op}}\) gives a simplicial object in
\(\mathcal{M}\)!

Another name for a monoid object in \(\mathcal{M}^{\mathrm{op}}\) is a
``comonoid object in \(\mathcal{M}\)''. Remember,
\(\mathcal{M}^{\mathrm{op}}\) is just like \(\mathcal{M}\) but with all
the arrows turned around. So if we've got a monoid object in
\(\mathcal{M}^{\mathrm{op}}\), it gives us a similar gadget in
\(\mathcal{M}\), but with all the arrows turned around. More precisely,
a comonoid object in \(\mathcal{M}\) is an object, say \(m\), with
``coproduct'' \[c\colon m\to m\otimes m\] and ``counit''
\[e\colon m\to 1\] morphisms, satisfying ``coassociativity'' and the
left and right ``counit laws''. You get these laws by taking
associativity and the left/right unit laws, writing them out as
commutative diagrams, and turning all the arrows around.

So: a comonoid object in a monoidal category \(\mathcal{M}\) gives a
simplicial object in \(\mathcal{M}\). Now let's see how this is related
to adjoint functors. Suppose we have an adjunction, so we have some
functors \[
  \begin{tikzcd}
    \mathcal{C} \rar[bend left=30,"L"]
    & \mathcal{D} \lar[bend left=30,"R"]
  \end{tikzcd}
\] and natural transformations \[
  \begin{aligned}
    e&\colon LR\Rightarrow 1_{\mathcal{D}}
  \\i&\colon 1_{\mathcal{C}} \Rightarrow RL
  \end{aligned}
\] satisfying the same equations I didn't write before.

Let \(\operatorname{Hom}(\mathcal{C},\mathcal{C})\) be the category
whose objects are functors from \(\mathcal{C}\) to itself and whose
morphisms are natural transformations between such functors. This is a
monoidal category, since we can compose functors from \(\mathcal{C}\) to
itself. In \protect\hyperlink{week92}{``Week 92''} I showed that
\(\operatorname{Hom}(\mathcal{C},\mathcal{C})\) has a monoid object in
it, namely \(RL\). The product for this monoid object is
\[R\cdot e\cdot L\colon RLRL \Rightarrow RL\] and the unit is
\[i\colon 1_{\mathcal{C}} \Rightarrow RL\] Folks often call this sort of
thing a ``monad''.

Similarly, \(\operatorname{Hom}(\mathcal{D},\mathcal{D})\) is a monoidal
category containing a comonoid object, namely \(LR\). The coproduct for
this comonoid object is \[L\cdot i\cdot R\colon LR \Rightarrow LRLR\]
and the counit is \[e\colon LR \Rightarrow 1_{\mathcal{D}}\] People call
this thing a ``comonad''. But what matters here is that we've seen this
comonoid object automatically gives us a simplicial object in
\(\operatorname{Hom}(\mathcal{D},\mathcal{D})\)! If we pick any object
\(d\) of \(\mathcal{D}\), we get a functor
\[\operatorname{Hom}(\mathcal{D},\mathcal{D})\to\mathcal{D}\] by taking
\[\operatorname{Hom}(\mathcal{D},\mathcal{D})\times\mathcal{D}\to \mathcal{D}\]
and plugging in \(d\) in the second argument. This functor lets us push
our simplicial object in \(\operatorname{Hom}(\mathcal{D},\mathcal{D})\)
forwards to a simplicial object in \(\mathcal{D}\). Voila!

\hypertarget{week119}{%
\section{April 13, 1998}\label{week119}}

I've been slacking off on This Week's Finds lately because I was busy
getting stuff done at Riverside so that I could visit the Center for
Gravitational Physics and Geometry here at Penn State with a fairly
clean slate. Indeed, sometimes my whole life seems like an endless
series of distractions designed to prevent me from writing This Week's
Finds. However, now I'm here and ready to have some fun\ldots.

Recently I've been trying to learn about grand unified theories, or
``GUTs''. These were popular in the late 1970s and early 1980s, when the
Standard Model of particle interactions had fully come into its own and
people were looking around for a better theory that would unify all the
forces and particles present in that model - in short, everything except
gravity.

The Standard Model works well but it's fairly baroque, so it's natural
to hope for some more elegant theory underlying it. Remember how it
goes:

\begin{longtable}[]{@{}lll@{}}
\caption*{Gauge bosons}\tabularnewline
\toprule
Electromagnetic force & Weak force & Strong force\tabularnewline
\midrule
\endfirsthead
\toprule
Electromagnetic force & Weak force & Strong force\tabularnewline
\midrule
\endhead
photon & W\textsubscript{+}, W\textsubscript{-}, Z & 8
gluons\tabularnewline
\bottomrule
\end{longtable}

\begin{longtable}[]{@{}llll@{}}
\caption*{Fermions}\tabularnewline
\toprule
\textbf{Leptons} & & \textbf{Quarks} &\tabularnewline
\midrule
\endfirsthead
\toprule
\textbf{Leptons} & & \textbf{Quarks} &\tabularnewline
\midrule
\endhead
electron & electron neutrino & down quark & up quark\tabularnewline
muon & muon neutrino & strange quark & charm quark\tabularnewline
tauon & tauon neutrino & bottom quark & top quark\tabularnewline
\bottomrule
\end{longtable}

\begin{longtable}[]{@{}c@{}}
\caption*{Higgs boson (not yet seen)}\tabularnewline
\toprule
?\tabularnewline
\endhead
\bottomrule
\end{longtable}

The strong, electromagnetic and weak forces are all described by
Yang--Mills fields, with the gauge group
\(\mathrm{SU}(3)\times\mathrm{SU}(2)\times\mathrm{U}(1)\). In what
follows I'll assume you know the rudiments of gauge theory, or at least
that you can fake it.

\(\mathrm{SU}(3)\) is the gauge group of the strong force, and its 8
generators correspond to the gluons.
\(\mathrm{SU}(2)\times\mathrm{U}(1)\) is the gauge group of the
electroweak force, which unifies electromagnetism and the weak force.
It's \emph{not} true that the generators of \(\mathrm{SU}(2)\)
corresponds to the W\textsubscript{+}, W\textsubscript{-} and Z while
the generator of \(\mathrm{U}(1)\) corresponds to the photon. Instead,
the photon corresponds to the generator of a sneakier \(\mathrm{U}(1)\)
subgroup sitting slantwise inside \(\mathrm{SU}(2)\times\mathrm{U}(1)\);
the basic formula to remember here is: \[Q = I_3 + Y/2\] where \(Q\) is
ordinary electric charge, \(I_3\) is the 3rd component of ``weak
isospin'', i.e.~the generator of \(\mathrm{SU}(2)\) corresponding to the
matrix \[
  \left(
    \begin{array}{cc}
      \frac12&0\\0&-\frac12
    \end{array}
  \right)
\] and \(Y\), ``hypercharge'', is the generator of the \(\mathrm{U}(1)\)
factor. The role of the Higgs particle is to spontaneously break the
\(\mathrm{SU}(2)\times\mathrm{U}(1)\) symmetry, and also to give all the
massive particles their mass. However, I don't want to talk about that
here; I want to focus on the fermions and how they form representations
of the gauge group
\(\mathrm{SU}(3)\times\mathrm{SU}(2)\times\mathrm{U}(1)\), because I
want to talk about how grand unified theories attempt to simplify this
picture --- at the expense of postulating more Higgs bosons.

The fermions come in 3 generations, as indicated in the chart above. I
want to explain how the fermions in a given generation are grouped into
irreducible representations of
\(\mathrm{SU}(3)\times\mathrm{SU}(2)\times\mathrm{U}(1)\). All the
generations work the same way, so I'll just talk about the first
generation. Also, every fermion has a corresponding antiparticle, but
this just transforms according to the dual representation, so I will
ignore the antiparticles here.

Before I tell you how it works, I should remind you that all the
fermions are, in addition to being representations of
\(\mathrm{SU}(3)\times\mathrm{SU}(2)\times\mathrm{U}(1)\), also
spin-\(1/2\) particles. The massive fermions --- the quarks and the
electron, muon and tauon --- are Dirac spinors, meaning that they can
spin either way along any axis. The massless fermions --- the neutrinos
--- are Weyl spinors, meaning that they always spin counterclockwise
along their axis of motion. This makes sense because, being massless,
they move at the speed of light, so everyone can agree on their axis of
motion! So the massive fermions have two helicity states, which we'll
refer to as ``left-handed'' and ``right-handed'', while the neutrinos
only come in a ``left-handed'' form.

(Here I am discussing the Standard Model in its classic form. I'm
ignoring any modifications needed to deal with a possible nonzero
neutrino mass. For more on Standard Model, neutrino mass and different
kinds of spinors, see \protect\hyperlink{week93}{``Week 93''}.)

Okay. The Standard Model lumps the left-handed neutrino and the
left-handed electron into a single irreducible representation of
\(\mathrm{SU}(3)\times\mathrm{SU}(2)\times\mathrm{U}(1)\):
\[(\nu_L,\mathrm{e}_L) \qquad\qquad (1,2,-1)\] This \(2\)-dimensional
representation is called \((1,2,-1)\), meaning that it's the tensor
product of the \(1\)-dimensional trivial rep of \(\mathrm{SU}(3)\), the
2-dimensional fundamental rep of \(\mathrm{SU}(2)\), and the
\(1\)-dimensional rep of \(\mathrm{U}(1)\) with hypercharge \(-1\).

Similarly, the left-handed up and down quarks fit together as:
\[(\mathrm{u}_L, \mathrm{u}_L, \mathrm{u}_L, \mathrm{d}_L, \mathrm{d}_L, \mathrm{d}_L) \qquad\qquad (3,2,1/3)\]
Here I'm writing both quarks 3 times since they also come in 3 color
states. In other words, this \(6\)-dimensional representation is the
tensor product of the \(3\)-dimensional fundamental rep of
\(\mathrm{SU}(3)\), the \(2\)-dimensional fundamental rep of
\(\mathrm{SU}(2)\), and the \(1\)-dimensional rep of \(\mathrm{U}(1)\)
with hypercharge \(1/3\). That's why we call this rep \((3,2,1/3)\).

(If you are familiar with the irreducible representations of
\(\mathrm{U}(1)\) you will know that they are usually parametrized by
integers. Here we are using integers divided by \(3\). The reason is
that people defined the charge of the electron to be \(-1\) before
quarks were discovered, at which point it turned out that the smallest
unit of charge was \(1/3\) as big as had been previously believed.)

The right-handed electron stands alone in a \(1\)-dimensional rep, since
there is no right-handed neutrino:
\[\mathrm{e}_R \qquad\qquad (1,1,-2).\] Similarly, the right-handed up
quark stands alone in a \(3\)-dimensional rep, as does the right-handed
down quark:
\[(\mathrm{u}_R,\mathrm{u}_R,\mathrm{u}_R) \qquad\qquad (3,1,4/3)\]
\[(\mathrm{d}_R,\mathrm{d}_R,\mathrm{d}_R) \qquad\qquad (3,1,-2/3)\]
That's it. If you want to study this stuff, try using the formula
\[Q = I_3 + Y/2\] to figure out the charges of all these particles. For
example, since the right-handed electron transforms in the trivial rep
of \(\mathrm{SU}(2)\), it has \(I_3 = 0\), and if you look up there
you'll see that it has \(Y = -2\). This means that its electric charge
is \(Q = -1\), as we already knew.

Anyway, we obviously have a bit of a mess on our hands! The Standard
Model is full of tantalizing patterns, but annoyingly complicated. The
idea of grand unified theories is to find a pattern lurking in all this
data by fitting the group
\(\mathrm{SU}(3)\times\mathrm{SU}(2)\times\mathrm{U}(1)\) into a larger
group. The smallest-dimensional ``simple'' Lie group that works is
\(\mathrm{SU}(5)\). Here ``simple'' is a technical term that eliminates,
for example, groups that are products of other groups --- these aren't
very ``unified''. Georgi and Glashow came up with their ``minimal''
\(\mathrm{SU}(5)\) grand unified theory in 1975. The idea is to stick
\(\mathrm{SU}(3)\times\mathrm{SU}(2)\) into \(\mathrm{SU}(5)\) in the
obvious diagonal way, leaving just enough room to cram in the
\(\mathrm{U}(1)\) if you are clever.

Now if you add up the dimensions of all the representations above you
get \(2 + 6 + 1 + 3 + 3 = 15\). This means we need to find a
\(15\)-dimensional representation of \(\mathrm{SU}(5)\) to fit all these
particles. There are various choices, but only one that really works
when you take all the physics into account. For a nice simple account of
the detective work needed to figure this out, see:

\begin{enumerate}
\def\labelenumi{\arabic{enumi})}
\tightlist
\item
  Edward Witten, ``Grand unification with and without supersymmetry'',
  in \emph{Introduction to supersymmetry in particle and nuclear
  physics}, edited by O. Castanos, A. Frank, L. Urrutia, Plenum Press,
  1984.
\end{enumerate}

I'll just give the answer. First we take the \(5\)-dimensional
fundamental representation of \(\mathrm{SU}(5)\) and pack fermions in as
follows:
\[(\mathrm{d}_R, \mathrm{d}_R, \mathrm{d}_R, \mathrm{e}^+_R, \bar{\nu}_R) \qquad\qquad 5 = (3,1,-2/3) + (1,2,-1)\]

Here \(\mathrm{e}^+_R\) is the right-handed positron and \(\bar{\nu}_R\)
is the right-handed antineutrino --- curiously, we need to pack some
antiparticles in with particles to get things to work out right. Note
that the first 3 particles in the above list, the 3 states of the
right-handed down quark, transform according to the fundamental rep of
\(\mathrm{SU}(3)\) and the trivial rep of \(\mathrm{SU}(2)\), while the
remaining two transform according to the trivial rep of
\(\mathrm{SU}(3)\) and the fundamental rep of \(\mathrm{SU}(2)\). That's
how it has to be, given how we stuffed
\(\mathrm{SU}(3)\times\mathrm{SU}(2)\) into \(\mathrm{SU}(5)\).

Note also that the charges of the 5 particles on this list add up to
zero. That's also how it has to be, since the generators of
\(\mathrm{SU}(5)\) are traceless. Note that the down quark must have
charge \(-1/3\) for this to work! In a sense, the \(\mathrm{SU}(5)\)
model says that quarks \emph{must} have charges in units of \(1/3\),
because they come in 3 different colors! This is pretty cool.

Then we take the \(10\)-dimensional representation of \(\mathrm{SU}(5)\)
given by the 2nd exterior power of the fundamental representation ---
i.e., antisymmetric \(5\times5\) matrices - and pack the rest of the
fermions in like this: \[
  \left(
    \begin{array}{ccccc}
      0 & \bar{\mathrm{u}}_L & \bar{\mathrm{u}}_L & \mathrm{u}_L & \mathrm{d}_L
    \\-\bar{\mathrm{u}_L} & 0 & \bar{\mathrm{u}_L} & \mathrm{u}_L & \mathrm{d}_L
    \\-\bar{\mathrm{u}_L} & -\bar{\mathrm{u}_L} & 0 & \mathrm{u}_L & \mathrm{d}_L
    \\-\mathrm{u}_L & -\mathrm{u}_L & -\mathrm{u}_L & 0 & \mathrm{e}^+_L
    \\-\mathrm{d}_L & -\mathrm{u}_L & -\mathrm{d}_L & \mathrm{e}^+_L & 0
    \end{array}
  \right)
  \qquad\qquad
  10 = (3,2,1/3) + (1,1,2) + (3,1,-4/3)
\] Here the \(\bar{\mathrm{u}}\) is the antiparticle of the up quark ---
again we've needed to use some antiparticles. However, you can easily
check that these two representations of \(\mathrm{SU}(5)\) together with
their duals account for all the fermions and their antiparticles.

The \(\mathrm{SU}(5)\) theory has lots of nice features. As I already
noted, it explains why the up and down quarks have charges \(2/3\) and
\(-1/3\), respectively. It also gives a pretty good prediction of
something called the Weinberg angle, which is related to the ratio of
the masses of the W and Z bosons. It also makes testable new
predictions! Most notably, since it allows quarks to turn into leptons,
it predicts that protons can decay --- with a halflife of somewhere
around \(10^{29}\) or \(10^{30}\) years. So people set off to look for
proton decay\ldots.

However, even when the \(\mathrm{SU}(5)\) model was first proposed, it
was regarded as slightly inelegant, because it didn't unify all the
fermions of a given generation in a \emph{single} irreducible
representation (together with its dual, for antiparticles). This is one
reason why people began exploring still larger gauge groups. In 1975
Georgi, and independently Fritzsch and Minkowski, proposed a model with
gauge group \(\mathrm{SO}(10)\). You can stuff \(\mathrm{SU}(5)\) into
\(\mathrm{SO}(10)\) as a subgroup in such a way that the 5- and
10-dimensional representations of \(\mathrm{SU}(5)\) listed above both
fit into a single \(16\)-dimensional rep of \(\mathrm{SO}(10)\), namely
the chiral spinor rep. Yes, 16, not 15 --- that wasn't a typo! The
\(\mathrm{SO}(10)\) theory predicts that in addition to the 15 states
listed above there is a 16th, corresponding to a right-handed neutrino!
I'm not sure yet how the recent experiments indicating a nonzero
neutrino mass fit into this business, but it's interesting.

Somewhere around this time, people noticed something interesting about
these groups we've been playing with. They all fit into the
``\(\mathrm{E}\) series''!

I don't have the energy to explain Dynkin diagrams and the
\(\mathrm{ABCDEFG}\) classification of simple Lie groups here, but
luckily I've already done that, so you can just look at
\protect\hyperlink{week62}{``Week 62''} ---
\protect\hyperlink{week65}{``Week 65''} to learn about that. The point
is, there is an infinite series of simple Lie groups associated to
rotations in real vector spaces --- the \(\mathrm{SO}(n)\) groups, also
called the \(\mathrm{B}\) and \(\mathrm{D}\) series. There is an
infinite series of them associated to rotations in complex vector spaces
--- the \(\mathrm{SU}(n)\) groups, also called the \(\mathrm{A}\)
series. And there is infintie series of them associated to rotations in
quaternionic vector spaces --- the \(\mathrm{Sp}(n)\) groups, also
called the \(\mathrm{C}\) series. And there is a ragged band of 5
exceptions which are related to the octonions, called \(mathrm{G}_2\),
\(\mathrm{F}_4\), \(\mathrm{E}_6\), \(\mathrm{E}_7\), and
\(\mathrm{E}_8\). I'm sort of fascinated by these --- see
\protect\hyperlink{week90}{``Week 90''},
\protect\hyperlink{week91}{``Week 91''}, and
\protect\hyperlink{week106}{``Week 106''} for more --- so I was
extremely delighted to find that the \(\mathrm{E}\) series plays a
special role in grand unified theories.

Now, people usually only talk about \(\mathrm{E}_6\), \(\mathrm{E}_7\),
and \(\mathrm{E}_8\), but one can work backwards using Dynkin diagrams
to define \(\mathrm{E}_5\), \(\mathrm{E}_4\), \(\mathrm{E}_3\),
\(\mathrm{E}_2\), and \(\mathrm{E}_1\). Let's do it! Thanks go to Allan
Adler and Robin Chapman for helping me understand how this works\ldots.

\(\mathrm{E}_8\) is a big fat Lie group whose Dynkin diagram looks like
this: \[
  \begin{tikzpicture}
    \draw[thick] (0,0) node{$\bullet$} to (1,0) node{$\bullet$} to (2,0) node{$\bullet$} to (3,0) node {$\bullet$} to (4,0) node {$\bullet$} to (5,0) node {$\bullet$} to (6,0) node {$\bullet$};
    \draw[thick] (2,0) to (2,1) node{$\bullet$};
  \end{tikzpicture}
\] If we remove the rightmost root, we obtain the Dynkin diagram of a
subgroup called \(\mathrm{E}_7\): \[
  \begin{tikzpicture}
    \draw[thick] (0,0) node{$\bullet$} to (1,0) node{$\bullet$} to (2,0) node{$\bullet$} to (3,0) node {$\bullet$} to (4,0) node {$\bullet$} to (5,0) node {$\bullet$};
    \draw[thick] (2,0) to (2,1) node{$\bullet$};
  \end{tikzpicture}
\] If we again remove the rightmost root, we obtain the Dynkin diagram
of a subgroup of \(\mathrm{E}_7\), namely \(\mathrm{E}_6\): \[
  \begin{tikzpicture}
    \draw[thick] (0,0) node{$\bullet$} to (1,0) node{$\bullet$} to (2,0) node{$\bullet$} to (3,0) node {$\bullet$} to (4,0) node {$\bullet$};
    \draw[thick] (2,0) to (2,1) node{$\bullet$};
  \end{tikzpicture}
\] This was popular as a gauge group for grand unified models, and the
reason why becomes clear if we again remove the rightmost root,
obtaining the Dynkin diagram of a subgroup we could call
\(\mathrm{E}_5\): \[
  \begin{tikzpicture}
    \draw[thick] (0,0) node{$\bullet$} to (1,0) node{$\bullet$} to (2,0) node{$\bullet$} to (3,0) node {$\bullet$};
    \draw[thick] (2,0) to (2,1) node{$\bullet$};
  \end{tikzpicture}
\] But this is really just good old \(\mathrm{SO}(10)\), which we were
just discussing! And if we yet again remove the rightmost root, we get
the Dynkin diagram of a subgroup we could call \(\mathrm{E}_4\): \[
  \begin{tikzpicture}
    \draw[thick] (0,0) node{$\bullet$} to (1,0) node{$\bullet$} to (2,0) node{$\bullet$};
    \draw[thick] (2,0) to (2,1) node{$\bullet$};
  \end{tikzpicture}
\]

\begin{verbatim}
      o      
      |      
o--o--o
\end{verbatim}

This is just \(\mathrm{SU}(5)\)! Let's again remove the rightmost root,
obtaining the Dynkin diagram for \(\mathrm{E}_3\). Well, it may not be
clear what counts as the rightmost root, but here's what I want to get
when I remove it: \[
  \begin{tikzpicture}
    \draw[thick] (0,0) node{$\bullet$} to (1,0) node{$\bullet$};
    \node at (2,1) {$\bullet$};
  \end{tikzpicture}
\] This is just \(\mathrm{SU}(3)\times\mathrm{SU}(2)\), sitting inside
\(\mathrm{SU}(5)\) in the way we just discussed! So for some mysterious
reason, the Standard Model and grand unified theories seem to be related
to the \(\mathrm{E}\) series!

We could march on and define \(\mathrm{E}_2\): \[
  \begin{tikzpicture}
    \node at (0,0) {$\bullet$};
    \node at (2,1) {$\bullet$};
  \end{tikzpicture}
\] which is just \(\mathrm{SU}(2)\times\mathrm{SU}(2)\), and
\(\mathrm{E}_1\): \[
  \begin{tikzpicture}
    \node at (0,0) {$\bullet$};
  \end{tikzpicture}
\] which is just \(\mathrm{SU}(2)\)\ldots{} but I'm not sure what's so
great about these groups.

By the way, you might wonder what's the real reason for removing the
roots in the order I did --- apart from getting the answers I wanted to
get --- and the answer is, I don't really know! If anyone knows, please
tell me. This could be an important clue.

Now, this stuff about grand unified theories and the \(\mathrm{E}\)
series is one of the reasons why people like string theory, because
heterotic string theory is closely related to \(\mathrm{E}_8\) (see
\protect\hyperlink{week95}{``Week 95''}). However, I must now tell you
the \emph{bad} news about grand unified theories. And it is \emph{very}
bad.

The bad news is that those people who went off to detect proton decay
never found it! It became clear in the mid-1980s that the proton
lifetime was at least \(10^{32}\) years or so, much larger than what the
\(\mathrm{SU}(5)\) theory most naturally predicts. Of course, if one is
desperate to save a beautiful theory from an ugly fact, one can resort
to desperate measures. For example, one can get the \(\mathrm{SU}(5)\)
model to predict very slow proton decay by making the grand unification
mass scale large. Unfortunately, then the coupling constants of the
strong and electroweak forces don't match at the grand unification mass
scale. This became painfully clear as better measurements of the strong
coupling constant came in.

Theoretical particle physics never really recovered from this crushing
blow. In a sense, particle physics gradually retreated from the goal of
making testable predictions, drifting into the wonderland of pure
mathematics\ldots{} first supersymmetry, then supergravity, and then
superstrings\ldots{} ever more elegant theories, but never yet a
verified experimental prediction. Perhaps we should be doing something
different, something better? Easy to say, hard to do! If we see a
superpartner at CERN, a lot of this ``superthinking'' will be vindicated
--- so I guess most particle physicists are crossing their fingers and
praying for this to happen.

The following textbook on grand unified theories is very nice,
especially since it begins with a review of the Standard Model:

\begin{enumerate}
\def\labelenumi{\arabic{enumi})}
\setcounter{enumi}{1}
\tightlist
\item
  Graham G. Ross, \emph{Grand Unified Theories}, Benjamin-Cummings,
  1984.
\end{enumerate}
\noindent
This one is a bit more idiosyncratic, but also good --- Mohapatra is
especially interested in theories where CP violation arises via
spontaneous symmetry breaking:

\begin{enumerate}
\def\labelenumi{\arabic{enumi})}
\setcounter{enumi}{2}
\tightlist
\item
  Ranindra N. Mohapatra,  \emph{Unification and Supersymmetry:the
  Frontiers of Quark-Lepton Physics}, Springer, Berlin, 1992.
\end{enumerate}
\noindent
I also found the following articles interesting:

\begin{enumerate}
\def\labelenumi{\arabic{enumi})}
\setcounter{enumi}{3}
\item
  D. V. Nanopoulos, ``Tales of the GUT age'', in \emph{Grand Unified
  Theories and Related Topics, proceedings of the 4th Kyoto Summer
  Institute}, World Scientific, Singapore, 1981.
\item
  P. Ramond, ``Grand unification'', in \emph{Grand Unified Theories and
  Related Topics, proceedings of the 4th Kyoto Summer Institute}, World
  Scientific, Singapore, 1981.
\end{enumerate}

Okay, now for some homotopy theory! I don't think I'm ever gonna get to
the really cool stuff\ldots{} in my attempt to explain everything
systematically, I'm getting worn out doing the preliminaries. Oh well,
on with it\ldots{} now it's time to start talking about loop spaces!
These are really important, because they tie everything together.
However, it takes a while to deeply understand their importance.

\begin{center}\rule{0.5\linewidth}{0.5pt}\end{center}

\hypertarget{homotopy_N}{\textbf{N.}}
\emph{The loop space of a topological space.} Suppose we
have a ``pointed space'' \(X\), that is, a topological space with a
distinguished point called the ``basepoint''. Then we can form the space
\(LX\) of all ``based loops'' in \(X\) --- loops that start and end at
the basepoint.

One reason why \(LX\) is so nice is that its homotopy groups are the
same as those of \(X\), but shifted: \[\pi_i(LX) = \pi_{i+1}(X).\]
Another reason \(LX\) is nice is that it's almost a topological group,
since one can compose based loops, and every loop has an ``inverse''.
However, one must be careful here! Unless one takes special care,
composition will only be associative up to homotopy, and the ``inverse''
of a loop will only be the inverse up to homotopy.

Actually we can make composition strictly associative if we work with
``Moore paths''. A Moore path in \(X\) is a continuous map
\[f\colon[0,T]\to X\] where \(T\) is an arbitrary nonnegative real
number. Given a Moore path \(f\) as above and another Moore path
\[g\colon[0,S]\to X\] which starts where \(f\) ends, we can compose them
in an obvious way to get a Moore path \[fg\colon[0,T+S]\to X\] Note that
this operation is associative ``on the nose'', not just up to homotopy.
If we define \(LX\) using Moore paths that start and end at the
basepoint, we can easily make \(LX\) into a topological monoid --- that
is, a topological space with a continuous associative product and a unit
element. (If you've read section \hyperlink{homotopy_L}{L} of
\protect\hyperlink{week117}{``Week 117''}) you'll know 
this is just a monoid
object in \(\mathsf{Top}\)!) In particular, the unit element of \(LX\)
is the path \(i\colon[0,0]\to X\) that just sits there at the basepoint
of \(X\).

\(LX\) is not a topological group, because even Moore paths don't have
strict inverses. But \(LX\) is \emph{close} to being a group. We can
make this fact precise in various ways, some more detailed than others.
I'm pretty sure one way to say it is this: the natural map from \(LX\)
to its ``group completion'' is a homotopy equivalence.

\begin{center}\rule{0.5\linewidth}{0.5pt}\end{center}

\hypertarget{homotopy_O}{\textbf{O.}}
\emph{The group completion of a topological monoid.} Let
\(\mathsf{TopMon}\) be the category of topological monoids and let
\(\mathsf{TopGp}\) be the category of topological groups. There is a
forgetful functor \[F\colon\mathsf{TopGp}\to\mathsf{TopMon}\] and this
has a left adjoint \[G\colon\mathsf{TopMon}\to\mathsf{TopGp}\] which
takes a topological monoid and converts it into a topological group by
throwing in formal inverses of all the elements and giving the resulting
group a nice topology. This functor \(G\) is called ``group completion''
and was first discussed by Quillen (in the simplicial context, in an
unpublished paper), and independently by Barratt and Priddy:

\begin{enumerate}
\def\labelenumi{\arabic{enumi})}
\setcounter{enumi}{5}
\tightlist
\item
  M. G. Barratt and S. Priddy, ``On the homology of non-connected
  monoids and their associated groups'', \emph{Comm. Math. Helv.}
  \textbf{47} (1972), 1--14.
\end{enumerate}
\noindent
For any topological monoid \(M\), there is a natural map from \(M\) to
\(F(G(M))\), thanks to the miracle of adjoint functors. This is the
natural map I'm talking about in the previous section!

\hypertarget{week120}{%
\section{May 6, 1998}\label{week120}}

Now that I'm hanging out with the gravity crowd at Penn State, I might
as well describe what's been going on here lately.

First of all, Ashtekar and Krasnov have written an expository account of
their work on the entropy of quantum black holes:

\begin{enumerate}
\def\labelenumi{\arabic{enumi})}
\tightlist
\item
  Abhay Ashtekar and Kirill Krasnov, ``Quantum geometry and black
  holes'', available as
  \href{https://arxiv.org/abs/gr-qc/9804039}{\texttt{gr-qc/9804039}}.
\end{enumerate}
\noindent
But if you prefer to see a picture of a quantum black hole without any
equations, try this:
\[\includegraphics[max width=0.35\linewidth]{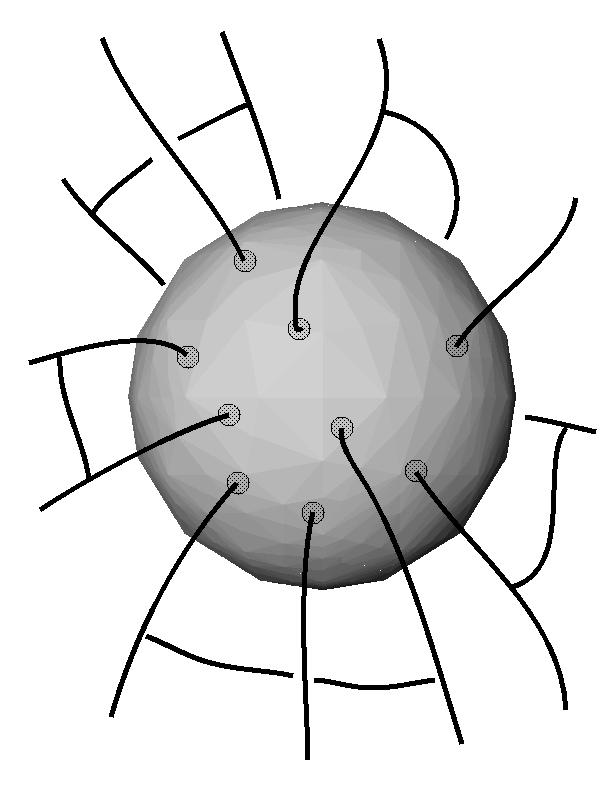}\] 
This shows a bunch of spin networks poking the horizon, giving it area
and curvature. Of course, this is just a theory.

Second, there's been a burst of new work studying quantum gravity in
terms of spin foams. A spin foam looks a bit like a bunch of soap suds
--- with the faces of the bubbles and the edges where the bubbles meet
labelled by spins \(j = 0, 1/2, 1, 3/2, \ldots\). Spin foams are an
attempt at a quantum description of the geometry of spacetime. If you
slice a spin foam with a hyperplane representing ``\(t = 0\)'' you get a
spin network: a graph with its edges and vertices labelled by spins.
Spin networks have been used in quantum gravity for a while now to
describe the geometry of space at a given time, so it's natural to hope
that they're a slice of something that describes the geometry of
spacetime.

As usual in quantum gravity, it's too early to tell if this approach
will work. As usual, it has lots of serious problems. But before going
into the problems, let me remind you how spin foams are supposed to
work.

To relate spin foams to more traditional ideas about spacetime, one can
consider spin foams living in a triangulated 4-manifold: one spin foam
vertex sitting in each \(4\)-simplex, one spin foam edge poking through
each tetrahedron, and one spin foam face intersecting each triangle.
Labelling the spin foam edges and faces with spins is supposed to endow
the triangulated 4-manifold with a ``quantum 4-geometry''. In other
words, it should let us compute things like the areas of the triangles,
the volumes of the tetrahedra, and the 4-volumes of the \(4\)-simplices.
There are some arguments going on now about the right way to do this,
but it's far from an arbitrary business: the interplay between group
representation theory and geometry says a lot about how it should go. In
the simplified case of \(3\)-dimensional spacetime, it's fairly well
understood --- the hard part, and the fun part, is getting it to work in
4 dimensions.

Assuming we can do this, the next trick is to compute an amplitude for
each spin foam vertex in a nice way, much as one computes amplitudes for
vertices of Feynman diagrams. A spin foam vertex is supposed to
represent an ``event'' --- if we slice the spin foam by a hyperplane we
get a spin network, and as we slide this slice ``forwards in time'', the
spin network changes its topology whenever we pass a spin foam vertex.
The amplitude for a vertex tells us how likely it is for this event to
happen. As usual in quantum theory, we need to take the absolute value
of an amplitude and square it to get a probability.

We also need to compute amplitudes for spin foam edges and faces, called
``propagators'', in analogy to the amplitudes one computes for the edges
of Feynman diagrams. Multiplying all the vertex amplitudes and
propagators for a given spin foam, one gets the amplitude for the whole
spin foam. This tells us how likely it is for the whole spin foam to
happen.

Barrett and Crane came up with a specific way to do all this stuff,
Reisenberger came up with a different way, I came up with a general
formalism for understanding this stuff, and now people are busy arguing
about the merits of different approaches. Here are some papers on the
subject --- I'll pick up where I left off in
\protect\hyperlink{week113}{``Week 113''}.

\begin{enumerate}
\def\labelenumi{\arabic{enumi})}
\setcounter{enumi}{1}
\tightlist
\item
  Louis Crane, David N. Yetter, ``On the classical limit of the balanced
  state sum'', available as
  \href{https://arxiv.org/abs/gr-qc/9712087}{\texttt{gr-qc/9712087}}.
\end{enumerate}
\noindent
The goal here is to show that in the limit of large spins, the amplitude
given by Barrett and Crane's formula approaches \[\exp(iS)\] where \(S\)
is the action for classical general relativity --- suitably discretized,
and in signature \(++++\). The key trick is to use an idea invented by
Regge in 1961.

Regge came up with a discrete analog of the usual formula for the action
in classical general relativity. His formula applies to a triangulated
4-manifold whose edges have specified lengths. In this situation, each
triangle has an ``angle deficit'' associated to it. It's easier to
visualize this two dimensions down, where each vertex in a triangulated
2-manifold has an angle deficit given by adding up angles for all the
triangles having it as a corner, and then subtracting \(2\pi\). No angle
deficit means no curvature: the triangles sit flat in a plane. The idea
works similarly in 4 dimensions. Here's Regge's formula for the action:
take each triangle in your triangulated 4-manifold, take its area,
multiply it by its angle deficit, and then sum over all the triangles.

Simple, huh? In the continuum limit, Regge's action approaches the
integral of the Ricci scalar curvature --- the usual action in general
relativity. For more see:

\begin{enumerate}
\def\labelenumi{\arabic{enumi})}
\setcounter{enumi}{2}
\tightlist
\item
  T. Regge, ``General relativity without coordinates'', \emph{Nuovo
  Cimento} \textbf{19} (1961), 558--571.
\end{enumerate}
\noindent
So, Crane and Yetter try to show that in the limit of large spins, the
Barrett--Crane spin foam amplitude approaches \(\exp(iS)\) where \(S\) is
the Regge action. There argument is interesting but rather sketchy.
Someone should try to fill in the details!

However, it's not clear to me that the large spin limit is physically
revelant. If spacetime is really made of lots of \(4\)-simplices
labelled by spins, the \(4\)-simplices have got to be quite small, so
the spins labelling them should be fairly small. It seems to me that the
right limit to study is the limit where you triangulate your 4-manifold
with a huge number of \(4\)-simplices labelled by fairly small spins.
After all, in the spin network picture of the quantum black hole, it
seems that spin network edges labelled by spin \(1/2\) contribute most
of the states (see \protect\hyperlink{week112}{``Week 112''}).

When you take a spin foam living in a triangulated 4-manifold and slice
it in a way that's compatible with the triangulation, the spin network
you get is a 4-valent graph. Thus it's not surprising that Barrett and
Crane's formula for vertex amplitudes is related to an invariant of
4-valent graphs with edges labelled by spins. There's already a branch
of math relating such invariants to representations of groups and
quantum groups, and their formula fits right in. Yetter has figured out
how to generalize this graph invariant to \(n\)-valent graphs with edges
labelled by spins, and he's also studied more carefully what happens
when one ``\(q\)-deforms'' the whole business --- replacing the group by
the corresponding quantum group. This should be related to quantum
gravity with nonzero cosmological constant, if all the mathematical
clues aren't lying to us. See:

\begin{enumerate}
\def\labelenumi{\arabic{enumi})}
\setcounter{enumi}{3}
\tightlist
\item
  David N. Yetter, ``Generalized Barrett--Crane vertices and invariants
  of embedded graphs'', available as
  \href{https://arxiv.org/abs/math.QA/9801131}{\texttt{math.QA/9801131}}.
\end{enumerate}
\noindent
Barrett has also given a nice formula in terms of integrals for the
invariant of 4-valent graphs labelled by spins. This is motivated by the
physics and illuminates it nicely:

\begin{enumerate}
\def\labelenumi{\arabic{enumi})}
\setcounter{enumi}{4}
\tightlist
\item
  John W. Barrett, ``The classical evaluation of relativistic spin
  networks'', available as
  \href{https://arxiv.org/abs/math.QA/9803063}{\texttt{math.QA/9803063}}.
\end{enumerate}
\noindent
Let me quote the abstract:

\begin{quote}
The evaluation of a relativistic spin network for the classical case of
the Lie group \(\mathrm{SU}(2)\) is given by an integral formula over
copies of \(\mathrm{SU}(2)\). For the graph determined by a
\(4\)-simplex this gives the evaluation as an integral over a space of
geometries for a \(4\)-simplex.
\end{quote}

Okay, so much for the good news. What about the bad news? To explain
this I need to get a bit more specific about Barrett and Crane's
approach.

Their approach is based on a certain way to describe the geometry of a
4-simplex. Instead of specifying lengths of edges as in the old Regge
approach, we specify bivectors for all its faces. Geometrically, a
bivector is just an ``oriented area element''; technically, the space of
bivectors is the dual of the space of \(2\)-forms. If we have a
\(4\)-simplex in R\^{}4 and we choose orientations for its triangular
faces, there's an obvious way to associate a bivector to each face. We
get 10 bivectors this way.

What constraints do these 10 bivectors satisfy? They can't be arbitrary!
First, for any four triangles that are all the faces of the same
tetrahedron, the corresponding bivectors must sum to zero. Second, every
bivector must be ``simple'' --- it must be the wedge product of two
vectors. Third, whenever two triangles are the faces of the same
tetrahedron, the sum of the corresponding bivectors must be simple.

It turns out that these constraints are almost but \emph{not quite
enough} to imply that 10 bivectors come from a \(4\)-simplex.
Generically, it there are four possibilities: our bivectors come from a
\(4\)-simplex, the \emph{negatives} of our bivectors come from a
\(4\)-simplex, their \emph{Hodge duals} come from a 4-simplex, or
\emph{the negatives of their Hodge duals} come from a 4-simplex.

If we ignore this and describe the \(4\)-simplex using bivectors
satisfying the three constraints above, and then quantize this
description, we get the picture of a ``quantum \(4\)-simplex'' that is
the starting-point for the Barrett--Crane model. But clearly it's
dangerous to ignore this problem.

Actually, I learned about this problem from Robert Bryant over on
sci.math.research, and I discussed it in my paper on spin foam models,
citing Bryant of course. Barrett and Crane overlooked this problem in
the first version of their paper, but now they recognize its importance.
Two papers have recently appeared which investigate it further:

\begin{enumerate}
\def\labelenumi{\arabic{enumi})}
\setcounter{enumi}{5}
\item
  Michael P. Reisenberger, ``Classical Euclidean general relativity from
  `left-handed area = right-handed area'\,'', available as
  \href{https://arxiv.org/abs/gr-qc/9804061}{\texttt{gr-qc/9804061}}.
\item
  Roberto De Pietri and Laurent Freidel, ``\(\mathfrak{so}(4)\)
  Plebanski Action and relativistic spin foam model'',
  available as
  \href{https://arxiv.org/abs/gr-qc/9804071}{\texttt{gr-qc/9804071}}.
\end{enumerate}
\noindent
These papers study classical general relativity formulated as a
constrained \(\mathrm{SO}(4)\) \(BF\) theory. The constraints needed
here are mathematically just the same as the constraints needed to
ensure that 10 bivectors come from the faces of an actual \(4\)-simplex!
This is part of the magic of this approach. But again, if one only
imposes the three constraints I listed above, it's not quite enough: one
gets fields that are either solutions of general relativity \emph{or}
solutions of three other theories! This raises the worry that the
Barrett--Crane model is a quantization, not exactly of general
relativity, but of general relativity mixed in with these extra
theories.

Here's another recent product of the Center for Classical and Quantum
Gravity here at Penn State:

\begin{enumerate}
\def\labelenumi{\arabic{enumi})}
\setcounter{enumi}{7}
\tightlist
\item
  Laurent Freidel and Kirill Krasnov, ``Discrete space-time volume for
  \(3\)-dimensional \(BF\) theory and quantum gravity'', 
  available as
  \href{https://arxiv.org/abs/hep-th/9804185}{\texttt{hep-th/9804185}}.
\end{enumerate}
\noindent
Freidel and Krasnov study the volume of a single \(3\)-simplex as an
observable in the context of the Turaev-Viro model --- a topological
quantum field theory which is closely related to quantum gravity in
spacetime dimension 3.

And here are some other recent papers on quantum gravity written by
folks who either work here at the CGPG or at least occasionally drift
through. I'll just quote the abstracts of these:

\begin{enumerate}
\def\labelenumi{\arabic{enumi})}
\setcounter{enumi}{8}
\tightlist
\item
  Ted Jacobson, ``Black hole thermodynamics today'', to appear in
  \emph{Proceedings of the Eighth Marcel Grossmann Meeting}, World
  Scientific, 1998.  Also available as
  \href{https://arxiv.org/abs/gr-qc/9801015}{\texttt{gr-qc/9801015}}.
\end{enumerate}

\begin{quote}
A brief survey of the major themes and developments of black hole
thermodynamics in the 1990's is given, followed by summaries of the
talks on this subject at MG8 together with a bit of commentary, and
closing with a look towards the future.
\end{quote}

\begin{enumerate}
\def\labelenumi{\arabic{enumi})}
\setcounter{enumi}{9}
\tightlist
\item
  Rodolfo Gambini and Jorge Pullin, ``Does loop quantum gravity imply
  \(\Lambda = 0\)?'', available as
  \href{https://arxiv.org/abs/gr-qc/9803097}{\texttt{gr-qc/9803097}}.
\end{enumerate}

\begin{quote}
We suggest that in a recently proposed framework for quantum gravity,
where Vassiliev invariants span the the space of states, the latter is
dramatically reduced if one has a non-vanishing cosmological constant.
This naturally suggests that the initial state of the universe should
have been one with \(\Lambda=0\).
\end{quote}

\begin{enumerate}
\def\labelenumi{\arabic{enumi})}
\setcounter{enumi}{10}
\tightlist
\item
  R. Gambini, O. Obregon and J. Pullin, ``Yang--Mills analogues of the
  Immirzi ambiguity'', available as
  \href{https://arxiv.org/abs/gr-qc/9801055}{\texttt{gr-qc/9801055}}.
\end{enumerate}

\begin{quote}
We draw parallels between the recently introduced `Immirzi ambiguity''
of the Ashtekar-like formulation of canonical quantum gravity and other
ambiguities that appear in Yang--Mills theories, like the \(\theta\)
ambiguity. We also discuss ambiguities in the Maxwell case, and
implication for the loop quantization of these theories.
\end{quote}

\begin{enumerate}
\def\labelenumi{\arabic{enumi})}
\setcounter{enumi}{11}
\tightlist
\item
  John Baez and Stephen Sawin, ``Diffeomorphism-invariant spin network
  states'', \emph{Jour. Funct. Analysis}, \textbf{158} (1998), 253--266.   
  Also available as
  \href{https://arxiv.org/abs/q-alg/9708005}{\texttt{q-alg/9708005}}.
\end{enumerate}

\begin{quote}
We extend the theory of diffeomorphism-invariant spin network states
from the real-analytic category to the smooth category. Suppose that
\(G\) is a compact connected semisimple Lie group and \(P\to M\) is a
smooth principal \(G\)-bundle. A `cylinder function' on the space of
smooth connections on \(P\) is a continuous complex function of the
holonomies along finitely many piecewise smoothly immersed curves in
\(M\). We construct diffeomorphism-invariant functionals on the space of
cylinder functions from `spin networks': graphs in \(M\) with edges
labeled by representations of \(G\) and vertices labeled by intertwining
operators. Using the `group averaging' technique of Ashtekar, Marolf,
Mourao and Thiemann, we equip the space spanned by these
`diffeomorphism-invariant spin network states' with a natural inner
product.
\end{quote}

Finally, here are two recent reviews of string theory and supersymmetry:

\begin{enumerate}
\def\labelenumi{\arabic{enumi})}
\setcounter{enumi}{12}
\item
  John H. Schwarz and Nathan Seiberg, ``String theory, supersymmetry,
  unification, and all that'', \emph{Reviews of Modern Physics} \textbf{71} (1999), S112.
  Also available as
  \href{https://arxiv.org/abs/hep-th/9803179}{\texttt{hep-th/9803179}}.
\item
  Keith R. Dienes and Christopher Kolda, ``Twenty open questions in
  supersymmetric particle physics'', available as
  \href{https://arxiv.org/abs/hep-ph/9712322}{\texttt{hep-ph/9712322}}.
\end{enumerate}

I'm afraid I'll slack off on my ``tour of homotopy theory'' this week. I
want to get to fun stuff like model categories and \(E_\infty\) spaces,
but it's turning out to be a fair amount of work to reach that goal!
That's what always happens with This Week's Finds: I start learning
about something and think ``oh boy, this stuff is great; I'll write it
up really carefully so that everyone can understand it,'' but then this
turns out to be so much work that by the time I'm halfway through I'm
off on some other kick.

\hypertarget{week121}{%
\section{May 15, 1998}\label{week121}}

This time I want to talk about higher-dimensional algebra and its
applications to topology. Marco Mackaay has just come out with a
fascinating paper that gives a construction of \(4\)-dimensional TQFTs
from certain ``monoidal \(2\)-categories''.

\begin{enumerate}
\def\labelenumi{\arabic{enumi})}
\tightlist
\item
  Marco Mackaay, ``Spherical \(2\)-categories and 4-manifold
  invariants'', available as
  \href{https://arxiv.org/abs/math.QA/9805030}{\texttt{math.QA/9805030}}.
\end{enumerate}

Beautifully, this construction is just a categorified version of Barrett
and Westbury's construction of \(3\)-dimensional topological quantum
field theories from ``monoidal categories''. Categorification --- the
process of replacing equations by isomorphisms --- is supposed to take
you up the ladder of dimensions. Here we are seeing it in action!

To prepare you understand Mackaay's paper, maybe I should explain the
idea of categorification. Since I recently wrote something about this, I
think I'll just paraphrase a bit of that. Some of this is already
familiar to long-time customers, so if you know it all already, just
skip it.

\begin{enumerate}
\def\labelenumi{\arabic{enumi})}
\setcounter{enumi}{1}
\tightlist
\item
  John Baez and James Dolan, ``Categorification'', in
  \emph{Higher Category Theory}, eds. Ezra Getzler and Mikhail Kapranov,
    \emph{Contemp.\ Math.} \textbf{230}, AMS, Providence, Rhode Island, 1998, 
     pp.\ 1--36.      Also available as
  \href{https://arxiv.org/abs/math.QA/9802029}{\texttt{math.QA/9802029}}.
\end{enumerate}

So, what's categorification? This tongue-twisting term, invented by
Louis Crane, refers to the process of finding category-theoretic analogs
of ideas phrased in the language of set theory, using the following
analogy between set theory and category theory:

\begin{longtable}{cc}
\toprule
\textbf{set theory} & \textbf{category theory} \tabularnewline
\midrule
\endhead
elements & objects \tabularnewline
equations between elements & isomorphisms between objects \tabularnewline
\midrule 
sets & categories \tabularnewline
functions between sets & functors between categories \tabularnewline
equations between functions & natural isomorphisms between functors \tabularnewline
\end{longtable}

Just as sets have elements, categories have objects. Just as there are
functions between sets, there are functors between categories.
Interestingly, the proper analog of an equation between elements is not
an equation between objects, but an isomorphism. More generally, the
analog of an equation between functions is a natural isomorphism between
functors.

For example, the category \(\mathsf{FinSet}\), whose objects are finite
sets and whose morphisms are functions, is a categorification of the set
\(\mathbb{N}\) of natural numbers. The disjoint union and Cartesian
product of finite sets correspond to the sum and product in
\(\mathbb{N}\), respectively. Note that while addition and
multiplication in \(\mathbb{N}\) satisfy various equational laws such as
commutativity, associativity and distributivity, disjoint union and
Cartesian product satisfy such laws \emph{only up to natural
isomorphism}. This is a good example of how equations between functions
get replaced by natural isomorphisms when we categorify.

If one studies categorification one soon discovers an amazing fact: many
deep-sounding results in mathematics are just categorifications of facts
we learned in high school! There is a good reason for this. All along,
we have been unwittingly ``decategorifying'' mathematics by pretending
that categories are just sets. We ``decategorify'' a category by
forgetting about the morphisms and pretending that isomorphic objects
are equal. We are left with a mere set: the set of isomorphism classes
of objects.

To understand this, the following parable may be useful. Long ago, when
shepherds wanted to see if two herds of sheep were isomorphic, they
would look for an explicit isomorphism. In other words, they would line
up both herds and try to match each sheep in one herd with a sheep in
the other. But one day, along came a shepherd who invented
decategorification. She realized one could take each herd and ``count''
it, setting up an isomorphism between it and some set of ``numbers'',
which were nonsense words like ``one, two, three,\ldots{}'' specially
designed for this purpose. By comparing the resulting numbers, she could
show that two herds were isomorphic without explicitly establishing an
isomorphism! In short, by decategorifying the category of finite sets,
the set of natural numbers was invented.

According to this parable, decategorification started out as a stroke of
mathematical genius. Only later did it become a matter of dumb habit,
which we are now struggling to overcome by means of categorification.
While the historical reality is far more complicated, categorification
really has led to tremendous progress in mathematics during the 20th
century. For example, Noether revolutionized algebraic topology by
emphasizing the importance of homology groups. Previous work had focused
on Betti numbers, which are just the dimensions of the rational homology
groups. As with taking the cardinality of a set, taking the dimension of
a vector space is a process of decategorification, since two vector
spaces are isomorphic if and only if they have the same dimension.
Noether noted that if we work with homology groups rather than Betti
numbers, we can solve more problems, because we obtain invariants not
only of spaces, but also of maps.

In modern lingo, the \(n\)th rational homology is a \emph{functor}
defined on the \emph{category} of topological spaces, while the \(n\)th
Betti number is a mere \emph{function}, defined on the \emph{set} of
isomorphism classes of topological spaces. Of course, this way of
stating Noether's insight is anachronistic, since it came before
category theory. Indeed, it was in Eilenberg and Mac Lane's subsequent
work on homology that category theory was born!

Decategorification is a straightforward process which typically destroys
information about the situation at hand. Categorification, being an
attempt to recover this lost information, is inevitably fraught with
difficulties. One reason is that when categorifying, one does not merely
replace equations by isomorphisms. One also demands that these
isomorphisms satisfy some new equations of their own, called ``coherence
laws''. Finding the right coherence laws for a given situation is
perhaps the trickiest aspect of categorification.

For example, a monoid is a set with a product satisfying the associative
law and a unit element satisfying the left and right unit laws. The
categorified version of a monoid is a ``monoidal category''. This is a
category \(\mathcal{C}\) with a product
\[\otimes\colon\mathcal{C}\times\mathcal{C}\to \mathcal{C}\] and unit
object \(1\). If we naively impose associativity and the left and right
unit laws as equational laws, we obtain the definition of a ``strict''
monoidal category. However, the philosophy of categorification suggests
instead that we impose them only up to natural isomorphism. Thus, as
part of the structure of a ``weak'' monoidal category, we specify a
natural isomorphism
\[a_{x,y,z}\colon (x \otimes y) \otimes z\to x \otimes (y \otimes z)\]
called the ``associator'', together with natural isomorphisms
\[l_x\colon 1 \otimes x\to x,\] \[r_x\colon x \otimes 1\to x.\] Using
the associator one can construct isomorphisms between any two
parenthesized versions of the tensor product of several objects.
However, we really want a \emph{unique} isomorphism. For example, there
are 5 ways to parenthesize the tensor product of 4 objects, which are
related by the associator as follows: \[
  \begin{tikzpicture}
    \node (mr) at (18:3) {$(x \otimes (y \otimes z)) \otimes w$};
    \node (t) at (90:3) {$((x \otimes y) \otimes z) \otimes w$};
    \node (ml) at (162:3) {$(x \otimes y) \otimes (z \otimes w)$};
    \node (bl) at (214:2.5) {$x \otimes (y \otimes (z \otimes w))$};
    \node (br) at (326:2.5) {$x \otimes ((y \otimes z) \otimes w)$};
    \draw[->] (t) to (ml);
    \draw[->] (ml) to (bl);
    \draw[->] (t) to (mr);
    \draw[->] (mr) to (br);
    \draw[->] (br) to (bl);
  \end{tikzpicture}
\] In the definition of a weak monoidal category we impose a coherence
law, called the ``pentagon identity'', saying that this diagram
commutes. Similarly, we impose a coherence law saying that the following
diagram built using \(a\), \(l\) and \(r\) commutes: \[
  \begin{tikzpicture}
    \node (mr) at (44:2) {$1\otimes(x\otimes 1)$};
    \node (t) at (136:2) {$(1\otimes x)\otimes1$};
    \node (ml) at (198:2) {$x\otimes1$};
    \node (bl) at (270:2) {$x$};
    \node (br) at (342:2) {$1\otimes x$};
    \draw[->] (t) to (ml);
    \draw[->] (ml) to (bl);
    \draw[->] (t) to (mr);
    \draw[->] (mr) to (br);
    \draw[->] (br) to (bl);
  \end{tikzpicture}
\] This definition raises an obvious question: how do we know we have
found all the right coherence laws? Indeed, what does ``right'' even
\emph{mean} in this context? Mac Lane's coherence theorem gives one
answer to this question: the above coherence laws imply that any two
isomorphisms built using \(a\), \(l\) and \(r\) and having the same
source and target must be equal.

Further work along these lines allow us to make more precise the sense
in which \(\mathbb{N}\) is a decategorification of \(\mathsf{FinSet}\).
For example, just as \(\mathbb{N}\) forms a monoid under either addition
or multiplication, \(\mathsf{FinSet}\) becomes a monoidal category under
either disjoint union or Cartesian product if we choose the isomorphisms
\(a\), \(l\), and \(r\) sensibly. In fact, just as \(\mathbb{N}\) is a
``rig'', satisfying all the ring axioms except those involving additive
inverses, \(\mathsf{FinSet}\) is what one might call a ``rig category''.
In other words, it satisfies the rig axioms up to natural isomorphisms
satisfying the coherence laws discovered by Kelly and Laplaza, who
proved a coherence theorem in this context.

Just as the decategorification of a monoidal category is a monoid, the
decategorification of any rig category is a rig. In particular,
decategorifying the rig category \(\mathsf{FinSet}\) gives the rig
\(\mathbb{N}\). This idea is especially important in combinatorics,
where the best proof of an identity involving natural numbers is often a
``bijective proof'': one that actually establishes an isomorphism
between finite sets.

While coherence laws can sometimes be justified retrospectively by
coherence theorems, certain puzzles point to the need for a deeper
understanding of the \emph{origin} of coherence laws. For example,
suppose we want to categorify the notion of ``commutative monoid''. The
strictest possible approach, where we take a strict monoidal category
and impose an equational law of the form \(x\otimes y = y\otimes x\), is
almost completely uninteresting. It is much better to start with a weak
monoidal category equipped with a natural isomorphism
\[B_{x,y}\colon x\otimes y\to y\otimes x\] called the ``braiding'' and
then impose coherence laws called ``hexagon identities'' saying that the
following two diagrams built from the braiding and the associator
commute: \[
  \begin{tikzcd}
    x\otimes(y\otimes z) \dar \rar
    & (y\otimes z)\otimes x
  \\(x\otimes y)\otimes z \dar
    & y\otimes (z\otimes x) \uar
  \\(y\otimes x)\otimes z \rar
    & y\otimes (x\otimes z) \uar
  \end{tikzcd}
\] \[
  \begin{tikzcd}
    (x\otimes y)\otimes z \dar \rar
    & z\otimes (x\otimes y)
  \\x\otimes (y\otimes z) \dar
    & (z\otimes x)\otimes y \uar
  \\y\otimes (z\otimes y) \rar
    & (x\otimes z)\otimes y \uar
  \end{tikzcd}
\] This gives the definition of a weak ``braided monoidal category''. If
we impose an additional coherence law saying that \(B_{x,y}\) is the
inverse of \(B_{y,x}\), we obtain the definition of a ``symmetric
monoidal category''. Both of these concepts are very important; which
one is ``right'' depends on the context. However, neither implies that
every pair of parallel morphisms built using the braiding are equal. A
good theory of coherence laws must naturally account for these facts.

The deepest insights into such puzzles have traditionally come from
topology. In homotopy theory it causes problems to work with spaces
equipped with algebraic structures satisfying equational laws, because
one cannot transport such structures along homotopy equivalences. It is
better to impose laws \emph{only up to homotopy}, with these homotopies
satisfying certain coherence laws, but again only up to homotopy, with
these higher homotopies satisfying their own higher coherence laws, and
so on. Coherence laws thus arise naturally in infinite sequences. For
example, Stasheff discovered the pentagon identity and a sequence of
higher coherence laws for associativity when studying the algebraic
structure possessed by a space that is homotopy equivalent to a loop
space. Similarly, the hexagon identities arise as part of a sequence of
coherence laws for spaces homotopy equivalent to double loop spaces,
while the extra coherence law for symmetric monoidal categories arises
as part of a sequence for spaces homotopy equivalent to triple loop
spaces. The higher coherence laws in these sequences turn out to be
crucial when we try to \emph{iterate} the process of categorification.

To \emph{iterate} the process of categorification, we need a concept of
``\(n\)-category'' --- roughly, an algebraic structure consisting of a
collection of objects (or ``0-morphisms''), morphisms between objects
(or ``1-morphisms''), \(2\)-morphisms between morphisms, and so on up to
\(n\)-morphisms. There are various ways of making this precise, and
right now there is a lot of work going on devoted to relating these
different approaches. But the basic thing to keep in mind is that the
concept of ``\((n+1)\)-category'' is a categorification of the concept
of ``\(n\)-category''. What were equational laws between \(n\)-morphisms
in an \(n\)-category are replaced by natural \((n+1)\)-isomorphisms,
which need to satisfy certain coherence laws of their own.

To get a feeling for how these coherence laws are related to homotopy
theory, it's good to think about certain special kinds of
\(n\)-category. If we have an \((n+k)\)-category that's trivial up to
but not including the k-morphism level, we can turn it into an
\(n\)-category by a simple reindexing trick: just think of its
\(j\)-morphisms as \((j-k)\)-morphisms! We call the \(n\)-categories we
get this way ``\(k\)-tuply monoidal \(n\)-categories''. Here is a little
chart of what they amount to for various low values of \(n\) and \(k\):

\vfill
\newpage

\begin{longtable}[]{@{}llll@{}}
\caption*{\(k\)-tuply monoidal \(n\)-categories}\tabularnewline
\toprule
\begin{minipage}[b]{0.26\columnwidth}\raggedright
\strut
\end{minipage} & \begin{minipage}[b]{0.21\columnwidth}\raggedright
\(n=0\)\strut
\end{minipage} & \begin{minipage}[b]{0.21\columnwidth}\raggedright
\(n=1\)\strut
\end{minipage} & \begin{minipage}[b]{0.21\columnwidth}\raggedright
\(n=2\)\strut
\end{minipage}\tabularnewline
\midrule
\endfirsthead
\toprule
\begin{minipage}[b]{0.26\columnwidth}\raggedright
\strut
\end{minipage} & \begin{minipage}[b]{0.21\columnwidth}\raggedright
\(n=0\)\strut
\end{minipage} & \begin{minipage}[b]{0.21\columnwidth}\raggedright
\(n=1\)\strut
\end{minipage} & \begin{minipage}[b]{0.21\columnwidth}\raggedright
\(n=2\)\strut
\end{minipage}\tabularnewline
\midrule
\endhead
\begin{minipage}[t]{0.26\columnwidth}\raggedright
\(k=0\)\strut
\end{minipage} & \begin{minipage}[t]{0.21\columnwidth}\raggedright
sets\strut
\end{minipage} & \begin{minipage}[t]{0.21\columnwidth}\raggedright
categories\strut
\end{minipage} & \begin{minipage}[t]{0.21\columnwidth}\raggedright
\(2\)-categories\strut
\end{minipage}\tabularnewline
\begin{minipage}[t]{0.26\columnwidth}\raggedright
\strut
\end{minipage} & \begin{minipage}[t]{0.21\columnwidth}\raggedright
\strut
\end{minipage} & \begin{minipage}[t]{0.21\columnwidth}\raggedright
\strut
\end{minipage} & \begin{minipage}[t]{0.21\columnwidth}\raggedright
\strut
\end{minipage}\tabularnewline
\begin{minipage}[t]{0.26\columnwidth}\raggedright
\(k=1\)\strut
\end{minipage} & \begin{minipage}[t]{0.21\columnwidth}\raggedright
monoids\strut
\end{minipage} & \begin{minipage}[t]{0.21\columnwidth}\raggedright
monoidal categories\strut
\end{minipage} & \begin{minipage}[t]{0.21\columnwidth}\raggedright
monoidal \(2\)-categories\strut
\end{minipage}\tabularnewline
\begin{minipage}[t]{0.26\columnwidth}\raggedright
\strut
\end{minipage} & \begin{minipage}[t]{0.21\columnwidth}\raggedright
\strut
\end{minipage} & \begin{minipage}[t]{0.21\columnwidth}\raggedright
\strut
\end{minipage} & \begin{minipage}[t]{0.21\columnwidth}\raggedright
\strut
\end{minipage}\tabularnewline
\begin{minipage}[t]{0.26\columnwidth}\raggedright
\(k=2\)\strut
\end{minipage} & \begin{minipage}[t]{0.21\columnwidth}\raggedright
commutative monoids\strut
\end{minipage} & \begin{minipage}[t]{0.21\columnwidth}\raggedright
braided monoidal categories\strut
\end{minipage} & \begin{minipage}[t]{0.21\columnwidth}\raggedright
braided monoidal \(2\)-categories\strut
\end{minipage}\tabularnewline
\begin{minipage}[t]{0.26\columnwidth}\raggedright
\strut
\end{minipage} & \begin{minipage}[t]{0.21\columnwidth}\raggedright
\strut
\end{minipage} & \begin{minipage}[t]{0.21\columnwidth}\raggedright
\strut
\end{minipage} & \begin{minipage}[t]{0.21\columnwidth}\raggedright
\strut
\end{minipage}\tabularnewline
\begin{minipage}[t]{0.26\columnwidth}\raggedright
\(k=3\)\strut
\end{minipage} & \begin{minipage}[t]{0.21\columnwidth}\raggedright
`` "\strut
\end{minipage} & \begin{minipage}[t]{0.21\columnwidth}\raggedright
symmetric monoidal categories\strut
\end{minipage} & \begin{minipage}[t]{0.21\columnwidth}\raggedright
weakly involutory monoidal \(2\)-categories\strut
\end{minipage}\tabularnewline
\begin{minipage}[t]{0.26\columnwidth}\raggedright
\strut
\end{minipage} & \begin{minipage}[t]{0.21\columnwidth}\raggedright
\strut
\end{minipage} & \begin{minipage}[t]{0.21\columnwidth}\raggedright
\strut
\end{minipage} & \begin{minipage}[t]{0.21\columnwidth}\raggedright
\strut
\end{minipage}\tabularnewline
\begin{minipage}[t]{0.26\columnwidth}\raggedright
\(k=4\)\strut
\end{minipage} & \begin{minipage}[t]{0.21\columnwidth}\raggedright
`` "\strut
\end{minipage} & \begin{minipage}[t]{0.21\columnwidth}\raggedright
`` "\strut
\end{minipage} & \begin{minipage}[t]{0.21\columnwidth}\raggedright
strongly involutory monoidal \(2\)-categories\strut
\end{minipage}\tabularnewline
\begin{minipage}[t]{0.26\columnwidth}\raggedright
\strut
\end{minipage} & \begin{minipage}[t]{0.21\columnwidth}\raggedright
\strut
\end{minipage} & \begin{minipage}[t]{0.21\columnwidth}\raggedright
\strut
\end{minipage} & \begin{minipage}[t]{0.21\columnwidth}\raggedright
\strut
\end{minipage}\tabularnewline
\begin{minipage}[t]{0.26\columnwidth}\raggedright
\(k=5\)\strut
\end{minipage} & \begin{minipage}[t]{0.21\columnwidth}\raggedright
`` "\strut
\end{minipage} & \begin{minipage}[t]{0.21\columnwidth}\raggedright
`` "\strut
\end{minipage} & \begin{minipage}[t]{0.21\columnwidth}\raggedright
`` "\strut
\end{minipage}\tabularnewline
\bottomrule
\end{longtable}

One reason James Dolan and I got so interested in this chart is the
``tangle hypothesis''. Roughly speaking, this says that
\(n\)-dimensional surfaces embedded in \((n+k)\)-dimensional space can
be described purely algebraically using the a certain special
``\(k\)-tuply monoidal \(n\)-category with duals''. If true, this
reduces lots of differential topology to pure algebra! It also helps you
understand the parameters \(n\) and \(k\): you should think of \(n\) as
``dimension'' and \(k\) as ``codimension''.

For example, take \(n = 1\) and \(k = 2\). Knots, links and tangles in
3-dimensional space can be described algebraically using a certain
``braided monoidal categories with duals''. This was the first
interesting piece of evidence for the tangle hypothesis. It has spawned
a whole branch of math called ``quantum topology'', which people are
trying to generalize to higher dimensions.

More recently, Laurel Langford tackled the case \(n = 2\), \(k = 2\).
She proved that \(2\)-dimensional knotted surfaces in \(4\)-dimensional
space can be described algebraically using a certain ``braided monoidal
\(2\)-category with duals''. These so-called ``2-tangles'' are
particularly interesting to me because of their relation to spin foam
models of quantum gravity, which are also all about surfaces in 4-space.
For references, see \protect\hyperlink{week103}{``Week 103''}. But if
you want to learn about more about this, you couldn't do better than to
start with:

\begin{enumerate}
\def\labelenumi{\arabic{enumi})}
\setcounter{enumi}{2}
\tightlist
\item
  J. S. Carter and M. Saito, \emph{Knotted Surfaces and Their Diagrams},
  American Mathematical Society, Providence, 1998.
\end{enumerate}

This is a magnificently illustrated book which will really get you able
to \emph{see} \(2\)-dimensional surfaces knotted in 4d space. At the end
it sketches the statement of Langford's result.

Another interesting thing about the above chart is that \(k\)-tuply
monoidal \(n\)-categories keep getting ``more commutative'' as \(k\)
increases, until one reaches \(k = n+2\), at which point things
stabilize. There is a lot of evidence suggesting that this
``stabilization hypothesis'' is true for all \(n\). Assuming it's true,
it makes sense to call a \(k\)-tuply monoidal \(n\)-category with
\(k\geqslant n+2\) a ``stable \(n\)-category''.

Now, where does homotopy theory come in? Well, here you need to look at
\(n\)-categories where all the \(j\)-morphisms are invertible for all
\(j\). These are called ``\(n\)-groupoids''. Using these, one can
develop a translation dictionary between \(n\)-category theory and
homotopy theory, which looks like this:

\begin{longtable}{c|c}
\toprule
\(\omega\)-groupoids & homotopy types \\
\(n\)-groupoids & homotopy \(n\)-types \\
\(k\)-tuply groupal \(\omega\)-groupoids & homotopy types of
\(k\)-fold loop spaces \\
\(k\)-tuply groupal \(n\)-groupoids & homotopy \(n\)-types of
\(k\)-fold loop spaces \\
\(k\)-tuply monoidal \(\omega\)-groupoids & homotopy types of
\(E_k\) spaces \\
\(k\)-tuply monoidal \(n\)-groupoids & homotopy \(n\)-types of
\(E_k\) spaces\\
stable \(\omega\)-groupoids & homotopy types of infinite loop
spaces \\
stable \(n\)-groupoids & homotopy \(n\)-types of infinite loop
spaces \\
\(\mathbb{Z}\)-groupoids & homotopy types of spectra \\
\bottomrule
\end{longtable}

The entries on the left-hand side are very natural from an algebraic
viewpoint; the entries on the right-hand side are things topologists
already study. We explain what all these terms mean in the paper, but
maybe I should say something about the first two rows, which are the
most basic in a way. A homotopy type is roughly a topological space ``up
to homotopy equivalence'', and an \(\omega\)-groupoid is a kind of
limiting case of an \(n\)-groupoid as \(n\) goes to infinity. If
infinity is too scary, you can work with homotopy \(n\)-types, which are
basically homotopy types with no interesting topology above dimension
\(n\). These should correspond to \(n\)-groupoids.

Using these basic correspondences we can then relate various special
kinds of homotopy types to various special kinds of
\(\omega\)-groupoids, giving the rest of the rows of the chart. Homotopy
theorists know a lot about the right-hand column, so we can use this to
get a lot of information about the left-hand column. In particular, we
can work out the coherence laws for \(n\)-groupoids, and --- this is the
best part, but the least understood --- we can then \emph{guess} a lot
of stuff about the coherence laws for \emph{general} \(n\)-categories.
In short, we are using homotopy theory to get our foot in the door of
\(n\)-category theory.

I should emphasize, though, that this translation dictionary is
partially conjectural. It gets pretty technical to say what exactly is
and is not known, especially since there's pretty rapid progress going
on. Even in the last few months there have been some interesting
developments. For example, Breen has come out with a paper relating
\(k\)-tuply monoidal \(n\)-categories to Postnikov towers and various
far-out kinds of homological algebra:

\begin{enumerate}
\def\labelenumi{\arabic{enumi})}
\setcounter{enumi}{3}
\tightlist
\item
  Lawrence Breen, ``Braided \(n\)-categories and
  \(\Sigma\)-structures'', Prepublications Matematiques de l'Universite
  Paris 13, 98-06, January 1998, to appear in the \emph{Proceedings of
  the Workshop on Higher Category Theory and Mathematical Physics at
  Northwestern University, Evanston, Illinois, March 1997}, eds.~Ezra
  Getzler and Mikhail Kapranov.
\end{enumerate}

Also, the following folks have also developed a notion of ``iterated
monoidal category'' whose nerve gives the homotopy type of a \(k\)-fold
loop space, just as the nerve of a category gives an arbitrary homotopy
type:

\begin{enumerate}
\def\labelenumi{\arabic{enumi})}
\setcounter{enumi}{4}
\tightlist
\item
  C. Balteanu, Z. Fiedorowicz, R. Schwaenzl, and R. Vogt, ``Iterated
  monoidal categories'', \emph{Adv. Math.} \textbf{176} (2003), 277--349. 
   Also available at
  \href{http://arxiv.org/abs/math.AT/9808082}{\texttt{math.AT/9808082}}.
\end{enumerate}

Anyway, in addition to explaining the relationship between
\(n\)-category theory and homotopy theory, Dolan's and my paper
discusses iterated categorifications of the very simplest algebraic
structures: the natural numbers and the integers. The natural numbers
are the free monoid on one generator; the integers are the free group on
one generator. We believe this is just the start of this chart listing various
algebraic structures and the free such structures on one generator:

\begin{longtable}{c|c}
\toprule
sets & the one-element set \\
monoids & the natural numbers \\
groups & the integers \\
\(k\)-tuply monoidal \(n\)-categories & the braid
\(n\)-groupoid in codimension \(k\) \\
\(k\)-tuply monoidal \(\omega\)-categories & the braid
\(\omega\)-groupoid in codimension \(k\) \\
stable \(n\)-categories & the braid \(n\)-groupoid in infinite
codimension \\
stable \(\omega\)-categories & the braid \(\omega\)-groupoid in
infinite codimension \\
\(k\)-tuply monoidal \(n\)-categories with duals & the
\(n\)-category of framed \(n\)-tangles in \(n+k\) dimensions
\\
stable \(n\)-categories with duals & the framed cobordism
\(n\)-category \\
\(k\)-tuply groupal \(n\)-groupoids & the homotopy \(n\)-tpye
of the \(k\)th loop space of \(S^k\) \\
\(k\)-tuply groupal \(\omega\)-groupoids & the homotopy type of
the \(k\)th loop space of \(S^k\) \\
stable \(\omega\)-groupoids & the homotopy type of the infinite
loop space \(S^\infty\) \\
\(\mathbb{Z}\)-groupoids & the sphere spectrum \\
\bottomrule
\end{longtable}

You may or may not know the guys on the right-hand side, but some of
them are very interesting and complicated, so it's really exciting that
they are all in some sense categorified and/or stabilized versions of
the integers and natural numbers.

Whew! There is more to say, but I'll just mention a few related papers
and then quit. If you're interested in \(n\)-categories you could also
check out ``the tale of \(n\)-categories'', starting in
\protect\hyperlink{week73}{``Week 73''}.

\begin{enumerate}
\def\labelenumi{\arabic{enumi})}
\setcounter{enumi}{5}
\tightlist
\item
   Martin Neuchl, \emph{Representation Theory of Hopf Categories}, 
  Ph.D. thesis, Department of Mathematics, University of Munich,
  1997. Available at \href{http://math.ucr.edu/home/baez/neuchl.ps}{\texttt{http://math.ucr.}} \href{http://math.ucr.edu/home/baez/neuchl.ps}{\texttt{edu/home/baez/neuchl.ps}} and  \href{http://math.ucr.edu/home/baez/neuchl.pdf}{\texttt{http://math.ucr.edu/home/baez/neuchl.pdf}}
\end{enumerate}

Just as the category of representations of a Hopf algebra gives a nice
monoidal category, the \(2\)-category of representations of a Hopf
category gives a nice monoidal \(2\)-category! Categorification strikes
again --- and this is perhaps our best hope for getting our hands on
the data needed to stick into Mackaay's machine and get concrete
examples of a 4d topological quantum field theories!

\begin{enumerate}
\def\labelenumi{\arabic{enumi})}
\setcounter{enumi}{6}
\tightlist
\item
  Jim Stasheff, ``Grafting Boardman's cherry trees to quantum field
  theory'', available as
  \href{https://arxiv.org/abs/math.AT/9803156}{\texttt{math.AT/9803156}}.
\end{enumerate}

Starting with Boardman and Vogt's work, and shortly thereafter that of
May, operads have become really important in homotopy theory, string
theory, and now \(n\)-category theory; this review article sketches some
of the connections.

\begin{enumerate}
\def\labelenumi{\arabic{enumi})}
\setcounter{enumi}{7}
\item
  Masoud Khalkhali, ``On cyclic homology of \(A_\infty\) algebras'',
   available as
  \href{https://arxiv.org/abs/math.QA/9805051}{\texttt{math.QA/9805051}}.

  Masoud Khalkhali, Homology of \(L_\infty\) algebras and cyclic
  homology, available as
  \href{https://arxiv.org/abs/math.QA/9805052}{\texttt{math.QA/9805052}}.
\end{enumerate}

An \(A_\infty\) algebra is an algebra that is associative \emph{up to an
associator} which satisfies the pentagon identity \emph{up to a
pentagonator} which satisfies it's own coherence law up to something, ad
infinitum. The concept goes back to Stasheff's work on \(A_\infty\)
spaces --- spaces with a homotopy equivalence to a space equipped with
an associative product. (These are the same thing as what I called
\(E_1\) spaces in the translation dictionary between \(n\)-groupoid
theory and homotopy theory.) But here it's been transported from Top
over to Vect. Similarly, an \(L_\infty\) algebra is a Lie algebra ``up
to an infinity of higher coherence laws''. Loday-Quillen and Tsygan
showed that that the Lie algebra homology of the algebra of stable
matrices over an associative algebra is isomorphic, as a Hopf algebra,
to the exterior algebra of the cyclic homology of the algebra. In the
second paper above, Khalkali gets the tools set up to extend this result
to the category of \(L_\infty\) algebras.

\hypertarget{week122}{%
\section{June 24, 1998}\label{week122}}

In summertime, academics leave the roost and fly hither and thither,
seeking conferences and conversations in far-flung corners of the world.
At the end of May, everyone started leaving the Center for Gravitational
Physics and Geometry: Lee Smolin for the Santa Fe Institute, Abhay
Ashtekar for Uruguay and Argentina, Kirill Krasnov for his native
Ukraine, and so on. It got so quiet that I could actually get some work
done, were it not for the fact that I, too, flew the coop: first for
Chicago, then Portugal, and then to one of the most isolated,
technologically backwards areas on earth: my parents' house. Connected
to cyberspace by only the thinnest of threads, writing new issues of
This Week's Finds became almost impossible\ldots.

I did, however, read some newsgroups, and by this means Jim Carr
informed me that an article on spin foam models of quantum gravity had
appeared in Science News. I can't resist mentioning it, since it quotes
me:

\begin{enumerate}
\def\labelenumi{\arabic{enumi})}
\tightlist
\item
  Ivars Peterson, ``Loops of gravity: calculating a foamy quantum
  space-time'', \emph{Science News}, June 13, 1998, Vol. \textbf{153},
  No.~24, 376--377.
\end{enumerate}

It gives a little history of loop quantum gravity, spin networks, and
the new burst of interest in spin foams. Nothing very technical --- but
good if you're just getting started. If you want something more
detailed, but still user-friendly, try Rovelli's new paper:

\begin{enumerate}
\def\labelenumi{\arabic{enumi})}
\setcounter{enumi}{1}
\tightlist
\item
  Carlo Rovelli and Peush Upadhya, ``Loop quantum gravity and quanta of
  space: a primer'', available as
  \href{https://arxiv.org/abs/gr-qc/9806079}{\texttt{gr-qc/9806079}}.
\end{enumerate}

I haven't read it yet, since I'm still in a rather low-tech portion of
the globe, but it gives simplified derivations of some of the basic
results of loop quantum gravity, like the formula for the eigenvalues of
the area operator. As explained in \protect\hyperlink{week110}{``Week
110''}, one of the main predictions of loop quantum gravity is that
geometrical observables such as the area of any surface take on a
discrete spectrum of values, much like the energy levels of a hydrogen
atom. At first the calculation of the eigenvalues of the area operator
seemed rather complicated, but by now it's well-understood, so Rovelli
and Upadhya are able to give a simpler treatment.

While I'm talking about the area operator, I should mention another
paper by Rovelli, in which he shows that its spectrum is not affected by
the presence of matter (or more precisely, fermions):

\begin{enumerate}
\def\labelenumi{\arabic{enumi})}
\setcounter{enumi}{2}
\tightlist
\item
  Carlo Rovelli and Merced Montesinos, ``The fermionic contribution to
  the spectrum of the area operator in nonperturbative quantum
  gravity'', available as
  \href{https://arxiv.org/abs/gr-qc/9806120}{\texttt{gr-qc/9806120}}.
\end{enumerate}

This is especially interesting because it fits in with other pieces of
evidence that fermions could simply be the ends of wormholes --- an old
idea of John Wheeler (see \protect\hyperlink{week109}{``Week 109''}).

I should also mention some other good review articles that have turned
up recently. Rovelli has written a survey comparing string theory, the
loop representation, and other approaches to quantum gravity, which is
very good because it points out the flaws in all these approaches, which
their proponents are usually all too willing to keep quiet about:

\begin{enumerate}
\def\labelenumi{\arabic{enumi})}
\setcounter{enumi}{3}
\tightlist
\item
  Carlo Rovelli, ``Strings, loops and others: a critical survey of the
  present approaches to quantum gravity''. Plenary lecture on quantum
  gravity at the \emph{GR15 conference, Pune, India}, available
  as \href{https://arxiv.org/abs/gr-qc/9803024}{\texttt{gr-qc/9803024}}.
\end{enumerate}

Also, Loll has written a review of approaches to quantum gravity that
assume spacetime is discrete. It does \emph{not} discuss the spin foam
approach, which is too new; instead it mainly talks about lattice
quantum gravity, the Regge calculus, and the dynamical triangulations
approach. In lattice quantum gravity you treat spacetime as a fixed
lattice, usually a hypercubical one, and work with discrete versions of
the usual fields appearing in general relativity. In the Regge calculus
you triangulate your \(4\)-dimensional spacetime --- i.e., chop it into
a bunch of \(4\)-dimensional simplices --- and use the lengths of the
edges of these simplices as your basic variables. (For more details see
\protect\hyperlink{week120}{``Week 120''}.) In the dynamical
triangulations approach you also triangulate spacetime, but not in a
fixed way --- you consider all possible triangulations. However, you
assume all the edges of all the simplices have the same length --- the
Planck length, say. Thus all the information about the geometry of
spacetime is in the triangulation itself --- hence the name ``dynamical
triangulations''. Everything becomes purely combinatorial - there are no
real numbers in our description of spacetime geometry anymore. This
makes the dynamical triangulations approach great for computer
simulations. Computer simulations of quantum gravity! Loll reports on
the results of a lot of these:

\begin{enumerate}
\def\labelenumi{\arabic{enumi})}
\setcounter{enumi}{4}
\tightlist
\item
  Renate Loll, ``Discrete approaches to quantum gravity in four
  dimensions'', available as
  \href{https://arxiv.org/abs/gr-qc/9805049}{\texttt{gr-qc/9805049}}.
\end{enumerate}

Here are some other good places to learn about the dynamical
triangulations approach to quantum gravity:

\begin{enumerate}
\def\labelenumi{\arabic{enumi})}
\setcounter{enumi}{6}
\item
  J. Ambj\o rn, ``Quantum gravity represented as dynamical
  triangulations'', \emph{Class. Quant. Grav.} \textbf{12} (1995)
  2079--2134.
\item
  J. Ambj\o rn, M. Carfora, and A. Marzuoli, \emph{The Geometry of
  Dynamical Triangulations}, Springer, Berlin, 1998. Also
  available as
  \href{https://arxiv.org/abs/hep-th/9612069}{\texttt{hep-th/9612069}}.
\end{enumerate}

I can't resist pointing out an amusing relationship between dynamical
triangulations and mathematical logic, which Ambj\o rn mentions in his
review article. In computer simulations using the dynamical
triangulations approach, one wants to compute the average of certain
quantities over all triangulations of a fixed compact manifold --- e.g.,
the \(4\)-dimensional sphere, \(S^4\). The typical way to do this is to
start with a particular triangulation and then keep changing it using
various operations --- ``Pachner moves'' --- that are guaranteed to
eventually take you from any triangulation of a compact
\(4\)-dimensional manifold to any other.

Now here's where the mathematical logic comes in. Markov's theorem says
there is no algorithm that can decide whether or not two triangulations
are triangulations of the same compact \(4\)-dimensional manifold.
(Technically, by ``the same'' I mean ``piecewise linearly
homeomorphic'', but don't worry about that!) If they \emph{are}
triangulations of the same manifold, blundering about using the Pachner
moves will eventually get you from one to the other, but if they are
\emph{not}, you may never know for sure.

On the other hand, \(S^4\) may be special. It's an open question whether
or not \(S^4\) is ``algorithmically detectable''. In other words, it's
an open question whether or not there's an algorithm that can decide
whether or not a triangulation is a triangulation of the
\(4\)-dimensional sphere.

Now, suppose \(S^4\) is \emph{not} algorithmically detectable. Then the
maximum number of Pachner moves it takes to get between two
triangulations of the 4-sphere must grow really fast: faster than any
computable function! After all, if it didn't, we could use this upper
bound to know when to give up when using Pachner moves to try to reduce
our triangulation to a known triangulation of \(S^4\). So there must be
``bottlenecks'' that make it hard to efficiently explore the set of all
triangulations of \(S^4\) using Pachner moves. For example, there must
be pairs of triangulations such that getting from one to other via
Pachner moves requires going through triangulations with a \emph{lot}
more \(4\)-simplices.

However, computer simulations using triangulations with up to 65,536
4-simplices have not yet detected such ``bottlenecks''. What's going on?
Well, maybe S\^{}4 actually \emph{is} algorithmically detectable. Or
perhaps it's not, but the bottlenecks only occur for triangulations that
have more than 65,536 \(4\)-simplices to begin with. Interestingly, one
dimension up, it's known that the \(5\)-dimensional sphere is \emph{not}
algorithmically detectable, so in this case bottlenecks \emph{must}
exist --- but computer simulations still haven't seen them.

I should emphasize that in addition to this funny computability stuff,
there is also a whole lot of interesting \emph{physics} coming out of
the dynamical triangulations approach to quantum gravity. Unfortunately
I don't have the energy to explain this now --- so read those review
articles, and check out that nice book by Ambj\o rn, Carfora and Marzuoli!

On another front\ldots{} Ambj\o rn and Loll, who are both hanging out at
the AEI these days, have recently teamed up to study causality in a
lattice model of \(2\)-dimensional Lorentzian quantum gravity:

\begin{enumerate}
\def\labelenumi{\arabic{enumi})}
\setcounter{enumi}{8}
\tightlist
\item
  J. Ambj\o rn and R. Loll, ``Non-perturbative Lorentzian quantum gravity,
  causality and topology change'', available as
  \href{https://arxiv.org/abs/hep-th/9805108}{\texttt{hep-th/9805108}}.
\end{enumerate}

I'll just quote the abstract:

\begin{quote}
We formulate a non-perturbative lattice model of two-dimensional
Lorentzian quantum gravity by performing the path integral over
geometries with a causal structure. The model can be solved exactly at
the discretized level. Its continuum limit coincides with the theory
obtained by quantizing 2d continuum gravity in proper-time gauge, but it
disagrees with 2d gravity defined via matrix models or Liouville theory.
By allowing topology change of the compact spatial slices (i.e.~baby
universe creation), one obtains agreement with the matrix models and
Liouville theory.
\end{quote}

And now for something completely different\ldots{}

I've been hearing rumbles off in the distance about some interesting
work by Kreimer relating renormalization, Feynman diagrams, and Hopf
algebras. A friendly student of Kreimer named Mathias Mertens handed me
a couple of the basic papers when I was in Portugal:

\begin{enumerate}
\def\labelenumi{\arabic{enumi})}
\setcounter{enumi}{9}
\item
  Dirk Kreimer, ``Renormalization and knot theory'', \emph{Journal of
  Knot Theory and its Ramifications}, \textbf{6} (1997), 479--581.
  Also available as
  \href{https://arxiv.org/abs/q-alg/9607022}{\texttt{q-alg/9607022}}.

  Dirk Kreimer, ``On the Hopf algebra structure of perturbative quantum
  field theories'', available as
  \href{https://arxiv.org/abs/q-alg/9707029}{\texttt{q-alg/9707029}}.
\end{enumerate}

I'm looking through them but I don't really understand them yet. The
basic idea seems to be something like this. In quantum field theory you
compute the probability for some reaction among particles by doing
integrals which correspond in a certain way to pictures called Feynman
diagrams. Often these integrals give infinite answers, which forces you
to do a trick called renormalization to cancel the infinities and get
finite answers. Part of why this trick works is that while your
integrals diverge, they usually diverge at a well-defined rate. For
example, you might get something asymptotic to a constant times
\(1/d^k\), where \(d\) is the spatial cutoff you put in to get a finite
answer. And the constant you get here can be explicitly computed. For
example, it often involves numbers like \(\zeta(n)\), where \(\zeta\) is
the Riemann zeta function, much beloved by number theorists:
\[\zeta(n) = \frac{1}{1^n} + \frac{1}{2^n} + \frac{1}{3^n} + \ldots\]
Kreimer noticed that if you take the Feynman diagram and do some tricks
to turn it into a drawing of a knot or link, the constant you get is
related in interesting ways to the topology of this knot or link! More
complicated knots or links give fancier constants, and there are all
sorts of suggestive patterns. He worked out a bunch of examples in the
first paper cited above, and since then people have worked out lots
more, which you can find in the references.

Apparently the secret underlying reason for these patterns comes from
the combinatorics of renormalization, which Kreimer was able to
summarize in a certain algebraic structure called a Hopf algebra. Hopf
algebras are important in both combinatorics and physics, so perhaps
this shouldn't be surprising. But there is still a lot of mysterious
stuff going on, at least as far as I can tell.

What's really intriguing about all this is \emph{which} quantum field
theories Kreimer was studying when he discovered this stuff: \emph{not}
topological quantum field theories like Chern--Simons theory, which
already have well-understood relationship to knot theory, but instead,
field theories that ordinary particle physicists have been thinking
about for decades, like quantum electrodynamics, \(\varphi^4\) theory in
4 dimensions, and \(\varphi^3\) theory in 6 dimensions --- field
theories where renormalization is a deadly serious business, thanks to
nasty problems like ``overlapping divergences''.

The idea that knot theory is relevant to \emph{these} field theories is
exciting but also somewhat puzzling, since they don't live in
3-dimensional spacetime the way Chern--Simons theory does. People
familiar with Chern--Simons theory have already been seeing fascinating
patterns relating knot theory, quantum field theory and number theory.
Is this new stuff related? Or is it something completely different?
Kreimer seems to think it's related.

According to Kirill Krasnov, the famous mathematician Alain Connes is
going around telling people to learn about this stuff. Apparently Connes
is now writing a paper on it with Kreimer, and it was Connes who got the
authors of this paper interested in the subject:

\begin{enumerate}
\def\labelenumi{\arabic{enumi})}
\setcounter{enumi}{10}
\tightlist
\item
  Thomas Krajewski and Raimar Wulkenhaar, ``On Kreimer's Hopf algebra
  structure of Feynman graphs'', available as
  \href{https://arxiv.org/abs/hep-th/9805098}{\texttt{hep-th/9805098}}.
\end{enumerate}

Since I haven't plunged in yet, I'll just quote the abstract:

\begin{quote}
We reinvestigate Kreimer's Hopf algebra structure of perturbative
quantum field theories. In Kreimer's original work, overlapping
divergences were first disentangled into a linear combination of
disjoint and nested ones using the Schwinger-Dyson equation. The linear
combination then was tackled by the Hopf algebra operations. We present
a formulation where the coproduct itself produces the linear
combination, without reference to external input.
\end{quote}

With any luck, mathematicians will study this stuff and finally
understand renormalization!

\hypertarget{week123}{%
\section{September 19, 1998}\label{week123}}

\begin{quote}
It all started out as a joke. Argument for argument's sake. Alison and
her infuriating heresies.

``A mathematical theorem,'' she'd proclaimed, "only becomes true when a
physical system tests it out: when the system's behaviour depends in
some way on the theorem being \emph{true} or \emph{false}.

It was June 1994. We were sitting in a small paved courtyard, having
just emerged from the final lecture in a one-semester course on the
philosophy of mathematics --- a bit of light relief from the hard grind
of the real stuff. We had fifteen minutes to to kill before meeting some
friends for lunch. It was a social conversation --- verging on mild
flirtation --- nothing more. Maybe there were demented academics,
lurking in dark crypts somewhere, who held views on the nature of
mathematical truth which they were willing to die for. But were were
twenty years old, and we \emph{knew} it was all angels on the head of a
pin.

I said, ``Physical systems don't create mathematics. Nothing
\emph{creates} mathematics --- it's timeless. All of number theory would
still be exactly the same, even if the universe contained nothing but a
single electron.''

Alison snorted. ``Yes, because even \emph{one electron}, plus a
space-time to put it in, needs all of quantum mechanics and all of
general relativity --- and all the mathematical infrastructure they
entail. One particle floating in a quantum vacuum needs half the major
results of group theory, functional analysis, differential geometry
---''

``OK, OK! I get the point. But if that's the case\ldots{} the events in
the first picosecond after the Big Bang would have `constructed' every
last mathematical truth required by \emph{any} physical system, all the
way to the Big Cruch. Once you've got the mathematics which underpins
the Theory of Everything\ldots{} that's it, that's all you ever need.
End of story.''

``But it's not. To \emph{apply} the Theory of Everything to a particular
system, you still need all the mathematics for dealing with \emph{that
system} --- which could include results far beyond the mathematics the
TOE itself requires. I mean, fifteen billion years after the Big Bang,
someone can still come along and prove, say\ldots{} Fermat's Last
Theorem.'' Andrew Wiles at Princeton had recently announced a proof of
the famous conjecture, although his work was still being scrutinised by
his colleagues, and the final verdict wasn't yet in. ``Physics never
needed \emph{that} before.''

I protested, ``What do you mean, `before'? Fermat's Last Theorem never
has --- and never will --- have anything to do with any branch of
physics.''

Alison smiled sneakily. ``No \emph{branch}, no. But only because the
class of physical systems whose behaviour depend on it is so ludicrously
specific: the brains of mathematicians who are trying to validate the
Wiles proof.''

``Think about it. Once you start trying to prove a theorem, then even if
the mathematics is so `pure' that it has no relevance to any other
object in the universe\ldots{} you've just made it relevant to
\emph{yourself}. You have to choose \emph{some} physical process to test
the theorem --- whether you use a computer, or a pen and paper\ldots{}
or just close your eyes and shuffle \emph{neurotransmitters}. There's no
such thing as a proof which doesn't rely on physical events, and whether
they're inside or outside your skull doesn't make them any less real.''
\end{quote}

And this is just the beginning\ldots{} the beginning of Greg Egan's tale
of an inconsistency in the axioms of arithmetic --- a ``topological
defect'' left over in the fabric of mathematics, much like the cosmic
strings or monopoles hypothesized by certain physicists thinking about
the early universe --- and the mathematicians who discover it and
struggle to prevent a large corporation from exploiting it for their own
nefarious purposes. This is the title story of his new collection,
``Luminous''.

I should also mention his earlier collection of stories, named after a
sophisticated class of mind-altering nanotechnologies, the
``axiomatics'', that affect specific beliefs of anyone who uses them:

\begin{enumerate}
\def\labelenumi{\arabic{enumi})}
\item
  Greg Egan, \emph{Axiomatic}, Orion Books, 1995.

  Greg Egan, \emph{Luminous}, Orion Books, 1998.
\end{enumerate}

Some of the stories in these volumes concern math and physics, such as
``The Planck Dive'', about some far-future explorers who send copies of
themselves into a black hole to study quantum gravity firsthand. One
nice thing about this story, from a pedant's perspective, is that Egan
actually works out a plausible scenario for meeting the technical
challenges involved --- with the help of a little 23rd-century
technology. Another nice thing is the further exploration of a world in
which everyone has long been uploaded to virtual ``scapes'' and can
easily modify and copy themselves --- a world familiar to readers of his
novel ``Diaspora'' (see \protect\hyperlink{week115}{``Week 115''}). But
what I really like is that it's not just a hard-science extravaganza;
it's a meditation on mortality. You can never really know what it's like
to cross an event horizon unless you do it\ldots.

Other stories focus on biotechnology and philosophical problems of
identity. The latter sort will especially appeal to everyone who liked
this book:

\begin{enumerate}
\def\labelenumi{\arabic{enumi})}
\setcounter{enumi}{1}
\tightlist
\item
  Daniel C. Dennett and Douglas R. Hofstadter, \emph{The Mind's I:
  Fantasies and Reflections on Self and Soul}, Bantam Books, New York, 1982.
\end{enumerate}

Among these, one of my favorite is called ``Closer''. How close can you
be to someone without actually \emph{being them}? Would temporarily
merging identities with someone you loved help you understand them
better? Luckily for you penny-pinchers out there, this particular story
is available free at the following website:

\begin{enumerate}
\def\labelenumi{\arabic{enumi})}
\setcounter{enumi}{2}
\tightlist
\item
  Greg Egan, \emph{Closer},
  \url{https://www.gregegan.net/MISC/CLOSER/Closer.html}
\end{enumerate}

Whoops! I'm drifting pretty far from mathematical physics, aren't I?
Self-reference has a lot to do with mathematical logic, but\ldots. To
gently drift back, let me point out that Egan has a website in which he
explains special and general relativity in a nice, nontechnical way:

\begin{enumerate}
\def\labelenumi{\arabic{enumi})}
\setcounter{enumi}{3}
\tightlist
\item
  Greg Egan, \emph{Foundations},
  \url{https://www.gregegan.net/FOUNDATIONS/}
\end{enumerate}

Also, here are some interesting papers:

\begin{enumerate}
\def\labelenumi{\arabic{enumi})}
\setcounter{enumi}{4}
\item
  Gordon L. Kane, ``Experimental evidence for more dimensions
  reported'', \emph{Physics Today}, May 1998, 13--16.

  Paul M. Grant, ``Researchers find extraordinarily high temperature
  superconductivity in bio-inspired nanopolymer'', \emph{Physics Today},
  May 1998, 17--19.

  Jack Watrous, ``Ribosomal robotics approaches critical experiments;
  government agencies watch with mixed interest'', \emph{Physics Today},
  May 1998, 21--23.
\end{enumerate}
\noindent
What these papers have in common is that they are all works of science
fiction, not science. They read superficially like straight science
reporting, but they are actually the winners of \emph{Physics Today}'s
``Physics Tomorrow'' essay contest!

For example, Grant writes:

\begin{quote}
Little's concept involved replacing the phonons --- characterized by
the Debye temperature --- with excitons, whose much higher
characteristic energies are on the order of 2 eV, or 23,000 K. If
excitons were to become the electron-pairing `glue', superconductors
with \(T_c\)'s as high as 500 K might be possible, even under weak
coupling conditions. Little even proposed a possible realization of the
idea: a structure composed of a conjugated polymer chain (polyene)
dressed with highly polarizable molecule (aromatics) as side groups.
Simply stated, the polyene chain would be a normal metal with a single
mobile electron per C-H molecular unit; electrons on separate units
would be paired by interacting with the exciton field on the polarizable
side groups.
\end{quote}

Actually, I think this part is perfectly true --- William A. Little
suggested this way to achieve high-temperature superconductivity back in
the 1960s. The science fiction part is just the description, later on in
Grant's article, of how Little's dream is actually achieved.

Okay, enough science fiction! Time for some real science! Quantum
gravity, that is. (Stop snickering, you skeptics\ldots.)

\begin{enumerate}
\def\labelenumi{\arabic{enumi})}
\setcounter{enumi}{5}
\tightlist
\item
  Laurent Freidel and Kirill Krasnov, ``Spin foam models and the
  classical action principle'', \emph{Adv. Theor. Math. Phys.} \textbf{2} (1999), 
  1183--1247.  Also available as
 \href{https://arxiv.org/abs/hep-th/9807092}{\texttt{hep-th/}}
  \href{https://arxiv.org/abs/hep-th/9807092}{\texttt{9807092}}.
\end{enumerate}
\noindent
I described the spin foam approach to quantum gravity in
\protect\hyperlink{week113}{``Week 113''}. But let me remind you how the
basic idea goes. A good way to get a handle on this idea is by analogy
with Feynman diagrams. In ordinary quantum field theory there is a
Hilbert space of states called ``Fock space''. This space has a basis of
states in which there are a specific number of particles at specific
positions. We can visualize such a state simply by imagining a bunch of
points in space, with labels to say which particles are which kinds:
electrons, quarks, and so on. One of the main jobs of quantum field
theory is to let us compute the amplitude for one such state to evolve
into another as time passes. Feynman showed that we can do it by
computing a sum over graphs in spacetime. These graphs are called
Feynman diagrams, and they represent ``histories''. For example, \[
  \begin{tikzpicture}[thick]
    \draw (0,0) to node[rotate=225,label={[label distance=-2mm]left:{u}}]{$\blacktriangle$} (1,-1);
    \draw (1,-1) to node[rotate=-225,label={[label distance=-3.2mm]left:{d}}]{$\blacktriangle$} (0,-2);
    \draw [decorate, decoration=snake] (1,-1) to node[label=above:{W}]{} (3,-1);
    \draw (4,0) to node[rotate=-225,label={[label distance=-2mm]right:{e}}]{$\blacktriangle$} (3,-1);
    \draw (3,-1) to node[rotate=225,label={[label distance=-2mm]right:{$\nu$}}]{$\blacktriangle$} (4,-2);
  \end{tikzpicture}
\] would represent a history in which an up quark emits a W boson and
turns into a down quark, with the W being absorbed by an electron,
turning it into a neutrino. Time passes as you march down the page.
Quantum field theory gives you rules for computing amplitudes for any
Feyman diagram. You sum these amplitudes over all Feynman diagrams
starting at one state and ending at another to get the total amplitude
for the given transition to occur.

Now, where do these rules for computing Feynman diagram amplitudes come
from? They are not simply postulated. They come from perturbation
theory. There is a general abstract formula for computing amplitudes in
quantum field theory, but it's not so easy to use this formula in
concrete calculations, except for certain very simple field theories
called ``free theories''. These theories describe particles that don't
interact at all. They are mathematically tractable but physically
uninteresting. Uninteresting, that is, \emph{except} as a starting-point
for studying the theories we \emph{are} interested in --- the so-called
``interacting theories''.

The trick is to think of an interacting theory as containing parameters,
called ``coupling constants'', which when set to zero make it reduce to
a free theory. Then we can try to expand the transition amplitudes we
want to know as a Taylor series in these parameters. As usual, computing
the coefficients of the Taylor series only requires us to to compute a
bunch of derivatives. And we can compute these derivatives using the
free theory! Typically, computing the \(n\)th derivative of some
transition amplitude gives us a bunch of integrals which correspond to
Feynman diagrams with n vertices.

By the way, this means you have to take the particles you see in Feynman
diagrams with a grain of salt. They don't arise purely from the
mathematics of the interacting theory. They arise when we
\emph{approximate} that theory by a free theory. This is not an idle
point, because we can take the same interacting theory and approximate
it by \emph{different} free theories. Depending on what free theory we
use, we may say different things about which particles our interacting
theory describes! In condensed matter physics, people sometimes use the
term ``quasiparticle'' to describe a particle that appears in a free
theory that happens to be handy for some problem or other. For example,
it can be helpful to describe vibrations in a crystal using ``phonons'',
or waves of tilted electron spins using ``spinons''. Condensed matter
theorists rarely worry about whether these particles ``really exist''.
The question of whether they ``really exist'' is less interesting than
the question of whether the particular free theory they inhabit provides
a good approximation for dealing with a certain problem. Particle
physicists, too, have increasingly come to recognize that we shouldn't
worry too much about which elementary particles ``really exist''.

But I digress! My point was simply to say that Feynman diagrams arise
from approximating interacting theories by free theories. The details
are complicated and in most cases nobody has ever succeeded in making
them mathematically rigorous, but I don't want to go into that here.
Instead, I want to turn to spin foams.

Everything I said about Feynman diagrams has an analogy in this approach
to quantum gravity. The big difference is that ordinary ``free
theories'' are formulated on a spacetime with a fixed metric - usually
Minkowski spacetime, with its usual flat metric. Attempts to approximate
quantum gravity by this sort of free theory failed dismally. Perhaps the
fundamental reason is that general relativity doesn't presume that
spacetime has a fixed metric --- au contraire, it's a theory in which
the metric is the main variable!

So the idea of Freidel and Krasnov is to approximate quantum graivty
with a very different sort of ``free theory'', one in which the metric
is a variable. The theory they use is called ``\(BF\) theory''. I said a
lot about \(BF\) theory in \protect\hyperlink{week36}{``Week 36''}, but
here the main point is simply that it's a topological quantum field
theory, or TQFT. A TQFT is a quantum field theory that does not presume
a fixed metric, but of a very simple sort, because it has no local
degrees of freedom. I very much like the idea that a TQFT might serve as
a novel sort of ``free theory'' for the purposes of studying quantum
gravity.

Everything that Freidel and Krasnov do is reminscent of familiar quantum
field theory, but also very different, because their starting-point is
\(BF\) theory rather than a free theory of a traditional sort. For
example, just as ordinary quantum field theory starts out with Fock
space, in the spin network approach to quantum gravity we start with a
nice simple Hilbert space of states. But this space has a basis
consisting, not of collections of 0-dimensional particles sitting in
space at specified positions, but of \(1\)-dimensional ``spin networks''
sitting in space. (For more on spin networks, see
\protect\hyperlink{week55}{``Week 55''} and
\protect\hyperlink{week110}{``Week 110''}.) And instead of using
\(1\)-dimensional Feynman diagrams to compute transition amplitudes, the
idea is now to use 2-dimensional gadgets called ``spin foams''. The
amplitudes for spin foams are easy to compute in \(BF\) theory, because
there are a lot of explicit formulas using the so-called ``Kauffman
bracket'', which is an easily computable invariant of spin networks. So
then the trick is to use this technology to compute spin foam amplitudes
for quantum gravity.

Now, I shouldn't give you the wrong impression here. There are lots of
serious problems and really basic open questions in this work, and the
whole thing could turn out to be fatally flawed somehow. Nonetheless,
something seems right about it, so I find it very interesting.

Anyway, on to some other papers. I'm afraid I don't have enough energy
for detailed descriptions, because I'm busy moving into a new house, so
I'll basically just point you at them\ldots.

\begin{enumerate}
\def\labelenumi{\arabic{enumi})}
\setcounter{enumi}{6}
\tightlist
\item
  Abhay Ashtekar, Alejandro Corichi and Jose A. Zapata, ``Quantum theory
  of geometry III: non-commutativity of Riemannian structures'',
   available as
  \href{https://arxiv.org/abs/gr-qc/9806041}{\texttt{gr-qc/9806041}}.
\end{enumerate}

This is the long-awaited third part of a series giving a mathematically
rigorous formalism for interpreting spin network states as ``quantum
3-geometries'', that is, quantum states describing the metric on
3-dimensional space together with its extrinsic curvature (as it sits
inside \(4\)-dimensional spacetime). Here's the abstract:

\begin{quote}
``The basic framework for a systematic construction of a quantum theory
of Riemannian geometry was introduced recently. The quantum versions of
Riemannian structures --- such as triad and area operators --- exhibit a
non-commutativity. At first sight, this feature is surprising because it
implies that the framework does not admit a triad representation. To
better understand this property and to reconcile it with intuition, we
analyze its origin in detail. In particular, a careful study of the
underlying phase space is made and the feature is traced back to the
classical theory; there is no anomaly associated with quantization. We
also indicate why the uncertainties associated with this
non-commutativity become negligible in the semi-classical regime.''
\end{quote}

In case you're wondering, the ``triad'' field is more or less what
mathematicians would call a ``frame field'' or ``soldering form'' ---
and it's the same as the ``B'' field in \(BF\) theory. It encodes the
information about the metric in Ashtekar's formulation to general
relativity.

Moving on to matters \(n\)-categorical, we have:

\begin{enumerate}
\def\labelenumi{\arabic{enumi})}
\setcounter{enumi}{7}
\tightlist
\item
  Andre Hirschowitz, Carlos Simpson, ``Descente pour les \(n\)-champs''
  (Descent for \(n\)-stacks),
  available as
  \href{https://arxiv.org/abs/math.AG/9807049}{\texttt{math.AG/9807049}}.
\end{enumerate}
\noindent
Apparently this provides a theory of ``\(n\)-stacks'', which are the
\(n\)-categorical generalization of sheaves. Ever since Grothendieck's
600-page letter to Quillen (see \protect\hyperlink{week35}{``Week
35''}), this has been the holy grail of \(n\)-category theory.
Unfortunately I haven't mustered sufficient courage to force my way
through 240 pages of French, so I don't really know the details!

For the following two \(n\)-category papers, exploring some themes close
to my heart, I'll just quote the abstracts:

\begin{enumerate}
\def\labelenumi{\arabic{enumi})}
\setcounter{enumi}{8}
\tightlist
\item
  Michael Batanin, ``Computads for finitary monads on globular sets'',
  available at
\href{https://citeseerx.ist.psu.edu/pdf/5ebfd5445a832eb9a104b2f9ef7ae9bc58741af6}{\texttt{https://citeseerx.ist.psu.edu/pdf/5ebfd5445a832eb9a104b2f9ef7ae9bc58}} 
\break
\href{https://citeseerx.ist.psu.edu/pdf/5ebfd5445a832eb9a104b2f9ef7ae9bc58741af6}
{\texttt{741af6}}
\end{enumerate}

\begin{quote}This work arose as a reflection on the foundation of higher dimensional
category theory. One of the main ingredients of any proposed definition
of weak \(n\)-category is the shape of diagrams (pasting scheme) we
accept to be composable. In a globular approach {[}due to Batanin{]}
each \(k\)-cell has a source and target \((k-1)\)-cell. In the opetopic
approach of Baez and Dolan and the multitopic approach of Hermida,
Makkai and Power each \(k\)-cell has a unique \((k-1)\)-cell as target
and a whole \((k-1)\)-dimensional pasting diagram as source. In the
theory of strict \(n\)-categories both source and target may be a
general pasting diagram.

The globular approach being the simplest one seems too restrictive to
describe the combinatorics of higher dimensional compositions. Yet, we
argue that this is a false impression. Moreover, we prove that this
approach is a basic one from which the other type of composable diagrams
may be derived. One theorem proved here asserts that the category of
algebras of a finitary monad on the category of n-globular sets is
\emph{equivalent} to the category of algebras of an appropriate monad on
the special category (of computads) constructed from the data of the
original monad. In the case of the monad derived from the universal
contractible operad this result may be interpreted as the equivalence of
the definitions of weak \(n\)-categories (in the sense of Batanin) based
on the `globular' and general pasting diagrams. It may be also
considered as the first step toward the proof of equivalence of the
different definitions of weak \(n\)-category.

We also develop a general theory of computads and investigate some
properties of the category of generalized computads. It turned out, that
in a good situation this category is a topos (and even a presheaf topos
under some not very restrictive conditions, the property firstly
observed by S. Schanuel and reproved by A. Carboni and P. Johnstone for
2-computads in the sense of Street).
\end{quote}

\begin{enumerate}
\def\labelenumi{\arabic{enumi})}
\setcounter{enumi}{9}
\tightlist
\item
  Tom Leinster, ``Structures in higher-dimensional category theory'', \break
  available as
  \href{https://arxiv.org/abs/math/0109021}{\texttt{math/0109021}}
\end{enumerate}

\begin{quote}
This is an exposition of some of the constructions which have arisen
in higher-dimensional category theory. We start with a review of the
general theory of operads and multicategories. Using this we give an
account of Batanin's definition of \(n\)-category; we also give an
informal definition in pictures. Next we discuss Gray-categories and
their place in coherence problems. Finally, we present various
constructions relevant to the opetopic definitions of \(n\)-category.
New material includes a suggestion for a definition of lax cubical
\(n\)-category; a characterization of small Gray-categories as the small
substructures of \(2\text{-}\mathsf{Cat}\); a conjecture on coherence
theorems in higher dimensions; a construction of the category of trees
and, more generally, of \(n\)-pasting diagrams; and an analogue of the
Baez--Dolan slicing process in the general theory of operads.
\end{quote}

Okay --- now for something completely different. In
\protect\hyperlink{week122}{``Week 122''} I said how Kreimer and Connes
have teamed up to write a paper relating Hopf algebras, renormalization,
and noncommutative geometry. Now it's out:

\begin{enumerate}
\def\labelenumi{\arabic{enumi})}
\setcounter{enumi}{10}
\tightlist
\item
  Alain Connes and Dirk Kreimer, ``Hopf algebras, renormalization and
  noncommutative geometry'', \emph{Comm. Math. Phys.} \textbf{199} (1998),             
  203--242.  Also available as
  \href{https://arxiv.org/abs/hep-th/9808042}{\texttt{hep-th/9808042}}.
\end{enumerate}
\noindent
Also, here's an introduction to Kreimer's work:

\begin{enumerate}
\def\labelenumi{\arabic{enumi})}
\setcounter{enumi}{11}
\tightlist
\item
  Dirk Kreimer, ``How useful can knot and number theory be for loop
  calculations?'', Talk given at the workshop \emph{Loops and Legs in
  Gauge Theories}, available as
  \href{https://arxiv.org/abs/hep-th/9807125}{\texttt{hep-th/}}
  \href{https://arxiv.org/abs/hep-th/9807125}{\texttt{9807125}}
\end{enumerate}

Switching over to homotopy theory and its offshoots\ldots{} when I
visited Dan Christensen at Johns Hopkins this spring, he introduced me
to all the homotopy theorists there, and Jack Morava gave me a paper
which really indicates the extent to which new-fangled ``quantum
topology'' has interbred with good old-fashioned homotopy theory:

\begin{enumerate}
\def\labelenumi{\arabic{enumi})}
\setcounter{enumi}{11}
\tightlist
\item
  Jack Morava, ``Quantum generalized cohomology'', available as
  \href{https://arxiv.org/abs/math.QA/9807058}{\texttt{math.QA/9807058}}.
\end{enumerate}
\noindent
Again, I'll just quote the abstract rather than venturing my own
summary:

\begin{quote}
We construct a ring structure on complex cobordism tensored with the
rationals, which is related to the usual ring structure as quantum
cohomology is related to ordinary cohomology. The resulting object
defines a generalized two-dimensional topological field theory taking
values in a category of spectra.
\end{quote}

Finally, Morava has a student who gave me an interesting paper on
operads and moduli spaces:

\begin{enumerate}
\def\labelenumi{\arabic{enumi})}
\setcounter{enumi}{12}
\tightlist
\item
  Satyan L. Devadoss, ``Tessellations of moduli spaces and the mosaic
  operad'', in \emph{Homotopy Invariant Algebraic Structures---A Conference 
  in Honor of J. Michael Boardman}, edited by J. P. Meyer, J. Morava and W. S.
   Wilson, \emph{Contemp.\ Math.} \textbf{239}.    Also available as
  \href{https://arxiv.org/abs/math.QA/9807010}{\texttt{math.QA/9807010}}.
\end{enumerate}

\begin{quote}
``We construct a new (cyclic) operad of `mosaics' defined by polygons
with marked diagonals. Its underlying (aspherical) spaces are the sets
\(M_{0,n}(\mathbb{R})\) of real points of the moduli space of punctured
Riemann spheres, which are naturally tiled by Stasheff associahedra. We
(combinatorially) describe them as iterated blow-ups and show that their
fundamental groups form an operad with similarities to the operad of
braid groups.''
\end{quote}

\begin{center}\rule{0.5\linewidth}{0.5pt}\end{center}

\begin{quote}
\emph{Some things are so serious that one can only jest about them.}

--- Niels Bohr
\end{quote}

\hypertarget{week124}{%
\section{October 23, 1998}\label{week124}}

I'm just back from Tucson, where I talked a lot with my friend Minhyong
Kim, who teaches at the math department of the University of Arizona. I
met Minhyong in 1986 when I was a postdoc and he was a grad student at
Yale. At the time, strings were all the rage. Having recently found 5
consistent superstring theories, many physicists were giddy with
optimism, some even suggesting that the Theory of Everything would be
completed before the turn of the century. A lot of mathematicians were
going along for the ride, delighted by the beautiful and intricate
mathematical infrastructure: conformal field theory, vertex operator
algebras, and so on. Minhyong was considering doing his thesis on one of
these topics, so we spent a lot of time talking about mathematical
physics.

However, he eventually decided to work with Serge Lang on arithmetic
geometry. This is a branch of algebraic geometry where you work over the
integers instead of a field --- especially important for Diophantine
equations. Personally, I was a bit disappointed. Perhaps it was because
I thought physics was more important than the decadent pleasures of pure
mathematics --- or perhaps it was because it made it much less likely
that we'd ever collaborate on a paper.

However, a lot of the math Minhyong learned when studying string theory
is also important in arithmetic geometry. An example is the theory of
elliptic curves. Roughly speaking, an elliptic curve is a torus formed
taking a parallelogram in the complex plane and identifying opposite
edges.

You might wonder why something basically doughnut-shaped is called an
elliptic curve! Let's clear that up right away. The ``elliptic'' part
comes from a relationship to elliptic functions, which generalize the
familiar trig functions from circles to ellipses. The ``curve'' part
comes from the fact that it takes one complex number \(z = x+iy\) to
describe your location on a surface with two real coordinates \((x,y)\),
so showoffs like to say that a torus is one-dimensional --- one
\emph{complex} dimension, that is! --- hence a ``curve''. In short, you
have to already understand elliptic curves to know why the heck they're
called elliptic curves.

Anyway, why are elliptic curves important? On the one hand, they show up
all over number theory, like in Wiles' proof of Fermat's last theorem.
On the other hand, in string theory, a string traces out a surface in
spacetime called the string worldsheet, and points on this surface are
conveniently described using a single complex number, so it's what those
showoffs call a ``curve'' --- and among the simplest examples are
elliptic curves!

If you're interested to see how Fermat's last theorem was reduced to a
problem about elliptic curves --- the so-called Shimura-Taniyama-Weil
conjecture --- you can look at the textbooks on elliptic curves listed
in \protect\hyperlink{week13}{``Week 13''}. But I won't say anything
about this, since I don't understand it. Instead, I want to talk about
how elliptic curves show up in string theory. For more on how these two
applications fit together, try:

\begin{enumerate}
\def\labelenumi{\arabic{enumi})}
\tightlist
\item
  Yuri I. Manin, ``Reflections on arithmetical physics'', in
  \emph{Conformal Invariance and String Theory}, eds.~Petre Dita and
  Vladimir Georgescu, Academic Press, 1989.
\end{enumerate}

Let me just quote the beginning:

\begin{quote}
The development of theoretical physics in the last quarter of the
twentieth century is guided by a very romantic system of values.
Aspiring to describe fundamental processes at the Planck scale,
physicists are bound to lose any direct connection with the observable
world. In this social context the sophisticated mathematics emerging in
string theory ceases to be only a technical tool needed to calculate
some measurable effects and becomes a matter of principle.

Today at least some of us are again nurturing an ancient Platonic
feeling that mathematical ideas are somehow predestined to describe the
physical world, however remote from reality their origins seem to be.

From this viewpoint one should perversely expect number theory to become
the most applicable branch of mathematics."
\end{quote}

I think this remark wisely summarizes both the charm and the dangers of
physics that relies more heavily on criteria of mathematical elegance
than of experimental verification.

Anyway, I don't want to get too deep into the theory of elliptic curves;
just enough so we see why the number 24 is so important in string
theory. You may remember that bosonic string theory works best in 26
dimensions (while the physically more important superstring theory,
which includes spin-\(1/2\) particles, works best in 10). Why is this
true? Well, there are various answers, but one is that if you think of
the string as wiggling in the 24 directions perpendicular to its own
2-dimensional surface --- two \emph{real} dimensions, that is! ---
various magical properties of the number 24 conspire to make things work
out.

What are these magical properties of the number 24? Well,
\[1^2 + 2^2 + 3^2 + \ldots + 24^2\] is itself a perfect square, and 24
is the only integer with this property besides silly ones like 0 and 1.
As described in \protect\hyperlink{week95}{``Week 95''}, this has some
very profound relationships to string theory. Unfortunately, I don't
know any way to deduce from this that bosonic string theory \emph{works
best} in 26 dimensions.

One reason bosonic string theory works best in 26 dimensions is that
\[1 + 2 + 3 + \ldots = -\frac{1}{12}\] and \(2 \times 12 = 24\). Of
course, this explanation is unsatisfactory in many ways. First of all,
you might wonder what the above equation means! Doesn't the sum
diverge???

Actually this is the \emph{least} unsatisfactory feature of the
explanation. Although the sum diverges, you can still make sense of it.
The Riemann zeta function is defined by
\[\zeta(s) = 1^{-s} + 2^{-s} + 3^{-s} + \ldots\] whenever the real part
of \(s\) is greater than \(1\), which makes the sum converge. But you
can analytically continue it to the whole complex plane, except for a
pole at \(1\). If you do this, you find that
\[\zeta(-1) = -\frac{1}{12}.\] Thus we may jokingly say that
\(1 + 2 + 3 + \ldots = -1/12\). But the real point is how the zeta
function shows up in string theory, and quantum field theory in general.
(It's also big in number theory.)

Unfortunately, the details quickly get rather technical; one has to do
some calculations and so on. That's the really unsatisfactory part. I
want something that clearly relates strings and the number 24, something
so simple even a child could understand it, and which, when you work out
all the implications, implies that bosonic string theory only makes
sense in 26 dimensions. I don't expect a child to be able to figure out
all the implications\ldots{} but I want the essence to be childishly
simple.

Here it is. Suppose the string worldsheet is an elliptic curve. Then we
can make it by taking a ``lattice'' of parallelograms in the complex
plane: \[
  \begin{tikzpicture}[scale=0.7]
    \draw[->] (-3,0) to (4,0) node[label=below:{$\Re(z)$}]{};
    \draw[->] (0,-3) to (0,4) node[label=left:{$\Im(z)$}]{};
    \foreach \m in {-1,0,1,2}
    {
      \foreach \n in {-1,0,1,2}
      {
        \node at ({\m*1.5-\n/3-0.2},{1.5*\n+\m-0.5}) {$\bullet$};
      }
    }
  \end{tikzpicture}
\] and identifying each point in each parallelogram with the
corresponding points on all the others. This rolls the plane up into a
torus. Now, two lattices are more symmetrical than the rest. One of them
is the square lattice: \[
  \begin{tikzpicture}[scale=0.7]
    \draw[->] (-2.5,0) to (4,0) node[label=below:{$\Re(z)$}]{};
    \draw[->] (0,-2) to (0,4) node[label=left:{$\Im(z)$}]{};
    \foreach \m in {-1,0,1,2}
    {
      \foreach \n in {-1,0,1,2}
      {
        \node at ({1.5*\m},{1.6*\n}) {$\bullet$};
      }
    }
  \end{tikzpicture}
\] which has 4-fold rotational symmetry. The other is the lattice with
lots of equilateral triangles in it: \[
  \begin{tikzpicture}[scale=0.7]
    \draw[->] (-3,0) to (4.5,0) node[label=below:{$\Re(z)$}]{};
    \draw[->] (0,-2) to (0,3.5) node[label=left:{$\Im(z)$}]{};
    \foreach \m in {-1,0,1,2}
    {
      \foreach \n in {-1,0,1,2}
      {
        \node at ({1.5*\m+0.75*\n},{1.33*\n}) {$\bullet$};
      }
    }
  \end{tikzpicture}
\] which has 6-fold rotational symmetry. The magic property of the
number 24, which makes string theory work so well in 26 dimensions, is
that \[4 \times 6 = 24\]

Okay, great. But if you're anything like me, at this point you're
wondering how the heck this actually helps. Why should string theory
care about these specially symmetrical lattices? And why should we
\emph{multiply} 4 and 6? So far everything I've said has been flashy but
insubstantial. Next week I'll fill in some of the details. Of course,
I'll need to turn up the sophistication level a notch or two.

In the meantime, you can read a bit more about this stuff in the
following article on Richard Borcherds, who won the Fields medal for his
work relating bosonic string theory, the Leech lattice in 24 dimensions,
and the Monster group:

\begin{enumerate}
\def\labelenumi{\arabic{enumi})}
\setcounter{enumi}{1}
\tightlist
\item
  W. Wayt Gibbs, ``Monstrous moonshine is true'', \emph{Scientific
  American}, November 1998, 40--41. Also available at
  \texttt{http://www.sciam.com/1998/1198issue/1198profile.html}.
\end{enumerate}

Gibbs asked me to come up with a simple explanation of the \(j\)
invariant for elliptic curves; you can \(j\)udge how well I succeeded.
For a more detailed attempt to do the same thing, see
\protect\hyperlink{week66}{``Week 66''}, which also has more references
on the Monster group. By the way, John McKay didn't actually make his
famous discovery relating the \(j\) invariant and Monster while reading
a 19th-century book on elliptic modular functions; he says ``It was du
Val's Elliptic Functions book in which \(j\) is expanded incorrectly as
a \(q\)-series --- very much a 20th century book.'' Apart from that, the
article seems accurate, as far as I can tell.

If you really want to understand how elliptic curves are related to
strings, you need to learn some conformal field theory. For that, try:

\begin{enumerate}
\def\labelenumi{\arabic{enumi})}
\setcounter{enumi}{2}
\tightlist
\item
  Phillippe Di Francesco, Pierre Mathieu, and David Senechal,
  \emph{Conformal Field Theory}, Springer, 1997.
\end{enumerate}

This is a truly wonderful tour of the subject. It's 890 pages long, but
it's designed to be readable by both mathematicians and physicists, so
you can look at the bits you want. It starts out with a 60-page
introduction to quantum field theory and a 30-page introduction to
statistical mechanics. The reason is that when we perform the
substitution called the ``Wick transform'': \[it/\hbar\mapsto k/T,\]
quantum field theory turns into statistical mechanics, and a nice
Lorentzian manifold may turn into a Riemannian manifold --- in other
words, ``spacetime'' turns into ``space''. And this gives conformal
field theory a double personality.

First, conformal field theory studies quantum field theories in 2
dimensions that are invariant under all conformal transformations ---
transformations that preserve angles but not necessarily lengths. These
are important in string theory because we can think of them as
transformations of the string worldsheet that preserve its complex
structure.

Secondly, if we do a Wick transform, these quantum field theories become
2-dimensional \emph{statistical mechanics} problems that are invariant
under all conformal transformations. This may seem an esoteric concern,
but thin films of material can often be treated as \(2\)-dimensional for
all practical purposes, and conformal invariance is typical at
``critical points'' --- boundaries between two or more phases for which
there is no latent heat, such as the boundary between the magnetized and
unmagnetized phases of a ferromagnet. In 2 dimensions, one can use
conformal field theory to thoroughly understand these critical points.

After this warmup, the book covers the fundamentals of conformal field
theory proper, including:

\begin{itemize}
\tightlist
\item
  the idea of conformal invariance (which is especially powerful in 2
  dimensions because then the group of conformal transformations is
  infinite-dimensional),
\item
  the free boson and fermion fields,
\item
  operator product expansions,
\item
  the Virasoro algebra (which is closely related to the Lie algebra of
  the group of conformal transformations, and has a representation on
  the Hilbert space of states of any conformal field theory),
\item
  minimal models (roughly, conformal field theories whose Hilbert space
  is built from finitely many irreducible representations of the
  Virasoro algebra),
\item
  the Coulomb-gas formalism (a way to describe minimal models in terms
  of the free boson and fermion fields),
\item
  modular invariance (the study of conformal field theory on tori ---
  this is where the elliptic curves start sneaking into the picture,
  dragging along with them the wonderful machinery of elliptic
  functions, theta functions, the Dedekind eta function, and so forth),
\item
  critical percolation (applying conformal field theory to systems where
  a substance is trying to ooze through a porous medium, with special
  attention paid to the critical point when the holes are \emph{just}
  big enough to let it ooze all the way through),
\item
  the \(2\)-dimensional Ising model (applying conformal field theory to
  ferromagnets, with special attention paid to the critical point when
  the temperature is \emph{just} low enough for ferromagnetism to set
  in)
\end{itemize}

By now we're at page 486. I'm getting tired just summarizing this thing!

Anyway, the book then turns to conformal field theories having Lie group
symmetries: in particular, the so-called Wess-Zumino-Witten or ``WZW''
models. Pure mathematicians are free to join here, even amateurs,
because we are now treated to a wonderful 78-page introduction to simple
Lie algebras, starting from scratch and working rapidly through all
sorts of fun stuff, skipping all the yucky proofs. Then we get a 54-page
introduction to affine Lie algebras, which are infinite- dimensional
generalizations of the simple Lie algebras, and play a crucial role in
string theory. Finally, we get a detailed 143-page course on WZW models
--- which are basically conformal field theories where your field takes
values in a Lie group --- and coset models --- where your field takes
values in a Lie group modulo a subgroup. It sounds like all minimal
models can be described as coset models, though I'm not quite sure.

Whew! Believe it or not, the authors plan a second volume! Anyway, this
is a wonderful book to have around. I was just about to buy a copy in
Chicago last spring --- on sale for a mere \$50 --- when I discovered
I'd lost my credit card. Sigh. The big ones always get away\ldots.

There are various formalisms for doing conformal field theory that
aren't covered in the above text. For example, the theory of ``vertex
operator algebras'', or ``vertex algebras'' is really popular among
mathematicians studying conformal field theory and the Monster group.

The standard definition of a vertex operator algebra is long and
complicated: it summarizes a lot of what you'd want a conformal field
theory to be like, but it's hard to learn to love it unless you already
know some \emph{other} approaches to conformal field theory. There's
another definition using operads that's much nicer, which will
eventually catch on --- some people complain that operads are too
abstract, but that's just hogwash. But anyway, there is a definite need
for more elementary texts on the subject. Here's one:

\begin{enumerate}
\def\labelenumi{\arabic{enumi})}
\setcounter{enumi}{3}
\tightlist
\item
  Victor Kac, \emph{Vertex Algebras for Beginners}, American
  Mathematical Society, University Lecture Series vol.~\textbf{10},
  1997.
\end{enumerate}

And then of course there is string theory proper. How do you learn that?
There's always the bible by Green, Schwarz and Witten (see
\protect\hyperlink{week118}{``Week 118''}), but a lot of stuff has
happened since that was written. Luckily, Joseph Polchinski has come out
with a ``new testament''; I haven't seen it yet but physicists say it's
very good:

\begin{enumerate}
\def\labelenumi{\arabic{enumi})}
\setcounter{enumi}{4}
\tightlist
\item
  Joseph Polchinski, \emph{String Theory}, two volumes, Cambridge U.\
  Press, 1998.
\end{enumerate}

There are also other textbooks, of course. Here's one that's free if you
print it out yourself:

\begin{enumerate}
\def\labelenumi{\arabic{enumi})}
\setcounter{enumi}{5}
\tightlist
\item
  E. Kiritsis, \emph{Introduction to Superstring Theory}, available as
  \href{https://arxiv.org/abs/hep-th/9709062}{\texttt{hep-th/9709062}}.
\end{enumerate}

For a more mathematical approach, you might want to try this when it
comes out:

\begin{enumerate}
\def\labelenumi{\arabic{enumi})}
\setcounter{enumi}{6}
\tightlist
\item
  \emph{Quantum Fields and Strings: A Course for Mathematicians},
  eds.~P. Deligne, P. Etinghof, D. Freed, L. Jeffrey, D. Kazhdan, D.
  Morrison and E. Witten, American Mathematical Society, to appear.
\end{enumerate}

Finally, when you get sick of all this new-fangled stuff and want to
read about the good old days when physicists predicted new particles
that actually wound up being \emph{observed}, you can turn to this book
about Dirac and his work:

\begin{enumerate}
\def\labelenumi{\arabic{enumi})}
\setcounter{enumi}{7}
\tightlist
\item
  Abraham Pais, Maurice Jacob, David I. Olive, and Michael F. Atiyah,
  \emph{Review of Paul Dirac: The Man and His Work}, Cambridge U.\ Press,
  1998.
\end{enumerate}

Also try this:

\begin{enumerate}
\def\labelenumi{\arabic{enumi})}
\setcounter{enumi}{8}
\tightlist
\item
  Michael Berry, ``Paul Dirac: the purest soul in physics'',
  \emph{Physics World}, February 1998, pp.~36--40.
\end{enumerate}

\begin{center}\rule{0.5\linewidth}{0.5pt}\end{center}

\begin{quote}
\emph{First, and above all for Dirac, the logic that led to the theory
was, although deeply sophisticated, in a sense beautifully simple. Much
later, when someone asked him (as many must have done before)} ``How did
you find the Dirac equation?'' \emph{he is said to have replied} ``I
found it beautiful.''

--- Michael Berry
\end{quote}

\hypertarget{week125}{%
\section{November 3, 1998}\label{week125}}

Last week I promised to explain some mysterious connections between
elliptic curves, string theory, and the number 24. I claimed that it all
boils down to the fact that there are two especially symmetric lattices
in the plane, namely the square lattice: \[
  \begin{tikzpicture}[scale=0.7]
    \draw[->] (-2.5,0) to (4,0) node[label=below:{$\Re(z)$}]{};
    \draw[->] (0,-2) to (0,4) node[label=left:{$\Im(z)$}]{};
    \foreach \m in {-1,0,1,2}
    {
      \foreach \n in {-1,0,1,2}
      {
        \node at ({1.5*\m},{1.6*\n}) {$\bullet$};
      }
    }
  \end{tikzpicture}
\] with 4-fold symmetry, and the hexagonal lattice: \[
  \begin{tikzpicture}[scale=0.7]
    \draw[->] (-3,0) to (4.5,0) node[label=below:{$\Re(z)$}]{};
    \draw[->] (0,-2) to (0,3.5) node[label=left:{$\Im(z)$}]{};
    \foreach \m in {-1,0,1,2}
    {
      \foreach \n in {-1,0,1,2}
      {
        \node at ({1.5*\m+0.75*\n},{1.33*\n}) {$\bullet$};
      }
    }
  \end{tikzpicture}
\] with 6-fold symmetry. Now it's time for me to start backing up those
claims.

First I need to talk a bit about lattices and
\(\mathrm{SL}(2,\mathbb{Z})\). As I explained in
\protect\hyperlink{week66}{``Week 66''}, a lattice in the complex plane
consists of all points that are integer linear combinations of two
complex numbers, say \(\omega_1\) and \(\omega_2\). However, we can
change these numbers without changing the lattice by letting \[
  \begin{aligned}
    \omega'_1 &= a\omega_1+b\omega_2
  \\\omega'_2 &= c\omega_1+d\omega_2
  \end{aligned}
\] where \[
  \left(
    \begin{array}{cc}
      a&b\\c&d
    \end{array}
  \right)
\] is a \(2\times2\) invertible matrix of integers whose inverse again
consists of integers. Usually it's good to require that our
transformation preserve the handedness of the basis
\((\omega_1,\omega_2)\), which means that this matrix should have
determinant \(1\). Such matrices form a group called
\(\mathrm{SL}(2,\mathbb{Z})\). In the context of elliptic curves it's
also called the ``modular group''.

Now associated to the square lattice is a special element of
\(\mathrm{SL}(2,\mathbb{Z})\) that corresponds to a 90 degree rotation.
Everyone calls it \(S\): \[
  S = \left(
    \begin{array}{cc}
      0&-1\\1&0
    \end{array}
  \right)
\] Associated to the hexagonal lattice is a special element of
\(\mathrm{SL}(2,\mathbb{Z})\) that corresponds to a 60 degree rotation.
Everyone calls it \(ST\): \[
  ST = \left(
    \begin{array}{cc}
      0&-1\\1&1
    \end{array}
  \right)
\] (See, there's a matrix they already call \(T\), and \(ST\) is the
product of \(S\) and that one.) Now, you may complain that the matrix
\(ST\) doesn't look like a rotation, but you have to be careful! What I
mean is, if you take the hexagonal lattice and pick a basis for it like
this: \[
  \begin{tikzpicture}[scale=0.7]
    \draw[->] (-3,0) to (4,0) node[label=below:{$\Re(z)$}]{};
    \draw[->] (0,-2) to (0,4) node[label=left:{$\Im(z)$}]{};
    \foreach \m in {-1,0,1,2}
    {
      \foreach \n in {-1,0,1,2}
      {
        \node (\m\n) at ({1.5*\m+0.75*\n},{1.33*\n}) {$\circ$};
      }
    }
    \node at (00) {$\bullet$};
    \node at (-0.3,-0.33) {\scriptsize$0$};
    \node at (10) {$\bullet$};
    \node at (1.5,-0.4) {\scriptsize$\omega_1$};
    \node at (01) {$\bullet$};
    \node at (0.75,0.9) {\scriptsize$\omega_2$};
  \end{tikzpicture}
\] then in \emph{this} basis the matrix \(ST\) represents a 60 degree
rotation.

So far this is pretty straightforward, but now come some surprises.
First, it turns out that \(\mathrm{SL}(2,\mathbb{Z})\) is
\emph{generated} by \(S\) and \(ST\). In other words, every \(2\times2\)
integer matrix with determinant \(1\) can be written as a product of a
bunch of copies of \(S\), \(ST\), and their inverses. Second, all the
relations satisfied by \(S\) and \(ST\) follow from these obvious ones:
\[
  \begin{aligned}
    S^4 &= 1
  \\(ST)^6 &= 1
  \end{aligned}
\] together with \[S^2 = (ST)^3\] which holds because both sides
describe a 180 degree rotation.

Right away this implies that \(\mathrm{SL}(2,\mathbb{Z})\) has a certain
inherent ``12-ness'' to it. Let me explain.
\(\mathrm{SL}(2,\mathbb{Z})\) is a nonabelian group --- this is how
someone with a Ph.D.~says that matrix multiplication doesn't commute ---
but suppose we abelianize it by imposing extra relations \emph{forcing}
commutativity. Then we get a group generated by \(S\) and \(ST\),
satisfying the above relations together with an extra one saying that
\(S\) and \(ST\) commute. This is the group \(\mathbb{Z}/12\), which has
12 elements!

This ``12-ness'' has a lot to do with the magic properties of the number
24 in string theory. But to see how this ``12-ness'' affects string
theory, we need to talk about elliptic curves a bit more. It will take
forever unless I raise the mathematical sophistication level a little.
So\ldots.

We can define an elliptic curve to be a torus \(\mathbb{C}/L\) formed by
taking the complex plane \(\mathbb{C}\) and modding out by a lattice L.
Since \(\mathbb{C}\) is an abelian group and \(L\) is a subgroup, this
torus is an abelian group, but in the theory of elliptic curves we
consider it not just as a group but also as a complex manifold. Thus two
elliptic curves \(\mathbb{C}/L\) and \(\mathbb{C}/L'\) are considered
isomorphic if there is a complex-analytic function from one to the other
that's also an isomorphism of groups. This happens precisely when there
is a nonzero number \(z\) such that \(zL = L'\), or in other words,
whenever \(L'\) is a rotated and/or dilated version of \(L\).

There's a wonderful space called the ``moduli space'' of elliptic
curves: each point on it corresponds to an isomorphism class of elliptic
curves. In physics, we think of each point in it as describing the
geometry of a torus-shaped string worldsheet. Thus in the path-integral
approach to string theory we need to integrate over this space, together
with a bunch of other moduli spaces corresponding to string worldsheets
with different topologies. All these moduli spaces are important and
interesting, but the moduli space of elliptic curves is a nice simple
example when you're first trying to learn this stuff. What does this
space look like?

Well, suppose we have an elliptic curve \(\mathbb{C}/L\). We can take
our lattice \(L\) and describe it in terms of a right-handed basis
\((\omega_1, \omega_2)\). For the purposes of classifying the describing
the elliptic curve up to isomorphism, it doesn't matter if we multiply
these basis elements by some number \(z\), so all that really matters is
the ratio \[\tau = \omega1/\omega2.\] Since our basis was right-handed,
\(\tau\) lives in the upper half-plane, which people like to call \(H\).

Okay, so now we have described our elliptic curve in terms of a complex
number \(\tau\) lying in \(H\). But the problem is, we could have chosen
a different right-handed basis for our lattice \(L\) and gotten a
different number \(\tau\). We've got to think about that. Luckily, we've
already seen how we can change bases without changing the lattice: we
just apply a matrix in \(\mathrm{SL}(2,\mathbb{Z})\), getting a new
basis \[
  \begin{aligned}
    \omega'_1 &= a\omega_1+b\omega_2
  \\\omega'_2 &= c\omega_1+d\omega_2
  \end{aligned}
\] This has the effect of changing \(\tau\) to
\[\tau' = \frac{a \tau + b}{c \tau + d}.\] If you don't see why, figure
it out --- you've gotta understand this to understand elliptic curves!

Anyway, two numbers \(\tau\) and \(\tau'\) describe isomorphic elliptic
curves if and only if they differ by the above sort of transformation.
So we've figured out the moduli space of elliptic curves: it's the
quotient space \(H/\mathrm{SL}(2,\mathbb{Z})\), where
\(\mathrm{SL}(2,\mathbb{Z})\) acts on \(H\) as above!

Now, the quotient space \(H/\mathrm{SL}(2,\mathbb{Z})\) is not a smooth
manifold, because while the upper halfplane \(H\) is a manifold and the
group \(\mathrm{SL}(2,\mathbb{Z})\) is discrete, the action of
\(\mathrm{SL}(2,\mathbb{Z})\) on \(H\) is not free: i.e., certain points
in \(H\) don't move when you hit them with certain elements of
\(\mathrm{SL}(2,\mathbb{Z})\).

If you don't see why this causes trouble, think about a simpler example,
like the group \(G = \mathbb{Z}/n\) acting as rotations of the complex
plane, \(\mathbb{C}\). Most points in the plane move when you rotate
them, but the origin doesn't. The quotient space \(\mathbb{C}/G\) is a
cone with its tip corresponding to the origin. It's smooth everywhere
except the tip, where it has a ``conical singularity''. The moral of the
story is that when we mod out a manifold by a group of symmetries, we
get a space with singularities corresponding to especially symmetrical
points in the original manifold.

So we expect that \(H/\mathrm{SL}(2,\mathbb{Z})\) has singularities
corresponding to points in \(H\) corresponding to especially symmetrical
lattices. These, of course, are our friends the square and hexagonal
lattices!

But let's be a bit more careful. First of all, \emph{nothing} in \(H\)
moves when you hit it with the matrix \(-1\). But that's no big deal: we
can just replace the group \(\mathrm{SL}(2,\mathbb{Z})\) by
\[\mathrm{PSL}(2,\mathbb{Z}) = \mathrm{SL}(2,\mathbb{Z})/\{\pm1\}\]
Since \(-1\) doesn't move \emph{any} points of \(H\), the action of
\(\mathrm{SL}(2,\mathbb{Z})\) on \(H\) gives an action of
\(\mathrm{PSL}(2,\mathbb{Z})\), and the moduli space of elliptic curves
is \(H/\mathrm{PSL}(2,\mathbb{Z})\).

Now most points in \(H\) aren't preserved by any element of
\(\mathrm{PSL}(2,\mathbb{Z})\). However, certain points are! The point
\[\tau = i\] corresponding to the square lattice, is preserved by \(S\)
and all its powers. And the point \[\tau = \exp(2\pi i/3)\]
corresponding to the hexagonal lattice, is preserved by \(ST\) and all
its powers. These give rise to two conical singularities in the moduli
space of elliptic curves. Away from these points, the moduli space is
smooth.

Lest you get the wrong impression, I should hasten to reassure you that
the moduli space is not all that complicated: it looks almost like the
complex plane! There's a famous one-to-one and onto function from the
moduli space to the complex plane: it's called the ``modular function''
and denoted by \(j\). So the moduli space is \emph{topologically} just
like the complex plane; the only difference is that it fails to be
\emph{smooth} at two points, where there are conical singularities.

This may seem a bit hard to visualize, but it's actually not too hard.
Here's one way. Start with the region in the upper half-plane outside
the unit circle and between the vertical lines \(x = -1/2\) and
\(x = 1/2\). It looks sort of like this: \[
  \begin{tikzpicture}[scale=0.7]
    \draw[->] (-2.4,0) to (2.4,0) node[label=right:{$\Re(z)$}]{};
    \draw[->] (0,-1) to (0,5) node[label=left:{$\Im(z)$}]{};
    \node at (1,0) {\scriptsize$\bullet$};
    \node at (1,-0.6) {\scriptsize$B'$};
    \node at (-1,0) {\scriptsize$\bullet$};
    \node at (-1,-0.6) {\scriptsize$B$};
    \node at (0,1) {\scriptsize$\bullet$};
    \node at (-0.2,0.7) {\scriptsize$A$};
    \draw[thick] (1,0) arc(0:180:1);
    \draw[thick] (1,0) to (1,3.5);
    \draw[thick,dashed] (1,3.5) to (1,4.5);
    \draw[thick] (-1,0) to (-1,3.5);
    \draw[thick,dashed] (-1,3.5) to (-1,4.5);
    \draw (0.4,0.9) to (1,1.2);
    \foreach \y in {0.5,1,1.5,2,2.5}
      \draw (-1,{\y+0.2}) to (1,{\y+1.2});
    \draw (-1,3.2) to (0,3.7);
    \draw[dashed] (0,3.7) to (1,4.2);
    \draw[dashed] (-1,3.7) to (0.8,4.6);
    \draw[dashed] (-1,4.2) to (-0.2,4.6);
  \end{tikzpicture}
\] Then glue the vertical line starting at \(B\) to the one starting at
\(B'\), and glue the arc \(AB\) to the arc \(AB'\). We get a space
that's smooth everywhere except at the points \(A\) and \(B = B'\),
where there are conical singularities. The total angle around the point
\(A\) is just 180 degrees --- half what it would be if the moduli space
were smooth there. The total angle around \(B\) is just 120 degrees ---
one third what it would be if the moduli space were smooth there.

The reason this works is that the region shown above is a ``fundamental
domain'' for the action of \(\mathrm{PSL}(2,\mathbb{Z})\) on \(H\). In
other words, every elliptic curve is isomorphic to one where the
parameter \(\tau\) lies in this region. The point \(A\) is where
\(\tau = i\), and the point B is where \(\tau = exp(2\pi i/3)\).

Now let's see where the ``12-ness'' comes into this picture. Minhyong
Kim explained this to me in a very nice way, but to tell you what he
said, I'll have to turn up the level of mathematical sophistication
another notch. (Needless to say, all the errors will be mine.)

So, I'll assume you know what a ``complex line bundle'' is --- this is
just another name for a \(1\)-dimensional complex vector bundle. Locally
a section of a complex line bundle looks a lot like a complex-valued
function, but this isn't true globally unless your line bundle is
trivial. If you aren't careful, sometimes you may \emph{think} you have
a function defined on a space, only to discover later that it's actually
a section of a line bundle. This sort of thing happens all the time in
physics. In string theory, when you're doing path integrals on moduli
space, you have to make sure that what you're integrating is really a
function! So it's important to understand all the line bundles on moduli
space.

Now, given any sort of space, we can form the set of all isomorphism
classes of line bundles over this space. This is actually an abelian
group, since when we tensor two line bundles we get another line bundle,
and when you tensor any line bundle with its dual, you get the trivial
line bundle, which plays the role of the multiplicative identity for
tensor products. This group is called the ``Picard group'' of your
space.

What's the Picard group of the moduli space of elliptic curves? Well,
when I said ``any sort of space'' I was hinting that there are all sorts
of spaces --- topological spaces, smooth manifolds, algebraic varieties,
and so on --- each one of which comes with its own particular notion of
line bundle. Thus, before studying the Picard group of moduli space we
need to decide what context we're going to work in! As a mere
\emph{topological space}, we've seen that the moduli space of elliptic
curves is indistinguishable from the plane, and every \emph{topological}
line bundle over the plane is trivial, so in \emph{this} context the
Picard group is the trivial group --- boring!

But the moduli space is actually much more than a mere topological
space. It's not a smooth manifold, but it's awfully close: it's the
quotient of the smooth manifold \(H\) by the discrete group
\(\mathrm{SL}(2,\mathbb{Z})\), and its singularities are pretty mild in
nature.

Somehow we should take advantage of this when defining the Picard group
of the moduli space. One way to do so involves the theory of ``stacks''.
Without getting into the details of this theory, let me just vaguely
sketch what it does for us here. For a much more careful treatment, with
more of an algebraic geometry flavor, try:

\begin{enumerate}
\def\labelenumi{\arabic{enumi})}
\tightlist
\item
  David Mumford, ``Picard groups of moduli problems'', in
  \emph{Arithmetical Algebraic Geometry}, ed.~O. F. G. Schilling, Harper
  and Row, New York, 1965.
\end{enumerate}

Suppose a discrete group \(G\) acts on a smooth manifold \(X\). A
``\(G\)-equivariant'' line bundle on \(X\) is a line bundle equipped
with an action of \(G\) that gets along with the action of \(G\) on
\(X\). If \(G\) acts freely on \(X\), a line bundle on \(X/G\) is the
same as a \(G\)-equivariant line bundle on \(X\). This isn't true when
the action of \(G\) on \(X\) isn't free. But we can still go ahead and
\emph{define} the Picard group of \(X/G\) to be the group of isomorphism
classes of \(G\)-equivariant line bundles on \(X\). Of course we should
say something to let people know that we're using this funny definition.
In our example, people call it the Picard group of the moduli
\emph{stack} of elliptic curves.

So what's this group, anyway?

Well, it turns out that you can get any
\(\mathrm{SL}(2,\mathbb{Z})\)-equivariant line bundle on \(H\), up to
isomorphism, by taking the trivial line bundle on \(H\) and using a
\(1\)-dimensional representation of \(\mathrm{SL}(2,\mathbb{Z})\) to say
how it acts on the fiber. So we just need to understand
\(1\)-dimensional representations of \(\mathrm{SL}(2,\mathbb{Z})\). The
set of isomorphism classes of these forms a group under tensor product,
and this is the group we're after.

Well, a \(1\)-dimensional representation of a group always factors
through the abelianization of that group. We saw the abelianization of
\(\mathrm{SL}(2,\mathbb{Z})\) was \(\mathbb{Z}/12\). But everyone knows
that the group of \(1\)-dimensional representations of \(\mathbb{Z}/n\)
is again \(\mathbb{Z}/n\) - this is called Pontryagin duality. So: the
Picard group of the moduli stack of elliptic curves is
\(\mathbb{Z}/12\).

So we see again an inherent ``12-ness'' built into the theory of
elliptic curves! You may be wondering how this makes the number 24 so
important in string theory. In particular, where does that extra factor
of 2 come from? I'll say a little more about this next Week. I may or
may not manage to tie together the loose ends!

You may also be wondering about ``stacks''. In this you're not alone.
There's an amusing passage about stacks in the following book:

\begin{enumerate}
\def\labelenumi{\arabic{enumi})}
\setcounter{enumi}{1}
\tightlist
\item
  Joe Harris and Ian Morrison, \emph{Moduli of Curves}, Springer,
  Berlin, 1998.
\end{enumerate}

They write:

\begin{quote}
``Of course, here I'm working with the moduli stack rather than with the
moduli space. For those of you who aren't familiar with stacks, don't
worry: basically, all it means is that I'm allowed to pretend that the
moduli space is smooth and that there's a universal family over it.''

Who hasn't heard these words, or their equivalent, spoken in a talk? And
who hasn't fantasized about grabbing the speaker by the lapels and
shaking him until he says what --- exactly --- he means by them? But
perhaps you're now thinking that all that is in the past, and that at
long last you're going to learn what a stack is and what they do.

Fat chance.
\end{quote}

Actually Mumford's paper cited above gives a nice introduction to the
theory of stacks without mentioning the dreaded word ``stack''.
Alternatively, you can try this:

\begin{enumerate}
\def\labelenumi{\arabic{enumi})}
\setcounter{enumi}{2}
\tightlist
\item
The Stacks Projects, available at \url{https://stacks.math.columbia.edu/}
\end{enumerate}

But let me just briefly say a bit about stacks and the moduli stack of
elliptic curves in particular. A stack is a weak sheaf of categories.
For this to make sense you must already know what a sheaf is! In the
simplest case, a sheaf over a topological space, the sheaf \(S\) gives
you a set \(S(U)\) for each open set U, and gives you a function
\(S(U,V)\colon S(U)\to S(V)\) whenever the open set \(U\) is contained
in the open set \(V\). These functions must satisfy some laws. The
notion of ``stack'' is just a categorification of this idea. That is, a
stack \(S\) over a topological space gives you a \emph{category}
\(S(U)\) for each open set \(U\), and gives you a \emph{functor}
\(S(U,V)\colon S(U)\to S(V)\). These functors satisfy the same laws as
before, but \emph{only up to specified natural isomorphism}. And these
natural isomorphisms must in turn satisfy some new laws of their own,
so-called coherence laws.

In the case at hand there's a stack over the moduli space of elliptic
curves. For any open set \(U\) in the moduli space, an object of
\(S(U)\) is a family of elliptic curves over \(U\), such that each
elliptic curve in the family sits over the point in moduli space
corresponding to its isomorphism class. Similarly, a morphism in
\(S(U)\) is a family of isomorphisms of elliptic curves. This allows us
to keep track of the fact that some elliptic curves have more
automorphisms than others! And it takes care of the funny stuff that
happens at the singular points in the moduli space.

By the way, this watered-down summary leaves out a lot of the algebraic
geometry that you usually see when people talk about stacks.

Finally, one more thing --- it looks like Kreimer and company are making
great progress on understanding renormalization in a truly elegant way.

\begin{enumerate}
\def\labelenumi{\arabic{enumi})}
\setcounter{enumi}{3}
\tightlist
\item
  D. J. Broadhurst and D. Kreimer, ``Renormalization automated by Hopf
  algebra'', available as
  \href{https://arxiv.org/abs/hep-th/9810087}{\texttt{hep-th/9810087}}.
\end{enumerate}

Let me quote the abstract:

\begin{quote}
It was recently shown that the renormalization of quantum field theory
is organized by the Hopf algebra of decorated rooted trees, whose
coproduct identifies the divergences requiring subtraction and whose
antipode achieves this. We automate this process in a few lines of
recursive symbolic code, which deliver a finite renormalized expression
for any Feynman diagram. We thus verify a representation of the operator
product expansion, which generalizes Chen's lemma for iterated
integrals. The subset of diagrams whose forest structure entails a
unique primitive subdivergence provides a representation of the Hopf
algebra \(H_R\) of undecorated rooted trees. Our undecorated Hopf
algebra program is designed to process the 24,213,878 BPHZ contributions
to the renormalization of 7,813 diagrams, with up to 12 loops. We
consider 10 models, each in 9 renormalization schemes. The two simplest
models reveal a notable feature of the subalgebra of Connes and
Moscovici, corresponding to the commutative part of the Hopf algebra
\(H_T\) of the diffeomorphism group: it assigns to Feynman diagrams
those weights which remove \(\zeta\) values from the counterterms of the
minimal subtraction scheme. We devise a fast algorithm for these
weights, whose squares are summed with a permutation factor, to give
rational counterterms.
\end{quote}

\hypertarget{week126}{%
\section{November 17, 1998}\label{week126}}

To round off some things I said in the previous two weeks, let me say a
bit more about string theory and Euler's mysterious equation
\[1 + 2 + 3 + \cdots = -\frac{1}{12}.\] 
For this I'll need to assume a
nodding acquaintance with quantum field theory.

There are two complementary ways to attack almost any problem in quantum
field theory: the Lagrangian approach, also known as ``path-integral
quantization'', and the Hamiltonian approach, also called ``canonical
quantization''. Let me describe string theory from both viewpoints. I'll
only talk about bosonic string theory, because my goal is to sketch why
it works best in \(26\)-dimensional spacetime, and because it's simpler
than superstring theory. Also, I'll only talk about closed strings.

Classically, such a string is simply a map from a closed surface into
spacetime. In the Lagrangian approach to quantization, we start by
choosing a formula for the action. We use the simplest possibility,
namely the \emph{area} of the surface. Of course, to define the area of
a surface in spacetime, we need the spacetime to have a metric. The
simplest thing is to work with \(n\)-dimensional Minkowski spacetime, so
let's do that.

We find the equations of motion of the string by extremizing the action.
These equations imply that if we watch the string in space as time
passes, it acts like collection of loops made of perfectly elastic
material. These loops vibrate, split and join as time passes.

It's perhaps a bit easier to see how the strings vibrate if we go over
to the Hamiltonian approach. This is a bit subtle, because string theory
has an enormous amount of ``gauge symmetry'' --- by which physicists
mean any symmetry that arises from the ability to switch between
different mathematical descriptions of what counts as the same physical
situation. There's a recipe to figure out the gauge symmetries of any
theory starting from the action. Applying this to string theory, it
turns out that two maps from a surface into spacetime count as
``physically the same'' if they differ only by a reparametrization of
the surface that's being mapped into spacetime.

When going over to the Hamiltonian approach, we have to deal with this
gauge symmetry. There are different ways to deal with it --- but we
can't just ignore it. Suppose we use the approach called ``lightcone
gauge-fixing''. This amounts to choosing a parametrization of our
surface so that the 2 coordinates on it are related in a simple way to 2
of the coordinates on \(n\)-dimensional Minkowski space. We can do this
because of the reparametrization gauge symmetry. But once we've done it,
we no longer have any more freedom to reparametrize our surface. In
short, we've squeezed all the juice out of our gauge symmetry: this is
what ``gauge-fixing'' is all about.

We started by studying a map from a surface \(S\) into \(n\)-dimensional
spacetime, which we can think of as field on \(S\) with \(n\)
components. However, in lightcone gauge, 2 components of this field are
given by simple formulas in terms of the rest. This lets us think of our
string as a field \(X\) with only \(n-2\) components. And when we do
this, it satisfies the simplest equation you could imagine! Namely, the
wave equation
\[\left(\frac{d^2}{dt^2}-\frac{d^2}{dx^2}\right) X(t,x) = 0.\] This is
same equation that describes an idealized violin string. The only
difference is that now, instead of a segment of violin string, we have a
bunch of closed loops of string. The energy, or Hamiltonian, is also
given by the usual wave equation Hamiltonian:
\[H = \frac12\int\left[\left(\frac{dX}{dt}\right)^2+\left(\frac{dX}{dx}\right)^2\right]dx.\]
The first term represents the kinetic energy of the string, while the
second represents its potential energy --- the energy it has due to
being stretched.

Henceforth I'll ignore the fact that loops of string can split or join,
and only talk about the vibrations of a single loop of string. Using the
linearity of the wave equation, we can decompose any solution of the
wave equation into sine waves moving in either direction --- so-called
``left-movers'' and ``right-movers'' --- together with a solution of the
form \[X(t,x) = A + Bt\] which describes the motion of the string's
center of mass. The left-movers and right-movers don't interact with
each other or with the center-of-mass motion, so we can learn a lot just
by studying the right-movers.

For starters, suppose the field \(X\) has just one component. Then the
right-moving vibrational modes look like
\[X(t,x) = A\sin(ik(t-x)) + B\cos(ik(t-x))\] with frequencies
\(k = 1,2,3,\ldots\). Abstractly, each of these vibrational modes is
just like a harmonic oscillator of frequency \(k\), so we can think of
the string as a big collection of harmonic oscillators.

Now suppose we quantize our string --- or more precisely, the
right-moving modes. By what we've said, this just amounts to quantizing
a bunch of harmonic oscillators, one of frequency \(k\) for each natural
number \(k\). This is great, since the harmonic oscillator is one of the
easiest physical systems to quantize!

As you may know, the quantum harmonic oscillator has discrete energy
levels with energies \(k/2, 3k/2, 5k/2,\ldots\). (Here I'm working in
units where \(\hbar = 1\); otherwise I'd need a factor of \(\hbar\).) In
particular, the energy of the lowest-energy state is called the
``zero-point energy'' or ``vacuum energy''. It usually doesn't hurt much
to subtract this off by redefining the Hamiltonian, but sometimes it's
important.

Now, what's the total zero-point energy of all the right-moving modes?
To figure this out, we add up the zero-point energy \(k/2\) for all
frequencies \(k = 1,2,3,\ldots\), obtaining
\[\frac{1 + 2 + 3 + \ldots }{2}.\] Of course this is divergent, but
there are lots of sneaky tricks for assigning values to divergent
series, so let's not be disheartened! Euler figured out such a trick for
calculating the sum \(1 + 2 + 3 +\ldots\) , and he got the value
\(-1/12\). If we momentarily assume this makes sense, then the total
zero-point energy works out to be \(-\frac{1}{24}\) More generally, if
we have a string in \(n\)-dimensional Minkowski spacetime, the field
\(X\) has \(n-2\) components, so the total zero-point energy is
\[-\frac{n-2}{24}\] Now, for other reasons, it turns out that string
theory works best when this zero-point energy is \(-1\). This is a bit
tricky to explain, but it has to do with the subtleties of gauge-fixing
in quantum field theory. Things that work nicely at the classical level
can easily screw up at the quantum level; in particular, symmetries of a
classical theory can be lost when you quantize. One has to really check
that the light-cone gauge fixing doesn't screw up the Lorentz-invariance
of string theory. It turns out that it \emph{does} screw it up unless
the zero-point energy of the right-movers is \(-1\). So bosonic string
theory works best when \[\frac{n-2}{24} = 1\] or in other words, when
\(n = 26\).

You really shouldn't take my word for this stuff! You can find more
details around pages 95--96 in volume 1 of the following book:

\begin{enumerate}
\def\labelenumi{\arabic{enumi})}
\tightlist
\item
  Michael B. Green, John H. Schwarz and Edward Witten, \emph{Superstring
  Theory}, two volumes, Cambridge U.\ Press, Cambridge, 1987.
\end{enumerate}

There's a lot I should say to fill in the details, but the most urgent
matter is to explain Euler's mysterious formula
\[1 + 2 + 3 + \ldots = -\frac{1}{12}\] As I said in
\protect\hyperlink{week124}{``Week 124''}, this is an example of zeta
function regularization. The Riemann zeta function is defined by
\[\zeta(s) = \frac{1}{1^s} + \frac{1}{2^s} + \frac{1}{3^s} + \ldots\]
when the sum converges, but it analytically continues to values of s
where the sum doesn't converge. If we do the analytic continuation, we
get \[\zeta(-1) = -\frac{1}{12}.\] Proving this rigorously is a bit of
work. One way is to use the ``functional equation'' for the Riemann zeta
function, which says that \[F(s) = F(1-s)\] where
\[F(s) = \pi^{-\frac{s}{2}}\Gamma\left(\frac{s}{2}\right)\zeta(s)\] and
\(\Gamma\) is the famous function with \(\Gamma(n) = (n-1)!\) for
\(n = 1,2,3,\ldots\) and \(\Gamma(s+1) = s \Gamma(s)\) for all \(s\).
Using \[\Gamma\left(\frac12\right) = \sqrt{\pi}\] and
\[\zeta(2) = \frac{\pi^2}{6},\] the functional equation implies
\(\zeta(-1) = -1/12\). But of course you have to prove the functional
equation! A nice exposition of this can be found in:

\begin{enumerate}
\def\labelenumi{\arabic{enumi})}
\setcounter{enumi}{1}
\tightlist
\item
  Neal Koblitz, \emph{Introduction to Elliptic Curves and Modular
  Forms}, 2nd edition, Springer, Berlin, 1993.
\end{enumerate}

I don't know Euler's original argument that \(\zeta(-1) = -1/12\).
However, Dan Piponi recently gave the following ``physicist's proof'' on
the newsgroup \texttt{sci.physics.research}. Let \(D\) be the
differentiation operator: \[D = \frac{d}{dx}\] Then Taylor's formula
says that translating a function to the left by a distance \(c\) is the
same as applying the operator \(e^{cD}\) to it, since
\[e^{cD} f = f + cf' + \left(\frac{c^2}{2!}\right)f'' + \ldots\] Using
some formal manipulations we obtain \[
  \begin{aligned}
    f(0) + f(1) + f(2) + \ldots
    & = ((1+e^D+e^{2D}+\ldots)f)(0)
  \\& = \left(\left(\frac{1}{1-e^D}\right)f\right)(0)
  \end{aligned}
\] or if \(F\) is an integral of \(f\), so that \(DF = f\),
\[f(0) + f(1) + f(2) + \ldots = \left(\left(\frac{D}{1 - e^D}\right) F\right)(0)\]
This formula can be made rigorous in certain contexts, but now we'll
throw rigor to the winds and apply it to the function \(f(x) = x\),
obtaining
\[0+1+2+3+\ldots = \left(\left(\frac{D}{1 - e^D}\right) F\right)(0)\]
where \[F(x) = \frac{x^2}{2}\] To finish the job, we work out the
beginning of the Taylor series for \(D/(1 - e^D)\). The coefficients of
this are closely related to the Bernoulli numbers, and this could easily
lead us into further interesting digressions, but all we need to know is
\[\frac{D}{1 - e^D} = -1 + \frac{D}{2} - \frac{D^2}{12} + \ldots\]
Applying this operator to \(F(x) = x^2/2\) and evaluating the result at
\(x=0\), the only nonzero term comes from the \(D^2\) term in the power
series, so we get
\[1 + 2 + 3 + \ldots = \left(\left(-\frac{D^2}{12}\right)F\right)(0) = -\frac{1}{12}\]
Voilà!

By the way, after he came up with this proof, Dan Piponi found an almost
identical proof in the following book:

\begin{enumerate}
\def\labelenumi{\arabic{enumi})}
\setcounter{enumi}{2}
\tightlist
\item
  G. H. Hardy, \emph{Divergent Series}, Chelsea Pub. Co., New York,
  1991.
\end{enumerate}

Now let me change gears. Besides the Riemann zeta function, there are a
lot of other special functions that show up in the study of elliptic
curves. Unsurprisingly, many of them are also important in string
theory. For example, consider the partition function of bosonic string
theory. What do I mean by a ``partition function'' here? Well, whenever
we have a quantum system with a Hamiltonian \(H\), its partition
function is defined to be
\[Z(\beta) = \operatorname{tr}(\exp(-\beta H))\] where \(\beta > 0\) is
the inverse temperature. This function is fundamental to statistical
mechanics, for reasons that I'm too lazy to explain here.

Before tackling the bosonic string, let's work out the partition
function for a quantum harmonic oscillator. To keep life simple, let's
subtract off the zero-point energy so the energy levels are \(0\),
\(k\), \(2k\), and so on. Mathematically, these energy levels are just
the eigenvalues of the harmonic oscillator Hamiltonian, \(H\). Thus the
eigenvalues of \(\exp(-\beta H)\) are \(1\), \(\exp(-\beta k)\),
\(\exp(-2\beta k)\), etc. The trace of this operator is just the sum of
its eigenvalues, so we get
\[Z(\beta) = 1 + \exp(-\beta k) + \exp(-2\beta k) + \ldots = \frac{1}{1 - \exp(-\beta k)}\]
This was first worked out by Planck, who assumed the harmonic oscillator
had discrete, evenly spaced energy levels and computed its partition
function as part of his struggle to understand the thermodynamics of the
electromagnetic field.

Okay, now let's do the bosonic string. To keep life simple we again
subtract off the zero-point energy. Also, we'll consider only the
right-moving modes, and we'll start by assuming the field \(X\)
describing the vibrations of the string has only one component. As we
saw before, the string then becomes the same as a collection of quantum
harmonic oscillators with frequencies \(k = 1, 2, 3, \ldots\). We've
seen that the oscillator with frequency \(k\) has partition function
\(1/(1 - \exp(-\beta k))\). To get the partition function of a quantum
system built from a bunch of noninteracting parts, you multiply the
partition functions of the parts (since the trace of a tensor product of
operators is the product of their traces). So the partition function of
our string is \[\prod_{k\in\mathbb{N}} \frac{1}{1 -\exp(-\beta k)}\] So
far, so good. But now suppose we take the zero-point energy into
account. We do this by subtracting \(1/24\) from the Hamiltonian of the
string, which has the effect of multiplying its partition function by
\(\exp(\beta/24)\). Thus we get
\[Z(\beta) = \exp(\beta/24) \prod \frac{1}{1 - \exp(-\beta k)}\] Lo and
behold: the reciprocal of the Dedekind eta function!

What's that, you ask? It's a very important function in the theory of
elliptic curves. People usually write it as a function of
\(q = \exp(-\beta)\), like this:
\[\eta(q) = q^{\frac{1}{24}} \prod (1 - q^k)\] But to see the relation
to elliptic curves we should switch variables yet again and write
\(q = \exp(2\pi i\tau)\). I already talked about this variable \(\tau\)
in \protect\hyperlink{week125}{``Week 125''}, where we were studying the
elliptic curve formed by curling up a parallelogram like this in the
complex plane: \[
  \begin{tikzpicture}[scale=0.7]
    \draw[->] (-2.5,0) to (4,0) node[label=below:{$\Re(z)$}]{};
    \draw[->] (0,-2) to (0,4) node[label=left:{$\Im(z)$}]{};
    \foreach \m in {-1,0,1,2}
    {
      \foreach \n in {-1,0,1,2}
      {
        \node (\m\n) at ({1.5*\m+0.75*\n},{1.33*\n}) {$\circ$};
      }
    }
    \node at (00) {$\bullet$};
    \node at (-0.3,-0.4) {$0$};
    \node at (10) {$\bullet$};
    \node at (1.5,-0.4) {$1$};
    \node at (01) {$\bullet$};
    \node at (0.75,0.9) {$\tau$};
    \node at (11) {$\bullet$};
    \node at (2.25,0.9) {$\tau+1$};
  \end{tikzpicture}
\] In physics, this elliptic curve is just one possibility for the shape
of a surface traced out by a string. The number \(1\) says how far the
surface goes in the \emph{space} direction before it loops around, and
the number \(\tau\) says how far it goes in the \emph{time} direction
before it loops around!

(The idea of ``looping around in time'' may seem bizarre, but it's very
important in physics. It turns out that studying the statistical
mechanics of a system at a given inverse temperature is the same as
studying Euclidean quantum field theory on a spacetime where time is
periodic with a given period. This idea is what relates the variables
\(\beta\) and \(\tau\).)

Now as I explained in \protect\hyperlink{week13}{``Week 13''}, the above
elliptic curve is not just an abstract torus-shaped thingie. We can also
think of it as the set of complex solutions of the following cubic
equation in two variables: \[y^2 = 4x^3 - g_2 x - g_3\] where the
numbers \(g_2\) and \(g_3\) are certain functions of \(\tau\). Moreover,
this equation defines an elliptic curve whenever the polynomial on the
right-hand side doesn't have repeated roots. So among other things,
elliptic curves are really just a way of studying cubic equations!

But when does \(4x^3 - g_2 x - g_3\) have repeated roots? Precisely when
the ``discriminant'' \[\Delta = g_2^3 - 27 g_3^2\] equals zero. This is
just the analog for cubics of the more familiar discriminant for
quadratic equations.

Now for the cool part: there's an explicit formula for the discriminant
in terms of the variable \(\tau\). And it involves the 24th power of the
Dedekind eta function! Here it is: \[\Delta = (2 \pi)^{12}\eta^{24}\] If
you haven't seen this before, it should seem \emph{amazing} that the
discriminant of a cubic equation can be computed using the 24th power of
a partition function that shows up in string theory. Of course that's
not how it went historically: Dedekind discovered his eta function long
before strings came along. What's really happening is that string theory
is exploiting special features of complex curves, and thus acquires some
of the ``24-ness'' of elliptic curves.

If I have the energy, next time I'll give you a better explanation of
why bosonic string theory works best in 26 dimensions, using some
special properties of the Dedekind eta function.

Meanwhile, if you want to see pictures of the Dedekind eta function,
together with some cool formulas it satisfies, try these:

\begin{enumerate}
\def\labelenumi{\arabic{enumi})}
\setcounter{enumi}{3}
\item
  Mathworld, ``Dedekind eta function'',
  \href{http://mathworld.wolfram.com/DedekindEtaFunction.html}{\texttt{{http://mathworld.wolfram.com/Dedekind}}}
  \href{http://mathworld.wolfram.com/DedekindEtaFunction.html}{\texttt{{EtaFunction.html}}}
\item
  Wikipedia, ``Dedekind eta function'',
  \href{http://en.wikipedia.org/wiki/Dedekind_eta_function}{\texttt{{http://en.wikipedia.org/wiki/Dedekind\_}}}
  \href{http://en.wikipedia.org/wiki/Dedekind_eta_function}{\texttt{{eta\_function}}}
\end{enumerate}

\begin{center}\rule{0.5\linewidth}{0.5pt}\end{center}

\begin{quote}
"Dear Sir,

I am very much gratified on perusing your letter of the 8th February
1913. I was expecting a reply from you similar to the one which a
Mathematics Professor at London wrote asking me to study carefully
Bromwich's Infinite Series and not fall into the pitfall of divergent
series. I have found a friend in you who views my labors
sympathetically. This is already some encouragement to me to proceed
with an onward course. I find in many a place in your letter rigourous
proofs are required and so on and you ask me to communicate the method
of proof. If I had given you my methods of proof I am sure you will
follow the London Professor. But as a fact I did not give him any proof
but made some assertions as the following under my new theory. I told
him that the sum of an infinite number of terms in the series
\(1+2+3+4+\ldots=-1/12\) under my theory. If I tell you this you will at
once point out to me the lunatic asylum as my goal."

--- Srinivasa Ramanujan, second letter to G. H. Hardy
\end{quote}

\hypertarget{week127}{%
\section{November 30, 1998}\label{week127}}

If you like \(\pi\), take a look at this book:

\begin{enumerate}
\def\labelenumi{\arabic{enumi})}
\tightlist
\item
  Lennart Berggren, Jonathan Borwein and Peter Borwein, \emph{\(\pi\): A
  Source Book}, Springer, Berlin, 1997.
\end{enumerate}

It's full of reprints of original papers about \(\pi\), from the Rhind
Papyrus right on up to the 1996 paper by Bailey, Borwein and Plouffe in
which they figured out how to compute a given hexadecimal digit of
\(\pi\) without computing the previous digits --- see
\protect\hyperlink{week66}{``Week 66''} for more about that. By the way,
Colin Percival has recently used this technique to compute the 5
trillionth binary digit of \(\pi\)! (It's either zero or one, I forget
which.) Percival is a 17-year old math major at Simon Fraser University,
and now he's leading a distributed computation project to calculate the
quadrillionth binary digit of \(\pi\). Anyone with a Pentium or faster
computer using Windows 95, 98, or NT can join. For more information,
see:

\begin{enumerate}
\def\labelenumi{\arabic{enumi})}
\setcounter{enumi}{1}
\tightlist
\item
  PiHex project,
  \url{http://wayback.cecm.sfu.ca/projects/pihex/}
\end{enumerate}

Anyway, the above book is \emph{full} of fun stuff, like a one-page
proof of the irrationality of \(\pi\) which uses only elementary
calculus, due to Niven, and the following weirdly beautiful formula due
to Euler, which unfortunately is not explained:
\[\frac{\pi}{2} = \frac32\cdot\frac56\cdot\frac76\cdot\frac{11}{10}\cdot\frac{13}{14}\cdot\frac{17}{18}\cdot\frac{19}{18}\cdot\ldots\]
Here the numerators are the odd primes, and the denominators are the
closest numbers of the form \(4n+2\).

Since I've been learning about elliptic curves and the like lately, I
was also interested to see a lot of relations between \(\pi\) and
modular functions. For example, the book has a reprint of Ramanujan's
paper ``Modular equations and approximations to \(\pi\)'', in which he
derives a bunch of bizarre formulas for \(\pi\), some exact but others
approximate, like this:
\[\pi \sim \frac{12}{\sqrt{190}}\ln((2\sqrt{2}+\sqrt{10})(3+\sqrt{10}))\]
which is good to 18 decimal places. The strange uses to which genius can
be put know no bounds!

Okay, now I'd like to wrap up my little story about why bosonic string
theory works best in 26 dimensions. This time I want to explain how you
do the path integral in string theory. Most of what I'm about to say
comes from some papers that my friend Minhyong Kim recommended to me
when I started pestering him about this stuff:

\begin{enumerate}
\def\labelenumi{\arabic{enumi})}
\setcounter{enumi}{2}
\item
  Jean-Benoit Bost, ``Fibres determinants, determinants regularises, et
  mesures sur les espaces de modules des courbes complexes'',
  \emph{Asterisque} \textbf{152}--\textbf{153} (1987), 113--149.
\item
  A. A. Beilinson and Y. I. Manin, ``The Mumford form and the Polyakov
  measure in string theory'', \emph{Comm. Math. Phys.} \textbf{107}
  (1986), 359--376.
\end{enumerate}

For a more elementary approach, try chapters IX and X.4 in this book:

\begin{enumerate}
\def\labelenumi{\arabic{enumi})}
\setcounter{enumi}{4}
\tightlist
\item
  Charles Nash, \emph{Differential Topology and Quantum Field Theory},
  Academic Press, New York, 1991.
\end{enumerate}

As I explained in \protect\hyperlink{week126}{``Week 126''}, a string
traces out a surface in spacetime called the ``string worldsheet''.
Let's keep life simple and assume the string worldsheet is a torus and
that spacetime is Euclidean \(\mathbb{R}^n\). Then to figure out the
expectation value of any physical observable, we need to calculate its
integral over the space of all maps from a torus to \(\mathbb{R}^n\).

To make this precise, let's use \(X\) to denote a map from the torus to
\(\mathbb{R}^n\). Then a physical observable will be some function
\(f(X)\), and its expectation value will be
\[\frac{1}{Z} \int f(X) \exp(-S(X)) dX\] Here \(S(X)\) is the action for
string theory, which is just the \emph{area} of the string worldsheet.
The number \(Z\) is a normalizing factor called the partition function:
\[Z = \int \exp(-S(X)) dX\] But there's a big problem here! As usual in
quantum field thoery, the space we're trying to integrate over is
infinite-dimensional, so the above integrals have no obvious meaning.
Technically speaking, the problem is that there's no Lebesgue measure
``\(dX\)'' on an infinite-dimensional manifold.

Mathematicians might throw up our hands in despair and give up at this
point. But physicists take a more pragmatic attitude: they just keep
massaging the problem, breaking rules here and there if necessary, until
they get something manageable. Physicists call this ``calculating the
path integral'', but from a certain viewpoint what they're really doing
is \emph{defining} the path integral, since it only has a precise
meaning after they're done.

In the case at hand, it was Polyakov who figured out the right massage:

\begin{enumerate}
\def\labelenumi{\arabic{enumi})}
\setcounter{enumi}{5}
\tightlist
\item
  A. M. Polyakov, ``Quantum geometry of bosonic strings'', \emph{Phys.
  Lett.} \textbf{B103} (1981), 207.
\end{enumerate}

He rewrote the above integral as a double integral: first an integral
over the space of metrics \(g\) on the torus, and then inside, for each
metric, an integral over maps \(X\) from the torus into spacetime. Of
course, any such map gives a metric on the torus, so this double
integral is sort of redundant. However, introducing this redundancy
turns out not to hurt. In fact, it helps!

To keep life simple, let's just talk about the partition function
\[Z = \int \exp(-S(X)) dX\] If we can handle this, surely we can handle
the integral with the observable \(f(X)\) in it. Polyakov's trick turns
the partition function into a double integral:
\[Z = \int ( \int \exp(-\langle X, \Delta X \rangle) dX) dg\] where
\(\Delta\) is the Laplacian on the torus and the angle brackets stand
for the usual inner product of \(\mathbb{R}^n\)-valued functions, both
defined using the metric \(g\).

At first glance Polyakov's trick may seem like a step backwards: now we
have two ill-defined integrals instead of one! However, it's actually a
step forward. Now we can do the inside integral by copying the formula
for the integral of a Gaussian of finitely many variables --- a standard
trick in quantum field theory. Ignoring an infinite constant that would
cancel later anyway, the inside integral works out to be:
\[(\operatorname{det}\Delta)^{-\frac12}\] But wait! The Laplacian here
is a linear operator on the vector space of \(\mathbb{R}^n\)-valued
functions on the torus. This is an infinite-dimensional vector space, so
we can't blithely talk about determinants the way we can in finite
dimensions. In finite dimensions, the determinant of a self-adjoint
matrix is the product of its eigenvalues. But the Laplacian has
infinitely many eigenvalues, which keep getting bigger and bigger. How
do we define the product of all its eigenvalues?

Of course the lowest eigenvalue of the Laplacian is zero, and you might
think that would settle it: the product of them all must be zero! But
that would make the above expression meaningless, and we are not going
to give up so easily. Instead, we will simply ignore the zero
eigenvalue! That way, we only have to face the product of all the
\emph{rest}.

(Actually there's something one can do which is slightly more careful
than simply ignoring the zero eigenvalue, but I'll talk about that
later.)

Okay, so now we have a divergent product to deal with. Well, in
\protect\hyperlink{week126}{``Week 126''} I used a trick called zeta
function regularization to make sense of a divergent sum, and we can use
that trick here to make sense of our divergent product. Suppose we have
a self-adjoint operator \(A\) with a discrete spectrum consisting of
positive eigenvalues. Then the ``zeta function'' of \(A\) is defined by:
\[\zeta(s) = \operatorname{tr}(A^{-s})\] To compute \(\zeta(s)\) we just
take all the eigenvalues of \(A\), raise them to the \(-s\) power, and
add them up. For example, if \(A\) has eigenvalues \(1,2,3,\ldots\),
then \(\zeta(s)\) is just the usual Riemann zeta function, which we
already talked about in \protect\hyperlink{week126}{``Week 126''}.

This stuff doesn't quite apply if \(A\) is the Laplacian on a compact
manifold, since one of its eigenvalues is zero. But we have already
agreed to throw out the zero eigenvalue, so let's do that when defining
\(\zeta(s)\) as a sum over eigenvalues. Then it turns out that the sum
converges when the real part of \(s\) is positive and large. Even
better, there's a theorem saying that \(\zeta(s)\) can be analytically
continued to \(s = 0\). Thus we can use the following trick to define
the determinant of the Laplacian.

Suppose that \(A\) is a self-adjoint matrix with positive eigenvalues.
Then \[\zeta(s) = \operatorname{tr}(\exp(-s \ln A))\] Differentiating
gives \[\zeta'(s) = -\operatorname{tr}(\ln A \exp(-s \ln A))\] and
setting \(s\) to zero we get \[\zeta'(0) = -\operatorname{tr}(\ln A).\]

But there's a nice little formula saying that
\[\operatorname{det}(A) = \exp(\operatorname{tr}(\ln A))\] so we get
\[\operatorname{det}(A) = \exp(-\zeta'(0)).\] Now we can use this
formula to \emph{define} the determinant of the Laplacian on a compact
manifold! Sometimes this is called a ``regularized determinant''.

So --- where are we? We used Polyakov's trick to write the partition
function of our torus-shaped string as
\[Z = \int ( \int \exp(-\langle X, \Delta X\rangle) dX) dg,\] then we
did the inside integral and got this:
\[Z = \int (\operatorname{det}\Delta)^{-\frac12} dg\] and then we
figured out a meaning for the determinant here.

What next? Well, since the Laplacian on \(\mathbb{R}^n\)-valued
functions is the direct sum of \(n\) copies of the Laplacian on
real-valued functions, we expect that
\[(\operatorname{det}\Delta)^{-\frac12} = (\operatorname{det}\mathrm{laplacian})^{-\frac{n}{2}}\]
where ``\(\mathrm{laplacian}\)'' stands for the Laplacian on ordinary
real-valued functions on the torus. One can actually check this
rigorously using the definition in terms of zeta functions. That's
reassuring: at least \emph{one} step of our calculation is rigorous! So
we get
\[Z = \int (\operatorname{det}\mathrm{laplacian})^{-\frac{n}{2}} dg\]
Great. But we are not out of the woods yet. We still have an integral
over the space of metrics to do --- another nasty infinite-dimensional
integral.

Time for another massage!

Look at this formula again:
\[Z = \int ( \int \exp(-\langle X, \Delta X\rangle) dX) dg\] The
Laplacian depends on the metric \(g\), and so does the inner product.
However, the combination \(\langle X, \Delta X\rangle\) depends only on
the ``conformal structure'' --- i.e., it doesn't change if we multiply
the metric by a position-dependent scale factor. It also doesn't change
under diffeomorphisms.

Now the space of conformal structures on a torus modulo diffeomorphisms
is something familiar: it's just the moduli space of elliptic curves! We
figured out what this space looks like in
\protect\hyperlink{week125}{``Week 125''}. It's finite-dimensional and
there's a nice way to integrate over it, called the Weil-Petersson
measure. So we can hope to replace the outside integral --- the integral
over metrics --- by an integral over this moduli space.

Indeed, we could hope that
\[Z \overset{\text{we hope!}}{=} \int (\int \exp(-\langle X, \Delta X\rangle) dX) d[g]\]
where now the outside integral is over moduli space and \(d[g]\) is the
Weil-Petersson measure. The hope, of course, is that the stuff on the
inside is well-defined as a function on moduli space.

Actually this hope is a bit naive. Even though
\(\langle X, \Delta X\rangle\) doesn't change if we rescale the metric,
the whole inside integral \[\int \exp(-\langle X, \Delta X\rangle) dX\]
\emph{does} change. This may seem odd, but remember, we did a lot of
hair-raising manipulations before we even got this integral to mean
anything! We basically wound up \emph{defining} it to be
\[(\operatorname{det}\Delta)^{-\frac12},\] and one can check that this
\emph{does} change when we rescale the metric. This problem is called
the ``conformal anomaly''.

Are we stuck? No! Luckily, there is \emph{another} problem, which
cancels this one when \(n = 26\). They say two wrongs don't make a
right, but with anomalies that's often the only way to get things to
work\ldots.

So what's this other problem? It's that we shouldn't just replace the
measure \(dg\) by the measure \(d[g]\) as I did in my naive formula for
the partition function. We need to actually figure out the relation
between them. Of course this is hard to do, because the measure \(dg\)
doesn't exist, rigorously speaking! Still, if we do a bunch more
hair-raising heuristic manipulations, which I will spare you, we can get
a formula relating \(dg\) and and \(d[g]\), and using this we get
\[Z = \int ( \int \exp(-\langle X, \Delta X\rangle) dX) f(g) d[g]\]
where \(f(g)\) is some function of the metric. There's a perfectly
explicit formula for this function, but your eyeballs would fall out if
I showed it to you. Anyway, the real point is that \emph{in 26
dimensions and only in 26 dimensions}, the integrand
\[( \int \exp(-\langle X, \Delta X\rangle) dX) f(g)\] is invariant under
rescalings of the metric (as well as diffeomorphisms). In other words,
the conformal anomaly in \[\int \exp(-\langle X, \Delta X\rangle) dX\]
is precisely canceled by a similar conformal anomaly in \(f(g)\), so
their product is a well-defined function on moduli space, so it makes
sense to integrate it against \(d[g]\). We can then go ahead and figure
out the partition function quite explicitly.

By now, if you're a rigorous sort of pure mathematician, you must be
suffering from grave doubts about the sanity of this whole procedure.
But physicists regard this chain of reasoning, especially the miraculous
cancellation of anomalies at the end, as a real triumph. And indeed,
it's far \emph{better} than \emph{most} of what happens in quantum field
theory!

I've heard publishers of science popularizations say that each equation
in a book diminishes its readership by a factor of \(2\). I don't know
if this applies to This Week's Finds, but normally I try very hard to
keep the equations to a minimum. This week, however, I've been very bad,
and if my calculations are correct, by this point I am the only one
reading this. So I might as well drop all pretenses of expository prose
and write in a way that only I can follow! The real reason I'm writing
this, after all, is to see if I understand this stuff.

Okay, so now I'd like to see if I understand how one explicitly
calculates this integral:
\[\int \exp(-\langle X, \mathrm{laplacian} X\rangle) dX\] Since we're
eventually going to integrate this (times some stuff) over moduli space,
we might as well assume the metric on our torus is gotten by curling up
the following parallelogram in the complex plane: \[
  \begin{tikzpicture}[scale=0.7]
    \draw[->] (-3,0) to (4,0) node[label=below:{$\Re(z)$}]{};
    \draw[->] (0,-2) to (0,4) node[label=left:{$\Im(z)$}]{};
    \foreach \m in {-1,0,1,2}
    {
      \foreach \n in {-1,0,1,2}
      {
        \node (\m\n) at ({1.5*\m+0.75*\n},{1.33*\n}) {$\circ$};
      }
    }
    \node at (00) {$\bullet$};
    \node at (-0.3,-0.4) {$0$};
    \node at (10) {$\bullet$};
    \node at (1.5,-0.4) {$1$};
    \node at (01) {$\bullet$};
    \node at (0.75,0.9) {$\tau$};
    \node at (11) {$\bullet$};
    \node at (2.25,0.9) {$\tau+1$};
  \end{tikzpicture}
\] There are at least two ways to do the calculation. One is to actually
work out the eigenvalues of the Laplacian on this torus and then do the
zeta function regularization to compute its determinant. Di Francesco,
Mathieu, and Senechal do this in the textbook I talked about in
\protect\hyperlink{week124}{``Week 124''}. They get
\[\int \exp(-\langle X, \mathrm{laplacian} X\rangle) dX = \frac{1}{\sqrt{\Im(\tau)} |\eta(\tau)|^2}\]
where ``\(\eta\)'' is the Dedekind eta function, regarded as function of
\(\tau\). But the calculation is pretty brutal, and it seems to me that
there should be a much easier way to get the answer. The left-hand side
is just the partition function for an massless scalar field on the
torus, and we basically did that back in
\protect\hyperlink{week126}{``Week 126''}. More precisely, we considered
just the right-moving modes and we got the following partition function:
\[\frac{1}{\eta(\tau)}\] How about the left-moving modes? Well, I'd
guess that their partition function is just the complex conjugate,
\[\frac{1}{\eta(\tau)^*}\] since right-movers correspond to holomorphic
functions and left-movers correspond to antiholomorphic functions in
this Euclidean picture. It's just a guess! And finally, what about the
zero-frequency mode? I have no idea. But we should presumably multiply
all three partition functions together to get the partition function of
the whole system --- that's how it usually works. And as you can see, we
\emph{almost} get the answer that Di Francesco, Mathieu, and Senechal
got. It would work out \emph{perfectly} if the partition function of the
zero-frequency mode were \(1/\sqrt{\Im(\tau)}\). By the way,
\(\Im(\tau)\) is just the \emph{area} of the torus.

As evidence that something like this might work, consider this: the
zero-frequency mode is presumably related to the zero eigenvalue of the
Laplacian. We threw that out when we defined the regularized determinant
of the Laplacian, but as I hinted, more careful calculations of
\[\int \exp(-\langle X, \mathrm{laplacian} X\rangle) dX\] don't just
ignore the zero eigenvalue. Instead, they somehow use it to get an extra
factor of \(1/\sqrt{\Im(\tau)}\). Admittedly, the calculations are not
particularly convincing: a more obvious guess would be that it gives a
factor of infinity. Di Francesco, Mathieu, and Senechal practically
admit that they \emph{need} this factor just to get modular invariance,
and that they'll do whatever it takes to get it. Nash just sticks in the
factor of \(1/\sqrt{\Im(\tau)}\), mutters something vague, and hurriedly
moves on.

Clearly the reason people want this factor is because of how the eta
function transforms under modular transformations. In
\protect\hyperlink{week125}{``Week 125''} I said that the group
\(\mathrm{PSL}(2,\mathbb{Z})\) is generated by two elements \(S\) and
\(T\), and if you look at the formulas there you'll see they act in the
following way on \(\tau\): \[
  \begin{aligned}
    S &\colon \tau \mapsto -\frac{1}{\tau}
  \\T &\colon \tau \mapsto \tau+1
  \end{aligned}
\] The Dedekind eta function satisfies \[
  \begin{aligned}
    \eta\left(-\frac{1}{\tau}\right) &= \left(\frac{\tau}{i}\right)^{\frac12}\eta(\tau)
  \\\eta(\tau+1) &= \exp\left(\frac{2\pi i}{24}\right)\eta\tau
  \end{aligned}
\] The second one is really easy to see from the definition; the first
one is harder. Anyway, using these facts it's easy to see that
\[\frac{1}{\sqrt{\Im(\tau)} |\eta(\tau)|^2}\] is invariant under
\(\mathrm{PSL}(2,\mathbb{Z})\), so it's really a function on moduli
space --- but only if that factor of \(1/\sqrt{\Im(\tau)}\) is in there!

Finally, I'd like to say something about why the conformal anomalies
cancel in 26 dimensions. When I began thinking about this stuff I was
hoping it'd be obvious from the transformation properties of the eta
function --- since they have that promising number ``24'' in them ---
but right now I do \emph{not} see anything like this going on. Instead,
it seems to be something like this: in the partition function
\[Z = \int (\int \exp(-\langle X, \Delta X\rangle) dX) f(g) d[g]\] the
mysterious function \(f\) is basically just the determinant of the
Laplacian on \emph{vector fields} on the torus. So ignoring those darn
zero eigenvalues the whole integrand here is
\[\operatorname{det}(\mathrm{laplacian})^{\frac{n}{2}} \operatorname{det}(\mathrm{laplacian}')\]
where ``\(\mathrm{laplacian}\)'' is the Laplacian on real-valued
functions and " \(\mathrm{laplacian'}\)" is the Laplacian on vector
fields. Now these determinants aren't well-defined functions on the
space of conformal structures; they're really sections of certain
``determinant bundles''. But in this situation, the determinant bundle
for the Laplacian on vector fields \emph{just so happens} to be the 13th
tensor power of the determinant bundle for the Laplacian on functions
--- so the whole expression above is a well-defined function on the
space of conformal structures, and thence on moduli space, precisely
when \(n = 26\)!!!

Now this ``just so happens'' cannot really be a coincidence. There
\emph{are} no coincidences in mathematics. That's why it pays to be
paranoid when you're a mathematician: nothing is random, everything fits
into a grand pattern, it's all just staring you in the face if you'd
only notice it. (Chaitin has convincingly argued otherwise using
G\"odel's theorem, and certainly some patterns in mathematics seem
``purely accidental'', but right now I'm just waxing rhapsodic,
expressing a feeling one sometimes gets\ldots.)

Indeed, look at the proof in Nash's book that one of these determinant
bundles is the 13th tensor power of the other --- I think this result is
due to Mumford, but Nash's proof is easy to read. What does he do? He
works out the first Chern class of both bundles using the index theorem
for families, and he gets something involving the Todd genus --- and the
Todd genus, as we all know, is defined using the same function
\[\frac{x}{1 - e^x} = -1 + \frac{x}{2} - \frac{x^2}{12} + \ldots\] 
that
we talked about in \protect\hyperlink{week126}{``Week 126''} when
computing the zero-point energy of the bosonic string! And yet again,
it's that darn \(-1/12\) in the power series expansion that makes
everything tick. That's where the \(13\) comes from! It's all an
elaborate conspiracy!

But of course the conspiracy is far grander than I've even begun to let
on. If we keep digging away at it, we're eventually led to nothing other
than\ldots.

\begin{center}
\textbf{monstrous moonshine}!!!
\end{center}
But I don't have the energy to talk about \emph{that} now. For more,
try:

\begin{enumerate}
\def\labelenumi{\arabic{enumi})}
\setcounter{enumi}{6}
\item
  Richard E.\ Borcherds, ``What is moonshine?'', talk given upon winning
  the Fields medal, available as
  \href{https://arxiv.org/abs/math.QA/9809110}{\texttt{math.QA/9809110}}.
\item
  Peter Goddard, ``The work of R.\ E.\ Borcherds'', available as
  \href{https://arxiv.org/abs/math.QA/9808136}{\texttt{math.QA/9808136}}.
\end{enumerate}

Okay, if you've actually read this far, you deserve a treat. First, try
this cartoon, which you'll see is quite relevant:

\begin{enumerate}
\def\labelenumi{\arabic{enumi})}
\setcounter{enumi}{8}
\tightlist
\item
  Cartoon by J. F. Cartier,
  \href{https://web.archive.org/web/19980701182649/http://www.physik.uni-frankfurt.de:80/~jr/gif/cartoon/cart0785.gif}{\texttt{https://web.archive.org/web/19980701182649/http:/}}  \href{https://web.archive.org/web/19980701182649/http://www.physik.uni-frankfurt.de:80/~jr/gif/cartoon/cart0785.gif}{\texttt{/www.physik.uni-frankfurt.de:80/\textasciitilde{}jr/gif/cartoon/cart0785.gif}}
\end{enumerate}

Second, let's calculate the determinant of an operator \(A\) whose
eigenvalues are the numbers \(1, 2, 3, \ldots\). You can think of this
operator as the Hamiltonian for the wave equation on the circle, where
we only keep the right-moving modes. As I already said, the zeta
function of this operator is the Riemann zeta function. This function
has \(\zeta'(0) = -\ln(2 \pi)/2\), so using our cute formula relating
determinants and zeta functions, we get
\[\operatorname{det}(A) = \exp(-\zeta'(0)) = (2\pi)^{\frac12}.\] Just
for laughs, if we pretend that the determinant of \(A\) is the product
of its eigenvalues as in the finite-dimensional case, we get:
\[1\cdot 2\cdot 3\cdot \ldots = (2 \pi)^{\frac12}\] or if you really
want to ham it up, \[\infty! = (2 \pi)^{\frac12}.\] Cute, eh? Dan Piponi
told me this, as well as some of the other things I've been talking
about. You can also find it in Bost's paper.

\begin{center}\rule{0.5\linewidth}{0.5pt}\end{center}

Notes and digressions:

\begin{itemize}
\item
  In all of the above, I put a minus sign into my Laplacian, so that it
  has nonnegative eigenvalues. This is common among erudite mathematical
  physics types, who like ``positive elliptic operators''.
\item
  The zeta function trick for defining the determinant of the Laplacian
  works for any positive elliptic operator on a compact manifold. A huge
  amount has been written about this trick. It's all based on the fact
  that the zeta function of a positive elliptic operator analytically
  continues to \(s = 0\). This fact was proved by Seeley:

  \begin{enumerate}
  \def\labelenumi{\arabic{enumi})}
  \setcounter{enumi}{9}
  \tightlist
  \item
    R. T. Seeley, ``Complex powers of an elliptic operator'',
    \emph{Proc. Symp. Pure Math.} \textbf{10} (1967), 288--307.
  \end{enumerate}
\item
  Why is the Polyakov action \(\langle X, \Delta X\rangle\) conformally
  invariant? Because the Laplacian has dimensions of
  \(1/\mathrm{length}^2\), while the integral used to define the inner
  product has dimensions of \(\mathrm{length}^2\), since the torus is
  \(2\)-dimensional. This is the magic of 2 dimensions! The path
  integral for higher-dimensional ``branes'' has not yet been made
  precise, because this magic doesn't happen there.
\item
  About Euler's weirdly beautiful formula for \(\pi\), Robin Chapman
  writes:

  \begin{quote}
  \[\frac{\pi}{2} = \frac32\cdot\frac56\cdot\frac76\cdot\frac{11}{10}\cdot\frac{13}{14}\cdot\frac{17}{18}\cdot\frac{19}{18}\cdot\ldots\tag{1}\]
  Using the Euler product for \(\zeta(2)\) gives
  \[\frac{\pi^2}{6} = \zeta(2) = 1+\frac{1}{2^2}+\frac{1}{3^2}+\ldots = \left(\frac{4}{3}\right)\left(\frac{9}{8}\right)\left(\frac{25}{24}\right)\cdots\left(\frac{p^2}{p^2-1}\right)\cdots\]
  and dropping the \(p=2\) term and dividing by (\(1\)) we see that
  (\(1\)) is equivalent to
  \[\frac{\pi}{4} = \left(\frac{3}{4}\right)\left(\frac{5}{4}\right)\left(\frac{7}{8}\right)\cdots\left(\frac{p}{p-\chi(p)}\right)\cdots \tag{2}\]
  where the numerators are odd primes and the denominators are their
  adjacent multiples of \(4\). Also \(\chi\) is the modulo \(4\)
  Dirichlet character. Now
  \[\frac{p}{p-\chi(p)} = 1+\frac{\chi(p)}{p}+\frac{\chi(p^2)}{p^2}+\ldots\]
  and if we multiply these formally the RHS of (\(2\)) becomes
  \[1-\frac13+\frac15-\frac17+\frac19-\ldots\] i.e., Gregory's series
  for \(\pi/4\). Admittedly it's not apparent that this formal
  manipulation is valid. However for Dirichlet \(L\)-functions the Euler
  product is valid at \(s = 1\). This requires some delicate analysis:
  for details see Landau's book on prime numbers or Davenport's
  \emph{Multiplicative Number Theory}.
  \end{quote}
\end{itemize}

\hypertarget{week128}{%
\section{January 4, 1999}\label{week128}}

This week I'd like to catch you up on the latest developments in quantum
gravity. First, a book that everyone can enjoy:

\begin{enumerate}
\def\labelenumi{\arabic{enumi})}
\tightlist
\item
  John Archibald Wheeler and Kenneth Ford, \emph{Geons, Black Holes, and
  Quantum Foam: A Life in Physics}, Norton, New York, 1998.
\end{enumerate}

This is John Wheeler's autobiography. If Wheeler's only contribution to
physics was being Bohr's student and Feynman's thesis advisor, that in
itself would have been enough. But he did much more. He played a crucial
role in the Manhattan project and the subsequent development of the
hydrogen bomb. He worked on nuclear physics, cosmic rays, muons and
other elementary particles. And he was also one of the earlier people to
get really excited about the more outlandish implications of general
relativity. For example, he found solutions of Einstein's equation that
correspond to regions of gravitational field held together only by their
own gravity, which he called ``geons''. He was not the first to study
black holes, but he was one of the first people to take them seriously,
and he invented the term ``black hole''. And the reason he is \emph{my}
hero is that he took seriously the challenge of reconciling general
relativity and quantum theory. Moreover, he recognized how radical the
ideas needed to accomplish this would be --- for example, the idea that
spacetime might not be truly be a continuum at short distance scales,
but instead some sort of ``quantum foam''.

Anyone interested in the amazing developments in physics during the 20th
century should read this book! Here is the story of how he first met
Feynman:

\begin{quote}
Dick Feynman, who had earned his bachelor's degree at MIT, showed up at
my office door as a brash and appealing twenty-one-year-old in the fall
of 1939 because, as a new student with a teaching assistantship, he had
been assigned to grade papers for me in my mechanics course. As we sat
down to talk about the course and his duties, I pulled out and placed on
the table between us a pocket watch. Inspired by my father's keenness
for time-and-motion studies, I was keeping track of how much time I
spent on teaching and teaching-related activities, how much on research,
and how much on departmental or university chores. This meeting was in
the category of teaching-related. Feynman may have been a little taken
aback by the watch but he was not one to be intimidated. He went out and
bought a dollar watch (as I learned later), so he would be ready for our
next meeting. When we got together again, I pulled out my watch and put
it on the table between us. Without cracking a smile, Feynman pulled out
his watch and put it on the table next to mine. His theatrical sense was
perfect. I broke down laughing, and soon he was laughing as hard as I,
until both of us had tears in our eyes. It took quite a while for us to
sober up and get on with our discussion. This set the tone for a
wonderful friendship that endured for the rest of his life.
\end{quote}

Next for something a wee bit more technical:

\begin{enumerate}
\def\labelenumi{\arabic{enumi})}
\setcounter{enumi}{1}
\tightlist
\item
  Steven Carlip, \emph{Quantum Gravity in 2+1 Dimensions}, Cambridge
  U.\ Press, Cambridge, 1998.
\end{enumerate}
\noindent
If you want to learn about quantum gravity in 2+1 dimensions this is the
place to start, because Carlip is the world's expert on this subject,
and he's pretty good at explaining things.

(By the way, physicists write ``2+1 dimensions'', not because they can't
add, but to emphasize that they are talking about 2 dimensions of space
and 1 dimension of time.)

Quantum gravity in 2+1 dimensions is just a warmup for what physicists
are really interested in --- quantum gravity in 3+1 dimensions. Going
down a dimension really simplifies things, because Einstein's equations
in 2+1 dimensions say that the energy and momentum flowing through a
given point of spacetime completely determine the curvature there,
unlike in higher dimensions. In particular, spacetime is \emph{flat} in
the vacuum in 2+1 dimensions, so there's no gravitational radiation.
Nonetheless, quantum gravity in 2+1 dimensions is very interesting, for
a number of reasons. Most importantly, we can solve the equations
exactly, so we can use it as a nice testing-ground for all sorts of
ideas people have about quantum gravity in 3+1 dimensions.

Quantum gravity is hard for various reasons, but most of all it's hard
because, unlike traditional quantum field theory, it's a
``background-free'' theory. What I mean by this is that there's no fixed
way of measuring times and distances. Instead, times and distances must
be measured with the help of the geometry of spacetime, and this
geometry undergoes quantum fluctuations. That throws most of our usual
methods for doing physics right out the window! Quantum gravity in 2+1
dimensions gives us, for the first time, an example of a background-free
theory where we can work out everything in detail.

Here's the table of contents of Carlip's book:

\begin{quote}
\begin{enumerate}
\def\labelenumi{\arabic{enumi}.}
\tightlist
\item
  Why (2+1)-dimensional gravity?
\item
  Classical general relativity in 2+1 dimensions
\item
  A field guide to the (2+1)-dimensional spacetimes
\item
  Geometric structures and Chern--Simons theory
\item
  Canonical quantization in reduced phase space
\item
  The connection representation
\item
  Operator algebras and loops
\item
  The Wheeler-DeWitt equation
\item
  Lorentzian path integrals
\item
  Euclidean path integrals and quantum cosmology
\item
  Lattice methods
\item
  The (2+1)-dimensional black hole
\item
  Next steps A. Appendix: The topology of manifolds B. Appendix:
  Lorentzian metrics and causal structure C. Appendix: Differential
  geometry and fiber bundles
\end{enumerate}
\end{quote}

And now for some stuff that's available online. First of all, anyone who
wants to keep up with research on gravity should remember to read
``Matters of Gravity''. I've talked about it before, but here's the
latest edition:

\begin{enumerate}
\def\labelenumi{\arabic{enumi})}
\setcounter{enumi}{2}
\tightlist
\item
  Jorge Pullin, editor, \emph{Matters of Gravity}, vol.~12, available at
  \href{https://arxiv.org/abs/gr-qc/9809031}{\texttt{gr-qc/9809031}}.
\end{enumerate}

There's a lot of good stuff in here. Quantum gravity buffs will
especially be interested in Gary Horowitz's article ``A nonperturbative
formulation of string theory?'' and Lee Smolin's ``Neohistorical
approaches to quantum gravity''. The curious title of Smolin's article
refers to \emph{new} work on quantum gravity involving a sum over
\emph{histories} --- or in other words, spin foam models.

Even if you can't go to a physics talk, these days you can sometimes
find it on the world-wide web. Here's one by John Barrett:

\begin{enumerate}
\def\labelenumi{\arabic{enumi})}
\setcounter{enumi}{3}
\tightlist
\item
  John W. Barrett, ``State sum models for quantum gravity'', Penn State
  relativity seminar, August 27, 1998, audio and text of transparencies
  available at
\href{https://web.archive.org/web/20010818150205/http://vishnu.nirvana.phys.psu.edu/online/Html/Seminars/Fall1998/Barrett/}{\texttt{https://}}
\href{https://web.archive.org/web/20010818150205/http://vishnu.nirvana.phys.psu.edu/online/Html/Seminars/Fall1998/Barrett/}{\texttt{web.archive.org/web/20010818150205/http://vishnu.nirvana.phys.psu.edu/}}
\href{https://web.archive.org/web/20010818150205/http://vishnu.nirvana.phys.psu.edu/online/Html/Seminars/Fall1998/Barrett/}{\texttt{online/Html/Seminars/Fall1998/Barrett/}}  
\end{enumerate}
\noindent
Barrett and Crane have a theory of quantum gravity, which I've also
worked on; I discussed it last in \protect\hyperlink{week113}{``Week
113''} and \protect\hyperlink{week120}{``Week 120''}. Before I describe
it I should warn the experts that this theory deals with Riemannian
rather than Lorentzian quantum gravity (though Barrett and Crane are
working on a Lorentzian version, and I hear Friedel and Krasnov are also
working on this). Also, it only deals with vacuum quantum gravity ---
empty spacetime, no matter.

In this theory, spacetime is chopped up into \(4\)-simplices. A
\(4\)-simplex is the \(4\)-dimensional analog of a tetrahedron. To
understand what I'm going to say next, you really need to understand
\(4\)-simplices, so let's start with them.

It's easy to draw a \(4\)-simplex. Just draw 5 dots in a kind of circle
and connect them all to each other! You get a pentagon with a pentagram
inscribed in it. This is a perspective picture of a \(4\)-simplex
projected down onto your \(2\)-dimensional paper. If you stare at this
picture you will see the \(4\)-simplex has 5 tetrahedra, 10 triangles,
10 edges and 5 vertices in it.

The shape of a \(4\)-simplex is determined by 10 numbers. You can take
these numbers to be the lengths of its edges, but if you want to be
sneaky you can also use the areas of its triangles. Of course, there are
some constraints on what areas you can choose for there to \emph{exist}
a 4-simplex having triangles with those areas. Also, there are some
choices of areas that fail to make the shape \emph{unique}: for one of
these bad choices, the \(4\)-simplex can flop around while keeping the
areas of all its triangles fixed. But generically, this non-uniqueness
doesn't happen.

In Barrett and Crane's theory, we chop spacetime into \(4\)-simplices
and describe the geometry of spacetime by specifying the area of each
triangle. But the geometry is ``quantized'', meaning that the area takes
a discrete spectrum of possible values, given by \[\sqrt{j(j+1)}\] where
the ``spin'' \(j\) is a number of the form \(0, 1/2, 1, 3/2, \ldots\).
This formula will be familiar to you if you've studied the quantum
mechanics of angular momentum. And that's no coincidence! The cool thing
about this theory of quantum gravity is that you can discover it just by
thinking a long time about general relativity and the quantum mechanics
of angular momentum, as long as you also make the assumption that
spacetime is chopped into \(4\)-simplices.

So: in Barrett and Crane's theory the geometry of spacetime is described
by chopping spacetime into \(4\)-simplices and labelling each triangle
with a spin. Let's call such a labelling a ``quantum 4-geometry''.
Similarly, the geometry of space is described by chopping space up into
tetrahedra and labelling each triangle with a spin. Let's call this a
``quantum 3-geometry''.

The meat of the theory is a formula for computing a complex number
called an ``amplitude'' for any quantum 4-geometry. This number plays
the usual role that amplitudes do in quantum theory. In quantum theory,
if you want to compute the probability that the world starts in some
state \(\psi\) and ends up in some state \(\psi'\), you just look at all
the ways the world can get from \(\psi\) to \(\psi'\), compute an
amplitude for each way, add them all up, and take the square of the
absolute value of the result. In the special case of quantum gravity,
the states are quantum 3-geometries, and the ways to get from one state
to another are quantum 4-geometries.

So, what's the formula for the amplitude of a quantum 4-geometry? It
takes a bit of work to explain this, so I'll just vaguely sketch how it
goes. First we compute amplitudes for each \(4\)-simplex and multiply
all these together. Then we compute amplitudes for each triangle and
multiply all these together. Then we multiply these two numbers.

(This is analogous to how we compute amplitudes for Feynman diagrams in
ordinary quantum field theory. A Feynman diagram is a graph whose edges
have certain labellings. To compute its amplitude, first we compute
amplitudes for each edge and multiply them all together. Then we compute
amplitudes for each vertex and multiply them all together. Then we
multiply these two numbers. One goal of work on ``spin foam models'' is
to more deeply understand this analogy with Feynman diagrams.)

Anyway, to convince oneself that this formula is ``good'', one would
like to relate it to other approaches to quantum gravity that also
involve \(4\)-simplices. For example, there is the Regge calculus, which
is a discretized version of \emph{classical} general relativity. In this
approach you chop spacetime into \(4\)-simplices and describe the shape
of each \(4\)-simplex by specifying the lengths of its edges. Regge
invented a formula for the ``action'' of such a geometry which
approaches the usual action for classical general relativity in the
continuum limit. I explained the formula for this ``Regge action'' in
\protect\hyperlink{week120}{``Week 120''}.

Now if everything were working perfectly, the amplitude for a
\(4\)-simplex in the Barrett--Crane model would be close to \(\exp(iS)\),
where \(S\) is the Regge action of that \(4\)-simplex. This would mean
that the Barrett--Crane model was really a lot like a path integral in
quantum gravity. Of course, in the Barrett--Crane model all we know is
the areas of the triangles in each \(4\)-simplex, while in the Regge
calculus we know the lengths of its edges. But we can translate between
the two, at least generically, so this is no big deal.

Recently, Barrett and Williams came up with a nice argument saying that
in the limit where the triangles have large areas, the amplitude for a
4-simplex in the Barrett--Crane theory is proportional, not to
\(\exp(iS)\), but to \(\cos(S)\):

\begin{enumerate}
\def\labelenumi{\arabic{enumi})}
\setcounter{enumi}{4}
\tightlist
\item
  John W. Barrett and Ruth M. Williams, ``The asymptotics of an
  amplitude for the \(4\)-simplex'', available as
  \href{https://arxiv.org/abs/gr-qc/9809032}{\texttt{gr-qc/9809032}}.
\end{enumerate}

This argument is not rigorous --- it uses a stationary phase
approximation that requires further justification. But Regge and Ponzano
used a similar argument to show the same sort of thing for quantum
gravity in 3 dimensions, and their argument was recently made rigorous
by Justin Roberts, with a lot of help from Barrett:

\begin{enumerate}
\def\labelenumi{\arabic{enumi})}
\setcounter{enumi}{5}
\tightlist
\item
  Justin Roberts, ``Classical \(6j\)-symbols and the tetrahedron'',
   available as \href{https://arxiv.org/abs/math-ph/9812013}{\texttt{math-ph/}} \href{https://arxiv.org/abs/math-ph/9812013}{\texttt{9812013}}.
\end{enumerate}
\noindent
So one expects that with work, one can make Barrett and Williams'
argument rigorous.

But what does it mean? Why does he get \(\cos(S)\) instead of
\(\exp(iS)\)? Well, as I said, the same thing happens one dimension down
in the so-called Ponzano--Regge model of \(3\)-dimensional Riemannian
quantum gravity, and people have been scratching their heads for decades
trying to figure out why. And by now they know the answer, and the same
answer applies to the Barrett--Crane model.

The problem is that if you describe \(4\)-simplex using the areas of its
triangles, you don't \emph{completely} know its shape. (See, I lied to
you before --- that's why you gotta read the whole thing.) You only know
it \emph{up to reflection}. You can't tell the difference between a
\(4\)-simplex and its mirror-image twin using only the areas of its
triangles! When one of these has Regge action \(S\), the other has
action \(-S\). The Barrett- Crane model, not knowing any better, simply
averages over both of them, getting
\[\frac12(\exp(iS) + \exp(-iS)) = \cos(S)\] So it's not really all that
bad; it's doing the best it can under the circumstances. Whether this is
good enough remains to be seen.

(Actually I didn't really \emph{lie} to you before; I just didn't tell
you my definition of ``shape'', so you couldn't tell whether
mirror-image 4-simplices should count as having the same shape.
Expository prose darts between the Scylla of overwhelming detail and the
Charybdis of vagueness.)

Okay, on to a related issue. In the Barrett--Crane model one describes a
quantum 4-geometry by labelling all the triangles with spins. This
sounds reasonable if you think about how the shape of a \(4\)-simplex is
almost determined by the areas of its triangles. But if you actually
examine the derivation of the model, it starts looking more odd. What
you really do is take the space of geometries of a \emph{tetrahedron}
embedded in \(\mathbb{R}^4\), and use a trick called geometric
quantization to get something called the ``Hilbert space of a quantum
tetrahedron in 4 dimensions''. You then build your \(4\)-simplices out
of these quantum tetrahedra.

Now the Hilbert space of a quantum tetrahedron has a basis labelled by
the eigenvalues of operators corresponding to the areas of its 4
triangular faces. In physics lingo, it takes 4 ``quantum numbers'' to
describe the shape of a quantum tetrahedron in 4 dimensions.

But classically, the shape of a tetrahedron is \emph{not} determined by
the areas of its triangles: it takes 6 numbers to specify its shape, not
just 4. So there is something funny going on.

At first some people thought there might be more states of the quantum
tetrahedron than the ones Barrett and Crane found. But Barbieri came up
with a nice argument suggesting that Barrett and Crane had really found
all of them:

\begin{enumerate}
\def\labelenumi{\arabic{enumi})}
\setcounter{enumi}{6}
\tightlist
\item
  Andrea Barbieri, ``Space of the vertices of relativistic spin
  networks'', available as
  \href{https://arxiv.org/abs/gr-qc/9709076}{\texttt{gr-qc/9709076}}.
\end{enumerate}

While convincing, this argument was not definitive, since it assumed
something plausible but not yet proven --- namely, that the ``\(6j\)
symbols don't have too many exceptional zeros''. Later, Mike
Reisenberger came up with a completely rigorous argument:

\begin{enumerate}
\def\labelenumi{\arabic{enumi})}
\setcounter{enumi}{7}
\tightlist
\item
  Michael P. Reisenberger, ``On relativistic spin network vertices'',
  available as
  \href{https://arxiv.org/abs/gr-qc/9809067}{\texttt{gr-qc/9809067}}.
\end{enumerate}

But while this settled the facts of the matter, it left open the
question of ``why'' --- why does it take \emph{6} numbers to describe
the shape of classical tetrahedron in 4 dimensions but only \emph{4}
numbers to describe the shape of a quantum one? John Barrett and I have
almost finished a paper on this, so I'll give away the answer.

Not surprisingly, the key is that in quantum mechanics, not all
observables commute. You only use the eigenvalues of \emph{commuting}
observables to label a basis of states. The areas of the quantum
tetrahedron's faces commute, and there aren't any other independent
commuting observables. It's a bit like how in classical mechanics you
can specify both the position and momentum of a particle, but in quantum
mechanics you can only specify one.

This isn't news, of course. And indeed, people knew perfectly well that
for this reason, it takes only \emph{5} numbers to describe the shape of
a quantum tetrahedron in 3 dimensions. The real puzzle was why it takes
even fewer numbers when your quantum tetrahedron lives in 4 dimensions!
It seemed bizarre that adding an extra dimension would reduce the number
of degrees of freedom! But it's true, and it's just a spinoff of the
uncertainty principle. Crudely speaking, in 4 dimensions the fact that
you know your tetrahedron lies in some hyperplane makes you unable to
know as much about its shape.

Here are some other talks available on the web:

\begin{enumerate}
\def\labelenumi{\arabic{enumi})}
\setcounter{enumi}{8}
\tightlist
\item
  Abhay Ashtekar, Chris Beetle and Steve Fairhurst, ``Mazatlan lectures
  on black holes'', slides available at
  \texttt{http://vishnu.nirvana.phys.psu.edu/online/Html/Conferences/Mazatlan/}
\end{enumerate}

These explain a new concept of ``nonrotating isolated horizon'' which
allow one to formulate and prove the zeroth and first laws of black hole
mechanics in a way that only refers to the geometry of spacetime near
the horizon. For more details try:

\begin{enumerate}
\def\labelenumi{\arabic{enumi})}
\setcounter{enumi}{9}
\tightlist
\item
  Abhay Ashtekar, Chris Beetle and S.\ Fairhurst, ``Isolated horizons: a
  generalization of black hole mechanics'', available as
  \href{https://arxiv.org/abs/gr-qc/9812065}{\texttt{gr-qc/9812065}}.
\end{enumerate}

This concept also serves as the basis for a forthcoming 2-part paper
where Ashtekar, Corichi, Krasnov and I compute the entropy of a quantum
black hole (see \protect\hyperlink{week112}{``Week 112''} for more on
this).

Finally, here are a couple more papers. I don't have time to say much
about them, but they're both pretty neat:

\begin{enumerate}
\def\labelenumi{\arabic{enumi})}
\setcounter{enumi}{10}
\tightlist
\item
  Matthias Arnsdorf and R.\ S.\ Garcia, ``Existence of spinorial states in
  pure loop quantum gravity'', available as
  \href{https://arxiv.org/abs/gr-qc/9812006}{\texttt{gr-qc/9812006}}.
\end{enumerate}

I'll just quote the abstract:

\begin{quote}
We demonstrate the existence of spinorial states in a theory of
canonical quantum gravity without matter. This should be regarded as
evidence towards the conjecture that bound states with particle
properties appear in association with spatial regions of non-trivial
topology. In asymptotically trivial general relativity the momentum
constraint generates only a subgroup of the spatial diffeomorphisms. The
remaining diffeomorphisms give rise to the mapping class group, which
acts as a symmetry group on the phase space. This action induces a
unitary representation on the loop state space of the Ashtekar
formalism. Certain elements of the diffeomorphism group can be regarded
as asymptotic rotations of space relative to its surroundings. We
construct states that transform non-trivially under a \(2\pi\) rotation:
gravitational quantum states with fractional spin.
\end{quote}

\begin{enumerate}
\def\labelenumi{\arabic{enumi})}
\setcounter{enumi}{13}
\tightlist
\item
  Steve Carlip, ``Black hole entropy from conformal field theory in any
  dimension'', available as
  \href{https://arxiv.org/abs/hep-th/9812013}{\texttt{hep-th/9812013}}.
\end{enumerate}

Again, here's the abstract:

\begin{quote}
When restricted to the horizon of a black hole, the `gauge' algebra of
surface deformations in general relativity contains a physically
important Virasoro subalgebra with a calculable central charge. The
fields in any quantum theory of gravity must transform under this
algebra; that is, they must admit a conformal field theory description.
With the aid of Cardy's formula for the asymptotic density of states in
a conformal field theory, I use this description to derive the
Bekenstein--Hawking entropy. This method is universal - it holds for any
black hole, in any dimension, and requires no details of quantum gravity
--- but it is also explicitly statistical mechanical, based on the
counting of microscopic states.
\end{quote}

On Thursday I'm flying to Schladming, Austria to attend a workshop on
geometry and physics organized by Harald Grosse and Helmut Gausterer.
Some cool physicists will be there, like Daniel Kastler and Julius Wess.
If I understand what they're talking about I'll try to explain it here.
Happy new year!

\begin{center}\rule{0.5\linewidth}{0.5pt}\end{center}

\textbf{Addendum:} Above I wrote:

\begin{quote}
Recently, Barrett and Williams came up with a nice argument saying that
in the limit where the triangles have large areas, the amplitude for a
\(4\)-simplex in the Barrett--Crane theory is proportional, not to
\(\exp(iS)\), but to \(\cos(S)\)\ldots.

This argument is not rigorous --- it uses a stationary phase
approximation that requires further justification. But similar argument
to show the same sort of thing for quantum gravity in 3 dimensions, and
their argument was recently made rigorous by Justin Roberts, with a lot
of help from Barrett\ldots.

So one expects that with work, one can make Barrett and Williams'
argument rigorous.
\end{quote}

In fact one can't make it rigorous: it's wrong! In the limit of large
areas the amplitude for a \(4\)-simplex in the Barrett--Crane model is
wildly different from \(\cos(S)\), or \(\exp(iS)\), or anything like
that. Dan Christensen, Greg Egan and I showed this in a couple of papers
that I discuss in \protect\hyperlink{week170}{``Week 170''} and
\protect\hyperlink{week198}{``Week 198''}. Our results were confirmed by
John Barrett, Chris Steel, Laurent Friedel and David Louapre.

By now --- I'm writing this in 2009 --- it's generally agreed that the
Barrett--Crane model is wrong and another model is better. To read about
this new model, see \protect\hyperlink{week280}{``Week 280''}.

\hypertarget{week129}{%
\section{February 15, 1999}\label{week129}}

For the last 38 years the Austrians have been having winter workshops on
nuclear and particle physics in a little Alpine ski resort town called
Schladming. This year it was organized by Helmut Gausterer and Hermann
Grosse, and the theme was ``Geometry and Quantum Physics'':

\begin{enumerate}
\def\labelenumi{\arabic{enumi})}
\tightlist
\item
  Geometry and Quantum Physics lectures, 38th Internationale
  Universit\"atswochen f\"ur Kern- und Teilchenphysik, available at
  \href{https://web.archive.org/web/19991114074855/http://physik.kfunigraz.ac.at/utp/iukt/iukt_99/iukt99-lect.html}{\texttt{https://web.archive.org/web/199911}}
 \href{https://web.archive.org/web/19991114074855/http://physik.kfunigraz.ac.at/utp/iukt/iukt_99/iukt99-lect.html}{\texttt{14074855/http://physik.kfunigraz.ac.at/utp/iukt/iukt\_99/iukt99-lect.}} \hfill \break \href{https://web.archive.org/web/19991114074855/http://physik.kfunigraz.ac.at/utp/iukt/iukt_99/iukt99-lect.html}{\texttt{html}}
\end{enumerate}

I was invited to give some talks about spin foam models, and the other
talks looked interesting, so I decided to leave my warm and sunny home
for the chilly north. I flew out to Salzburg in early January and took a
train to Schladming from there. Jet-lagged and exhausted, I almost slept
through my train stop, but I made it and soon collapsed into my hotel
bed.

The next day I alternately slept and prepared my talks. The workshop
began that evening with a speech by Helmut Grosse, a speech by the town
mayor, and a reception featuring music by a brass band. The last two
struck me as a bit unusual --- there's something peculiarly Austrian
about drinking beer and discussing quantum gravity over loud oompah
music! This was also the first conference I've been to that featured
skiing and bowling competitions.

Anyway, there were a number of 4-hour minicourses on different subjects,
which should eventually appear as articles in this book:

\begin{enumerate}
\def\labelenumi{\arabic{enumi})}
\setcounter{enumi}{1}
\tightlist
\item
  \emph{Geometry and Quantum Physics}, proceedings of the 38th Internationale
  Universit\"atswochen f\"ur Kern- und Teilchenphysik, Schladming,
  Austria, Jan.~9-16, 1999, eds.~H. Gausterer, H. Grosse and L. Pittner,
  Lecture Notes in Physics \textbf{543}, Springer, Berlin, 2000.
\end{enumerate}
\noindent
Right now they exist in the form of lecture notes:

\begin{itemize}
\tightlist
\item
  Anton Alekseev: ``Symplectic and noncommutative geometry of systems
  with symmetry''
\item
  John Baez: ``Spin foam models of quantum gravity''
\item
  Cesar Gomez: ``Duality and D-branes''
\item
  Daniel Kastler: ``Noncommutative geometry and fundamental physical
  interactions''
\item
  John Madore: ``An introduction to noncommutative geometry''
\item
  Rudi Seiler: ``Geometric properties of transport in quantum Hall
  systems''
\item
  Julius Wess: ``Physics on noncommutative spacetime structures''
\end{itemize}

All these talks were about different ways of combining quantum theory
and geometry. Quantum theory is so strange that ever since its invention
there has been a huge struggle to come to terms with it at all levels.
It took a while for it to make its full impact in pure mathematics, but
now you can see it happening all over: there are lots of papers on
quantum topology, quantum geometry, quantum cohomology, quantum groups,
quantum logic\ldots{} even quantum set theory! There are even some
fascinating attempts to apply quantum mechanics to unsolved problems in
number theory like the Riemann hypothesis\ldots{} will they bear fruit?
And if so, what does this mean about the world? Nobody really knows yet;
we're in a period of experimentation - a bit of a muddle.

I don't have the energy to summarize all these talks so I'll concentrate
on part of Alekseev's --- just a tiny smidgen of it, actually! But
first, let me just quickly say a word about each speaker's topic.

Alekseev talked about some ideas related to the stationary phase
approximation. This is one of the main tools linking classical mechanics
to quantum mechanics. It's a trick for approximately computing the
integral of a function of the form \(\exp(iS(x))\) knowing only \(S(x)\)
and its 2nd derivative at points where its first derivative vanishes. In
physics, people use it to compute path integrals in the semiclassical
limit where what matters most is paths near the classical trajectories.
Alekseev discussed problems where the stationary phase approximation
gives the exact answer. There's a wonderful thing called the
Duistermaat-Heckman formula which says that this happens in certain
situations with circular symmetry. There are also generalizations to
more complicated symmetry groups. These are related to `equivariant
cohomology' --- more on that later.

I talked about the spin foam approach to quantum gravity. I've already
discussed this in \protect\hyperlink{week113}{``Week 113''},
\protect\hyperlink{week114}{``Week 114''},
\protect\hyperlink{week120}{``Week 120''}, and
\protect\hyperlink{week128}{``Week 128''}, so there's no need to say
more here.

Cesar Gomez gave a wonderful introduction to string theory, starting
from scratch and rapidly working up to T-duality and D-branes. The idea
behind T-duality is very simple and pretty. Basically, if you have
closed strings living in a space with one dimension curled up into a
circle of radius \(R\), there is a symmetry that involves replacing
\(R\) by \(1/R\) and switching two degrees of freedom of the string,
namely the number of times it winds around the curled-up direction and
its momentum in the curled-up direction. Both these numbers are
integers. D-branes are something that shows up when you consider the
consequences of this symmetry for \emph{open} strings.

String theory is rather conservative in that, at least until recently,
it usually treated spacetime as a manifold with a fixed geometry and
only applied quantum mechanics to the description of the strings
wiggling around \emph{in} spacetime. In spin foam models, by contrast,
spacetime itself is modelled quantum-mechanically as a kind of
higher-dimensional version of a Feynman diagram. There are also other
ideas about how to treat spacetime quantum-mechanically. One of them is
to treat the coordinates on spacetime as noncommuting variables. In this
approach, called noncommutative geometry, the uncertainty principle
limits our ability to simultaneously know all the coordinates of a
particle's position, giving spacetime a kind of quantum ``fuzziness''.
Personally I don't find noncommutative geometry convincing as a theory
of physical spacetime, because there are no clues that spacetime
actually has this sort of fuzziness. But I find it quite interesting as
mathematics.

Daniel Kastler talked about Alain Connes' theories of physics based on
noncommutative geometry. He discussed both the original Connes--Lott
version of the Standard Model and newer theories that include gravity.
Kastler is a real character! As usual, his talks lauded Connes to the
heavens and digressed all over the map in a frustrating but entertaining
manner. Throughout the conference, he kept us well-fed with anecdotes,
bringing back the aura of heroic bygone days. A random example: Pauli
liked to work long into the night --- so when a student asked ``Could I
meet you at your office at 9 a.m.?'' he replied ``No, I can't possibly
stay that late''.

One nice idea mentioned by Kastler came from this paper:

\begin{enumerate}
\def\labelenumi{\arabic{enumi})}
\setcounter{enumi}{2}
\tightlist
\item
  Alain Connes, ``Noncommutative geometry and reality'', \emph{J. Math.
  Phys.} \textbf{36} (1995), 6194.
\end{enumerate}
\noindent
The idea is to equip spacetime with extra curled-up dimensions shaped
like the quantum group \(\mathrm{SU}_q(2)\) where \(q\) is a 3rd root of
unity. A quantum group is actually a kind of noncommutative algebra, but
using Connes' ideas you can think of it as a kind of ``space''. If you
mod out this particular algebra by its nilradical, you get the algebra
\(M_1(\mathbb{C})\oplus M_2(\mathbb{C})\oplus M_3(\mathbb{C})\), where
\(M_n(\mathbb{C})\) is the algebra of \(n\times n\) complex matrices.
This has a tantalizing relation to the gauge group of the Standard
Model, namely \(\mathrm{U}(1)\times\mathrm{SU}(2)\times\mathrm{SU}(3)\).

John Madore also spoke about noncommutative geometry, but more on the
general theory and less on the applications to physics. He concentrated
on the notion of a ``differential calculus'' --- a structure you can
equip an algebra with in order to do differential geometry thinking of
it as a kind of ``space''.

Julius Wess also spoke on noncommutative geometry, focussing on a
\(q\)-deformed version of quantum mechanics. The process of
``\(q\)-deformation'' is something you can do not only to groups like
\(\mathrm{SU}(2)\) but also other spaces. You get noncommutative
algebras, and these often have nice differential calculi that let you go
ahead and do noncommutative geometry. Wess had a nice humorous way of
defusing tense situations. When one questioner pointedly asked him
whether the material he was presenting was useful in physics or merely a
pleasant game, he replied ``That's a very good question. I will try to
answer that later. For now you're just like students in calculus: you
don't know why you're learning all this stuff\ldots.'' And when Kastler
and other mathematicians kept hassling him over whether an operator was
self-adjoint or merely hermitian, he begged for mercy by saying ``I
would like to be a physicist. That was my dream from the beginning.''

Anyway, I hope that from these vague descriptions you get some sense of
the ferment going on in mathematical physics these days. Everyone agrees
that quantum theory should change our ideas about geometry. Nobody
agrees on how.

Now let me turn to Alekseev's talk. In addition to describing his own
work, he explained many things I'd already heard about. But he did it so
well that I finally understood them! Let me talk about one of these
things: equivariant deRham cohomology. For this, I'll assume you know
about deRham cohomology, principal bundles, connections and curvature.
So I assume you know that given a manifold \(M\), we can learn a lot
about its topology by looking at differential forms on \(M\) and
figuring out the space of closed \(p\)-forms modulo exact ones --- the
so-called \(p\)th deRham cohomology of \(M\). But now suppose that some
Lie group \(G\) acts on \(M\) in a smooth way. What can differential
forms tell us about the topology of this group action?

All sorts of things! First suppose that \(G\) acts freely on \(M\) ---
meaning that \(gx\) is different from \(x\) for any point \(x\) of \(M\)
and any element \(g\) of \(G\) other than the identity. Then the
quotient space \(M/G\) is a manifold. Even better, the map \(M\to M/G\)
gives us a principal \(G\)-bundle with total space \(M\) and base space
\(M/G\).

Can we figure out the deRham cohomology of \(M/G\)? Of course if we were
smart enough we could do it by working out \(M/G\) and then computing
its cohomology. But there's a sneakier way to do it using the
differential forms on \(M\). The map \(M\to M/G\) lets us pull back any
form on \(M/G\) to get a form on \(M\). This lets us think of forms on
\(M/G\) as forms on \(M\) satisfying certain equations --- people call
them ``basic'' differential forms because they come from the base space
\(M/G\).

What are these equations? Well, note that each element \(v\) of the Lie
algebra of \(G\) gives a vector field on \(M\), which I'll also call
\(v\). This give two operations on the differential forms on \(M\): the
Lie derivative \(L_v\) and the interior product \(i_v\). It's easy to
see that any basic differential form is annihilated by these operations
for all \(v\). The converse is true too! So we have some nice equations
describing the basic forms.

If we now take the space of closed basic \(p\)-forms modulo the exact
basic \(p\)-forms, we get the deRham cohomology of \(M/G\)! This lets us
study the topology of \(M/G\) using differential forms on \(M\). It's
very convenient.

If the action of \(G\) on \(M\) isn't free, the quotient space \(M/G\)
might not be a manifold. This doesn't stop us from defining ``basic''
differential forms on \(M\) just as before. We can also define some
cohomology groups by taking the closed basic \(p\)-forms modulo the
exact ones. But topologists know from long experience that another
approach is often more useful. Group actions that aren't free are
touchy, sensitive creatures --- a real nuisance to work with. Luckily,
when you have an action that's not free, you can tweak it slightly to
make it free. This involves ``puffing up'' the space that the group acts
on --- replacing it by a bigger space that the group acts on freely.

For example, suppose you have a group \(G\) acting on a one-point space.
Unless \(G\) is trivial, this action isn't free. In fact, it's about as
far from free as you can get! But we can ``puff it up'' and get a space
called \(EG\). Like the one-point space, \(EG\) is contractible, but
\(G\) acts freely on it. Actually there are various spaces with these
two properties, and it doesn't much matter which one we use --- people
call them all \(EG\). People call the quotient space \(EG/G\) the
``classifying space'' of \(G\), and they denote it by \(BG\).

More generally, suppose we have \emph{any} action of \(G\) on a manifold
\(M\). How can we puff up \(M\) to get a space on which \(G\) acts
freely? Simple: just take its product with \(EG\). Since \(G\) acts on
\(M\) and \(EG\), it acts on the product \(M\times EG\) in an obvious
way. Since \(G\) acts freely on \(EG\), its action on \(M\times EG\) is
free. And since \(EG\) is contractible, the space \(M\times EG\) is a
lot like \(M\), at least as far as topology goes. More precisely, it has
the same homotopy type!

Actually the last 2 paragraphs can be massively generalized at no extra
cost. There's no need for \(G\) to be a Lie group or for \(M\) to be a
manifold. \(G\) can be any topological group and \(M\) can be any
topological space! But since I want to talk about \emph{deRham}
cohomology, I don't need this extra generality here.

Anyway, now we know the right substitute for the quotient space \(M/G\)
when the action of \(G\) on \(M\) isn't free: it's the quotient space
\((M\times EG)/G\).

So now let's figure out how to compute the \(p\)th deRham cohomology of
\((M\times EG)/G\). Since \(G\) acts freely on \(M\times EG\), this
should be just the closed basic \(p\)-forms on \(M\times EG\) modulo the
exact ones, where ``basic'' is defined as before. In fact this is true.
We call the resulting space the \(p\)th ``equivariant deRham
cohomology'' of the space \(M\). It's a kind of well-behaved substitute
for the deRham cohomology of \(M/G\) in the case when \(M/G\) isn't a
manifold.

There's only one slight problem: the space \(EG\) is very big, so it's
not easy to deal with differential forms on \(M\times EG\)!

You'll note that I didn't say much about what \(EG\) looks like. All I
said is that it's some contractible space on which \(G\) acts freely. I
didn't even say it was a manifold, so it's not even obvious that
``differential forms on \(EG\)'' makes sense! If you are smart you can
choose your space \(EG\) so that it's a manifold. However, you'll
usually need it to be infinite-dimensional.

Differential forms make perfect sense on infinite-dimensional manifolds,
but they can be a bit tiresome when we're trying to do explicit
calculations. Luckily there is a small subalgebra of the differential
forms on \(EG\) that's sufficient for the purpose of computing
equivariant cohomology! This is called the ``Weil algebra'', \(WG\).

To guess what this algebra is, let's just list all the obvious
differential forms on \(EG\) that we can think of. Well, I guess none of
them are obvious unless we know a few more facts! First of all, since
the action of \(G\) on \(EG\) is free, the quotient map \(EG\to BG\)
gives us a principal \(G\)-bundle with total space \(EG\) and base space
\(BG\). This bundle is very interesting. It's called the ``universal''
principal \(G\)-bundle. The reason is that any other principal
\(G\)-bundle is a pullback of this one.

(I guess I'm upping the sophistication level again here: I'm assuming
you know how to pull back bundles!)

Even better, if we choose our space \(EG\) so that it's a manifold, then
there is a god-given connection on the bundle \(EG\to BG\), and any
other principal \(G\)-bundle \emph{with connection} is a pullback of
this one.

(And now I'm assuming you know how to pull back connections! However,
this pullback stuff is not necessary in what follows, so just ignore it
if you like.)

Okay, so how can we get a bunch of differential forms on \(EG\) just
using the fact that it's the total space of a \(G\)-bundle equipped with
a connection?

Well, whenever we have a \(G\)-bundle \(E\to B\), we can think of a
connection on it as a \(1\)-form on \(E\) taking values in the Lie
algebra of \(G\). Let's see what differential forms on \(E\) this gives
us! Let's call the connection \(A\). If we pick a basis of the Lie
algebra, we can take the components of \(A\) in this basis, and we get a
bunch of \(1\)-forms \(A_i\) on \(E\). We also get a bunch of
\(2\)-forms \(dA_i\). We also get a bunch of \(2\)-forms
\(A_i\wedge A_j\). And so on.

In general, we can form all possible linear combinations of wedge
products of the \(A_i\)'s and the \(dA_i\)'s. We get a big fat algebra.
In the case when our bundle is \(EG\to BG\), equipped with its god-given
connection, we define this algebra to be the Weil algebra, \(WG\)!

Great. But let's try to define \(WG\) in a purely algebraic way, so we
can do computations with it more easily. We're starting out with the
1-forms \(A_i\) and taking all linear combinations of wedge products of
them and their exterior derivatives. There are in fact no relations
except the obvious ones, so \(WG\) is just ``the supercommutative
differential graded algebra freely generated by the variables \(A_i\)''.
Note: all the mumbo-jumbo about supercommutative differential graded
algebras is a way of mentioning the \emph{obvious} relations.

Warning: people don't usually describe the Weil algebra quite this way.
They usually seem describe it in terms of the connection \(1\)-forms and
curvature \(2\)-forms. However, the curvature is related to the
connection by the formula \(F = dA + A\wedge A\), and if you use this
you can go from the usual description of the Weil algebra to mine --- I
think.

(Actually, people often describe the Weil algebra as an algebra
generated by a bunch of things of degree 1 and a bunch of things of
degree 2, without telling you that the things of degree 1 are secretly
components of a connection \(1\)-form and the things of degree 2 are
secretly components of a curvature \(2\)-form! That's why I'm telling
you all this stuff --- so that if you ever study this stuff you'll have
a better chance of seeing what's going on behind all the murk.)

Okay, so here is the upshot. Say we want to compute the equivariant
deRham cohomology of some manifold \(M\) on which \(G\) acts. In other
words, we want to compute the deRham cohomology of \((M\times EG)/G\).
On the one hand, we can start with the differential forms on
\(M\times EG\), figure out the ``basic'' \(p\)-forms, and take the space
of closed basic \(p\)-forms modulo exact ones. But remember: up to
details of analysis, the algebra of differential forms on \(M\times EG\)
is just the tensor product of the algebra of forms on \(M\) and the
algebra of forms on \(EG\). And we have this nice small ``substitute''
for the algebra of forms on \(EG\), namely the Weil algebra \(WG\). So
let's take the algebra of differential forms on \(M\) and just tensor it
with \(WG\). We get a differential graded algebra with Lie derivative
operations \(L_v\) and interior product operations \(i_v\) defined on
it. We then proceed as before: we take the space of closed basic
elements of degree \(p\) modulo exact ones. Voila! This is something one
can actually compute, with sufficient persistence. And it gives the same
answer, at least when \(G\) is connected and simply connected.

There are all sorts of other things to say. For example, if we take the
simplest posssible case, namely when \(M\) is a single point, this gives
a nice trick for computing the deRham cohomology of \(EG/G = BG\). Guys
in this cohomology ring are called ``characteristic classes'', and
they're really important in physics. Since any principal \(G\)-bundle is
a pullback of \(EG\to BG\), and cohomology classes pull back, these
characteristic classes give us cohomology classes in the base space of
any principal \(G\)-bundle --- thus helping us classify \(G\)-bundles.
But if I started explaining this now, we'd be here all night.

Also sometime I should say more about how to construct \(EG\).

\hypertarget{week130}{%
\section{February 27, 1999}\label{week130}}

All sorts of cool stuff is happening in physics --- and I don't mean
mathematical physics, I mean real live experimental physics! I feel
slightly guilty for not mentioning it on This Week's Finds. Let me
atone.

Here's the big news in a nutshell: we may have been wrong about four
fundamental constants of nature. We thought they were zero, but maybe
they're not! I'm talking about the masses of the neutrinos and the
cosmological constant.

Let's start with neutrinos.

There are three kinds of neutrinos: electron, muon, and \(\tau\)
neutrinos. They are closely akin to the charged particles whose names
they borrow --- the electron, muon and \(\tau\) --- but unlike those
particles they are electrically neutral and very light. They are rather
elusive, since they interact only via the weak force and gravity. I'm
sure you've all heard how a neutrino can easily make it through hundreds
of light years of lead without being absorbed.

But despite their ghostly nature, neutrinos play a very real role in
physics, since radioactive decay often involves a neutron turning into a
proton while releasing an electron and an electron antineutrino. (In
fact, Pauli proposed the existence of neutrinos in 1930 to account for a
little energy that went missing in this process. They were only directly
observed in 1956.) Similarly, in nuclear fusion, a proton may become a
neutron while releasing a positron and an electron neutrino. For
example, when a type II supernova goes off, it emits so many neutrinos
that if you're anywhere nearby, they'll kill you before anything else
gets to you! Indeed, in 1987 a supernova in the Large Magellanic Cloud,
about 100,000 light years away, was detected by four separate neutrino
detectors.

I said neutrinos were ``very light'', but just how light? So far most
work has only given upper bounds. In the 1980s, the Russian ITEP group
claimed to have found a nonzero mass for the electron neutrino, but this
was subsequently blamed on problems with their apparatus. As of now,
laboratory experiments give upper bounds of 4.4 eV for the electron
neutrino mass, 0.17 MeV for the muon neutrino, and 18 MeV for the
\(\tau\) neutrino. By contrast, the electron's mass is 0.511 MeV, the
muon's is 106 MeV, and the \(\tau\)'s is a whopping 1771 MeV.

For this reason, the conventional wisdom used to be that neutrinos were
massless. After all, the electron neutrino is definitely far lighter
than any known particle except the photon --- which is massless. The
larger upper bounds on the other neutrino's masses are mainly due to the
greater difficulty in doing the experiments.

Having neutrinos be massless would also nicely explain their most
stunning characteristic, namely that they're only found in a left-handed
form. What I mean by this is that they spin counterclockwise when viewed
head-on as they come towards you. It turns out that this violation of
left-right symmetry comes fairly easily to massless particles, but only
with more difficulty to massive ones. The reason is simple: massless
particles move at the speed of light, so you can't outrun them. Thus
everyone, regardless of their velocity, agrees on what it means for such
a particle to be spinning one way or another as it comes towards them.
This is not the case for a massive particle!

There was, however, a fly in the ointment. Since the sun is powered by
fusion, it should emit lots of neutrinos. In fact, the standard solar
model predicts that here on earth we are bombarded by 60 billion solar
neutrinos per square centimeter per second! So in the late 1960s, a team
led by Ray Davis set out to detect these neutrinos by putting a tank of
100,000 gallons of perchloroethylene down into a gold mine in Homestake,
South Dakota. Lots of different nuclear reactions are going on in the
sun, producing neutrinos of different energies. The Homestake experiment
can only detect the most energetic ones --- those produced when boron-8
decays into beryllium-8. These neutrinos have enough energy to turn
chlorine-37 in the tank into argon-37. Being a noble gas, the argon can
be separated out and measured. This is not easy --- one only expects
about 4 atoms of argon a day! So the experiment required extreme care
and went on for decades.

They only saw about a quarter as many neutrinos as expected.

Of course, with an experiment as delicate as this, there are always many
possibilities for error, including errors in the standard solar model.
So a Japanese group decided to use a tank of 2,000 tons of water in a
mine in Kamioka to look for solar neutrinos. This ``Kamiokande''
experiment used photomultiplier tubes to detect the Cherenkov radiation
formed by electrons that happen to be hit by neutrinos. Again it was
sensitive only to high-energy neutrinos.

After 5 years, they started seeing signs of a correlation between
sunspot activity and their neutrino count. Interesting. But more
interesting still, they didn't see as many neutrinos as expected. Only
about half as many, in fact.

Starting in the 1990s, various people began to build detectors that
could detect lower-energy neutrinos --- including those produced in the
dominant fusion reactions powering the sun. For this it's good to use
gallium-71, which turns to germanium-71 when bombarded by neutrinos. The
GALLEX detector in Italy uses 30 tons of gallium in the form of gallium
chloride dissolved in water. The SAGE detector, located in a tunnel in
the Caucasus mountains, uses 60 tons of molten metallic gallium. This
isn't quite as scary as it sounds, because gallium has a very low
melting point --- it melts in your hand! But still, of course, these
experiments are very difficult.

Again, these experiments didn't see as many neutrinos as expected.

By this point, the theorists had worked themselves into a full head of
steam trying to account for the missing neutrinos. Currently the most
popular theory is that some of the electron neutrinos have turned into
muon and \(\tau\) neutrinos by the time they reach earth. These other
neutrinos would be not be registered by our detectors.

Folks call this hypothetical process ``neutrino oscillation''. For it to
happen, the neutrinos need to have a nonzero mass. After all, a massless
particle moves at the speed of light, so it doesn't experience any
passage of time --- thanks to relativistic time dilation. Only particles
with mass can become something else while they are whizzing along
minding their own business.

If in fact you posit a small mass for the neutrinos, oscillations happen
automatically as long as the ``mass eigenstates'' are different from the
``flavor eigenstates''. By ``flavor'' we mean whether the neutrino is an
electron, muon or \(\tau\) neutrino. For simplicity, imagine that the
state of a neutrino at rest is given by a vector whose 3 components are
the amplitudes for it to be one of the three different flavors. If all
but one of these components are zero we have a neutrino with a definite
flavor --- a ``flavor eigenstate''. On the other hand, the energy of a
particle at rest is basically just its mass. Thus in the present context
the energy of the neutrino is described by a \(3\times3\) self-adjoint
matrix \(H\), the ``Hamiltonian'', whose eigenvectors are called ``mass
eigenstates''. These may or may not be the same as the flavor
eigenstates! Schroedinger's equation says that any state \(\psi\) of the
neutrino evolves as follows: \[\frac{d\psi}{dt} = -iH\psi.\] Thus if
\(\psi\) starts out being a mass eigenstate it stays a mass eigenstate.
But if it starts out being a flavor eigenstate, it won't stay a flavor
eigenstate --- unless the mass and flavor eigenstates coincide! Instead,
it will oscillate.

I bet you were wondering when the math would start. Don't worry, there
won't be much this time.

Anyway, for other particles, like quarks, it's well-known that the mass
and flavor eigenstates \emph{don't} coincide. So we shouldn't be
surprised at neutrino oscillations, at least if neutrinos actually have
nonzero mass.

Actually things are more complicated than I'm letting on. In addition to
oscillating in empty space, it's possible that neutrinos oscillate
\emph{more} as they are passing through the sun itself, thanks to
something called the MSW effect --- named after Mikheyev, Smirnov and
Wolfenstein. And there are two different ways for neutrinos to have
mass, depending on whether they are Dirac spinors or Majorana spinors
(see \protect\hyperlink{week93}{``Week 93''}).

But I don't want to get caught up in theoretical nuances here! I want to
talk about experiments, and I haven't even gotten to the new stuff yet
--- the stuff that's getting everybody \emph{really} confused!

First of all, there's now some laboratory evidence for neutrino
oscillations coming from the Liquid Scintillator Neutrino Detector at
Los Alamos. What these folks do is let positively charged pions decay
into antimuons and muon neutrinos. Then they check to see if any muon
neutrinos become electron neutrinos. They claim that they do! They also
claim to see evidence of muon antineutrinos becoming electron
antineutrinos.

Secondly, and more intriguing still, there are a bunch of experiments
involving atmospheric neutrinos: Super-Kamiokande, Soudan 2, IMB, and
MACRO. You see, when cosmic rays smack into the upper atmosphere, they
produce all sorts of particles, including electron and muon neutrinos
and their corresponding antineutrinos. Cosmic ray experts think they
know how many of each sort of neutrino should be produced. But the
experimenters down on the ground are seeing different numbers!

Again, this could be due to neutrino oscillations. But what's
\textbf{really} cool is that the numbers seem to depend on where the
neutrinos are coming from: from the sky right above the detector, from
right below the detector --- in which case they must have come all the
way through the earth --- or whatever. Neutrinos coming from different
directions take different amounts of time to get from the upper
atmosphere to the detector. Thus an obvious explanation for the
experimental results is that we're actually seeing the oscillation
process \textbf{as it takes place}.

If this is true, we can try to get detailed information about the
neutrino mass matrix from the numbers these experiments are measuring!

And this is exactly what people have been doing. But they're finding
something very strange. If all the experiments are right, and nobody is
making any mistakes, it seems that NO choice of neutrino mass matrix
really fits all the data! To fit all the data, folks need to do
something drastic --- like posit a 4th kind of neutrino!

Now, it's no light matter to posit another neutrino. The known neutrinos
couple to the weak force in almost identical ways. This allows one to
create equal amounts of neutrino-antineutrino pairs of all 3 flavors by
letting Z bosons decay --- the Z being the neutral carrier of the weak
force. When a Z boson seemingly decays into ``nothing'', we can safely
bet that it has decayed into a neutrino- antineutrino pair. In 1989, an
elegant and famous experiment at CERN showed that Z bosons decay into
``nothing'' at exactly the rate one would expect if there were 3 flavors
of neutrino. Thus there can only be extra flavors of neutrino if they
are very massive, if they couple very differently to the weak force, or
if some other funny business is going on.

Now, electron or muon neutrinos are unlikely to oscillate into a
\emph{very massive} sort of neutrino --- basically because of energy
conservation. So if we want an extra neutrino to explain the
experimental results we find ourselves stuck with, it'll have to be one
that couples to the weak force very differently from the ones we know. A
simple, but drastic, possibility is that it not interact via the weak
force at all! Folks call this a ``sterile'' neutrino.

Now, sterile neutrinos would blow a big hole in the Standard Model, much
more so than plain old \emph{massive} neutrinos. So things are getting
very interesting.

Wilczek recently wrote a nice easy-to-read paper describing arguments
that \emph{massive} neutrinos fit in quite nicely with the possibility
that the Standard Model is just part of a bigger, better theory --- a
``Grand Unified Theory''. I sketched the basic ideas of the
\(\mathrm{SU}(5)\) and \(\mathrm{SO}(10)\) grand unified theories in
\protect\hyperlink{week119}{``Week 119''}. Recall that in the
\(\mathrm{SU}(5)\) theory, the left-handed parts of all fermions of a
given generation fit into two irreducible representations of
\(\mathrm{SU}(5)\) --- a 5-dimensional rep and a \(10\)-dimensional rep.
For example, for the first generation, the \(5\)-dimensional rep
consists of the left-handed down antiquark (which comes in 3 colors),
the left-handed electron, and the left-handed electron neutrino. The
\(10\)-dimensional rep consists of the left-handed up quark, down quark,
and up antiquark (which come in 3 colors each), together with the
left-handed positron.

In the \(\mathrm{SO}(10)\) theory, all these particles AND ONE MORE fit
into a single 16-dimensional irreducible representation of
\(\mathrm{SO}(10)\). What could this extra particle be?

Well, since this extra particle transforms trivially under
\(\mathrm{SU}(5)\), it must not feel the electromagnetic, weak or strong
force! Thus it's tempting to take this missing particle to be the
left-handed electron antineutrino. Of course, we don't see such a
particle --- we only see antineutrinos that spin clockwise. But if
neutrinos are massive Dirac spinors there must be such a particle, and
having it not feel the electromagnetic, weak or strong force would
nicely explain \emph{why} we don't see it.

Grotz and Klapdor consider this possibility in their book on the weak
interaction (see below), but unfortunately, it seems this theory would
make the electron neutrino have a mass of about 5 MeV --- much too big!
Sigh. So Wilczek, following the conventional wisdom, assumes the missing
particle is very massive --- he calls it the ``N''. And he summarizes
some arguments that this massive particle could help give the neutrinos
very small masses, via something called the ``seesaw mechanism''.
Unfortunately I don't have the energy to describe this now, so for more
you should look at his paper (referred to below).

To wrap up, let me just say one final thing about the cosmic
significance of the neutrino. Massive neutrinos could account for some
of the ``missing mass'' that cosmologists are worrying about. So there's
an indirect connection between the neutrino mass and the cosmological
constant! The cosmological constant is essentially the energy density of
the vacuum. It was long assumed to be zero, but now there are some
glimmerings of evidence that it's not. In fact, some people are quite
convinced that it's not. The fate of the universe hangs in the
balance\ldots.

Unfortunately I am too tired now to say much more about this. So let me
just give you a nice easy starting-point:

\begin{enumerate}
\def\labelenumi{\arabic{enumi})}
\tightlist
\item
  \emph{Special Report: Revolution in Cosmology}, Scientific American,
  January 1999. Includes the articles ``Surveying space-time with
  supernovae'' by Craig J. Horgan, Robert P. Kirschner and Nicholoas B.
  Suntzeff, ``Cosmological antigravity'' by Lawrence M. Krauss, and
  ``Inflation in a low-density universe'' by Martin A. Bucher and David
  N. Spergel.
\end{enumerate}

How can you learn more about neutrinos? It can't hurt to start here:

\begin{enumerate}
\def\labelenumi{\arabic{enumi})}
\setcounter{enumi}{1}
\tightlist
\item
  Nikolas Solomey, \emph{The Elusive Neutrino}, Scientific American
  Library, 1997.
\end{enumerate}
\noindent
If you want to dig in deeper, you need to learn about the weak force,
since we've only seen neutrinos via their weak interaction with other
particles. The following book is a great place to start:

\begin{enumerate}
\def\labelenumi{\arabic{enumi})}
\setcounter{enumi}{2}
\tightlist
\item
  K. Grotz and H. V. Klapdor, \emph{The Weak Interaction in Nuclear,
  Particle and Astrophysics}, Adam Hilger, Bristol, 1990.
\end{enumerate}
\noindent
Then you'll be ready for this book, which examines every aspect of
neutrinos in detail --- complete with copies of historical papers:

\begin{enumerate}
\def\labelenumi{\arabic{enumi})}
\setcounter{enumi}{3}
\tightlist
\item
  Klaus Winter, ed., \emph{Neutrino Physics}, Cambridge U. Press,
  Cambridge, 1991.
\end{enumerate}
\noindent
And then, if you want to study the physics of \emph{massive}
neutrinos, you should try this:

\begin{enumerate}
\def\labelenumi{\arabic{enumi})}
\setcounter{enumi}{4}
\tightlist
\item
  Felix Boehm and Petr Vogel, \emph{Physics of Massive Neutrinos},
  Cambridge U. Press, Cambridge, 1987.
\end{enumerate}
\noindent
But neutrino physics is moving fast, and lots of the new stuff hasn't
made its way into books yet, so you should also look at other stuff. For
links to lots of great neutrino websites, including websites for most of
the experiments I mentioned, try:

\begin{enumerate}
\def\labelenumi{\arabic{enumi})}
\setcounter{enumi}{5}
\tightlist
\item
  ``The neutrino oscillation industry'',
  \url{https://www.hep.anl.gov/ndk/hypertext/}
\end{enumerate}
\noindent
For some recent general overviews, try these:

\begin{enumerate}
\def\labelenumi{\arabic{enumi})}
\setcounter{enumi}{6}
\item
  Paul Langacker, ``Implications of neutrino mass'',
  \href{https://www.physics.upenn.edu/~pgl/neutrino/jhu/jhu.html}{\texttt{https://www.physics.upenn.}}
  \href{https://www.physics.upenn.edu/~pgl/neutrino/jhu/jhu.html}{\texttt{edu/\textasciitilde{}pgl/neutrino/jhu/jhu.html}}
\item
  Boris Kayser, ``Neutrino mass: where do we stand, and where are we
  going?'', available as
  \href{https://arxiv.org/abs/hep-ph/9810513}{\texttt{hep-ph/9810513}}.
\end{enumerate}
\noindent
For information on various experiments, try these:

\begin{enumerate}
\def\labelenumi{\arabic{enumi})}
\setcounter{enumi}{8}
\item
  GALLEX collaboration, ``GALLEX solar neutrino observations: complete
  results for GALLEX II'', \emph{Phys. Lett.} B357 (1995), 237--247.

  ``Final results of the CR-51 neutrino source experiments in GALLEX'',
  \emph{Phys. Lett.} \textbf{B420} (1998), 114--126.

  ``GALLEX solar neutrino observations: results for GALLEX IV'',
  \emph{Phys. Lett.} \textbf{B447} (1999), 127--133.
\item
  SAGE collaboration, ``Results from SAGE'', \emph{Phys. Lett.}
  \textbf{B328} (1994), 234--248.

  ``The Russian-American gallium experiment (SAGE) CR neutrino source
  measurement'', \emph{Phys. Rev.~Lett.} \textbf{77} (1996), 4708--4711.
\item
  LSND collaboration, ``Evidence for neutrino oscillations from muon
  decay at rest'', \emph{Phys. Rev.} \textbf{C54} (1996) 2685--2708.
  Also available as
  \href{https://arxiv.org/abs/nucl-ex/9605001}{\texttt{nucl-ex/9605001}}.

  ``Evidence for anti-muon-neutrino \(\to\) anti-electron-neutrino
  oscillations from the LSND experiment at LAMPF'', \emph{Phys.
  Rev.~Lett.} \textbf{77} (1996), 3082--3085.  Also available as
  \href{https://arxiv.org/abs/nucl-ex/9605003}{\texttt{nucl-ex/9605003}}.

  ``Evidence for muon-neutrino \(\to\) electron-neutrino oscillations
  from LSND'', \emph{Phys. Rev.~Lett.} \textbf{81} (1998), 1774--1777.
  Also available as
  \href{https://arxiv.org/abs/nucl-ex/9709006}{\texttt{nucl-ex/9709006}}.

  ``Results on muon-neutrino \(\to\) electron-neutrino oscillations from
  pion decay in flight'', \emph{Phys. Rev.} \textbf{C58} (1998),
  2489--2511.
\item
  Super-Kamiokande collaboration, ``Evidence for oscillation of
  atmospheric neutrinos'', \emph{Phys. Rev.~Lett.} \textbf{81} (1998),
  1562--1567.  Also available as
  \href{https://arxiv.org/abs/hep-ex/9807003}{\texttt{hep-ex/9807003}}.
\item
  MACRO collaboration, ``Measurement of the atmospheric neutrino-induced
  upgoing muon flux'', \emph{Phys. Lett.} \textbf{B434} (1998),
  451--457.  Also available as
  \href{https://arxiv.org/abs/hep-ex/9807005}{\texttt{hep-ex/9807005}}.
\item
  IMB collaboration, ``A search for muon-neutrino oscillations with the
  IMB detector'', \emph{Phys. Rev.~Lett.} \textbf{69} (1992),
  1010--1013.
\end{enumerate}
\noindent
For a fairly model-independent attempt to figure out something about
neutrino masses from the latest crop of experiments, see:

\begin{enumerate}
\def\labelenumi{\arabic{enumi})}
\setcounter{enumi}{15}
\tightlist
\item
  V. Barger, T. J. Weiler, and K. Whisnant, ``Inferred 4.4 eV upper
  limits on the muon- and tau-neutrino masses''.  Also available as
  \href{https://arxiv.org/abs/hep-ph/9808367}{\texttt{hep-ph/9808367}}.
\end{enumerate}
\noindent
For a nice summary of the data, and an argument that it's evidence for
the existence of a sterile neutrino, see:

\begin{enumerate}
\def\labelenumi{\arabic{enumi})}
\setcounter{enumi}{16}
\tightlist
\item
  David O. Caldwell, ``The status of neutrino mass''.  Also available
  as
  \href{https://arxiv.org/abs/hep-ph/9804367}{\texttt{hep-ph/9804367}}.
\end{enumerate}
\noindent
For a very readable argument that massive neutrinos are evidence for a
supersymmetric \(\mathrm{SO}(10)\) grand unified theory, see
\noindent
\begin{enumerate}
\def\labelenumi{\arabic{enumi})}
\setcounter{enumi}{17}
\tightlist
\item
  Frank Wilczek, ``Beyond the Standard Model: this time for real''.  Also available as
  \href{https://arxiv.org/abs/hep-ph/9809509}{\texttt{hep-ph/9809509}}.
\end{enumerate}
\noindent
Finally, with all these cracks developing in the Standard Model, it's
nice to think again about the rise of the Standard Model. The following
book is packed with the reminiscences of many theorists and
experimentalists involved in developing this wonderful theory of
particles and forces, including Bjorken, 't Hooft, Veltman, Susskind,
Polyakov, Richter, Iliopoulos, Gell-Mann, Weinberg, Lederman, Goldhaber,
Cronin, and Kobayashi:

\begin{enumerate}
\def\labelenumi{\arabic{enumi})}
\setcounter{enumi}{18}
\tightlist
\item
  Lilian Hoddeson, Laurie Brown, Michael Riordan and Max Dresden, eds.,
  \emph{The Rise of the Standard Model: Particle Physics in the 1960s
  and 1970s}, Cambridge U.\ Press, Cambridge, 1997.
\end{enumerate}
\noindent
It's a must for anyone with an interest in the history of physics!

\hypertarget{week131}{%
\section{March 7, 1999}\label{week131}}

I've been thinking more about neutrinos and their significance for grand
unified theories. The term ``grand unified theory'' sounds rather
pompous, but in its technical meaning it refers to something with
limited goals: a quantum field theory that attempts to unify all the
forces \emph{except} gravity. This limitation lets you pretend spacetime
is flat.

The heyday of grand unified theories began in the mid-1970s, shortly
after the triumph of the Standard Model. As you probably know, the
Standard Model is a quantum field theory describing all known particles
and all the known forces except gravity: the electromagnetic, weak and
strong forces. The Standard Model treats the electromagnetic and weak
forces in a unified way --- so one speaks of the ``electroweak'' force
--- but it treats the strong force seperately.

In 1975, Georgi and Glashow invented a theory which fit all the known
particles of each generation into two irreducible representations of
\(\mathrm{SU}(5)\). Their theory had some very nice features: for
example, it unified the strong force with the electroweak force, and it
explained why quark charges come in multiples of \(1/3\). It also made
some new predictions, most notably that protons decay with a halflife of
something like \(10^{29}\) or \(10^{30}\) years. Of course, it's
slightly inelegant that one needs \emph{two} irreducible representations
of \(\mathrm{SU}(5)\) to account for all the particles of each
generation. Luckily \(\mathrm{SU}(5)\) fits inside \(\mathrm{SO}(10)\)
in a nice way, and Georgi used this to concoct a slightly bigger theory
where all 15 particles of each generation, AND ONE MORE, fit into a
single irreducible representation of \(\mathrm{SO}(10)\). I described
the mathematics of all this in \protect\hyperlink{week119}{``Week
119''}, so I won't do so again here.

What's the extra particle? Well, when you look at the math, one obvious
possibility is a right-handed neutrino. As I explained last week, the
existence of a right-handed neutrino would make it easier for neutrinos
to have mass. And this in turn would allow ``oscillations'' between
neutrinos of different generations --- possibly explaining the
mysterious shortage of electron neutrinos that we see coming from the
sun.

This ``solar neutrino deficit'' had already been seen by 1975, so
everyone got very excited about grand unified theories. The order of the
day was: look for neutrino oscillations and proton decay!

A nice illustration of the mood of the time can be found in a talk
Glashow gave in 1980:

\begin{enumerate}
\def\labelenumi{\arabic{enumi})}
\tightlist
\item
  Sheldon Lee Glashow, ``The new frontier'', in \emph{First Workshop on
  Grand Unification}, eds.~Paul H. Frampton, Sheldon L. Glashow and Asim
  Yildiz, Math Sci Press, Brookline Massachusetts, 1980, pp.~3--8.
\end{enumerate}

I'd like to quote some of his remarks because it's interesting to
reflect on what has happened in the intervening two decades:

\begin{quote}
Pions, muons, positrons, neutrons and strange particles were found
without the use of accelerators. More recently, most developments in
elementary particle physics depended upon these expensive artificial
aids. Science changes quickly. A time may come when accelerators no
longer dominate our field: not yet, but perhaps sooner than some may
think.

Important discoveries await the next generation of accelerators. QCD and
the electroweak theory need further confirmation. We need know how b
quarks decay. The weak interaction intermediaries must be seen to be
believed. The top quark (or perversions needed by topless theories)
lurks just out of range. Higgs may wait to be found. There could well be
a fourth family of quarks and leptons. There may even be unanticipated
surprises. We need the new machines.
\end{quote}

Of course we now know how the b (or ``bottom'') quark decays, we've seen
the t (or ``top'') quark, we've seen the weak interaction
intermediaries, and we're quite sure there is not a fourth generation of
quarks and leptons. There have been no unanticipated surprises.
Accelerators grew ever more expensive until the U.S. Congress withdrew
funding for the Superconducting Supercollider in 1993. The Higgs is
still waiting to be found or proved nonexistent. Experiments at CERN
should settle that issue by 2003 or so.

\begin{quote}
On the other hand, we have for the first time an apparently correct
\emph{theory} of elementary particle physics. It may be, in a sense,
phenomenologically complete. It suggests the possibility that there are
no more surprises at higher energies, at least for energies that are
remotely accessible. Indeed, PETRA and ISR have produced no surprises.
The same may be true for PEP, ISABELLE, and the TEVATRON. Theorists do
expect novel higher-energy phenomena, but only at absurdly inacessible
energies. Proton decay, if it is found, will reinforce the belief in the
great desert extending from \(100\) GeV to the unification mass of
\(10^{14}\) GeV. Perhaps the desert is a blessing in disguise. Ever
larger and more costly machines conflict with dwindling finances and
energy reserves. All frontiers come to an end.

You may like this scenario or not; it may be true or false. But, it is
neither impossible, implausible, or unlikely. And, do not despair nor
prematurely lament the death of particle physics. We have a ways to go
to reach the desert, with exotic fauna along the way, and even the
desolation of a desert can be interesting. The end of the high-energy
frontier in no ways implies the end of particle physics. There are many
ways to skin a cat. In this talk I will indicate several exciting lines
of research that are well away from the high-energy frontier. Important
results, perhaps even extraordinary surprises, await us. But, there is
danger on the way.

The passive frontier of which I shall speak has suffered years of benign
neglect. It needs money and manpower, and it must compete for this with
the accelerator establishment. There is no labor union of physicists who
work at accelerators, but sometimes it seems there is. It has been
argued that plans for accelerator construction must depend on the
``needs'' of the working force: several thousands of dedicated
high-energy experimenters. This is nonsense. Future accelerators must be
built in accordance with scientific, not demographic, prioriries. The
new machines are not labor-intensive, must not be forced to be so. Not
all high energy physicsts can be accomodated at the new machines. The
high-energy physicist has no guaranteed right to work at an accelerator,
he has not that kind of job security. He must respond to the challenge
of the passive frontier .
\end{quote}

Of course, the collapse of the high-energy physics job market and the
death of the Superconducting Supercollider give these words a certain
poignancy. But what is this ``passive frontier'' Glashow mentions? It
means particle physics that doesn't depend on very high energy particle
accelerators. He lists a number of options:

\begin{enumerate}
\def\labelenumi{\Alph{enumi})}
\item
  CP phenomenology. The Standard Model is not symmetrical under
  switching particles and their antiparticles --- called ``charge
  conjugation'', or ``C''. Nor is it symmetrical under switching left
  and right --- called ``parity'', or ``P''. It's almost, but not quite,
  symmetrical under the combination of both operations, called ``CP''.
  Violation of CP symmetry is evident in the behavior of the neutral
  kaon. Glashow suggests looking for CP violation in the form of a
  nonzero magnetic dipole moment for the neutron. As far as I know, this
  has still not been seen.
\item
  New kinds of stable matter. Glashow proposes the search for new stable
  particles as ``an ambitious and risky field of scientific endeavor''.
  People have looked and haven't found anything.
\item
  Neutrino masses and neutrino oscilllations. Glashow claims that
  ``neutrinos should have masses, and should mix''. He now appears to be
  right. It took almost 20 years for the trickle of experimental results
  to become the lively stream we see today, but it happened. He urges
  ``Let us not miss the next nearby supernova!'' Luckily we did not.
\item
  Astrophysical neutrino physics. In addition to solar neutrinos and
  neutrinos from supernovae, there are other interesting connections
  between neutrinos and astrophysics. The background radiation from the
  big bang should contain neutrinos as well as the easier-to-see
  photons. More precisely, there should be about 100 neutrinos of each
  generation per cubic centimeter of space, thanks to this effect. These
  ``relic neutrinos'' have not been seen, but that's okay: by now they
  would be too low in energy to be easily detected. Glashow notes that
  if neutrinos had a nonzero mass, relic neutrinos could contribute
  substantially to the total density of the universe. The heaviest
  generation weighing 30 eV or so might be enough to make the universe
  eventually recollapse! On the other hand, for neutrinos to be
  gravitationally bound to galaxies, they'd need to be at least 20 eV or
  so.
\item
  Magnetic monopoles. Most grand unified theories predict the existence
  of magnetic monopoles due to ``topological defects'' in the Higgs
  fields. Glashow urges people to look for these. This has been done,
  and they haven't been seen.
\item
  Proton decay. As Glashow notes, proton decay would be the ``king of
  the new frontier''. Reflecting the optimism of 1980, he notes that
  ``to some, it is a foregone conclusion that proton decay is about to
  be seen by experiments now abuilding''. But alas, people looked very
  hard and did not find it! This killed the \(\mathrm{SU}(5)\) theory.
  Many people switched to supersymmetric theories, which are more
  compatible with very slow proton decay. But with the continuing lack
  of new experiments to explain, enthusiasm for grand unified theories
  gradually evaporated, and theoretical particle physics took refuge in
  the elegant abstractions of string theory.
\end{enumerate}

But now, 20 years later, interest in grand unified theories seems to be
coming back. We have a rich body of mysterious experimental results
about neutrino oscillations. Somebody should explain them!

On a slightly different note, one of my little side hobbies is to study
the octonions and dream about how they might be useful in physics. One
place they show up is in the \(\mathrm{E}_6\) grand unified theory ---
the next theory up from the \(\mathrm{SO}(10)\) theory. I said a bit
about this in \protect\hyperlink{week119}{``Week 119''}, but I just
bumped into another paper on it in the same conference proceedings that
Glashow's paper appears is:

\begin{enumerate}
\def\labelenumi{\arabic{enumi})}
\setcounter{enumi}{1}
\tightlist
\item
  Feza Gursey, ``Symmetry breaking patterns in \(\mathrm{E}_6\)'', in
  \emph{First Workshop on Grand Unification}, eds.~Paul H. Frampton,
  Sheldon L. Glashow and Asim Yildiz, Math Sci Press, Brookline
  Massachusetts, 1980, pp.~39--55.
\end{enumerate}

He says something interesting that I want to understand someday ---
maybe someone can explain why it's true. He says that \(\mathrm{E}_6\)
is a ``complex'' group, \(\mathrm{E}_7\) is a ``pseudoreal'' group, and
\(\mathrm{E}_8\) is a ``real'' group. This use of terminology may be
nonstandard, but what he means is that \(\mathrm{E}_6\) admits complex
representations that are not their own conjugates, \(\mathrm{E}_7\)
admits complex reps that are their own conjugates, and that all complex
reps of \(\mathrm{E}_8\) are complexifications of real ones (and hence
their own conjugates). This should have something to do with the
symmetry of the Dynkin diagram of \(\mathrm{E}_6\).

Octonions are also prominent in string theory and in the grand unified
theories proposed by my friends Geoffrey Dixon and Tony Smith --- see
\protect\hyperlink{week59}{``Week 59''},
\protect\hyperlink{week91}{``Week 91''}, and
\protect\hyperlink{week104}{``Week 104''}. I'll probably say more about
this someday\ldots.

The reason I'm interested in neutrinos is that I want to learn what
evidence there is for laws of physics going beyond the Standard Model
and general relativity. This is also why I'm trying to learn a bit of
astrophysics. The new hints of evidence for a nonzero cosmological
constant, the missing mass problem, the large-scale structure of the
universe, and even the puzzling \(\gamma\)-ray bursters --- they're all
food for thought along these lines.

The following book caught my eye since it looked like just what I need
--- an easy tutorial in the latest developments in cosmology:

\begin{enumerate}
\def\labelenumi{\arabic{enumi})}
\setcounter{enumi}{2}
\tightlist
\item
  Greg Bothun, \emph{Modern Cosmological Observations and Problems},
  Taylor \& Francis, London, 1998.
\end{enumerate}

On reading it, some of the remarks about particle physics made me
unhappy. For example, Bothun says ``the observed entropy \(S\) of the
universe, as measured by the ratio of baryons to photons, is
\(\sim5\times10^{-10}\).'' But as Ted Bunn explained to me, the entropy
is actually correlated to the ratio of photons to baryons --- the
reciprocal of this number. Bothun also calls the kinetic energy density
of the field postulated in inflationary cosmology, ``essentially an
entropy field that currently drives the uniform expansion and cooling of
the universe''. This makes no sense to me. There are also a large number
of typos, the most embarrassing being ``virilizing'' for
``virializing''.

But there's a lot of good stuff in this book! The author's specialty is
large-scale structure, and I learned a lot about that. Just to set the
stage, recall that the Milky Way has a disc about 30 kiloparsecs in
diameter and contains roughly 100 or 200 billion stars. But our galaxy
is one of a cluster of about 20 galaxies, called the Local Group. In
addition to our galaxy and the Large and Small Magellanic Clouds which
orbit it, this contains the Andromeda Galaxy (also known as M31),
another spiral galaxy called M33, and a bunch of dwarf irregular
galaxies. The Local Group is about a megaparsec in radius.

This is typical. Galaxies often lie in clusters which are a few
megaparsecs in radius, containing from a handful to hundreds of big
galaxies. Some famous nearby clusters include the Virgo cluster (about
20 megaparsecs away) and the Coma cluster (about 120 megaparsecs away).
Thousands of clusters have been cataloged by Abell and collaborators.

And then there are superclusters, each typically containing 3--10
clusters in an elongated ``filament'' about 50 megaparsecs in diameter.
I don't mean to make this sound more neat than it actually is, because
nobody is very sure about intergalactic distances, and the structures
themselves are rather messy. But there are various discernible patterns.
For example, superclusters tend to occur at the edges of larger roundish
``voids'' which have few galaxies in them. These voids are very large,
about 100 or 200 megaparsecs across. In general, galaxies tend to be
falling into denser regions and moving away from the voids. For example,
the Milky Way is falling towards the center of the Local Supercluster at
about 300 kilometers per second, and the Local Supercluster is also
falling towards the next nearest one --- the Hydra-Centaurus
Supercluster --- at about 300 kilometers per second.

Now, if the big bang theory is right, all this stuff was once very
small, and the universe was much more homogeneous. Obviously gravity
tends to amplify inhomogeneities. The problem is to understand in a
quantitative way how these inhomogeneities formed as the universe grew.

Here are a couple of other books that I'm finding useful --- they're a
bit more mathematical than Bothun's. I'm trying to stick to new books
because this subject is evolving so rapidly:

\begin{enumerate}
\def\labelenumi{\arabic{enumi})}
\setcounter{enumi}{3}
\item
  Jayant V. Narlikar, \emph{Introduction to Cosmology}, Cambridge U.\
  Press, Cambridge, 1993.
\item
  Peter Coles and Francesco Lucchin, \emph{Cosmology: The Origin and
  Evolution of Cosmic Structure}, Wiley, New York, 1995.
\end{enumerate}

While I was looking around, I also bumped into the following book on
black holes:

\begin{enumerate}
\def\labelenumi{\arabic{enumi})}
\setcounter{enumi}{5}
\tightlist
\item
  Sandip K. Chakrabarti, ed., \emph{Observational Evidence for Black Holes
  in the Universe}, Kluwer, Dordrecht, 1998.
\end{enumerate}
\noindent
It mentioned some objects I'd never heard of before. I want to tell you
about them, just because they're so cool!

\begin{itemize}
\item
  \emph{X-ray novae}: First, what's a nova? Well, remember that a white
  dwarf is a small, dense, mostly burnt-out star. When one member of a
  binary star is a white dwarf, and the other dumps some of its gas on
  this white dwarf, the gas can undergo fusion and emit a huge burst of
  energy --- as much as 10,000 times what the sun emits in a year. To
  astronomers it may look like a new star is suddenly born --- hence the
  term ``nova''. But not all novae emit mostly visible light --- some
  emit X-rays or even \(\gamma\)-rays. A ``X-ray nova'' is an X-ray
  source that suddenly appears in a few days and then gradually fades
  away in a year or less. Many of these are probably neutron stars
  rather than white dwarfs. But a bunch are probably black holes!
\item
  \emph{Blazars}: A ``blazar'' is a galactic nucleus that's shooting out
  a jet of hot plasma almost directly towards us, exhibiting rapid
  variations in power. Like quasars and other active galactic nuclei,
  these are probably black holes sucking in stars and other stuff and
  forming big accretion disks that shoot out jets of plasma from both
  poles.
\item
  \emph{Mega masers}: A laser is a source of coherent light caused by
  stimulated emission --- a very quantum-mechanical gadget. A maser is
  the same sort of thing but with microwaves. In fact, masers were
  invented before lasers --- they are easier to make because the
  wavelength is longer. In galaxies, clouds of water vapor, hydroxyl,
  silicon monoxide, methanol and other molecules can form enormous
  natural masers. In our galaxy the most powerful such maser is W49N,
  which has a power equal to that of the Sun. But recently, still more
  powerful masers have been bound in other galaxies, usually associated
  with active galactic nuclei. These are called ``mega masers'' and they
  are probably powered by black holes. The first mega maser was
  discovered in 1982; it is a hydroxyl ion maser in the galaxy IC4553,
  with a luminosity about 1000 times that of our sun. Subsequently
  people have found a bunch of water mega masers. The most powerful so
  far is in TXFS2226-184 --- it has a luminosity of about 6100 times
  that of the Sun!
\end{itemize}

\begin{center}\rule{0.5\linewidth}{0.5pt}\end{center}

\textbf{Addendum:} Here is something from Allen Knutson in response to
my remark that \(\mathrm{E}_6\) has complex representations that aren't
their own conjugates. I hoped that this is related to the symmetry of
the Dynkin diagram of \(\mathrm{E}_6\), and Allen replied:

\begin{quote}
It does. The automorphism \(G\to G\) that exchanges representations with
their duals, the Cartan involution, may or may not be an inner
automorphism. The group of outer automorphisms of \(G\) (\(G\) simple)
is iso to the diagram automorphism group. So no diagram auts, means the
Cartan involution is inner, means all reps are iso to their duals,
i.e.~possess invariant bilinear forms.
\end{quote}

\begin{quote}
(Unfortunately it's not iff --- the \(D_n\)'s alternate between whether
the Cartan involution is inner, much as their centers alternate between
\(\mathbb{Z}_4\) and \(\mathbb{Z}_2^2\).)
\end{quote}

\begin{quote}
Any rep either has complex trace sometimes, or a real, or a quaternionic
structure, morally because of Artin-Wedderburn applied to the real group
algebra of \(G\). Given a rep one can find out which by looking at the
integral over \(G\) of \(\operatorname{Tr}(g^2)\), which comes out
\(0\), \(1\), or \(-1\) (respectively). This is the ``Schur indicator''
and can be found in Serre's LinReps of Finite Groups.

Allen K.
\end{quote}

\hypertarget{week132}{%
\section{April 2, 1999}\label{week132}}

Today I want to talk about \(n\)-categories and quantum gravity again.
For starters let me quote from a paper of mine about this stuff:

\begin{enumerate}
\def\labelenumi{\arabic{enumi})}
\tightlist
\item
  John Baez, ``Higher-dimensional algebra and Planck-scale physics'', in
   \emph{Physics Meets Philosophy at the Planck Scale},
  eds.~Craig Callender and Nick Huggett, Cambridge U.\ Press, Cambridge,
  2001, pp.\ 177--195.   Also available as
  \href{https://arxiv.org/abs/gr-qc/9902017}{\texttt{gr-qc/9902017}}.
\end{enumerate}
\noindent
By the way, this book should be pretty fun to read --- it'll contain
papers by both philosophers and physicists, including a bunch who have
already graced the pages of This Week's Finds, like Barbour, Isham,
Rovelli, Unruh, and Witten. I'll say more about it when it comes out.

Okay, here are some snippets from this paper. It starts out talking
about the meaning of the Planck length, why it may be important in
quantum gravity, and what a theory of quantum gravity should be like:

\begin{quote}
Two constants appear throughout general relativity: the speed of light
\(c\) and Newton's gravitational constant \(G\). This should be no
surprise, since Einstein created general relativity to reconcile the
success of Newton's theory of gravity, based on instantaneous action at
a distance, with his new theory of special relativity, in which no
influence travels faster than light. The constant \(c\) also appears in
quantum field theory, but paired with a different partner: Planck's
constant \(\hbar\). The reason is that quantum field theory takes into
account special relativity and quantum theory, in which \(\hbar\) sets
the scale at which the uncertainty principle becomes important.

It is reasonable to suspect that any theory reconciling general
relativity and quantum theory will involve all three constants \(c\),
\(G\), and \(\hbar\). Planck noted that apart from numerical factors
there is a unique way to use these constants to define units of length,
time, and mass. For example, we can define the unit of length now called
the `Planck length' as follows: \[L = \sqrt{\frac{\hbar G}{c^3}}\] This
is extremely small: about \(1.6\times 10^{-35}\) meters. Physicists have
long suspected that quantum gravity will become important for
understanding physics at about this scale. The reason is very simple:
any calculation that predicts a length using only the constants \(c\),
\(G\) and \(\hbar\) must give the Planck length, possibly multiplied by
an unimportant numerical factor like \(2\pi\).

For example, quantum field theory says that associated to any mass \(m\)
there is a length called its Compton wavelength, \(L_C\), such that
determining the position of a particle of mass \(m\) to within one
Compton wavelength requires enough energy to create another particle of
that mass. Particle creation is a quintessentially
quantum-field-theoretic phenomenon. Thus we may say that the Compton
wavelength sets the distance scale at which quantum field theory becomes
crucial for understanding the behavior of a particle of a given mass. On
the other hand, general relativity says that associated to any mass
\(m\) there is a length called the Schwarzschild radius, \(L_S\), such
that compressing an object of mass \(m\) to a size smaller than this
results in the formation of a black hole. The Schwarzschild radius is
roughly the distance scale at which general relativity becomes crucial
for understanding the behavior of an object of a given mass. Now,
ignoring some numerical factors, we have \[L_C = \frac{\hbar}{mc}\] and
\[L_S = \frac{Gm}{c^2}.\] These two lengths become equal when \(m\) is
the Planck mass. And when this happens, they both equal the Planck
length!

At least naively, we thus expect that both general relativity and
quantum field theory would be needed to understand the behavior of an
object whose mass is about the Planck mass and whose radius is about the
Planck length. This not only explains some of the importance of the
Planck scale, but also some of the difficulties in obtaining
experimental evidence about physics at this scale. Most of our
information about general relativity comes from observing heavy objects
like planets and stars, for which \(L_S\gg L_C\). Most of our
information about quantum field theory comes from observing light
objects like electrons and protons, for which \(L_C\gg L_S\). The Planck
mass is intermediate between these: about the mass of a largish cell.
But the Planck length is about \(10^{-20}\) times the radius of a
proton! To study a situation where both general relativity and quantum
field theory are important, we could try to compress a cell to a size
\(10^{-20}\) times that of a proton. We know no reason why this is
impossible in principle. But we have no idea how to actually accomplish
such a feat.

There are some well-known loopholes in the above argument. The
`unimportant numerical factor' I mentioned above might actually be very
large, or very small. A theory of quantum gravity might make testable
predictions of dimensionless quantities like the ratio of the muon and
electron masses. For that matter, a theory of quantum gravity might
involve physical constants other than \(c\), \(G\), and \(\hbar\). The
latter two alternatives are especially plausible if we study quantum
gravity as part of a larger theory describing other forces and
particles. However, even though we cannot prove that the Planck length
is significant for quantum gravity, I think we can glean some wisdom
from pondering the constants \(c\), \(G\), and \(\hbar\) --- and more
importantly, the physical insights that lead us to regard these
constants as important.

What is the importance of the constant \(c\)? In special relativity,
what matters is the appearance of this constant in the Minkowski metric
\[ds^2 = c^2 dt^2 - dx^2 - dy^2 - dz^2\] which defines the geometry of
spacetime, and in particular the lightcone through each point. Stepping
back from the specific formalism here, we can see several ideas at work.
First, space and time form a unified whole which can be thought of
geometrically. Second, the quantities whose values we seek to predict
are localized. That is, we can measure them in small regions of
spacetime (sometimes idealized as points). Physicists call such
quantities `local degrees of freedom'. And third, to predict the value
of a quantity that can be measured in some region \(R\), we only need to
use values of quantities measured in regions that stand in a certain
geometrical relation to \(R\). This relation is called the `causal
structure' of spacetime. For example, in a relativistic field theory, to
predict the value of the fields in some region \(R\), it suffices to use
their values in any other region that intersects every timelike path
passing through \(R\). The common way of summarizing this idea is to say
that nothing travels faster than light. I prefer to say that a good
theory of physics should have \emph{local degrees of freedom propagating
causally}.

In Newtonian gravity, \(G\) is simply the strength of the gravitational
field. It takes on a deeper significance in general relativity, where
the gravitational field is described in terms of the curvature of the
spacetime metric. Unlike in special relativity, where the Minkowski
metric is a `background structure' given a priori, in general relativity
the metric is treated as a field which not only affects, but also is
affected by, the other fields present. In other words, the geometry of
spacetime becomes a local degree of freedom of the theory.
Quantitatively, the interaction of the metric and other fields is
described by Einstein's equation \[G_{ab} = 8 \pi G T_{ab}\] where the
Einstein tensor \(G_{ab}\) depends on the curvature of the metric, while
the stress-energy tensor \(T_{ab}\) describes the flow of energy and
momentum due to all the other fields. The role of the constant \(G\) is
thus simply to quantify how much the geometry of spacetime is affected
by other fields. Over the years, people have realized that the great
lesson of general relativity is that a good theory of physics should
contain no geometrical structures that affect local degrees of freedom
while remaining unaffected by them. Instead, all geometrical structures
--- and in particular the causal structure --- should themselves be
local degrees of freedom. For short, one says that the theory should be
background-free.

The struggle to free ourselves from background structures began long
before Einstein developed general relativity, and is still not complete.
The conflict between Ptolemaic and Copernican cosmologies, the dispute
between Newton and Leibniz concerning absolute and relative motion, and
the modern arguments concerning the `problem of time' in quantum gravity
--- all are but chapters in the story of this struggle. I do not have
room to sketch this story here, nor even to make more precise the
all-important notion of `geometrical structure'. I can only point the
reader towards the literature, starting perhaps with the books by
Barbour and Earman, various papers by Rovelli, and the many references
therein.

Finally, what of \(\hbar\)? In quantum theory, this appears most
prominently in the commutation relation between the momentum \(p\) and
position \(q\) of a particle: \[pq - qp = -i \hbar,\] together with
similar commutation relations involving other pairs of measurable
quantities. Because our ability to measure two quantities simultaneously
with complete precision is limited by their failure to commute,
\(\hbar\) quantifies our inability to simultaneously know everything one
might choose to know about the world. But there is far more to quantum
theory than the uncertainty principle. In practice, \(\hbar\) comes
along with the whole formalism of complex Hilbert spaces and linear
operators.

There is a widespread sense that the principles behind quantum theory
are poorly understood compared to those of general relativity. This has
led to many discussions about interpretational issues. However, I do not
think that quantum theory will lose its mystery through such
discussions. I believe the real challenge is to better understand why
the mathematical formalism of quantum theory is precisely what it is.
Research in quantum logic has done a wonderful job of understanding the
field of candidates from which the particular formalism we use has been
chosen. But what is so special about this particular choice? Why, for
example, do we use complex Hilbert spaces rather than real or
quaternionic ones? Is this decision made solely to fit the experimental
data, or is there a deeper reason? Since questions like this do not yet
have clear answers, I shall summarize the physical insight behind
\(\hbar\) by saying simply that a good theory of the physical universe
should be a \emph{quantum theory} --- leaving open the possibility of
eventually saying something more illuminating.

Having attempted to extract the ideas lying behind the constants \(c\),
\(G\), and \(\hbar\), we are in a better position to understand the task
of constructing a theory of quantum gravity. General relativity
acknowledges the importance of \(c\) and \(G\) but idealizes reality by
treating \(\hbar\) as negligibly small. From our discussion above, we
see that this is because general relativity is a background-free
classical theory with local degrees of freedom propagating causally. On
the other hand, quantum field theory as normally practiced acknowledges
\(c\) and \(\hbar\) but treats \(G\) as negligible, because it is a
background-dependent quantum theory with local degrees of freedom
propagating causally.

The most conservative approach to quantum gravity is to seek a theory
that combines the best features of general relativity and quantum field
theory. To do this, we must try to find a \emph{background-free quantum
theory with local degrees of freedom propagating causally}. While this
approach may not succeed, it is definitely worth pursuing. Given the
lack of experimental evidence that would point us towards fundamentally
new principles, we should do our best to understand the full
implications of the principles we already have!

From my description of the goal one can perhaps see some of the
difficulties. Since quantum gravity should be background-free, the
geometrical structures defining the causal structure of spacetime should
themselves be local degrees of freedom propagating causally. This much
is already true in general relativity. But because quantum gravity
should be a quantum theory, these degrees of freedom should be treated
quantum-mechanically. So at the very least, we should develop a quantum
theory of some sort of geometrical structure that can define a causal
structure on spacetime.
\end{quote}

Then I talk about topological quantum field theories, which are
background-free quantum theories \emph{without} local degrees of
freedom, and what we have learned from them. Basically what we've
learned is that there's a deep analogy between the mathematics of
spacetime (e.g. differential topology) and the mathematics of quantum
theory. This is interesting because in background-free quantum theories
we expect that spacetime, instead of serving as a ``stage'' which
events play out, actually becomes part of the play of events itself ---
and must itself be described using quantum theory. So it's very
interesting to try to connect the concepts of spacetime and quantum
theory. The analogy goes like this:

\begin{longtable}[]{@{}ll@{}}
\toprule
\begin{minipage}[b]{0.56\columnwidth}\raggedright
Differential topology\strut
\end{minipage} & \begin{minipage}[b]{0.38\columnwidth}\raggedright
Quantum theory\strut
\end{minipage}\tabularnewline
\midrule
\endhead
\begin{minipage}[t]{0.56\columnwidth}\raggedright
\((n-1)\)-dimensional manifold (space)\strut
\end{minipage} & \begin{minipage}[t]{0.38\columnwidth}\raggedright
Hilbert space (states)\strut
\end{minipage}\tabularnewline
\begin{minipage}[t]{0.56\columnwidth}\raggedright
cobordism between \((n-1)\)-dimensional manifolds (spacetime)\strut
\end{minipage} & \begin{minipage}[t]{0.38\columnwidth}\raggedright
operator (process)\strut
\end{minipage}\tabularnewline
\begin{minipage}[t]{0.56\columnwidth}\raggedright
composition of cobordisms\strut
\end{minipage} & \begin{minipage}[t]{0.38\columnwidth}\raggedright
composition of operators\strut
\end{minipage}\tabularnewline
\begin{minipage}[t]{0.56\columnwidth}\raggedright
identity cobordism\strut
\end{minipage} & \begin{minipage}[t]{0.38\columnwidth}\raggedright
identity operator\strut
\end{minipage}\tabularnewline
\bottomrule
\end{longtable}

And if you know a little category theory, you'll see what we have here
are two categories: the category of cobordisms and the category of
Hilbert spaces. A topological quantum field theory is a functor from the
first to the second\ldots.

Okay, now for some other papers:

\begin{enumerate}
\def\labelenumi{\arabic{enumi})}
\setcounter{enumi}{1}
\item
  Geraldine Brady and Todd H. Trimble, ``A string diagram calculus for
  predicate logic'', and C. S. Peirce's system Beta, available at
  \texttt{http://people.cs.uchicago.edu/\textasciitilde{}brady}

  Geraldine Brady and Todd H. Trimble, ``A categorical interpretation of
  Peirce's propositional logic Alpha'', \emph{Jour.\ Pure Appl.\ Alg.}
  \textbf{149} (2000), 213--239.

  Geraldine Brady and Todd H. Trimble, ``The topology of relational
  calculus''.
\end{enumerate}

Charles Peirce is a famously underappreciated American philosopher who
worked in the late 1800s. Among other things, like being the father of
pragmatism, he is also one of the fathers of higher-dimensional algebra.
As you surely know if you've read me often enough, part of the point of
higher-dimensional algebra is to break out of ``linear thinking''. By
``linear thinking'' I mean the tendency to do mathematics in ways that
are easily expressed in terms of \(1\)-dimensional strings of symbols.
In his work on logic, Peirce burst free into higher dimensions. He
developed a way of reasoning using diagrams that he called ``existential
graphs''. Unfortunately this work by Peirce was never published! One
reason is that existential graphs were difficult and expensive to print.
As a result, his ideas languished in obscurity.

By now it's clear that higher-dimensional algebra is useful in physics:
examples include Feynman diagrams and the spin networks of Penrose. The
theory of \(n\)-categories is beginning to provide a systematic language
for all these techniques. So it's worth re-evaluating Peirce's work and
seeing how it fits into the picture. And this is what the papers by
Brady and Trimble do!

\begin{enumerate}
\def\labelenumi{\arabic{enumi})}
\setcounter{enumi}{2}
\tightlist
\item
  J.\ Scott Carter, Louis H.\ Kauffman, and Masahico Saito, ``Structures
  and diagrammatics of four dimensional topological lattice field
  theories'',  \emph{Adv. Math.}\textbf{146} (1999), 39--100.  Also available as
  \href{https://arxiv.org/abs/math.GT/9806023}{\texttt{math.GT/9806023}}.
\end{enumerate}
\noindent
We can get \(3\)-dimensional topological quantum field theories from
certain Hopf algebras. As I described in
\protect\hyperlink{week38}{``Week 38''}, Crane and Frenkel made the
suggestion that by categorifying this construction we should get
\(4\)-dimensional TQFTs from certain Hopf categories. This paper makes
the suggestion precise in a certain class of examples! Basically these
are categorified versions of the Dijkgraaf--Witten theory.

\begin{enumerate}
\def\labelenumi{\arabic{enumi})}
\setcounter{enumi}{3}
\tightlist
\item
  J. Scott Carter, Daniel Jelsovsky, Selichi Kamada, Laurel Langford and
  Masahico Saito, ``Quandle cohomology and state-sum invariants of
  knotted curves and surfaces'', available as
  \href{https://arxiv.org/abs/math.GT/9903135}{\texttt{math.GT/9903135}}.
\end{enumerate}
\noindent
Yet another attack on higher dimensions! This one gets invariants of
2-links --- surfaces embedded in \(\mathbb{R}^4\) --- from the
cohomology groups of ``quandles''. I don't really understand how this
fits into the overall scheme of higher-dimensional algebra yet. They
show their invariant distinguishes between the 2-twist spun trefoil (a
certain sphere knotted in \(\mathbb{R}^4\) and the same sphere with the reversed
orientation.

\begin{enumerate}
\def\labelenumi{\arabic{enumi})}
\setcounter{enumi}{4}
\tightlist
\item
  Tom Leinster, ``Structures in higher-dimensional category theory'',
  available as  hfill \break \href{https://arxiv.org/abs/math/0109021}{\texttt{math/0109021}}.
\end{enumerate}

This is a nice tour of ideas in higher-dimensional algebra. Right now
one big problem with the subject is that there are lots of approaches
and not a clear enough picture of how they fit together. Leinster's
paper is an attempt to start seeing how things fit together.

\begin{enumerate}
\def\labelenumi{\arabic{enumi})}
\setcounter{enumi}{5}
\tightlist
\item
  Claudio Hermida, ``Higher-dimensional multicategories'', slides of a
  lecture given in 1997.
\end{enumerate}
\noindent
This talk presents some of the work by Makkai, Power and Hermida on
their definition of \(n\)-categories. For more on their work see
\protect\hyperlink{week107}{``Week 107''}.

\begin{enumerate}
\def\labelenumi{\arabic{enumi})}
\setcounter{enumi}{6}
\tightlist
\item
  Carlos Simpson, ``On the Breen--Baez--Dolan stabilization hypothesis for
  Tamsamani's weak \(n\)-categories'', available as
  \href{https://arxiv.org/abs/math.CT/9810058}{\texttt{math.CT/9810058}}.
\end{enumerate}
\noindent
For quite a while now James Dolan and I have been talking about
something we call the ``stabilization hypothesis''. I gave an
explanation of this in \protect\hyperlink{week121}{``Week 121''}, but
briefly, it says that the \(n\)th column of the following chart (which
extends to infinity in both directions) stabilizes after \(2n+2\) rows:

\vfill
\newpage

\begin{longtable}[]{@{}llll@{}}
\caption*{\(k\)-tuply monoidal \(n\)-categories}\tabularnewline
\toprule
\begin{minipage}[b]{0.26\columnwidth}\raggedright
\strut
\end{minipage} & \begin{minipage}[b]{0.21\columnwidth}\raggedright
\(n=0\)\strut
\end{minipage} & \begin{minipage}[b]{0.21\columnwidth}\raggedright
\(n=1\)\strut
\end{minipage} & \begin{minipage}[b]{0.21\columnwidth}\raggedright
\(n=2\)\strut
\end{minipage}\tabularnewline
\midrule
\endfirsthead
\toprule
\begin{minipage}[b]{0.26\columnwidth}\raggedright
\strut
\end{minipage} & \begin{minipage}[b]{0.21\columnwidth}\raggedright
\(n=0\)\strut
\end{minipage} & \begin{minipage}[b]{0.21\columnwidth}\raggedright
\(n=1\)\strut
\end{minipage} & \begin{minipage}[b]{0.21\columnwidth}\raggedright
\(n=2\)\strut
\end{minipage}\tabularnewline
\midrule
\endhead
\begin{minipage}[t]{0.26\columnwidth}\raggedright
\(k=0\)\strut
\end{minipage} & \begin{minipage}[t]{0.21\columnwidth}\raggedright
sets\strut
\end{minipage} & \begin{minipage}[t]{0.21\columnwidth}\raggedright
categories\strut
\end{minipage} & \begin{minipage}[t]{0.21\columnwidth}\raggedright
\(2\)-categories\strut
\end{minipage}\tabularnewline
\begin{minipage}[t]{0.26\columnwidth}\raggedright
\strut
\end{minipage} & \begin{minipage}[t]{0.21\columnwidth}\raggedright
\strut
\end{minipage} & \begin{minipage}[t]{0.21\columnwidth}\raggedright
\strut
\end{minipage} & \begin{minipage}[t]{0.21\columnwidth}\raggedright
\strut
\end{minipage}\tabularnewline
\begin{minipage}[t]{0.26\columnwidth}\raggedright
\(k=1\)\strut
\end{minipage} & \begin{minipage}[t]{0.21\columnwidth}\raggedright
monoids\strut
\end{minipage} & \begin{minipage}[t]{0.21\columnwidth}\raggedright
monoidal categories\strut
\end{minipage} & \begin{minipage}[t]{0.21\columnwidth}\raggedright
monoidal \(2\)-categories\strut
\end{minipage}\tabularnewline
\begin{minipage}[t]{0.26\columnwidth}\raggedright
\strut
\end{minipage} & \begin{minipage}[t]{0.21\columnwidth}\raggedright
\strut
\end{minipage} & \begin{minipage}[t]{0.21\columnwidth}\raggedright
\strut
\end{minipage} & \begin{minipage}[t]{0.21\columnwidth}\raggedright
\strut
\end{minipage}\tabularnewline
\begin{minipage}[t]{0.26\columnwidth}\raggedright
\(k=2\)\strut
\end{minipage} & \begin{minipage}[t]{0.21\columnwidth}\raggedright
commutative monoids\strut
\end{minipage} & \begin{minipage}[t]{0.21\columnwidth}\raggedright
braided monoidal categories\strut
\end{minipage} & \begin{minipage}[t]{0.21\columnwidth}\raggedright
braided monoidal \(2\)-categories\strut
\end{minipage}\tabularnewline
\begin{minipage}[t]{0.26\columnwidth}\raggedright
\strut
\end{minipage} & \begin{minipage}[t]{0.21\columnwidth}\raggedright
\strut
\end{minipage} & \begin{minipage}[t]{0.21\columnwidth}\raggedright
\strut
\end{minipage} & \begin{minipage}[t]{0.21\columnwidth}\raggedright
\strut
\end{minipage}\tabularnewline
\begin{minipage}[t]{0.26\columnwidth}\raggedright
\(k=3\)\strut
\end{minipage} & \begin{minipage}[t]{0.21\columnwidth}\raggedright
`` "\strut
\end{minipage} & \begin{minipage}[t]{0.21\columnwidth}\raggedright
symmetric monoidal categories\strut
\end{minipage} & \begin{minipage}[t]{0.21\columnwidth}\raggedright
weakly involutory monoidal \(2\)-categories\strut
\end{minipage}\tabularnewline
\begin{minipage}[t]{0.26\columnwidth}\raggedright
\strut
\end{minipage} & \begin{minipage}[t]{0.21\columnwidth}\raggedright
\strut
\end{minipage} & \begin{minipage}[t]{0.21\columnwidth}\raggedright
\strut
\end{minipage} & \begin{minipage}[t]{0.21\columnwidth}\raggedright
\strut
\end{minipage}\tabularnewline
\begin{minipage}[t]{0.26\columnwidth}\raggedright
\(k=4\)\strut
\end{minipage} & \begin{minipage}[t]{0.21\columnwidth}\raggedright
`` "\strut
\end{minipage} & \begin{minipage}[t]{0.21\columnwidth}\raggedright
`` "\strut
\end{minipage} & \begin{minipage}[t]{0.21\columnwidth}\raggedright
strongly involutory monoidal \(2\)-categories\strut
\end{minipage}\tabularnewline
\begin{minipage}[t]{0.26\columnwidth}\raggedright
\strut
\end{minipage} & \begin{minipage}[t]{0.21\columnwidth}\raggedright
\strut
\end{minipage} & \begin{minipage}[t]{0.21\columnwidth}\raggedright
\strut
\end{minipage} & \begin{minipage}[t]{0.21\columnwidth}\raggedright
\strut
\end{minipage}\tabularnewline
\begin{minipage}[t]{0.26\columnwidth}\raggedright
\(k=5\)\strut
\end{minipage} & \begin{minipage}[t]{0.21\columnwidth}\raggedright
`` "\strut
\end{minipage} & \begin{minipage}[t]{0.21\columnwidth}\raggedright
`` "\strut
\end{minipage} & \begin{minipage}[t]{0.21\columnwidth}\raggedright
`` "\strut
\end{minipage}\tabularnewline
\bottomrule
\end{longtable}

Carlos Simpson has now made this hypothesis precise and proved it using
Tamsamani's definition of \(n\)-categories! And he did it using the same
techniques that Graeme Segal used to study \(k\)-fold loop
spaces\ldots{} exploiting the relation between \(n\)-categories and
homotopy theory. This makes me really happy.

\begin{enumerate}
\def\labelenumi{\arabic{enumi})}
\setcounter{enumi}{7}
\tightlist
\item
  Mark Hovey, \emph{Model Categories}, American Mathematical Society
  Mathematical Surveys and Monographs, vol.~\textbf{63}, Providence,
  Rhode Island, 1999.  Also available as
  \url{http://www.math.uiuc.edu/K-theory/0278/index.html}
\end{enumerate}
\noindent
Speaking of that kind of thing, the technique of model categories is
really important for homotopy theory and \(n\)-categories, and this book
is a really great place to learn about it.

\begin{enumerate}
\def\labelenumi{\arabic{enumi})}
\setcounter{enumi}{8}
\tightlist
\item
  Frank Quinn, ``Group-categories and their field theories'', \emph{Geom. 
 Topol. Monogr.} \textbf{2} (1999), 407--453.  Also available as
  \href{https://arxiv.org/abs/math.GT/9811047}{\texttt{math.GT/9811047}}.
\end{enumerate}
\noindent
This one is about the algebra behind certain topological quantum field
theories. I'll just quote the abstract:

\begin{quote}
A group-category is an additively semisimple category with a monoidal
product structure in which the simple objects are invertible. For
example in the category of representations of a group, \(1\)-dimensional
representations are the invertible simple objects. This paper gives a
detailed exploration of ``topological quantum field theories'' for
group-categories, in hopes of finding clues to a better understanding of
the general situation. Group-categories are classified in several ways
extending results of Fr\"ohlich and Kerler. Topological field theories
based on homology and cohomology are constructed, and these are shown to
include theories obtained from group-categories by Reshetikhin--Turaev
constructions. Braided-commutative categories most naturally give
theories on 4-manifold thickenings of 2-complexes; the usual 3-manifold
theories are obtained from these by normalizing them (using results of
Kirby) to depend mostly on the boundary of the thickening. This is
worked out for group-categories, and in particular we determine when the
normalization is possible and when it is not.
\end{quote}

\begin{enumerate}
\def\labelenumi{\arabic{enumi})}
\setcounter{enumi}{9}
\tightlist
\item
  Sjoerd Crans, ``A tensor product for Gray-categories'', \emph{Theory
  and Applications of Categories}, Vol. \textbf{5}, 1999, No.~2, pp
  12--69, available at
  \href{http://www.tac.mta.ca/tac/volumes/1999/n2/abstract.html}{\texttt{http://www.tac.mta.ca/}}  \href{http://www.tac.mta.ca/tac/volumes/1999/n2/abstract.html}{\texttt{tac/volumes/1999/n2/abstract.html}}
\end{enumerate}
\noindent
A Gray-category is what some people call a semistrict \(3\)-category:
not as general as a weak \(3\)-category, but general enough.
Technically, Gray-categories are defined as categories enriched over the
category of \(2\)-categories equipped with a tensor product invented by
John Gray. To define semistrict \(4\)-categories one might similarly try
to equip Gray-categories with a suitable tensor product. And this is
what Crans is studying. Let me quote the abstract:

\begin{quote}
In this paper I extend Gray's tensor product of \(2\)-categories to a
new tensor product of Gray-categories. I give a description in terms of
generators and relations, one of the relations being an `'interchange''
relation, and a description similar to Gray's description of his tensor
product of \(2\)-categories. I show that this tensor product of
Gray-categories satisfies a universal property with respect to
quasi-functors of two variables, which are defined in terms of
lax-natural transformations between Gray-categories. The main result is
that this tensor product is part of a monoidal structure on Gray-Cat,
the proof requiring interchange in an essential way. However, this does
not give a monoidal \{(bi)closed\} structure, precisely because of
interchange And although I define composition of lax-natural
transformations, this composite need not be a lax-natural transformation
again, making Gray-Cat only a partial Gray-Cat-cateegory.
\end{quote}

\hypertarget{week133}{%
\section{April 23, 1999}\label{week133}}

I'd like to start with a long quote from a paper by Ashtekar:

\begin{enumerate}
\def\labelenumi{\arabic{enumi})}
\tightlist
\item
  Abhay Ashtekar, ``Quantum mechanics of geometry'', available
  as \href{https://arxiv.org/abs/gr-qc/9901023}{\texttt{gr-qc/9901023}}.
\end{enumerate}

\begin{quote}
During his Goettingen inaugural address in 1854, Riemann suggested that
the geometry of space may be more than just a fiducial, mathematical
entity serving as a passive stage for physical phenomena, and may in
fact have direct physical meaning in its own right. General relativity
provided a brilliant confirmation of this vision: curvature of space now
encodes the physical gravitational field. This shift is profound. To
bring out the contrast, let me recall the situation in Newtonian
physics. There, space forms an inert arena on which the dynamics of
physical systems --- such as the solar system --- unfolds. It is like a
stage, an unchanging backdrop for all of physics. In general relativity,
by contrast, the situation is very different. Einstein's equations tell
us that matter curves space. Geometry is no longer immune to change. It
reacts to matter. It is dynamical. It has ``physical degrees of
freedom'' in its own right. In general relativity, the stage disappears
and joins the troupe of actors! Geometry is a physical entity, very much
like matter.

Now, the physics of this century has shown us that matter has
constituents and the \(3\)-dimensional objects we perceive as solids are
in fact made of atoms. The continuum description of matter is an
approximation which succeeds brilliantly in the macroscopic regime but
fails hopelessly at the atomic scale. It is therefore natural to ask: Is
the same true of geometry? If so, what is the analog of the `atomic
scale?' We know that a quantum theory of geometry should contain three
fundamental constants of Nature, \(c\), \(G\), \(\hbar\), the speed of
light, Newton's gravitational constant and Planck's constant. Now, as
Planck pointed out in his celebrated paper that marks the beginning of
quantum mechanics, there is a unique combination,
\[L = \sqrt{\frac{\hbar G}{c^3}},\] of these constants which has
dimension of length. (\(L \sim 10^{-33}\,\mathrm{cm}\).) It is now
called the Planck length. Experience has taught us that the presence of
a distinguished scale in a physical theory often marks a potential
transition; physics below the scale can be very different from that
above the scale. Now, all of our well-tested physics occurs at length
scales much bigger than \(L\). In this regime, the continuum picture
works well. A key question then is: Will it break down at the Planck
length? Does geometry have constituents at this scale? If so, what are
its atoms? Its elementary excitations? Is the space-time continuum only
a `coarse-grained' approximation? Is geometry quantized? If so, what is
the nature of its quanta?

To probe such issues, it is natural to look for hints in the procedures
that have been successful in describing matter. Let us begin by asking
what we mean by quantization of physical quantities. Take a simple
example --- the hydrogen atom. In this case, the answer is clear: while
the basic observables --- energy and angular momentum --- take on a
continuous range of values classically, in quantum mechanics their
eigenvalues are discrete; they are quantized. So, we can ask if the same
is true of geometry. Classical geometrical quantities such as lengths,
areas and volumes can take on continuous values on the phase space of
general relativity. Are the eigenvalues of corresponding quantum
operators discrete? If so, we would say that geometry is quantized and
the precise eigenvalues and eigenvectors of geometric operators would
reveal its detailed microscopic properties.

Thus, it is rather easy to pose the basic questions in a precise
fashion. Indeed, they could have been formulated soon after the advent
of quantum mechanics. Answering them, on the other hand, has proved to
be surprisingly difficult. The main reason, I believe, is the inadequacy
of standard techniques. More precisely, to examine the microscopic
structure of geometry, we must treat Einstein gravity quantum
mechanically, i.e., construct at least the basics of a quantum theory of
the gravitational field. Now, in the traditional approaches to quantum
field theory, one \emph{begins} with a continuum, background geometry.
To probe the nature of quantum geometry, on the other hand, we should
\emph{not} begin by assuming the validity of this picture. We must let
quantum gravity decide whether this picture is adequate; the theory
itself should lead us to the correct microscopic model of geometry.

With this general philosophy, in this article I will summarize the
picture of quantum geometry that has emerged from a specific approach to
quantum gravity. This approach is non-perturbative. In perturbative
approaches, one generally begins by assuming that space-time geometry is
flat and incorporates gravity --- and hence curvature --- step by step
by adding up small corrections. Discreteness is then hard to unravel.

{[}Footnote: The situation can be illustrated by a harmonic oscillator:
While the exact energy levels of the oscillator are discrete, it would
be very difficult to ``see'' this discreteness if one began with a free
particle whose energy levels are continuous and then tried to
incorporate the effects of the oscillator potential step by step via
perturbation theory.{]}

In the non-perturbative approach, by contrast, there is no background
metric at all. All we have is a bare manifold to start with. All fields
--- matter as well as gravity/geometry --- are treated as dynamical from
the beginning. Consequently, the description can not refer to a
background metric. Technically this means that the full diffeomorphism
group of the manifold is respected; the theory is generally covariant.

As we will see, this fact leads one to Hilbert spaces of quantum states
which are quite different from the familiar Fock spaces of particle
physics. Now gravitons --- the three dimensional wavy undulations on a
flat metric --- do not represent fundamental excitations. Rather, the
fundamental excitations are \emph{one} dimensional. Microscopically,
geometry is rather like a polymer. Recall that, although polymers are
intrinsically one dimensional, when densely packed in suitable
configurations they can exhibit properties of a three dimensional
system. Similarly, the familiar continuum picture of geometry arises as
an approximation: one can regard the fundamental excitations as `quantum
threads' with which one can `weave' continuum geometries. That is, the
continuum picture arises upon coarse-graining of the semi-classical
`weave states'. Gravitons are no longer the fundamental mediators of the
gravitational interaction. They now arise only as approximate notions.
They represent perturbations of weave states and mediate the
gravitational force only in the semi-classical approximation. Because
the non-perturbative states are polymer-like, geometrical observables
turn out to have discrete spectra. They provide a rather detailed
picture of quantum geometry from which physical predictions can be made.

The article is divided into two parts. In the first, I will indicate how
one can reformulate general relativity so that it resembles gauge
theories. This formulation provides the starting point for the quantum
theory. In particular, the one-dimensional excitations of geometry arise
as the analogs of ``Wilson loops'' which are themselves analogs of the
line integrals \(\exp(i\smallint A\cdot dl)\) of electromagnetism. In
the second part, I will indicate how this description leads us to a
quantum theory of geometry. I will focus on area operators and show how
the detailed information about the eigenvalues of these operators has
interesting physical consequences, e.g., to the process of Hawking
evaporation of black holes.
\end{quote}

I feel like quoting more, but I'll resist. It's a nice semi-technical
introduction to loop quantum gravity --- a very good place to start if
you know some math and physics but are just getting started on the
quantum gravity business.

Next, here are some papers by younger folks working on loop quantum
gravity:

\begin{enumerate}
\def\labelenumi{\arabic{enumi})}
\setcounter{enumi}{1}
\item
  Fotini Markopoulou, ``The internal description of a causal set: What
  the universe looks like from the inside'', \emph{Comm.\  Math.\ Phys.}
   \textbf{211} (2000), 559--583.  Also available as
  \href{https://arxiv.org/abs/gr-qc/9811053}{\texttt{gr-qc/9811053}}.

  Fotini Markopoulou, ``Quantum causal histories'', \emph{Class.\ Quant.\
  Grav.} \textbf{17} (2000), 2059--2072. Also available
  as
  \href{https://arxiv.org/abs/hep-th/9904009}{\texttt{hep-th/9904009}}.
\end{enumerate}
\noindent
Fotini Markopoulou is perhaps the first person to take the issue of
causality really seriously in loop quantum gravity. In her earlier work
with Lee Smolin (see \protect\hyperlink{week99}{``Week 99''} and
\protect\hyperlink{week114}{``Week 114''}) she proposed a way to equip
an evolving spin network (or what I'd call a spin foam) with a partial
order on its vertices, representing a causal structure. In these papers
she is further developing these ideas. The first one uses topos theory!
It's good to see brave young physicists who aren't scared of using a
little category theory here and there to make their ideas precise.
Personally I feel confused about causality in loop quantum gravity --- I
think we'll have to muck around and try different things before we find
out what works. But Markopoulou's work is the main reason I'm even
\emph{daring} to think about these issues\ldots.

\begin{enumerate}
\def\labelenumi{\arabic{enumi})}
\setcounter{enumi}{2}
\tightlist
\item
  Seth A. Major, ``Embedded graph invariants in Chern--Simons theory'',
 \emph{Nucl.\ Phys. B} \textbf{550} (1999), 531--560.  Also
  available as
  \href{https://arxiv.org/abs/hep-th/9810071}{\texttt{hep-th/9810071}}.
\end{enumerate}
\noindent
In This Week's Finds I've already mentioned Seth Major has worked with
Lee Smolin on \(q\)-deformed spin networks in quantum gravity (see
\protect\hyperlink{week72}{``Week 72''}). There is a fair amount of
evidence, though as yet no firm proof, that \(q\)-deforming your spin
networks corresponds to introducing a nonzero cosmological constant. The
main technical problem with \(q\)-deformed spin networks is that they
require a ``framing'' of the underlying graph. Here Major tackles that
problem\ldots.

And now for something completely different, arising from a thread on
sci.physics.research started by Garrett Lisi. What's the gauge group of
the Standard Model? Everyone will tell you it's
\(\mathrm{U}(1)\times\mathrm{SU}(2)\times\mathrm{SU}(3)\), but as Marc
Bellon pointed out, this is perhaps not the most accurate answer. Let me
explain why and figure out a better answer.

Every particle in the Standard Model transforms according to some
representation of
\(\mathrm{U}(1)\times\mathrm{SU}(2)\times\mathrm{SU}(3)\), but some
elements of this group act trivially on all these representations. Thus
we can find a smaller group which can equally well be used as the gauge
group of the Standard Model: the quotient of
\(\mathrm{U}(1)\times\mathrm{SU}(2)\times\mathrm{SU}(3)\) by the
subgroup of elements that act trivially.

Let's figure out this subgroup! To do so we need to go through all the
particles and figure out which elements of
\(\mathrm{U}(1)\times\mathrm{SU}(2)\times\mathrm{SU}(3)\) act trivially
on all of them.

Start with the gauge bosons. In any gauge theory, the gauge bosons
transform in the adjoint representation, so the elements of the gauge
group that act trivially are precisely those in the \emph{center} of the
group. \(\mathrm{U}(1)\) is abelian so its center is all of
\(\mathrm{U}(1)\). Elements of \(\mathrm{SU}(n)\) that lie in the center
must be diagonal. The \(n\times n\) diagonal unitary matrices with
determinant \(1\) are all of the form \(\exp(2\pi i/n)\), and these form
a subgroup isomorphic to \(\mathbb{Z}/n\). It follows that the center of
\(\mathrm{U}(1)\times\mathrm{SU}(2)\times\mathrm{SU}(3)\) is
\(\mathrm{U}(1)\times\mathbb{Z}/2\times\mathbb{Z}/3\).

Next let's look at the other particles. If you forget how these work,
see \protect\hyperlink{week119}{``Week 119''}. For the fermions, it
suffices to look at those of the first generation, since the other two
generations transform exactly the same way. First of all, we have the
left-handed electron and neutrino: \[(\nu_L, \mathrm{e}_L)\] These form
a \(2\)-dimensional representation. This representation is the tensor
product of the irreducible rep of \(\mathrm{U}(1)\) with hypercharge
\(-1\), the isospin-\(1/2\) rep of \(\mathrm{SU}(2)\), and the trivial
rep of \(\mathrm{SU}(3)\).

A word about notation! People usually describe irreducible reps of
\(\mathrm{U}(1)\) by integers. For historical reasons, hypercharge comes
in integral multiples of \(1/3\). Thus to get the appropriate integer we
need to multiply the hypercharge by \(3\). Also, the group
\(\mathrm{SU}(2)\) here is associated, not to spin in the sense of
angular momentum, but to something called ``weak isospin''. That's why
we say ``isospin-\(1/2\) rep'' above. Mathematically, though, this is
just the usual spin-\(1/2\) representation of \(\mathrm{SU}(2)\).

Next we have the left-handed up and down quarks, which come in 3 colors
each:
\[(\mathrm{u}_L, \mathrm{u}_L, \mathrm{u}_L, \mathrm{d}_L, \mathrm{d}_L, \mathrm{d}_L)\]
This \(6\)-dimensional representation is the tensor product of the
irreducible rep of \(\mathrm{U}(1)\) with hypercharge \(1/3\), the
isospin-\(1/2\) rep of \(\mathrm{SU}(2)\), and the fundamental rep of
\(\mathrm{SU}(3)\).

That's all the left-handed fermions. Note that they all transform
transform according to the isospin-\(1/2\) rep of \(\mathrm{SU}(2)\) ---
we call them ``isospin doublets''. The right-handed fermions all
transform according to the isospin-\(0\) rep of \(\mathrm{SU}(2)\) ---
they're ``isospin singlets''. First we have the right-handed electron:
\[\mathrm{e}_R\] This is the tensor product of the irreducible rep of
\(\mathrm{U}(1)\) with hypercharge \(-2\), the isospin-\(0\) rep of
\(\mathrm{SU}(2)\), and the trivial rep of \(\mathrm{SU}(3)\). Then
there are the right-handed up quarks:
\[(\mathrm{u}_R, \mathrm{u}_R, \mathrm{u}_R)\] which form the tensor
product of the irreducible rep of \(\mathrm{U}(1)\) with hypercharge
\(4/3\), the isospin-\(0\) rep of \(\mathrm{SU}(2)\), and the
fundamental rep of \(\mathrm{SU}(3)\). And then there are the
right-handed down quarks: \[(\mathrm{d}_R, \mathrm{d}_R, \mathrm{d}_R)\]
which form the tensor product of the irreducible rep of
\(\mathrm{U}(1)\) with hypercharge \(2/3\), the isospin-\(0\) rep of
\(\mathrm{SU}(2)\), and the \(3\)-dimensional fundamental rep of
\(\mathrm{SU}(3)\).

Finally, besides the fermions, there is the --- so far unseen --- Higgs
boson: \[(\mathrm{H}_+, \mathrm{H}_0)\] This transforms according to the
tensor product of the irreducible rep of \(\mathrm{U}(1)\) with
hypercharge \(1\), the isospin-\(1/2\) rep of \(\mathrm{SU}(2)\), and
the 1-dimensional trivial rep of \(\mathrm{SU}(3)\).

Okay, let's see which elements of
\(\mathrm{U}(1)\times\mathbb{Z}/2\times\mathbb{Z}/3\) act trivially on
all these representations! Note first that the generator of
\(\mathbb{Z}/2\) acts as multiplication by \(1\) on the isospin singlets
and \(-1\) on the isospin doublets. Similarly, the generator of
\(\mathbb{Z}/3\) acts as multiplication by \(1\) on the leptons and
\(\exp(2\pi i/3)\) on the quarks. Thus everything in
\(\mathbb{Z}/2\times\mathbb{Z}/3\) acts as multiplication by some sixth
root of unity. So to find elements of
\(\mathrm{U}(1)\times\mathbb{Z}/2\times\mathbb{Z}/3\) that act
trivially, we only need to consider guys in \(\mathrm{U}(1)\) that are
sixth roots of unity.

To see what's going on, we make a little table using the information
I've described:

\begin{longtable}[]{@{}lccc@{}}
\toprule
\begin{minipage}[b]{0.03\columnwidth}\raggedright
\strut
\end{minipage} & \begin{minipage}[b]{0.30\columnwidth}\centering
Action of \(\exp(\pi i/3)\) in \(\mathrm{U}(1)\)\strut
\end{minipage} & \begin{minipage}[b]{0.23\columnwidth}\centering
Action of \(-1\) in \(\mathrm{SU}(2)\)\strut
\end{minipage} & \begin{minipage}[b]{0.33\columnwidth}\centering
Action of \(\exp(2\pi i/3)\) in \(\mathrm{SU}(3)\)\strut
\end{minipage}\tabularnewline
\midrule
\endhead
\begin{minipage}[t]{0.03\columnwidth}\raggedright
\(\mathrm{e}_L\)\strut
\end{minipage} & \begin{minipage}[t]{0.30\columnwidth}\centering
\(-1\)\strut
\end{minipage} & \begin{minipage}[t]{0.23\columnwidth}\centering
\(-1\)\strut
\end{minipage} & \begin{minipage}[t]{0.33\columnwidth}\centering
\(1\)\strut
\end{minipage}\tabularnewline
\begin{minipage}[t]{0.03\columnwidth}\raggedright
\(\nu_L\)\strut
\end{minipage} & \begin{minipage}[t]{0.30\columnwidth}\centering
\(-1\)\strut
\end{minipage} & \begin{minipage}[t]{0.23\columnwidth}\centering
\(-1\)\strut
\end{minipage} & \begin{minipage}[t]{0.33\columnwidth}\centering
\(1\)\strut
\end{minipage}\tabularnewline
\begin{minipage}[t]{0.03\columnwidth}\raggedright
\(\mathrm{u}_L\)\strut
\end{minipage} & \begin{minipage}[t]{0.30\columnwidth}\centering
\(\exp(\pi i/3)\)\strut
\end{minipage} & \begin{minipage}[t]{0.23\columnwidth}\centering
\(-1\)\strut
\end{minipage} & \begin{minipage}[t]{0.33\columnwidth}\centering
\(\exp(2\pi i/3)\)\strut
\end{minipage}\tabularnewline
\begin{minipage}[t]{0.03\columnwidth}\raggedright
\(\mathrm{d}_L\)\strut
\end{minipage} & \begin{minipage}[t]{0.30\columnwidth}\centering
\(\exp(\pi i/3)\)\strut
\end{minipage} & \begin{minipage}[t]{0.23\columnwidth}\centering
\(-1\)\strut
\end{minipage} & \begin{minipage}[t]{0.33\columnwidth}\centering
\(\exp(2\pi i/3)\)\strut
\end{minipage}\tabularnewline
\begin{minipage}[t]{0.03\columnwidth}\raggedright
\(\mathrm{e}_R\)\strut
\end{minipage} & \begin{minipage}[t]{0.30\columnwidth}\centering
\(1\)\strut
\end{minipage} & \begin{minipage}[t]{0.23\columnwidth}\centering
\(1\)\strut
\end{minipage} & \begin{minipage}[t]{0.33\columnwidth}\centering
\(1\)\strut
\end{minipage}\tabularnewline
\begin{minipage}[t]{0.03\columnwidth}\raggedright
\(\mathrm{u}_R\)\strut
\end{minipage} & \begin{minipage}[t]{0.30\columnwidth}\centering
\(\exp(4\pi i/3)\)\strut
\end{minipage} & \begin{minipage}[t]{0.23\columnwidth}\centering
\(1\)\strut
\end{minipage} & \begin{minipage}[t]{0.33\columnwidth}\centering
\(\exp(2\pi i/3)\)\strut
\end{minipage}\tabularnewline
\begin{minipage}[t]{0.03\columnwidth}\raggedright
\(\mathrm{d}_R\)\strut
\end{minipage} & \begin{minipage}[t]{0.30\columnwidth}\centering
\(\exp(4\pi i/3)\)\strut
\end{minipage} & \begin{minipage}[t]{0.23\columnwidth}\centering
\(1\)\strut
\end{minipage} & \begin{minipage}[t]{0.33\columnwidth}\centering
\(\exp(2\pi i/3)\)\strut
\end{minipage}\tabularnewline
\begin{minipage}[t]{0.03\columnwidth}\raggedright
H\strut
\end{minipage} & \begin{minipage}[t]{0.30\columnwidth}\centering
\(-1\)\strut
\end{minipage} & \begin{minipage}[t]{0.23\columnwidth}\centering
\(-1\)\strut
\end{minipage} & \begin{minipage}[t]{0.33\columnwidth}\centering
\(1\)\strut
\end{minipage}\tabularnewline
\bottomrule
\end{longtable}

And we look for patterns!

See any?

The most important one for our purposes is that if we multiply all three
numbers in each row, we get \(1\).

This means that the element \((\exp(\pi i/3),-1,\exp(2\pi i/3))\) in
\(\mathrm{U}(1)\times\mathrm{SU}(2)\times\mathrm{SU}(3)\) acts trivially
on all particles. This element generates a subgroup isomorphic to
\(\mathbb{Z}/6\). If you think a bit harder you'll see there are no
\emph{other} patterns that would make any \emph{more} elements of
\(\mathrm{U}(1)\times\mathrm{SU}(2)\times\mathrm{SU}(3)\) act trivially.
And if you think about the relation between charge and hypercharge,
you'll see this pattern has a lot to do with the fact that quark charges
in multiples of \(1/3\), while leptons have integral charge. There's
more to it than that, though\ldots.

Anyway, the ``true'' gauge group of the Standard Model --- i.e., the
smallest possible one --- is not
\(\mathrm{U}(1)\times\mathrm{SU}(2)\times\mathrm{SU}(3)\), but the
quotient of this by the particular \(\mathbb{Z}/6\) subgroup we've just
found. Let's call this group \(G\).

There are two reasons why this might be important. First, Marc Bellon
pointed out a nice way to think about \(G\): it's the subgroup of
\(\mathrm{U}(2)\times\mathrm{U}(3)\) consisting of elements \((g,h)\)
with \[(\operatorname{det} g)(\operatorname{det} h) = 1.\] If we embed
\(\mathrm{U}(2)\times\mathrm{U}(3)\) in \(\mathrm{U}(5)\) in the obvious
way, then this subgroup \(G\) actually lies in \(\mathrm{SU}(5)\),
thanks to the above equation. And this is what people do in the
\(\mathrm{SU}(5)\) grand unified theory. They don't actually stuff all
of \(\mathrm{U}(1)\times\mathrm{SU}(2)\times\mathrm{SU}(3)\) into
\(\mathrm{SU}(5)\), just the group \(G\)! For more details, see
\protect\hyperlink{week119}{``Week 119''}. Better yet, try this book
that Brett McInnes recommended to me:

\begin{enumerate}
\def\labelenumi{\arabic{enumi})}
\setcounter{enumi}{3}
\tightlist
\item
  Lochlainn O'Raifeartaigh, \emph{Group Structure of Gauge Theories},
  Cambridge U.\ Press, Cambridge, 1986.
\end{enumerate}

Second, this magical group \(G\) has a nice action on a
\(7\)-dimensional manifold which we can use as the fiber for a
\(11\)-dimensional Kaluza--Klein theory that mimics the Standard Model in
the low-energy limit. The way to get this manifold is to take
\(S^3\times S^5\) sitting inside \(C^2\times C^3\) and mod out by the
action of \(\mathrm{U}(1)\) as multiplication by phases. The group \(G\)
acts on \(C^2\times C^3\) in an obvious way, and using this it's easy to
see that it acts on \((C^2\times C^3)/\mathrm{U}(1)\).

I'm not sure where to read more about this, but you might try:

\begin{enumerate}
\def\labelenumi{\arabic{enumi})}
\setcounter{enumi}{4}
\item
  Edward Witten, ``Search for a realistic Kaluza--Klein theory'',
  \emph{Nucl. Phys.} \textbf{B186} (1981), 412--428.

  Edward Witten, ``Fermion quantum numbers in Kaluza--Klein theory'', in
  \emph{Shelter Island II, Proceedings: Quantum Field Theory and the
  Fundamental Problems of Physics}, ed.~T. Appelquist et al, MIT Press,
  1985, pp.~227--277.
\item
  Thomas Appelquist, Alan Chodos and Peter G.O. Freund, editors,
  \emph{Modern Kaluza--Klein Theories}, Addison-Wesley, Menlo Park,
  California, 1987.
\end{enumerate}

\hypertarget{week134}{%
\section{June 8, 1999}\label{week134}}

My production of ``This Week's Finds'' has slowed to a trickle as I've
been struggling to write up a bunch of papers. Deadlines, deadlines! I
hate deadlines, but when you write things for other people, or with
other people, that's what you get into. I'll do my best to avoid them in
the future. Now I'm done with my chores and I want to have some fun.

I spent last weekend with a bunch of people talking about quantum
gravity in a hunting lodge by a lake in Minnowbrook, New York:

\begin{enumerate}
\def\labelenumi{\arabic{enumi})}
\tightlist
\item
  \emph{Minnowbrook Symposium on Space-Time Structure}, program and
  transparencies of talks available at
  \href{https://web.archive.org/web/19991008165301/http://www.phy.syr.edu/research/he_theory/minnowbrook/}{\texttt{https://web.archive.org/web/19991008165301/http://www.}}
\href{https://web.archive.org/web/19991008165301/http://www.phy.syr.edu/research/he_theory/minnowbrook/}{\texttt{phy.syr.edu/research/he\_theory/minnowbrook/}}
\end{enumerate}
\noindent
The idea of this get-together, organized by Kameshwar Wali and some
other physicists at Syracuse University, was to bring together people
working on string theory, loop quantum gravity, noncommutative geometry,
and various discrete approaches to spacetime. People from these
different schools of thought don't talk to each other as much as they
should, so this was a good idea. People gave lots of talks, asked lots
of tough questions, argued, and learned what each other were doing. But
I came away with a sense that we're quite far from understanding quantum
gravity: every approach has obvious flaws.

One big problem with string theory is that people only know how to study
it on a spacetime with a fixed background metric. Even worse, things are
poorly understood except when the metric is static --- that is, roughly
speaking, when geometry of space does not change with the passage of
time.

For example, people understand a lot about string theory on spacetimes
that are the product of Minkowski spacetime and a fixed Calabi--Yau
manifold. There are lots of Calabi--Yau manifolds, organized in
continuous multi-parameter families called moduli spaces. This suggests
the idea that the geometry of the Calabi--Yau manifold could change with
time. This idea is lurking behind a lot of interesting work. For
example, Brian Greene gave a nice talk on ``mirror symmetry''. Different
Calabi--Yau manifolds sometimes give the same physics; these are called
``mirror manifolds''. Because of this, a curve in one moduli space of
Calabi--Yau manifolds can be physically equivalent to a curve in some
other moduli space, which sometimes lets you continue the curve beyond a
singularity in the first moduli space. Physicists like to think of these
curves as representing spacetime geometries where the Calabi--Yau
manifold changes with time. The problem is, there's no fully worked out
version of string theory that allows for a time-dependent Calabi--Yau
manifold!

There's a good reason for this: one shouldn't expect anything so simple
to make sense, except in the ``adiabatic approximation'' where things
change very slowly with time. The product of Minkowski spacetime with a
fixed Calabi--Yau manifold is a solution of the \(10\)-dimensional
Einstein equations, and this is part of why this kind of spacetime
serves as a good background for string theory. But we do not get a
solution if the geometry of the Calabi--Yau manifold varies from point to
point in Minkowski spacetime --- except in the adiabatic approximation.

There are also problems with ``unitarity'' in string theory when the
geometry of space changes with time. This is already familiar from
ordinary quantum field theory on curved spacetime. In quantum field
theory, people usually like to describe time evolution using unitary
operators on a Hilbert space of states. But this approach breaks down
when the geometry of space changes with time. People have studied this
problem in detail, and there seems to be no completely satisfactory way
to get around it. No way, that is, except the radical step of ceasing to
treat the geometry of spacetime as a fixed ``background''. In other
words: stop doing quantum field theory on spacetime with a
pre-established metric, and invent a background-free theory of quantum
gravity! But this is not so easy --- see
\protect\hyperlink{week132}{``Week 132''} for more on what it would
entail.

Apparently this issue is coming to the attention of string theorists now
that they are trying to study their theory on non-static background
metrics, such as anti-de Sitter spacetime. Indeed, someone at the
conference said that a bunch of top string theorists recently got
together to hammer out a strategy for where string theory should go
next, but they got completely stuck due to this problem. I think this is
good: it means string theorists are starting to take the foundational
issues of quantum gravity more seriously. These issues are deep and
difficult.

However, lest I seem to be picking on string theory unduly, I should
immediately add that all the other approaches have equally serious
flaws. For example, loop quantum gravity is wonderfully background-free,
but so far it is almost solely a theory of kinematics, rather than
dynamics. In other words, it provides a description of the geometry of
\emph{space} at the quantum level, but says little about
\emph{spacetime}. Recently people have begun to study dynamics with the
help of ``spin foams'', but we still can't compute anything well enough
to be sure we're on the right track. So, pessimistically speaking, it's
possible that the background-free quality of loop quantum gravity has
only been achieved by simplifying assumptions that will later prevent us
from understanding dynamics.

Alain Connes expressed this worry during Abhay Ashtekar's talk, as did
Arthur Jaffe afterwards. Technically speaking, the main issue is that
loop quantum gravity assumes that unsmeared Wilson loops are sensible
observables at the kinematical level, while in other theories, like
Yang--Mills theory, one always needs to smear the Wilson loops. Of course
these other theories aren't background-free, so loop quantum gravity
probably \emph{should} be different. But until we know that loop quantum
gravity really gives gravity (or some fancier theory like supergravity)
in the large-scale limit, we can't be sure it should be different in
this particular way. It's a legitimate worry\ldots{} but only time will
tell!

I could continue listing approaches and their flaws, including Connes'
own approach using noncommutative geometry, but let me stop here. The
only really good news is that different approaches have \emph{different}
flaws. Thus, by comparing them, one might learn something!

Some more papers have come out recently which delve into the
philosophical aspects of this muddle:

\begin{enumerate}
\def\labelenumi{\arabic{enumi})}
\setcounter{enumi}{1}
\item
  Carlo Rovelli, ``Quantum spacetime: what do we know?'', in
  \emph{Physics Meets Philosophy at the Planck Scale}, eds.~Craig
  Callender and Nick Huggett, Cambridge U.\ Press, Cambridge, 2001.  
  Also available as
  \href{https://arxiv.org/abs/gr-qc/9903045}{\texttt{gr-qc/9903045}}.
\item
  J. Butterfield and C. J. Isham, ``Spacetime and the philosophical
  challenge of quantum gravity'', in \emph{Physics Meets
  Philosophy at the Planck Scale}, eds.~Craig Callender and Nick
  Huggett, Cambridge U.\ Press, Cambridge, 2001. Preprint available as
  \href{https://arxiv.org/abs/gr-qc/9903072}{\texttt{gr-qc/9903072}}.
\end{enumerate}

Rovelli's paper is a bit sketchy, but it outlines ideas which I find
very appealing --- I always find him to be very clear-headed about the
conceptual issues of quantum gravity. I found the latter paper a bit
frustrating, because it lays out a wide variety of possible positions
regarding quantum gravity, but doesn't make a commitment to any one of
them. However, this is probably good when one is writing to an audience
of philosophers: one should explain the problems instead of trying to
sell them on a particular claimed solution, because the proposed
solutions come and go rather rapidly, while the problems remain. Let me
quote the abstract:

\begin{quote}
We survey some philosophical aspects of the search for a quantum theory
of gravity, emphasising how quantum gravity throws into doubt the
treatment of spacetime common to the two `ingredient theories' (quantum
theory and general relativity), as a \(4\)-dimensional manifold equipped
with a Lorentzian metric. After an introduction, we briefly review the
conceptual problems of the ingredient theories and introduce the
enterprise of quantum gravity. We then describe how three main research
programmes in quantum gravity treat four topics of particular
importance: the scope of standard quantum theory; the nature of
spacetime; spacetime diffeomorphisms, and the so-called problem of time.
By and large, these programmes accept most of the ingredient theories'
treatment of spacetime, albeit with a metric with some type of quantum
nature; but they also suggest that the treatment has fundamental
limitations. This prompts the idea of going further: either by
quantizing structures other than the metric, such as the topology; or by
regarding such structures as phenomenological. We discuss this in
Section 5.
\end{quote}

Now let me mention a few more technical papers that have come out in the
last few months:

\begin{enumerate}
\def\labelenumi{\arabic{enumi})}
\setcounter{enumi}{3}
\tightlist
\item
  John Baez and John W. Barrett, ``The quantum tetrahedron in 3 and 4
  dimensions'', \emph{Adv. Theor. Math. Phys.} \textbf{3} (1999), 815--850.
   Also available as
  \href{https://arxiv.org/abs/gr-qc/gr-qc/9903060}{\texttt{gr-qc/9903060}}.
\end{enumerate}
\noindent
The idea here is to form a classical phase whose points represent
geometries of a tetrahedron in 3 or 4 dimensions, and then apply
geometric quantization to obtain a Hilbert space of states. These
Hilbert spaces play an important role in spin foam models of quantum
gravity. The main goal of the paper is to explain why the quantum
tetrahedron has fewer degrees of freedom in 4 dimensions than in 3
dimensions. Let me quote from the introduction:

\begin{quote}
State sum models for quantum field theories are constructed by giving
amplitudes for the simplexes in a triangulated manifold. The simplexes
are labelled with data from some discrete set, and the amplitudes depend
on this labelling. The amplitudes are then summed over this set of
labellings, to give a discrete version of a path integral. When the
discrete set is a finite set, then the sum always exists, so this
procedure provides a bona fide definition of the path integral.

State sum models for quantum gravity have been proposed based on the Lie
algebra \(\mathfrak{so}(3)\) and its \(q\)-deformation. Part of the
labelling scheme is then to assign irreducible representations of this
Lie algebra to simplexes of the appropriate dimension. Using the
\(q\)-deformation, the set of irreducible representations becomes
finite. However, we will consider the undeformed case here as the
geometry is more elementary.

Irreducible representations of \(\mathfrak{so}(3)\) are indexed by a
non-negative half-integers \(j\) called spins. The spins have different
interpretations in different models. In the Ponzano--Regge model of
\(3\)-dimensional quantum gravity, spins label the edges of a
triangulated 3-manifold, and are interpreted as the quantized lengths of
these edges. In the Ooguri--Crane--Yetter state sum model, spins label
triangles of a triangulated 4-manifold, and the spin is interpreted as
the norm of a component of the B-field in a BF Lagrangian. There is also
a state sum model of \(4\)-dimensional quantum gravity in which spins
label triangles. Here the spins are interpreted as areas.

Many of these constructions have a topologically dual formulation. The
dual 1-skeleton of a triangulated surface is a trivalent graph, each of
whose edges intersect exactly one edge in the original triangulation.
The spin labels can be thought of as labelling the edges of this graph,
thus defining a spin network. In the Ponzano--Regge model, transition
amplitudes between spin networks can be computed as a sum over
labellings of faces of the dual 2-skeleton of a triangulated 3-manifold.
Formulated this way, we call the theory a `spin foam model'.

A similar dual picture exists for \(4\)-dimensional quantum gravity. The
dual 1-skeleton of a triangulated 3-manifold is a 4-valent graph each of
whose edges intersect one triangle in the original triangulation. The
labels on the triangles in the 3-manifold can thus be thought of as
labelling the edges of this graph. The graph is then called a
`relativistic spin network'. Transition amplitudes between relativistic
spin networks can be computed using a spin foam model. The path integral
is then a sum over labellings of faces of a 2-complex interpolating
between two relativistic spin networks.

In this paper we consider the nature of the quantized geometry of a
tetrahedron which occurs in some of these models, and its relation to
the phase space of geometries of a classical tetrahedron in 3 or 4
dimensions. Our main goal is to solve the following puzzle: why does the
quantum tetrahedron have fewer degrees of freedom in 4 dimensions than
in 3 dimensions? This seeming paradox turns out to have a simple
explanation in terms of geometric quantization. The picture we develop
is that the four face areas of a quantum tetrahedron in four dimensions
can be freely specified, but that the remaining parameters cannot, due
to the uncertainty principle.
\end{quote}

Naively one would expect the quantum tetrahedron to have the same number
of degrees of freedom in 3 and 4 dimensions (since one is considering
tetrahedra mod rotations). However, quantum mechanics is funny about
these things! For example, the Hilbert space of two spin-\(1/2\)
particles whose angular momenta point in opposite directions is smaller
than the Hilbert space of a single spin-\(1/2\) particle, even though
classically you might think both systems have the same number of degrees
of freedom. In fact a very similar thing happens for the quantum
tetrahedron in 3 and 4 dimensions.

\begin{enumerate}
\def\labelenumi{\arabic{enumi})}
\setcounter{enumi}{4}
\tightlist
\item
  Abhay Ashtekar, Alejandro Corichi and Kirill Krasnov, ``Isolated
  horizons: the classical phase space'', \emph{Adv. Theor. Math. Phys.} \textbf{3} 
   (1999), 419--478.  Also available as
  \href{https://arxiv.org/abs/gr-qc/9905089}{\texttt{gr-qc/9905089}}.
\end{enumerate}

This paper explains in more detail the classical aspects of the
calculation of the entropy of a black hole in loop quantum gravity (see
\protect\hyperlink{week112}{``Week 112''} for a description of this
calculation). Let me quote the abstract:

\begin{quote}
A Hamiltonian framework is introduced to encompass non-rotating (but
possibly charged) black holes that are ``isolated'' near future
time-like infinity or for a finite time interval. The underlying
space-times need not admit a stationary Killing field even in a
neighborhood of the horizon; rather, the physical assumption is that
neither matter fields nor gravitational radiation fall across the
portion of the horizon under consideration. A precise notion of
non-rotating isolated horizons is formulated to capture these ideas.
With these boundary conditions, the gravitational action fails to be
differentiable unless a boundary term is added at the horizon. The
required term turns out to be precisely the Chern--Simons action for the
self-dual connection. The resulting symplectic structure also acquires,
in addition to the usual volume piece, a surface term which is the
Chern--Simons symplectic structure. We show that these modifications
affect in subtle but important ways the standard discussion of
constraints, gauge and dynamics. In companion papers, this framework
serves as the point of departure for quantization, a statistical
mechanical calculation of black hole entropy and a derivation of laws of
black hole mechanics, generalized to isolated horizons. It may also have
applications in classical general relativity, particularly in the
investigation of analytic issues that arise in the numerical studies of
black hole collisions.
\end{quote}

The following are some review articles on spin networks, spin foams and
the like:

\begin{enumerate}
\def\labelenumi{\arabic{enumi})}
\setcounter{enumi}{5}
\item
  Roberto De Pietri, ``Canonical `loop' quantum gravity and spin foam
  models", in \emph{Recent Developments in General Relativity} Springer, 
  Berlin, 2000, p.\ 43--61.  Also available as
  \href{https://arxiv.org/abs/gr-qc/9903076}{\texttt{gr-qc/9903076}}.
\item
  Seth Major, ``A spin network primer'',  \emph{Amer. Jour.
  Phys.} \textbf{67} (1999), 972--980.  Available as
  \href{https://arxiv.org/abs/gr-qc/9905020}{\texttt{gr-qc/9905020}}.
\item
  Seth Major, ``Operators for quantized directions'', \emph{Class. Quant.
 Grav.} \textbf{16} (1999), 3859--3877.  Also available
  as \href{https://arxiv.org/abs/gr-qc/9905019}{\texttt{gr-qc/9905019}}.
\item
  John Baez, ``An introduction to spin foam models of \(BF\) theory and
  quantum gravity'', in \emph{Geometry and Quantum Physics}, eds.~Helmut
  Gausterer and Harald Grosse, Lecture Notes in Physics,
  Springer, Berlin, 2000, pp.~25--93. Preprint available as
  \href{https://arxiv.org/abs/gr-qc/9905087}{\texttt{gr-qc/9905087}}.
\end{enumerate}

By the way, Barrett and Crane have come out with a paper sketching a
spin foam model for Lorentzian (as opposed to Riemannian) quantum
gravity:

\begin{enumerate}
\def\labelenumi{\arabic{enumi})}
\setcounter{enumi}{9}
\tightlist
\item
  John Barrett and Louis Crane, ``A Lorentzian signature model for
  quantum general relativity'', \emph{Class. Quant. Grav.} \textbf{17} (2000),     
   3101--3118.  Also available as
  \href{https://arxiv.org/abs/gr-qc/9904025}{\texttt{gr-qc/9904025}}.
\end{enumerate}

However, this model is so far purely formal, because it involves
infinite sums that probably diverge. We need to keep working on this!
Now that I'm getting a bit of free time, I want to tackle this issue.
Meanwhile, Iwasaki has come out with an alternative spin foam model of
Riemannian quantum gravity:

\begin{enumerate}
\def\labelenumi{\arabic{enumi})}
\setcounter{enumi}{10}
\tightlist
\item
  Junichi Iwasaki, ``A surface theoretic model of quantum gravity'',
  available as
  \href{https://arxiv.org/abs/gr-qc/9903112}{\texttt{gr-qc/9903112}}.
\end{enumerate}
\noindent
Alas, I don't really understand this model yet. 

Finally, to wrap things up, something completely different:

\begin{enumerate}
\def\labelenumi{\arabic{enumi})}
\setcounter{enumi}{11}
\tightlist
\item
  Richard E. Borcherds, ``Quantum vertex algebras'', available
  as
  \href{https://arxiv.org/abs/math.QA/9903038}{\texttt{math.QA/9903038}}.
\end{enumerate}
\noindent
I like how the abstract of this paper starts: ``The purpose of this
paper is to make the theory of vertex algebras trivial''. Good! Trivial
is not bad, it's good. Anything one understands is automatically
trivial.

\hypertarget{week135}{%
\section{July 31, 1999}\label{week135}}

Well, darn it, now I'm too busy running around to conferences to write
This Week's Finds! First I went to Vancouver, then to Santa Barbara, and
for almost a month now I've been in Portugal, bouncing between Lisbon
and Coimbra. But let me try to catch up\ldots.

From June 16th to 19th, Steve Savitt and Steve Weinstein of the
University of British Columbia held a workshop designed to get
philosophers and physicists talking about the conceptual problems of
quantum gravity:

\begin{enumerate}
\def\labelenumi{\arabic{enumi})}
\tightlist
\item
  \emph{Toward a New Understanding of Space, Time and Matter}, workshop
  home page at \url{http://axion.physics.ubc.ca/Workshop/}
\end{enumerate}

After a day of lectures by Chris Isham, John Earman, Lee Smolin and
myself, we spent the rest of the workshop sitting around in a big room
with a beautiful view of Vancouver Bay, discussing various issues in a
fairly organized way. For example, Chris Isham led a discussion on
``What is a quantum theory?'' in which he got people to question the
assumptions underlying quantum physics, and Simon Saunders led one on
``Quantum gravity: physics, metaphysics or mathematics?'' in which we
pondered the scientific and sociological implications of the fact that
work on quantum gravity is motivated more by desire for consistency,
clarity and mathematical elegance than the need to fit new experimental
data.

It's pretty clear that understanding quantum gravity will make us
rethink some fundamental concepts --- the question is, which ones? By
the end of the conference, almost every basic belief or concept relevant
to physics had been held up for careful scrutiny and found questionable.
Space, time, causality, the real numbers, set theory --- you name it! It
was a bit unnerving --- but it's good to do this sort of thing now and
then, to prevent hardening of the mental arteries, and it's especially
fun to do it with a big bunch of physicists and philosophers. However, I
must admit that I left wanting nothing more than to do lots of grungy
calculations in order to bring myself back down to earth --- relatively
speaking, of course.

I particularly enjoyed Chris Isham's talk about topos theory because it
helped me understand one way that topos theory could be applied to
quantum theory. I've tended to regard topoi as ``too classical'' for
quantum theory, because while the internal logic of a topos is
intuitionistic (the principle of exclude middle may fail), it's still
not very quantum. For example, in a topos the operation ``and'' still
distributes over ``or'', and vice versa, while failure of this sort of
distributivity is a hallmark of quantum logic. If you don't know what I
mean, try these books, in rough order of increasing difficulty:

\begin{enumerate}
\def\labelenumi{\arabic{enumi})}
\setcounter{enumi}{1}
\item
  David W. Cohen, \emph{An Introduction to Hilbert Space and Quantum
  Logic}, Springer, Berlin, 1989.
\item
  C. Piron, \emph{Foundations of Quantum Physics}, W. A. Benjamin,
  Reading, Massachusetts, 1976.
\item
  C. A. Hooker, editor, \emph{The Logico-algebraic Approach to Quantum
  Mechanics}, two volumes, D. Reidel, Boston, 1975-1979.
\end{enumerate}

Perhaps even more importantly, topoi are Cartesian! What does this mean?
Well, it means that we can define a ``product'' of any two objects in a
topos. That is, given objects \(a\) and \(b\), there's an object
\(a\times b\) equipped with morphisms \[p\colon a\times b\to a\] and
\[q\colon a\times b\to b\] called ``projections'', satisfying the
following property: given morphisms from some object \(c\) into \(a\)
and \(b\), say \[f\colon c\to a\] and \[g\colon c\to b\] there's a
unique morphism \(f\times g\colon c\to a\times b\) such that if we
follow it by \(p\) we get \(f\), and if we follow it by \(q\) we get
\(g\). This is just an abstraction of the properties of the usual
Cartesian product of sets, which is why we call a category ``Cartesian''
if any pair of objects has a product.

Now, it's a fun exercise to check that in a Cartesian category, every
object has a morphism \[\Delta\colon a\to a\times a\] called the
``diagonal'', which when composed with either of the two projections
from \(a\times a\) to a gives the identity. For example, in the topos of
sets, the diagonal morphism is given by \[\Delta(x) = (x,x)\] We can
think of the diagonal morphism as allowing ``duplication of
information''. This is not generally possible in quantum mechanics:

\begin{enumerate}
\def\labelenumi{\arabic{enumi})}
\setcounter{enumi}{4}
\tightlist
\item
  William Wooters and Wocjciech Zurek, ``A single quantum cannot be
  cloned'', \emph{Nature} \textbf{299} (1982), 802--803.
\end{enumerate}

The reason is that in the category of Hilbert spaces, the tensor product
is not a product in the above sense! In particular, given a Hilbert
space \(H\), there isn't a natural diagonal operator
\[\Delta\colon H\to H tensor H\] and there aren't even natural
projection operators from \(H\otimes H\) to \(H\). As pointed out to me
by James Dolan, this non-Cartesianness of the tensor product gives
quantum theory much of its special flavor. Besides making it impossible
to ``clone a quantum'', it's closely related to how quantum theory
violates Bell's inequality, because it means we can't think of an
arbitrary state of a two-part quantum system as built by tensoring
states of both parts.

Anyway, this has made me feel for a while that topos theory isn't
sufficiently ``quantum'' to be useful in understanding the peculiar
special features of quantum physics. However, after Isham and I gave our
talks, someone pointed out to me that one can think of a topological
quantum field theory as a presheaf of Hilbert spaces over the category
\(\mathsf{nCob}\) whose morphisms are \(n\)-dimensional cobordisms. Now,
presheaves over any category form a topos, so this means we should be
able to think of a topological quantum field theory as a ``Hilbert space
object'' in the topos of presheaves over \(\mathsf{nCob}\). From this
point of view, the peculiar ``quantumness'' of topological quantum field
theory comes from it being a Hilbert space object, while its peculiar
``variability'' --- i.e., the fact that it assigns a different Hilbert
space to each \((n-1)\)-dimensional manifold representing space ---
comes from the fact that it's an object in a topos. (Topoi are known for
being very good at handling things like ``variable sets''.) I'm not sure
how useful this is, but it's worth keeping in mind.

While I'm talking about quantum logic, let me raise a puzzle concerning
the Kochen-Specker theorem. Remember what this says: if you have a
Hilbert space \(H\) with dimension more than \(2\), there's no map \(F\)
from self-adjoint operators on \(H\) to real numbers with the following
properties:

\begin{enumerate}
\def\labelenumi{\alph{enumi})}
\tightlist
\item
  For any self-adjoint operator \(A\), \(F(A)\) lies in the spectrum of
  \(A\),
\end{enumerate}

and

\begin{enumerate}
\def\labelenumi{\alph{enumi})}
\setcounter{enumi}{1}
\tightlist
\item
  For any continuous \(f\colon\mathbb{R}\to\mathbb{R}\),
  \(f(F(A)) = F(f(A))\).
\end{enumerate}

This means there's no sensible consistent way of thinking of all
observables as simultaneously having values in a quantum system!

Okay, the puzzle is: what happens if the dimension of \(H\) equals
\(2\)? I don't actually know the answer, so I'd be glad to hear it if
someone can figure it out!

By the way, I once wanted to do an undergraduate research project on
mathematical physics with Kochen. He asked me if I knew the spectral
theorem, I said ``no'', and he said that in that case there was no point
in me trying to work with him. I spent the next summer reading Reed and
Simon's book on Functional Analysis and learning lots of different
versions of the spectral theorem. I shudder to think that perhaps this
is why I spent years studying analysis before eventually drifting
towards topology and algebra. But no: now that I think about it, I was
already interested in analysis at the time, since I'd had a wonderful
real analysis class with Robin Graham.

Okay, now let me say a bit about the next conference I went to. From
June 22nd to 26th there was a conference on ``Strong Gravitational
Fields'' at the Institute for Theoretical Physics at U. C. Santa
Barbara. This finished up a wonderful semester-long program by Abhay
Ashtekar, Gary Horowitz and Jim Isenberg:

\begin{enumerate}
\def\labelenumi{\arabic{enumi})}
\setcounter{enumi}{5}
\tightlist
\item
  \emph{Classical and Quantum Physics of Strong Gravitational Fields},
  program homepage with transparencies and audio files of talks at
  \href{https://online.kitp.ucsb.edu/online/gravity99/}{\texttt{https://online.kitp.ucsb.edu/}} \href{https://online.kitp.ucsb.edu/online/gravity99/}    
   {\texttt{online/gravity99/}}
\end{enumerate}
\noindent
Like the whole program, the conference covered a wide range of topics
related to gravity: string theory and loop quantum gravity,
observational and computational black hole physics, and \(\gamma\)-ray
bursters. I can't summarize all this stuff; since I usually spend a lot
of talking about quantum gravity here, let me say a bit about other
things instead.

John Friedman gave an interesting talk on gravitational waves from
unstable neutron stars. When a pulsar is young, like about 5000 years
old, it typically rotates about its axis once every 16 milliseconds or
so. A good example is N157B, a pulsar in the Large Magellanic Cloud.
Using the current spindown rate one can extrapolate and guess that
pulsars have about a 5-millisecond period at birth. It's interesting to
think about what makes a newly formed neutron star slow down. Theorists
have recently come up with a new possible mechanism: namely, a new sort
of gravitational-wave-driven instability of relativistic stars that
could force newly formed slow down to a 10-millisecond period. It's very
clever: the basic idea is that if a star is rotating very fast, a
rotational mode that rotates slower than the star will gravitationally
radiate \emph{positive} angular momentum, but such modes carry
\emph{negative} angular momentum, since they rotate slower than the
star. If you think about it carefully, you'll see this means that
gravitational radiation should tend to amplify such modes! I asked for a
lowbrow analog of this mechanism and it turns out that a similar sort of
thing is at work in the formation of water waves by the wind --- with
linear momentum taking the place of angular momentum. Anyway, it's not
clear that this process really ever has a chance to happen, because it
only works when the neutron star is not too hot and not too cold, but
it's pretty cool.

Richard Price gave a nice talk on computer simulation of black hole
collisions. Quantitatively understanding the gravitational radiation
emitted in black hole and neutron star collisions is a big business
these days --- it's one of the NSF's ``grand challenge'' problems. The
reason is that folks are spending a lot of money building gravitational
wave detectors like LIGO:

\begin{enumerate}
\def\labelenumi{\arabic{enumi})}
\setcounter{enumi}{6}
\item
  LIGO project home page, \url{https://www.ligo.caltech.edu/}
\end{enumerate}
\noindent
and they need to know exactly what to look for. Now, head-on collisions
are the easiest to understand, since one can simplify the calculation
using axial symmetry. Unfortunately, it's not very likely that two black
holes are going to crash into each other head-on. One really wants to
understand what happens when two black holes spiral into each other.
There are two extreme cases: the case of black holes of equal mass, and
the case of a very light black hole of mass falling into a heavy one.

The latter case is 95\% understood, since we can think of the light
black hole as a ``test particle'' --- ignoring its effect on the heavy
one. The light one slowly spirals into the heavy one until it reaches
the innermost stable orbit, and then falls in. We can use the theory of
a relativistic test particle falling into a black hole to understand the
early stages of this process, and use black hole perturbation theory to
study the ``ringdown'' of the resulting single black hole in the late
stages of the process. (By ``ringdown'' I mean the process whereby an
oscillating black hole settles down while emitting gravitational
radiation.) Even the intermediate stages are manageable, because the
radiation of the small black hole doesn't have much effect on the big
one.

By contrast, the case of two black holes of equal mass is less well
understood. We can treat the early stages, where relativistic effects
are small, using a post-Newtonian approximation, and again we can treat
the late stages using black hole perturbation theory. But things get
complicated in the intermediate stage, because the radiation of each
hole greatly effects the other, and there is no real concept of
``innermost stable orbit'' in this case. To make matters worse, the
intermediate stage of the process is exactly the one we really want to
understand, because this is probably when most of the gravitational
waves are emitted!

People have spent a lot of work trying to understand black hole
collisions through number-crunching computer calculation, but it's not
easy: when you get down to brass tacks, general relativity consists of
some truly scary nonlinear partial differential equations. Current work
is bedeviled by numerical instability and also the problem of simulating
enough of a region of spacetime to understand the gravitational
radiation being emitted. Fans of mathematical physics will also realize
that gauge-fixing is a major problem. There is a lot of interest in
simplifying the calculations through ``black hole excision'': anything
going on inside the event horizon can't affect what happens outside, so
if one can get the computer to \emph{find} the horizon, one can forget
about simulating what's going on inside! But nobody is very good at
doing this yet\ldots{} even using the simpler concept of ``apparent
horizon'', which can be defined locally. So there is some serious work
left to be done!

(For more details on both these talks, go to the conference website and
look at the transparencies.)

I also had some interesting talks with people about black hole entropy,
some of which concerned a new paper by Steve Carlip. I'm not really able
to do justice to the details, but it seems important\ldots.

\begin{enumerate}
\def\labelenumi{\arabic{enumi})}
\setcounter{enumi}{8}
\tightlist
\item
  Steve Carlip, ``Entropy from conformal field theory at Killing
  horizons'', \emph{Class.\ Quant.\ Grav.} \textbf{16} (1999), 3327--3348. 
   Also available as
  \href{https://arxiv.org/abs/gr-qc/9906126}{\texttt{gr-qc/9906126}}.
\end{enumerate}
\noindent
Let me just quote the abstract:

\begin{quote}
On a manifold with boundary, the constraint algebra of general
relativity may acquire a central extension, which can be computed using
covariant phase space techniques. When the boundary is a (local) Killing
horizon, a natural set of boundary conditions leads to a Virasoro
subalgebra with a calculable central charge. Conformal field theory
methods may then be used to determine the density of states at the
boundary. I consider a number of cases --- black holes, Rindler space,
de Sitter space, Taub--NUT and Taub--Bolt spaces, and dilaton gravity ---
and show that the resulting density of states yields the expected
Bekenstein--Hawking entropy. The statistical mechanics of black hole
entropy may thus be fixed by symmetry arguments, independent of details
of quantum gravity.
\end{quote}

There was also a lot of talk about ``isolated horizons'', a concept that
plays a fundamental role in certain treatments of black holes in loop
quantum gravity:

\begin{enumerate}
\def\labelenumi{\arabic{enumi})}
\setcounter{enumi}{9}
\item
  Abhay Ashtekar, Christopher Beetle, and Stephen Fairhurst, ``Mechanics
  of isolated horizons'', \emph{Class. Quant. Grav.} \textbf{17} (2000),
   253--298.  Also available as
  \href{https://arxiv.org/abs/gr-qc/9907068}{\texttt{gr-qc/9907068}}.
\item
  Jerzy Lewandowski, ``Spacetimes admitting isolated horizons'',
  \emph{Class.\ Quant.\ Grav.} \textbf{17} (2000), L53--L59.  Also available as
  \href{https://arxiv.org/abs/gr-qc/9907058}{\texttt{gr-qc/9907058}}.
\end{enumerate}
\noindent
For more on isolated horizons try the references in
\protect\hyperlink{week128}{``Week 128''}.

Finally, on a completely different note, I've recently seen some new
papers related to the McKay correspondence --- see
\protect\hyperlink{week65}{``Week 65''} if you don't know what
\emph{that} is! I haven't read them yet, but I just want to remind
myself that I should, so I'll list them here:

\begin{enumerate}
\def\labelenumi{\arabic{enumi})}
\setcounter{enumi}{11}
\item
  John McKay, ``Semi-affine Coxeter--Dynkin graphs and
  \(G\subseteq\mathrm{SU}_2(\mathbb{C})\)'', available as
  \href{https://arxiv.org/abs/math.QA/9907089}{\texttt{math.QA/9907089}}.
\item
  Igor Frenkel, Naihuan Jing and Weiqiang Wang, ``Vertex representations
  via finite groups and the McKay correspondence'',  	\emph{Internat.\ Math.\ 
   Res.\ Notices} \textbf{4} (2000), 195--222.   Also available as
  \href{https://arxiv.org/abs/math.QA/9907166}{\texttt{math.QA/9907166}}.

   Igor Frenkel, Naihuan Jing and Weiqiang Wang, ``Quantum vertex 
   representations via finite groups and the McKay
  correspondence'', \emph{Comm\. Math.\ Phys.} \textbf{211} (2000), 
   365--393.  Also available as
  \href{https://arxiv.org/abs/math.QA/9907175}{\texttt{math.QA/9907175}}.
\end{enumerate}

Next time I want to talk about the big category theory conference in
honor of Mac Lane's 90th birthday! Then I'll be pretty much caught up on
the conferences\ldots.

\begin{center}\rule{0.5\linewidth}{0.5pt}\end{center}

Robert Israel's answer to my puzzle about the Kochen-Specker theorem:

\begin{quote}
It's not true in dimension \(2\). Note that for a self-adjoint
\(2\times2\) matrix \(A\), any \(f(A)\) is of the form \(a A + b I\) for
some real scalars \(a\) and \(b\) (this is easy to see if you
diagonalize \(A\)). The self-adjoint matrices that are not multiples of
\(I\) split into equivalence classes, where \(A\) and \(B\) are
equivalent if \(B = a A + b I\) for some scalars \(a, b\) (\(a \ne 0\)).
Pick a representative \(A\) from each equivalence class, choose \(F(A)\)
as one of the eigenvalues of \(A\), and then
\(F(a A + b I) = a F(A) + b.\) Of course, \(F(b I) = b\). Then \(F\)
satisfies the two conditions.
\end{quote}

\begin{quote}
The reason this doesn't work in higher dimensions is that in higher
dimensions you can have two self-adjoint matrices \(A\) and \(B\) which
don't commute, \(F(A) = G(B)\) for some functions \(F\) and \(G\), and
\(F(A)\) is not a multiple of \(I\).
\end{quote}

\begin{quote}
Robert Israel \hfill \break 
Department of Mathematics \hfill \break
University of British Columbia Vancouver, BC, Canada V6T 1Z2
\end{quote}

\hypertarget{week136}{%
\section{August 21, 1999}\label{week136}}

I spent most of last month in Portugal, spending time with Roger Picken
at the Instituto Superior Tecnico in Lisbon and attending the category
theory school and conference in Coimbra, which was organized by Manuela
Sobral:

\begin{enumerate}
\def\labelenumi{\arabic{enumi})}
\tightlist
\item
  Category Theory 99 website, with abstracts of talks,
  \url{http://www.mat.uc.pt/~ct99/}
\end{enumerate}
\noindent
The conference was a big deal this year, because it celebrated the 90th
birthday of Saunders Mac Lane, who with Samuel Eilenberg invented
category theory in 1945. Mac Lane was there and in fine fettle. He gave
a nice talk about working with Eilenberg, and after the banquet in his
honor, he even sang a song about Riemann while wrapped in a black cloak!

(In case you're wondering, the cloak was contributed by some musicians.
In Coimbra, the folks who play fado music tend to wear black cloaks. A
few days ago we'd seen them serenade a tearful old man and then wrap him
in a cloak, so one of our number suggested that they try this trick on
Mac Lane. Far from breaking into tears, he burst into song.)

The conference was exquisitely well-organized, packed with top category
theorists, and stuffed with so many cool talks I scarcely know where to
begin describing them\ldots{} I'll probably say a bit about a random
sampling of them next time, and the proceedings will appear in a special
issue of the Journal of Pure and Applied Algebra honoring Mac Lane's
90th birthday, so keep your eye out for that if you're interested. The
school featured courses by Cristina Pedicchio, Vaughan Pratt, and some
crazy mathematical physicist who thinks the laws of physics are based on
\(n\)-categories. The notes can be found in the following book:

\begin{enumerate}
\def\labelenumi{\arabic{enumi})}
\setcounter{enumi}{1}
\tightlist
\item
  \emph{School on Category Theory and Applications}, Coimbra, July
  13-17, 199, Textos de Matematica Serie B No.~\textbf{21}, Departamento
  De Matematica da Universidade de Coimbra. Contains: ``\(n\)-Categories''
  by John Baez, ``Algebraic theories'' by M. Cristina Pedicchio, and
  ``Chu Spaces: duality as a common foundation for computation and
  mathematics'' by Vaughan Pratt.
\end{enumerate}

Pedicchio's course covered various generalizations of Lawvere's
wonderful concept of an algebraic theory. Recall from
\protect\hyperlink{week53}{``Week 53''} that we can think of a category
\(\mathcal{C}\) with extra properties or structure as a kind of
``theory'', and functors \(F\colon \mathcal{C}\to\mathsf{Set}\)
preserving this structure as ``models'' of the theory. For example, a
``finite products theory'' \(\mathcal{C}\) is just a category with
finite products. In this case, a model is a functor
\(F\colon\mathcal{C}\to\mathsf{Set}\) preserving finite products, and a
morphism of models is a natural transformation between such functors.
This gives us a category \(\mathsf{Mod}(\mathcal{C})\) of models of
\(\mathcal{C}\).

To understand what this really means, let's restrict attention the
simplest case, when all the objects in \(\mathcal{C}\) are products of a
given object \(x\). In this case Pedicchio calls \(\mathcal{C}\) an
``algebraic theory''. A model \(F\) is then really just a set \(F(x)\)
together with a bunch of \(n\)-ary operations coming from the morphisms
in \(\mathcal{C}\), satisfying equational laws coming from the equations
between morphisms in \(\mathcal{C}\). Any sort of algebraic gadget
that's just a set with a bunch of \(n\)-ary operations satisfying
equations can be described using a theory of this sort. For example:
monoids, groups, abelian groups, rings\ldots{} and so on. We can
describe any of these using a suitable algebraic theory, and in each
case, the category \(\mathsf{Mod}(\mathcal{C})\) will be the category of
these algebraic gadgets.

Now, what I didn't explain last time I discussed this was the notion of
theory-model duality. Fans of ``duality'' in all its forms are sure to
like this! There's a functor
\[R\colon\mathsf{Mod}(\mathcal{C})\to\mathsf{Set}\] which carries each
model \(F\) to the set \(F(x)\). We can think of this as a functor which
forgets all the operations of our algebraic gadget and remembers only
the underlying set. Now, if you know about adjoint functors (see
\protect\hyperlink{week77}{``Week
77''}--\protect\hyperlink{week79}{``Week 79''}), this should immediately
make you want to find a left adjoint for \(R\), namely a functor
\[L\colon\mathsf{Set}\to\mathsf{Mod}(\mathcal{C})\] sending each set to
the ``free'' algebraic gadget on this set. Indeed, such a left adjoint
exists!

Given this pair of adjoint functors we can do all sorts of fun stuff. In
particular, we can talk about the category of ``finitely generated free
models'' of our theory. The objects here are objects of
\(\mathsf{Mod}(\mathcal{C})\) of the form \(L(S)\) where \(S\) is a
finite set, and the morphisms are the usual morphisms in
\(\mathsf{Mod}(\mathcal{C})\). Let me call this category
\(\mathsf{fgFree}\text{-}\mathsf{Mod}(\mathcal{C})\).

Now for the marvelous duality theorem:
\(\mathsf{fgFree}\text{-}\mathsf{Mod}(\mathcal{C})\) is equivalent to
the opposite of the category \(\mathcal{C}\). In other words, you can
reconstruct an algebraic theory from its category of finitely generated
free algebras in the simplest manner imaginable: just reversing the
direction of all the morphisms! This is so nice I won't explain why it's
true\ldots{} I don't want to deprive you of the pleasure of looking at
some simple examples and seeing for yourself how it works. For example,
take the theory of groups, and figure out how every operation appearing
in the definition of ``group'' corresponds to a homomorphism between
finitely generated free groups.

There are lots of other interesting questions related to theory-model
duality. For example: what kinds of categories arise as categories of
models of an algebraic theory? Pedicchio calls these ``algebraic
categories'', and she told us some nice theorems characterizing them.
Or: given the category of free models of an algebraic theory, can you
fatten it up to get the category of \emph{all} models? Pedicchio
mentioned a process called ``exact completion'' that does the job. Or:
starting from just the category of models of a theory, can you tell
which are the free models? Alas, I don't know the answer to this\ldots{}
but I'm sure people do.

Even better, all of this can be generalized immensely, to theories of a
more flexible sort than the ``algebraic theories'' I've been talking
about so far. For example, we can study ``essentially algebraic
theories'', which are just categories with finite limits. Given one of
these, say \(\mathcal{C}\), we define a model to be a functor
\(F\colon\mathcal{C}\to\mathsf{Set}\) preserving finite limits. This
allows one to study algebraic structures with partially defined
operations. I already gave an example in
\protect\hyperlink{week53}{``Week 53''} --- there's a category with
finite limits called ``the theory of categories'', whose models are
categories! One can work out theory-model duality in this bigger
context, where it's called Gabriel-Ulmer duality:

\begin{enumerate}
\def\labelenumi{\arabic{enumi})}
\setcounter{enumi}{2}
\tightlist
\item
  P. Gabriel and F. Ulmer, \emph{Lokal praesentierbare Kategorien},
  Lecture Notes in Mathematics \textbf{221}, Springer, Berlin, 1971.
\end{enumerate}

But this stuff goes far beyond that, and Pedicchio led us at a rapid
pace all the way up to the latest work. A lot of the basic ideas here
came from Lawvere's famous thesis on algebraic semantics, so it was nice
to see him attending these lectures, and even nicer to hear that 26
years after he wrote it, his thesis is about to be published:

\begin{enumerate}
\def\labelenumi{\arabic{enumi})}
\setcounter{enumi}{3}
\tightlist
\item
  William Lawvere, \emph{Functorial Semantics of Algebraic Theories},
  Ph.D.~thesis, University of Columbia, 1963.   Publlished in 
  \href{Reprints in Theory and Applications of Categories} \textbf{5} (2004), 
   pp. 1--121.  Available at \href{http://www.tac.mta.ca/tac/reprints/articles/5/tr5abs.html}{\text{http://www.tac.mta.ca/tac/reprints/articles/5/tr5abs.html}}
\end{enumerate}

It was also nice to find
out that Lawvere and Schanuel are writing a book on ``objective number
theory''\ldots{} which will presumably be more difficult, but hopefully
not less delightful, than their wonderful introduction to category
theory for people who know \emph{nothing} about fancy mathematics:

\begin{enumerate}
\def\labelenumi{\arabic{enumi})}
\setcounter{enumi}{4}
\tightlist
\item
  William Lawvere and Steve Schanuel, \emph{Conceptual Mathematics: A
  First Introduction to Categories}, Cambridge U.\ Press, Cambridge, 1997.
\end{enumerate}
\noindent
This is the book to give to all your friends who are wondering what
category theory is about and want to learn a bit without too much pain.
If you've read this far and understood what I was talking about, you
must have such friends! If you \emph{didn't} understand what I was
talking about, read this book!

By the way, Lawvere told me that he started out wanting to do physics,
and wound up doing his thesis on algebraic semantics when he started to
trying to formalize what a physical theory was. It's interesting that
the modern notion of ``topological quantum field theory'' is very much
modelled after Lawvere's ideas, but with symmetric monoidal categories
with duals replacing the categories with finite products which Lawvere
considered! I guess he was just ahead of his time. In fact, he has
returned to physics in more recent years - but that's another story.

Okay, let me change gears now\ldots.

Some \(n\)-category gossip. Ross Street has a student who has defined a
notion of semistrict \(n\)-category up to \(n = 5\), and Sjoerd Crans
has defined semistrict \(n\)-categories (which he calls ``teisi'') for
\(n\) up to \(6\). However, the notion still seems to resist definition
for general \(n\), which prompted my pal Lisa Raphals to compose the
following limerick:

\begin{quote}
A theoretician of ``\(n\)'' \hfill\break
Considered conditions on when \hfill\break  
Some mathematicians \hfill\break 
Could find definitions \hfill\break
For \(n\) even greater than ten.  \hfill
\end{quote}

Interestingly, work on weak \(n\)-categories seems to be proceeding at a
slightly faster clip --- they've gotten to \(n = \infty\) already. In
fact, during the conference Michael Batanin came up to me and said that
a fellow named Penon had published a really terse definition of weak
\(\omega\)-categories that seems equivalent to Batanin's own (see
\protect\hyperlink{week103}{``Week 103''}) --- at least after some minor
tweaking. Batanin was quite enthusiastic and said he plans to write a
paper about this stuff.

Later, when I went to Cambridge England, Tom Leinster gave a talk
summarizing Penon's definition:

\begin{enumerate}
\def\labelenumi{\arabic{enumi})}
\setcounter{enumi}{5}
\tightlist
\item
  Jacques Penon, ``Approache polygraphique des \(\infty\)-categories non
  strictes'', in \emph{Cahiers Top. Geom. Diff.} \textbf{40} (1999),
  31--79.
\end{enumerate}
\noindent
It seems pretty cool, so I'd like to tell you what Leinster said ---
using his terminology rather than Penon's (which of course is in
French). To keep this short I'm going to assume you know a reasonable
amount of category theory.

First of all, a ``reflexive globular set'' is a collection of sets and
functions like this: \[
  \begin{tikzcd}
    X_0
        \rar["i" description]
    & X_1
        \lar[shift right=5,"s" description]
        \lar[shift left=5,"t" description]
        \rar["i" description]
    & X_2
        \lar[shift right=5,"s" description]
        \lar[shift left=5,"t" description]
        \rar["i" description]
    & \ldots
        \lar[shift right=5,"s" description]
        \lar[shift left=5,"t" description]
  \end{tikzcd}
\] going on to infinity, satisfying these equations: \[
  \begin{gathered}
    s(s(x)) = s(t(x))
  \\t(s(x)) = t(t(x))
  \\s(i(x)) = t(i(x)) = x.
  \end{gathered}
\] We call the elements of \(X_n\) ``\(n\)-cells'', and call \(s(x)\)
and \(t(x)\) the ``source'' and ``target'' of the \(n\)-cell \(x\),
respectively. \(If s(x) = a\) and \(t(x) = b\), we think of \(x\) as
going from \(a\) to \(b\), and write \(x\colon a\to b\).

If we left out all the stuff about the maps \(i\) we would simply have a
``globular set''. These are important in \(n\)-category theory because
strict \(\omega\)-categories, and also Batanin's weak
\(\omega\)-categories, are globular sets with extra structure. This also
true of Penon's definition, but he starts right away with ``reflexive''
globular sets, which have these maps \(i\) that are a bit like the
degeneracies in the definition of a simplicial set (see
\protect\hyperlink{week115}{``Week 115''}). In Penon's definition
\(i(x)\) plays the role of an ``identity \(n\)-morphism'', so we also
write \(i(x)\) as \(1_x\colon x\to x\).

Let \(\mathsf{RGlob}\) be the category of reflexive globular sets, where
morphisms are defined in the obvious way. (In other words,
\(\mathsf{RGlob}\) is a presheaf category --- see
\protect\hyperlink{week115}{``Week 115''} for an explanation of this
notion.)

In this setup, the usual sort of strict \(\omega\)-category may be
defined as a reflexive globular set \(X\) together with various
``composition'' operations that allow us to compose \(n\)-cells \(x\)
and \(y\) whenever \(t^j(x)=s^j(x)\), obtaining an \(n\)-cell
\[x \circ_j y\] We get one such composition operation for each \(n\) and
each \(j\) such that \(1\leqslant j\leqslant n\). We impose some obvious
axioms of two sorts:

\begin{enumerate}
\def\labelenumi{\Alph{enumi})}
\item
  axioms determining the source and target of a composite; and
\item
  strict associativity, unit and interchange laws.
\end{enumerate}

I'll assume you know these axioms or can fake it. (If you read the
definition of strict \(2\)-category in \protect\hyperlink{week80}{``Week
80''}, perhaps you can get an idea for what kinds of axioms I'm talking
about.)

Now, strict \(\omega\)-categories are great, but we need to weaken this
notion. So, first Penon defines an ``\(\omega\)-magma'' to be something
exactly like a strict \(\omega\)-category but without the axioms of type
B. You may recall that a ``magma'' is defined by Bourbaki to be a set
with a binary operation satisyfing no laws whatsoever --- the primeval
algebraic object! An \(\omega\)-magma is just as lawless, and a lot
bigger and meaner.

Strict \(\omega\)-categories are too strict: all laws hold as equations.
\(\omega\)-magmas are too weak: no laws hold at all! How do we get what
we want?

We define a category \(Q\) whose objects are quadruples
\((M,p,\mathcal{C},[\cdot,\cdot])\) where:

\begin{itemize}
\tightlist
\item
  \(M\) is an \(\omega\)-magma
\item
  \(\mathcal{C}\) is a strict \(\omega\)-category
\item
  \(p\colon M\to\mathcal{C}\) is a morphism of \(\omega\)-magmas (i.e.~a
  morphism of reflexive globular sets strictly preserving all the
  \(\omega\)-magma operations)
\item
  \([\cdot,\cdot]\) is a way of lifting equations between
  \(n\)-morphisms in the image of the projection \(p\) to
  \((n+1)\)-morphisms in \(M\). More precisely: given \(n\)-cells
  \[f,g\colon a\to b\] in \(M\) such that \(p(f) = p(g)\), we have an
  \((n+1)\)-cell \[[f,g]\colon f\to g\] in \(M\) such that
  \(p([f,g]) = 1_{p(f)} = 1_{p(g)}\). We require that \([f,f] = 1_f\).
\end{itemize}

A morphism in \(Q\) is defined to be the obvious thing: a morphism
\(f\colon M\to M'\) of \(\omega\)-magmas and a morphism
\(f\colon\mathcal{C}\to\mathcal{C}'\) of strict \(\omega\)-categories,
strictly preserving all the structure in sight.

Okay, now we define a functor \[U\colon Q\to\mathsf{RGlob}\] by
\[\mathrm{U}(M,p,C,[\cdot,\cdot]) = M\] where we think of \(M\) as just
a reflexive globular set. Penon proves that \(U\) has a left adjoint
\[F\colon\mathsf{RGlob}\to Q\] This adjunction defines a monad
\[T\colon \mathsf{RGlob}\to \mathsf{RGlob}\] and Penon defines a ``weak
\(\omega\)-category'' to be an algebra of this monad.

(See \protect\hyperlink{week92}{``Week 92''} and
\protect\hyperlink{week118}{``Week 118''} for how you get monads from
adjunctions. Alas, I think I haven't gotten around to explaining the
concept of an algebra of a monad! So much to explain, so little time!)

Now, if you know some category theory and think a while about this, you
will see that in a weak \(\omega\)-category defined this way, all the
laws like associativity hold \emph{up to equivalence}, with the
equivalences satisfying the necessary coherence laws \emph{up to
equivalence}, and so ad infinitum. Crudely speaking, the lifting
\([\cdot,\cdot]\) is what turns equations into \(n\)-morphisms. To get a
feeling for how this work, you have to figure out what the left adjoint
\(F\) looks like. Penon works this out in detail in the second half of
his paper.

\hypertarget{week137}{%
\section{September 4, 1999}\label{week137}}

Now I'm in Cambridge England, chilling out with the category theorists,
so it makes sense for me to keep talking about category theory. I'll
start with some things people discussed at the conference in Coimbra
(see last week).

\begin{enumerate}
\def\labelenumi{\arabic{enumi})}
\tightlist
\item
  Michael M\"uger, ``Galois theory for braided tensor categories and the
  modular closure'', \emph{Adv. Math.} \textbf{150} (2000), 151--201.  
   Also available as
  \href{https://arxiv.org/abs/math.CT/9812040}{\texttt{math.CT/9812040}}.
\end{enumerate}
\noindent
A braided monoidal category is simple algebraic gadget that captures a
bit of the essence of \(3\)-dimensionality in its rawest form. It has a
bunch of ``objects'' which we can draw a labelled dots like this: \[
  \begin{tikzpicture}
    \node[label=above:{$x$}] at (0,0) {$\bullet$};
  \end{tikzpicture}
\] So far this is just 0-dimensional. Next, given a bunch of objects we
get a new object, their ``tensor product'', which we can draw by setting
the dots side by side. So, for example, we can draw \(x\otimes y\) like
this: \[
  \begin{tikzpicture}
    \node[label=above:{$x$}] at (0,0) {$\bullet$};
    \node[label=above:{$y$}] at (1,0) {$\bullet$};
  \end{tikzpicture}
\] This is \(1\)-dimensional. But in addition we have, for any pair of
objects \(x\) and \(y\), a bunch of ``morphisms'' \(f\colon x\to y\). We
can draw a morphism from a tensor product of objects to some other
tensor product of objects as a picture like this: \[
  \begin{tikzpicture}
    \node[label=above:{$x$}] at (0,0) {$\bullet$};
    \node[label=above:{$y$}] at (1,0) {$\bullet$};
    \node[label=above:{$z$}] at (2,0) {$\bullet$};
    \draw[thick] (0,0) to (0,-1.5);
    \draw[thick] (1,0) to (1,-1.5);
    \draw[thick] (2,0) to (2,-1.5);
    \draw[thick,rounded corners] (-0.2,-1.5) rectangle ++(2.4,-1);
    \node at (1,-2) {$f$};
    \draw[thick] (0.5,-2.5) to (0.5,-4);
    \draw[thick] (1.5,-2.5) to (1.5,-4);
    \node[label=below:{$u$}] at (0.5,-4) {$\bullet$};
    \node[label=below:{$v$}] at (1.5,-4) {$\bullet$};
  \end{tikzpicture}
\] This picture is \(2\)-dimensional. In addition, we require that for
any pair of objects \(x\) and \(y\) there is a ``braiding'', a special
morphism from \(x\otimes y\) to \(y\otimes x\). We draw it like this: \[
  \begin{tikzpicture}
    \node[label=above:{$x$}] at (0,0) {$\bullet$};
    \node[label=above:{$y$}] at (1,0) {$\bullet$};
    \begin{knot}
      \strand[thick] (1,0)
        to [out=down,in=up] (0,-2);
      \strand[thick] (0,0)
        to [out=down,in=up] (1,-2);
    \end{knot}
    \node[label=below:{$y$}] at (0,-2) {$\bullet$};
    \node[label=below:{$x$}] at (1,-2) {$\bullet$};
  \end{tikzpicture}
\] With this crossing of strands, the picture has become
\(3\)-dimensional!

We also require that we can ``compose'' a morphism \(f\colon x\to y\)
and a morphism \(g\colon y\to z\) and get a morphism
\(fg\colon x\to z\). We draw this by sticking one picture on top of each
other.   I'll draw a fancy example where all the objects
in question are themselves tensor products of other objects: \[
  \begin{tikzpicture}
    \node[label=above:{$\phantom{|}x\phantom{|}$}] at (0,0) {$\bullet$};
    \node[label=above:{$\phantom{|}y\phantom{|}$}] at (1,0) {$\bullet$};
    \node[label=above:{$\phantom{|}z\phantom{|}$}] at (2,0) {$\bullet$};
    \draw[thick] (0,0) to (0,-1);
    \draw[thick] (1,0) to (1,-1);
    \draw[thick] (2,0) to (2,-1);
    \draw[thick,rounded corners] (-0.2,-1) rectangle ++(2.4,-1);
    \node at (1,-1.5) {$f$};
    \draw[thick] (0.5,-2) to (0.5,-3);
    \draw[thick] (1.5,-2) to (1.5,-3);
    \draw[thick,rounded corners] (-0.2,-3) rectangle ++(2.4,-1);
    \node at (1,-3.5) {$g$};
    \draw[thick] (0,-4) to (0,-5);
    \draw[thick] (0.66,-4) to (0.66,-5);
    \draw[thick] (1.33,-4) to (1.33,-5);
    \draw[thick] (2,-4) to (2,-5);
    \node[label=below:{$\phantom{|}a\phantom{|}$}] at (0,-5) {$\bullet$};
    \node[label=below:{$\phantom{|}b\phantom{|}$}] at (0.66,-5) {$\bullet$};
    \node[label=below:{$\phantom{|}c\phantom{|}$}] at (1.33,-5) {$\bullet$};
    \node[label=below:{$\phantom{|}d\phantom{|}$}] at (2,-5) {$\bullet$};
  \end{tikzpicture}
\]

Finally, we require that the tensor product, braiding and composition
satisfy a bunch of axioms. I won't write these down because I already
did so in \protect\hyperlink{week121}{``Week 121''}, but the point is
that they all make geometrical sense --- or more precisely, topological
sense --- in terms of the above pictures.

The pictures I've drawn should make you think about knots and tangles
and circuit diagrams and Feynman diagrams and all sorts of things like
that --- and it's true, you can understand all these things very
elegantly in terms of braided monoidal categories! Sometimes it's nice
to throw in another rule: \[
  \begin{tikzpicture}
    \node[label=above:{$\phantom{|}x\phantom{|}$}] at (0,0) {$\bullet$};
    \node[label=above:{$\phantom{|}y\phantom{|}$}] at (1,0) {$\bullet$};
    \begin{knot}
      \strand[thick] (1,0)
        to [out=down,in=up] (0,-2);
      \strand[thick] (0,0)
        to [out=down,in=up] (1,-2);
    \end{knot}
    \node[label=below:{$\phantom{|}y\phantom{|}$}] at (0,-2) {$\bullet$};
    \node[label=below:{$\phantom{|}x\phantom{|}$}] at (1,-2) {$\bullet$};
    \node at (2,-1) {$=$};
    \begin{scope}[shift={(3,0)}]
      \node[label=above:{$\phantom{|}x\phantom{|}$}] at (0,0) {$\bullet$};
      \node[label=above:{$\phantom{|}y\phantom{|}$}] at (1,0) {$\bullet$};
      \begin{knot}
        \strand[thick] (0,0)
          to [out=down,in=up] (1,-2);
        \strand[thick] (1,0)
          to [out=down,in=up] (0,-2);
      \end{knot}
      \node[label=below:{$\phantom{|}y\phantom{|}$}] at (0,-2) {$\bullet$};
      \node[label=below:{$\phantom{|}x\phantom{|}$}] at (1,-2) {$\bullet$};
    \end{scope}
  \end{tikzpicture}
\] where we cook up the second picture using the inverse of the
braiding. This rule is good when you don't care about the difference
between overcrossings and undercrossings. If this rule holds we say our
braided monoidal category is ``symmetric''. Topologically, this rule
makes sense when we study \(4\)-dimensional or higher-dimensional
situations, where we have enough room to untie all knots. For example,
the traditional theory of Feynman diagrams is based on symmetric
monoidal categories (like the category of representations of the
Poincar\'e group), and it works very smoothly in \(4\)-dimensional
spacetime.

But \(3\)-dimensional spacetime is a bit different. For example, when we
interchange two identical particles, it really makes a difference
whether we do it like this: \[
  \begin{tikzpicture}
    \node[label=above:{$\phantom{|}x\phantom{|}$}] at (0,0) {$\bullet$};
    \node[label=above:{$\phantom{|}y\phantom{|}$}] at (1,0) {$\bullet$};
    \begin{knot}
      \strand[thick] (1,0)
        to [out=down,in=up] (0,-2);
      \strand[thick] (0,0)
        to [out=down,in=up] (1,-2);
    \end{knot}
    \node[label=below:{$\phantom{|}y\phantom{|}$}] at (0,-2) {$\bullet$};
    \node[label=below:{$\phantom{|}x\phantom{|}$}] at (1,-2) {$\bullet$};
  \end{tikzpicture}
\] or like this: \[
  \begin{tikzpicture}
    \node[label=above:{$\phantom{|}x\phantom{|}$}] at (0,0) {$\bullet$};
    \node[label=above:{$\phantom{|}y\phantom{|}$}] at (1,0) {$\bullet$};
    \begin{knot}
      \strand[thick] (0,0)
        to [out=down,in=up] (1,-2);
      \strand[thick] (1,0)
        to [out=down,in=up] (0,-2);
    \end{knot}
    \node[label=below:{$\phantom{|}y\phantom{|}$}] at (0,-2) {$\bullet$};
    \node[label=below:{$\phantom{|}x\phantom{|}$}] at (1,-2) {$\bullet$};
  \end{tikzpicture}
\] 
Thus in 3d spacetime, besides bosons and fermions, we have other
sorts of particles that act differently when we interchange them ---
sometimes people call them ``anyons'', and sometimes people talk about
``exotic statistics''.

Now let me dig into some more technical aspects of the picture.

Starting with Reshetikhin and Turaev, people have figured out how to use
braided monoidal categories to construct topological quantum field
theories in \(3\)-dimensional spacetime. But they can't do it starting
from any old braided monoidal category, because quantum field theory has
a lot to do with Hilbert spaces. So usually they start from a special
sort called a ``modular tensor category''. This is a kind of hybrid of a
braided monoidal category and a Hilbert space.

In fact, apart from one technical condition --- which is at the heart of
M\"uger's work --- we can get the definition of a modular tensor category
by taking the definition of ``Hilbert space'', categorifying it once to
get the definition of ``2-Hilbert space'', and then throwing in a tensor
product and braiding that are compatible with this structure.

It's amazing that by such abstract conceptual methods we come up with
almost precisely what's needed to construct topological quantum field
theories in 3 dimensions! It's a great illustration of the power of
category theory. It's almost like getting something for nothing! But
I'll resist the temptation to tell you the details, since
\protect\hyperlink{week99}{``Week 99''} explains a bunch of it, and the
rest is in here:

\begin{enumerate}
\def\labelenumi{\arabic{enumi})}
\setcounter{enumi}{1}
\tightlist
\item
  John Baez, ``Higher-dimensional algebra II: 2-Hilbert spaces'',
  \emph{Adv.\ Math.} \textbf{127} (1997), 125--189. Also available as
  \href{https://arxiv.org/abs/q-alg/9609018}{\texttt{q-alg/9609018}}.
\end{enumerate}
\noindent
In this paper I call a 2-Hilbert space with a compatible tensor product
a ``2-H*-algebra", and if it also has a compatible braiding, I call
it a ``braided 2-H*-algebra''. This terminology is bit clunky, but for
consistency I'll use it again here.

Okay, great: we \emph{almost} get the definition of modular tensor
category by elegant conceptual methods. But there is one niggling but
crucial technical condition that remains! There are lots of different
ways to state this condition, but M\"uger proves they're equivalent to
the following very elegant one.

Let's define the ``center'' of a braided monoidal category to be the
category consisting of all objects \(x\) such that \[
  \begin{tikzpicture}
    \node[label=above:{$\phantom{|}x\phantom{|}$}] at (0,0) {$\bullet$};
    \node[label=above:{$\phantom{|}y\phantom{|}$}] at (1,0) {$\bullet$};
    \begin{knot}
      \strand[thick] (1,0)
        to [out=down,in=up] (0,-2);
      \strand[thick] (0,0)
        to [out=down,in=up] (1,-2);
    \end{knot}
    \node[label=below:{$\phantom{|}y\phantom{|}$}] at (0,-2) {$\bullet$};
    \node[label=below:{$\phantom{|}x\phantom{|}$}] at (1,-2) {$\bullet$};
    \node at (2,-1) {$=$};
    \begin{scope}[shift={(3,0)}]
      \node[label=above:{$\phantom{|}x\phantom{|}$}] at (0,0) {$\bullet$};
      \node[label=above:{$\phantom{|}y\phantom{|}$}] at (1,0) {$\bullet$};
      \begin{knot}
        \strand[thick] (0,0)
          to [out=down,in=up] (1,-2);
        \strand[thick] (1,0)
          to [out=down,in=up] (0,-2);
      \end{knot}
      \node[label=below:{$y$}] at (0,-2) {$\bullet$};
      \node[label=below:{$x$}] at (1,-2) {$\bullet$};
    \end{scope}
  \end{tikzpicture}
\] for all \(y\), and all morphisms between such objects. The center of
a braided monoidal category is obviously a symmetric monoidal category.
The term ``center'' is supposed to remind you of the usual center of a
monoid --- the elements that commute with all the others. And indeed,
both kinds of center are special cases of a general construction that
pushes you down the columns of the ``periodic table'':

\vfill
\newpage

\begin{longtable}[]{@{}llll@{}}
\caption*{\(k\)-tuply monoidal \(n\)-categories}\tabularnewline
\toprule
\begin{minipage}[b]{0.26\columnwidth}\raggedright
\strut
\end{minipage} & \begin{minipage}[b]{0.21\columnwidth}\raggedright
\(n=0\)\strut
\end{minipage} & \begin{minipage}[b]{0.21\columnwidth}\raggedright
\(n=1\)\strut
\end{minipage} & \begin{minipage}[b]{0.21\columnwidth}\raggedright
\(n=2\)\strut
\end{minipage}\tabularnewline
\midrule
\endfirsthead
\toprule
\begin{minipage}[b]{0.26\columnwidth}\raggedright
\strut
\end{minipage} & \begin{minipage}[b]{0.21\columnwidth}\raggedright
\(n=0\)\strut
\end{minipage} & \begin{minipage}[b]{0.21\columnwidth}\raggedright
\(n=1\)\strut
\end{minipage} & \begin{minipage}[b]{0.21\columnwidth}\raggedright
\(n=2\)\strut
\end{minipage}\tabularnewline
\midrule
\endhead
\begin{minipage}[t]{0.26\columnwidth}\raggedright
\(k=0\)\strut
\end{minipage} & \begin{minipage}[t]{0.21\columnwidth}\raggedright
sets\strut
\end{minipage} & \begin{minipage}[t]{0.21\columnwidth}\raggedright
categories\strut
\end{minipage} & \begin{minipage}[t]{0.21\columnwidth}\raggedright
\(2\)-categories\strut
\end{minipage}\tabularnewline
\begin{minipage}[t]{0.26\columnwidth}\raggedright
\strut
\end{minipage} & \begin{minipage}[t]{0.21\columnwidth}\raggedright
\strut
\end{minipage} & \begin{minipage}[t]{0.21\columnwidth}\raggedright
\strut
\end{minipage} & \begin{minipage}[t]{0.21\columnwidth}\raggedright
\strut
\end{minipage}\tabularnewline
\begin{minipage}[t]{0.26\columnwidth}\raggedright
\(k=1\)\strut
\end{minipage} & \begin{minipage}[t]{0.21\columnwidth}\raggedright
monoids\strut
\end{minipage} & \begin{minipage}[t]{0.21\columnwidth}\raggedright
monoidal categories\strut
\end{minipage} & \begin{minipage}[t]{0.21\columnwidth}\raggedright
monoidal \(2\)-categories\strut
\end{minipage}\tabularnewline
\begin{minipage}[t]{0.26\columnwidth}\raggedright
\strut
\end{minipage} & \begin{minipage}[t]{0.21\columnwidth}\raggedright
\strut
\end{minipage} & \begin{minipage}[t]{0.21\columnwidth}\raggedright
\strut
\end{minipage} & \begin{minipage}[t]{0.21\columnwidth}\raggedright
\strut
\end{minipage}\tabularnewline
\begin{minipage}[t]{0.26\columnwidth}\raggedright
\(k=2\)\strut
\end{minipage} & \begin{minipage}[t]{0.21\columnwidth}\raggedright
commutative monoids\strut
\end{minipage} & \begin{minipage}[t]{0.21\columnwidth}\raggedright
braided monoidal categories\strut
\end{minipage} & \begin{minipage}[t]{0.21\columnwidth}\raggedright
braided monoidal \(2\)-categories\strut
\end{minipage}\tabularnewline
\begin{minipage}[t]{0.26\columnwidth}\raggedright
\strut
\end{minipage} & \begin{minipage}[t]{0.21\columnwidth}\raggedright
\strut
\end{minipage} & \begin{minipage}[t]{0.21\columnwidth}\raggedright
\strut
\end{minipage} & \begin{minipage}[t]{0.21\columnwidth}\raggedright
\strut
\end{minipage}\tabularnewline
\begin{minipage}[t]{0.26\columnwidth}\raggedright
\(k=3\)\strut
\end{minipage} & \begin{minipage}[t]{0.21\columnwidth}\raggedright
`` "\strut
\end{minipage} & \begin{minipage}[t]{0.21\columnwidth}\raggedright
symmetric monoidal categories\strut
\end{minipage} & \begin{minipage}[t]{0.21\columnwidth}\raggedright
weakly involutory monoidal \(2\)-categories\strut
\end{minipage}\tabularnewline
\begin{minipage}[t]{0.26\columnwidth}\raggedright
\strut
\end{minipage} & \begin{minipage}[t]{0.21\columnwidth}\raggedright
\strut
\end{minipage} & \begin{minipage}[t]{0.21\columnwidth}\raggedright
\strut
\end{minipage} & \begin{minipage}[t]{0.21\columnwidth}\raggedright
\strut
\end{minipage}\tabularnewline
\begin{minipage}[t]{0.26\columnwidth}\raggedright
\(k=4\)\strut
\end{minipage} & \begin{minipage}[t]{0.21\columnwidth}\raggedright
`` "\strut
\end{minipage} & \begin{minipage}[t]{0.21\columnwidth}\raggedright
`` "\strut
\end{minipage} & \begin{minipage}[t]{0.21\columnwidth}\raggedright
strongly involutory monoidal \(2\)-categories\strut
\end{minipage}\tabularnewline
\begin{minipage}[t]{0.26\columnwidth}\raggedright
\strut
\end{minipage} & \begin{minipage}[t]{0.21\columnwidth}\raggedright
\strut
\end{minipage} & \begin{minipage}[t]{0.21\columnwidth}\raggedright
\strut
\end{minipage} & \begin{minipage}[t]{0.21\columnwidth}\raggedright
\strut
\end{minipage}\tabularnewline
\begin{minipage}[t]{0.26\columnwidth}\raggedright
\(k=5\)\strut
\end{minipage} & \begin{minipage}[t]{0.21\columnwidth}\raggedright
`` "\strut
\end{minipage} & \begin{minipage}[t]{0.21\columnwidth}\raggedright
`` "\strut
\end{minipage} & \begin{minipage}[t]{0.21\columnwidth}\raggedright
`` "\strut
\end{minipage}\tabularnewline
\bottomrule
\end{longtable}

I described this in \protect\hyperlink{week74}{``Week 74''} and
\protect\hyperlink{week121}{``Week 121''}, so I won't do so again. My
point here is really just that lots of this \(3\)-dimensional stuff is
part of a bigger picture that applies to all different dimensions. For
more details, including a description of the center construction, try:

\begin{enumerate}
\def\labelenumi{\arabic{enumi})}
\setcounter{enumi}{2}
\tightlist
\item
  John Baez and James Dolan, ``Categorification'', in \emph{Higher
  Category Theory}, eds.~Ezra Getzler and Mikhail Kapranov, 
  \emph{Contemp.\ Math.}~\textbf{230}, AMS, Providence, 1998, pp.~1--36. 
   Also available at
  \href{https://arxiv.org/abs/math.QA/9802029}{\texttt{math.QA/9802029}}.
\end{enumerate}

Anyway, M\"uger's elegant characterization of a modular tensor category
amounts to this: it's a braided 2-H*-algebra whose center is
``trivial''. This means that every object in the center is a direct sum
of copies of the object \(1\) --- the unit for the tensor product.

M\"uger does a lot more in his paper that I won't describe here, and he
also said a lot of interesting things in his talk about the general
concept of center. For example, the center of a monoidal category is a
braided monoidal category. In particular, if you take the center of a
2-H*-algebra you get a braided 2-H*-algebra. But what if you then
take this braided 2-H*-algebra and look at \emph{its} center? Well, it
turns out to be ``trivial'' in the above sense!

There's a bit of overlap between M\"uger's paper and this one:

\begin{enumerate}
\def\labelenumi{\arabic{enumi})}
\setcounter{enumi}{3}
\tightlist
\item
  A. Bruguieres, ``Categories premodulaires, modularisations et
  invariants des varietes de dimension 3'', available as 
\href{https://imag.umontpellier.fr/~bruguieres/docs/modular.pdf}{\texttt{https://imag.umontpellier.fr/\textasciitilde{}bruguieres/}}
\href{https://imag.umontpellier.fr/~bruguieres/docs/modular.pdf}{\texttt{docs/modular.pdf}}
\end{enumerate}
\noindent
One especially important issue they both touch upon is this: if you have
a braided 2-H*-algebra, is there any way to mess with it slightly
to get a modular tensor category? The answer is yes. Thus we can really
get a topological quantum field theory from any braided 2-H*-algebra.
But this raises another question: can we describe this topological
quantum field theory directly, without using the modular tensor
category? The answer is again yes! For details see:

\begin{enumerate}
\def\labelenumi{\arabic{enumi})}
\setcounter{enumi}{4}
\tightlist
\item
  Stephen Sawin, ``Jones--Witten invariants for nonsimply-connected Lie
  groups and the geometry of the Weyl alcove'', available as
  \href{https://arxiv.org/abs/math.QA/9905010}{\texttt{math.QA/9905010}}.
\end{enumerate}

This paper uses this machinery to get topological quantum field theories
related to Chern--Simons theory. People have thought about this a lot,
ever since Reshetikhin and Turaev, but the really great thing about this
paper is that it handles the case when the gauge group isn't
simply-connected. This introduces a lot of subtleties which previous
papers touched upon only superficially. Sawin works it out much more
thoroughly by an analysis of subsets of the Weyl alcove that are closed
under tensor product. It's very pretty, and reading it is very good
exercise if you want to learn more about representations of quantum
groups.

Now, I said that a lot of this is part of a bigger picture that works in
higher dimensions. However, a lot of this higher-dimensional stuff
remains very mysterious. Here are two cool papers that make some
progress in unlocking these mysteries:

\begin{enumerate}
\def\labelenumi{\arabic{enumi})}
\setcounter{enumi}{5}
\item
  Marco Mackaay, ``Finite groups, spherical \(2\)-categories, and
  4-manifold invariants'', available as
  \href{https://arxiv.org/abs/math.QA/9903003}{\texttt{math.QA/9903003}}.
\item
  Mikhail Khovanov, ``A categorification of the Jones polynomial'',
  available as \hfill \break
  \href{https://arxiv.org/abs/math.QA/9908171}{\texttt{math.QA/9908171}}.
\end{enumerate}

Marco Mackaay spoke about his work in Coimbra, and I had grilled him
about it in Lisbon beforehand, so I think I understand it pretty well.
Basically what he's doing is categorifying the \(3\)-dimensional
topological quantum field theories studied by Dijkgraaf and Witten to
get \(4\)-dimensional theories. It fits in very nicely with his earlier
work described in \protect\hyperlink{week121}{``Week 121''}.

People have been trying to categorify the magic of quantum groups for
quite some time now, and Khovanov appears to have made a good step in
that direction by describing the Jones polynomial of a link as the
``graded Euler characteristic'' of a chain complex of graded vector
spaces. Since graded Euler characteristic is a generalization of the
dimension of a vector space, and taking the dimension is a process of
decategorification (i.e.~vector spaces are isomorphic iff they have the
same dimension), Khovanov's chain complex can be thought of as a
categorified version of the Jones polynomial.

I would like to understand better the relation between Khovanov's work
and the work of Crane and Frenkel on categorifying quantum groups (see
\protect\hyperlink{week58}{``Week 58''}). For this, I guess I should
read the following papers:

\begin{enumerate}
\def\labelenumi{\arabic{enumi})}
\setcounter{enumi}{7}
\item
  J. Bernstein, I. Frenkel and M. Khovanov, ``A categorification of the
  Temperley--Lieb algebra and Schur quotients of \(U(\mathfrak{sl}_2)\)
  by projective and Zuckerman functors'', \emph{Selecta
  Mathematica}, New Series \textbf{5} (1999), 199--241.  Also available as \href{https://arxiv.org/abs/math/0002087}{\texttt{math/0002087}}.
\item
  Mikhail Khovanov, \emph{Graphical Calculus, Canonical Bases and
  Kazhdan--Lusztig Theory}, Ph.D.~thesis, Department of Mathematics,
 Yale University, 1997.  Available at \href{https://www.math.columbia.edu/~khovanov/research/thesis.pdf}{\texttt{https://www.math.columbia.edu/~khovanov/research/thesis.pdf}}
\end{enumerate}

\hypertarget{week138}{%
\section{September 12, 1999}\label{week138}}

I haven't been going to the Newton Institute much during my stay in
Cambridge, even though it's right around the corner from where I live.
There's always interesting math and physics going on at the Newton
Institute, and this summer they had some conferences on cosmology, but
I've been trying to get away from it all for a while. Still, I couldn't
couldn't resist the opportunity to go to James Hartle's 60th birthday
party, which was held there on September 2nd.

Hartle is famous for his work on quantum gravity and the foundations of
quantum mechanics, so some physics bigshots came and gave talks. First
Chris Isham spoke on applications of topos theory to quantum mechanics,
particularly in relation to Hartle's work on the so-called ``decoherent
histories'' approach to quantum mechanics, which he developed with
Murray Gell-Mann. Then Roger Penrose spoke on his ideas of
gravitationally induced collapse of the wavefunction. Then Gary Gibbons
spoke on the geometry of quantum mechanics. All very nice talks!

Finally, Stephen Hawking spoke on his work with Hartle. This talk was
the most personal in nature: Hawking interspersed technical descriptions
of the papers they wrote together with humorous reminiscences of their
get-togethers in Santa Barbara and elsewhere, including a long trip in a
Volkswagen beetle with Hawking's wheelchair crammed into the back seat.

I think Hawking said he wrote 4 papers with Hartle. The first really
important one was this:

\begin{enumerate}
\def\labelenumi{\arabic{enumi})}
\tightlist
\item
  James Hartle and Stephen Hawking, ``Path integral derivation of black
  hole radiance'', \emph{Phys. Rev.} \textbf{D13} (1976), 2188.
\end{enumerate}

As I explained in \protect\hyperlink{week111}{``Week 111''}, Hawking
wrote a paper in 1975 establishing a remarkable link between black hole
physics and thermodynamics. He showed that a black hole emits radiation
just as if it had a temperature inversely proportional to its mass.
However, this paper was regarded with some suspicion at the time, not
only because the result was so amazing, but because the calculation
involved modes of the electromagnetic field of extremely short
wavelengths near the event horizon --- much shorter than the Planck
length.

For this reason, Hartle and Hawking decided to redo the calculation
using path integrals --- a widely accepted technique in particle
physics. Hawking's background was in general relativity, so he wasn't
too good at path integrals; Hartle had more experience with particle
physics and knew how to do that kind of stuff.

(By now, of course, Hawking can do path integrals quicker than most
folks can balance their checkbook. This was a while ago.)

This wasn't straightforward. In particle physics people usually do
calculations assuming spacetime is flat, so Hartle and Hawking needed to
adapt the usual path-integral techniques to the case when spacetime
contains a black hole. The usual trick in path integrals is to replace
the time variable t by an imaginary number, then do the calculation, and
then analytically continue the answer back to real times. This isn't so
easy when there's a black hole around!

For starters, you have to analytically continue the Schwarzschild
solution (the usual metric for a nonrotating black hole) to imaginary
values of the time variable. When you do this, something curious
happens: you find that the Schwarzschild solution is periodic in the
imaginary time direction. And the period is proportional to the black
hole's mass.

Now, if you are good at physics, you know that doing quantum field
theory calculations where imaginary time is periodic with period \(1/T\)
is the same as doing statistical mechanics calculations where the
temperature is \(T\). So right away, you see that a black hole acts like
it has a temperature inversely proportional to its mass!

(In case you're worried, I'm using units where \(\hbar\), \(c\), \(G\),
and \(k\) are equal to \(1\).)

Anyway, that's how people think about the Hartle-Hawking paper these
days. I haven't actually read it, so my description may be a bit
anachronistic. Things usually look simpler and clearer in retrospect.

The other really important paper by Hartle and Hawking is this one:

\begin{enumerate}
\def\labelenumi{\arabic{enumi})}
\setcounter{enumi}{1}
\tightlist
\item
  James Hartle and Stephen Hawking, ``Wavefunction of the universe'',
  \emph{Phys. Rev.} \textbf{D28} (1983), 2960.
\end{enumerate}

In quantum mechanics, we often describe the state of a physical system
by a wavefunction --- a complex-valued function on the classical
configuration space. If quantum mechanics applies to the whole universe,
this naturally leads to the question: what's the wavefunction of the
universe? In the above paper, Hartle and Hawking propose an answer.

Now, it might seem a bit overly ambitious to guess the wavefunction of
the entire universe, since we haven't even seen the whole thing yet. And
indeed, if someone claims to know the wavefunction of the whole
universe, you might think they were claiming to know everything that has
happened or will happen. Which naturally led Gell-Mann to ask Hartle:
``If you know the wavefunction of the universe, why aren't you rich
yet?''

But the funny thing about quantum theory is that, thanks to the
uncertainty principle, you can know the wavefunction of the universe,
and still be completely clueless as to which horse will win at the races
tomorrow, or even how many planets orbit the sun.

That will either make sense to you, or it won't, and I'm not sure
anything \emph{short} I might write will help you understand it if you
don't already. A full explanation of this business would lead me down
paths I don't want to tread just now --- right into that morass they
call ``the interpretation of quantum mechanics''.

So instead of worrying too much about what it would \emph{mean} to know
the wavefunction of the universe, let me just explain Hartle and
Hawking's formula for it. Mind you, this formula may or may not be
correct, or even rigorously well-defined --- there's been a lot of
argument about it in the physics literature. However, it's pretty cool,
and definitely worth knowing.

Here things get a wee bit more technical. Suppose that space is a
3-sphere, say \(X\). The classical configuration space of general
relativity is the space of metrics on \(X\). The wavefunction of the
universe should be some complex-valued function on this classical
configuration space. And here's Hartle and Hawking's formula for it:
\[\psi(q) = \int_{g|X=q} \exp(-S(g)/\hbar)dg.\] Now you can wow your
friends by writing down this formula and saying ``Here's the
wavefunction of the universe!''

But, what does it mean?

Well, the integral is taken over Riemannian metrics \(g\) on a 4-ball
whose boundary is \(X\), but we only integrate over metrics that
restrict to a given metric \(q\) on \(X\) --- that's what I mean by
writing \(g|X = q\). The quantity \(S(g)\) is the Einstein-Hilbert
action of the metric \(g\) --- in other words, the integral of the Ricci
scalar curvature of \(g\) over the 4-ball. Finally, of course, \(\hbar\)
is Planck's constant.

The idea is that, formally at least, this wavefunction is a solution of
the Wheeler-DeWitt equation, which is the basic equation of quantum
gravity (see \protect\hyperlink{week43}{``Week 43''}).

The measure ``\(dg\)'' is, unfortunately, ill-defined! In other words,
one needs to use lots of clever tricks to extract physics from this
formula, as usual for path integrals. But one can do it, and Hawking and
others have spent a lot of time ever since 1983 doing exactly this. This
led to a subject called ``quantum cosmology''.

I should add that there are lots of ways to soup up the basic
Hartle-Hawking formula. If we have other fields around besides gravity,
we just throw them into the action in the action in the obvious way and
integrate over them too. If our manifold \(X\) representing space is not
a 3-sphere, we can pick some other 4-manifold having it as boundary. If
we can't make up our mind which 4-manifold to use, we can try a ``sum
over topologies'', summing over all 4-manifolds with \(X\) as boundary.
We can do this even when \(X\) is a 3-sphere, actually --- but it's a
bit controversial whether we should, and also whether the sum converges.

Well, there's a lot more to say, like what the physical interpretation
of the Hartle-Hawking formula is, and what predicts. It's actually quite
cool --- in a sense, it says that the universe tunnelled into being out
of nothingness! But that sounds like a bunch of nonsense --- the sort of
fluff they write on the front of science magazines to sell copies. To
really explain it takes quite a bit more work. And unfortunately, it's
just about dinner-time, so I want to stop now.

Anyway, it was an interesting birthday party.

\hypertarget{week139}{%
\section{September 19, 1999}\label{week139}}

Last time I described some of the talks at James Hartle's 60th birthday
celebration at the Newton Institute. But I also met some people at that
party that I'd been wanting to talk to. There's a long story behind
this, so if you don't mind, I'll start at the beginning\ldots.

A while ago, Phillip Helbig, one of the two moderators of
\texttt{sci.physics.research} who do astrophysics, drew my attention to
an interesting paper:

\begin{enumerate}
\def\labelenumi{\arabic{enumi})}
\tightlist
\item
  Vipul Periwal, ``Cosmological and astrophysical tests of quantum
  gravity'', available as
  \href{https://arxiv.org/abs/astro-ph/9906253}{\texttt{astro-ph/9906253}}.
\end{enumerate}
\noindent
The basic idea behind this is that quantum gravity effects could cause
deviations from Newton's inverse square law at large distance scales,
and that these deviations might explain various puzzles in astrophysics,
like the ``missing mass problem'' and the possibly accelerating
expansion of the universe.

This would be great, because it might not only help us understand these
astrophysics puzzles, but also help solve the big problem with quantum
gravity, namely the shortage of relevant experimental data.

But of course one needs to read the fine print before getting too
excited about ideas like this!

Following the argument in Periwal's paper requires some familiarity with
the renormalization group, since that's what people use to study how
``constants'' like the charge of the electron or Newton's gravitational
constant depend on the distance scale at which you measure them --- due
to quantum effects. Reading the paper, I immediately became frustrated
with my poor understanding of the renormalization group. It's really
important, so I decided to read more about it and explain it in the
simplest possible terms on \texttt{sci.physics.research} --- since to
understand stuff, I like to try to explain it.

In the process, I found this book very helpful:

\begin{enumerate}
\def\labelenumi{\arabic{enumi})}
\tightlist
\item
  Michael E. Peskin and Daniel V. Schroeder, \emph{An Introduction to
  Quantum Field Theory}, Addison-Wesley, Reading, Massachusetts, 1995.
\end{enumerate}

The books I'd originally learned quantum field theory from didn't
incorporate the modern attitude towards renormalization, due to Kenneth
Wilson --- the idea that quantum field theory may not ultimately be true
at very short distance scales, but that's okay, because if we assume
it's a good approximation at pretty short distance scales, it becomes a
\emph{better} approximation at \emph{larger} distance scales. This is
especially important when you're thinking about quantum gravity, where
godawful strange stuff may be happening at the Planck length. Peskin and
Schroeder explain this idea quite well. For my own sketchy summary, try
this:

\begin{enumerate}
\def\labelenumi{\arabic{enumi})}
\setcounter{enumi}{1}
\tightlist
\item
  John Baez, Renormalization made easy, available at
  \href{http://math.ucr.edu/home/baez/renormalization.html}{\texttt{http://math.ucr.edu/home/}}  \href{http://math.ucr.edu/home/baez/renormalization.html}{\texttt{baez/renormalization.html}}
\end{enumerate}
\noindent
I deliberately left out as much math as possible, to concentrate on the
basic intuition.

Thus fortified, I returned to Periwal's paper, and it made a bit more
sense. Let me describe the main idea: how we might expect Newton's
gravitational constant to change with distance.

So, suppose we have any old quantum field theory with a coupling
constant \(G\) in it. In fact, \(G\) will depend on the length scale at
which we measure it. But using Planck's constant and the speed of light
we can translate length into \(1/\mathrm{momentum}\). This allows us to
think of \(G\) as a function of momentum. Roughly speaking, when you
shoot particles at each other at higher momenta, they come closer
together before bouncing off, so measuring a coupling constant at a
higher momentum amounts to measuring at a shorter distance scale.

The equation describing how \(G\) depends on the momentum \(p\) is
called the ``Callan--Symanzik equation''. In general it looks like this:
\[\frac{dG}{d(\ln p)} = \beta(G)\] But all the fun starts when we use
our quantum field theory to calculate the right hand side, which is
called --- surprise! --- the ``beta function'' of our theory. Typically
we get something like this:
\[\frac{dG}{d(\ln p)} = (n-d)G+aG_2+bG_3+\ldots\] Here \(n\) is the
dimension of spacetime and \(d\) is a number called the ``upper critical
dimension''. You see, it's fun when possible to think of our quantum
field theory as defined in a spacetime of arbitrary dimension, and then
specialize to the case at hand. I'll show you how work out \(d\) in a
minute. It's harder to work out the numbers \(a\), \(b\), and so on ---
for this, you need to do some computations using the quantum field
theory in question.

What does the Callan--Symanzik equation really mean? Well, for starters
let's neglect the higher-order terms and suppose that
\[\frac{dG(p)}{d(\ln p)} = (n-d)G\] This says \(G\) is proportional to
\(p^{n-d}\). There are 3 cases:

\begin{enumerate}
\def\labelenumi{\Alph{enumi})}
\item
  When \(n < d\), our coupling constant gets \emph{smaller} at higher
  momentum scales, and we say our theory is ``superrenormalizable''.
  Roughly, this means that at larger and larger momentum scales, our
  theory looks more and more like a ``free field theory'' --- one where
  particles don't interact at all. This makes superrenormalizable
  theories easy to study by treating them as a free field theory plus a
  small perturbation.
\item
  When \(n > d\), our coupling constant gets \emph{larger} at higher
  momentum scales, and we say our theory is ``nonrenormalizable''. Such
  theories are hard to study using perturbative calculations in free
  field theory.
\item
  When \(n = d\), we are right on the brink between the two cases above.
  We say our theory is ``renormalizable'', but we really have to work
  out the next term in the beta function to see if the coupling constant
  grows or shrinks with increasing momentum.
\end{enumerate}

Consider the example of general relativity. We can figure out the upper
critical dimension using a bit of dimensional analysis and handwaving.
Let's work in units where Planck's constant and the speed of light are
\(1\). The Lagrangian is the Ricci scalar curvature divided by
\(8\pi G\), where \(G\) is Newton's gravitational constant. We need to
get something dimensionless when we integrate the Lagrangian over
spacetime to get the action, since we exponentiate the action when doing
path integrals in quantum field theory. Curvature has dimensions of
\(1/\mathrm{length}^2\), so when spacetime has dimension \(n\), \(G\)
must have dimensions of \(\mathrm{length}^{n-2}\).

This means that if you are a tiny little person with a ruler \(X\) times
smaller than mine, Newton's constant will seem \(X^{n-2}\) times bigger
to you. But measuring Newton's constant at a length scale that's \(X\)
times smaller is the same as measuring it at a momentum scale that's
\(X\) times bigger. We already solved the Callan--Symanzik equation and
saw that when we measure \(G\) at the momentum scale \(p\), we get an
answer proportional to \(p^{n-d}\). We thus conclude that \(d = 2\).

(If you're a physicist, you might enjoy finding the holes in the above
argument, and then plugging them.)

This means that quantum gravity is nonrenormalizable in 4 dimensions.
Apparently gravity just keeps looking stronger and stronger at shorter
and shorter distance scales. That's why quantum gravity has
traditionally been regarded as hard --- verging on hopeless.

However, there is a subtlety. We've been ignoring the higher-order terms
in the beta function, and we really shouldn't!

This is obvious for renormalizable theories, since when \(n = d\), the
beta function looks like \[\frac{dG}{d(\ln p)} = aG_2+bG_3+\ldots\] so
if we ignore the higher-order terms, we are ignoring the whole
right-hand side! To see the effect of these higher-order terms let's
just consider the simple case where \[\frac{dG}{d(\ln p)} = aG_2\] If
you solve this you get \[G = \frac{c}{1-ac\ln p}\] where \(c\) is a
positive constant. What does this mean? Well, if \(a < 0\), it's obvious
even before solving the equation that \(G\) slowly \emph{decreases} with
increasing momentum. In this case we say our theory is ``asymptotically
free''. For example, this is true for the strong force in the Standard
Model, so in collisions at high momentum quarks and gluons act a lot
like free particles. (For more on this, try
\protect\hyperlink{week94}{``Week 94''}.)

On the other hand, if \(a > 0\), the coupling constant \(G\)
\emph{increases} with increasing momentum. To make matters worse, it
becomes \emph{infinite} at a sufficiently high momentum! In this case
we say our theory has a ``Landau pole'', and we cluck our tongues
disapprovingly, because it's not a good thing. For example, this is what
happens in quantum electrodynamics when we don't include the weak force.
Of course, one should really consider the effect of even higher-order
terms in the beta function before jumping to conclusions. However,
particle physicists generally feel that among renormalizable field
theories, the ones with \(a < 0\) are good, and the ones with \(a > 0\)
are bad.

Okay, now for the really fun part. Perturbative quantum gravity in 2
dimensions is not only renormalizable (because this is the upper
critical dimension), it's also asympotically free! Thus in \(n\)
dimensions, we have \[\frac{dG}{d(\ln p) = (n-2)G+aG_2+\ldots}\] where
\(a < 0\). If we ignore the higher-order terms which I have written as
``\(\ldots\)'', this implies something very interesting for quantum
gravity in 4 dimensions. Suppose that at low momenta \(G\) is small.
Then the right-hand side is dominated by the first term, which is
positive. This means that as we crank up the momentum scale, \(G\) keeps
getting bigger. This is what we already saw about nonrenormalizable
theories. But after a while, when \(G\) gets big, the second term starts
mattering more --- and it's negative. So the growth of \(G\) starts
slowing!

In fact, it's easy to see that as we keep cranking up the momentum,
\(G\) will approach the value for which \[\frac{dG}{d(\ln p)} = 0.\] We
call this value an ``ultraviolet stable fixed point'' for the
gravitational constant. Mathematically, what we've got is a flow in the
space of coupling constants, and an ultraviolet stable fixed point is
one that attracts nearby points as we flow in the direction of higher
momenta. This particular kind of ultraviolet stable fixed point ---
coming from an asymptotically free theory in dimensions above its upper
critical dimension --- is called a ``Wilson--Fisher fixed point''.

So: perhaps quantum gravity is saved from an ever-growing Newton's
constant at small distance scales by a Wilson--Fisher fixed point! But
before we break out the champagne, note that we neglected the
higher-order terms in the beta function in our last bit of reasoning.
They can still screw things up. For example, if
\[\frac{dG}{d(\ln p)} = (n-2)G+aG_2+bG_3\] and \(b\) is positive, there
will not be a Wilson--Fisher fixed point when the dimension \(n\) gets
too far above \(2\). Is \(4\) too far above \(2\)? Nobody knows for
sure. We can't really work out the beta function exactly. So, as usual
in quantum gravity, things are a bit iffy.

However, Periwal cites the following paper as giving numerical evidence
for a Wilson--Fisher fixed point in quantum gravity:

\begin{enumerate}
\def\labelenumi{\arabic{enumi})}
\setcounter{enumi}{2}
\tightlist
\item
  Herbert W. Hamber and Ruth M. Williams, ``Newtonian potential in
  quantum Regge gravity'', \emph{Nucl. Phys.} \textbf{B435} (1995),
  361--397.
\end{enumerate}

And he draws some startling conclusions from the existence of this fixed
point. He says it should have consequences for the missing mass problem
and the value of the cosmological constant! However, I found it hard to
follow his reasoning, so I decided to track down some of the references
--- starting with the above paper.

Now, Ruth Williams works at Cambridge University, so I was not surprised
to find her at Hartle's party. She was busy talking to John Barrett, who
also does quantum gravity, up at Nottingham University. I arranged to
stop by her office, get a copy of her paper, and have her explain it to
me. I also arranged to visit John in Nottingham and have him explain his
work with Louis Crane on Lorentzian spin foam models --- but more about
that next week!

Anyway, here's how the Hamber--Williams paper goes, very roughly. They
simulate quantum gravity by chopping up a \(4\)-dimensional torus into
\(16\times16\times16\times16\) hypercubes, chopping each hypercube into
24 \(4\)-simplices in the obvious way, and then doing a Monte Carlo
calculation of the path integral using the Regge calculus, which is a
discretized version of general relativity suited to triangulated
manifolds (see \protect\hyperlink{week119}{``Week 119''} for details).
Their goal was to work out how Newton's constant varies with distance.
They did it by calculating correlations between Wilson loops that wrap
around the torus. They explain how you can deduce Newton's constant from
this information, but I don't have the energy to describe that here.
Anyway, they claim that Newton's constant varies with distance as one
would expect if there was a Wilson--Fisher fixed point!

(It's actually more complicated that this because besides Newton's
constant, there is also another coupling constant in their theory: the
cosmological constant. And of course this is very important for
potential applications to astrophysics.)

Unfortunately, I'm still mystified about a large number of things. Let
me just mention two. First, Hamber and Williams consider values of \(G\)
which are \emph{greater} than the Wilson--Fisher fixed point. Since this
is an ultraviolet stable fixed point, such values of \(G\) flow
\emph{down} to the fixed point as we crank up the momentum scale. Or in
other words, in this regime Newton's constant gets \emph{bigger} with
increasing distances. At least to my naive brain, this sounds nice for
explaining the missing mass problem. But the funny thing is, this regime
is utterly different from the regime where \(G\) is close to zero ---
namely, \emph{less} than the Wilson--Fisher fixed point. I thought all
the usual perturbative quantum gravity calculations were based on the
assumption that at macroscopic distance scales \(G\) is small, and flows
up to the fixed point as we crank up the momentum scale! Are these folks
claiming this picture is completely wrong? I'm confused.

Another puzzle is that Periwal thinks Newton's constant will start to
grow at distance scales roughly comparable to the radius of the universe
(or more precisely, the Hubble length). But it looks like Hamber and
Williams say their formula for \(G\) as a function of momentum holds at
\emph{short} distance scales.

I guess I need to read more stuff, starting perhaps with Weinberg's old
paper on quantum gravity and the renormalization group:

\begin{enumerate}
\def\labelenumi{\arabic{enumi})}
\setcounter{enumi}{3}
\tightlist
\item
  Steven Weinberg, ``Ultraviolet divergences in quantum theories of
  gravitation'', in \emph{General Relativity: an Einstein Centenary
  Survey}, eds.~Stephen Hawking and Werner Israel, Cambridge U.\ Press,
  Cambridge, 1979.
\end{enumerate}
\noindent
and then perhaps turning to his paper on the cosmological constant:

\begin{enumerate}
\def\labelenumi{\arabic{enumi})}
\setcounter{enumi}{4}
\tightlist
\item
  Steven Weinberg, ``The cosmological constant problem'',
  \emph{Rev.~Mod. Phys.} \textbf{61} (1989), 1.
\end{enumerate}
\noindent
and some books on the renormalization group and quantum gravity:

\begin{enumerate}
\def\labelenumi{\arabic{enumi})}
\setcounter{enumi}{5}
\item
  Claude Itzykson and Jean-Michel Drouffe, \emph{Statistical Field
  Theory}, two volumes, Cambridge U.\ Press, Cambridge, 1989.
\item
  Jean Zinn-Justin, \emph{Quantum Field Theory and Critical Phenomena},
  Oxford U. Press, Oxford, 1993.
\item
  Jan Ambj\o rn, Bergfinnur Durhuus, and Thordur Jonsson, \emph{Quantum
  Geometry: A Statistical Field Theory Approach}, Cambridge U.\ Press,
  Cambridge, 1997.
\end{enumerate}
\noindent
I should also think more about this recent paper, which claims to find a
phase transition in a toy model of quantum gravity where one does the
path integral over a special class of metrics --- namely those with 2
Killing vector fields.

\begin{enumerate}
\def\labelenumi{\arabic{enumi})}
\setcounter{enumi}{8}
\tightlist
\item
  Viqar Husain and Sebastian Jaimungal, ``Phase transition in quantum
  gravity'',  \emph{AIP Conference Proceedings}
 \textbf{493}, American Institute of Physics, Woodbury, New York,
  1999, pp.\ 238--242. Also available as
  \href{https://arxiv.org/abs/gr-qc/9908056}{\texttt{gr-qc/9908056}}.
\end{enumerate}
\noindent
But if anyone can help me clear up these issues, please let me know!

Okay, enough of that. Another person I met at the party was Roger
Penrose! Later I visited him in Oxford. Though recently retired, he
still holds monthly meetings at his house in the country, attended by a
bunch of young mathematicians and physicists. At the one I went to, the
discussion centered around Penrose's forthcoming book. The goal of this
book is to explain modern physics to people who know only a little math,
but are willing to learn more. A nice thing about it is that it treats
various modern physics fads without the uncritical adulation that mars
many popularizations. In particular, when I visited, he was busy writing
a chapter on inflationary cosmology, so he talked about a bunch of
problems with that theory, and cosmology in general.

I've never been sold on inflation, since it relies on fairly speculative
aspects of grand unified theories (or GUTs), so most of these problems
merely amused me. Theorists take a certain wicked glee in seeing someone
else's theory in trouble. However, one of these problems concerned the
Standard Model, and this hit closer to home. Penrose made the standard
observation that the most distant visible galaxies in opposite
directions have not had time to exchange information --- at least not
since the time of recombination, when the initial hot fireball cooled
down enough to become transparent. But if the symmetry between the
electromagnetic and weak forces is spontaneously broken only when the
Higgs field cools down enough to line up, as the Standard Model
suggests, this raises the danger that the Higgs field could wind up
pointing in different directions in different patches of the visible
universe! --- since these different ``domains'' would not yet have had
time to expand to the point where a single one fills the whole visible
universe. But we don't see such domains --- or more precisely, we don't
see the ``domain walls'' one would expect at their boundaries.

Of course, inflation is an attempt to deal with similar problems, but
inflation is posited to happen at GUT scale energies, too soon (it
seems) to solve \emph{this} problem, which happens when things cool down
to the point where the electroweak symmetry breaks.

Again, if anyone knows anything about this, I'd love to hear about it.

\hypertarget{week140}{%
\section{October 16, 1999}\label{week140}}

Let's start with something fun: biographies!

\begin{enumerate}
\def\labelenumi{\arabic{enumi})}
\item
  Norman Macrae, \emph{John von Neumann: The Scientific Genius Who
  Pioneered the Modern Computer, Game Theory, Nuclear Deterrence and
  Much More}, American Mathematical Society, Providence, Rhode Island,
  1999.
\item
  Steve Batterson, \emph{Stephen Smale: The Mathematician Who Broke the
  Dimension Barrier}, American Mathematical Society, Providence, Rhode
  Island, 2000.
\end{enumerate}

Von Neumann might be my candidate for the best mathematical physicist of
the 20th century. His work ranged from the ultra-pure to the
ultra-applied. At one end: his work on axiomatic set theory. At the
other: designing and building some of the first computers to help design
the hydrogen bomb --- which was so applied, it got him in trouble at the
Institute for Advanced Studies! But there's so much stuff in between:
the mathematical foundations of quantum mechanics (von Neumann algebras,
the Stone-Von Neumann theorem and so on), ergodic theory, his work on
Hilbert's fifth problem, the Manhattan project, game theory, the theory
of self-reproducing cellular automata\ldots. You may or may not like
him, but you can't help being awed. Hans Bethe, no dope himself, said of
von Neumann that ``I always thought his brain indicated that he belonged
to a new species, an evolution beyond man''. The mathematician Polya
said ``Johnny was the only student I was ever afraid of.'' Definitely an
interesting guy.

While von Neumann is one of those titans that dominated the first half
of the 20th century, Smale is more typical of the latter half --- he
protested the Vietnam war, and now he even has his own web page!

\begin{enumerate}
\def\labelenumi{\arabic{enumi})}
\setcounter{enumi}{2}
\tightlist
\item
  Stephen Smale's web page,
  \url{http://www.math.berkeley.edu/~smale/}
\end{enumerate}

He won the Fields medal in 1966 for his work on differential topology.
Some of his work is what you might call ``pure'': figuring out how to
turn a sphere inside out without any crinkles, proving the Poincar\'e
conjecture in dimensions 5 and above, stuff like that. But a lot of it
concerns dynamical systems: cooking up strange attractors using the
horseshoe map, proving that structural stability is not generic, and so
on --- long before the recent hype about chaos theory began! More
recently he's also been working on economics, game theory, and the
relation of computational complexity to algebraic geometry.

Now for some papers on spin networks and spin foams:

\begin{enumerate}
\def\labelenumi{\arabic{enumi})}
\setcounter{enumi}{3}
\tightlist
\item
  Roberto De Pietri, Laurent Freidel, Kirill Krasnov, and Carlo Rovelli,
  ``Barrett--Crane model from a Boulatov--Ooguri field theory over a
  homogeneous space'', \emph{Nucl.\ Phys.\ B} \textbf{574} (2000), 
     785--806.  Also available as
  \href{https://arxiv.org/abs/hep-th/9907154}{\texttt{hep-th/9907154}}.
\end{enumerate}
\noindent
The Barrett--Crane model is a very interesting theory of quantum gravity.
I've described it already in \protect\hyperlink{week113}{``Week 113''},
\protect\hyperlink{week120}{``Week 120''} and
\protect\hyperlink{week128}{``Week 128''}, so I won't go into much
detail about it --- I'll just plunge in\ldots.

The original Barrett--Crane model involved a fixed triangulation of
spacetime. One can also cook up versions where you sum over
triangulations. In some ways the most natural is to sum over all ways of
taking a bunch of \(4\)-simplices and gluing them face-to-face until no
faces are left free. Some of these ways give you manifolds; others
don't. In this paper, the authors show that this ``sum over
triangulations'' version of the Barrett--Crane model can be thought of as
a quantum field theory on a product of 4 copies of the 3-sphere. Weird,
huh?

But it's actually not so weird. The space of complex functions on the
\((n-1)\)-sphere is naturally a representation of \(\mathrm{SO}(n)\).
But there's another way to think of this representation. Consider an
triangle in \(\mathbb{R}^n\). We can associate vectors to two of its
edges, say \(v\) and \(w\), and form the wedge product of these vectors
to get a bivector \(v\wedge w\). This bivector describes the area
element associated to the triangle. If we pick an orientation for the
triangle, this bivector is uniquely determined. Now, a bivector of the
form \(v\wedge w\) is called ``simple''. The space of simple bivectors
is naturally a Poisson manifold --- i.e., we can define Poisson brackets
of functions on this space --- so we can think of it as a ``classical
phase space''. Using geometric quantization, we can quantize this
classical phase space and get a Hilbert space. Since rotations act on
the classical phase space, they act on this Hilbert space, so it becomes
a representation of \(\mathrm{SO}(n)\). And this representation is
isomorphic to the space of complex functions on the \((n-1)\)-sphere!

Thus, we can think of a complex function on the \((n-1)\)-sphere as a
``quantum triangle'' in \(\mathbb{R}^n\), as long as we really just care
about the area element associated to the triangle. One can develop this
analogy in detail and make it really precise. In particular, one can
describe a ``quantum tetrahedron'' in \(\mathbb{R}^n\) as a collection
of 4 quantum triangles satisfying some constraints that say the fit
together into a tetrahedron. These quantum tetrahedra act almost like
ordinary tetrahedra when they are large, but when the areas of their
faces becomes small compared to the square of the Planck length, the
uncertainty principle becomes important: you can't simultaneously know
everything about their geometry with perfect precision.

Let me digress for a minute and sketch the history of this stuff. The
quantum tetrahedron in 3 dimensions was invented by Barbieri --- see
\protect\hyperlink{week110}{``Week 110''}. It quickly became important
in the study of spin foam models. Then Barrett and I systematically
worked out how to construct the quantum tetrahedron in 3 and 4
dimensions using geometric quantization --- see
\protect\hyperlink{week134}{``Week 134''}. Subsequently, Freidel and
Krasnov figured out how to generalize this stuff to higher dimensions:

\begin{enumerate}
\def\labelenumi{\arabic{enumi})}
\setcounter{enumi}{4}
\item
  Laurent Freidel, Kirill Krasnov and Raymond Puzio, ``BF description of
  higher-dimensional gravity'', \emph{Adv. Theor. Math. Phys.} \textbf{3} 
  (1999), 1289--1324.  Also available as
  \href{https://arxiv.org/abs/hep-th/9901069}{\texttt{hep-th/9901069}}.
\item
  Laurent Freidel and Kirill Krasnov, ``Simple spin networks as Feynman
  graphs'', \emph{Jour. Math. Phys.} \textbf{41} (2000), 1681--1690.
   Available as
  \href{https://arxiv.org/abs/hep-th/9903192}{\texttt{hep-th/9903192}}.
\end{enumerate}

So much for history --- now back to business. So far I've told you that
the state of a ``quantum triangle'' in 4 dimensions is given by a
complex function on the 3-sphere. And I've told you that a ``quantum
tetrahedron'' is a collection of 4 quantum triangles satisfying some
constraints. More precisely, let \[H = L^2(S^3)\] be the Hilbert space
for a quantum triangle in 4 dimensions. Then the Hilbert space for a
quantum tetrahedron is a certain subspace \(T\) of
\(H\otimes H\otimes H\otimes H\), where ``\(\otimes\)'' denotes the
tensor product of Hilbert spaces. Concretely, we can think of states in
\(T\) as complex functions on the product of 4 copies of \(S^3\). These
complex functions need to satisfy some constraints, but let's not worry
about those\ldots.

Now let's ``second quantize'' the Hilbert space \(T\). This is physics
jargon for making a Hilbert space out of the algebra of polynomials on
\(T\) --- usually called the ``Fock space'' on \(T\). As usual, there
are two pictures of states in this Fock space: the ``field'' picture and
the ``particle'' picture. On the one hand, they are states of a quantum
field theory on the product of 4 copies of \(S^3\). But on the other
hand, they are states of an arbitrary collection of quantum tetrahedra
in 4 dimensions. In other words, we've got ourselves a quantum field
theory whose ``elementary particles'' are quantum tetrahedra!

The idea of the de Pietri-Freidel-Krasnov-Rovelli paper is to play these
two pictures off each other and develop a new way of thinking about the
Barrett--Crane model. The main thing these guys do is write down a
Lagrangian with some nice properties. Throughout quantum field theory,
one of the big ideas is to start with a Lagrangian and use it to compute
the amplitudes of Feynman diagrams. A Feynman diagram is a graph with
edges corresponding to ``particles'' and vertices corresponding to
``interactions'' where a bunch of particles turns into another bunch of
particles.

But in the present context, the so-called ``particles'' are really
quantum tetrahedra! Thus the trick is to write down a Lagrangian giving
Feynman diagrams with 5-valent vertices. If you do it right, these
5-valent vertices correspond exactly to ways that 5 quantum tetrahedra
can fit together as the 5 faces of a \(4\)-simplex! Let's call such a
thing a ``quantum \(4\)-simplex''. Then your Feynman diagrams correspond
exactly to ways of gluing together a bunch of quantum \(4\)-simplices
face-to-face. Better yet, if you set things up right, the amplitude for
such a Feynman diagram exactly matches the amplitude that you'd compute
for a triangulated manifold using the Barrett--Crane model!

In short, what we've got here is a quantum field theory whose Feynman
diagrams describe ``quantum geometries of spacetime'' --- where
spacetime is not just a fixed triangulated manifold, but any possible
way of gluing together a bunch of \(4\)-simplices face-to-face.

Sounds great, eh? So what are the problems? I'll just list three. First,
we don't know that the ``sum over Feynman diagrams'' converges in this
theory. In fact, it probably does not --- but there are probably ways to
deal with this. Second, the model is Riemannian rather than Lorentzian:
we are using the rotation group when we should be using the Lorentz
group. Luckily this is addressed in a new paper by Barrett and Crane.
Third, we aren't very good at computing things with this sort of model
--- short of massive computer simulations, it's tough to see what it
actually says about physics. In my opinion this is the most serious
problem: we should either start doing computer simulations of spin foam
models, or develop new analytical techniques for handling them --- or
both!

Now, this ``new paper by Barrett and Crane'' is actually not brand new.
It's a revised version of something they'd already put on the net:

\begin{enumerate}
\def\labelenumi{\arabic{enumi})}
\setcounter{enumi}{6}
\tightlist
\item
  John Barrett and Louis Crane, ``A Lorentzian signature model for
  quantum general relativity'', \emph{Class. Quant. Grav.} \textbf{17} (2000),     
   3101--3118.  Also available as
  \href{https://arxiv.org/abs/gr-qc/9904025}{\texttt{gr-qc/}} 
  \href{https://arxiv.org/abs/gr-qc/9904025}{\texttt{9904025}}.
\end{enumerate}
\noindent
However, it's much improved. When I went up to Nottingham at the end of
the summer, I had Barrett explain it to me. I learned all sorts of cool
stuff about representations of the Lorentz group. Unfortunately, I don't
now have the energy to explain all that stuff. I'll just say this:
everything I said above generalizes to the Lorentzian case. The main
difference is that we use the \(3\)-dimensional hyperboloid
\[H^3 = \{t^2 - x^2 - y^2 - z^2 = 1\}\] wherever we'd been using the
3-sphere \[S^3 = \{t^2 + x^2 + y^2 + z^2 = 1\}.\] It's sort of obvious
in retrospect, but it's nice that it works out so neatly!

Okay, here are some more papers on spin networks and spin foams. Since
I'm getting lazy, I'll just quote the abstracts:

\begin{enumerate}
\def\labelenumi{\arabic{enumi})}
\setcounter{enumi}{7}
\tightlist
\item
  Sameer Gupta, ``Causality in spin foam models'', \emph{Phys.\ Rev.\ D}
 \textbf{61} (2000), 064014.  Also available as
  \href{https://arxiv.org/abs/gr-qc/9908018}{\texttt{gr-qc/9908018}}.
\end{enumerate}

\begin{quote}
We compute Teitelboim's causal propagator in the context of canonical
loop quantum gravity. For the Lorentzian signature, we find that the
resultant power series can be expressed as a sum over branched, colored
two-surfaces with an intrinsic causal structure. This leads us to define
a general structure which we call a ``causal spin foam''. We also
demonstrate that the causal evolution models for spin networks fall in
the general class of causal spin foams.
\end{quote}

\begin{enumerate}
\def\labelenumi{\arabic{enumi})}
\setcounter{enumi}{8}
\tightlist
\item
  Matthias Arnsdorf and Sameer Gupta, ``Loop quantum gravity on
  non-compact spaces'', \emph{Nucl.\ Phys.\ B} \textbf{577} (2000),
   529--546. Also available as
  \href{https://arxiv.org/abs/gr-qc/9909053}{\texttt{gr-qc/9909053}}.
\end{enumerate}

\begin{quote}
We present a general procedure for constructing new Hilbert spaces for
loop quantum gravity on non-compact spatial manifolds. Given any fixed
background state representing a non-compact spatial geometry, we use the
Gel'fand-Naimark-Segal construction to obtain a representation of the
algebra of observables. The resulting Hilbert space can be interpreted
as describing fluctuation of compact support around this background
state. We also give an example of a state which approximates classical
flat space and can be used as a background state for our construction.
\end{quote}

\begin{enumerate}
\def\labelenumi{\arabic{enumi})}
\setcounter{enumi}{9}
\tightlist
\item
  Seth A. Major, ``Quasilocal energy for spin-net gravity'',  \emph{Class.\ 
  Quant.\  Grav.\ } \break \textbf{17} (2000), 1467--1487.  Also
  available as
  \href{https://arxiv.org/abs/gr-qc/9906052}{\texttt{gr-qc/9906052}}.
\end{enumerate}

\begin{quote}
The Hamiltonian of a gravitational system defined in a region with
boundary is quantized. The classical Hamiltonian, and starting point for
the regularization, is required by functional differentiablity of the
Hamiltonian constraint. The boundary term is the quasilocal energy of
the system and becomes the ADM mass in asymptopia. The quantization is
carried out within the framework of canonical quantization using spin
networks. The result is a gauge invariant, well-defined operator on the
Hilbert space induced from the state space on the whole spatial
manifold. The spectrum is computed. An alternate form of the operator,
with the correct naive classical limit, but requiring a restriction on
the Hilbert space, is also defined. Comparison with earlier work and
several consequences are briefly explored.
\end{quote}

\begin{enumerate}
\def\labelenumi{\arabic{enumi})}
\setcounter{enumi}{10}
\tightlist
\item
  C. Di Bartolo, R. Gambini, J. Griego, J. Pullin, ``Consistent
  canonical quantization of general relativity in the space of Vassiliev
  knot invariants'', \emph{Phys.\ Rev.\ Lett.\ }\textbf{84} (200), 2314--2317.  Also available as
  \href{https://arxiv.org/abs/gr-qc/9909063}{\texttt{gr-qc/9909063}}.
\end{enumerate}

\begin{quote}
We present a quantization of the Hamiltonian and diffeomorphism
constraint of canonical quantum gravity in the spin network
representation. The novelty consists in considering a space of
wavefunctions based on the Vassiliev knot invariants. The constraints
are finite, well defined, and reproduce at the level of quantum
commutators the Poisson algebra of constraints of the classical theory.
A similar construction can be carried out in \(2+1\) dimensions leading to
the correct quantum theory.
\end{quote}

\begin{enumerate}
\def\labelenumi{\arabic{enumi})}
\setcounter{enumi}{11}
\tightlist
\item
  John Baez, ``Spin foam perturbation theory'', in \emph{Diagrammatic Morphisms and Applications}, eds.\ David Radford, Fernando Souza, and David Yetter, \emph{Contemp. Math.} \textbf{318}, American Mathematical Society, Providence, Rhode Island, 2003, pp.\ 9--21.  Also available as
  \href{https://arxiv.org/abs/gr-qc/9910050}{\texttt{gr-qc/9910050}}.
\end{enumerate}

\begin{quote}
We study perturbation theory for spin foam models on triangulated
manifolds. Starting with any model of this sort, we consider an
arbitrary perturbation of the vertex amplitudes, and write the evolution
operators of the perturbed model as convergent power series in the
coupling constant governing the perturbation. The terms in the power
series can be efficiently computed when the unperturbed model is a
topological quantum field theory. Moreover, in this case we can
explicitly sum the whole power series in the limit where the number of
top-dimensional simplices goes to infinity while the coupling constant
is suitably renormalized. This `dilute gas limit' gives spin foam models
that are triangulation-independent but not topological quantum field
theories. However, we show that models of this sort are rather trivial
except in dimension 2.
\end{quote}

\hypertarget{week141}{%
\section{October 26, 1999}\label{week141}}

How can you resist a book with a title like ``Inconsistent
Mathematics''?

\begin{enumerate}
\def\labelenumi{\arabic{enumi})}
\tightlist
\item
  Chris Mortensen, \emph{Inconsistent Mathematics}, Kluwer, Dordrecht,
  1995.
\end{enumerate}

Ever since G\"odel showed that all sufficiently strong systems formulated
using the predicate calculus must either be inconsistent or incomplete,
most people have chosen what they perceive as the lesser of two evils:
accepting incompleteness to save mathematics from inconsistency. But
what about the other option?

This book begins with the startling sentence: ``The following idea has
recently been gaining support: that the world is or might be
inconsistent.'' As we reel in shock, Mortensen continues:

\begin{quote}
Let us consider set theory first. The most natural set theory to adopt
is undoubtedly one which has unrestricted set abstraction (also known as
naive comprehension). This is the natural principle which declares that
to every property there is a unique set of things having the property.
But, as Russell showed, this leads rapidly to the contradiction that the
the Russell set {[}the set of all sets that do not contain themselves as
a member{]} both is and is not a member of itself. The overwhelming
majority of logicians took the view that this contradiction required a
weakening of unrestricted abstraction in order to ensure a consistent
set theory, which was in turn seen as necessary to provide a consistent
foundation for mathematics. But all ensuing attempts at weakening set
abstraction proved to be in various ways ad hoc. Da Costa and Routley
both suggested instead that the Russell set might be dealt with more
naturally in an inconsistent but nontrivial set theory (where triviality
means that every sentence is provable).
\end{quote}

An inconsistent but nontrivial logical system is called
\emph{paraconsistent}. But it's not so easy to create such systems. To
keep an inconsistency from infecting the whole system and making it
trivial, we need to drop the rule of classical logic which says that
``\(A\) and \(\operatorname{not}(A)\) implies \(B\)'' for all
propositions \(A\) and \(B\). Unfortunately, this rule is built into the
propositional calculus from the very start!

So, we need to revise the propositional calculus.

One way to do it is to abandon ``material implication'' --- the form of
implication where you can calculate the truth value of ``\(P\) implies
\(Q\)'' from those of \(P\) and \(Q\) using the following truth table:

\begin{longtable}[]{@{}ccc@{}}
\toprule
\(P\) & \(Q\) & \(P\implies Q\)\tabularnewline
\midrule
\endhead
T & T & T\tabularnewline
T & F & F\tabularnewline
F & T & T\tabularnewline
F & F & T\tabularnewline
\bottomrule
\end{longtable}

With material implication, a false statement implies \emph{every}
statement, so any inconsistency is fatal. But in real life, if we
discover we have one inconsistent belief, we don't conclude we can fly
and go jump off a building! Material implication is really just our best
attempt to define implication using truth tables with 2 truth values:
true and false. So it's not surprising that logicians have investigated
other forms of implication.

One obvious approach is to use more truth values, like ``true'',
``false'', and ``don't know''. There's a long history of work on such
multi-valued logics.

Another approach, initiated by Anderson and Belnap, is called
``relevance logic''. In relevance logic, ``\(P\) implies \(Q\)'' can
only be true if there is a conceptual connection between \(P\) and
\(Q\). So if \(B\) has nothing to do with \(A\), we don't get ``\(A\)
and \(\operatorname{not}(A)\) implies \(B\)''.

This book describes a logical system called ``RQ'' --- relevance logic
with quantifiers. It also describes a system called ``R\#'', which is a
version of the Peano axioms of arithmetic based on RQ instead of the
usual predicate calculus. Following the work of Robert Meyer, it proves
that R\# is nontrivial in the sense described above. Moreover, this
proof can be carried out R\# itself! However, you can carry out the
proof of G\"odel's 2nd incompleteness theorem in R\#, so R\# cannot prove
itself consistent.

To paraphrase Mortensen: ``But this is not really a puzzle. The
explanation is that relevant and other paraconsistent logics turn on
making a distinction between inconsistency and triviality, the former
being weaker than the latter; whereas classical logical cannot make this
distinction. For what the present author's intuitions are worth, these
do seem to be different concpets. Thus for R\#, consistency cannot be
proved by finitistic means by G\"odel's second theorem, whereas
nontriviality can be shown. Since Peano arithmetic collapses this
distinction, both kinds of consistency are infected by the same
unprovability.''

Mortensen also mentions another approach to get rid of ``\(A\) and
\(\operatorname{not}(A)\) implies \(B\)'' without getting rid of
material implication. This is to get rid of the rule that ``\(A\) and
\(\operatorname{not}(A)\)'' is false! He calls this ``Brazilian logic''.
Presumably this is not because your average Brazilian thinks this way,
but because the inventor of this approach, Da Costa, is Brazilian.

Brazilian logic sounds very bizarre at first, but in fact it's just the
dual of intuitionistic logic, where you drop the rule that ``\(A\) or
\(\operatorname{not}(A)\)'' is true. Intuitionistic logic is nicely
modeled by open sets in a topological space: ``and'' is intersection,
``or'' is union, and ``not'' is the interior of the complement.
Similarly, Brazilian logic is modeled by closed sets. In intuitionistic
logic we allow a slight gap between A and \(\operatorname{not}(A)\); in
Brazilian logic we allow a slight overlap.

In short, this book is full of fascinating stuff. Lots of passages are
downright amusing at first, like this:

\begin{quote}
{[}\ldots{]} there have been calls recently for inconsistent calculus,
appealing to the history of calculus in which inconsistent claims
abound, especially about infinitesimals (Newton, Leibniz, Bernoulli,
l'Hospital, even Cauchy). However, inconsistent calculus has resisted
development.
\end{quote}

But you always have to remember that the author is interested in
theories which, though inconsistent, are still paraconsistent. And I
think he really makes a good case for his claim that inconsistent
mathematics is worth studying --- even if our ultimate goal is to
\emph{avoid} inconsistency!

Okay, now let me switch gears drastically and say a bit about ``exotic
spheres'' --- smooth manifolds that are homeomorphic but not
diffeomorphic to the \(n\)-sphere with its usual smooth structure.
People on \texttt{sci.physics.research} have been talking about this
stuff lately, so it seems like a good time for a mini-essay on the
subject. Also, my colleague Fred Wilhelm works on the geometry of exotic
spheres, and he just gave a talk on it here at U. C. Riverside, so I
should pass along some of his wisdom while I still remember it.

First, recall the ``Hopf bundle''. It's easy to describe starting with
the complex numbers. The unit vectors in \(\mathbb{C}^2\) form the
sphere \(S^3\). The unit complex numbers form a group under
multiplication. As a manifold this is just the circle \(S^1\), but as a
group it's better known as \(\mathrm{U}(1)\). You can multiply a unit
vector by a unit complex number and get a new unit vector, so \(S^1\)
acts on \(S^3\). The quotient space is the complex projective space
\(\mathbb{CP}^1\), which is just the sphere \(S^2\). So what we've got
here is fiber bundle: \[S^1 \to S^3\to S^2 = \mathbb{CP}^1\] with fiber
\(S^1\), total space \(S^3\) and base space \(S^2\). This is the Hopf
bundle. It's famous because the map from the total space to the base was
the first example of a topologically nontrivial map from a sphere to a
sphere of lower dimension. In the lingo of homotopy theory, we say it's
the generator of the group \(\pi_3(S^2)\).

Now in \protect\hyperlink{week106}{``Week 106''} I talked about how we
can mimic this construction by replacing the complex numbers with any
other division algebra. If we use the real numbers we get a fiber bundle
\[S^0\to S^1\to \mathbb{RP}^1 = S^1\] where \(S^0\) is the group of unit
real numbers, better known as \(\mathbb{Z}/2\). This bundle looks like
the edge of a Moebius strip. If we use the quaternions we get a more
interesting fiber bundle: \[S^3\to S^7\to \mathbb{HP}^1 = S^4\] where
\(S^3\) is the group of unit quaternions, better known as
\(\mathrm{SU}(2)\). We can even do something like this with the
octonions, and we get a fiber bundle
\[S^7\to S^{15}\to \mathbb{OP}^1 = S^8\] but now \(S^7\), the unit
octonions, doesn't form a group --- because the octonions aren't
associative.

Anyway, it's the quaternionic version of the Hopf bundle that serves as
the inspiration for Milnor's construction of exotic 7-spheres. These
exotic 7-spheres are actually total spaces of \emph{other} bundles with
fiber \(S^3\) and base space \(S^4\). The easiest way to get your hands
on these bundles is to take \(S^4\), chop it in half along the equator,
put a trivial \(S^3\)-bundle over each hemisphere, and then glue these
together. To glue these bundles together we need a way to attach the
fibers over each point x of the equator. In other words, for each point
\(x\) in the equator of \(S^4\) we need a map \[f_x\colon S^3\to S^3\]
which should vary smoothly with \(x\). But the equator of \(S^4\) is
just \(S^3\), and \(S^3\) is a group --- the unit quaternions --- so we
can take \[f_x(y) = x^n y x^m\] for any pair of integers \((n,m)\).

This gives us a bunch of \(S^3\)-bundles over \(S^4\). The total space
\(X(n,m)\) of any one of these bundles is obviously a smooth
\(7\)-dimensional manifold. But when is it homeomorphic to the 7-sphere?
And when is it \emph{diffeomorphic} to the 7-sphere with its usual
smooth structure?

Well, first we use some Morse theory. You can learn a lot about the
topology of a smooth manifold if you have a ``Morse function'' on the
manifold: a smooth real-valued function all of whose critical points are
nondegenerate. If you don't believe me, read this book:

\begin{enumerate}
\def\labelenumi{\arabic{enumi})}
\setcounter{enumi}{1}
\tightlist
\item
  John Milnor, \emph{Morse Theory}, Princeton U. Press, Princeton, 1960.
\end{enumerate}
\noindent
When \(n + m = 1\) there's a Morse function on \(X(n,m)\) with only two
critical points --- a maximum and a minimum. This implies that
\(X(n,m)\) is homeomorphic to a sphere!

Once we know that \(X(n,m)\) is homeomorphic to \(S^7\), we have to
decide when it's diffeomorphic to \(S^7\) with its usual smooth
structure. This is the hard part. Notice that \(X(n,m)\) is the unit
sphere bundle of a vector bundle over \(S^4\) whose fiber is the
quaternions. We can understand a bunch about \(X(n,m)\) using the
characteristic classes of this vector bundle. In particular, we can
compute the Euler number and the Pontrjagin number of this vector
bundle. Using the Euler number we can show that \(X(n,m)\) is
homeomorphic to a sphere \emph{only} if \(n + m = 1\) --- you can't
really do this using Morse theory. But more importantly, using the
Pontrjagin number, we can show that in this case \(X(n,m)\) is
diffeomorphic to \(S^7\) with its usual smooth structure if and only if
\((n-m)^2 = 1 \mod 7\). Otherwise it's ``exotic''.

For the details of the above argument you can try the following book:

\begin{enumerate}
\def\labelenumi{\arabic{enumi})}
\setcounter{enumi}{2}
\tightlist
\item
  B. A. Dubrovin, A. T. Fomenko and S. P. Novikov, \emph{Modern Geometry
  --- Methods and Applications, Part III: Introduction to Homology
  Theory}, Springer Graduate Texts \textbf{125}, Springer, Berlin, 1990.
\end{enumerate}
\noindent
or the original paper:

\begin{enumerate}
\def\labelenumi{\arabic{enumi})}
\setcounter{enumi}{3}
\tightlist
\item
  John Milnor, ``On manifolds homeomorphic to the 7-sphere'', \emph{Ann.\
  Math.} \textbf{64} (1956), 399--405.
\end{enumerate}

Now, with quite a bit more work, you can show that smooth structures on
the \(n\)-sphere form an group under connected sum --- the operation of
chopping out a small hole in two spheres and gluing them together ---
and you can show that this group is \(\mathbb{Z}/28\) for \(n = 7\).
This means that if we consider two smooth structures on the 7-sphere the
same when they're related by an \emph{orientation-preserving}
diffeomorphism, we get exactly 28 kinds. Unfortunately we don't get all
of them by the above explicit construction. For more details, see:

\begin{enumerate}
\def\labelenumi{\arabic{enumi})}
\setcounter{enumi}{4}
\tightlist
\item
  M. Kervaire and J. Milnor, ``Groups of homotopy spheres I'',
  \emph{Ann. Math.} \textbf{77} (1963), 504--537.
\end{enumerate}

By the way, part II of the above paper doesn't exist! Instead, you
should read this:

\begin{enumerate}
\def\labelenumi{\arabic{enumi})}
\setcounter{enumi}{5}
\tightlist
\item
  J. Levine, ``Lectures on groups of homotopy spheres'', in
  \emph{Algebraic and Geometric Topology}, Springer Lecture Notes in
  Mathematics \textbf{1126}, Springer, Berlin, 1985, pp.~62--95.
\end{enumerate}

Anyway, if you're wondering why I'm talking so much about exotic
7-spheres, instead of lower-dimensional examples that are easier to
visualize, check out this table of groups of smooth structures on the
\(n\)-sphere:

\begin{longtable}[]{@{}ll@{}}
\toprule
\(n\) & group of smooth structures on the \(n\)-sphere\tabularnewline
\midrule
\endhead
\(0\) & \(1\)\tabularnewline
\(1\) & \(1\)\tabularnewline
\(2\) & \(1\)\tabularnewline
\(3\) & \(1\)\tabularnewline
\(4\) & ?\tabularnewline
\(5\) & \(1\)\tabularnewline
\(6\) & \(1\)\tabularnewline
\(7\) & \(\mathbb{Z}/28\)\tabularnewline
\(8\) & \(\mathbb{Z}/2\)\tabularnewline
\(9\) &
\(\mathbb{Z}/2\times\mathbb{Z}/2\times\mathbb{Z}/2\)\tabularnewline
\(10\) & \(\mathbb{Z}/6\)\tabularnewline
\(11\) & \(\mathbb{Z}/992\)\tabularnewline
\(12\) & \(1\)\tabularnewline
\(13\) & \(\mathbb{Z}/3\)\tabularnewline
\(14\) & \(\mathbb{Z}/2\)\tabularnewline
\(15\) & \(\mathbb{Z}/8128\times\mathbb{Z}/2\)\tabularnewline
\(16\) & \(\mathbb{Z}/2\)\tabularnewline
\(17\) &
\(\mathbb{Z}/2\times\mathbb{Z}/2\times\mathbb{Z}/2\times\mathbb{Z}/2\)\tabularnewline
\(18\) & \(\mathbb{Z}/8\times\mathbb{Z}/2\)\tabularnewline
\bottomrule
\end{longtable}

Dimension 7 is the simplest interesting case --- except perhaps for
dimension 4, where the answer is unknown! The ``smooth Poincar\'e
conjecture'' says that there's only one smooth structure on the
4-sphere, but this remains a conjecture\ldots.

As you can see, there are lots of exotic 11-spheres. In fact, this is
relevant to string theory! You can get an \(n\)-sphere with any possible
smooth structure by taking two \(n\)-dimensional balls and gluing them
together along their boundary using some orientation-preserving
diffeomorphism \[f\colon S^{n-1}\to S^{n-1}.\] Orientation-preserving
diffeomorphisms like this form a group called
\(\mathrm{Diff}_+(S^{n-1})\). Using the above trick, it turns out that
the group of smooth structures on the n-sphere is isomorphic to the
group of \emph{connected components} of \(\mathrm{Diff}_+(S^{n-1})\). So
the existence of exotic 11-spheres means that there are lots of ``exotic
diffeomorphisms'' of the 10-sphere!

Now, string theory lives in 10 dimensions, and one wants certain
quantities to be invariant under orientation-preserving diffeomorphisms
of spacetime --- otherwise you say the theory has ``gravitational
anomalies''. First you have to check this for ``small diffeomorphisms''
of spacetime, that is, those connected to the identity map by a
continuous path. But then you have to check it for ``large
diffeomorphisms'' --- those living in different connected components of
the diffeomorphism group. When spacetime is a 10-sphere, this means you
need to check diffeomorphism invariance for all 991 components of
\(\mathrm{Diff}_+(S^{n-1})\) besides the component containing the
identity. These components correspond to exotic 11-spheres!

Witten did this in the following paper:

\begin{enumerate}
\def\labelenumi{\arabic{enumi})}
\setcounter{enumi}{6}
\tightlist
\item
  Edward Witten, ``Global gravitational anomalies'', \emph{Comm. Math.
  Phys.} \textbf{100} (1985), 197--229.
\end{enumerate}
\noindent
This may be the first paper about exotic spheres in physics.

There are other interesting things to do with an exotic sphere. One is
to put a metric on it and look at its curvature. The sphere with its
usual ``round'' metric is very symmetrical and has positive curvature
everywhere. There are various meanings of ``positive curvature'', but
the round sphere has positive curvature in all possible ways! One kind
of curvature is ``sectional curvature''. In general, it's hard to find
compact manifolds other than the sphere with its usual smooth structure
that have metrics with everywhere positive sectional curvature. Gromoll
and Meyer found an exotic 7-sphere with a metric having
\emph{nonnegative} sectional curvature:

\begin{enumerate}
\def\labelenumi{\arabic{enumi})}
\setcounter{enumi}{7}
\tightlist
\item
  Detlef Gromoll and Wolfgang Meyer, ``An exotic sphere with nonnegative
  sectional curvature'', \emph{Ann. Math.} \textbf{100} (1974),
  401--406.
\end{enumerate}

The construction isn't terribly hard so let me describe it. First, start
with the group \(\mathrm{Sp}(2)\), consisting of \(2\times2\) unitary
quaternionic matrices (see \protect\hyperlink{week64}{``Week 64''}). As
always with compact Lie groups, this has a metric that's invariant under
right and left translations, and this metric is unique up to a constant
scale factor. The group of unit quaternions acts as metric-preserving
maps (aka ``isometries'') of \(\mathrm{Sp}(2)\) in the following way:
let the quaternion \(q\) map \[
  \left(
    \begin{array}{cc}
      a&b\\c&d
    \end{array}
  \right)
\] to \[
  \left(
    \begin{array}{cc}
      qaq^{-1}&qb\\qcq^{-1}&qd
    \end{array}
  \right)
\] The quotient space is an exotic 7-sphere, and it inherits a metric
with nonnegative sectional curvature.

Now, since compact manifolds with positive sectional curvature are tough
to find, you might wonder if this exotic 7-sphere can be given a metric
with \emph{positive} sectional curvature. And the answer is: almost! It
can be given a metric having positive sectional curvature except on a
set of measure zero. This was recently proved by Wilhelm:

\begin{enumerate}
\def\labelenumi{\arabic{enumi})}
\setcounter{enumi}{8}
\tightlist
\item
  Frederick Wilhelm, ``An exotic sphere with positive curvature almost
  everywhere'', \emph{Jour.\ Geometric Analysis} \textbf{11} (2001),  
   519--560.
\end{enumerate}

It's also an interesting theorem, due to Hitchin, that for any \(n > 0\)
there exist exotic spheres of dimensions \(8n+1\) and \(8n+2\) having no
metric of positive scalar curvature:

\begin{enumerate}
\def\labelenumi{\arabic{enumi})}
\setcounter{enumi}{7}
\tightlist
\item
  Nigel Hitchin, ``Harmonic spinors'', \emph{Adv. Math.} \textbf{14}
  (1974), 1--55.
\end{enumerate}
\noindent
So some exotic spheres are not so as ``round'' as you might think! In
fact, 3 of the exotic spheres in 10 dimensions cannot be given a metric
such that the connected component of the isometry group is bigger than
\(\mathrm{U}(1)\times\mathrm{U}(1)\), so these are quite ``bumpy''. This
follows from results of Reinhard Schultz, who happens to be the
department chair here:

\begin{enumerate}
\def\labelenumi{\arabic{enumi})}
\setcounter{enumi}{8}
\tightlist
\item
  Reinhard Schultz, ``Circle actions on homotopy spheres bounding
  plumbing manifolds'', \emph{Proc. A.M.S.} \textbf{36} (1972),
  297--300.
\end{enumerate}

There's a lot more to say about exotic spheres, but let me just briefly
mention two things. First, there are cool connections between exotic
spheres and higher-dimensional knot theory. If you want a small taste of
this stuff, try:

\begin{enumerate}
\def\labelenumi{\arabic{enumi})}
\setcounter{enumi}{9}
\tightlist
\item
  Louis Kauffman, \emph{Knots and Physics}, World Scientific, Singapore,
  1991.
\end{enumerate}
\noindent
Look in the index under ``exotic spheres''.

Second, people have computed the effect of exotic 7-spheres on quantum
gravity path integrals in 7 dimensions:

\begin{enumerate}
\def\labelenumi{\arabic{enumi})}
\setcounter{enumi}{10}
\tightlist
\item
  Kristin Schleich and Donald Witt, Exotic spaces in quantum gravity,
  \emph{Class. Quant. Grav.} \textbf{16} (1999), 2447--2469. Also
  available as
  \href{https://arxiv.org/abs/gr-qc/9903086}{\texttt{gr-qc/9903086}}.
\end{enumerate}
\noindent
I'm not sure exotic spheres are \emph{really} relevant to physics, but
it would be cool, so I'm glad some people are trying to establish
connections.

Okay, that's enough for exotic spheres, at least for now! I've got a few
more things here that I just want to mention\ldots.

I've been learning a bit about Calabi--Yau manifolds and mirror symmetry
in string theory lately. The basic idea is that string theory on
different spacetime manifolds can be physically equivalent. I don't know
enough to want to try to explain this stuff yet, but here are some place
to look in case you're interested:

\begin{enumerate}
\def\labelenumi{\arabic{enumi})}
\setcounter{enumi}{11}
\item
  Claire Voisin, \emph{Mirror Symmetry}, American Mathematical Society,
  Providence, Rhode Island, 1999.
\item
  David A. Cox and Sheldon Katz, \emph{Mirror Symmetry and Algebraic
  Geometry}, American Mathematical Society, Providence, Rhode Island,
  1999.
\item
  Shing-Tung Yau, editor, \emph{Mirror Symmetry I}, American
  Mathematical Society, Providence, Rhode Island, 1998.

  Brian Green and Shing-Tung Yau, editors, \emph{Mirror Symmetry II},
  American Mathematical Society, Providence, Rhode Island, 1997.

  Duong H. Phong, Luc Vinet and Shing-Tung Yau, editors, \emph{Mirror
  Symmetry III}, American Mathematical Society, Providence, Rhode Island, 1999.
\end{enumerate}
\noindent
So far I'm mainly trying to learn really basic stuff, and for this, the
following lectures are proving handy:

\begin{enumerate}
\def\labelenumi{\arabic{enumi})}
\setcounter{enumi}{15}
\tightlist
\item
  P. Candelas, ``Lectures on complex manifolds'', in \emph{Superstrings
  '87}, eds.~L. Alvarez-Gaume et al, World Scientific, Singapore, 1988,
  pp.~1--88.
\end{enumerate}

On a different note, the American Mathematical Society has come out with
some good-looking books on surgery theory --- the process of making new
manifolds from old by cutting and pasting. I've got these on my reading
list, so if anyone wants to buy me a Christmas present, here's what you
should get:

\begin{enumerate}
\def\labelenumi{\arabic{enumi})}
\setcounter{enumi}{16}
\item
  Robert E. Gompf and Andras I Stipsicz, \emph{4-Manifolds and Kirby
  Calculus}, Amderican Mathematical Society, Providence, Rhode Island, 1999.
\item
  C. T. C. Wall and A. A. Ranicki, \emph{Surgery on Compact Manifolds},
  2nd edition, American Mathematical Society, Providence, Rhode Island, 1999.
\end{enumerate}

Finally, there's some cool stuff going on with operads that I haven't
been able to keep up with. Let me quote the abstracts:

\begin{enumerate}
\def\labelenumi{\arabic{enumi})}
\setcounter{enumi}{18}
\tightlist
\item
  Alexander A. Voronov, ``Homotopy Gerstenhaber algebras'', 
  in \emph{In Conf\'erence Mosh\'e Flato 1999: Quantization, Deformations, 
 and Symmetries}, Volume II, Springer, Berlin, 2000, pp.\ 307--331.
  Also available as
  \href{https://arxiv.org/abs/math.QA/9908040}{\texttt{math.QA/9908040}}.
\end{enumerate}

\begin{quote}
The purpose of this paper is to complete Getzler--Jones' proof of
Deligne's Conjecture, thereby establishing an explicit relationship
between the geometry of configurations of points in the plane and the
Hochschild complex of an associative algebra. More concretely, it is
shown that the \(B_\infty\)-operad, which is generated by multilinear
operations known to act on the Hochschild complex, is a quotient of a
certain operad associated to the compactified configuration spaces.
Different notions of homotopy Gerstenhaber algebras are discussed: one
of them is a \(B_\infty\)-algebra, another, called a homotopy
\(G\)-algebra, is a particular case of a \(B_\infty\)-algebra, the
others, a \(G_\infty\)-algebra, an \(E^1\)-bar-algebra, and a weak
\(G_\infty\)-algebra, arise from the geometry of configuration spaces.
Corrections to the paper math.QA/9602009 of Kimura, Zuckerman, and the
author related to the use of a nonextant notion of a homotopy
Gerstenhaber algebra are made.
\end{quote}

\begin{enumerate}
\def\labelenumi{\arabic{enumi})}
\setcounter{enumi}{19}
\tightlist
\item
  Maxim Kontsevich, ``Operads and motives in deformation quantization'',
  \emph{Lett. Math. Phys.} \textbf{48} (1999), 35--72.  Also
  available as
  \href{https://arxiv.org/abs/math.QA/9904055}{\texttt{math.QA/9904055}}.
\end{enumerate}

\begin{quote}
It became clear during last 5-6 years that the algebraic world of
associative algebras (abelian categories, triangulated categories, etc)
has many deep connections with the geometric world of two-dimensional
surfaces. One of the manifestations of this is Deligne's conjecture
(1993) which says that on the cohomological Hochschild complex of any
associative algebra naturally acts the operad of singular chains in the
little discs operad. Recently D. Tamarkin discovered that the operad of
chains of the little discs operad is formal, i.e.~it is homotopy
equivalent to its cohomology. From this fact and from Deligne's
conjecture follows almost immediately my formality result in deformation
quantization. I review the situation as it looks now. Also I conjecture
that the motivic Galois group acts on deformation quantizations, and
speculate on possible relations of higher-dimensional algebras and of
motives to quantum field theories.
\end{quote}

\begin{enumerate}
\def\labelenumi{\arabic{enumi})}
\setcounter{enumi}{20}
\tightlist
\item
  James E. McClure and Jeffrey H. Smith, ``A solution of Deligne's
  conjecture'', available as
  \href{https://arxiv.org/abs/math.QA/9910126}{\texttt{math.QA/9910126}}
\end{enumerate}

\begin{quote}
Deligne asked in 1993 whether the Hochschild cochain complex of an
associative ring has a natural action by the singular chains of the
little 2-cubes operad. In this paper we give an affirmative answer to
this question. We also show that the topological Hochschild cohomology
spectrum of an associative ring spectrum has an action of an operad
equivalent to the little 2-cubes.
\end{quote}

\begin{center}\rule{0.5\linewidth}{0.5pt}\end{center}

My original table of groups of smooth structures on spheres had some
mistakes in it which were corrected by Linus Kramer, Marco Mackaay, Tony
Smith and Frederick Wilhelm. In fact, the table in the book by Dubrovin,
Fomenko and Novikov differs from the table in Kervaire and Milnor's
paper! The table above comes from Kervaire and Milnor, taking advantage
of some subsequent work in dimension 3 and also some work of Brumfield
which nailed down the groups in dimensions 9 and 17 --- see below for
more information.

The paper by Kervaire and Milnor has a cool formula for the \emph{order}
of the group of smooth structures on the \((4n-1)\)-sphere for
\(n > 1\). It's: \[\frac{2^{2n-4} (2^{2n-1}-1) P(4n-1) B(n) a(n)}{n}\]
where:

\begin{itemize}
\tightlist
\item
  \(P(k)\) is the order of the \(k\)th stable homotopy group of spheres
\item
  \(B(k)\) is the \(k\)th Bernoulli number, in the sequence
  \(1/6, 1/30, 1/42, 1/30, 5/66, 691/2730, 7/6, \ldots\)
\item
  \(a(k)\) is 1 or 2 according to whether k is even or odd.
\end{itemize}

\begin{center}\rule{0.5\linewidth}{0.5pt}\end{center}

Here are some remarks by Linus Kramer on exotic spheres in dimensions 9
and 17, which he posted to \texttt{sci.math.research} in response to a
question of mine. Kervaire and Milnor said the group of exotic spheres
in dimension 9 was \((\mathbb{Z}/2)^3\) or
\(\mathbb{Z}/2\times\mathbb{Z}/4\), and the group in dimension 17 was
\((\mathbb{Z}/2)^4\) or \((\mathbb{Z}/2)^2\times\mathbb{Z}/4\). Linus
writes:

\begin{quote}
The list by Kervaire and Milnor seems to be correct; in dimension 9, the
group is \((\mathbb{Z}/2)^3\), and in dimension 17 it's
\((\mathbb{Z}/2)^4\). This follows from the results of Brumfield
{[}\emph{Mich. Math. J.} \bf{17}{]} stated on the first page of his paper, plus the
list of the first stable homotopy groups of spheres, and the properties
of Adams' J-homomorphism \(J\colon\pi_n(SO)\to\pi^s_n\).

There is an exact sequence
\[0\to bP_k\to\Gamma_{k-1}\to\pi^s_{k-1}/\operatorname{im}(J)\to 0\]
provided that \(k+3\) is not a power of \(2\) (\(\Gamma_{k-1}\) is the
group we are looking for). Now for \(k-1=8,9,10,17\), we have
\(k+3=12,13,14,21\), and these are not powers of \(2\). Now
\(bP_k=0,\mathbb{Z}/2,0,\mathbb{Z}/\theta_{16}\) for \(k=9,10,11,16\).
So for \(k=9,11\), the group \(\Gamma_{k-1}\) is the same as
\(\pi^s_{k-1}/\operatorname{im}(J)\), and for \(k=10,16\) use
Brumfield's result that then the group is
\(\Gamma_{k-1} = \mathbb{Z}/2+(\pi^s_{k-1}/\operatorname{im}(J))\).

Hope I didn't make a mistake while chasing through all these exact
sequences\ldots{} Of course, the result relies also on some tables,
namely the first stable homotopy groups of spheres.

Linus Kramer
\end{quote}

\hypertarget{week142}{%
\section{December 5, 1999}\label{week142}}

I was recently infected by a meme --- a self-propagating pattern of
human behavior. Now I want to pass it on to you! I like this particular
meme because it's so simple. It's even simpler than the parasites
described on my webpage:

\begin{enumerate}
\def\labelenumi{\arabic{enumi})}
\tightlist
\item
  John Baez, ``Subcellular life forms'', available at
  \href{http://math.ucr.edu/home/baez/subcellular.html}{\texttt{http://math.ucr.edu/home/baez/}}   \href{http://math.ucr.edu/home/baez/subcellular.html}{\texttt{subcellular.html}} 
\end{enumerate}

I wrote this webpage when I was trying to understand some of the
simplest self-reproducing entities: viruses, viroids, virusoids,
plasmids, prions, and various forms of junk DNA. Viroids are especially
simple. Unlike a virus, a viroid doesn't even have a protein coat: it's
just a naked RNA molecule! So instead of actively breaking into the host
cell, it must passively wait to be absorbed. Then somehow it hijacks the
machinery of the cell nucleus to reproduce itself. Theodore Diener
discovered the first viroid in 1971: the potato spindle tuber viroid,
which makes potatoes abnormally long and sometimes cracked. At first
people doubted the possibility of a life form smaller than a virus. But
by now the complete molecular structure of this viroid has been worked
out. It consists of only 359 nucleotides --- or in other words, about
12,000 atoms!

But since a meme relies on the complex apparatus of human culture to
reproduce itself, it can get away with being even simpler than a viroid.
It can even be the simplest sort of thing of all: an abstract
mathematical structure defined by a short list of axioms!

A good example is the game of tic-tac-toe. It's not very interesting,
but it's just interesting enough to keep propagating itself through
human children, who are highly susceptible to the charm of simple games.
Most children soon develop an immunity to tic-tac-toe, just like measles
and mumps --- but only after passing it on to some other child.

Unfortunately, the meme that infected me is a lot harder to shake,
because it's a lot more interesting. I'm talking about the game of Go.

This game is played on a \(19\times19\) square grid. Each player starts
with a large supply of stones --- black for the first player, white for
the second. They take turns putting a stone on a grid point. A group of
stones of one color ``dies'' and is removed from the board when it is
surrounded by stones of the other color. More precisely, we say a stone
is ``dead'' when none of its nearest neighbors of the same color have
nearest neighbors of the same color which have nearest neighbors of the
same color which\ldots{} have nearest neighbors that are still vacant
grid points.

There are also two subsidary rules, designed to keep silly things from
happening.

First, you are not allowed to put a stone someplace where it will
immediately die, \emph{unless} doing so immediately kills one or more of
the other player's stones --- in which case their stones die, and yours
lives.

Second, if putting down your stone kills a stone of the other player,
but they could immediately put that stone back and kill yours, leading
to an infinite loop, we say that ``ko'' has occurred. In this case, the
other player is forbidden from putting their stone back right away.

How do you win? Simply put, the goal is to end up with as much
``territory'' as possible. Territory includes grid points occupied by
stones of your color, and also vacant grid points that the other player
could not occupy without their stones eventually dying. (In practice, Go
players do not fight to the bitter end, so territory also includes
stones of the other color that are ``doomed to die''.)

That's basically it!

The cool part is that starting from these simple rules, a whole world of
strategy unfolds, full of specific tricks --- but also quite general
philosophical lessons about ``power'', ``territory'', and ``threat''. In
a good game, both players start by efficiently marking out some
territory, putting stones down in a widely separated way that looks
random to the beginner, but in fact is delicately balanced between being
too conservative and too ambitious. The midgame starts when both players
start trying to surround each other and threaten to kill stones. But be
careful: threatening to kill stones can be better than actually killing
them, and the difference between ``surrounding'' and ``being
surrounded'' is rather subtle! The endgame comes when territory is
almost fully demarcated, with only a few squabbles around the edges. The
endgame game proves unexpectedly difficult for beginners, since one can
snatch defeat from the jaws of victory even at this stage.

A well-developed Go game is said to be like a work of art, with all
opposing forces neatly balanced in a harmonious pattern. As a
mathematical physicist, it reminds me of the Ising model at a phase
transition, when there are as many black grid points as white ones, and
arbitrarily large clusters of both colors. Perhaps there's even a real
relation to the theory of ``self-organized criticality'', in which a
system spontaneously works its way to the brink of a phase transition.

People say Go was developed in China between 4 and 6 thousand years ago.
Its early history is obscure, but it is said to have started, not as a
game, but as a tool for divination and the teaching of military
strategy. I'm no expert, but to me Go seems like a nice illustration of
yin-yang philosophy --- the idea that the dynamic complexity of the
universe arises from the dialectic interplay of binary opposites. For a
good introduction to what I'm talking about, you can't beat the \emph{I Ching}
--- the ``Classic of Changes'', a Chinese divination text compiled in
the 9th century B.C., but containing material that probably dates back
at least a few centuries earlier. This book describes the significance
of 64 ``hexagrams'', which are patterns built from 6 bits of
information, like this:
\[
  \begin{tikzpicture}  
    \draw[ultra thick] (-1,0.75) to (1,0.75);
    \draw[ultra thick] (-1,0.45) to (-0.33,0.45);  \draw[ultra thick] (0.33,0.45) to (1,0.45);
    \draw[ultra thick] (-1,0.15) to (1,0.15);
    \draw[ultra thick] (-1,-0.15) to (-0.33,-0.15); \draw[ultra thick] (0.33,-0.15) to (1,-0.15);
    \draw[ultra thick] (-1,-0.45) to (1,-0.45);
    \draw[ultra thick] (-1,-0.75) to (1,-0.75);
  \end{tikzpicture}
\]
The idea that complex patterns can be described using bits was borrowed
from the Chinese by Leibniz, who invented the concept of binary
arithmetic and dreamt of a purely mechanical approach to logic based on
simple rules. Now, of course, these ideas dominate modern technology! So
perhaps it's not surprising that Go still holds an attraction for many
mathematicians and physicists.

In fact, I bet some you are smirking and wondering why I didn't learn Go
much earlier! The reason is that I've always avoided playing games,
except for the ``great game'' of mathematical physics. I only tried
playing Go the weekend before last, while visiting my friend Bruce Smith
up in San Rafael after giving a talk on quantum tetrahedra at Stanford.
Bruce explained Go to me and showed me how it was philosophically
interesting. But most importantly, he showed me a computer program that
plays Go. Computers aren't great at Go, but they're good enough to beat
an amateur like me, so they're good to learn from at first, and for some
reason I prefer to play a computer than another person --- perhaps
because computers don't gloat.

The computer program I played against is called ``GNU Go''. You can
download it free from the internet, thanks to the Free Software
Foundation:

\begin{enumerate}
\def\labelenumi{\arabic{enumi})}
\setcounter{enumi}{1}
\tightlist
\item
  GNU Go, \url{http://www.gnu.org/software/gnugo/}
\end{enumerate}
\noindent
You can adjust the size of the board and also the handicap --- the
number of stones you get right away when you start. To use this program
in a UNIX environment you need an interface program called ``cgoban'',
which is also free:

\begin{enumerate}
\def\labelenumi{\arabic{enumi})}
\setcounter{enumi}{2}
\tightlist
\item
  CGoban,
  \url{http://www.inetarena.com/~wms/comp/cgoban/}
\end{enumerate}
\noindent
On Windows you can use an interface available from the GNU Go webpage.

For more information on Go start here:

\begin{enumerate}
\def\labelenumi{\arabic{enumi})}
\setcounter{enumi}{3}
\tightlist
\item
  American Go Association, \url{http://www.usgo.org/}
\end{enumerate}
\noindent
You can find lots of go books listed at this website. Personally I found
these books to be a nice introduction to the game, but they may be hard
to find:

\begin{enumerate}
\def\labelenumi{\arabic{enumi})}
\setcounter{enumi}{4}
\tightlist
\item
  The Nihon Kiin, \emph{Go: The World's Most Fascinating Game}, two
  volumes, Sokosha Printing Co., Tokyo, 1973.
\end{enumerate}
\noindent
When you get more advanced, there are a lot of books to read, with fun
titles like ``Get Strong at Invading'', ``Reducing Territorial
Frameworks'', and ``Utilizing Outward Influence''. It pays to study
``joseki'', or openings:

\begin{enumerate}
\def\labelenumi{\arabic{enumi})}
\setcounter{enumi}{5}
\tightlist
\item
  Ishida Yoshio, \emph{Dictionary of Basic Joseki}, 3 volumes, Ishi
  Press International, San Jose, California, 1977.
\end{enumerate}
\noindent
It's also good to study ``tsume-go'', or ``life and death problems'',
where you figure out which player can win in various configurations. A
mathematician would call this the ``local'' analysis of Go:

\begin{enumerate}
\def\labelenumi{\arabic{enumi})}
\setcounter{enumi}{6}
\tightlist
\item
  Cho Chikun, \emph{All About Life and Death}, two volumes, Ishi Press
  International, San Jose, California, 1993.
\end{enumerate}
\noindent
Ishi Press puts out a lot of other books on Go, but I haven't been able
to get ahold of them yet. I'm sort of fascinated by one that talks about
a difficult abstract concept called ``thickness'', since I suspect this
is a global rather than local concept:

\begin{enumerate}
\def\labelenumi{\arabic{enumi})}
\setcounter{enumi}{7}
\tightlist
\item
  Ishidea Yoshio, \emph{All About Thickness: Understanding Moyo and
  Influence}, Ishi Press International, San Jose, California.
\end{enumerate}
\noindent
If you want to get mathematical about Go endgames, try this:

\begin{enumerate}
\def\labelenumi{\arabic{enumi})}
\setcounter{enumi}{8}
\tightlist
\item
  Elwyn Berlekamp and David Wolfe, \emph{Mathematical Go: Chilling Gets
  the Last Point}, A. K. Peters, Wellesley, Massachusetts, 1994.
\end{enumerate}
\noindent
If you want to get computational, try this:

\begin{enumerate}
\def\labelenumi{\arabic{enumi})}
\setcounter{enumi}{9}
\tightlist
\item
  Markus Enzenberger, \emph{The integration of a priori knowledge into a
  Go playing neural network},
  \url{http://www.cgl.ucsf.edu/go/Programs/neurogo-html/NeuroGo.html}
\end{enumerate}
\noindent
If instead you prefer to curl up with a good novel based on a game of
Go, try this:

\begin{enumerate}
\def\labelenumi{\arabic{enumi})}
\setcounter{enumi}{10}
\tightlist
\item
  Yasunari Kawabata, \emph{The Master of Go}, trans. Edward G.
  Seidensticker, Knopf, New York, 1972.
\end{enumerate}
\noindent
On a different note, here are two good editions of the I Ching:

\begin{enumerate}
\def\labelenumi{\arabic{enumi})}
\setcounter{enumi}{11}
\item
  \emph{The I Ching or Book of Changes}, trans. Richard Wilhelm and Cary
  F. Baynes, Princeton U. Press, Princeton, 1969.

  \emph{The Classic of Changes: A New Translation of the I Ching as
  Interpreted by Wang Bi}, trans. Richard John Lynn, Columbia U. Press,
  1994.
\end{enumerate}
\noindent
Okay. Enough culture --- time for some math!

I was invited to Stanford University by David Carlton, who works on
modular forms, and I found out from him and his friends that the
Shimura--Taniyama--Weil conjecture has been proved! This might have been a
nice scoop for This Week's Finds, but by now it's appeared in the
\emph{Notices of the AMS}, so everyone knows about it:

\begin{enumerate}
\def\labelenumi{\arabic{enumi})}
\setcounter{enumi}{12}
\tightlist
\item
  Henri Darmon, ``A proof of the full Shimura--Taniyama--Weil conjecture
  is announced'', \emph{Notices of the American Mathematical Society},
  \textbf{46} no. 11 (December 1999), 1397--1401.  Also available as
\href{https://www.ams.org/notices/199911/comm-darmon.pdf}{\texttt{https://www.ams.org/notices/199911/comm-darmon.}} \href{https://www.ams.org/notices/199911/comm-darmon.pdf}{\texttt{pdf}}
\end{enumerate}

Andrew Wiles proved part of this conjecture in order to prove Fermat's
Last Theorem, but the conjecture is actually much more interesting than
Fermat's Last Theorem, and a proof of the whole thing was announced this
summer by Breuil, Conrad, Diamond and Taylor.

What does the conjecture say?

Well, first you have to know a bit about elliptic curves. An ``elliptic
curve'' is the space of solutions of an equation like this:
\[y^2 = x^3 + ax + b\] They come up naturally in string theory, and I've
talked about them already in \protect\hyperlink{week13}{``Week 13''} and
\protect\hyperlink{week124}{``Week
124''}--\protect\hyperlink{week127}{``Week 127''}. If all the variables
in sight are complex numbers, an elliptic curve looks like a torus, but
number theorists like to consider the case where the coefficients a and
b are rational. By a simple change of variables you can then get the
coefficients to be integers. Then it makes sense to work modulo a prime
number \(p\): in other words, to think of all the variables as living in
the field of integers \(\mod p\), better known as \(\mathbb{Z}/p\). If
you're smart, you can tell if an elliptic curve \(\mod p\) is
``singular'' or not: being nonsingular is like being a smooth manifold.
People say an elliptic curve has ``good reduction at \(p\)'' if it's
nonsingular \(\mod p\). For any given elliptic curve, this is true
except for finitely many primes.

Any elliptic curve \(E\) has finitely many points \(\mod p\). Let's call
the number of points \(N(E,p)\) and set \[a(E,p) = p - N(E,p).\] If this
list of numbers satisfies a certain condition, which I'll describe in a
minute, we say our elliptic curve is ``modular''. The
Shimura--Taniyama--Weil conjecture states that all elliptic curves are
modular.

Okay, so what does ``modular'' mean? Well, for this we need a little
digression on modular forms. In \protect\hyperlink{week125}{``Week
125''} I described the moduli space of elliptic curves, which is the
space of all different shapes an elliptic curve can have. I showed that
this space was \(H/\mathrm{SL}(2,\mathbb{Z})\), where \(H\) is the upper
half of the complex plane and \(\mathrm{SL}(2,\mathbb{Z})\) is the group
of \(2\times2\) integer matrices with determinant \(1\). A modular form
is basically just a holomorphic section of some line bundle over the
moduli space of elliptic curves. But if this sounds too high-tech, don't
be scared! We can also think of it as an analytic function on the upper
half-plane that transforms in a nice way under the action of
\(\mathrm{SL}(2,\mathbb{Z})\). Remember, any matrix \[
  \left(
    \begin{array}{cc}
      a&b\\c&d
    \end{array}
  \right)
\] in \(\mathrm{SL}(2,\mathbb{Z})\) acts on the upper half-plane as
follows: \[\tau \mapsto \frac{a \tau + b}{c \tau + d}\] For an analytic
function \(f\colon H\to\mathbb{C}\) to be a ``modular form of weight
\(k\)'', it must transform as follows:
\[f\left(\frac{a \tau + b}{c \tau + d}\right) = (c \tau + d)^k f(\tau)\]
for some integer \(k\). We also require that \(f\) satisfy some growth
conditions as \(\tau\to\infty\), so we can expand it as a Taylor series
\[f(\tau) = \sum a_n q^n\] where \[q = \exp(2 \pi i \tau)\] is a
variable that equals \(0\) when \(\tau = \infty\). The nicest modular
forms are the ``cusp forms'', which have \(a_0 = 0\), and thus vanish at
\(\tau = \infty\).

Next, we can straightforwardly generalize everything I just said if we
replace \(\mathrm{SL}(2,\mathbb{Z})\) by various subgroups thereof.
(This amounts to studying holomorphic sections of line bundles over some
moduli space of elliptic curves \emph{equipped with extra structure}.)
For example, we can use the subgroup \(\Gamma_0(N)\) consisting of those
matrices in \(\mathrm{SL}(2,\mathbb{Z})\) whose lower-left entries are
divisible by \(N\). If we use this group instead of
\(\mathrm{SL}(2,\mathbb{Z})\), we get what are called modular forms of
``level \(N\)''. We define ``weight'' of such a modular form just as
before, and ditto for ``cusp forms''.

And now we can say what it means for an elliptic curve to be modular! We
say an elliptic curve \(E\) is ``modular'' if for some \(N\) there's a
weight \(2\) level \(N\) cusp form \[f(\tau) = \sum a_n q^n\] normalized
so that \(a_1 = 1\), with the property that \[a_p = a(E,p)\] for all
primes \(p\) at which \(E\) has good reduction.

So now you know what the Shimura--Taniyama--Weil conjecture says: all
elliptic curves are modular! It's not obvious that this implies Fermat's
Last Theorem, but it does, thanks to a trick invented by Gerhard Frey.

There turn out to be fascinating but mysterious relationships between
the Shimura--Taniyama--Weil conjecture, something called the Langlands
program, and topological quantum field theory:

\begin{enumerate}
\def\labelenumi{\arabic{enumi})}
\setcounter{enumi}{13}
\tightlist
\item
  Mikhail Kapranov, ``Analogies between the Langlands correspondence and
  topological quantum field theory'', in \emph{Functional Analysis on
  the Eve of the 21st Century}, Volume I, Birkhauser, Boston, 1995,
  pp.~119--151.
\end{enumerate}
\noindent
For this reason --- and others --- it's not so surprising that David
Carlton and some of his buddies are interested in \(n\)-categories. In
fact, Carlton caught a small error in the definition of \(n\)-categories
due to James Dolan and myself --- it turns out that the number ``\(1\)''
should be the number ``\(2\)'' at one particular place in the
definition! Anyone who can spot a problem like that is friend of mine.

Even better, Carlton is now interested in understanding the
\((n+1)\)-category of all \(n\)-categories, which is crucial for really
doing anything with \(n\)-categories. Makkai has a new paper on this
subject, and I realize now that I've never mentioned this paper on This
Week's Finds, so let me conclude by quoting the abstract. It's pretty
long and detailed, and probably only of interest to \(n\)-category
addicts\ldots.

\begin{enumerate}
\def\labelenumi{\arabic{enumi})}
\setcounter{enumi}{14}
\tightlist
\item
  M. Makkai, The multitopic \(\omega\)-category of all multitopic
  \(\omega\)-categories, available at
  \url{www.math.mcgill.ca/makkai/mltomcat04/mltomcat04.pdf}
\end{enumerate}

\begin{quote}
``The paper gives two definitions: that of ``multitopic
\(\omega\)-category" and that of ``the (large) multitopic set of all
(small) multitopic \(\omega\)-categories''. It also announces the
theorem that the latter is a multitopic \(\omega\)-category. (The proof
of the theorem will be contained in a sequel to this paper.)

The work has two direct sources. One is the paper {[}H/M/P{]} (for the
references, see at the end of this abstract) in which, among others, the
concept of ``multitopic set'' was introduced. The other is the present
author's work on FOLDS, First Order Logic with Dependent Sorts. The
latter was reported on in {[}M2{]}. A detailed account of the work on
FOLDS is in {[}M3{]}. For the understanding of the present paper, what
is contained in {[}M2{]} suffices. In fact, section 1 of the present
paper gives the definitions of all that's needed in this paper; so,
probably, there won't be even a need to consult {[}M2{]}.

The concept of multitopic set, the main contribution of {[}H/M/P{]},
was, in turn, inspired by the work of J. Baez and J. Dolan {[}B/D{]}.
Multitopic sets are a variant of opetopic sets of loc. cit. The name
``multitopic set'' refers to multicategories, a concept originally due
to J. Lambek {[}L{]}, and given an only moderately generalized
formulation in {[}H/M/P{]}. The earlier ``opetopic set'' of {[}B/D{]} is
based on a concept of operad. I should say that the exact relationship
of the two concepts (``multitopic set'' and ``opetopic set'') is still
not clarified. The main aspect in which the theory of multitopic sets is
in a more advanced state than that of opetopic sets is that, in
{[}H/M/P{]}, there is an explicitly defined category \(\mathsf{Mlt}\) of
\textbf{multitopes}, with the property that the category of multitopic
sets is equivalent to the category of \(\mathsf{Set}\)-valued functors
on \(\mathsf{Mlt}\), a result given a detailed proof in {[}H/M/P{]}. The
corresponding statement on opetopic sets and opetopes is asserted in
{[}B/D{]}, but the category of opetopes is not described. In this paper,
the category of multitopes plays a basic role.

Multitopic sets and multitopes are described in section 2 of this paper;
for a complete treatment, the paper {[}H/M/P{]} should be consulted.

The indebtedness of the present work to the work of Baez and Dolan goes
further than that of {[}H/M/P{]}. The second ingredient of the
Baez/Dolan definition, after ``opetopic set'', is the concept of
``universal cell''. The Baez/Dolan definition of weak \(n\)-category
achieves the remarkable feat of specifying the composition structure by
universal properties taking place in an opetopic set. In particular, a
(weak) opetopic (higher-dimensional) category is an opetopic set with
additional properties ( but with no additional data), the main one of
the additional properties being the existence of sufficiently many
universal cells. This is closely analogous to the way concepts like
``elementary topos'' are specified by universal properties: in our
situation, ``multitopic set'' plays the ``role of the base'' played by
``category'' in the definition of ``elementary topos''. In {[}H/M/P{]},
no universal cells are defined, although it was mentioned that their
definition could be supplied without much difficulty by imitating
{[}B/D{]}. In this paper, the ``universal (composition) structure'' is
supplied by using the concept of FOLDS-equivalence of {[}M2{]}.

In {[}M2{]}, the concepts of ``FOLDS-signature'' and
``FOLDS-equivalence'' are introduced. A (FOLDS-) signature is a category
with certain special properties. For a signature \(L\) , an
\textbf{\(L\)-structure} is a \(\mathsf{Set}\)-valued functor on \(L\).
To each signature \(L\), a particular relation between two variable
\(L\)-structures, called \(L\)-equivalence, is defined. Two
\(L\)-structures M, N, are \(L\)-equivalent iff there is a so-called
\(L\)-equivalence span \(M\leftarrow P\rightarrow N\) between them;
here, the arrows are ordinary natural transformations, required to
satisfy a certain property called ``fiberwise surjectivity''.

The slogan of the work {[}M2{]}, {[}M3{]} on FOLDS is that \emph{all
meaningful properties of \(L\)-structures are invariant under
\(L\)-equivalence}. As with all slogans, it is both a normative
statement (``you should not look at properties of \(L\)-structures that
are not invariant under \(L\)-equivalence''), and a statement of fact,
namely that the ``interesting'' properties of \(L\)-structures are in
fact invariant under \(L\)-equivalence. (For some slogans, the
``statement of fact'' may be false.) The usual concepts of
``equivalence'' in category theory, including the higher dimensional
ones such as ``biequivalence'', are special cases of \(L\)-equivalence,
upon suitable, and natural, choices of the signature \(L\); {[}M3{]}
works out several examples of this. Thus, in these cases, the slogan
above becomes a tenet widely held true by category theorists. I claim
its validity in the generality stated above.

The main effort in {[}M3{]} goes into specifying a language, First Order
Logic with Dependent Sorts, and showing that the first order properties
invariant under \(L\)-equivalence are precisely the ones that can be
defined in FOLDS. In this paper, the language of FOLDS plays no role.
The concepts of ``FOLDS-signature'' and ``FOLDS-equivalence'' are fully
described in section 1 of this paper.

The definition of \textbf{multitopic \(\omega\)-category} goes, in
outline, as follows. For an arbitrary multitope \(\Sigma\) of dimension
\(\geqslant2\), for a multitopic set \(S\), for a pasting diagram
\(\alpha\) in \(S\) of shape the domain of \(\Sigma\) and a cell \(a\)
in \(S\) of the shape the codomain of \(\Sigma\), such that \(a\) and
\(\alpha\) are parallel, we define what it means to say that a is a
\textbf{composite} of \(\alpha\). First, we define an auxiliary FOLDS
signature \(L\langle \Sigma\rangle\) extending \(\mathsf{Mlt}\), the
signature of multitopic sets. Next, we define structures
\(S\langle a\rangle\) and \(S\langle \alpha\rangle\), both of the
signature \(L\langle \Sigma\rangle\), the first constructed from the
data \(S\) and \(a\) , the second from \(S\) and \(\alpha\), both
structures extending \(S\) itself. We say that \(a\) is a composite of
\(\alpha\) if there is a FOLDS-equivalence-span \(E\) between
\(S\langle a\rangle\) and \(S\langle \alpha\rangle\) that restricts to
the identity equivalence-span from \(S\) to \(S\). Below, I'll refer to
\(E\) as an \textbf{equipment} for a being a composite of \(\alpha\). A
multitopic set is a \textbf{mulitopic \(\omega\)-category} iff every
pasting diagram \(\alpha\) in it has at least one composite.

The analog of the universal arrows in the Baez/Dolan style definition is
as follows. A \textbf{universal arrow} is defined to be an arrow of the
form \(b\colon\alpha\to a\) where \(a\) is a composite of \(\alpha\) via
an equipment \(E\) that relates \(b\) with the identity arrow on \(a\):
in turn, the identity arrow on \(a\) is any composite of the empty
pasting diagram of dimension \(\dim(a)+1\) based on \(a\). Note that the
main definition does \emph{not} go through first defining ``universal
arrow''.

A new feature in the present treatment is that it aims directly at weak
\(\omega\)-categories; the finite dimensional ones are obtained as
truncated versions of the full concept. The treatment in {[}B/D{]}
concerns finite dimensional weak categories. It is important to
emphasize that a multitopic \(\omega\)-category is still just a
multitopic set with additional properties, but with no extra data.

The definition of ``multitopic \(\omega\)-category'' is given is section
5; it uses sections 1, 2 and 4, but not section 3.

The second main thing done in this paper is the definition of
\(\mathsf{MltOmegaCat}\). This is a particular large multitopic set. Its
definition is completed only by the end of the paper. The 0-cells of
\(\mathsf{MltOmegaCat}\) are the samll multitopic \(\omega\)-categories,
defined in section 5. Its 1-cells, which we call 1-transfors (thereby
borrowing, and altering the meaning of, a term used by Sjoerd Crans
{[}Cr{]}) are what stand for ``morphisms'', or ``functors'', of
multitopic \(\omega\)-categories. For instance, in the \(2\)-dimensional
case, multitopic \(2\)-categories correspond to ordinary bicategories by
a certain process of ``cleavage'', and the 1-transfors correspond to
homomorphisms of bicategories {[}Be{]}. There are \(n\)-dimensional
transfors for each \(n\) in \(\mathbb{N}\). For each multitope (that is,
``shape'' of a higher dimensional cell) \(PI\), we have the
\textbf{\(PI\)-transfors}, the cells of shape \(PI\) in
\(\mathsf{MltOmegaCat}\).

For each fixed multitope \(PI\), a \(PI\)-transfor is a
\emph{\(PI\)-colored multitopic set} with additional properties.
``\(PI\)-colored multitopic sets'' are defined in section 3; when \(PI\)
is the unique zero-dimensional multitope, \(PI\)-colored multitopic sets
are the same as ordinary multitopic sets. Thus, the definition of a
transfor of an arbitrary dimension and shape is a generalization of that
of ``multitopic \(\omega\)-category''; the additional properties are
also similar, they being defined by FOLDS-based universal properties.
There is one new element though. For \(\dim(PI)\geqslant2\) , the
concept of \(PI\)-transfor involves a universal property which is an
\(\omega\)-dimensional, FOLDS-style generalization of the concept of
right Kan-extension (right lifting in the terminology used by Ross
Street). This is a ``right-adjoint'' type universal property, in
contrast to the ``left-adjoint'' type involved in the concept of
composite (which is a generalization of the usual tensor product in
modules).

The main theorem, stated but not proved here, is that
\(\mathsf{MltOmegaCat}\) is a multitopic \(\omega\)-category.

The material in this paper has been applied to give formulations of
\(\omega\)-dimensional versions of various concepts of homotopy theory;
details will appear elesewhere.

References:

{[}B/D{]} J. C. Baez and J. Dolan, ``Higher-dimensional algebra III.
\(n\)-categories and the algebra of opetopes''. \emph{Advances in
Mathematics} \textbf{135} (1998), 145--206.

{[}Be{]} J. Benabou, ``Introduction to bicategories''. In: \emph{Lecture
Notes in Mathematics} \textbf{47}, Springer, Berlin, 1967, pp.\ 1--77.

{[}Cr{]} S. Crans, ``Localizations of transfors''. \emph{Macquarie
Mathematics Reports} no. \textbf{98/237}.

{[}H/M/P{]} C. Hermida, M. Makkai and J. Power, ``On weak higher
dimensional categories I'', \emph{Jour.\ Pure Appl.\ Alg.}, \textbf{154} (2000), 221--246.

{[}L{]} J. Lambek, ``Deductive systems and categories II''.
\emph{Lecture Notes in Mathematics} \textbf{86} Springer,
Berlin, 1969, pp.\ 76--122

{[}M2{]} M. Makkai, ``Towards a categorical foundation of mathematics''.
In: \emph{Logic Colloquium '95} (J. A. Makowski and E. V. Ravve,
editors) \emph{Lecture Notes in Logic} \textbf{11}, Springer, Berlin, 1998.

{[}M3{]} M. Makkai, ``First Order Logic with Dependent Sorts'',
available as {\rm \href{https://www.math.mcgill.ca/makkai/folds/foldsinpdf/FOLDS.pdf}{https://}}
{\rm \href{https://www.math.mcgill.ca/makkai/folds/foldsinpdf/FOLDS.pdf}{www.math.mcgill.ca/makkai/folds/foldsinpdf/FOLDS.pdf}}
\end{quote}

\hypertarget{week143}{%
\section{December 29, 1999}\label{week143}}

Since this is the last Week of the millennium, I'll make sure to pack it
full of retrospectives and prognostications. But I'd like to start with
an update on something I discussed a while back.

By the way, please don't give me flak about how the millennium starts in
2001. I use the CE or Common Era system, which starts counting at the
year zero, not the AD system, which starts at the year one because it
was invented in 526 CE by Dennis the Diminutive, long before the number
zero caught on. In my opinion, the real ``millennium bug'' is that
anyone is still using the antiquated AD system!

Anyway\ldots.

In \protect\hyperlink{week73}{``Week 73''}, I mentioned a theory about
why molecules important in biology tend to come in a consistent
chirality, or handedness. For example, there's lots of dextrose in
nature --- this being the right-handed form of the sugar sucrose --- but
not much of its left-handed counterpart, levulose. It's no surprise that
\emph{one or the other} would dominate, but you might guess that
\emph{which one} was just an accident of history. After all, there's no
fundamental difference between right and left, right?

Or is there? Actually there is: the weak nuclear force distinguishes
between the two! So some people theorized that very slight differences
in energy levels, due to the weak force, favor the formation of
left-handed amino acids and right-handed sugars --- which is what we see
in nature.

Recently people have found evidence for a somewhat different version of
this theory:

\begin{enumerate}
\def\labelenumi{\arabic{enumi})}
\tightlist
\item
  Robert F. Service, ``Does life's handedness come from within?'',
  \emph{Science} \textbf{286} (November 12, 1999), 1282--1283.
\end{enumerate}

When radioactive atoms decay via the weak force, the electrons they
shoot off tend to have a left-handed spin. Could this affect the
handedness of molecules or crystals that happen to be forming in the
vicinity? Sodium chlorate is a chemical that can form both left-handed
and right-handed crystals, so researchers took some solutions of the
stuff and let them crystallize while blasting them with electrons formed
by the decay of radioactive strontium. Sure enough, this biased the
handedness of the crystals! Blasting the stuff with right-handed
positrons favored formation of crystals of the opposite handedness. The
strangest part was that the effect was even bigger than expected.

I still think the whole business is pretty iffy --- after all, the flux
of radiation in this experiment was a lot bigger than what we normally
see on earth. But it would sure be neat if the origin of chirality in
biology was related to the deeper mystery of chirality in particle
physics.

Okay, now for a little retrospective. Don't worry --- I won't list the
top ten developments in mathematical physics of the last millennium!
Instead, I just want to recommend two papers. First, an old paper by
Poincar\'e:

\begin{enumerate}
\def\labelenumi{\arabic{enumi})}
\setcounter{enumi}{1}
\tightlist
\item
  Henri Poincar\'e, ``The present and future of mathematical physics'',
  \emph{Bull. Amer. Math. Soc.} \textbf{12} (1906), 240--260. Reprinted
  as part of a retrospective issue of the \emph{Bull. Amer. Math.
  Soc}., \textbf{37} (2000), 25--38.   Also available at
  \href{https://www.ams.org/bull/2000-37-01/S0273-0979-99-00801-0/S0273-0979-99-00801-0.pdf}{\texttt{https://www.}}
\href{https://www.ams.org/bull/2000-37-01/S0273-0979-99-00801-0/S0273-0979-99-00801-0.pdf}{\texttt{ams.org/bull/2000-37-01/S0273-0979-99-00801-0/S0273-0979-99-00801-0.pdf}}
\end{enumerate}

This article is based on a speech he gave in 1904. After a fascinating
review of the development of mathematical physics, he makes some
accurate predictions about quantum mechanics and special relativity ---
but closes on a conservative note:

\begin{quote}
In what direction we are going to expand we are unable to foresee.
Perhaps it is the kinetic theory of gases that will forge ahead and
serve as a model for the others. In that case, the facts that appeared
simple to us at first will be nothing more than the resultants of a very
large number of elementary facts which the laws of probability alone
would induce to work toward the same end. A physical law would then
assume an entirely new aspect; it would no longer be merely a
differential equation, it would assume the character of a statistical
law.

Perhaps too we shall have to construct an entirely new mechanics, which
we can only just get a glimpse of, where, the inertia increasing with
the velocity, the velocity of light would be a limit beyond which it
would be impossible to go. The ordinary, simpler mechanics would remain
a first approximation since it would be valid for velocities that are
not too great, so that the old dynamics would be found in the new. We
should have no reason to regret that we believed in the older
principles, and indeed since the velocities that are too great for the
old formulas will always be exceptional, the safest thing to do in
practice would be to act as though we continued to believe in them. They
are so useful that a place should be saved for them. To wish to banish
them altogether would be to deprive oneself of a valuable weapon. I
hasten to say, in closing, that we are not yet at that pass, and that
nothing proves as yet that they will not come out of the fray victorious
and intact.
\end{quote}

The same issue of the AMS Bulletin also has a lot of other interesting
papers and book reviews from the last century, by folks like Birkhoff,
Einstein and Weyl.

The second paper I recommend is a new one by Rovelli:

\begin{enumerate}
\def\labelenumi{\arabic{enumi})}
\setcounter{enumi}{2}
\tightlist
\item
  Carlo Rovelli, ``The century of the incomplete revolution: searching
  for general relativistic quantum field theory'', \emph{J. Math. Phys.} 
  \textbf{41} (2000), 3776--3800.   Also available as
  \href{https://arxiv.org/abs/hep-th/9910131}{\texttt{hep-th/9910131}}.
\end{enumerate}

Let me just quote the abstract:

\begin{quote}
In fundamental physics, this has been the century of quantum mechanics
and general relativity. It has also been the century of the long search
for a conceptual framework capable of embracing the astonishing features
of the world that have been revealed by these two ``first pieces of a
conceptual revolution''. I discuss the general requirements on the
mathematics and some specific developments towards the construction of
such a framework. Examples of covariant constructions of (simple)
generally relativistic quantum field theories have been obtained as
topological quantum field theories, in nonperturbative zero-dimensional
string theory and its higher dimensional generalizations, and as spin
foam models. A canonical construction of a general relativistic quantum
field theory is provided by loop quantum gravity. Remarkably, all these
diverse approaches have turn out to be related, suggesting an intriguing
general picture of general relativistic quantum physics.
\end{quote}

Now for the prognostications. Since we should never forget that the
towering abstractions of mathematical physics are ultimately tested by
experiment, I'd like to talk about some interesting physics
\emph{experiments} that are coming up in the next millennium. These days
more and more interesting information about physics is coming from
astronomy, so I'll concentrate on work that lies on this interface.

In \protect\hyperlink{week80}{``Week 80''} I talked about how Gravity
Probe B will try to detect an effect of general relativity called
``frame-dragging'' caused by the earth's rotation. I also talked about
how LIGO --- the Laser Interferometric Gravitational Wave Observatory
--- will try to detect gravitational waves:

\begin{enumerate}
\def\labelenumi{\arabic{enumi})}
\setcounter{enumi}{3}
\tightlist
\item
  LIGO homepage, \url{http://www.ligo.caltech.edu/}
\end{enumerate}

If all works as planned, LIGO should be great for studying the final
death spirals of binary black holes and/or neutron stars. When it starts
taking data sometime around 2002, it should be able to detect the final
``chirp'' of gravitational radiation produced a pair of inspiralling
neutron stars in the Virgo Cluster, a cluster of galaxies about 15
megaparsecs away. Such an event would distort the spacetime metric
\emph{here} by only about \(1\) part in \(10^{21}\). This is why LIGO
needs to compare oscillations in the lengths of two arms of an
interferometer, each \(4\) kilometers long, with an accuracy of
\(10^{-16}\) centimeters: about one hundred-millionth of the diameter of
a hydrogen atom. To do this will require some \emph{very} clever tricks
to reduce noise.

As the experiment continues, they intend to improve the sensitivity
until it can detect distortions in the metric of only \(1\) part in
\(10^{22}\), and second-generation detectors should get to \(1\) part in
\(10^{23}\). At that point, we should be able to detect neutron star
``chirps'' from a distance of \(200\) megaparsecs. Events of this sort
should happen once or twice a year.

Since it's crucial to rule out spurious signals, LIGO will have two
detectors, one in Livingston, Louisiana and one in Hanford, Washington.
This should also allow us to tell where the gravitational waves are
coming from. And there are other gravitational wave detection projects
underway too! France and Germany are collaborating on a laser
interferometer called VIRGO, with arms 3 kilometers long, to be built in
Cascina, Italy:

\begin{enumerate}
\def\labelenumi{\arabic{enumi})}
\setcounter{enumi}{4}
\tightlist
\item
  VIRGO homepage, \url{https://www.virgo-gw.eu/}
\end{enumerate}
\noindent
Germany and Great Britain are collaborating on a 600-meter-long one
called GEO 600, to be built south of Hannover:

\begin{enumerate}
\def\labelenumi{\arabic{enumi})}
\setcounter{enumi}{5}
\tightlist
\item
  GEO 600 homepage, \url{https://www.geo600.org/}
\end{enumerate}
\noindent
The Japanese are working on one called TAMA 300, which is a 300-
meter-long warmup for a planned kilometer-long interferometer:

\begin{enumerate}
\def\labelenumi{\arabic{enumi})}
\setcounter{enumi}{6}
\tightlist
\item
  TAMA 300 homepage, \url{http://gwpo.nao.ac.jp/en/}
\end{enumerate}
\noindent
In addition, the Brazilian GRAVITON project is building something called
the Einstein Antenna, which uses mechanical resonance rather than
interferometry. The basic principle goes back to Joseph Weber's original
bar detectors, which tried to sense the vibrations of a 2-meter-long
aluminum cylinder induced by gravitational waves. But the design
involves lots of hot new technology: SQUIDS, buckyballs, and the like:

\begin{enumerate}
\def\labelenumi{\arabic{enumi})}
\setcounter{enumi}{7}
\tightlist
\item
  GRAVITON homepage,
  \href{https://web.archive.org/web/19991108211006/http://www.das.inpe.br/graviton/project.html}{\texttt{https://web.archive.org/web/19991108211006/}}
 \href{https://web.archive.org/web/19991108211006/http://www.das.inpe.br/graviton/project.html}{\texttt{http://www.das.inpe.br/graviton/project.html}}
\end{enumerate}

There are also other gravitational wave detectors being built\ldots{}
but ultimately, the really best ones will probably be built in outer
space. There are two good reasons for this. First, outer space is big:
when you're trying to detect very small distortions of the geometry of
spacetime, it helps to measure the distance between quite distant
points. Second, outer space is free of seismic noise and most other
sources of vibration. This is why people are working on the LISA project
--- the Laser Interferometer Space Antenna:

\begin{enumerate}
\def\labelenumi{\arabic{enumi})}
\setcounter{enumi}{8}
\tightlist
\item  NASA's homepage on the LISA project, \url{http://lisa.nasa.gov/}

  Wikipedia, Laser Interferometer Space Antenna, \href{https://en.wikipedia.org/wiki/Laser_Interferometer_Space_Antenna}{\texttt{https://en.wikipedia.org/wiki/}}
\href{https://en.wikipedia.org/wiki/Laser_Interferometer_Space_Antenna}{\texttt{Laser\_Interferometer\_Space\_Antenna}}
\end{enumerate}

The idea is to orbit 3 satellites in an equilateral triangle with sides
5 million kilometers long, and constantly measure the distance between
them to an accuracy of a tenth of an angstrom --- \(10^{-11}\) meters
--- using laser interferometry. (A modified version of the plan would
use 6 satellites.) The big distances would make it possible to detect
gravitational waves with frequencies of 0.0001 to 0.1 hertz, much lower
than the frequencies for which the ground-based detectors are optimized.
The plan involves a really cool technical trick to keep the satellites
from being pushed around by solar wind and the like: each satellite will
have a free-falling metal cube floating inside it, and if the satellite
gets pushed to one side relative to this mass, sensors will detect this
and thrusters will push the satellite back on course.

I don't think LISA has been funded yet, but if all goes well, it may fly
within 10 years or so. Eventually, a project called LISA 2 might be
sensitive enough to detect gravitational waves left over from the early
universe --- the gravitational analogue of the cosmic microwave
background radiation!

The microwave background radiation tells us about the universe when it
was roughly \(10^5\) years old, since that's when things cooled down
enough for most of the hydrogen to stop being ionized, making it
transparent to electromagnetic radiation. In physics jargon, that's when
electromagnetic radiation ``decoupled''. But the gravitational
background radiation would tell us about the universe when it was
roughly \(10^{-38}\) seconds old, since that's when gravitational
radiation decoupled. This figure could be way off due to physics we
don't understand yet, but anyway, we're talking about a window into the
\emph{really} early universe.

Actually, Mark Kamionkowski of Caltech has theorized that the European
Space Agency's ``Planck'' satellite may detect subtle hints of the
gravitational background radiation through its tendency to polarize the
microwave background radiation. You probably heard how COBE, the Cosmic
Background Explorer, detected slight anisotropies in the microwave
background radiation. Now people are going to redo this with much more
precision: while COBE had an angular resolution of 7 degrees, Planck
will have a resolution of 4 arcminutes. They hope to launch it in 2007:

\begin{enumerate}
\def\labelenumi{\arabic{enumi})}
\setcounter{enumi}{9}
\tightlist
\item
  Wikipedia, Planck (spacecraft), 
  \href{https://en.wikipedia.org/wiki/Planck_(spacecraft)}{\texttt{https://en.wikipedia.org/wiki/Planck\_}} \hfill \break
\href{https://en.wikipedia.org/wiki/Planck_(spacecraft)}{\texttt{(spacecraft)}}
\end{enumerate}

What else is coming up? Well, gravity people should be happy about the
new satellite-based X-ray telescopes, since these should be great for
looking at black holes. In July 1999, NASA launched one called
``Chandra''. (This is the nickname of Subrahmanyan Chandrasekhar, who
won the Nobel prize in 1983 for his work on stellar evolution, neutron
stars, black holes, and closed-form solutions of general relativity.)
The first pictures from Chandra are already coming out --- check out
this website:

\begin{enumerate}
\def\labelenumi{\arabic{enumi})}
\setcounter{enumi}{10}
\tightlist
\item
  Chandra homepage, \url{http://chandra.harvard.edu/}
\end{enumerate}
\noindent
On December 10th, the Europeans launched XMM, the ``X-ray Multi-Mirror
Mission'':

\begin{enumerate}
\def\labelenumi{\arabic{enumi})}
\setcounter{enumi}{11}
\tightlist
\item
  XMM homepage, \url{https://en.wikipedia.org/wiki/XMM-Newton}
\end{enumerate}
\noindent
This is a set of three X-ray telescopes that will have lower angular
resolution than Chandra, but 5--15 times more sensitivity. It'll also be
able to study X-ray spectra, thanks to a diffraction grating that
spreads the X-rays out by wavelength. And in January, the Japanese plan
to launch ASTRO-E, designed to look at shorter wavelength X-rays:

\begin{enumerate}
\def\labelenumi{\arabic{enumi})}
\setcounter{enumi}{12}
\tightlist
\item
  Wikipedia, Suzaku: ASTRO-E, \href{https://en.wikipedia.org/wiki/Suzaku_(satellite)#ASTRO-E}{\texttt{https://en.wikipedia.org/wiki/Suzaku\_}} \hfill \break
 \href{https://en.wikipedia.org/wiki/Suzaku_(satellite)#ASTRO-E}{\texttt{(satellite)\#ASTRO-E}}
\end{enumerate}

Taken together, this new generation of X-ray telescopes should tell us a
lot about the dynamics of the rapidly changing accretion disks of black
holes, where infalling gas and dust spirals in and heats up to the point
of emitting X-rays. They may also help us better understand the X-ray
afterglow of \(\gamma\)-ray bursters. As you probably have heard, these
rascals make ordinary supernovae look like wet firecrackers! Some folks
think they're caused when a supernova creates a black hole. But nobody
is sure.

Peering further into the future, here's a nice article about new
projects people are dreaming up to study physics using astronomy:

\begin{enumerate}
\def\labelenumi{\arabic{enumi})}
\setcounter{enumi}{13}
\tightlist
\item
  Robert Irion, ``Space becomes a physics lab'', \emph{Science} \textbf{286} (1999),
  2060--2062.
\end{enumerate}

In 2005 folks plan to launch GLAST, the Gamma-Ray Large Area Space
Telescope, designed to study \(\gamma\)-ray bursters and the like, and
also the Alpha Magnetic Spectrometer, designed to search for antimatter
in space. But there are also a bunch of interesting projects that are
still basically just a twinkle in someone's eye\ldots.

For example: OWL, the Orbiting Wide-Angle Light Collector, a pair of
satellites that would trace the paths of super-high-energy cosmic rays
through the earth's atmosphere. As I explained in
\protect\hyperlink{week81}{``Week 81''}, people have seen cosmic rays
with ridiculously high energies, like 320 Eev --- the energy of a
1-kilogram rock moving at 10 meters per second, all packed into one
particle. OWL would orbit the earth, watch these things, and figure out
where the heck they're coming from.

Or how about this: The Dark Matter Telescope! This would use
gravitational lensing to chart the ``dark matter'' which seems to
account for a good percentage of the mass in the universe --- if, of
course, dark matter really exists.

\begin{enumerate}
\def\labelenumi{\arabic{enumi})}
\setcounter{enumi}{14}
\tightlist
\item
  Wikipedia, \href{https://en.wikipedia.org/wiki/Vera_C._Rubin_Observatory}{\texttt{https://en.wikipedia.org/wiki/Vera\_C.\_Rubin\_Observatory}}
\end{enumerate}

Anyway, there should be a lot of exciting experiments coming up. But as
usual, the really exciting stuff will be the stuff we can't predict.

\hypertarget{week144}{%
\section{January 21, 2000}\label{week144}}

Since this is the first Week of the new millennium, I'd like to start
with a peek into the future. Not just the next hundred or thousand
years, either --- I'm sick of short-term planning. No, I'd like to talk
about the next few \emph{billion} years.

As you've probably all heard, if we don't do anything about it, the Sun
will turn into a red giant in about 5 billion years. If we get our act
together, we should have plenty of time to deal with this problem. But
when planning for the far future, it's dangerous to be too parochial!
Events outside our solar system can also affect us. For example, a
nearby supernova could be a real bummer. It wouldn't be the first time:
it seems that about 340,000 years ago there was one only 180 lightyears
away. At this distance it would have been as bright as a full moon, and
its X-rays and \(\gamma\)-rays would have stripped off the Earth's ozone
layer pretty badly for a while. A closer one could be a lot worse.

How we do know about this supernova? It's an interesting story. We
happen to live in a region of space called the Local Bubble, about 300
lightyears across, in which the interstellar gas is hotter and 5 to 10
times less dense than the surrounding stuff. People wondered about the
origin of this bubble until they studied a pulsar called Geminga about
300 lightyears away from us. Pulsars are rapidly spinning neutron stars
formed by supernovae, and by studying their spin rate and the rate their
spin is slowing down, you can guess when they were formed. Geminga turns
out to be about 340,000 years old. It's moving away from us at a known
rate, so back then it would have been 180 lightyears away --- in just
about the right place for a supernova to have created a shock wave
forming the Local Bubble.

I don't know the best place to read about the Local Bubble, but this
sounds promising:

\begin{enumerate}
\def\labelenumi{\arabic{enumi})}
\tightlist
\item
  M. J. Freyberg and J. Trumper, eds., \emph{The Local Bubble and
  Beyond}, Springer Lecture
  Notes in Physics \textbf{506}, Springer, Berlin, 1998.
\end{enumerate}

Looking further afield, we should also watch out for the health of the
Milky Way as a whole:

\begin{enumerate}
\def\labelenumi{\arabic{enumi})}
\setcounter{enumi}{1}
\item
  Robert Irion, ``A crushing end for our galaxy'', \emph{Science}
  \textbf{287} (2000), 62--64.
\item
  Roland Buser, ``The formation and early evolution of the Milky Way
  galaxy'', \emph{Science} \textbf{287} (2000), 69--74.
\end{enumerate}
\noindent
It now appears that the Milky Way, like most big spiral galaxies, was
built up by a gradual merger of smaller clouds of stars and gas. And it
seems this process is not finished. In 1994 people found a small galaxy
orbiting the Milky Way, almost hidden behind the dense dust clouds in
the galactic center. Called the Sagittarius dwarf galaxy, it is only
about \(1/1000\)th the mass of the Milky Way. Its eccentric orbit about
our galaxy is strewn with stars pulled away from it by tidal forces, and
it may have already passed through the outer parts of our galaxy's disk
several times. It may not survive the next pass, due in about 750
million years.

But that's not all. The Large and Small Magellanic Clouds, visible to
the naked eye in the Southern Hemisphere, are also dwarf galaxies
orbiting ours. They are considerably larger than the Sagittarius dwarf
galaxy. And they're not just orbiting the Milky Way: they are gradually
spiralling in and getting torn apart. If nothing interrupts this
process, they'll crash into our galaxy in about 10 billion years. When
when this happens, the shock waves from colliding gas should create
enough new stars to make our galaxy shine about 25\% brighter for the
next several hundred million years! This could prove quite a nuisance in
these parts.

But again, we should not make the mistake of parochialism: dangers from
afar may prove more urgent than those in our neighborhood. The dwarf
galaxies near us are nothing compared to Andromeda. This spiral galaxy
is twice the size of ours, about 2.5 million light years away, and
clearly visible from the Northern Hemisphere. Unfortunately, it's also
heading towards us at a speed of 140 kilometers per second! As it comes
closer, gravitational attraction will speed it up, so it may hit our
galaxy --- or at least come close --- in only 3 billion years. If this
happens, the two galaxies will first whiz past or through each other,
but then their gravitational attraction will pull them back together,
and after 1 or 2 billion more years they should coalesce into a single
big elliptical galaxy. Direct hits between individal stars are unlikely,
but many existing stars will be hurled out into intergalactic space, and
many new stars will be born as gas clouds collide.

You may think that I'm joking when I speak of planning ahead for such
events, but I'm not. We have plenty of time, so it's not very urgent ---
but it's not too soon to start thinking about these things. And if you
think it's hopelessly beyond our powers to deal with a collision of
galaxies, please remember that 3 billion years ago we were single-celled
organisms. With any luck, our abilities 3 billion years from now should
compare to our present abilities as our present abilities compare to
those of microorganisms! And if life on Earth screws up and dies out,
well, there are plenty of other planets out there.

By the way, while we're discussing matters galactic, remember how last
Week I said that the X-ray telescope Chandra has recently started taking
data? Well, the interesting news is already coming in! For a long time
people have wondered about the origin of the ``X-ray background
radiation'': a diffuse X-ray glow that covers the whole sky. On
Thursday, astronomers using Chandra discovered that most of this
radiation actually comes from about 70 million individual point sources!
Apparently, many of these are supermassive black holes at the center of
galaxies. There's already a lot of evidence for such black holes ---
which seem to power quasars and other active galactic nuclei --- but
it's delightful to find them in such large numbers. It might even be
taken as evidence for Smolin's hypothesis that the universe is optimized
for black hole production thanks to a process of Darwinian evolution
(see \protect\hyperlink{week31}{``Week 31''} and
\protect\hyperlink{week33}{``Week 33''} for details).

For more, try this:

\begin{enumerate}
\def\labelenumi{\arabic{enumi})}
\setcounter{enumi}{3}
\tightlist
\item
  ``Chandra resolves cosmic X-ray glow and finds mysterious new
  sources'', available at
 \href{http://chandra.harvard.edu/press/00_releases/press_011400bg.html}
  {\texttt{http://chandra.harvard.edu/press/00\_releases/press\_011400bg.html}}
\end{enumerate}
\noindent
You should also check out the Chandra website for nice new pictures of
the black holes at the center of the Milky Way and Andromeda.

Okay\ldots{} I've been sort of goofing off in the last few Weeks, but
now I want to return to some hardcore mathematics. In particular, I want
to talk about \(n\)-categories and homotopy theory, so I'm going to pick
up ``The Tale of \(n\)-Categories'' roughly where I left off in
\protect\hyperlink{week100}{``Week 100''}, and start connecting it to
the little introduction to homotopy theory I gave from
\protect\hyperlink{week115}{``Week 115''} to
\protect\hyperlink{week119}{``Week 119''}.

As I've said many times, the goal of \(n\)-category theory is to
eliminate equations from mathematics, or at least to be able to postpone
pretending that isomorphisms are equations for as long as you like. I've
repeatedly described the practical benefits of this, so I won't bother
doing so again --- I'll assume you're convinced of it!

To achieve this goal, an \(n\)-category is supposed to be some sort of
algebraic structure with objects, morphisms between objects,
\(2\)-morphisms between morphisms, and so on up to and including
\(n\)-morphisms, with various ways of composing all these guys. The idea
is then that we should never assert that two \(j\)-morphisms are
\emph{equal} except for \(j = n\). Instead, we should just specify an
\emph{equivalence} between them. An ``equivalence'' is a bit like an
isomorphism, but it's defined recursively from the top down. An
\(n\)-morphism is an equivalence iff it's an isomorphism, that is, iff
it's invertible. But for \(j < n\), a \(j\)-morphism is an equivalence
if it's invertible \emph{up to equivalence}.

There are various competing definitions of \(n\)-category at present,
but the key idea behind all the definitions of \emph{weak}
\(n\)-category is that the ways of composing \(j\)-morphisms should
satisfy associativity and all the other usual laws only up to
equivalence. For example, suppose we have some morphisms
\(a\colon w\to x\), \(b\colon w\to y\) and \(c\colon y\to z\) in a
\(1\)-category. Then associativity holds ``on the nose'', i.e., as an
equation: \[a(bc) = (ab)c.\] In a \(1\)-category there is no opportunity
for ``weakening'' this law. But in a weak \(2\)-category, associativity
holds only up to equivalence. In other words, we have an invertible
\(2\)-morphism called the ``associator''
\[A_{a,b,c}: (ab)c \Rightarrow a(bc)\] taking the part of the above
equation.

But there's a catch: when we replace equational laws by equivalences
this way, the equivalences need to satisfy laws of their own, or it
becomes impossible to work with them. These laws are called ``coherence
laws''. For the associator, the necessary coherence law is called the
pentagon equation. It says that this diagram commutes: 
\[
 \begin{tikzpicture}
\scalebox{1}{
    \node (mr) at (2,2) {$(a (b c)) d$};
    \node (t) at (0,3.5) {$((a b) c) d$};
    \node (ml) at (-2,2) {$(a b) (c d)$};
    \node (bl) at (-1.3,0) {$a (b (c d))$};
    \node (br) at (1.3,0) {$a ((b c) d)$};
    \draw[double,double equal sign distance,-implies] (t) to (ml);
    \draw[double,double equal sign distance,-implies] (ml) to (bl);
    \draw[double,double equal sign distance,-implies] (t) to (mr);
    \draw[double,double equal sign distance,-implies] (mr) to (br);
    \draw[double,double equal sign distance,-implies] (br) to (bl);
}
  \end{tikzpicture}
\]
I haven't labelled the double arrows here, but they are all
\(2\)-morphisms built from the associator in obvious ways\ldots{}
obvious if you know about \(2\)-categories, at least. The pentagon
equation says that the two basic ways of going from \(((ab)c)d\) to
\(a(b(cd))\) by rebracketing are equal to each other. But in fact,
Mac Lane's ``coherence theorem'' says that given the pentagon equation,
you can rebracket composites of arbitrarily many morphisms using the
associator over and over to your heart's content, and you'll never get
into trouble: all the ways of going from one bracketing to another are
equal.

In a weak \(3\)-category, the pentagon equation is replaced by a
\(3\)-morphism called the ``pentagonator''. This in turn satisfies a new
coherence law of its own, which I can't easily draw for you, because
doing so requires a \(3\)-dimensional diagram in the shape of a
polyhedron with 14 vertices, called the ``associahedron''.

As you might fear, this process never stops: there's an infinite list of
``higher coherence laws'' for associativity, which can be represented as
higher-dimensional associahedra. They were discovered by James Stasheff
around 1963. Here are the original papers:

\begin{enumerate}
\def\labelenumi{\arabic{enumi})}
\setcounter{enumi}{4}
\tightlist
\item
  James Stasheff, ``Homotopy associativity of H-spaces I'', \emph{Trans.
  Amer. Math. Soc.} \textbf{108} (1963), 275--292.

James Stasheff, ``Homotopy associativity of H-spaces II'', \emph{Trans.
Amer. Math. Soc.} \textbf{108} (1963), 293--312.
\end{enumerate}
\noindent
Personally, I find his book a lot easier to read:

\begin{enumerate}
\def\labelenumi{\arabic{enumi})}
\setcounter{enumi}{5}
\tightlist
\item
  James Stasheff, \emph{H-spaces from a Homotopy Point of View},
  Springer Lecture Notes in Mathematics \textbf{161}, Springer,
  Berlin, 1970.
\end{enumerate}

There's a wealth of interesting combinatorics lurking in the
associahedra. To my shame, I realize that I've never discussed this
stuff, so I'd better say a bit about it. Then next Week I'll return to
my real goal, which is to explain how you can use homotopy theory to
understand coherence laws. With any luck, I'll get around to telling you
all sorts of wonderful stuff about Postnikov towers, the cohomology of
Eilenberg--Mac Lane spaces, and so on. We'll see. So much math, so little
time\ldots.

Okay, here's how you build an associahedron.

First I'll describe the vertices, because they're very simple: the
correspond to all the ways of bracketing a string of \(n\) letters.
Well, that's a bit vague, so I'll do an example. Suppose \(n = 4\). Then
we get 5 bracketings:

\begin{itemize}
\tightlist
\item
  \(((ab)c)d\)
\item
  \((a(bc))d\)
\item
  \((ab)(cd)\)
\item
  \(a((bc)d)\)
\item
  \(a(b(cd))\)
\end{itemize}

These are exactly the vertices of the pentagon I drew earlier! And this
how it always works: the bracketings of \(n\) letters are the vertices
of the \((n-2)\)-dimensional associahedron. This should not be
surprising, since associativity is all about bracketing.

More precisely, we're interested in the bracketings of \(n\) letters
that correspond to binary planar trees with \(n\) leaves. For example,
when \(n = 4\): \[
  \begin{gathered}
    \begin{tikzpicture}[scale=0.7]
      \node[label=above:{$a$}] at (0,0) {};
      \node[label=above:{$b$}] at (1,0) {};
      \node[label=above:{$c$}] at (2,0) {};
      \node[label=above:{$d$}] at (3,0) {};
      \draw[thick] (0,0) to (1.5,-3);
      \draw[thick] (1,0) to (0.5,-1);
      \draw[thick] (2,0) to (1,-2);
      \draw[thick] (3,0) to (1.5,-3);
      \node at (4,-1.5) {$((ab)c)d$};
    \end{tikzpicture}
    \quad
    \begin{tikzpicture}[scale=0.7]
      \node[label=above:{$a$}] at (0,0) {};
      \node[label=above:{$b$}] at (1,0) {};
      \node[label=above:{$c$}] at (2,0) {};
      \node[label=above:{$d$}] at (3,0) {};
      \draw[thick] (0,0) to (1.5,-3);
      \draw[thick] (1,0) to (1.5,-1);
      \draw[thick] (2,0) to (1,-2);
      \draw[thick] (3,0) to (1.5,-3);
      \node at (4,-1.5) {$(a(bc))d$};
    \end{tikzpicture}
  \\\begin{tikzpicture}[scale=0.7]
      \node[label=above:{$a$}] at (0,0) {};
      \node[label=above:{$b$}] at (1,0) {};
      \node[label=above:{$c$}] at (2,0) {};
      \node[label=above:{$d$}] at (3,0) {};
      \draw[thick] (0,0) to (1.5,-3);
      \draw[thick] (1,0) to (0.5,-1);
      \draw[thick] (2,0) to (2.5,-1);
      \draw[thick] (3,0) to (1.5,-3);
      \node at (4,-1.5) {$(ab)(cd)$};
    \end{tikzpicture}
  \\\begin{tikzpicture}[scale=0.7]
      \node[label=above:{$a$}] at (0,0) {};
      \node[label=above:{$b$}] at (1,0) {};
      \node[label=above:{$c$}] at (2,0) {};
      \node[label=above:{$d$}] at (3,0) {};
      \draw[thick] (0,0) to (1.5,-3);
      \draw[thick] (1,0) to (2,-2);
      \draw[thick] (2,0) to (1.5,-1);
      \draw[thick] (3,0) to (1.5,-3);
      \node at (4,-1.5) {$a((bc)d)$};
    \end{tikzpicture}
    \quad
    \begin{tikzpicture}[scale=0.7]
      \node[label=above:{$a$}] at (0,0) {};
      \node[label=above:{$b$}] at (1,0) {};
      \node[label=above:{$c$}] at (2,0) {};
      \node[label=above:{$d$}] at (3,0) {};
      \draw[thick] (0,0) to (1.5,-3);
      \draw[thick] (1,0) to (2,-2);
      \draw[thick] (2,0) to (2.5,-1);
      \draw[thick] (3,0) to (1.5,-3);
      \node at (4,-1.5) {$a(b(cd))$};
    \end{tikzpicture}  
  \end{gathered}
\]
\noindent
We can think of these trees as recording the \emph{process} of
multiplying \(n\) things, with time marching down the page.

How many binary planar trees with \(n\) leaves are there, anyway? Well,
the answer is called the \((n-1)\)st Catalan number. These numbers were
first discovered by Euler, but they're named after Eugene Catalan, who
discovered their relation to binary trees. Here they are, starting from
the 0th one:
\[1, 1, 2, 5, 14, 42, 132, 429, 1430, 4862, 16796, 58786, 208012, 742900,
\ldots\] The \(n\)th Catalan number is also the number of ways of taking
a regular \((n+2)\)-gon and chopping it into triangles by connecting the
vertices by line segments that don't cross each other. It's also the
number of ways of getting from a street corner in Manhattan to another
street corner that's \(n\) blocks north and \(n\) blocks east, always
driving north or east, but making sure that at no stage have you gone a
greater total distance north than east. Get it? No? Maybe a picture will
help! When \(n = 3\), there are 5 ways: \[
  \begin{gathered}
    \begin{tikzpicture}
      \foreach \x in {0,1,2,3}
        \foreach \y in {0,1,2,3}
          \node (\x\y) at (\x,\y) {$\bullet$};
      \draw[ultra thick,-latex] (00) edge (10) (10) edge (20) (20) edge (30)
        (30) edge (31) (31) edge (32) (32) edge (33);
    \end{tikzpicture}
    \qquad\qquad
    \begin{tikzpicture}
      \foreach \x in {0,1,2,3}
        \foreach \y in {0,1,2,3}
          \node (\x\y) at (\x,\y) {$\bullet$};
      \draw[ultra thick,-latex] (00) edge (10) (10) edge (20)
        (20) edge (21)
        (21) edge (31)
        (31) edge (32) (32) edge (33);
    \end{tikzpicture}
  \\[2em]
    \begin{tikzpicture}
      \foreach \x in {0,1,2,3}
        \foreach \y in {0,1,2,3}
          \node (\x\y) at (\x,\y) {$\bullet$};
      \draw[ultra thick,-latex] (00) edge (10) (10) edge (20)
        (20) edge (21) (21) edge (22)
        (22) edge (32)
        (32) edge (33);
    \end{tikzpicture}
  \\[2em]
    \begin{tikzpicture}
      \foreach \x in {0,1,2,3}
        \foreach \y in {0,1,2,3}
          \node (\x\y) at (\x,\y) {$\bullet$};
      \draw[ultra thick,-latex] (00) edge (10)
        (10) edge (11)
        (11) edge (21) (21) edge (31)
        (31) edge (32) (32) edge (33);
    \end{tikzpicture}
    \qquad\qquad
    \begin{tikzpicture}
      \foreach \x in {0,1,2,3}
        \foreach \y in {0,1,2,3}
          \node (\x\y) at (\x,\y) {$\bullet$};
      \draw[ultra thick,-latex] (00) edge (10)
        (10) edge (11)
        (11) edge (21)
        (21) edge (22)
        (22) edge (32)
        (32) edge (33);
    \end{tikzpicture}
  \end{gathered}
\]

I leave it as a puzzle for you to understand why all these things are
counted by the Catalan numbers. If you want to see nicer pictures of all
these things, go here:

\begin{enumerate}
\def\labelenumi{\arabic{enumi})}
\setcounter{enumi}{6}
\tightlist
\item
  Robert M. Dickau, ``Catalan numbers'',
  \url{https://www.robertdickau.com/catalan.html}
\end{enumerate}
\noindent
For more problems whose answer involves the Catalan numbers, try this:

\begin{enumerate}
\def\labelenumi{\arabic{enumi})}
\setcounter{enumi}{7}
\tightlist
\item
  Kevin Brown, ``The meanings of Catalan numbers'',
  \href{https://web.archive.org/web/20171109040813/http://mathforum.org/kb/message.jspa?messageID=22219}{\texttt{https://web.archive.org/}}
\href{https://web.archive.org/web/20171109040813/http://mathforum.org/kb/message.jspa?messageID=22219}{\texttt{20171109040813/http://mathforum.org/kb/message.jspa?messageID=22219}}
\end{enumerate}

To figure out a formula for the Catalan numbers, we can use the
technique of generating functions:

\begin{enumerate}
\def\labelenumi{\arabic{enumi})}
\setcounter{enumi}{8}
\tightlist
\item
  Herbert Wilf, \emph{Generatingfunctionology}, Academic Press, New York,
  1994. Also available at
  \url{http://www.cis.upenn.edu/~wilf/}
\end{enumerate}

Briefly, the idea is to make up a power series \(T(x)\) where the
coefficient of \(x^n\) is the number of \(n\)-leaved binary trees. Since
by some irritating accident of history people call this the \((n-1)\)st
Catalan number, we have: \[T(x) = \sum_{n\in\mathbb{N}} C_{n-1} x^n\] We
can do this trick whenever we're counting how many structures of some
sort we can put on an \(n\)-element set. Nice operations on structures
correspond to nice operations on formal power series. Using this
correspondence we can figure out the function \(T\) and then do a Taylor
expansion to determine the Catalan numbers. Instead of explaining the
theory of how this all works, I'll just demonstrate it as a kind of
magic trick.

So: what is a binary tree? It's either a binary tree with one leaf (the
degenerate case) or a pair of binary trees stuck together. Now let's
translate this fact into an equation: \[T = x + T^2\] Huh? Well, in this
game ``plus'' corresponds to ``or'', ``times'' corresponds to ``and'',
and the power series ``\(x\)'' is the generating function for binary
trees with one leaf. So this equation really just says ``a binary tree
equals a binary tree with one leaf or a binary tree and a binary tree''.

Next, let's solve this equation for \(T\). It's just a quadratic
equation, so any high school student can solve it:
\[T = \frac{1-\sqrt{1-4x}}{2}.\] Now if we do a Taylor expansion we get
\[T = x + x^2 + 2x^3 + 5x^4 + 14x^5 + 42x^6 + \ldots\] Lo and behold ---
the Catalan numbers! If we're a bit smarter and use the binomial theorem
and mess around a bit, we get a closed-form formula for the Catalan
numbers: \[C_n = \frac{\binom{2n}{n}}{n + 1}\] Neat, huh? If you want to
understand the category-theoretic foundations of this trick, read about
Joyal's concept of ``species''. This makes precise the notion of a
``structure you can put on a finite set''. For more details, try:

\begin{enumerate}
\def\labelenumi{\arabic{enumi})}
\setcounter{enumi}{9}
\item
  Andre Joyal, ``Une theorie combinatoire des series formelles'',
  \emph{Adv. Math.} \textbf{42} (1981), 1--82.
\item
  F. Bergeron, G. Labelle, and P. Leroux, \emph{Combinatorial Species
  and Tree-Like Structures}, Cambridge U.\ Press, Cambridge, 1998.
\end{enumerate}

Anyway, now we know how many vertices the associahedron has. But what
about all the higher-dimensional faces of the associahedron? There's a
lot to say about this, but it's basically pretty simple: all the faces
of the \((n-2)\)-dimensional associahedron correspond to planar trees
with \(n\) leaves. It gets a little tricky to draw using ASCII, so I'll
just do the case \(n = 3\). The \(1\)-dimensional associahedron is the
unit interval, and in terms of trees it looks like this: \[
  \begin{tikzpicture}[scale=0.7]
    \begin{scope}
      \node[label=above:{$a$}] at (0,0) {};
      \node[label=above:{$b$}] at (1,0) {};
      \node[label=above:{$c$}] at (2,0) {};
      \draw[thick] (0,0) to (1,-2);
      \draw[thick] (1,0) to (0.5,-1);
      \draw[thick] (2,0) to (1,-2);
    \end{scope}
    \begin{scope}[shift={(3.5,0)}]
      \node[label=above:{$a$}] at (0,0) {};
      \node[label=above:{$b$}] at (1,0) {};
      \node[label=above:{$c$}] at (2,0) {};
      \draw[thick] (0,0) to (1,-2);
      \draw[thick] (1,0) to (1,-2);
      \draw[thick] (2,0) to (1,-2);
    \end{scope}
    \begin{scope}[shift={(7,0)}]
      \node[label=above:{$a$}] at (0,0) {};
      \node[label=above:{$b$}] at (1,0) {};
      \node[label=above:{$c$}] at (2,0) {};
      \draw[thick] (0,0) to (1,-2);
      \draw[thick] (1,0) to (1.5,-1);
      \draw[thick] (2,0) to (1,-2);
    \end{scope}
    \draw[ultra thick] (1,-3) to (8,-3);
    \node at (1,-3) {$\bullet$};
    \node at (8,-3) {$\bullet$};
  \end{tikzpicture}
\] Over at the left end of the interval we have the binary tree
corresponding to first composing \(a\) and \(b\), and then composing the
result with \(c\). At the right end, we have the binary tree
corresponding to first composing \(b\) and \(c\), and then composing the
result with \(a\). In the middle we have a ternary tree that corresponds
to \emph{simultaneously} composing \(a\), \(b\), and \(c\).

Actually, we can think of any point in the \((n-2)\)-dimensional
associahedron as an \(n\)-leaved tree whose edges have certain specified
lengths, so as you slide your finger across the \(1\)-dimensional
associahedron above, you can imagine the left-hand tree continuously
``morphing'' into the right-hand one. In this way of thinking, each
point of the associahedron corresponds to a particular \(n\)-ary
operation: a way of composing \(n\) things. To make this precise one
must use the theory of ``operads''. The theory of operads is really the
royal road to understanding \(n\)-categories, coherence laws, and their
relation to homotopy theory\ldots. But here, alas, I must stop.

\begin{center}\rule{0.5\linewidth}{0.5pt}\end{center}

Footnote --- If you want to know more about the deep inner meaning of
the Catalan numbers, try this:

\begin{enumerate}
\def\labelenumi{\arabic{enumi})}
\setcounter{enumi}{11}
\tightlist
\item
  Richard P. Stanley, \emph{Enumerative Combinatorics}, volume 2,
  Cambridge U.\ Press, Cambridge, 1999, pp.~219--229.
\end{enumerate}
\noindent
It lists 66 different combinatorial interpretations of these numbers! As
an exercise, it urges you to prove that they all work, ideally by
finding 4290 ``simple and elegant'' bijections between the various sets
being counted.

(Thanks go to my pal Bill Schmitt for mentioning this reference.)

\hypertarget{week145}{%
\section{February 9, 2000}\label{week145}}

I know I promised to talk about homotopy theory and \(n\)-categories,
but I've gotten sidetracked into thinking about projective planes, so
I'll talk about that this Week and go back to the other stuff later.
Sorry, but if I don't talk about what intrigues me at the instant I'm
writing this stuff, I can't get up the energy to write it.

So:

There are many kinds of geometry. After Euclidean geometry, one of the
first to become popular was projective geometry. Projective geometry is
the geometry of perspective. If you draw a picture on a piece of paper
and view it from a slant, distances and angles in the picture will get
messed up --- but lines will still look like lines. This kind of
transformation is called a ``projective transformation''. Projective
geometry is the study of those aspects of geometry that are preserved by
projective transformations.

Interestingly, \(2\)-dimensional projective geometry has some curious
features that don't show up in higher dimensions. To explain this, I
need to tell you about projective planes.

I talked a bit about projective planes in
\protect\hyperlink{week106}{``Week 106''}. The basic idea is to take the
ordinary plane and add some points at infinity so that every pair of
distinct lines intersects in exactly one point. Lines that were parallel
in the ordinary plane will intersect at one of the points at infinity.
This simplifies the axioms of projective geometry.

But what exactly do I mean by ``the ordinary plane''? Well, ever since
Descartes, most people think of the plane as \(\mathbb{R}^2\), which
consists of ordered pairs of real numbers. But algebraists also like to
use \(\mathbb{C}^2\), consisting of ordered pairs of complex numbers.
For that matter, you could take any field \(\mathbb{F}\) --- like the
rational numbers, or the integers modulo a prime --- and use
\(\mathbb{F}^2\). Algebraic geometers call this sort of thing an
``affine plane''.

A projective plane is a bit bigger than an affine plane. For this, start
with the \(3\)-dimensional vector space \(\mathbb{F}^3\). Then define
the projective plane over \(\mathbb{F}\), denoted \(\mathbb{FP}^2\), to
be the space of lines through the origin in \(\mathbb{F}^3\). You can
show the projective plane is the same as the affine plane together with
extra points, which play the role of ``points at infinity''.

In fact, you can generalize this a bit --- you can make sense of the
projective plane over \(\mathbb{F}\) whenever \(\mathbb{F}\) is a
division ring! A division ring is a like a field, but where
multiplication isn't necessarily commutative. The best example is the
quaternions. In \protect\hyperlink{week106}{``Week 106''} I talked about
the real, complex and quaternionic projective planes, their symmetry
groups, and their relation to quantum mechanics. Here's a good book
about this stuff, emphasizing the physics applications:

\begin{enumerate}
\def\labelenumi{\arabic{enumi})}
\tightlist
\item
  V. S. Varadarajan, \emph{Geometry of Quantum Mechanics},
  Springer, Berlin, 1985.
\end{enumerate}

So far, so good. But there's another approach to projective planes
that's even more general. This approach goes back to Euclidean geometry:
it's based on a list of axioms. In this approach, a projective plane
consists of a set of ``points'', a set of ``lines'', and a relation
which tells us whether or not a given point ``lies on'' a given line.
I'm putting quotes around all these words, because in this approach they
are undefined terms. All we get to work with are the following axioms:

\begin{enumerate}
\def\labelenumi{\Alph{enumi})}
\tightlist
\item
  Given two distinct points, there exists a unique line that both points
  lie on.
\item
  Given two distinct lines, there exists a unique point that lies on
  both lines.
\item
  There exist four points, no three of which lie on the same line.
\item
  There exist four lines, no three of which have the same point lying on
  them.
\end{enumerate}

Actually we can leave out either axiom C) or axiom D) --- the rest of
the axioms will imply the one we leave out. It's a nice little exercise
to convince yourself of this. I put in both axioms just to make it
obvious that this definition of projective plane is ``self-dual''. In
other words, if we switch the words ``point'' and ``line'' and switch
who lies on who, the definition stays the same!

Duality is one of the great charms of the theory of projective planes:
whenever you prove any theorem, you get another one free of charge with
the roles of points and lines switched, thanks to duality. There are
lots of different kinds of ``duality'' in mathematics, but this is
probably the grand-daddy of them all.

Now, it's easy to prove that starting from any division ring
\(\mathbb{F}\), we get a projective plane \(\mathbb{FP}^2\) satisfying
the above axioms. The fun part is to try to go the other way! Starting
from a structure satisfying the above axioms, can you cook up a division
ring that it comes from?

Well, starting from a projective plane, you can try to recover a
division ring as follows. Pick a line and throw out one point --- and
call that point ``the point at infinity''. What's left is an ``affine
line'' - let's call it \(L\). Let's try to make \(L\) into a division
ring. To do this, we first need to pick two different points in \(L\),
which we call \(0\) and \(1\). Then we need to cook up rules for adding
and multiplying points on \(L\).

For this, we use some tricks invented by the ancient Greeks!

This should not be surprising. After all, those dudes thought about
arithmetic in very geometrical ways. How can you add points on a line
using the geometry of the plane? Just ask any ancient Greek, and here's
what they'll say:

First pick a line \(L'\) that's parallel to \(L\) --- meaning that \(L\)
and \(L'\) intersect only at the point at infinity. Then pick a line
\(M\) that intersects \(L\) at the point \(0\) and \(L'\) at some point
which we call \(0'\). We get a picture like this: \[
  \begin{tikzpicture}[scale=1.3]
    \node[label={[label distance=-2mm]below left:{$0$}}] at (0,0) {$\bullet$};
    \node[label={[label distance=-2mm]below left:{$0'$}}] at (0.5,2) {$\bullet$};
    \draw[thick] (-0.125,-0.5) to node[fill=white]{$M$} (0.625,2.5);
    \draw[thick] (-0.5,0) to node[fill=white]{$L$} (4,0);
    \draw[thick] (0,2) to node[fill=white]{$L'$} (4.5,2);
  \end{tikzpicture}
\] Then, to add two points \(x\) and \(y\) on \(L\), draw this picture:
\[
  \begin{tikzpicture}[scale=1.3]
    \node[label={[label distance=-2mm]below left:{$0$}}] at (0,0) {$\bullet$};
    \node[label={[label distance=-2mm]below left:{$0'$}}] at (0.5,2) {$\bullet$};
    \draw[thick] (-0.125,-0.5) to node[fill=white]{$M$} (0.375,1.5) to (0.625,2.5);
    \draw[thick] (-0.5,0) to (2.5,0) to node[fill=white]{$L$} (4,0);
    \draw[thick] (0,2) to (3,2) to node[fill=white]{$L'$} (4.5,2);
    \node[label={[label distance=-2mm]below left:{$x$}}] at (0.7,0) {$\bullet$};
    \draw[thick] (0.5,2) to node[fill=white]{$N$} (0.75,-0.5);
    \node at (2.5,2) {$\bullet$};
    \node[label={[label distance=-2mm]below left:{$y$}}] at (2,0) {$\bullet$};
    \node[label={[label distance=-2mm]below left:{$z$}}] at (2.7,0) {$\bullet$};
    \draw[thick] (1.875,-0.5) to node[fill=white]{$M'$} (2.625,2.5);
    \draw[thick] (2.5,2) to node[fill=white]{$N'$} (2.75,-0.5);
  \end{tikzpicture}
\] In other words, draw a line \(M'\) parallel to \(M\) through the
point \(y\), draw a line \(N\) through \(x\) and \(0'\), and draw a line
\(N'\) parallel to \(N\) and going through the point where \(M'\) and
\(L'\) intersect. \(L\) and \(N'\) intersect at the point called
\(z\)\ldots{} and we define this point to be \(x + y\)!

This is obviously the right thing, because the two triangles in the
picture are congruent.

What about multiplication? Well, first draw a line \(L'\) that
intersects our line \(L\) only at the point \(0\). Then draw a line
\(M\) from the point \(1\) to some point \(1'\) that's on \(L'\) but not
on \(L\): \[
  \begin{tikzpicture}
    \draw[thick] (-1,0) to (3,0) node[label=right:{$L$}]{};
    \draw[thick] (0,0) to (1.5,2) node[label=right:{$L'$}]{};
    \draw[thick] (1,-1) node[label=right:{$M$}]{} to (1,1.33);
    \node[label={[label distance=-2mm]below left:{$0$}}] at (0,0){$\bullet$};
    \node[label={[label distance=-2mm]below left:{$1$}}] at (1,0){$\bullet$};
    \node[label={[label distance=-2mm]above left:{$1'$}}] at (1,1.33){$\bullet$};
  \end{tikzpicture}
\] Then, to multiply \(x\) and \(y\), draw this picture: \[
  \begin{tikzpicture}
    \draw[thick] (-1,0) to (6,0) node[label={[label distance=-2mm]right:{$L$}}]{};
    \draw[thick] (0,0) to node[label=above:{$L'$}]{} (3,4);
    \draw[thick] (1,-1) node[label={[label distance=-2mm]right:{$M$}}]{} to (1,1.33) to (1,1.75);
    \draw[thick] (2.5,-1) node[label={[label distance=-2mm]right:{$M'$}}]{} to (2.5,3.33) node{$\bullet$} to (2.5,4);
    \draw[thick] (2,-0.66) node[label={[label distance=-2mm]below:{$N$}}]{} to (1.66,0) node[label={[label distance=-3mm]above right:{$x$}}]{$\bullet$} to (1,1.33) to (0.66,2);
    \draw[thick] (4.66,-0.66) node[label={[label distance=-2mm]below:{$N'$}}]{} to (4.33,0) node[label={[label distance=-3mm]above right:{$z$}}]{$\bullet$} to (2.5,3.33) to (2.17,4);
    \node[label={[label distance=-2mm]below left:{$0$}}] at (0,0){$\bullet$};
    \node[label={[label distance=-2mm]below left:{$1$}}] at (1,0){$\bullet$};
    \node[label={[label distance=-3mm]above right:{$y$}}] at (2.5,0){$\bullet$};
    \node[label={[label distance=-2mm]left:{$1'$}}] at (1,1.33){$\bullet$};
  \end{tikzpicture}
\] In other words, draw a line \(N\) though \(1'\) and \(x\), draw a
line \(M'\) parallel to \(M\) through the point \(y\), and draw a line
\(N'\) parallel to \(N\) through the point where \(L'\) and \(M'\)
intersect. \(L\) and \(N'\) intersect at the point called \(z\)\ldots{}
and we define this point to be \(xy\)!

This is obviously the right thing, because the triangle containing the
points \(1\) and \(x\) is similar to the triangle containing \(y\) and
\(z\).

So, now that we've cleverly figured out how to define addition and
multiplication starting from a projective plane, we can ask: do we get a
division ring?

And the answer is: not necessarily. It's only true if our projective
plane is ``Desarguesian''. This is a special property named after an old
theorem about the real projective plane, proved by Desargues. A
projective plane is Desarguesian if Desargues' theorem holds for this
plane.

But wait --- there's an even more basic question we forgot to ask!
Namely: was our ancient Greek method of defining addition and
multiplication independent of the choices we made? We needed to pick
some points and lines to get things going. If you think about it hard,
these choices boil down to picking four points, no three of which lie on
a line --- exactly what axiom C) guarantees we can do.

Alas, it turns out that in general our recipe for addition and
multiplication really depends on \emph{how} we chose these four points.
But if our projective plane is Desarguesian, it does not!

In fact, if we stick to Desarguesian projective planes, everything works
very smoothly. For any division ring \(\mathbb{F}\) the projective plane
\(\mathbb{FP}^2\) is Desarguesian. Conversely, starting with a
Desarguesian projective plane, we can use the ancient Greek method to
cook up a division ring \(\mathbb{F}\). Best of all, these two
constructions are inverse to each other --- at least up to isomorphism.

At this point you should be pounding your desk and yelling ``Great ---
but what does ``Desarguesian" mean? I want the nitty-gritty details!''

Okay\ldots.

Given a projective plane, define a ``triangle'' to be three points that
don't lie on the same line. Now suppose you have two triangles \(xyz\)
and \(x'y'z'\). The sides of each triangle determine three lines, say
\(LMN\) and \(L'M'N'\). If we're really lucky, the line through \(x\)
and \(x'\), the line through \(y\) and \(y'\), and the line through
\(z\) and \(z'\) will all intersect at the same point. We say that our
projective plane is ``Desarguesian'' if whenever this happens, something
else happens: the intersection of \(L\) and \(L'\), the intersection of
\(M\) and \(M'\), and the intersection of \(N\) and \(N'\) all lie on
the same line.

If you have trouble visualizing what I just said, take a look at this
webpage, which also gives a proof of Desargues' theorem for the real
projective plane:

\begin{enumerate}
\def\labelenumi{\arabic{enumi})}
\setcounter{enumi}{1}
\tightlist
\item
  Roger Mohr and Bill Trigs, ``Desargues' Theorem'', available at
  \href{https://web.archive.org/web/20010210185727/http://spigot.anu.edu.au/people/samer/Research/Doc/ECV_Tut_Proj_Geom/node25.html}{\texttt{https://web.archive.}}
  \href{https://web.archive.org/web/20010210185727/http://spigot.anu.edu.au/people/samer/Research/Doc/ECV_Tut_Proj_Geom/node25.html}{\texttt{org/web/20010210185727/http://spigot.anu.edu.au/people/samer/Research/}}
  \href{https://web.archive.org/web/20010210185727/http://spigot.anu.edu.au/people/samer/Research/Doc/ECV_Tut_Proj_Geom/node25.html}{\texttt{Doc/ECV\_Tut\_Proj\_Geom/node25.html}}
\end{enumerate}

Desargues' theorem is a bit complicated, but one cool thing is that its
converse is its dual. This is easy to see if you stare at it: ``three
lines intersecting at the same point'' is dual to ``three points lying
on the the same line''. Even cooler, Desargues' theorem implies its own
converse! Thus the property of being Desarguesian is self-dual.

Another nice fact about Desarguesian planes concerns collineations. A
``collineation'' is a map from a projective plane to itself that
preserves all lines. Collineations form a group, and this group acts on
the set of all ``quadrangles'' --- a quadrangle being a list of four
points, no three of which lie on a line. Axiom C) says that every
projective plane has at least one quadrangle. It turns out that if a
projective plane is Desarguesian, the group of collineations acts
transitively on the set of quadrangles: given any two quadrangles,
there's a collineation carrying one to the other. This is the reason why
the ancient Greek trick for adding and multiplying doesn't depend on the
choice of quadrangle when our projective plane is Desarguesian!

An even more beautiful fact about Desarguesian planes concerns their
relation to higher dimensions. Just as we defined projective planes
through a list of axioms, we can also define projective spaces of any
dimension \(n = 1,2,3,\ldots\). The simplest example is
\(\mathbb{FP}^n\) --- the space of lines through the origin in
\(\mathbb{F}^{n+1}\), where \(\mathbb{F}\) is some division ring. The
neat part is that when \(n > 2\), this is the \emph{only} example.
Moreover, any projective plane sittting inside one of these
higher-dimensional projective spaces is automatically Desarguesian! So
the non-Desarguesian projective planes are really freaks of dimension 2.

All this is very nice. But there are some obvious further questions,
namely: what's special about projective planes that actually come from
\emph{fields}, and what can we say about non-Desarguesian projective
planes?

The key to the first question is an old theorem proved by the last of
the great Greek geometers, Pappus, in the 3rd century CE. It turns out
that in any projective plane coming from a field, the Pappus Theorem
holds. Conversely, any projective plane satisfying the Pappus Theorem
comes from a unique field. We call such projective planes ``Pappian''.

The Pappus theorem will be too scary if I explain it using only words,
so I'll tell you to look at a picture instead. The fun thing about this
picture is that you can move the points around with your
mouse and see how things change:

\begin{enumerate}
\def\labelenumi{\arabic{enumi})}
\setcounter{enumi}{2}
\tightlist
\item
  Pappus' theorem,
  \url{https://graemewilkin.github.io/Geometry/Pappus.html}
\end{enumerate}

Now, what about the non-Desarguesian projective planes? If we try to get
a division ring from an \emph{arbitrary} projective plane, we fail
miserably. However, we can still define addition and multiplication
using the tricks described above. These operations depend crucially our
choice of a quadrangle. But if we list all the axioms these operations
satisfy, we get the definition of an algebraic gadget called a ``ternary
ring''.

They're called ``ternary rings'' because they're usually described in
terms of a ternary operation that generalizes \(xy + z\). But the
precise definition is too depressing for me give here. It's a classic
example of what James Dolan calls ``centipede mathematics'', where you
take a mathematical concept and see how many legs you can pull off
before it can no longer walk. A ternary ring is like a division ring
that can just barely limp along on its last legs.

I'm not a big fan of centipede mathematics, but there is one really nice
example of a ternary ring that isn't a division ring. Namely, the
octonions! These are almost a division ring, but their multiplication
isn't associative.

I already talked about the octonions in
\protect\hyperlink{week59}{``Week 59''},
\protect\hyperlink{week61}{``Week 61''},
\protect\hyperlink{week104}{``Week 104''} and
\protect\hyperlink{week105}{``Week 105''}. In
\protect\hyperlink{week106}{``Week 106''}, I explained how you can
define \(\mathbb{OP}^2\), the projective plane over the octonions. This
is the best example of a non-Desarguesian projective plane. One reason
it's so great is that that its group of collineations is
\(\mathrm{E}_6\). \(\mathrm{E}_6\) is one of the five ``exceptional
simple Lie groups'' --- mysterious and exciting things that deserve all
the study they can get!

Next I want to talk about the relation between projective geometry and
the \emph{other} exceptional Lie groups, but first let me give you some
references. To start, here's a great book on projective planes and all
the curious centipede mathematics they inspire:

\begin{enumerate}
\def\labelenumi{\arabic{enumi})}
\setcounter{enumi}{3}
\tightlist
\item
  Frederick W. Stevenson, \emph{Projective Planes}, W. H. Freeman and
  Company, San Francisco, 1972.
\end{enumerate}
\noindent
You'll learn all about nearfields, quasifields, Moufang loops, Cartesian
groups, and so on.

Much of the same material is covered in these lectures by Hall, which
are unfortunately a bit hard to find:

\begin{enumerate}
\def\labelenumi{\arabic{enumi})}
\setcounter{enumi}{4}
\tightlist
\item
  Marshall Hall, \emph{Projective Planes and Other Topics}, California
  Institute of Technology, Pasadena, 1954.
\end{enumerate}
\noindent
For a more distilled introduction to the same stuff, try the last
chapter of Hall's book on group theory:

\begin{enumerate}
\def\labelenumi{\arabic{enumi})}
\setcounter{enumi}{5}
\tightlist
\item
  Marshall Hall, \emph{The Theory of Groups}, Macmillan, New York, 1959.
\end{enumerate}
\noindent
If you're only interested in Desarguesian projective planes, try this:

\begin{enumerate}
\def\labelenumi{\arabic{enumi})}
\setcounter{enumi}{6}
\tightlist
\item
  Robin Hartshorne, \emph{Foundations of Projective Geometry}, Benjamin,
  New York, 1967.
\end{enumerate}
\noindent
In particular, this book gives a nice account of the collineation group
in the Desarguesian case. The punchline is simple to state, so I'll tell
you. Suppose \(\mathbb{F}\) is a division ring. Then the collineation
group of \(\mathbb{FP}^2\) is generated by two obvious subgroups:
\(\mathrm{PGL}(3,\mathbb{F})\) and the automorphism group of
\(\mathbb{F}\). The intersection of these two subgroups is the group of
inner automorphisms of \(\mathbb{F}\).

If the above references are too intense, try this leisurely, literate
introduction to the subject first:

\begin{enumerate}
\def\labelenumi{\arabic{enumi})}
\setcounter{enumi}{7}
\tightlist
\item
  Daniel Pedoe, \emph{An Introduction to Projective Geometry},
  Macmillan, New York, 1963.
\end{enumerate}
\noindent
And you're really interested in the \emph{finite} projective planes, you
can try this reference, which assumes very little knowledge of algebra:

\begin{enumerate}
\def\labelenumi{\arabic{enumi})}
\setcounter{enumi}{8}
\tightlist
\item
  A. Adrian Albert and Reuben Sandler, \emph{An Introduction to Finite
  Projective Planes}, Holt, Rinehart and Winston, New York, 1968.
\end{enumerate}
\noindent
For a nice online introduction to projective geometry over the real
numbers and its applications to image analysis, try this:

\begin{enumerate}
\def\labelenumi{\arabic{enumi})}
\setcounter{enumi}{9}
\tightlist
\item
  Roger Mohr and Bill Triggs, ``Projective geometry for image
  analysis'', available at
  \url{https://hal.inria.fr/inria-00548361/}
\end{enumerate}

Finally, for interesting relations between projective geometry and
exceptional Lie groups, try this:

\begin{enumerate}
\def\labelenumi{\arabic{enumi})}
\setcounter{enumi}{10}
\tightlist
\item
  J. M. Landsberg and L. Manivel, ``The projective geometry of
  Freudenthal's magic square'', \emph{Jour.\ Alg. }\textbf{239}, 
  477--512. Also available as
  \href{https://arxiv.org/abs/math.AG/9908039}{\texttt{math.AG/9908039}}.
\end{enumerate}
\noindent
The Freudenthal--Tits magic square is a strange way of describing most of
the exceptional Lie groups in terms of the real numbers, complex
numbers, quaternions and octonions. In the usual way of describing it,
you start with two of these division algebras, say \(\mathbb{F}\) and
\(\mathbb{F}'\). Then let \(J(\mathbb{F})\) be the space of \(3\times3\)
self-adjoint matrices with coefficients in \(\mathbb{F}\). This is a
Jordan algebra with the product \(xy + yx\). As mentioned in
\protect\hyperlink{week106}{``Week 106''}, Jordan algebras have a lot to
do with projective planes. In particular, the nontrivial projections in
\(J(\mathbb{F})\) correspond to the 1- and \(2\)-dimensional subspaces
of \(\mathbb{F}^3\), and thus to the points and lines in the projective
plane \(\mathbb{FP}^2\).

Next, let \(J_0(\mathbb{F})\) be the subspace of \(J(\mathbb{F})\)
consisting of the \emph{traceless} self-adjoint matrices. Also, let
\(\Im(\mathbb{F}')\) be the space of pure imaginary element of
\(\mathbb{F}'\). Finally, let the magic Lie algebra
\(M(\mathbb{F},\mathbb{F}')\) be given by
\[M(K,K') = \operatorname{Der}(J(K)) \oplus J_0(K) \otimes \Im(K') \oplus \operatorname{Der}(K')\]
Here \(\oplus\) stands for direct sum, \(\otimes\) stands for tensor
product, and \(\mathrm{Der}\) stands for the space of derivations of the
algebra in question. It's actually sort of tricky to describe how to
make \(M(\mathbb{F},\mathbb{F}')\) into a Lie algebra, and I'm sort of
tired, so I'll wimp out and tell you to read this stuff:

\begin{enumerate}
\def\labelenumi{\arabic{enumi})}
\setcounter{enumi}{11}
\item
  Hans Freudenthal, ``Lie groups in the foundations of geometry'',
  \emph{Adv. Math.} \textbf{1} (1964) 143.
\item
  Jacques Tits, ``Algebres alternatives, algebres de Jordan et algebres
  de Lie exceptionelles'', \emph{Proc. Colloq. Utrecht},
  vol.~\textbf{135}, 1962.
\item
  R. D. Schafer, \emph{Introduction to Non-associative Algebras},
  Academic Press, New York, 1966.
\end{enumerate}

By the way, the paper by Freudenthal is a really mind-bending mix of Lie
theory and axiomatic projective geometry, definitely worth looking at.
Anyway, if you do things right you get the following square of Lie
algebras \(M(\mathbb{F},\mathbb{F}')\):

\begin{longtable}[]{@{}lcccc@{}}
\toprule
\begin{minipage}[b]{0.09\columnwidth}\raggedright
\strut
\end{minipage} & \begin{minipage}[b]{0.19\columnwidth}\centering
\(\mathbb{F}=\mathbb{R}\)\strut
\end{minipage} & \begin{minipage}[b]{0.19\columnwidth}\centering
\(\mathbb{F}=\mathbb{C}\)\strut
\end{minipage} & \begin{minipage}[b]{0.19\columnwidth}\centering
\(\mathbb{F}=\mathbb{H}\)\strut
\end{minipage} & \begin{minipage}[b]{0.19\columnwidth}\centering
\(\mathbb{F}=\mathbb{O}\)\strut
\end{minipage}\tabularnewline
\midrule
\endhead
\begin{minipage}[t]{0.09\columnwidth}\raggedright
\(\mathbb{F}'=\mathbb{R}\)\strut
\end{minipage} & \begin{minipage}[t]{0.19\columnwidth}\centering
\(\mathrm{A}_1\)\strut
\end{minipage} & \begin{minipage}[t]{0.19\columnwidth}\centering
\(\mathrm{A}_2\)\strut
\end{minipage} & \begin{minipage}[t]{0.19\columnwidth}\centering
\(\mathrm{C}_3\)\strut
\end{minipage} & \begin{minipage}[t]{0.19\columnwidth}\centering
\(\mathrm{F}_4\)\strut
\end{minipage}\tabularnewline
\begin{minipage}[t]{0.09\columnwidth}\raggedright
\(\mathbb{F}'=\mathbb{C}\)\strut
\end{minipage} & \begin{minipage}[t]{0.19\columnwidth}\centering
\(\mathrm{A}_2\)\strut
\end{minipage} & \begin{minipage}[t]{0.19\columnwidth}\centering
\(\mathrm{A}_2\oplus\mathrm{A}_2\)\strut
\end{minipage} & \begin{minipage}[t]{0.19\columnwidth}\centering
\(\mathrm{A}_5\)\strut
\end{minipage} & \begin{minipage}[t]{0.19\columnwidth}\centering
\(\mathrm{E}_6\)\strut
\end{minipage}\tabularnewline
\begin{minipage}[t]{0.09\columnwidth}\raggedright
\(\mathbb{F}'=\mathbb{H}\)\strut
\end{minipage} & \begin{minipage}[t]{0.19\columnwidth}\centering
\(\mathrm{C}_3\)\strut
\end{minipage} & \begin{minipage}[t]{0.19\columnwidth}\centering
\(\mathrm{A}_5\)\strut
\end{minipage} & \begin{minipage}[t]{0.19\columnwidth}\centering
\(\mathrm{B}_6\)\strut
\end{minipage} & \begin{minipage}[t]{0.19\columnwidth}\centering
\(\mathrm{E}_7\)\strut
\end{minipage}\tabularnewline
\begin{minipage}[t]{0.09\columnwidth}\raggedright
\(\mathbb{F}'=\mathbb{O}\)\strut
\end{minipage} & \begin{minipage}[t]{0.19\columnwidth}\centering
\(\mathrm{F}_4\)\strut
\end{minipage} & \begin{minipage}[t]{0.19\columnwidth}\centering
\(\mathrm{E}_6\)\strut
\end{minipage} & \begin{minipage}[t]{0.19\columnwidth}\centering
\(\mathrm{E}_7\)\strut
\end{minipage} & \begin{minipage}[t]{0.19\columnwidth}\centering
\(\mathrm{E}_8\)\strut
\end{minipage}\tabularnewline
\bottomrule
\end{longtable}

Here \(\mathbb{R}\), \(\mathbb{C}\), \(\mathbb{H}\) and \(\mathbb{O}\)
stand for the reals, complexes, quaternions and octonions. If you don't
know what all the Lie algebras in the square are, check out
\protect\hyperlink{week64}{``Week 64''}. (I should admit that the above
square is not very precise, because I don't say which real forms of the
Lie algebras in question are showing up.)

The first fun thing about this square is that \(\mathrm{F}_4\),
\(\mathrm{E}_6\), \(\mathrm{E}_7\) and \(\mathrm{E}_8\) are four of the
five exceptional simple Lie algebras --- and the fifth one,
\(\mathrm{G}_2\), is just \(\operatorname{Der}(\mathbb{O})\). So all the
exceptional Lie algebras are related to the octonions! And the second
fun thing about this square is that it's symmetrical along the diagonal,
even though this is not at all obvious from the definition. This is what
makes the square truly ``magic''.

I don't really understand the magic square, but it's on my to-do list.
That's one reason I'm glad there's a new paper out that describes the
magic square in a way that makes its symmetry manifest:

\begin{enumerate}
\def\labelenumi{\arabic{enumi})}
\setcounter{enumi}{14}
\tightlist
\item
  C. H. Barton and A. Sudbery, ``Magic squares of Lie algebras'',
  available as \hfill \break
  \href{https://arxiv.org/abs/math.RA/0001083}{\texttt{math.RA/0001083}}.
\end{enumerate}

It also generalizes the magic square in a number of directions. But what
I really want is for the connection between projective planes, division
algebras, exceptional Lie groups and the magic square to becomes truly
\emph{obvious} to me. I'm a long way from that point.

\begin{center}\rule{0.5\linewidth}{0.5pt}\end{center}

Here's an interesting email from David Broadhurst about the failure of
the Pappus theorem in quaternionic projective space:

\begin{quote}
Date: Fri, 3 Mar 2000 20:50:03 GMT \hfill \break
From: David Broadhurst \hfill \break
Subject: Paul Dirac and projective geometry

John:

Shortly before his death I spent a charming afternoon with Paul Dirac.
Contrary to his reputation, he was most forthcoming.

Among many things, I recall this: Dirac explained that while trained
as an engineer and known as a physicist, his aesthetics were mathematical.
He said (as I can best recall, nearly 20 years on): At a young age,
I fell in love with projective geometry.  I always wanted to use to use 
it in physics, but never found a place for it.

Then someone told him that the difference between complex and quaternionic
QM had been characterized as the failure of theorem in classical projective 
geometry.

Dirac's face beamed a lovely smile: ``Ah," he said, ``it was just such a thing
that I hoped to do".

I was reminded of this when backtracking to your 
\protect\hyperlink{week106}{``Week 106''}, today.

Best,

David
\end{quote}

\begin{center}\rule{0.5\linewidth}{0.5pt}\end{center}

\begin{quote}
\emph{The reader should not attempt to form a mental picture of a closed
straight line.}

--- Frank Ayres, Jr., \emph{Projective Geometry}
\end{quote}

\hypertarget{week146}{%
\section{March 11, 2000}\label{week146}}

\begin{quote}
Paper in white the floor of the room, and rule it off in one- foot
squares. Down on one's hands and knees, write in the first square a set
of equations conceived as able to govern the physics of the universe.
Think more overnight. Next day put a better set of equations into square
two. Invite one's most respected colleagues to contribue to other
squares. At the end of these labors, one has worked oneself out into the
doorway. Stand up, look back on all those equations, some perhaps more
hopeful than others, raise one's finger commandingly, and give the order
``Fly!'' Not one of those equations will put on wings, take off, or fly.
Yet the universe ``flies''.

Some principle uniquely right and compelling must, when one knows it, be
also so compelling that it is clear the universe is built, and must be
built, in such and such a way, and that it could not be otherwise. But
how can one disover that principle?
\end{quote}

John Wheeler was undoubtedly the author of these words, which appear
near the end of Misner, Thorne and Wheeler's textbook \emph{Gravitation},
published in 1972. Since then, more people than ever before in the
history of the world have tried their best to find this uniquely
compelling principle. A lot of interesting ideas, but no success yet.

But what if Wheeler was wrong? What if there is \emph{not} a uniquely
compelling principle or set of equations that governs our universe? For
example, what if \emph{all} equations govern universes? In other words,
what if all mathematical structures have just as much ``physical
existence'' (whatever that means!) as those describing our universe do?
Many of them will not contain patterns we could call awareness or
intelligence, but some will, and these would be ``seen from within'' as
``universes'' by their inhabitants. In this scenario, there's nothing
special about \emph{our} universe except that we happen to be in this
one.

In other words: what if there is ultimately no difference between
mathematical possibility and physical existence? This may seem crazy,
but personally I believe that most alternatives, when carefully
pondered, turn out to be even \emph{more} crazy.

Of course, it's fun to think about such ideas and difficult to get
anywhere with them. But tonight when I was nosing around the web, I
found out that someone has already developed and published this idea:

\begin{enumerate}
\def\labelenumi{\arabic{enumi})}
\item
  Max Tegmark, ``Is the `theory of everything' merely the ultimate
  ensemble theory?", \emph{Ann. Phys.} \textbf{270} (1998), 1--51.
  Also available as
  \href{https://arxiv.org/abs/gr-qc/9704009}{\texttt{gr-qc/9704009}}.

  Max Tegmark, ``Welcome to my crazy universe'', available at
  \href{https://space.mit.edu/home/tegmark/toe.html}{\texttt{https://space.mit.}}
  \href{https://space.mit.edu/home/tegmark/toe.html}{\texttt{edu/home/tegmark/toe.html}}

\item
  Marcus Chown, ``Anything goes'', \emph{New Scientist} \textbf{158}
  (1998) 26--30.  Also available at
  \url{https://space.mit.edu/home/tegmark/toe_press.html}
\end{enumerate}
\noindent
As far as I can tell, most of the time Max Tegmark is a perfectly
respectable physicist at the University of Pennsylvania; he works on the
cosmic microwave background radiation, the large-scale structure of the
universe (superclusters and the like), and type 1A supernovae. But he
has written a fascinating paper on the above hypothesis, which I urge
you to read. It's less far-out than you may think.

Okay, now on to quantum gravity. Jan Ambj\o rn and Renate Loll have teamed
up to work on discrete models of spacetime geometry, with an emphasis on
the Lorentzian geometry of triangulated manifolds. Much more has been
done over on the Riemannian side of things, so it's high time to focus
more energy on the physically realistic Lorentzian case. Of course, if
the metric is fixed you can often use a trick called ``Wick rotation''
to turn results about quantum field theory on Riemannian spacetime into
results for Lorentzian spacetime. But it's never been clear that this
works when the geometry of spacetime is a variable --- and quantized,
for that matter. So we need both more work on Wick rotation in this
context and also work that goes straight for the jugular: the Lorentzian
context.

Here are some of their papers:

\begin{enumerate}
\def\labelenumi{\arabic{enumi})}
\setcounter{enumi}{3}
\item
  J. Ambj\o rn, J. Correia, C. Kristjansen, and R. Loll, ``On the relation
  between Euclidean and Lorentzian 2d quantum gravity'', \emph{Phys.\ Lett.\ B}
  \textbf{475} (2000), 24--32.  Also available as
  \href{https://arxiv.org/abs/hep-th/9912267}{\texttt{hep-th/9912267}}.

  J. Ambj\o rn, J. Jurkiewicz and R. Loll, ``Lorentzian and Euclidean
  quantum gravity --- analytical and numerical results'', in \emph{M-Theory and Quantum Geometry}, edited by L.\ Thorlacius and T.\ Jonsson, NATO Science Series \textbf{556}, Springer, Berlin, 381--450.  Also  available as
  \href{https://arxiv.org/abs/hep-th/0001124}{\texttt{hep-th/0001124}}.

  J. Ambj\o rn, J. Jurkiewicz and R. Loll, ``A non-perturbative Lorentzian
  path integral for gravity'', \emph{Phys.\ Rev.\ Lett.} \textbf{85} (2000) 
   924--927. Also available as
  \href{https://arxiv.org/abs/hep-th/0002050}{\texttt{hep-th/0002050}}.
\end{enumerate}
\noindent
The last paper is especially interesting to me, since it tackles the
problem of defining a path integral for 3+1-dimensional Lorentzian
quantum gravity. They describe a path integral where you first slice
spacetime like a salami using surfaces of constant time, and then pack
each slice with simplices having edges with specified lengths --- the
edges being spacelike within each surface, and timelike when they go
from one surface to the next. They allow the number of simplices in each
slice to be variable. Actually they focus on the 2+1-dimensional case,
but they say the 3+1-dimensional case works similarly, and I actually
trust them enough to believe them about this --- especially since
nothing they do relies on the fact that 2+1-dimensional gravity lacks
local degrees of freedom. They can Wick-rotate this picture and get a
time evolution operator that's self-adjoint and positive, just like
you'd expect of an operator of the form \(\exp(-tH)\).

Speaking of Wick rotations in quantum gravity, here's another paper to
think about:

\begin{enumerate}
\def\labelenumi{\arabic{enumi})}
\setcounter{enumi}{4}
\tightlist
\item
  Abhay Ashtekar, Donald Marolf, Jose Mourao and Thomas Thiemann,
  ``Osterwalder--Schrader reconstruction and diffeomorphism invariance'',
  \emph{Class.\ Quant.\ Grav.} \textbf{17} (2000), 4919--4940.  Also available as
  \href{https://arxiv.org/abs/quant-ph/9904094}{\texttt{quant-ph/}}
 \href{https://arxiv.org/abs/quant-ph/9904094}{\texttt{9904094}}.
\end{enumerate}
\noindent
The Osterwalder--Shrader theorem is the result people use when they want
to \emph{rigorously} justify Wick rotations. Here these authors
generalize it so that it applies to a large class of background-free
field theories --- perhaps even quantum gravity! It turns out not to be
hard, once you go about it properly. Quite a surprise.

I've been working with Ashtekar and Krasnov for a couple of years now on
computing the entropy of black holes using loop quantum gravity. I
talked about this in \protect\hyperlink{week112}{``Week 112''}, right
after we came out with a short paper sketching the calculation. Now
we're almost done with the detailed paper. In the meantime, Ashtekar has
written a couple of pedagogical accounts explaining the basic idea. I
mentioned one he wrote with Krasnov in
\protect\hyperlink{week120}{``Week 120''}, and here's another:

\begin{enumerate}
\def\labelenumi{\arabic{enumi})}
\setcounter{enumi}{5}
\tightlist
\item
  Abhay Ashtekar, ``Interface of general relativity, quantum physics and
  statistical mechanics: some recent developments'', \emph{Annalen Phys.}      
  \textbf{9} (2000), 178--198.  Also available as
  \href{https://arxiv.org/abs/gr-qc/9910101}{\texttt{gr-qc/9910101}}.
\end{enumerate}

Let me just quote the abstract --- I can't bear to talk about this any
more until the actual paper is finished:

\begin{quote}
The arena normally used in black holes thermodynamics was recently
generalized to incorporate a broad class of physically interesting
situations. The key idea is to replace the notion of stationary event
horizons by that of `isolated horizons.' Unlike event horizons, isolated
horizons can be located in a space-time quasi-locally. Furthermore, they
need not be Killing horizons. In particular, a space-time representing a
black hole which is itself in equilibrium, but whose exterior contains
radiation, admits an isolated horizon. In spite of this generality, the
zeroth and first laws of black hole mechanics extend to isolated
horizons. Furthermore, by carrying out a systematic, non-perturbative
quantization, one can explore the quantum geometry of isolated horizons
and account for their entropy from statistical mechanical
considerations. After a general introduction to black hole
thermodynamics as a whole, these recent developments are briefly
summarized.
\end{quote}

There have also been a number of papers working out the details of the
classical notion of ``isolated horizon'' --- I've mentioned some
already, but let me just list them all here:

\begin{enumerate}
\def\labelenumi{\arabic{enumi})}
\setcounter{enumi}{6}
\item
  Abhay Ashtekar, Alejandro Corichi, and Kirill Krasnov, ``Isolated
  horizons: the classical phase space'', \emph{Adv. Theor. Math. 
  Phys.}\textbf{3} (1999), 419--478.  Also available as
  \href{https://arxiv.org/abs/gr-qc/9905089}{\texttt{gr-qc/9905089}}.

  Abhay Ashtekar, Christopher Beetle, and Stephen Fairhurst, ``Mechanics
  of isolated horizons'', \emph{Class. Quant. Grav.} \textbf{17} (2000)
  253--298.  Also available as
  \href{https://arxiv.org/abs/gr-qc/9907068}{\texttt{gr-qc/9907068}}.

  Abhay Ashtekar and Alejandro Corichi, ``Laws governing isolated
  horizons: inclusion of dilaton couplings'', \emph{Class. Quant. Grav.}   
  \textbf{17} (2000), 1317--1332.   Also available as
  \href{https://arxiv.org/abs/gr-qc/9910068}{\texttt{gr-qc/9910068}}.

  Jerzy Lewandowski, ``Space-times admitting isolated horizons'',
  \emph{Class. Quant. Grav.} \textbf{17} (2000), L53--L59.
  Also available as
  \href{https://arxiv.org/abs/gr-qc/9907058}{\texttt{gr-qc/9907058}}.
\end{enumerate}

Lewandowski's paper is important because it gets serious about studying
rotating isolated horizons --- this makes me feel a lot more optimistic
that we'll eventually be able to extend the entropy calculation to
rotating black holes (so far it's just done for the nonrotating case).

Okay, now let me turn my attention to spin foams. Last month,
Reisenberger and Rovelli came out with a couple of papers that push
forward the general picture of spin foams as Feynman diagrams,
generalizing the old work of Boulatov and Ooguri, and the newer work of
De Pietri et al.~Again, I'll just quote the abstracts\ldots.

\begin{enumerate}
\def\labelenumi{\arabic{enumi})}
\setcounter{enumi}{7}
\tightlist
\item
  Michael Reisenberger and Carlo Rovelli, ``Spin foams as Feynman
  diagrams'', in \emph{A Relativistic Spacetime Odyssey: Experiments and  
  Theoretical Viewpoints on General Relativity and Quantum Gravity}, edited
  by Ignazio Ciufolini, Daniele Dominici and Luca Lusanna, World Scientific,
 Singapore, pp.\ 431--448.  Also available as
  \href{https://arxiv.org/abs/gr-qc/0002083}{\texttt{gr-qc/0002083}}.
\end{enumerate}

\begin{quote}
It has been recently shown that a certain non-topological spin foam
model can be obtained from the Feynman expansion of a field theory over
a group. The field theory defines a natural ``sum over triangulations'',
which removes the cutoff on the number of degrees of freedom and
restores full covariance. The resulting formulation is completely
background independent: spacetime emerges as a Feynman diagram, as it
did in the old two-dimensional matrix models. We show here that any spin
foam model can be obtained from a field theory in this manner. We give
the explicit form of the field theory action for an arbitrary spin foam
model. In this way, any model can be naturally extended to a sum over
triangulations. More precisely, it is extended to a sum over
2-complexes.
\end{quote}

\begin{enumerate}
\def\labelenumi{\arabic{enumi})}
\setcounter{enumi}{8}
\tightlist
\item
  Michael Reisenberger and Carlo Rovelli, ``Spacetime as a Feynman
  diagram: the connection formulation'', available as
  \href{https://arxiv.org/abs/gr-qc/0002095}{\texttt{gr-qc/0002095}}.
\end{enumerate}

\begin{quote}
Spin foam models are the path integral counterparts to loop quantized
canonical theories. In the last few years several spin foam models of
gravity have been proposed, most of which live on finite simplicial
lattice spacetime. The lattice truncates the presumably infinite set of
gravitational degrees of freedom down to a finite set. Models that can
accomodate an infinite set of degrees of freedom and that are
independent of any background simplicial structure, or indeed any a
priori spacetime topology, can be obtained from the lattice models by
summing them over all lattice spacetimes. Here we show that this sum can
be realized as the sum over Feynman diagrams of a quantum field theory
living on a suitable group manifold, with each Feynman diagram defining
a particular lattice spacetime. We give an explicit formula for the
action of the field theory corresponding to any given spin foam model in
a wide class which includes several gravity models. Such a field theory
was recently found for a particular gravity model {[}De Pietri et al,
\href{https://arxiv.org/abs/hep-th/9907154}{\texttt{hep-th/9907154}}{]}.
Our work generalizes this result as well as Boulatov's and Ooguri's
models of three and four dimensional topological field theories, and
ultimately the old matrix models of two dimensional systems with
dynamical topology. A first version of our result has appeared in a
companion paper
{[}\href{https://arxiv.org/abs/gr-qc/0002083}{\texttt{gr-qc/0002083}}{]}:
here we present a new and more detailed derivation based on the
connection formulation of the spin foam models.
\end{quote}

I'm completely biased, but I think this is the way to go in quantum
gravity\ldots{} we need to think more about the Lorentzian side of
things, like Barrett and Crane have been doing, but these spin foam
models are so darn simple and elegant I can't help but think there's
something right about them --- especially when you see the sum over
triangulations pop out automatically from the Feynman diagram expansion
of the relevant path integral.

There's also been some good work on the relation between canonical
quantum gravity and Vassiliev invariants. The idea is to use this class
of knot invariants as a basis for a Hilbert space of diffeomorphism-
invariant states --- a tempting alternative to the Hilbert space having
spin networks as a basis. Maybe everything will start making sense when
we see how these two choices fit together. But anyway, these papers
tackle the crucial issue of the Hamiltonian constraint using this
Vassiliev approach, and get results startlingly similar to those
obtained by Thiemann using the spin network approach:

\begin{enumerate}
\def\labelenumi{\arabic{enumi})}
\setcounter{enumi}{9}
\item
  Cayetano Di Bartolo, Rodolfo Gambini, Jorge Griego, and Jorge Pullin,
  ``Consistent canonical quantization of general relativity in the space
  of Vassiliev invariants'', available as
  \href{https://arxiv.org/abs/gr-qc/9909063}{\texttt{gr-qc/9909063}}.

  ``Canonical quantum gravity in the Vassiliev invariants arena: I.
  Kinematical structure'', available as
  \href{https://arxiv.org/abs/gr-qc/9911009}{\texttt{gr-qc/9911009}}.

  ``Canonical quantum gravity in the Vassiliev invariants arena: II.
  Constraints, habitats and consistency of the constraint algebra'',
   available as
  \href{https://arxiv.org/abs/gr-qc/9911010}{\texttt{gr-qc/9911010}}.
\end{enumerate}

Finally, Martin Bojowald has written a couple of papers applying the
loop approach to quantum cosmology. The idea is to apply loop
quantization to a ``minisuperspace'' --- a phase space describing only
those solutions of general relativity that have a certain large symmetry
group.

\begin{enumerate}
\def\labelenumi{\arabic{enumi})}
\setcounter{enumi}{10}
\item
  Martin Bojowald, ``Loop Quantum Cosmology I: Kinematics'', 
  available as
  \href{https://arxiv.org/abs/gr-qc/9910103}{\texttt{gr-qc/9910103}}.

  Martin Bojowald, ``Loop Quantum Cosmology II: Volume Operators'',
  \href{https://arxiv.org/abs/gr-qc/9910104}{\texttt{gr-qc/9910104}}.
\end{enumerate}

\hypertarget{week147}{%
\section{May 20, 2000}\label{week147}}

Various books are coming out to commemorate the millennium\ldots.
describing the highlights of the math we've done so far, and laying out
grand dreams for the future. The American Mathematical Society has come
out with one:

\begin{enumerate}
\def\labelenumi{\arabic{enumi})}
\tightlist
\item
  \emph{Mathematics: Frontiers and Perspectives}, edited by Vladimir
  Arnold, Michael Atiyah, Peter Lax and Barry Mazur, AMS, Providence,
  Rhode Island, 2000.
\end{enumerate}
\noindent
This contains 30 articles by bigshots like Chern, Connes, Donaldson,
Jones, Lions, Manin, Mumford, Penrose, Smale, Vafa, Wiles and Witten. I
haven't actually read it yet, but I want to get ahold of it.

Springer Verlag is coming out with one, too:

\begin{enumerate}
\def\labelenumi{\arabic{enumi})}
\setcounter{enumi}{1}
\tightlist
\item
  \emph{Mathematics Unlimited: 2001 and Beyond}, edited by Bjorn
  Engquist and Wilfried Schmid, Springer Verlag, New York, 2000.
\end{enumerate}
\noindent
It should appear in the fall.

I don't know what the physicists are doing along these lines. The
American Physical Society has a nice timeline of 20th century physics on
their website:

\begin{enumerate}
\def\labelenumi{\arabic{enumi})}
\setcounter{enumi}{2}
\tightlist
\item
  The American Physical Society: ``A Century of Physics'', available at
  \href{https://web.archive.org/web/19990508143827/http://timeline.aps.org/APS/home_HighRes.html}{\texttt{https://web.archive.org/web/19990508143827/http://timeline.aps.org/APS/home\_HighRes.html}}
\end{enumerate}
\noindent
But I don't see anything about books.

One reason I haven't been doing many This Week's Finds lately is that
I've been buying and then moving into a new house. Another is that James
Dolan and I have been busily writing our own millennial pontifications,
which will appear in the Springer Verlag book:

\begin{enumerate}
\def\labelenumi{\arabic{enumi})}
\setcounter{enumi}{3}
\tightlist
\item
  John Baez and James Dolan, \emph{From finite sets to Feynman
  diagrams}.  Available as
  \href{http://arXiv.org/abs/math.QA/0004133}{\texttt{math.QA/0004133}}
\end{enumerate}
\noindent
So let me talk about this stuff a bit\ldots.

As usual, the underlying theme of this paper is categorification. I've
talked about this a lot already --- e.g.~in
\protect\hyperlink{week121}{``Week 121''} --- so I'll assume you
remember that when we categorify, we use the following analogy 
to take interesting \emph{equations} and see them as shorthand for even
more interesting \emph{isomorphisms}:

\vfill
\newpage

\begin{longtable}{cc}
\toprule
\textbf{set theory} & \textbf{category theory} \tabularnewline
\midrule
\endhead
elements & objects \tabularnewline
equations between elements & isomorphisms between objects \tabularnewline
\midrule 
sets & categories \tabularnewline
functions between sets & functors between categories \tabularnewline
equations between functions & natural isomorphisms between functors \tabularnewline
\end{longtable}

To take a simple example, consider the laws of basic arithmetic, like
\(a+b = b+a\) or \(a(b+c) = ab+ac\). We usually think of these as
\emph{equations} between \emph{elements} of the \emph{set} of natural
numbers. But really they arise from \emph{isomorphisms} between
\emph{objects} of the \emph{category} of finite sets.

For example, if we have finite sets \(a\) and \(b\), and we use \(a+b\)
to denote their disjoint union, then there is a natural isomorphism
between \(a+b\) and \(b+a\). Moreover, this isomorphism is even sort of
interesting! For example, suppose we use \(1\) to denote a set
consisting of one dot, and \(2\) to denote a set of two dots. Then the
natural isomorphism between \(1+2\) and \(2+1\) can be visualized as the
process of passing one dot past two, like this: \[
  \begin{tikzpicture}[knot gap=7]
    \draw[knot] (1,0) to (0,-3);
    \draw[knot] (1.5,0) to (0.5,-3);
    \draw[knot] (0,0) to (1.5,-3);
    \node at (0,0) {$\bullet$};
    \node at (1,0) {$\bullet$};
    \node at (1.5,0) {$\bullet$};
    \node at (0,-3) {$\bullet$};
    \node at (0.5,-3) {$\bullet$};
    \node at (1.5,-3) {$\bullet$};
  \end{tikzpicture}
\] This may seem like an excessively detailed ``picture proof'' that 1+2
indeed equals 2+1, perhaps suitable for not-too-bright kindergarteners.
But in fact it's just a hop, skip and a jump from here to heavy-duty
stuff like the homotopy groups of spheres. I sketched how this works in
\protect\hyperlink{week102}{``Week 102''} so I won't do so again here.
The point is, after we categorify, elementary math turns out to be
pretty powerful!

Now, let me make this idea of ``categorifying the natural numbers'' a
bit more precise. Let \(\mathsf{FinSet}\) stand for the category whose
objects are finite sets and whose morphisms are functions between these.
If we ``decategorify'' this category by forming the set of isomorphism
classes of objects, we get \(\mathbb{N}\), the natural numbers.

All the basic arithmetic operations on \(\mathbb{N}\) come from
operations on \(\mathsf{FinSet}\). I've already noted how addition comes
from disjoint union. Disjoint union is a special case of what category
theorists call the ``coproduct'', which makes sense for a bunch of
categories --- see \protect\hyperlink{week99}{``Week 99''} for the
general definition. Similarly, multiplication comes from the Cartesian
product of finite sets, which is a special case of what category
theorists call the ``product''. To get the definition of a product, you
just take the definition of a coproduct and turn all the arrows around.
There are also nice category-theoretic interpretations of the numbers
\(0\) and \(1\), and all the basic laws of addition and multiplication.
Exponentiation too!

Combinatorists have lots of fun thinking about how to take equations in
\(\mathbb{N}\) and prove them using explicit isomorphisms in
\(\mathsf{FinSet}\) --- they call such a proof a ``bijective proof''. To
read more about this, try:

\begin{enumerate}
\def\labelenumi{\arabic{enumi})}
\setcounter{enumi}{4}
\item
  James Propp and David Feldman, ``Producing new bijections from old'',
  \emph{Adv. Math.} \textbf{113} (1995), 1--44. Also available at
  \href{http://faculty.uml.edu/jpropp/cancel.pdf}{\texttt{http://faculty.uml.edu/jpropp/cancel.}} \href{http://faculty.uml.edu/jpropp/cancel.pdf}{\texttt{pdf}}
\item
  John Conway and Peter Doyle, ``Division by three''. Available at
  \href{https://arxiv.org/abs/math/0605779}{math/0605779}.
\end{enumerate}

The latter article studies this question: if I give you an isomorphism
between \(3x\) and \(3y\), can you construct a isomorphism between \(x\)
and \(y\)? Here of course \(x\) and \(y\) are finite sets, \(3\) is any
3-element set, and multiplication means Cartesian product. Of course you
can prove an isomorphism exists, but can you \emph{construct} one in a
natural way --- i.e., without making any random choices? The history of
this puzzle turns out to be very interesting. But I don't want to give
away the answer! See if you can do it or not.

Anyway, having categorified the natural numbers, we might be inclined to
go on and categorify the integers. Can we do it? In other words: can we
find something like the category of finite sets that includes ``sets
with a negative number of elements''? There turns out be an interesting
literature on this subject:

\begin{enumerate}
\def\labelenumi{\arabic{enumi})}
\setcounter{enumi}{6}
\item
  Daniel Loeb, ``Sets with a negative number of elements'', \emph{Adv.
  Math.\ } \textbf{91} (1992), 64--74.
\item
  S. Schanuel, Negative sets have Euler characteristic and
  dimension, in \emph{Category Theory}, 
   Lecture Notes in Mathematics \textbf{1488}, Springer,
  Berlin, 1991, pp.~379--385.
\item
  James Propp, ``Exponentiation and Euler measure'', available as
  \href{http://www.arXiv.org/abs/math.CO/0204009}{\texttt{math.CO/0204009}}.
\item
  Andre Joyal, ``Regle des signes en algebre combinatoire'',
  \emph{Comptes Rendus Mathematiques de l'Academie des Sciences, La
  Societe Royale du Canada}, \textbf{VII} (1985), 285--290.
\end{enumerate}
\noindent
See also \protect\hyperlink{week102}{``Week 102''} for more\ldots.

But I don't want to talk about negative sets right now! Instead, I want
to talk about \emph{fractional} sets. It may seem odd to tackle division
before subtraction, but historically, the negative numbers were invented
quite a bit \emph{after} the nonnegative rational numbers. Apparently
half an apple is easier to understand than a negative apple! This
suggests that perhaps `sets with fractional cardinality' are simpler
than `sets with negative cardinality'.

The key is to think carefully about the meaning of division. The usual
way to get half an apple is to chop one into ``two equal parts''. Of
course, the parts are actually \emph{not equal} --- if they were,
there would be only one part! They are merely \emph{isomorphic}. This
suggests that categorification will be handy.

Indeed, what we really have is a \(\mathbb{Z}/2\) symmetry group acting
on the apple which interchanges the two isomorphic parts. In general, if
a group \(G\) acts on a set \(S\), we can ``divide'' the set by the
group by taking the quotient \(S/G\), whose points are the orbits of the
action. If \(S\) and \(G\) are finite and \(G\) acts freely on \(S\),
this construction really corresponds to division, since the cardinality
\(|S/G|\) is equal to \(|S|/|G|\). However, it is crucial that the
action be free.

For example, why is \(6/2 = 3\)? We can take a set \(S\) consisting of
six dots in a row:
\[\bullet\quad\bullet\quad\bullet\quad\bullet\quad\bullet\quad\bullet\]
let \(G = \mathbb{Z}/2\) act freely by reflections, and identify all the
elements in each orbit to obtain 3-element set \(S/G\). Pictorially,
this amounts to folding the set \(S\) in half, so it is not surprising
that \(|S/G| = |S|/|G|\) in this case. Unfortunately, if we try a
similar trick starting with a 5-element set:
\[\bullet\quad\bullet\quad\bullet\quad\bullet\quad\bullet\] it fails
miserably! We don't obtain a set with 2.5 elements, because the group
action is not free: the point in the middle gets mapped to itself. So to
define ``sets with fractional cardinality'', we need a way to count the
point in the middle as ``half a point''.

To do this, we should first find a better way to define the quotient of
\(S\) by \(G\) when the action fails to be free. Following the policy of
replacing equations by isomorphisms, let us define the ``weak quotient''
\(S//G\) to be the category with elements of \(S\) as its objects, with
a morphism \(g\colon s\to s'\) whenever \(g(s) = s'\), and with
composition of morphisms defined in the obvious way.

Next, let us figure out a good way to define the ``cardinality'' of a
category. Pondering the examples above leads us to the following recipe:
for each isomorphism class of objects we pick a representative \(x\) and
compute the reciprocal of the number of automorphisms of this object;
then we sum over isomorphism classes.

It is easy to see that with this definition, the point in the middle of
the previous picture gets counted as `half a point' because it has two
automorphisms, so we get a category with cardinality 2.5. In general,
\[|S//G| = |S|/|G|\] whenever \(G\) is a finite group acting on a finite
set \(S\). This formula is a simplified version of `Burnside's lemma',
so-called because it is due to Cauchy and Frobenius. Burnside's lemma
gives the cardinality of the ordinary quotient. But the weak quotient is
nicer, so Burnside's lemma simplifies when we use weak quotients.

Now, the formula for the cardinality of a category makes sense even for
some categories that have infinitely many objects --- all we need is for
the sum to make sense. So let's try to compute the cardinality of the
category of finite sets! Since any \(n\)-element set has \(n!\)
automorphisms (i.e.~permutations), we get following marvelous formula:
\[|\mathsf{FinSet}| = e\] This turns out to explain lots of things about
the number \(e\).

Now, a category all of whose morphisms are isomorphisms is called a
``groupoid''. Any category \(\mathcal{C}\) has an underlying groupoid
\(\mathcal{C}_0\) with the same objects but only the isomorphisms as
morphisms. The cardinality of a category \(\mathcal{C}\) always equals
that of its underlying groupoid \(\mathcal{C}_0\). This suggests that
this notion should really be called "groupoid cardinality. If you're a
fan of \(n\)-categories, this suggests that we should generalize the
concept of cardinality to n-groupoids, or even \(\omega\)-groupoids. And
luckily, we don't need to understand \(\omega\)-groupoids very well to
try our hand at this! \(\omega\)-groupoids are supposed to be an
algebraic way of thinking about topological spaces up to homotopy. Thus
we just need to invent a concept of the `cardinality' of a topological
space which has nice formal properties and which agrees with the
groupoid cardinality in the case of homotopy 1-types. In fact, this is
not hard to do. We just need to use the homotopy groups \(\pi_k(X)\) of
the space \(X\).

So: let's define the ``homotopy cardinality'' of a topological space
\(X\) to be the alternating product
\[|X| = |\pi_1(X)|^{-1} \cdot |\pi_2(X)| \cdot |\pi_3(X)|^{-1} \cdot \,\, \ldots\]
when \(X\) is connected and the product converges; if \(X\) is not
connected, let's define its homotopy cardinality to be the sum of the
homotopy cardinalities of its connected components, when the sum
converges. We call spaces with well-defined homotopy cardinality
``tame''. The disjoint union or Cartesian product of tame spaces is
again tame, and we have \[
  \begin{aligned}
    |X + Y| &= |X| + |Y| ,
  \\|X \times Y| &= |X| \times |Y|
  \end{aligned}
\] just as you would hope.

Even better, homotopy cardinality gets along well with fibrations, which
we can think of as `twisted products' of spaces. Namely, if
\[F\to X\to B\] is a fibration and the base space B is connected, we
have \[|X| = |F| \times |B|\] whenever two of the three spaces in
question are tame (which implies the tameness of the third).

As a fun application of this fact, recall that any topological group
\(G\) has a ``classifying space'' \(BG\), meaning a space with a
principal \(G\)-bundle over it \[G\to EG\to BG\] whose total space \(EG\) is
contractible. I described how to construct the classifying space in
\protect\hyperlink{week117}{``Week 117''}, at least in the case of a
discrete group \(G\), but I didn't say much about why it's so great. The
main reason it's great is that \emph{any} \(G\)-bundle over \emph{any}
space is a pullback of the bundle \(EG\) over \(BG\). But right now,
what I want to note is that since \(EG\) is contractible it is tame, and
\(|EG| = 1\). Thus \(G\) is tame if and only if \(BG\) is, and
\[|BG| = 1 / |G|\] so we can think of \(BG\) as the ``reciprocal'' of
\(G\)!

This idea is already lurking behind the usual approach to ``equivariant
cohomology''. Suppose \(X\) is a space on which the topological group
\(G\) acts. When the action of \(G\) on \(X\) is free, it is fun to
calculate cohomology groups (and other invariants) of the quotient space
\(X/G\). When the action is not free, this quotient can be very
pathological, so people usually replace it by the ``homotopy quotient''
\(X//G\), which is defined as \((EG \times X)/G\). This is like the
ordinary quotient but with equations replaced by homotopies. And there
is a fibration \[X\to X//G\to BG,\] so when \(X\) and \(G\) are tame we
have \[|X//G| = |X| \times |BG| = |X|/|G|\] just as you would hope!

Now in the paper, Jim and I go on to talk about how all these ideas can
be put to use to give a nice explanation of the combinatorics of Feynman
diagrams. But I don't want to explain all that stuff here --- then you
wouldn't need to read the paper! Instead, I just want to point out
something mysterious about homotopy cardinality.

Homotopy cardinality is formally very similar to Euler characteristic.
The Euler characteristic \(\chi(X)\) is given by the alternating sum
\[\chi(X) = \dim(H_0(X)) - \dim(H_1(X)) + \dim(H_2(X)) - \ldots\]
whenever the sum converges, where \(H_n(X)\) is a vector space over the
rational numbers called the \(n\)th rational homology group of \(X\).
Just as for homotopy cardinality, we have \[
  \begin{aligned}
    \chi(X + Y) &= \chi(X) + \chi(Y),
  \\\chi(X \times Y) &= \chi(X) \times \chi(Y)
  \end{aligned}
\] and more generally, whenever \[F\to X\to B\] is a fibration and the
base space \(B\) is connected, we have
\[\chi(X) = \chi(F)\times \chi(B)\] whenever any two of the spaces have
well-defined Euler characteristic, which implies that the third does too
(unless I'm confused).

So Euler characteristic is a lot like homotopy cardinality. But not many
spaces have \emph{both} well-defined homotopy cardinality \emph{and}
well-defined Euler characteristic. So they're like Jekyll and Hyde ---
you hardly ever see them in the same place at the same time, so you
can't tell if they're really the same guy.

But there are some weird ways to try to force the issue and compute both
quantities for certain spaces. For example, suppose \(G\) is a finite
group. Then we can build \(BG\) starting from a simplicial set with 1
nondegenerate 0-simplex, \(|G|-1\) nondegenerate \(1\)-simplices,
\((|G|-1)^2\) nondegenerate 2-simplices, and so on. If there were only
finitely many nondegenerate simplices of all dimensions, we could
compute the Euler characteristic of this space as the alternating sum of
the numbers of such simplices. So let's try doing that here! We get:
\[\chi(BG) = 1 - (|G|-1) + (|G|-1)^2 - \ldots\] Of course the sum
diverges, but if we go ahead and use the geometric formula anyway, we
get \[\chi(BG) = 1/|G|\] which matches our previous (rigorous) result
that \[|BG| = 1/|G|.\] So maybe they're the same after all! There are
similar calculations like this in James Propp's paper ``Exponentiation
and Euler characteristic'', referred to above\ldots{} though he uses a
slightly different notion of Euler characteristic, due to Schanuel.
Clearly something interesting is going on with these ``divergent Euler
characteristics''. For appearances of this sort of thing in physics,
see:

\begin{enumerate}
\def\labelenumi{\arabic{enumi})}
\setcounter{enumi}{10}
\tightlist
\item
  Matthias Blau and George Thompson, ``\(N = 2\) topological gauge
  theory, the Euler characteristic of moduli spaces, and the Casson
  invariant'', \emph{Comm. Math. Phys.} \textbf{152} (1993), 41--71.
\end{enumerate}
\noindent
and the references therein. (I discussed this paper a bit in
\protect\hyperlink{week51}{``Week 51''}.)

However, there are still challenging tests to the theory that homotopy
cardinality and Euler characteristic are secretly the same. Here's a
puzzle due to James Dolan. Consider a Riemann surface of genus
\(g > 1\). Such a surface has Euler characteristic \(2 - 2g\), but such
a surface also has vanishing homotopy groups above the first, which
implies that it's \(BG\) for \(G\) equal to its fundamental group. If
homotopy cardinality and Euler characteristic were the same, this would
imply \[|G| = \frac{1}{|BG|} = \frac{1}{\chi(S)} = \frac{1}{2 - 2g}\]
But the fundamental group \(G\) is infinite! What's going on?

Well, I'm actually sort of glad that \(1/(2 - 2g)\) is \emph{negative}.
Sometimes a divergent series of positive integers can be cleverly summed
up to give a negative number. The simplest example is the geometric
series \[1 + 2 + 4 + 8 + 16 + \ldots = \frac{1}{1 - 2} = -1\] but in
\protect\hyperlink{week126}{``Week 126''} I talked about a more
sophisticated example that is very important in string theory:
\[1 + 2 + 3 + 4 + 5 + \ldots = \zeta(-1) = -\frac{1}{12}\] So maybe some
similar trickery can be used to count the elements of \(G\) and get a
divergent sum that can be sneakily evaluated to obtain \(1/(2-2g)\). Of
course, even if we succeed in doing this, the skeptics will rightly
question the significance of such tomfoolery. But there is sometimes a
lot of profound truth lurking in these bizarre formal manipulations, and
sometimes if you examine what's going on carefully enough, you can
discover cool stuff.

To wrap up, let me mention an interesting paper on the foundations of
categorification:

\begin{enumerate}
\def\labelenumi{\arabic{enumi})}
\setcounter{enumi}{11}
\tightlist
\item
  Claudio Hermida, ``From coherent structures to universal properties'', \emph{Jour.\ Pure Appl.\ Alg.} \textbf{165} (2001), 7--61.  Also
  available as \href{https://arxiv.org/abs/math/0006161}{\texttt{math/0006161}}.
\end{enumerate}
\noindent
and also two papers about \(2\)-groupoids and topology:

\begin{enumerate}
\def\labelenumi{\arabic{enumi})}
\setcounter{enumi}{12}
\tightlist
\item
   K. A. Hardie, K. H. Kamps, R. W. Kieboom, ``A homotopy \(2\)-groupoid of
   a Hausdorff space'', \emph{Appl. Categ. Structures} \textbf{8}, (2000), 209--234.

  K. A. Hardie, K. H. Kamps, R. W. Kieboom, ``A homotopy bigroupoid of a
  topological space'', \emph{Appl. Categ. Structures} \textbf{9}, (2001), 311--327.
\end{enumerate}
\noindent
I would talk about these if I had the energy, but it's already past my
bed-time. Good night!

\begin{center}\rule{0.5\linewidth}{0.5pt}\end{center}

\textbf{Addenda:} Toby Bartels had some interesting things to say about
this issue of This Week's Finds. Here is my reply, which quotes some of
his remarks\ldots.

\begin{center}\rule{0.5\linewidth}{0.5pt}\end{center}

\begin{quote}
\begin{enumerate}
\def\labelenumi{\arabic{enumi})}
\setcounter{enumi}{2}
\tightlist
\item
  The American Physical Society: ``A Century of Physics'', available at
  {\rm \href{https://web.archive.org/web/19990508143827/http://timeline.aps.org/APS/home_HighRes.html}{\texttt{https://web.archive.org/web/19990508143827/http://timeline.aps.org/APS/home\_HighRes.html}}}
\end{enumerate}
\end{quote}

\begin{quote}
{\bf TB: }I like how they make the famous picture of Buzz Aldrin, the one
that everyone thinks is a picture of Neil Armstrong, into a picture of
Neil Armstrong after all: ``Here he is reflected in Buzz Aldrin's
visor.''.
\end{quote}

\begin{quote}
{\bf JB: }Heh. Sounds like something a doting grandmother would say!
\end{quote}

\begin{center}\rule{0.5\linewidth}{0.5pt}\end{center}

\begin{quote}
\begin{enumerate}
\def\labelenumi{\arabic{enumi})}
\setcounter{enumi}{4}
\tightlist
\item
 John Conway and Peter Doyle, ``Division by three''. Available at
  {\rm \href{https://arxiv.org/abs/math/0605779}{math/0605779}}.
\end{enumerate}

This article studies this question: if I give you an isomorphism
between \(3x\) and \(3y\), can you construct a isomorphism between \(x\)
and \(y\)?
\end{quote}

\begin{quote}
{\bf TB: }The answer must be something that won't work if 3 is replaced
by an infinite cardinal. That said, I can't even figure out how to
divide by 2! If I take the 3 copies of \(X\) or \(Y\) and put them on
top of each other, I get a finite, 2-coloured, 3-valent, nonsimple,
undirected graph. I remember from combinatorics that the 2 colours of a
finite, 2-coloured, simple, undirected graph of fixed valency are
equipollent, but I can't remember the bijective proof. (Presumably it
can be adopted to nonsimple graphs.)
\end{quote}

\begin{quote}
{\bf JB: }It's a tricky business. Let me quote from the above article:

\vskip 0.5em
\textbf{History}

A proof that it is possible to divide by two was presented by Bernstein
in his Inaugural Dissertation of 1901, which appeared in Mathematische
Annallen in 1905; Bernstein also indicated how to extend his results to
division by any finite \(n\), but we are not aware of anyone other than
Bernstein himself who ever claimed to understand this argument. In 1922
Sierpinski published a simpler proof of division by two, and he worked
hard to extend his method to division by three, but never succeeded.

In 1927 Lindenbaum and Tarski announced, in an infamous paper that
contained statements (without proof) of 144 theorems of set theory, that
Lindenbaum had found a proof of division by three. Their failure to give
any hint of a proof must have frustrated Sierpinski, for it appears that
twenty years later he still did not know how to divide by three.
Finally, in 1949, in a paper `dedicated to Professor Waclaw Sierpinski
in celebration of his forty years as teacher and scholar', Tarski
published a proof. In this paper, Tarski explained that unfortunately he
couldn't remember how Lindenbaum's proof had gone, except that it
involved an argument like the one Sierpinski had used in dividing by
two, and another lemma, due to Tarski, which we will describe below.
Instead of Lindenbaum's proof, he gave another.

Now when we began the investigations reported on here, we were aware
that there was a proof in Tarski's paper, and Conway had even pored over
it at one time or another without achieving enlightenment. The problem
was closely related to the kind of question John had looked at in his
thesis, and it was also related to work that Doyle had done in the field
of bijective combinatorics. So we decided that we were going to figure
out what the heck was going on. Without too much trouble we figured out
how to divide by two. Our solution turned out to be substantially
equivalent to that of Sierpinski, though the terms in which we will
describe it below will not much resemble Sierpinski's. We tried and
tried and tried to adapt the method to the case of dividing by three,
but we kept getting stuck at the same point in the argument. So finally
we decided to look at Tarski's paper, and we saw that the lemma Tarski
said Lindenbaum had used was precisely what we needed to get past the
point we were stuck on! So now we had a proof of division by three that
combined an argument like that Sierpinski used in dividing by two with
an appeal to Tarski's lemma, and we figured we must have hit upon an
argument very much like that of Lindenbaum's. This is the solution we
will describe here: Lindenbaum's argument, after 62 years.
\end{quote}

\begin{center}\rule{0.5\linewidth}{0.5pt}\end{center}

\begin{quote}
So: let's define the ``homotopy cardinality'' of a topological space
\(X\) to be the alternating product
\(|X| = \prod_{i>0} |\pi_i(X)|^{(-1)^i}\) when \(X\) is connected and
the product converges.
\end{quote}

\begin{quote}
{\bf TB: }What about divergence to \(0\)? If \(\pi_i(X)\) is infinite for
some odd \(i\) but no even \(i\), can we say \(|X|\) is \(0\)?
\end{quote}

\begin{quote}
{\bf JB: }Well, we can, but we might regret it later. In a sense \(0\) is
no better than \(\infty\) when one is doing products, so if you allow
\(0\) as a legitimate value for a homotopy cardinality, you should allow
\(\infty\), but if you allow both, you get in trouble when you try to
multiply them. This dilemma is familiar from the case of infinite sums
(where \(+\infty\) and \(-\infty\) are the culprits), and the resolution
seems to be either:

\begin{itemize}
\tightlist
\item
  disallow both \(0\) and \(\infty\) as legitimate answers for the above
  product,
\end{itemize}

or

\begin{itemize}
\tightlist
\item
  allow both but then be extra careful when stating your theorems so
  that you don't run into problems.
\end{itemize}
\end{quote}

\begin{center}\rule{0.5\linewidth}{0.5pt}\end{center}

\begin{quote}
As a fun application of this fact, recall that any topological group
\(G\) has a ``classifying space'' \(BG\), meaning a space with a
principal \(G\)-bundle over it \(G\to EG\to BG\) whose total space
\(EG\) is contractible. I described how to construct the classifying
space in \protect\hyperlink{week117}{``Week 117''}, at least in the case
of a discrete group \(G\), but I didn't say much about why it's so
great. The main reason it's great is that any \(G\)-bundle over any
space is a pullback of the bundle \(EG\) over \(BG\). But right now,
what I want to note is that since \(EG\) is contractible it is tame, and
\(|EG| = 1\). Thus \(G\) is tame if and only if \(BG\) is, and
\(|BG| = 1 / |G|\), so we can think of \(BG\) as the reciprocal' of
\(G\)!
\end{quote}

\begin{quote}
{\bf TB: }On the other hand, \(G\) is already a kind of reciprocal of itself. If \(G\)
is a discrete group, it's a topological space with
\(|G|_\mathrm{homotopy} = |G|_\mathrm{set}\). But \(G\) is also a groupoid
with 1 object, and \(|G|_\mathrm{groupoid} = 1 / |G|_\mathrm{set}\). So,
\(|G|_\mathrm{homotopy} |G|_\mathrm{groupoid} = 1\).
\end{quote}

\begin{quote}
{\bf JB: }Believe it or not, you are reinventing \(BG\)! A groupoid can
be reinterpreted as a space with vanishing homotopy groups above the
first, and if you do this to the groupoid \(G\), you get \(BG\).

More generally:

Recall that we can take a pointed space \(X\) and form a pointed space
\(LX\) of loops in \(X\) that start and end at the basepoint. This
clearly has \[\pi_{n+1}(X) = \pi_n(LX)\] so if \(X\) is connected and
tame we'll have \[|LX| = 1/|X|\] Now with a little work you can make
\(LX\) (or a space homotopy-equivalent to it!) into a topological group
with composition of loops as the product.

And then it turns out that \(BLX\) is homotopy equivalent to \(X\) when
\(X\) is connected. Conversely, given a topological group \(G\), \(LBG\)
is homotopy equivalent to \(G\).

So what we're seeing is that topological groups and connected pointed
spaces are secretly the same thing, at least from the viewpoint of
homotopy theory. In topology, few things are as important as this fact.

But what's really going on here? Well, to go from a topological group
\(G\) to a connected pointed space, you have to form \(BG\), which has
all the same homotopy groups but just pushed up one notch:
\[\pi_{n+1}(BG) = \pi_n(G)\] And to go from a connected point space
\(X\) to a topological group, you have to form \(LX\), which has all the
same homotopy groups but just pushed down one notch:
\[\pi_{n-1}(LX) = \pi_n(X)\] This is actually the trick you are playing,
in slight disguise.

And the real point is that a 1-object \(\omega\)-groupoid can be
reinterpreted as an \(\omega\)-groupoid by forgetting about the object
and renaming all the \(j\)-morphisms ``\((j-1)\)-morphisms''.

See? When you finally get to the bottom of it, this ``\(BG\)'' business
is just a silly reindexing game!!! Of course no textbook can admit this
openly --- partially because they don't talk about \(\omega\)-groupoids.
\end{quote}

\begin{center}\rule{0.5\linewidth}{0.5pt}\end{center}

\begin{quote}
So Euler characteristic is a lot like homotopy cardinality. But not many
spaces have both well-defined homotopy cardinality and well-defined
Euler characteristic. So they're like Jekyll and Hyde --- you hardly
ever see them in the same place at the same time, so you can't tell if
they're really the same guy.
\end{quote}

\begin{quote}
{\bf TB: }So, are they ever both defined but different?
\end{quote}

\begin{quote}
{\bf JB: }I don't recall any examples where they're both finite, but
different. I know very few cases where they're both finite! How about
the point? How about the circle? How about the 2-sphere? I leave you to
ponder these cases.
\end{quote}

\begin{center}\rule{0.5\linewidth}{0.5pt}\end{center}

\begin{quote}
However, there are still challenging tests to the theory that homotopy
cardinality and Euler characteristic are secretly the same. Here's a
puzzle due to James Dolan. Consider a Riemann surface of genus
\(g > 1\). Such a surface has Euler characteristic \(2 - 2g\), but such
a surface also has vanishing homotopy groups above the first, which
implies that it's \(BG\) for \(G\) equal to its fundamental group. If
homotopy cardinality and Euler characteristic were the same, this would
imply \[|G| = 1/|BG| = 1/\chi(S) = 1/(2 - 2g)\].

But the fundamental group \(G\) is infinite! What's going on?
\end{quote}

\begin{quote}
{\bf TB: }This doesn't seem too surprising. \(1/(2 - 2g)\) is also
infinite. Just use the geometric series in reverse:
\[1/(2 - 2g) = (1/2) \sum_i g^i,\]

which diverges since \(g > 1\).
\end{quote}

\begin{quote}
{\bf JB: }Well, what I really want is a way of counting elements of the
fundamental group of the surface \(S\) which gives me a divergent sum
that I can cleverly sum up to get \(1/(2 - 2g)\).
\end{quote}

\begin{center}\rule{0.5\linewidth}{0.5pt}\end{center}

Later, my wish above was granted by Laurent Bartholdi and Danny
Ruberman! People have already figured out how to count the number of
elements in the fundamental group of a Riemann surface, resum, and get
\(1/(2 - 2g)\) in a nice way. Here are two references:

\begin{enumerate}
\def\labelenumi{\arabic{enumi})}
\setcounter{enumi}{13}
\item
  William J. Floyd and Steven P. Plotnick, ``Growth functions on
  Fuchsian groups and the Euler characteristic'', \emph{Invent. Math.}
  \textbf{88} (1987), 1--29.
\item
  R. I. Grigorchuk, ``Growth functions, rewriting systems and Euler
  characteristic'', \emph{Mat. Zametki} \textbf{58} (1995), 653--668,
  798.
\end{enumerate}

You can read more about Euler characteristic and homotopy cardinality
here:

\begin{enumerate}
\def\labelenumi{\arabic{enumi})}
\setcounter{enumi}{15}
\tightlist
\item
  John Baez, ``Euler characteristic versus homotopy cardinality'',
  lecture at the \emph{Fields Institute Program on Applied Homotopy
  Theory}, September 20, 2003. Available at
  \url{http://www.math.ucr.edu/home/baez/cardinality/}
\end{enumerate}

\begin{center}\rule{0.5\linewidth}{0.5pt}\end{center}

\begin{quote}
\emph{The imaginary expression \(\sqrt{-a}\) and the negative expression
\(-b\) have this resemblance, that either of them occurring as the solution of a problem indicates some inconsistency or absurdity.  As far as real meaning is concerned, both are imaginary, since \(0 - a\) is as inconceivable as \(\sqrt{-a}\).}

--- Augustus De Morgan, 1831
\end{quote}

\hypertarget{week148}{%
\section{June 5, 2000}\label{week148}}

Last week I talked about some millennium-related books. This week, some
millennial math problems! In 1900, at the second International Congress
of Mathematicians, Hilbert posed a famous list of 23 problems. No one
individual seems to have the guts to repeat that sort of challenge now.
But the newly-founded Clay Mathematics Institute, based in Cambridge
Massachusetts and run by Arthur Jaffe, has just laid out a nice list of
7 problems:

\begin{enumerate}
\def\labelenumi{\arabic{enumi})}
\tightlist
\item
  Clay Mathematics Institute, ``Millennium Prize Problems'',
  \href{http://www.claymath.org/millennium-problems}{\texttt{http://www.claymath.}}   \href{http://www.claymath.org/millennium-problems}{\texttt{org/millennium-problems}}
\end{enumerate}

There is a 1 million dollar prize for each one! Unlike most of Hilbert's
problems, these weren't cooked up specially for the occasion: they have
already proved their merit by resisting attack for some time.

Here they are:

\vskip 1em

\textbf{1. P = NP?}

This is the newest problem on the list and the easiest to explain. An
algorithm is ``polynomial-time'' if the time it takes to run is bounded
by some polynomial in the length of the input data. This is a crude but
easily understood condition to decide whether an algorithm is fast
enough to be worth bothering with. A ``nondeterministic
polynomial-time'' algorithm is one that can \emph{check} a purported
solution to a problem in an amount of time bounded by some polynomial in
the input data. All algorithms in P are in NP, but how about the
converse? Is P = NP? Stephen Cook posed this problem in 1971 and it's
still open. It seems unlikely to be true --- a good candidate for a
counterexample is the problem of factoring integers --- but nobody has
\emph{proved} that it's false. This is the most practical question of
the lot, because if the answer were ``yes'', there's a chance that one
could use this result to quickly crack all the current best encryption
schemes.

\vskip 1em
\textbf{2. The Poincaré conjecture}

Spheres are among the most fundamental topological spaces, but spheres
hold many mysteries. For example: is every \(3\)-dimensional manifold
with the same homotopy type as a 3-sphere actually homeomorphic to a
3-sphere? Or for short: are homotopy 3-spheres really 3-spheres?
Poincaré posed this puzzle in 1904 shortly after he knocked down an
easier conjecture of his by finding 3-manifolds with the same homology
groups as 3-spheres that weren't really 3-spheres. The
higher-dimensional analogues of Poincaré's question have all been
settled in the affirmative --- Smale, Stallings and Wallace solved it in
dimensions 5 and higher, and Freedman later solved the subtler
4-dimensional case --- but the \(3\)-dimensional case is still unsolved.
This is an excellent illustration of a fact that may seem surprising at
first: many problems in topology are toughest in fairly low dimensions!
The reason is that there's less ``maneuvering room''. The last couple
decades have seen a burst of new ideas in low-dimensional topology ---
this has been a theme of This Week's Finds ever since it started --- but
the Poincaré conjecture remains uncracked.

\vskip 1em
\textbf{3. The Birch--Swinnerton-Dyer conjecture}

This is a conjecture about elliptic curves, and indirectly, number
theory. For a precise definition of an elliptic curve I'll refer you to
\protect\hyperlink{week13}{``Week 13''} and
\protect\hyperlink{week125}{``Week 125''}, but basically, it's a
torus-shaped surface described by an algebraic equation like this:
\[y^2 = x^3 + ax + b\] Any elliptic curve is naturally an abelian group,
and the points on it with rational coordinates form a finitely generated
subgroup. When are there infinitely many such rational points? In 1965,
Birch and Swinnerton-Dyer conjectured a criterion involving something
called the ``\(L\)-function'' of the elliptic curve. The \(L\)-function
\(L(s)\) is an elegant encoding of how many solutions there are to the
above equation modulo \(p\), where \(p\) is any prime. The
Birch-Swinnerton-Dyer conjecture says that \(L(1) = 0\) if and only if
the elliptic curve has infinitely many rational points. More generally,
it says that the order of the zero of \(L(s)\) at \(s = 1\) equals the
rank of the group of rational points on the elliptic curve (that is, the
rank of the free abelian summand of this group.) A solution to this
conjecture would shed a lot of light on Diophantine equations, one of
which goes back to at least the 10th century --- namely, the problem of
finding which integers appear as the areas of right triangles all of
whose sides have lengths equal to rational numbers.

\vskip 1em
\textbf{4. The Hodge conjecture}

This question is about algebraic geometry and topology. A ``projective
nonsingular complex algebraic variety'' is basically a compact smooth
manifold described by a bunch of homogeneous complex polynomial
equations. Such a variety always has even dimension, say \(2n\). We can
take the DeRham cohomology of such a variety and break it up into parts
\(H^{p,q}\) labelled by pairs \((p,q)\) of integers between \(0\) and
\(n\), using the fact that every function is a sum of a holomorphic and
an antiholomorphic part. Sitting inside the DeRham cohomology is the
rational cohomology, The rational guys inside \(H^{p,p}\) are called
``Hodge forms''. By Poincaré duality any closed analytic subspace of our
variety defines a Hodge form --- this sort of Hodge form is called an
algebraic cycle. The Hodge conjecture, posed in 1950 states: every Hodge
form is a rational linear linear combination of algebraic cycles. It's
saying that we can concretely realize a bunch of cohomology classes
using closed analytic subspaces sitting inside our variety.

\vskip 1em
\textbf{5. Existence and mass gap for Yang--Mills theory}

One of the great open problems of modern mathematical physics is whether
the Standard Model of particle physics is mathematically consistent.
It's not even known whether ``pure'' Yang--Mills theory --- uncoupled to
fermions or the Higgs --- is a well-defined quantum field theory with
reasonable properties. To make this question precise, people have
formulated various axioms for a quantum field theory, like the so-called
``Haag-Kastler axioms''. The job of constructive quantum field theory is
to mathematically study questions like whether we can construct
Yang--Mills theory in such a way that it satisfies these axioms. But one
really wants to know more: at the very least, existence of Yang--Mills
theory coupled to fermions, together with a ``mass gap'' --- i.e., a
nonzero minimum mass for the particles formed as bound states of the
theory (like protons are bound states of quarks).

\vskip 1em
\newpage
\textbf{6. Existence and smoothness for the Navier--Stokes equations}

The Navier--Stokes equations are a set of partial differential equations
describing the flow of a viscous incompressible fluid. If you start out
with a nice smooth vector field describing the flow of some fluid, it
will often get complicated and twisty as turbulence develops. Nobody
knows whether the solution exists for all time, or whether it develops
singularities and becomes undefined after a while! In fact, numerical
evidence hints at the contrary. So one would like to know whether
solutions exist for all time and remain smooth --- or at least find
conditions under which this is the case. Of course, the Navier--Stokes
equations are only an approximation to the actual behavior of fluids,
since it idealizes them as a continuum when they are actually made of
molecules. But it's important to understand whether and how the
continuum approximation breaks down as turbulence develops.

\vskip 1em
\textbf{7. The Riemann hypothesis}

For \(\Re(s) > 1\) the Riemann zeta function is defined by
\[\zeta(s) = \frac{1}{1^s} + \frac{1}{2^s} + \frac{1}{3^s} + \ldots\]
But we can extend it by analytic continuation to most of the complex
plane --- it has a pole at \(s = 1\). The zeta function has a bunch of
zeros in the ``critical strip'' where \(\Re(s)\) is between \(0\) and
\(1\). In 1859, Riemann conjectured that all such zeros have real part
equal to \(1/2\). This conjecture has lots of interesting ramifications
for things like the distribution of prime numbers. By now, more than a
billion zeros in the critical strip have been found to have real part
\(1/2\); it has also been shown that ``most'' such zeros have this
property, but the Riemann hypothesis remains open.

If you solve one of these conjectures and win a million dollars because
you read about it here on This Week's Finds, please put me in your will.

Okay, now on to some other stuff.

This week was good for me in two ways. First of all, Ashtekar, Krasnov
and I finally finished a paper on black hole entropy that we've been
struggling away on for over 3 years. I can't resist talking about this
paper at length, since it's such a relief to be done with it. Second,
the guru of \(n\)-category theory, Ross Street, visited Riverside and
explained a bunch of cool stuff to James Dolan and me. I may talk about
this next time.

\begin{enumerate}
\def\labelenumi{\arabic{enumi})}
\setcounter{enumi}{1}
\tightlist
\item
  Abhay Ashtekar, John Baez and Kirill Krasnov, ``Quantum geometry of
  isolated horizons and black hole entropy'', \emph{Adv.\ Theor.\ Math.\ Phys.}
   \textbf{4}  (2001), 1--94.   Also available at
  \href{https://arxiv.org/abs/gr-qc/0005126}{\texttt{gr-qc/0005126}}.
\end{enumerate}
\noindent
I explained an earlier version of this paper in
\protect\hyperlink{week112}{``Week 112''}, but now I want to give a more
technical explanation. So:

The goal of this paper is to understand the geometry of black holes in a
way that takes quantum effects into account, using the techniques of
loop quantum gravity. We do not consider the region near the
singularity, which is poorly understood. Instead, we focus on the
geometry of the event horizon, since we wish to compute the entropy of a
black hole by counting the microstates of its horizon.

Perhaps I should say a word about why we want to do this. As explained
in \protect\hyperlink{week111}{``Week 111''}, Bekenstein and Hawking
found a formula relating the entropy \(S\) of a black hole to the area
\(A\) of its event horizon. It is very simple: \[S = A/4\] in units
where the speed of light, Newton's constant, Boltzmann's constant and
Planck's constant equal \(1\). Now, in quantum statistical mechanics,
the entropy of a system in thermal equilibrium is roughly the logarithm
of the number \(N\) of microstates it can occupy: \[S = \ln N.\] This is
exactly right when all the microstates have the same energy. Thus we
expect that a black hole of area \(A\) has about \[N = \exp(A/4)\]
microstates. For a solar-mass black hole, that's about \(\exp(10^{76})\)
microstates! Any good theory of quantum gravity must explain what these
microstates are. Since their number is related to the event horizon's
area, it is natural to guess that they're related to the geometry of the
event horizon. But how?

It's clear that everything will work perfectly if each little patch of
the event horizon with area \(4 \ln(2)\) has exactly 2 states. I think
Wheeler was the first to take this seriously enough to propose a toy
model where each such patch stores one bit of information, making the
black hole into something sort of like an enormous hard drive:

\begin{enumerate}
\def\labelenumi{\arabic{enumi})}
\setcounter{enumi}{2}
\tightlist
\item
  John Wheeler, ``It from bit'', in \emph{Sakharov Memorial Lecture on
  Physics}, Volume 2, eds.~L. Keldysh and V. Feinberg, Nova Science, New
  York, 1992.
\end{enumerate}
\noindent
Of course, this idea sounds a bit nutty. However, the quantum state of a
spinor contains exactly one bit of information, and loop quantum gravity
is based on the theory of spinors, so it's not as crazy as it might
seem\ldots. Still, there are some, ahem, \emph{details} to be worked
out!

So let's work them out.

The first step is to understand the classical mechanics of a black hole
in a way that allows us to apply the techniques of loop quantum gravity.
In other words, we want to describe a classical phase space for our
black hole. This step was done in a companion paper:

\begin{enumerate}
\def\labelenumi{\arabic{enumi})}
\setcounter{enumi}{3}
\tightlist
\item
  Abhay Ashtekar, Alejandro Corichi and Kirill Krasnov, ``Isolated
  horizons: the classical phase space'', \emph{Adv.\ Theor.\
   Math.\ Phys.} \textbf{3} (2000), 418--471.  Also
  available as
  \href{https://arxiv.org/abs/gr-qc/9905089}{\texttt{gr-qc/9905089}}.
\end{enumerate}

The idea is to consider the region of spacetime outside the black hole
and assume that its boundary is a cylinder of the form
\(\mathbb{R}\times S^2\). We demand that this boundary is an ``isolated
horizon'' --- crudely speaking, a surface that light cannot escape from,
with no matter or gravitational radiation falling in for the stretch of
time under consideration. To make this concept precise we need to impose
some boundary conditions on the metric and other fields at the horizon.
These are most elegantly described using Penrose's spinor formalism for
general relativity, as discussed in \protect\hyperlink{week109}{``Week
109''}. With the help of these boundary conditions, we can start with
the usual Lagrangian for general relativity, turn the crank, and work
out a description of the phase space for an isolated black hole.

If we temporarily ignore the presence of matter, a point in this phase
space describes the metric and extrinsic curvature of space outside the
black hole at a given moment of time. Technically, we do this using an
\(\mathrm{SU}(2)\) connection \(A\) together with an
\(\mathfrak{su}(2)\)-valued \(2\)-form \(E\). You can think of these as
analogous to the vector potential and electric field in
electromagnetism. As usual, they need to satisfy some constraints coming
from Einstein's equations for general relativity. They also need to
satisfy boundary conditions coming from the definition of an isolated
horizon.

Since the black hole is shaped like a ball, the boundary conditions hold
on a 2-sphere that I'll call the ``horizon 2-sphere''. One thing the
boundary conditions say is that on the horizon 2-sphere, the
\(\mathrm{SU}(2)\) connection \(A\) is completely determined by a
\(\mathrm{U}(1)\) connection, say \(W\). This \(\mathrm{U}(1)\)
connection is really important, because it describes the intrinsic
geometry of the horizon 2-sphere. Here's a good way to think about it:
first you restrict the spacetime metric to the horizon 2-sphere, and
then you work out the Levi-Civita connection of this metric on the
2-sphere. Finally, since loop quantum gravity is based on the parallel
transport of spinors, you work out the corresponding connection for
spinors on the 2-sphere, which is a \(\mathrm{U}(1)\) connection. That's
\(W\)!

The boundary conditions also say that on the horizon 2-sphere, the \(E\)
field is proportional to the curvature of \(W\). So on the horizon
2-sphere, \emph{all} the fields are determined by \(W\). This is even
true when we take the presence of matter into account. When we quantize,
it'll be the microstates of this field \(W\) that give rise to the black
hole entropy. Since \(W\) is just a technical way of describing the
shape of the horizon 2-sphere, this means that the black hole entropy
arises from the many slightly different possible shapes that the horizon
can have.

But I'm getting ahead of myself here! We haven't quantized yet; we're
just talking about the classical phase space for an isolated black hole.

The most unusual feature of this phase space is that its symplectic
structure is a sum of two terms. First, there is the usual integral over
space at a given time, which makes the \(E\) field canonically conjugate
to the \(A\) field away from the horizon 2-sphere. But then there is a
boundary term: an integral over the horizon 2-sphere. This gives the
geometry of the horizon a life of its own, which ultimately accounts for
the black hole entropy. Not surprisingly, this boundary term involves
the \(\mathrm{U}(1)\) connection \(W\). In fact, this boundary term is
just the symplectic structure for \(\mathrm{U}(1)\) Chern--Simons theory
on the 2-sphere! It's the simplest thing you can write down:
\[\omega(\delta W, \delta W') = \frac{k}{2\pi} \int \delta W\wedge \delta W'\]
Here \(\omega\) is the \(\mathrm{U}(1)\) Chern--Simons symplectic
structure; we're evaluating it on two tangent vectors to the space of
\(\mathrm{U}(1)\) connections on the 2-sphere, which we call
\(\delta W\) and \(\delta W'\). These are the same as \(1\)-forms, so we
can wedge them and integrate the result over the 2-sphere. The number
\(k\) is some constant depending on the area of the black hole\ldots{}
but more about that later!

I guess this Chern--Simons stuff needs some background to fully
appreciate. I have been talking about it for a long time here on This
Week's Finds. The quantum version of Chern--Simons theory is a
3-dimensional quantum field theory that burst into prominence thanks to
Witten's work relating it to the Jones polynomial, which is an invariant
of knots. At least heuristically, you can calculate the Jones polynomial
by doing a path integral in \(\mathrm{SU}(2)\) Chern--Simons theory. It
also turns out that Chern--Simons theory is deeply related to quantum
gravity in 3d spacetime. For quite a while, various people have hoped
that Chern--Simons theory was important for quantum gravity in 4d
spacetime, too --- see for example \protect\hyperlink{week56}{``Week
56''} and \protect\hyperlink{week57}{``Week 57''}. However, there have
been serious technical problems in most attempts to relate Chern--Simons
theory to physically realistic problems in 4d quantum gravity. I think
we may finally be straightening out some of these problems! But the
ironic twist is that we're using \(\mathrm{U}(1)\) Chern--Simons theory,
which is really very simple compared to the sexier \(\mathrm{SU}(2)\)
version. For example, \(\mathrm{U}(1)\) Chern--Simons theory also gives a
knot invariant, but it's basically just the self-linking number. And the
math of \(\mathrm{U}(1)\) Chern--Simons theory goes back to the 1800s ---
it's really just the mathematics of ``theta functions''.

As a historical note, I should add that the really nice derivation of
the Chern--Simons boundary term in the symplectic structure for isolated
black holes was found in a paper written \emph{after} the one I
mentioned above:

\begin{enumerate}
\def\labelenumi{\arabic{enumi})}
\setcounter{enumi}{4}
\tightlist
\item
  Abhay Ashtekar, Chris Beetle and Steve Fairhurst, ``Mechanics of
  isolated horizons'', \emph{Class.\ Quant.\ Grav.} \textbf{17}
  (2000), 253--298. Available at
  \href{https://arxiv.org/abs/gr-qc/9907068}{\texttt{gr-qc/9907068}}.
\end{enumerate}

Originally, everyone thought that to make the action differentiable as a
function of the fields, you had to add a boundary term to the usual
action for general relativity, and that this boundary term was
responsible for the boundary term in the symplectic structure. This
seemed a bit ad hoc. Of course, you need to differentiate the action to
get the field equations, so it's perfectly sensible to add an extra term
if that's what you need to make the action differentiable, but still you
wonder: where did the extra term COME FROM?

Luckily, Ashtekar and company eventually realized that while you
\emph{can} add an extra term to the action, you don't really \emph{need}
to. By cleverly using the ``isolated horizon'' boundary conditions, you
can show that the usual action for general relativity is already
differentiable without any extra term, and that it yields the
Chern--Simons boundary term in the symplectic structure.

Okay: we've got a phase space for an isolated black hole, and we've got
the symplectic structure on this phase space. Now what?

Well, now we should quantize this phase space! It's a bit complicated,
but thanks to the two-part form of the symplectic structure, it
basically breaks up into two separate problems: quantizing the A field
and its canonical conjugate \(E\) outside the horizon 2-sphere, and
quantizing the \(W\) field on this 2-sphere. The first problem is
basically just the usual problem of loop quantum gravity --- people know
a lot about that. The second problem is basically just quantizing
\(\mathrm{U}(1)\) Chern--Simons theory --- people know even \emph{more}
about that! But then you have to go back and put the two pieces
together. For that, it's crucial that on the horizon, the \(E\) field is
proportional to the curvature of the connection \(W\).

So: what do quantum states in the resulting theory look like? I'll
describe a basis of states for you\ldots.

Outside the black hole, they are described by spin networks. I've
discussed these in \protect\hyperlink{week110}{``Week 110''} and
elsewhere, but let me just recall the basics. A spin network is a graph
whose edges are labelled by irreducible representations of
\(\mathrm{SU}(2)\), or in other words spins \(j = 0, 1/2, 1, \ldots\).
Their vertices are labelled as well, but that doesn't concern us much
here. What matters more is that the spin network edges can puncture the
horizon 2-sphere. And it turns out that each puncture must be labelled
with a number \(m\) chosen from the set \[\{-j, -j+1, \ldots, j-1, j\}\]
These numbers \(m\) determine the state of the geometry of the horizon
2-sphere.

What do these numbers \(j\) and \(m\) really \emph{mean?} Well, they should be
vaguely familiar if you've studied the quantum mechanics of angular
momentum. The same math is at work here, but with a rather different
interpretation. Spin network edges represent quantized flux lines of the
gravitational \(E\) field. When a spin network edge punctures the
horizon 2-sphere, it contributes \emph{area} to the 2-sphere: a
spin-\(j\) edge contributes an area equal to \[8\pi\gamma\sqrt{j(j+1)}\]
for some constant \(\gamma\).

But due to the boundary conditions relating the \(E\) field to the
curvature of the connection \(W\), each spin network edge also
contributes \emph{curvature} to the horizon 2-sphere. In fact, this
2-sphere is flat except where a spin network edge punctures it; at the
punctures it has cone-shaped singularities. You can form a cone by
cutting out a wedge-shaped slice from a piece of paper and reattaching
the two new edges, and the shape of this cone is described by the
``deficit angle'' --- the angle of the wedge you removed. In this black
hole business, a puncture labelled by the number \(m\) gives a conical
singularity with a deficit angle equal to \[\frac{4\pi m}{k}\] where
\(k\) is the same constant appearing in the formula for the Chern-
Simons symplectic structure.

I guess now it's time to explain these mysterious constants! First of
all, \(\gamma\) is an undetermined dimensionless constant, usually
called the ``Immirzi parameter'' because it was first discovered by
Fernando Barbero. This parameter sets the scale at which area is
quantized! Of course, the formula for the area contributed by a
spin-\(j\) edge: \[8\pi\gamma\sqrt{j(j+1)}\] also has a factor of the
Planck area lurking in it, which you can't see because I've set \(c\),
\(G\), and \(\hbar\) to \(1\). That's not surprising. What's surprising
is the appearance of the Barbero-Immirzi parameter. So far, loop quantum
gravity cannot predict the value of this parameter from first
principles.

Secondly, the number \(k\), called the ``level'' in Chern--Simons theory,
is given by \[k = \frac{A}{4\pi\gamma}\] Okay, that's all for my quick
description of what we get when we quantize the phase space for an
isolated black hole. I didn't explain how the quantization procedure
actually \emph{works} --- it involves all sorts of fun stuff about theta
functions and so on. I just told the final result.

Now for the entropy calculation. Here we ask the following question:
``given a black hole whose area is within \(\varepsilon\) of \(A\), what
is the logarithm of the number of microstates compatible with this
area?'' This should be the entropy of the black hole --- and it won't
depend much on the number \(\varepsilon\), so long as its on the Planck
scale.

To calculate the entropy, first we work out all the ways to label
punctures by spins \(j\) so that the total area comes within
\(\varepsilon\) of \(A\). For any way to do this, we then count the
allowed ways to pick numbers \(m\) describing the intrinsic curvature of
the black hole surface. Then we sum these up and take the logarithm.

What's the answer? Well, I'll do the calculation for you now in a really
sloppy way, just to sketch how it goes. To get as many ways to pick the
numbers \(m\) as possible, we should concentrate on states where most of
the spins \(j\) labelling punctures equal \(1/2\). If \emph{all} these
spins equal \(1/2\), each puncture contributes an area
\[8 \pi \gamma \sqrt{j(j+1)} = 4 \pi \gamma \sqrt{3}\] to the horizon
2-sphere. Since the total area is close to \(A\), this means that there
are about \(A/(4 \pi \gamma \sqrt{3})\) punctures. Then for each
puncture we can pick \(m = -1/2\) or \(m = 1/2\). This gives
\[N = 2^{A/4\ \pi\ \gamma\ \sqrt{3}}\] ways to choose the \(m\) values.
If this were \emph{exactly} right, the entropy of the black hole would
be \[S = \ln N = \left(\frac{\ln 2}{4 \pi \gamma \sqrt{3}}\right) A\]
Believe it or not, this crude estimate asymptotically approaches the
correct answer as \(A\) approaches infinity. In other words, when the
black hole is in its maximum-entropy state, the vast majority of the
spin network edges poking through the horizon are spin-\(1/2\) edges.

So, what have we seen? Well, we've seen that the black hole entropy is
(asymptotically!) proportional to the area, just like Bekenstein and
Hawking said. That's good. But we don't get the Bekenstein--Hawking
formula \[S = A/4\] because there is an undetermined parameter in our
formula --- the Barbero-Immirzi parameter. That's bad. However, our
answer will match the Bekenstein--Hawking formula if we take
\[\gamma = \frac{\ln 2}{\pi\sqrt{3}}\] If we do this, we no longer have
that annoying undetermined constant floating around in loop quantum
gravity. In fact, we can say that we've determined the ``quantum of
area'' --- the smallest possible unit of area. That's good. And then
it's almost true that in our model, each little patch of the black hole
horizon with area \(4\ln(2)\) contains a single bit of information ---
since a spin-\(1/2\) puncture has area \(4\ln(2)\), and the angle
deficit at a puncture labelled with spin \(1/2\) can take only 2 values,
corresponding to \(m = -1/2\) and \(m = 1/2\). Of course, there are also
punctures labelled by higher values of \(j\), but the \(j = 1/2\)
punctures dominate the count of the microstates.

Of course, one might object to this procedure on the following grounds:
``You've been ignoring matter thus far. What if you include, say,
electromagnetic fields in the game? This will change the calculation,
and now you'll probably need a different value of \(\gamma\) to match
the Bekenstein--Hawking result!''

However, this is not true: we can redo the calculation including
electromagnetism, and the same \(\gamma\) works. That's sort of nice.

There are a lot of interesting comparisons between our way of computing
black hole entropy and the ways its done in string theory, and a lot of
other things to say, too but for that, you'll have to read the
paper\ldots{} I'm worn out now!

\begin{center}\rule{0.5\linewidth}{0.5pt}\end{center}

\textbf{Addenda:} I thank Herman Rubin and Lieven Marchand for some
corrections of errors I made while describing the Riemann hypothesis and
P = NP conjecture. I also thank J. Maurice Rojas for pointing out that
Steve Smale was an individual who \emph{did} have the guts to pose a
list of math problems for the 21st century, back in 1998. This appears
in:

\begin{enumerate}
\def\labelenumi{\arabic{enumi})}
\setcounter{enumi}{5}
\tightlist
\item
  Stephen Smale, ``Mathematical problems for the next century'',
  \emph{Mathematical Intelligencer} \textbf{20} (1998), 7--15.  
  
    Wikipedia, ``Smale's problems'', available at
   \url{https://en.wikipedia.org/wiki/Smale's_problems}.
\end{enumerate}
\noindent
I believe this also appears in the book edited by Arnold mentioned at
the beginning of \protect\hyperlink{week147}{``Week 147''}.

\begin{center}\rule{0.5\linewidth}{0.5pt}\end{center}

\begin{quote}
\emph{\ldots{} for beginners engaging in research, a most difficult
feature to grasp is that of quality --- that is, the depth of a problem.
Sometimes authors work courageously and at length to arrive at results
which they believe to be significant and which experts believe to be
shallow. This can be explained by the analogy of playing chess. A master
player can dispose of a beginner with ease no matter how hard the latter
tries. The reason is that, even though the beginner may have planned a
good number of moves ahead, by playing often the master has met many
similar and deeper problems; he has read standard works on various
aspects of the game so that he can recall many deeply analyzed
positions. This is the same in mathematical research. We have to play
often with the masters (that is, try to improve on the results of famous
mathematicians); we must learn the standard works of the game (that is,
the ``well-known'' results). If we continue like this our progress
becomes inevitable.}

--- Hua Loo-Keng, \emph{Introduction to Number Theory}
\end{quote}

\hypertarget{week149}{%
\section{June 12, 2000}\label{week149}}

Elliptic cohomology sits at the intersection of several well-travelled
mathematical roads. It boasts fascinating connections with homotopy
theory, string theory, elliptic curves, modular forms, and the
mysterious ubiquity of the number 24. This makes it very fascinating,
but also a bit intimidating to anyone who is not already an expert on
all these subjects.

Is \emph{anyone} actually an expert on all these subjects? Perhaps
Graeme Segal is! After all, he became famous for his work on homotopy
theory, he \emph{invented} the axioms of conformal field theory ---
borrowing lots of ideas from string theory, of course --- and I'm sure
he mastered the theory of elliptic curves one weekend when he was a kid.
So to learn about elliptic cohomology, one should really start here:

\begin{enumerate}
\def\labelenumi{\arabic{enumi})}
\tightlist
\item
  Graeme Segal, ``Elliptic cohomology'', \emph{Asterisque}
  \textbf{161--162} (1988), 187--201.
\end{enumerate}

Another good reference is this proceedings of a conference held at
Princeton in 1986:

\begin{enumerate}
\def\labelenumi{\arabic{enumi})}
\setcounter{enumi}{1}
\tightlist
\item
  Peter S. Landweber, editor, \emph{Elliptic Curves and Modular Forms in
  Algebraic Topology}, Springer Lecture Notes in Mathematics
  \textbf{1326}, Springer, Berlin, 1988.
\end{enumerate}
\noindent
This book is also helpful:

\begin{enumerate}
\def\labelenumi{\arabic{enumi})}
\setcounter{enumi}{2}
\tightlist
\item
  Charles B. Thomas, \emph{Elliptic Cohomology}, Kluwer, Dordrecht,
  1999.
\end{enumerate}
\noindent
though I'm afraid it's a bit long on details and short on the big
picture and physics intuition. For that, you might have to try this:

\begin{enumerate}
\def\labelenumi{\arabic{enumi})}
\setcounter{enumi}{3}
\tightlist
\item
  Edward Witten, ``Elliptic genera and quantum field theory'',
  \emph{Comm.\ Math.\ Phys.} \textbf{109} (1987), 525--536.
\end{enumerate}
\noindent
Also try this book, if you can get ahold of it:

\begin{enumerate}
\def\labelenumi{\arabic{enumi})}
\setcounter{enumi}{4}
\tightlist
\item
  Friedrich Hirzebruch, Thomas Berger and Rainer Jung, \emph{Manifolds
  and Modular Forms}, translated by Peter S. Landweber, Vieweg,
  Braunschweig (a publication of the Max Planck Institute for
  Mathematics in Bonn), 1992.
\end{enumerate}

Now to have a snowball's chance in hell of understanding elliptic
cohomology, you need to understand complex oriented cohomology theories.
So I have to start by telling you what \emph{those} are. This will be
sort of a crash course in algebraic topology. By the time I'm done with
that, I'll probably be too worn out to talk about elliptic cohomology
--- but at least I'll have laid the groundwork.

So: what's a ``generalized cohomology theory''?

This is a gadget that eats a topological space \(X\) and spits out a
sequence of abelian groups \(h^n(X)\). To be a generalized cohomology
theory, this gadget must satisfy a bunch of axioms called the
Eilenberg--Steenrod axioms. The most basic example is so-called ordinary
cohomology, so when you're first learning this stuff the main motivation
for the Eilenberg--Steenrod axioms is that they're all satisfied by
ordinary cohomology. But there are lots of other examples: various
flavors of K-theory, cobordism theory, and so on. Eventually, you learn
that underlying any generalized cohomology theory there is a list of
spaces \(E(n)\) such that \[h^n(X) = [X, E(n)]\] where the right-hand
side is the set of homotopy classes of maps from \(X\) to \(E(n)\). We
say this list of spaces \(E(n)\) ``represents'' the generalized
cohomology theory. Moreover, these spaces fit together to form a
``spectrum'', meaning that the space of based loops in \(E(n)\) is
\(E(n-1)\). It follows that each space \(E(n)\) is an infinite loop
space: a space of loops in a space of loops in a space of loops
in\ldots{} where you can go on as far as you like.

Conversely, given an infinite loop space \(E(0)\), we can use it to cook
up a spectrum and thus a generalized cohomology theory. So generalized
cohomology theories, spectra and infinite loop spaces are almost the
same thing.

But what's so important about them?

Well, secretly an infinite loop space is nothing but a homotopy
theorist's version of an abelian group. A bit more technically, we could
call it a ``homotopy coherent abelian group''. By this I mean a space
with a continuous binary operation satisfying all the usual laws for an
abelian group \emph{up to homotopy}, where these homotopies satisfy all
the nice laws you can imagine \emph{up to homotopy}, and so on ad
infinitum. In the context of homotopy theory, this is almost as good as
an abelian group. Pretty much anything a normal mathematician can do
with an abelian group, a homotopy theorist can do with an infinite loop
space!

For example, normal mathematicians often like to take an abelian group
and equip it with an extra operation called ``multiplication'' that
makes it into a \emph{ring}. Homotopy theorists like to do the same for
infinite loop spaces. But of course, the homotopy theorists only demand
that the ring axioms hold \emph{up to homotopy}, where the homotopies
satisfy a bunch of nice laws \emph{up to homotopy}, and so on. Usually
they do this in the context of spectra rather than infinite loop spaces
--- a distinction too technical for me to worry about here! --- so they
call this sort of thing a ``ring spectrum''. Similarly, corresponding to
a commutative ring, the homotopy theorists have a notion called an
``\(E_\infty\) ring spectrum''. The word ``\(E_\infty\)'' is just a
funny way of saying that the commutative law holds up to homotopy, with
the homotopies satisfying a bunch of laws up to homotopy, etcetera.

If you start with a ring spectrum, the corresponding cohomology theory
will have products. In other words, the cohomology groups \(h^n(X)\) of
any space \(X\) will fit together to form a graded ring called
\(h^*(X)\) --- the star stands for a little blank where you can stick in
any number ``\(n\)''. And if your ring spectrum is an \(E_\infty\) ring
spectrum, \(h^*(X)\) will be graded-commutative. This is what happens in
most of really famous generalized cohomology theories. For example, the
ordinary cohomology of a space is actually a graded-commutative ring
with a product called the ``cap product'', and similar things are true
for the most popular flavors of K-theory and cobordism theory.

Of course, it's quite a bit of work to make all this stuff precise:
people spent a lot of energy on it back in the 1970's. But it's very
beautiful, so everybody should learn it. For the details, try:

\begin{enumerate}
\def\labelenumi{\arabic{enumi})}
\setcounter{enumi}{5}
\item
  J. Adams, \emph{Infinite Loop Spaces}, Princeton U. Press, Princeton,
  1978.
\item
  J. Adams, \emph{Stable Homotopy and Generalized Homology}, 
  U. Chicago Press, Chicago, 1974.
\item
  J. P. May, \emph{The Geometry of Iterated Loop Spaces}, Lecture Notes
  in Mathematics \textbf{271}, Springer, Berlin, 1972.  Also available at 
 \href{http://www.math.uchicago.edu/~may/BOOKS/geom_iter.pdf}{\texttt{http://www.math.uchicago.edu/\~may/}} \href{http://www.math.uchicago.edu/~may/BOOKS/geom_iter.pdf}{\texttt{geom\_iter.pdf}}
\item
  J. P. May, F. Quinn, N. Ray and J. Tornehave, \emph{\(E_\infty\) Ring
  Spaces and \(E_\infty\) Ring Spectra}, Lecture Notes in Mathematics
  \textbf{577}, Springer, Berlin, 1977.  Also available at \href{http://www.math.uchicago.edu/~may/BOOKS/e_infty.pdf}{\texttt{http:/}} \href{http://www.math.uchicago.edu/~may/BOOKS/e_infty.pdf}{\texttt{/www.math.uchicago.edu/\~may/BOOKS/e\_infty.pdf}} 
\item
  G. Carlsson and R. Milgram, ``Stable homotopy and iterated loop
  spaces'', in \emph{Handbook of Algebraic Topology}, edited by Ioan M. James,
  North-Holland, Amsterdam, 1995.
\end{enumerate}

Now, there's a particularly nice class of generalized cohomology
theories called ``complex oriented cohomology theories''. Elliptic
cohomology is one of these, so to understand elliptic cohomology you
first have to study these guys a bit. Instead of just giving you the
definition, I'll lead up to it rather gradually\ldots.

Let's start with the integers, \(\mathbb{Z}\). These form an abelian
group under addition, so by what I said above they are a pitifully
simple special case of an infinite loop space. So there's some space
with a basepoint called \(K(\mathbb{Z},1)\) such that the space of all
based loops in \(K(\mathbb{Z},1)\) is \(\mathbb{Z}\).

Be careful here: I'm now using the word ``is'' the way homotopy
theorists do! I really mean the space of based loops in
\(K(\mathbb{Z},1)\) is \emph{homotopy equivalent} to \(\mathbb{Z}\). But
since we're doing homotopy theory, that's good enough.

Okay: so there's a space \(K(\mathbb{Z},1)\) such that the space of all
based loops in \(K(\mathbb{Z},1)\) is \(\mathbb{Z}\). Similarly, there's
a space \(K(\mathbb{Z},2)\) such that the space of all based loops in
\(K(\mathbb{Z},2)\) is \(K(\mathbb{Z},1)\). And so on\ldots{} that's
what it means to say that \(\mathbb{Z}\) is an infinite loop space.

These spaces \(K(\mathbb{Z},n)\) are called ``Eilenberg--Mac Lane
spaces'', and they fit together to form a spectrum called the
Eilenberg--Mac Lane spectrum. Since it's built using only the integers,
this is the simplest, nicest spectrum in the world. Thus the generalized
cohomology theory it represents has got to be something simple and nice.
And it is: it's just ordinary cohomology!

But what do the spaces \(K(\mathbb{Z},n)\) actually look like?

Well, for starters, \(K(\mathbb{Z},0)\) is just \(\mathbb{Z}\), by
definition.

\(K(\mathbb{Z},1)\) is just the circle, \(S^1\). You can check that the
space of based loops in \(S^1\) is homotopy equivalent to \(\mathbb{Z}\)
--- the key is that such loops are classified up to homotopy by an
integer called the \emph{winding number}. In quantum physics,
\(K(\mathbb{Z},1)\) usually goes by the name \(\mathrm{U}(1)\) --- the
group of unit complex numbers, or ``phases''.

\(K(\mathbb{Z},2)\) is a bit more complicated: it's infinite-dimensional
complex projective space, \(\mathbb{CP}^\infty\)! I talked a bunch about
projective spaces in \protect\hyperlink{week106}{``Week 106''}. There I
only talked about finite-dimensional ones like \(\mathbb{CP}^n\), but
you can define \(\mathbb{CP}^\infty\) as a ``direct limit'' of these as
\(n\) approaches \(\infty\), using the fact that \(\mathbb{CP}^n\) sits
inside \(\mathbb{CP}^{n+1}\) as a subspace. Alternatively, you can take
your favorite complex Hilbert space \(H\) with countably infinite
dimension and form the space of all \(1\)-dimensional subspaces in
\(H\). This gives a slightly fatter version of \(\mathbb{CP}^\infty\),
but it's homotopy equivalent, and it's a very natural thing to study if
you're a physicist: it's just the space of all ``pure states'' of the
quantum system whose Hilbert space is \(H\).

How about \(K(\mathbb{Z},3)\)? Well, I don't know a nice geometrical
description of this one. And this really pisses me off! There should be
some nice way to think of \(K(\mathbb{Z},3)\) as some sort of
infinite-dimensional manifold. What is it? Does anyone know? Jean-Luc
Brylinski raised this question at the Conference on Higher Category
Theory and Physics in 1997, and it's been bugging me ever since. From
the work of Brylinski which I summarized in
\protect\hyperlink{week25}{``Week 25''}, it's clear that a good answer
should shed light on stuff like quantum theory and string theory.
Basically, the point is that the integers, the group \(\mathrm{U}(1)\),
and infinite-dimensional complex projective space are all really
important in quantum theory. This is perhaps more obvious for the latter
two spaces --- the integers are so basic that it's hard to see what's so
``quantum-mechanical'' about them. However, since each of these spaces
is just the loop space of the next, they're all part of tightly linked
sequence\ldots{} and I want to know what comes next!

But I'm digressing. I really want to focus on \(K(\mathbb{Z},2)\), or in
other words, infinite-dimensional complex projective space. Note that
there's an obvious complex line bundle over this space. Remember, each
point in \(\mathbb{CP}^\infty\) is really a \(1\)-dimensional subspace
in some Hilbert space \(H\). So we can use these \(1\)-dimensional
subspaces as the fibers of a complex line bundle over
\(\mathbb{CP}^\infty\), called the ``canonical bundle''. I'll call this
line bundle \(L\).

The complex line bundle \(L\) is important because it's ``universal'':
all the rest can be obtained from this one! More precisely, suppose we
have any topological space \(X\) and any map
\[f\colon X\to \mathbb{CP}^\infty\] Then we can form a complex line
bundle over \(X\) whose fiber over any point \(x\) is just the fiber of
\(L\) over the point \(f(x)\). This bundle is called the ``pullback'' of
\(L\) by the map \(f\). And the really cool part is that \emph{any}
complex line bundle over \emph{any} space \(X\) is isomorphic to the
pullback of \(L\) by some map! Even better, two such line bundles are
isomorphic if and only if the maps \(f\) defining them are homotopic!
This reduces the study of many questions about complex line bundles to
the study of this guy \(L\).

For example, suppose we want to classify complex line bundles over any
space \(X\). From what I just said, this task is equivalent to the task
of classifying homotopy classes of maps
\[f\colon X\to \mathbb{CP}^\infty.\] But remember,
\(\mathbb{CP}^\infty\) is the Eilenberg-Maclane space
\(K(\mathbb{Z},2)\), and the Eilenberg-Maclane spectrum represents
ordinary cohomology! So
\[[X, \mathbb{CP}^\infty] = [X, K(\mathbb{Z},2)] = H^2(X)\] where
H\^{}2(X) stands for the 2nd ordinary cohomology group of \(X\). So the
following things are really the same:

\begin{itemize}
\tightlist
\item
  isomorphism classes of complex line bundles over \(X\)
\item
  homotopy classes of maps from \(X\) to \(\mathbb{CP}^\infty\)
\item
  elements of the ordinary cohomology group \(H^2(X)\).
\end{itemize}

This is great, because it gives us three different viewpoints to play
with. In particular, \(H^2(X)\) is easy to compute --- anyone who has
taken a basic course on algebraic topology can do it. But the other two
viewpoints are nice and geometrical, so they let us do things with
\(H^2(X)\) that we might not have thought of otherwise.

So now you know this: if you hand me a complex line bundle over \(X\), I
can cook up an element of \(H^2(X)\). People call this the ``first Chern
class'' of the line bundle. If you hand me two complex line bundles, I
can tell if they're isomorphic by seeing if their first Chern classes
are equal. Conversely, if you hand me any element of \(H^2(X)\), I can
cook up a complex line bundle over \(X\) whose first Chern class is that
element.

Of course, I haven't really explained \emph{how} I cook up all these
things. To learn that, you need to study this stuff a bit more.

But let's consider a couple of examples. Suppose \(X\) is the 2-sphere
\(S^2\). Since \[H^2(S^2) = \mathbb{Z}\] this means that first Chern
class of a line bundle over \(S^2\) is secretly just an integer. People
call this the ``first Chern number'' of the line bundle. The first
physicist to get excited about this was Dirac, who bumped into this idea
when thinking about magnetic monopoles and charge quantization. Dirac
didn't know about complex line bundles and Chern classes --- he was just
studying the change of phase of an electrically charged particle as you
move it around in the magnetic field produced by a monopole! But later,
the physicist Yang met the mathematician Chern and translated Dirac's
work into the language of line bundles. See

\begin{enumerate}
\def\labelenumi{\arabic{enumi})}
\setcounter{enumi}{10}
\tightlist
\item
  C. N. Yang, ``Magnetic monopoles, fiber bundles and gauge field'', in
  \emph{Selected Papers, 1945--1980, with Commentary}, W. H. Freeman and
  Company, San Francisco, 1983.
\end{enumerate}
\noindent
for the full story.

Next let's try a curiously self-referential example. It should be fun to
classify complex line bundles on \(\mathbb{CP}^\infty\), since this is
where the universal one lives! So let's take \(X = \mathbb{CP}^\infty\).
Since \(\mathbb{CP}^\infty\) is \(K(\mathbb{Z},2)\), a little abstract
nonsense shows that it's ordinary 2nd cohomology group is
\(\mathbb{Z}\):
\[H^2(\mathbb{CP}^\infty) = [\mathbb{CP}^\infty, \mathbb{CP}^\infty] = \mathbb{Z}.\]
This means that the first Chern class of a complex line bundle over
\(\mathbb{CP}^\infty\) is secretly just an integer. But what's the first
Chern class of the universal complex line bundle, \(L\)? Well, this
bundle is the pullback of itself via the \emph{identity} map
\[1\colon \mathbb{CP}^\infty\to \mathbb{CP}^\infty\] and this map
corresponds to the element \(1\) in
\([\mathbb{CP}^\infty, \mathbb{CP}^\infty] = \mathbb{Z}\). So the first
Chern class of \(L\) is \(1\). See how tautologous this argument is? It
sounds like it's saying something profound, but once you understand it,
it's really just saying \(1 = 1\).

The first Chern class of the universal bundle \(L\) is really important,
so let's call it \(c\). It's important because it's universal: it gives
us a nice way to think of the first Chern class of \emph{any} complex
line bundle. Up to isomorphism, any complex line bundle over any space
\(X\) comes from some map \[f\colon X\to \mathbb{CP}^\infty\] so to
compute the first Chern class of this line bundle, we can just work out
\(f^*(c)\), where \[f^*\colon H^2(\mathbb{CP}^\infty)\to H^2(X)\] is the
map induced by f.~If you don't see why this is true, think about it a
while --- it's just a big fat tautology!

The ideas we've been discussing raise some obvious questions. For
example, \(H^2(X)\) isn't just a set: it's an abelian group. We already
knew this from our basic course in algebraic topology, and now we also
know another explanation: \(\mathbb{CP}^\infty\) is an infinite loop
space, so it's like an abelian group for the purposes of homotopy
theory. In fact, this particular infinite loop space actually \emph{is}
an abelian group. Maps from anything into an abelian group form an
abelian group, which makes \[H^2(X) = [X, \mathbb{CP}^\infty]\] into an
abelian group. But now you're dying to know: what exactly do the product
map
\[m\colon\mathbb{CP}^\infty\times\mathbb{CP}^\infty\to \mathbb{CP}^\infty\]
and the inverse map \[i\colon\mathbb{CP}^\infty\to \mathbb{CP}^\infty\]
look like? And what does all this mean for the set of isomorphism
classes of complex line bundles on \(X\)? It's an abelian group --- but
what are products and inverses like in this abelian group?

Well, I won't answer the first question here: there's a very nice
explicit answer, and you can describe it in terms of particles and
antiparticles running around on the Riemann sphere, but it would be too
much of a digression to talk about it here. To learn more, study the
``Thom--Dold theorem'' and also some stuff about ``configuration spaces''
in topology:

\begin{enumerate}
\def\labelenumi{\arabic{enumi})}
\setcounter{enumi}{11}
\tightlist
\item
  Dusa McDuff, ``Configuration spaces of positive and negative
  particles'', \emph{Topology} \textbf{14} (1975), 91--107.
\end{enumerate}

The second question is much easier: the set of isomorphism classes of
complex line bundles on a space \(X\) becomes an abelian group with
\emph{tensor product} of line bundles as the product. Taking the
\emph{dual} of a line bundle gives the inverse in this group.

Putting these ideas together, we get a nice description of tensoring
line bundles in terms of the product
\[m\colon \mathbb{CP}^\infty\times\mathbb{CP}^\infty\to \mathbb{CP}^\infty\]
which I can explain even without saying what the product looks like.
Suppose I have two line bundles on \(X\) and I want to tensor them. I
might as well assume they are pullbacks of the universal bundle \(L\) by
some maps \[f\colon X\to \mathbb{CP}^\infty,\]
\[g\colon X\to \mathbb{CP}^\infty.\] It follows from what we've seen
that to tensor these bundles, I can just form the map
\[fg: X\to \mathbb{CP}^\infty\] given as the composite
\[X\xrightarrow{(f,g)}\mathbb{CP}^\infty\times\mathbb{CP}^\infty\xrightarrow{m}\mathbb{CP}^\infty\]
and then take the pullback of \(L\) by \(fg\).

In other words: since the canonical line bundle on
\(\mathbb{CP}^\infty\) is universal, \(\mathbb{CP}^\infty\) knows
everything there is to know about complex line bundles. In particular,
it knows everything there is to know about \emph{tensoring} complex line
bundles: the operation of tensoring is encoded in the \emph{product} on
\(\mathbb{CP}^\infty\). Similarly, the operation of taking the
\emph{dual} of a complex line bundle is encoded in the \emph{inverse}
operation \[i\colon \mathbb{CP}^\infty\to \mathbb{CP}^\infty.\] Now, if
you've absorbed everything I just said --- or better yet, if you already
knew it! --- you are ready for the definition of a ``complex oriented
cohomology theory''. Basically, it's a generalized cohomology theory
where all this stuff about line bundles and the first Chern class works
almost like it does in ordinary cohomology.

Suppose we have a generalized cohomology theory; let's see the
conditions under which it's ``complex oriented''.

For starters, it must come from a ring spectrum, so that \(h^*(X)\) is a
graded ring for any space \(X\). Whenever we're in this situation, it's
interesting to take \(X\) to be a single point: we get a ring \(R\)
called the ``cohomology ring of a point''. This has a god-given element
of degree 0, namely the multiplicative unit \(1\). In any generalized
cohomology theory we have an isomorphism
\[h^{n+k}(S^n) = h^k(\mathrm{point})\] and taking the god-given element
\(1\) in \(h^0(\mathrm{point})\), this gives gives a special element in
\(h^n(S^n)\) called the ``orientation'' of \(S^n\). Now note that
\(S^2\) is the same thing as \(\mathbb{CP}^1\), so that it sits inside
\(\mathbb{CP}^\infty\): \[S^2\to \mathbb{CP}^\infty\] This gives a map
\[h^2(\mathbb{CP}^\infty)\to h^2(S^2)\] We say our generalized
cohomology is ``complex oriented'' if there is an element \(c\) in
\(h^2(\mathbb{CP}^\infty)\) that maps to the orientation of \(S^2\)
under the above map, and changes sign under the inverse map
\[i^*\colon h^2(\mathbb{CP}^\infty)\to h^2(\mathbb{CP}^\infty).\] For
example, ordinary cohomology is complex oriented, where we take \(c\) to
be the first Chern class of the universal complex line bundle! This
follows from all the stuff I've said so far.

But lots of other generalized cohomology theories are complex oriented,
too. The most famous ones are complex K-theory and complex cobordism
theory. In fact, complex cobordism theory is the ``universal'' complex
oriented cohomology theory --- it's the most informative of the whole
lot. All the rest are like watered-down versions of this one. Ordinary
cohomology is the most watered-down of all. Complex K-theory is a bit
less watered-down. And elliptic cohomology is still less watered-down!

But what \emph{is} elliptic cohomology?

Well, I might or might not get around to talking about this next Week.
I've learned the hard way not to \emph{promise} to talk about things in
future issues: my mind is too scattered to be able to stick to one
subject in a predictable manner. For example, last week I hinted that
I'd talk about Ross Street's work on \(n\)-categories, but now I've
spent so much time blabbing about this other stuff that I don't think
I'll get around to it. Let me just list the papers he gave me:

\begin{enumerate}
\def\labelenumi{\arabic{enumi})}
\setcounter{enumi}{12}
\item
  Ross H. Street, ``The petit topos of globular sets'', \emph{Jour.\ Pure 
 Appl.\ Alg.\ }\textbf{154} (2000), 299--315.  Also available at 
 \url{https://core.ac.uk/display/82254278}
 
\item
  Ross H. Street and Michael Batanin, ``The universal property of the
  multitude of trees'',  Also available at \url{https://core.ac.uk/display/82645668}

\item
  Michael Batanin, ``Shuffle polytopes, cooperative games and
  \(2\)-dimensional coherence for higher dimensional categories''.  Apparently
  not available.
\end{enumerate}

The folks down in Sydney are making great progress on understanding
\(n\)-categories from the globular point of view, and the importance of
the category of \emph{trees} has become quite clear. You can think of
trees as generalized natural numbers, and then Batanin's operads are a
natural generalization of the usual operads, which have operations
taking a natural number's worth of arguments. The trees describe ways of
glomming a bunch of globes together to get a new globe. I wish I had
time to explain this better! But it takes a while, and it really
requires some pictures.

\begin{center}\rule{0.5\linewidth}{0.5pt}\end{center}

Footnote:

{[}1{]} Almost, but not quite: if I hand you the infinite loop space
\(E(0)\), you can only recover one connected component of the infinite
loop space \(E(1)\), namely the component containing the basepoint. So
there is more information in a spectrum than there is in an infinite
loop space. A spectrum is a sequence of infinite loop spaces where the
based loops in \(E(n)\) form the space \(E(n-1)\); starting from a
single infinite loop space we can cook up a spectrum, but it will be a
spectrum of a special sort, called a ``connective'' spectrum, where the
spaces \(E(n)\) are connected for \(n > 0\).

Given a spectrum we can define the generalized cohomology groups
\(H^n(X)\) even when \(n\) is negative, via:
\[H^n(X) = \lim_{k\to\infty} [\Sigma^k(X), E(n+k)]\] where
\(\Sigma^k(X)\) denotes the \(k\)-fold suspension of \(X\). If the
spectrum is connective, these groups will vanish when \(n\) is negative.
A good example of a connective spectrum is the spectrum for ordinary
cohomology (the Eilenberg--Mac Lane spectrum). A good example of a
nonconnective spectrum is the spectrum for real or complex K-theory.

\begin{center}\rule{0.5\linewidth}{0.5pt}\end{center}

\begin{quote}
\emph{This therefore, is mathematics: she reminds you of the invisible
forms of the soul; she gives light to her own discoveries; she awakens
the mind and purifies the intellect; she brings light to our intrinsic
ideas; she abolishes oblivion and ignorance which are ours by birth.}

--- Proclus, \emph{Commentary on the First Book of Euclid's \emph{Elements}}
\end{quote}

\hypertarget{week150}{%
\section{June 18, 2000}\label{week150}}

First I'd like to say some stuff about Lagrange points. Then I'll
continue talking about complex oriented cohomology theories.

Euler and Lagrange won the Paris Academy Prize in 1722 for their work on
the orbit of the moon. The essay that Lagrange submitted for this prize,
``Essai sur le probleme des trois corps'', contained some interesting
results on the 3-body problem. Among other things, he studied the case
where a lighter body was revolving about a heavier one in a circular
orbit, and found all the places where a third much lighter body could
sit ``motionless'' with respect to the other two. There are 5 such
places, and they're now called the ``Lagrange points'' L1 -- L5.

For example: imagine the moon going around the earth in a circular
orbit. Then there are 5 Lagrange points where we can put a satellite. 3
of these are unstable equilibria. They lie on the line through the earth
and moon. L1 is between the earth and moon, L2 is in the same direction
as the moon but a bit further out, and L3 is opposite the moon. The
other 2 Lagrange points are stable equilibria. L4 lies in the moon's
orbit 60 degrees ahead of the moon, while L5 lies 60 degrees behind the
moon.

Here's a nice picture of the Lagrange points:

\begin{enumerate}
\def\labelenumi{\arabic{enumi})}
\tightlist
\item
   NASA, ``Lagrange points'',
  \href{https://map.gsfc.nasa.gov/mission/observatory_l2.html}{\texttt{https://map.gsfc.nasa.gov/mission/observatory\_l2.}}
\href{https://map.gsfc.nasa.gov/mission/observatory_l2.html}{\texttt{html}}
\end{enumerate}
\noindent
In general, the points L4 and L5 will be stable equilibria as long as
the heavy body (e.g.~the earth in the above example) is more than 25
times as massive as the intermediate-sized one (e.g.~the moon). And in
case you're wondering, this magic number is just an approximation to the
exact figure, which is
\[\frac{25}{2}\left(1+\sqrt{1-\frac{4}{625}}\right) = 24.95993579437711227887\ldots\]
Here's a proof that the Lagrange points work as advertised, including a
derivation of the above number:

\begin{enumerate}
\def\labelenumi{\arabic{enumi})}
\setcounter{enumi}{1}
\tightlist
\item
  Neil J. Cornish, ``The Lagrange points'', available at
  \href{https://map.gsfc.nasa.gov/ContentMedia/lagrange.pdf}{\texttt{{https://map.gsfc.nasa.gov/}}}   \href{https://map.gsfc.nasa.gov/ContentMedia/lagrange.pdf}{\texttt{{ContentMedia/lagrange.pdf}}}
\end{enumerate}

Now, Lagrange did his calculation as a mathematical exercise and didn't
believe it was relevant to the actual solar system. But he was wrong
about that! The stable Lagrange points L4 and L5 are quite interesting,
because stuff tends to accumulate there.

For example, people have found over six hundred asteroids called
``Trojans'' at the stable Lagrange points of Jupiter's orbit around the
sun. The first to be discovered was 588 Achilles, back in 1906 --- the
number here meaning that it was the 588th asteroid found. In general,
the Trojans at L4 are named after Greek soldiers in the Trojan war,
while those at L5 are named after actual Trojans --- soldiers from the
city of Troy! You can see the Trojans quite clearly in this picture of
the asteroid belt:

\begin{enumerate}
\def\labelenumi{\arabic{enumi})}
\setcounter{enumi}{2}
\tightlist
\item
  Map of inner solar system, available at \href{https://commons.wikimedia.org/wiki/File:InnerSolarSystem-en.png}{\texttt{https://commons.wikimedia.org/wiki/}} \href{https://commons.wikimedia.org/wiki/File:InnerSolarSystem-en.png}{\texttt{File:InnerSolarSystem-en.png}}
\end{enumerate}

The asteroid 5261 Eureka is a ``Martian Trojan'', occupying the L5 point
of Mars' orbit around the sun. A second Martian Trojan was discovered at
L5 in 1998, but it doesn't have a name yet. For more information, try:

\begin{enumerate}
\def\labelenumi{\arabic{enumi})}
\setcounter{enumi}{3}
\tightlist
\item
Wikipeda, Trojan (celestial body), \href{https://en.wikipedia.org/wiki/Trojan_(celestial_body)}{\texttt{https://en.wikipedia.org/wiki/Trojan\_}} \hfill
\break \href{https://en.wikipedia.org/wiki/Trojan_(celestial_body)}{\texttt{(celestial\_body)}}
\end{enumerate}

There may also be a few small asteroids at the Lagrange points of Venus
and Earth's orbits around the sun. Does anyone know more about this? My
brief search revealed only some information about the curious asteroid
3753 Cruithne, which is a companion asteroid of the Earth. But 3753
Cruithne is not at a Lagrange point! Instead, it moves in a very
complicated spiralling horse-shoe shaped orbit relative to the earth.
For a beautiful explanation with pictures by the discoverers of this
phenomenon, see:

\begin{enumerate}
\def\labelenumi{\arabic{enumi})}
\setcounter{enumi}{4}
\tightlist
\item
  Paul Wiegert, Kimmo Innanen and Seppo Mikkola, Near-earth asteroid
  3753 Cruithne --- Earth's curious companion,
  \url{http://www.astro.queensu.ca/~wiegert/}
\end{enumerate}

For over 150 years, astronomers have been searching for other satellites
of the Earth besides the big one I see out my window right now. There
have been a lot of false alarms, with people even giving names to the
satellites they thought they discovered, like ``Kleinchen'' and
``Lilith''. For the fascinating story of these ``second moons'' and
other hypothetical planets, see:

\begin{enumerate}
\def\labelenumi{\arabic{enumi})}
\setcounter{enumi}{5}
\tightlist
\item
  Paul Schlyter, ``Hypothetical planets'',
 \href{https://www.astro.auth.gr/ANTIKATOPTRISMOI/nineplanets/nineplanets/hypo.html}{\texttt{https://www.astro.auth.gr/}} \hfill \break \href{https://www.astro.auth.gr/ANTIKATOPTRISMOI/nineplanets/nineplanets/hypo.html}{\texttt{ANTIKATOPTRISMOI/nineplanets/nineplanets/hypo.html}}
\end{enumerate}

Unfortunately, none of these second moons were real. But in the 1960s,
people discovered dust clouds at the stable Lagrange points of the
Moon's orbit around the Earth! They are about 4 times as big as the
Moon, but they are not very substantial.

What else lurks at Lagrange points?

Well, Saturn has a moon called Dione, and 60 degrees ahead of Dione,
right at the Lagrange point L4, there is a tiny moon called Helene.
Here's a picture of Helene:

\begin{enumerate}
\def\labelenumi{\arabic{enumi})}
\setcounter{enumi}{6}
\tightlist
\item
  Astronomy picture of the day: ``Dione's Lagrange moon Helene'',
  \href{http://antwrp.gsfc.nasa.gov/apod/ap951010.html}{texttt{http://antwrp.gsfc.nasa.gov/apod/ap951010.html}}   \href{http://antwrp.gsfc.nasa.gov/apod/ap951010.html}{texttt{http://antwrp.gsfc.nasa.gov/apod/ap951010.html}} 
\end{enumerate}

Isn't she cute? There's also a small moon called Telesto at the L4 point
of Saturn's moon Tethys, and one called Calypso at the L5 point of
Tethys:

\begin{enumerate}
\def\labelenumi{\arabic{enumi})}
\setcounter{enumi}{7}
\tightlist
\item
  Bill Arnett, Introduction to the nine planets: Tethys, Telesto and
  Calypso,
  \href{http://astrolink.mclink.it/nineplanets/tethys.html}{\texttt{http://}}   \href{http://astrolink.mclink.it/nineplanets/tethys.html}{\texttt{/astrolink.mclink.it/nineplanets/tethys.html}}
\end{enumerate}
\noindent
Does anyone know other natural occupants of Lagrange points?

For a long time crackpots and science fiction writers have talked about
a ``Counter-Earth'', complete with its own civilization, located at the
L3 point of the Earth's orbit around the Sun --- exactly where we can
never see it! I think satellite explorations have ruled out this
possibility by now, but since L3 is an unstable Lagrange point, it was
never very likely to begin with.

On the other hand, fans of space exploration have long dreamt of setting
up a colony at the stable Lagrange points of the Moon's orbit around the
earth --- right in those dust clouds, I guess. But now people are
putting artificial satellites at the Lagrange points of the Earth's
orbit around the Sun! For example, L1 is the home of ``SOHO'': the Solar
and Heliospheric Observatory. Sitting between the Earth and Sun gives
SOHO a nice clear view of sunspots, solar flares, and the solar wind.
It's not stable, but it can exert a bit of thrust now and then to stay
put. For more information and some pretty pictures, try:

\begin{enumerate}
\def\labelenumi{\arabic{enumi})}
\setcounter{enumi}{8}
\tightlist
\item
  SOHO website, \url{https://soho.nascom.nasa.gov/}
\end{enumerate}

In April 2001, NASA plans to put a satellite called ``MAP'' at the
Lagrange point L2. MAP is the Microwave Anisotropy Probe, which will
study anisotropies in the cosmic microwave background. These are
starting to be a really interesting window into the early history of the
universe. For more, see:

\begin{enumerate}
\def\labelenumi{\arabic{enumi})}
\setcounter{enumi}{9}
\tightlist
\item
  MAP website, \url{http://map.gsfc.nasa.gov/}
\end{enumerate}

Finally, for quantum analogues of the Lagrange points, see:

\begin{enumerate}
\def\labelenumi{\arabic{enumi})}
\setcounter{enumi}{10}
\tightlist
\item
  T. Uzer, Ernestine A. Lee, David Farrelly, and Andrea F. Brunello,
  Synthesis of a classical atom: wavepacket analogues of the Trojan
  asteroids, \emph{Contemp. Phys.} \textbf{41} (2000), 1--14. 
\end{enumerate}
\noindent
Okay, enough about Lagrange points! Now I want to talk a bit more about
complex oriented cohomology theories. Last time I left off at the
definition. So let me start with a little review, and then plunge ahead.

A generalized cohomology theory assigns to each space \(X\) a bunch of
groups \(h^n(X)\), one for each integer \(n\). We impose some axioms
that make them work very much like ordinary cohomology. However, when
\(X\) is a point, we no longer require that \(h^n(X)\) is trivial for
nonzero \(n\). It may not seem like much, but it turns out to make a big
difference! There are all sorts of very interesting examples.

Like what?

Well, there's K-theory, which lets us study a space by looking at vector
bundles over that space and its iterated suspensions. We get various
flavors of K-theory from various kinds of vector bundle: real K-theory,
complex K-theory, quaternionic K-theory, and so on. There's even a sort
of K-theory invented by Atiyah that's based on Clifford algebras, called
``KR theory''. For more about all these, try:

\begin{enumerate}
\def\labelenumi{\arabic{enumi})}
\setcounter{enumi}{11}
\item
  Dale Husemoller, \emph{Fibre Bundles}, Springer, Berlin,
  1975.
\item
  H. Blaine Lawson, Jr.~and Marie-Louise Michelsohn, \emph{Spin
  Geometry}, Princeton U.\ Press, Princeton, 1989.
\end{enumerate}

Then there's cobordism theory, which lets us study a space by looking at
manifolds mapped into that space and its iterated suspensions. To be
honest, this is actually how \emph{bordism} theory works --- this being
a generalized \emph{homology} theory. But every generalized homology
theory goes hand-in-hand with a generalized cohomology theory, and if
you understand one you understand the other, at least in
principle\ldots. Anyway, there are various flavors of cobordism theory
corresponding to various kinds of extra structure you can put on a
manifold: piecewise-linear cobordism theory, smooth cobordism theory,
oriented cobordism theory, spin cobordism theory, complex cobordism
theory, symplectic cobordism theory, stable homotopy theory, and so on.
For more about these, try:

\begin{enumerate}
\def\labelenumi{\arabic{enumi})}
\setcounter{enumi}{13}
\tightlist
\item
  Robert E. Stong, \emph{Notes on Cobordism Theory}, Princeton
  U.\ Press, Princeton, 1968.
\end{enumerate}

Finally, there are generalized cohomology theories inspired more by
algebra than by geometry, which only hardcore homotopy theorists seem to
understand, like Morava K-theory and Brown-Peterson theory.

To round off this little tour, I guess I should add that there are lots
of maps going between different generalized cohomology theories! As I
explained in \protect\hyperlink{week149}{``Week 149''}, each generalized
cohomology \(h^n\) corresponds to a ``spectrum'': a list of spaces, each
being the loop space of the next. Spectra form a category, and given a
map between spectra we get a map between their generalized cohomology
theories. So we shouldn't study these one of at a time: it's better to
play around with lots at once! For some important examples of the stuff
you can do this way, try:

\begin{enumerate}
\def\labelenumi{\arabic{enumi})}
\setcounter{enumi}{14}
\item
  P. E. Conner and E. E. Floyd, \emph{The Relation of Cobordism to
  K-theories}, Lecture Notes in Mathematics \textbf{28},
  Springer, Berlin, 1966.
\item
  Douglas C. Ravenel, \emph{Complex Cobordism and Stable Homotopy Groups
  of Spheres}, Academic Press, New York, 1986.
\end{enumerate}
\noindent
and also the books on generalized cohomology listed in
\protect\hyperlink{week149}{``Week 149''}.

Anyway, to bring order to this zoo, it's nice to pick out some of the
special features of \emph{ordinary} cohomology theory and study the
generalized cohomology theories that share these features.

For example, in ordinary cohomology, the cohomology groups of a space
fit together to form a graded ring. If a generalized cohomology theory
is like this, we call it ``multiplicative''. Whenever this is the case,
we get a graded ring \(R\) called the ``coefficient ring'' of our
theory, which is simply the cohomology ring of the one-point space. For
ordinary cohomology theory the coefficient ring is just \(\mathbb{Z}\),
but for other theories it can be very interesting and complicated. By
easy abstract nonsense, the cohomology ring of any space is an algebra
over the coefficient ring.

Another nice feature of ordinary cohomology is the first Chern class.
Whenever you have a complex line bundle over a space \(X\), you get an
invariant called its ``first Chern class'' which lives in \(H^2(X)\),
and this invariant is sufficiently powerful to completely classify such
line bundles. Last week I described the \emph{universal} line bundle
over infinite-dimensional complex projective space:
\[L\to \mathbb{CP}^\infty\] and showed how the first Chern class of
\emph{any} line bundle comes from the first Chern class of this one,
which I called \(c\).

If a generalized cohomology theory is multiplicative and there's an
element \(c\) of \(h^2(\mathbb{CP}^\infty)\) that acts like the first
Chern class of the universal line bundle, we call the theory ``complex
oriented''. Of course, to make this precise we need to isolate the key
features of the first Chern class and abstract them. I did this in
\protect\hyperlink{week149}{``Week 149''}, so I won't do it again here.
Instead, I'll just say a bit about what we can \emph{do} with a complex
oriented cohomology theory.

For starters, we can use the element \(c\) to get an invariant of
complex line bundles --- a kind of generalized version of the first
Chern class. To do this, just remember from
\protect\hyperlink{week149}{``Week 149''} that \(\mathbb{CP}^\infty\) is
the classifying space for complex line bundles. In other words,
\emph{any} line bundle over \emph{any} space \(X\) is isomorphic to a
pullback of the universal line bundle by some map
\[f\colon X\to\mathbb{CP}^\infty.\] Thus, given a line bundle we can
find such a map \(f\) and use it to pull back the element \(c\) to get
an element of \(h^2(X)\). This is exactly like the usual first Chern
class of our line bundle, except now we're using a generalized
cohomology theory instead of ordinary cohomology.

Can we get any \emph{other} invariants of complex line bundles from
\emph{other} elements of the cohomology of \(\mathbb{CP}^\infty\)? Not
really: in any complex oriented cohomology theory, the cohomology ring
of \(\mathbb{CP}^\infty\) is just the algebra of formal power series in
the element \(c\): \[h^*(\mathbb{CP}^\infty) = R[[c]]\] where \(R\) is
the coefficient ring.

However, we can get other invariants of complex \emph{vector} bundles,
which are analogous to the higher Chern classes. In fact, we can mimic a
huge amount of the usual theory of characteristic classes in the context
of a complex oriented cohomology theory. I'd love to talk about this,
but it would be a digression from my main goal, which is to make
elliptic cohomology at least halfway comprehensible to the amateur.

So: last time I mentioned that \(\mathbb{CP}^\infty\) is an abelian
topological group, with a multiplication map
\[m\colon \mathbb{CP}^\infty \times \mathbb{CP}^\infty\to \mathbb{CP}^\infty\]
and an inverse map \[i\colon \mathbb{CP}^\infty\to \mathbb{CP}^\infty.\]
And I explained how these represent the operations of \emph{tensoring}
two line bundles and taking the \emph{dual} of a line bundle,
respectively. Now let's see what we can do with these maps when we have
a complex oriented cohomology theory. First of all, since cohomology is
contravariant we get homomorphisms
\[m^*\colon h^*(\mathbb{CP}^\infty)\to h^*(\mathbb{CP}^\infty \times \mathbb{CP}^\infty)\]
and \[i^*\colon h^*(\mathbb{CP}^\infty)\to h^*(\mathbb{CP}^\infty)\] But
as I already said, we have \[h^*(\mathbb{CP}^\infty) = R[[c]]\] and
similarly we have
\[h^*(\mathbb{CP}^\infty \times \mathbb{CP}^\infty) = R[[c]] \otimes R[[c]]\]
where the product ``\(\otimes\)'' on the right side is a slightly
fattened-up version of the usual tensor product of algebras over \(R\).
So we really have homomorphisms
\[m^*\colon R[[c]]\to R[[c]] \otimes R[[c]]\]
\[i^*\colon R[[c]]\to R[[c]]\] which satisfy all the usual axioms for
the product and inverse in an abelian group --- but turned around
backwards.

Folks who like Hopf algebras will immediately note that \(R[[c]]\) is
like a Hopf algebra. A nice way to form a Hopf algebra is to take the
algebra of functions on a group and use the product and inverse in the
group to give this algebra extra operations called the ``coproduct'' and
``antipode''. We're doing the same thing here, except that we're using
formal power series in one variable \(c\) instead of functions of one
variable. So folks call \(R[[c]]\) a ``formal group law''.

In short: complex oriented cohomology theories give formal group laws!

Lest this seem overly abstract and unmotivated, remember that it's just
a way of talking about what happens to the generalized ``first Chern
class'' when we tensor line bundles. In a vague but useful way, we can
visualize guys in \(R[[c]]\) as formal power series on the line, where
the line has been equipped with some abelian group structure, at least
right near the origin. This group structure is what yields the coproduct
\[m^*\colon R[[c]]\to R[[c]] \otimes R[[c]]\] and antipode
\[i^*\colon R[[c]]\to R[[c]]\] But the real point is that this group
structure tells us how to compute the generalized ``first Chern class''
of a tensor product of line bundles starting from both of their
generalized first Chern classes.

Some examples may help. In ordinary cohomology theory, when we tensor
two line bundles, we just \emph{add} their first Chern classes. So in
this case, we've got the line made into a group using addition, and
\(R[[c]]\) becomes a formal group law called the ``additive formal group
law''.

Another famous example is complex K-theory. In this theory, when we
tensor two line bundles, we basically just \emph{multiply} their
generalized first Chern classes. That's because cohomology classes in
K-theory are just equivalence classes of vector bundles, and multiplying
them just \emph{means} tensoring them. So in this case, \(R[[c]]\)
becomes the ``multiplicative'' formal group law. Of course, some
fiddling around is required, because we don't usually think of the
multiplicative unit \(1\) as the origin of the line\ldots{} but we can
if we want.

What are some other examples? Well, complex cobordism theory is one.
This corresponds to the ``universal'' formal group law: a formal group
law that has a unique homomorphism to any other one! In this case \(R\)
is quite big: it's called the ``Lazard ring''. And this universal aspect
of complex cobordism theory makes it the king of all complex oriented
generalized cohomology theories. I really should spend about 10 pages
explaining it to you in detail, but I won't\ldots.

\ldots{} because I want to finally say a word about elliptic cohomology!

I've talked a lot about elliptic curves in
\protect\hyperlink{week13}{``Week 13''},
\protect\hyperlink{week125}{``Week 125''} and
\protect\hyperlink{week126}{``Week 126''}, so I get to assume you know
about them: an elliptic curve is a \(1\)-dimensional compact complex
manifold which is also an abelian group. As such, any elliptic curve
gives a formal group law. And thus we can wonder if this formal group
law comes from a complex oriented cohomology theory\ldots{} and it does!
And this is elliptic cohomology!

Now this is just the beginning of a long story: there's much more to
say, and I don't have the energy to say it here\ldots{} but I'll just
tantalize you with some of the high points. Since elliptic curves can be
thought of as the worldsheets of strings, there are a bunch of
interesting relationships between string theory and elliptic cohomology.
In addition to the references I've gave you in
\protect\hyperlink{week149}{``Week 149''}, see for example:

\begin{enumerate}
\def\labelenumi{\arabic{enumi})}
\setcounter{enumi}{16}
\tightlist
\item
  Hirotaka Tamanoi, \emph{Elliptic Genera and Vertex Operator
  Super-Algebras}, Springer Lecture Notes in Mathematics \textbf{1704},
  Springer, Berlin, 1999.
\end{enumerate}
\noindent
We also get lots of nice relationships to the theory of modular forms.
But the one thing we apparently don't have \emph{yet} is a very neat
\emph{geometrical} description of elliptic cohomology! There's got to be
one, since so much nice geometry is involved\ldots{} but what is it?

(For an answer to this question, see \protect\hyperlink{week197}{``Week
197''}.)

\begin{center}\rule{0.5\linewidth}{0.5pt}\end{center}

\textbf{Addendum:} For the proof that the complex cobordism theory
corresponds to the universal formal group law, read Quillen's original
paper on the subject:

\begin{enumerate}
\def\labelenumi{\arabic{enumi})}
\setcounter{enumi}{17}
\tightlist
\item
  Daniel Quillen, ``On the formal group laws of unoriented and complex
  cobordism theory'', \emph{Bull. Amer. Math. Soc.} \textbf{75} (1969),
  1293--1298. Also available as
  \href{http://projecteuclid.org/euclid.bams/1183530915}{\texttt{http://}}
  \href{http://projecteuclid.org/euclid.bams/1183530915}{\texttt{projecteuclid.org/euclid.bams/1183530915}}
\end{enumerate}

\begin{center}\rule{0.5\linewidth}{0.5pt}\end{center}

\begin{quote}
\emph{In any field, find the strangest thing and then explore it.}

--- John Wheeler
\end{quote}

\end{document}